\documentclass[12pt]{amsart}
\usepackage[utf8]{inputenc}
\usepackage[english]{babel}
\usepackage{subcaption}
\usepackage{enumerate}
\usepackage{amsfonts}
\usepackage{comment}
\usepackage{amsmath,amssymb}
\usepackage[makeroom]{cancel}
\usepackage{stmaryrd}
\usepackage{tikz-cd} 
\usepackage{enumerate}
\usepackage{empheq}
\usepackage{enumitem}
\usepackage{graphicx}
\usepackage{amsthm}
\usepackage{wrapfig}
\usepackage{hyperref}
\usepackage{chngcntr}
\counterwithin{figure}{section}
\newtheorem{theorem}{Theorem}
\newtheorem*{theorem*}{Theorem}
\newtheorem{prop}[theorem]{Proposition}
\newtheorem{conj}[theorem]{Conjecture}
\newtheorem{coro}{Corollary}[theorem]
\newtheorem{lemm}[theorem]{Lemma}
\theoremstyle{definition}
\newtheorem{defi/}[theorem]{Definition}
\newenvironment{defi}
  {%
   \pushQED{\qed}\begin{defi/}}
  {\popQED\end{defi/}}
\newtheorem*{defi*}{Definition}
\newtheorem{ques}[theorem]{Question}
\newtheorem{thdef}[theorem]{Theorem-Definition}
\newtheorem{cor-definition}[theorem]{Corollary-Definition}
\newtheorem{remadef}[theorem]{Remark-Definition}
\newtheorem{rema}[theorem]{Remark}
\newtheorem{conv}[theorem]{Convention}
\newtheorem*{conv*}{Convention}

\numberwithin{theorem}{section}

\newcommand{\inte}[1]{\overset{\circ}{#1}}
\newcommand{\clos}[1]{\overline{#1}}

\newcommand*\widefbox[1]{\fbox{\hspace{2em}#1\hspace{2em}}}

\newcommand{\quotient}[2]{{\raisebox{.2em}{$#1$}\left/\raisebox{-.2em}{$#2$}\right.}}

\numberwithin{equation}{section}
\makeatletter
\def\l@section#1{%
  \@tocline{1}{1em}{1.5em}{2.3em}{}{\bfseries\ignorespaces#1\unskip\hfill }%
}

\def\l@subsection#1{%
  \@tocline{2}{0.3em}{3em}{3.5em}{}{\ignorespaces#1\unskip\hfill }%
}

\def\l@subsubsection#1{%
  \@tocline{3}{0pt}{4em}{5em}{}{\ignorespaces#1\unskip\hfill}%
}
\makeatother

\begin{document}

\title{Markovian families for pseudo-Anosov flows} 
\author[I. Iakovoglou]{Ioannis Iakovoglou}
\address{IMJ-PRG, Sorbonne Université \\ \vspace{0.1cm} Supported by
 the ERC project 818737 Emergence of wild differentiable dynamical systems}
\email{ioannis.iakovoglou@ens-lyon.fr}

\begin{abstract}
Generalizing the classification approach described in  \cite{preprint} for transitive Anosov flows in dimension 3, in this paper we describe a method for classifying (not necessarily transitive) pseudo-Anosov flows on 3-manifolds up to orbital equivalence. 

To every pseudo-Anosov flow $\Phi$ (with no 1-prongs) on $M^3$ is associated a bifoliated plane $\mathcal{P}$ endowed with an action of $\pi_1(M)$. It is known that the previous action characterizes $\Phi$ up to orbital equivalence and admits infinitely many Markovian families (i.e. collections of rectangles in $\mathcal{P}$ generalizing the notion of Markov partition for group actions on the plane). Our goal in this paper consists in showing that : 

\begin{itemize}
\item if $\mathcal{R}$ is a Markovian family of $\Phi$, the number of orbits of rectangles of $\mathcal{R}$ and their pattern of intersection can be encoded by a finite combinatorial object, called a geometric type, which describes completely $\Phi$ up to Dehn-Goodman-Fried surgeries on a specific finite set $\Gamma$ of periodic orbits of $\Phi$
\item our previous choices of surgeries on $\Gamma$ can be read as sequences of rectangles in $\mathcal{R}$ and can be encoded by finite combinatorial objects, called cycles
\item a geometric type with cycles of $\mathcal{R}$  describes the original flow $\Phi$ up to orbital equivalence
\end{itemize} 

Several of the above results will be stated and proven in a slightly more general setting involving strong Markovian actions on the plane. Finally, due to the lack of bibliographic references on pseudo-Anosov flows in dimension 3, in the first part of the paper we provide an introduction to pseudo-Anosov flow theory containing several useful results for our classification approach together with their proofs.  
\end{abstract}
\maketitle
\vspace{-0.5cm}
\vspace{0.05cm}
\normalsize

\section{Introduction}
\subsection{General setting}
Thanks to its relation with the study of 3-manifolds and their foliations (see for instance \cite{Mosher1}, \cite{Mosher2}, \cite{Fe2} and more recently \cite{Fe3}) and parallelly to the recent interest in the project of classification of Anosov flows on $3$-manifolds, the problem of classification of pseudo-Anosov flows in dimension 3 has been nowadays progressively gaining in interest and importance.

To this day different approaches have been used in order to address the question of classification of pseudo-Anosov flows in dimension 3: 
\begin{enumerate}
    \item By examining the action of the fundamental group on the bifoliated plane of the flow and by using the topology of the ambient manifold: 
        \begin{itemize}
            
            \item the authors of \cite{BarbotFenley} show among others that the only pseudo-Anosov flow on a $3$-manifold with virtually solvable fundamental group is a suspension Anosov flow and that the only pseudo-Anosov flow on a Seifert manifold is up to a finite cover a geodesic flow on the unit tangent bundle of a hyperbolic surface
            \item the authors of \cite{flowsingraphmanifolds} classify the totally periodic pseudo-Anosov flows in graph manifolds
        \end{itemize}
    \item By examining the action of the fundamental group on the circle at infinity $\mathbb{S}_{\infty}$ of the bifoliated plane of a pseudo-Anosov flow, it is proven in \cite{circleatinfinity} that the orbit equivalence class of most pseudo-Anosov flows in dimension 3 -except from those whose bifoliated planes contain a tree of lozenges- is completely determined by the elements of the fundamental group acting on  $\mathbb{S}_{\infty}$ with fixed points. 
\end{enumerate}

Despite the variety of approaches towards a classification of pseudo-Anosov flows in dimension 3, the question of classification remains still open. In fact, till this moment little is known about pseudo-Anosov flows on toroidal manifolds containing simultaneously Seifert and atoroidal pieces. It is for this reason that we would like in this paper to adapt the classification approach developed in \cite{preprint} in the case pseudo-Anosov flows, thus providing a classification method that is a priori independant from the underlying topology of the manifold. This classification method has already been shown to provide a complete classification  of structurally stable diffeomorphisms on surfaces (see \cite{Asterisque} and \cite{Beguin}) and also to admit several applications in the study of Anosov flows on 3-manifolds. For instance, thanks to the results of \cite{preprint}, we have recently provided:
\begin{itemize}
    \item a way to recursively construct the fundamental groups of all 3-manifolds admitting transitive Anosov flows (see \cite{fundgroup})
    \item a classification of all totally periodic Anosov flows in graph manifolds (see \cite{totperiodicioannis}) 
    \item a necessary and sufficient condition for a group action on the plane to be realised as the action of the fundamental group on the bifoliated plane of an Anosov flow (see \cite{realisableactions})
\end{itemize}

\subsection{Markovian families and geometric types}
Fix $M$ a closed, orientable, smooth 3-manifold, $\Phi$ a  pseudo-Anosov flow on $M$ (with no 1-prongs) and $F^s$, $F^u$ its stable and unstable foliations.  

In \cite{Gabai}, Gabai and Oertel show that the universal cover of $M$ is homeomorphic to $\mathbb{R}^3$. Furthermore, by generalizing the techniques used in \cite{Fe1} and \cite{Ba1}, one can show that the space of orbits of  $\widetilde{\Phi}$, the lift of $\Phi$ on the universal cover of $M$, is  homeomorphic to $\mathbb{R}^2$ (see Theorem  \ref{thdef.bifoliatedplane}). The space of orbits of $\widetilde{\Phi}$ is therefore a plane endowed with the natural quotient of the lift of the stable and unstable manifolds of $\Phi$ on $\mathbb{R}^3$. In other words, the flow $\Phi$ is naturally associated to a plane $\mathcal{P}$ endowed with pair of transverse singular foliations $\mathcal{F}^s$ and $\mathcal{F}^u$. We call $(\mathcal{P},\mathcal{F}^s,\mathcal{F}^u)$ the \emph{bifoliated plane of $\Phi$}. 

Even though the flow $\widetilde{\Phi}$ no longer exists on $\mathcal{P}$, the bifoliated plane of $\Phi$ in naturally endowed with dynamics. The fundamental group of $M$ acts on $\widetilde{M}=\mathbb{R}^3$ by preserving the orbits of the lifted flow; hence the action of $\pi_1(M)$ on $\widetilde{M}$ descends to an action  on $\mathcal{P}$. By adapting the arguments used in \cite{Ba1}, it is known that the bifoliated plane $(\mathcal{P},\mathcal{F}^s,\mathcal{F}^u)$ together with the previous action of $\pi_1(M)$ describes completely the flow $\Phi$ up to orbital equivalence (see Theorem \ref{t.barbottheor}). Therefore,  a pseudo-Anosov flow in dimension 3 up to orbital equivalence can be thought of as a group action on a two dimensional plane. 

The previous fact was one of our main motivations for defining the notion of \emph{Markovian family}, an object generalizing the notion of Markov partition for a group action on a bifoliated plane: 

\begin{defi*}
A \emph{Markovian family of rectangles} in $\mathcal{P}$ is a set of mutually distinct rectangles (see Definition \ref{d.markovfamily}) $(R_i)_{i \in I}$ covering $\mathcal{P}$ such that 
\begin{enumerate}
\item $(R_i)_{i \in I}$ is the union of a finite number of orbits of rectangles of the action by $\pi_1(M)$ 
\item For every two distinct rectangles $R_i,R_j$ in $(R_i)_{i \in I}$, if $\overset{\circ}{R_i } \cap \overset{\circ}{R_j} \neq \emptyset$, then $R_i \cap R_j$ is a non-trivial horizontal subrectangle of $R_i$ (or $R_j$ resp.) and a non-trivial vertical subrectangle of $R_j$ (or $R_i$ resp.) 
\item Take any point $x\in \mathcal{P}$ and $U$ any compact neighborhood of $x$ in $\mathcal{P}$. The neighborhood $U$ can be covered by a finite number of rectangles of the family $(R_i)_{i \in I}$ 
\end{enumerate} 
\end{defi*}
For a more detailed definition, see Definition \ref{d.markovfamily}.

A Markovian family has many points in common with a Markov partition. More specifically, after our introduction to pseudo-Anosov flow theory, we will begin the main part of this paper by showing that, as a consequence of the existence of Markov partitions for all pseudo-Anosov flows : 

 \vspace{0.20cm}
\noindent
\textbf{Proposition.}\textit{ Any pseudo-Anosov flow in dimension 3 admits infinitely many Markovian families.}
 \vspace{0.20cm}
 
Furthermore, just as in the case of Markov partitions, one can define the notion of first return map for a Markovian family $\mathcal{R}$ and codify its action on the rectangles of $\mathcal{R}$ by a finite combinatorial object, called a \emph{geometric type}:

\begin{defi*}[Informal]
Take $n \in \mathbb{N}^{*}$ and  $(h_1,v_1),...,(h_n,v_n) \in {(\mathbb{N}^*)}^2$ such that $\sum_i h_i =\sum_i v_i$. Consider a bijection $\phi$ of the two sets $$\mathcal{H}=\lbrace (i,j)| i\in\llbracket 1,n \rrbracket, j \in \llbracket 1,h_i \rrbracket  \rbrace  \text{ and } \mathcal{V}=\lbrace (i,j)| i\in\llbracket 1,n \rrbracket, j \in \llbracket 1,v_i \rrbracket  \rbrace$$ and $u$ a function from the set  $\mathcal{H}$ to $\lbrace -1,+1\rbrace$. The data $G:=(n,(h_i)_{i \in \llbracket 1,n \rrbracket}, (v_i)_{i\in \llbracket 1,n \rrbracket}, \mathcal{H}, \mathcal{V},\phi, u)$ will be called a \emph{geometric type}. 

We will say that the geometric type $(n',(h'_i)_{i \in \llbracket 1,n' \rrbracket}, (v'_i)_{i\in \llbracket 1,n' \rrbracket}, \mathcal{H}', \mathcal{V}',\phi', u')$ is \emph{equivalent} to $G$ when 
\begin{itemize}
    \item $n=n'$
    \item up to reindexation $(h_i,v_i)=(h'_i,v'_i)$ for every $i\in \llbracket 1,n\rrbracket$
    \item there exist two bijections $H_h: \mathcal{H}\rightarrow \mathcal{H}'$ and $H_v: \mathcal{V}\rightarrow \mathcal{V}'$ such that for every $i$  
    \begin{itemize}
        \item either $H_h(i, j)= (i,j)$ for every $j\in \llbracket 1, h_i\rrbracket$ or $H_h(i, j)= (i,h_i-j)$ for every $j\in \llbracket 1, h_i\rrbracket$
        \item either $H_v(i, j)= (i,j)$ for every $j\in \llbracket 1, v_i\rrbracket$ or $H_v(i, j)= (i,v_i-j)$ for every $j\in \llbracket 1, v_i\rrbracket$
        \item $H_h$ and $H_v$ respect $\phi, \phi',u,u'$ 
    \end{itemize}
\end{itemize} (see Figure \ref{f.examplegeometrictype} and Definitions \ref{d.geometrictype} and \ref{d.equivalentgeomtypes} for more details). 
\end{defi*}
Even though the definitions of a geometric type and of the equivalence of geometric types are purely combinatorial, both of these notions admit simple geometric interpretations that will be provided in Section \ref{s.preliminaries}. Furthermore, using the above definition of equivalence, it is easy to see that an equivalence class of geometric types is necessarily finite and that checking whether or not two geometric types belong in the same equivalence class can be easily done  algorithmically. Finally, equivalence classes of geometric types naturally appear inside every pseudo-Anosov flow: 

 \textbf{Theorem-Definition A.}\textit{
Let $M$ be an orientable and closed 3-manifold and $\Phi$ a
pseudo-Anosov flow on $M$. To any Markovian family $\mathcal{R}$ of $\Phi$ we can canonically associate a finite set of pairwise equivalent geometric types, called the geometric types of $\mathcal{R}$ or the geometric types associated to $\mathcal{R}$.}     
\vspace{0.20cm}
 
For a more precise statement, see Definition \ref{d.geometrictypemarkovfamily} and Theorem \ref{t.associatemarkovfamiliestogeometrictype}. The above result will be proven in a slightly more general context involving actions on the plane admitting Markovian families. 

A natural question arising from Theorem-Definition A is whether or not a geometric type of a Markovian family of $\Phi$ contains enough information to describe $\Phi$ up to orbital equivalence. We will prove that the answer to this question is negative. More specifically, as in the case of Anosov flows, by adapting our proof in \cite{preprint}, we will prove that a geometric type associated to a Markovian family $\mathcal{R}$ of $\Phi$ describes completely $\Phi$ up to specific Dehn-Goodman-Fried surgeries (see Section \ref{s.surgeries} for a definition). In order to do that, we will first establish that: 

 \vspace{0.20cm}
\noindent
\textbf{Proposition. }\textit{Let $M$ be an orientable and closed 3-manifold, $\Phi$ a pseudo-Anosov flow on $M$ and $\mathcal{P}$ its bifoliated plane. Denote by $\widetilde{\Phi}$ the lift of $\Phi$ on the universal cover of $M$. Consider $\mathcal{R}$ any Markovian family  of $\Phi$ and $\Gamma$ the set of points  $x\in \mathcal{P}$ with the following two properties: 
\begin{enumerate}
    \item $x$ corresponds to an orbit of $\widetilde{\Phi}$ that projects to a periodic orbit of $\Phi$ in $M$
    \item $x$ does not intersect the interior of some rectangle in $\mathcal{R}$
\end{enumerate}
We have that $\Gamma$ is a $\pi_1(M)$-invariant and non-empty set, which is finite up to the action of $\pi_1(M)$.}
 \vspace{0.20cm}

In other words, to any Markovian family $\mathcal{R}$ of $\Phi$  we can canonically associate a finite number of periodic orbits of $\Phi$ that ``never intersect the interiors of the rectangles of $\mathcal{R}$''.  The previous set of periodic orbits will be called the set of \emph{boundary periodic orbits of $\mathcal{R}$}. See Definition \ref{d.boundaryperiodic} and  Proposition \ref{p.boundaryperiodicpoints} for more details. Once the previous proposition is proved, we will explain how one can adapt the techniques in \cite{preprint} in order to prove that: 

\vspace{0.20cm}
\textbf{Theorem B}
\textit{Let $\Phi_1,\Phi_2$ be two pseudo-Anosov flows on the orientable and closed manifolds $M_1^3, M_2^3$, $\mathcal{R}_1,\mathcal{R}_2$ two Markovian families in their bifoliated planes $\mathcal{P}_1,\mathcal{P}_2$ and $\Gamma_1^{M_1},\Gamma_2^{M_2}$ their associated boundary periodic orbits in $M_1$ and $M_2$. If the equivalence classes of geometric types associated to $\mathcal{R}_1$ and $\mathcal{R}_2$ are the same, then $M_1-\Gamma_1^{M_1} \approx M_2-\Gamma_2^{M_2}$ and $\Phi_1$ (up to orbital equivalence) can be obtained from $\Phi_2$ by performing a finite number of Dehn-Goodman-Fried surgeries on $\Gamma_2^{M_2}$. }
\vspace{0.20cm}

Once again the previous result will be shown in a slightly more general context involving Markovian actions on the plane. One of the difficulties of the above proof relies on the fact that we need to compare two pseudo-Anosov flows up to surgeries. In order to achieve this, we will make use of the \emph{bifoliated plane up to surgeries of $\Phi$} and of the following generalization of a classical theorem by Barbot (see Theorem \ref{t.barbottheor}): 

\begin{defi*}[Informal] Consider $\Gamma^M$ a finite set of periodic orbits of $\Phi$, $\widetilde{\Gamma}$ their lifts on the universal cover of $M$ and $\Gamma$ their projections on $\mathcal{P}$. We define the \emph{bifoliated plane $\clos{\mathcal{P}}$ of $\Phi$ up to surgeries on $\Gamma$} as the universal cover of $\mathcal{P}-\Gamma$ together with some points at infinity. 
\end{defi*}
See Section \ref{s.constructionPbar} for more details. 

\vspace{0.20cm}
\textbf{Theorem C}
\textit{Let $\Phi$ be a pseudo-Anosov flow on the orientable and closed manifold $M^3$ and $\Gamma$ a finite set of periodic orbits of $\Phi$. The bifoliated plane $\clos{\mathcal{P}}$ of $\Phi$ up to surgeries on $\Gamma$ can be endowed with two transverse singular foliations and an action of $\pi_1(M-\Gamma)$ by homeomorphisms. Together with this action and those foliations, $\clos{\mathcal{P}}$ describes completely the flow $\Phi$ up to orbital equivalence and up to Dehn-Goodman-Fried surgeries on $\Gamma$.}
\vspace{0.20cm}

See Theorem \ref{t.generalbarbot} for more details. According to Theorem B, a geometric type of $\Phi$ gives enough information in order to identify the original pseudo-Anosov flow with the exception of the tubular neighborhoods of finitely many periodic orbits of the flow. By including in the geometric type some additional information that will allow us to reconstruct the previous neighborhoods, we can transform the geometric type into a finite combinatorial  invariant up to orbital equivalence: 

\begin{defi*}[Informal]
    A geometric type $(R_1,...,R_n,(h_i)_{i \in \llbracket 1,n \rrbracket}, (v_i)_{i\in \llbracket 1,n \rrbracket}, \mathcal{H}, \mathcal{V},\phi, u)$ endowed with a finite number of sequences of rectangles in $\mathcal{H} \cup \mathcal{V} \cup (R_i)_{i \in \llbracket 1,n \rrbracket}$ will be called a \emph{geometric type with cycles}. 
\end{defi*}
See Definition \ref{d.geomtypecycles} for more details. Using geometric types with cycles and our proofs of Theorem-Definition A and Theorem B, we will show the following refined versions of the previous theorems: 

\vspace{0.20cm}
\textbf{Theorem-Definition D}\textit{
Let $M$ be an orientable and closed 3-manifold and $\Phi$ a 
pseudo-Anosov flow on $M$. To any Markovian family $\mathcal{R}$ of $\Phi$ we can canonically associate a finite set of pairwise equivalent geometric types with cycles, called the geometric types with cycles of $\mathcal{R}$ or the geometric types with cycles associated to $\mathcal{R}$.}  
\vspace{0.20cm}

See Definition \ref{d.geomtypecyclesmarkov} and Theorem \ref{t.geometrictypeswithcyclesclass} for more details. 

\textbf{Theorem E}\textit{
Let $\Phi_1,\Phi_2$ be two pseudo-Anosov flows on the orientable and closed manifolds $M_1^3,M_2^3$ and $\mathcal{R}_1,\mathcal{R}_2$ two Markovian families in $\mathcal{P}_1,\mathcal{P}_2$. If the equivalence classes of geometric types with cycles associated to  $\mathcal{R}_1,\mathcal{R}_2$ are the same, then $\Phi_1$ is orbitally equivalent to $\Phi_2$.}
\vspace{0.20cm}

Once again the two previous results will be shown in a more general context involving actions on the plane that leave invariant some Markovian family. 

\textit{Acknowledgments}. We would like to address a special thanks to Fran\c{c}ois B\'eguin for sharing his intuition around this work and for his very useful advice and ideas concerning the presentation of our results. The completion of this work was made possible thanks to support of the ERC project 818737: Emergence of wild differentiable dynamical systems. 

\newpage
\tableofcontents
\newpage

\section{An introduction to pseudo-Anosov flow theory}
In this first section of this paper, due to the shortage of bibliographic references concerning pseudo-Anosov flows on 3-manifolds, we decided to provide a short  (and rather incomplete) introduction to pseudo-Anosov flow theory in dimension 3. 

Several proofs included in this section are direct adaptations of their equivalents for Anosov flows in dimension 3. Nevertheless, the inclusion of this introductory part was made necessary for us, because of the existence of several definitions for a pseudo-Anosov flow in the literature, some of which are not yet known (by the author at least) to be equivalent in the context of topological non-transitive pseudo-Anosov flows in dimension 3.

\subsection{The definition of a pseudo-Anosov homeomorphism}

Consider $\mathcal{F}^2_h,\mathcal{F}^2_v$ the foliations by horizontal and vertical lines on $\mathbb{R}^2=\mathbb{C}$, $\mathcal{D}_2$ the euclidean open square $\{z\in\mathbb{C}||\text{Re}(z)|<1 \text{, } |{Im}(z)|<1\}$ and $\pi_p: \mathbb{C}\rightarrow \mathbb{C}$ the map associating $z$ to $z^p$, where $p\in \mathbb{N}^*$. Let $\mathcal{D}_1:= \pi_2(\mathcal{D}_2)$ and $\mathcal{D}_p:= \pi_p^{-1}(\mathcal{D}_1)$ for any $p\geq 3$. 

\begin{figure}[h]

  \begin{minipage}[ht]{0.4\textwidth}
    \centering 
     \vspace{0.8cm}
    \hspace{-0.5cm}
    \includegraphics[width=0.9\textwidth]{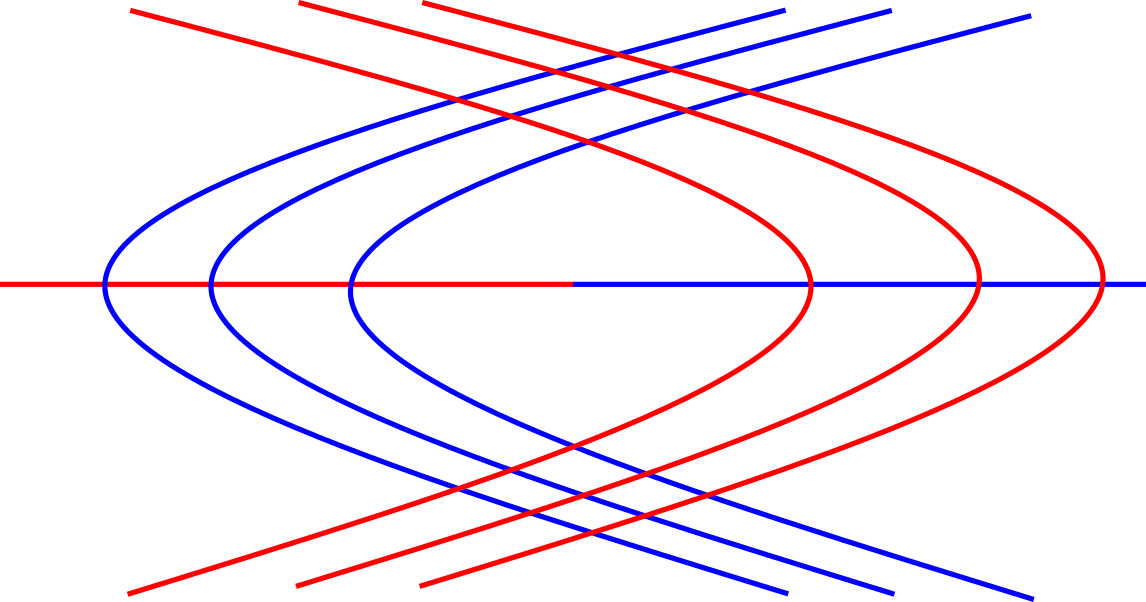}
  \hspace{-1cm}

  \end{minipage}
 \begin{minipage}[ht]{0.4\textwidth}
 \centering
    \includegraphics[width=0.8\textwidth]{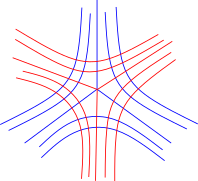}
    \hspace{-1cm}
    
  \end{minipage}
\caption{On the right, the foliations $\mathcal{F}^1_h,\mathcal{F}^1_v$. On the left, $\mathcal{F}^3_h,\mathcal{F}^3_v$ }
\label{f.prongsingularities}
  \end{figure}

The set of images of the leaves of $\mathcal{F}^2_h,\mathcal{F}^2_v$ by $\pi_2$ define two singular foliations $\mathcal{F}^1_h,\mathcal{F}^1_v$ (see Figure \ref{f.prongsingularities}). Similarly, for every $p\geq 3$, the set of pre-images of the leaves $\mathcal{F}^1_h,\mathcal{F}^1_v$ by $\pi_p$, define two singular foliations, say $\mathcal{F}^p_h,\mathcal{F}^p_v$ (see Figure \ref{f.prongsingularities}).  

\begin{defi}\label{d.singularfolisurfaces}
    Consider $\Sigma$ a (not necessarily closed) surface with empty boundary.  
    
    We will say that $\mathcal{G}$ is \emph{a singular foliation of $\Sigma$} if it is a partition of $\Sigma$ into sets, called \emph{leaves}, satisfying the following property: for every $x\in \Sigma$ there exists $U_x$ a neighborhood of $x$ in $\Sigma$, $k\geq 1$ and $h:U_x\rightarrow \mathcal{D}_k$ a homeomorphism such that $h(x)=0\in \mathbb{C}$ and $$h(\mathcal{G}\cap U_x)=\mathcal{F}_h^k\cap \mathcal{D}_k$$

    \noindent In other words, we ask that $h$ sends the restriction of the leaves of $\mathcal{G}$ on $U_x$ to the restriction of the leaves of $\mathcal{F}_h^k$ on $\mathcal{D}_k$. 

    \vspace{0.2cm}

   Furthermore, if $k\neq 2$, then we will call $x$ a \emph{$k-$prong singularity of $\mathcal{G}$} or more simply \emph{a singularity of $\mathcal{G}$}. Any other point of $\Sigma$ will be called a \emph{regular point of $\mathcal{G}$}. We will denote by $\text{Sing}(\mathcal{G})$ the set of singularities of $\mathcal{G}$. Similarly, we will say that a leaf of $\mathcal{G}$ is \emph{singular} if it contains a point in $\text{Sing}(\mathcal{G})$ and any other leaf of $\mathcal{G}$ will be called \emph{regular}. 
   
   Finally, if $x$ is a singular point in a leaf $L$ of $\mathcal{G}$, then the closure of any connected component of $L-\{x\}$ will be called a \emph{prong} of $x$.

 \end{defi}

The above definition of singular foliation is rather restrictive, as it forces all the singular leaves to be homeomorphic to the union of countably many segments intersecting along regular points or $k$-prong singularities. Even more, our previous definition disallows the existence of several other types of singularities that can arise in two dimensional foliations (elliptic singularities, source/sink types singularities, etc). Despite its restrictive character, by convention, unless explicitely mentionned otherwise, whenever we will refer in this paper to a singular foliation, we will mean a foliation satisfying Definition \ref{d.singularfolisurfaces}. 

Furthermore, it is not difficult to see using the previous definition that: 
\begin{remadef}\label{r.singularfolisurfaces}
    If $\mathcal{G}$ is a singular foliation on the surface $\Sigma$: 

    \begin{itemize}
        \item $\text{Sing}(\mathcal{G})$ forms a discrete set in $\Sigma$, thus the set of singular leaves of $\mathcal{G}$ is at most countable
        \item a singular leaf of $\mathcal{G}$ can a priori contain more than one singular point
        \item outside $\text{Sing}(\mathcal{G})$, the singular foliation $\mathcal{G}$ defines a $C^0$ (regular) foliation on $\Sigma-\text{Sing}(\mathcal{G})$. We will say that $\mathcal{G}$ is \emph{orientable} (resp. \emph{transversely orientable}) if its restriction on $M-\text{Sing}(\mathcal{G})$ is orientable (resp. transversely orientable)
    \end{itemize}
\end{remadef}

\begin{defi}\label{d.transversefolisurfaces}
   Consider $\mathcal{F}, \mathcal{G}$ two singular foliations on a surface $\Sigma$ with empty boundary. We will say that $\mathcal{F}, \mathcal{G}$ are \emph{transverse} if for every point $x\in \Sigma$ there exists $U_x$ a neighborhood of $x$ in $\Sigma$, $k\geq 1$ and $h:U_x\rightarrow \mathcal{D}_k$ a homeomorphism such that $h(x)=0$, $h(\mathcal{F}\cap U_x)=\mathcal{F}_h^k\cap \mathcal{D}_k$ and $h(\mathcal{G}\cap U_x)=\mathcal{F}_v^k\cap \mathcal{D}_k$. 
\end{defi}
\begin{defi}\label{d.pseudohomeo}

Let $\Sigma$ be a closed surface. A homeomorphism $f: \Sigma \rightarrow \Sigma$ is \emph{pseudo-Anosov} if there exists a pair of $f$-invariant singular transverse foliations $F^s$ and $F^u$ on $\Sigma$ with no $1$-prong singularties such that: 
\begin{itemize}
    \item $F^s$ (resp. $F^u$) is equipped with a transverse non-atomic (Borel) measure $\mu^s$ (resp. $\mu^u$) with full support assigning finite measure to any compact arc transverse to $F^s$ (resp. $F^u$)
    \item there exists $\lambda>1$ such that $f_{\star}(\mu^s)=\lambda \mu^s$ and $f_{\star}(\mu^u)=\lambda^{-1} \mu^u$ 
\end{itemize} 

If $F^s, F^u$ are allowed to have one prongs, we will call $f$ a \emph{pseudo-Anosov homeomorphism with one prongs}.
\end{defi}

Notice that as a result of the above definition, the pseudo-Anosov property of a homeomorphism is invariant by $C^0$-conjugation. 

Several different definitions of a pseudo-Anosov homeomorphism exist in the literature. Following our previous notations, some authors (see for instance \cite{GerberKatok}) require in the definition of a pseudo-Anosov homeomorphism (with no one prongs) that:
\begin{itemize}
    \item for every regular point $x$ of $F^s,F^u$ the map $h: U_x\rightarrow \mathcal{D}_2$ be a $C^{\infty}$ diffeomorphism such that $h_{*}(\mu^s), h_{*}(\mu^u)$ coincide with the Lebesgue measures $|\text{Im}(dz)|, |\text{Re}(dz)|$
    \item for every $p$-prong singularity $x$ of $F^s,F^u$ (with $p\geq 3$) the map $h: U_x\rightarrow \mathcal{D}_p$ be a $C^{\infty}$ diffeomorphism such that $ h_{*}(\mu^s),  h_{*}(\mu^u)$ coincide with the measures $|\text{Im}(dz^{p/2})|, |\text{Re}(dz^{p/2})|$
\end{itemize}

We would like to point out at this moment that:
\begin{prop}\label{p.gerberkatokdefinitionisequiv}
    Up to conjugation, any pseudo-Anosov homeomorphism $f$ in the sense of Definition \ref{d.pseudohomeo} has the above properties. 
\end{prop}
Before beginning the proof of this proposition, let us first give two useful definitions: 
\begin{defi}\label{d.standardpolygon}
Consider $k\geq 2$ and the previously defined foliations on the plane  $\mathcal{F}_{h,v}^k$. A closed topological disk $P$ in the plane will be called a \emph{polygon} if it is bounded by segments of $\mathcal{F}_{h}^k$ and $\mathcal{F}_{v}^k$ that alternate between them (see Figure \ref{f.polygons}). We will say that $P$ is a \emph{standard polygon} if one of two following conditions is verified:
    \begin{itemize}
        \item $P$ contains no singularity in its interior and its boundary consists of $2$ segments of $\mathcal{F}_{h}^k$ and $2$ segments of $\mathcal{F}_{v}^k$. In this case, $P$ is trivially bifoliated and is called a \emph{rectangle}.
        \item $P$ contains in its interior a $k$-prong singularity ($k\geq 3$) and its boundary consists of $k$ segments of $\mathcal{F}_{h}^k$ and $k$ segments of $\mathcal{F}_{v}^k$. In this case, we will call $P$ a \emph{standard $2k$-gon}.
    \end{itemize}

Consider now $\mathcal{F}, \mathcal{G}$ two singular transverse foliations on a surface $\Sigma$ with empty boundary. A closed topological disk $P$ in $\Sigma$ will be called a \emph{polygon} (resp. \emph{standard polygon}, \emph{rectangle}, \emph{standard $2k$-gon}) if there exists $U$ a neighborhood of $P$, $k\geq 2$, $V$ a neighborhood of $0\in \mathbb{C}$ and a homeomorphism $h:U\rightarrow V$ such that : 
\begin{itemize}
    \item $h$ sends the restriction of the leaves of $\mathcal{F}$ (resp. $ \mathcal{G}$) on $U$ to the restriction of the leaves of $\mathcal{F}_{h}^k$  (resp. $\mathcal{F}_{v}^k$) on $V$ 
    \item the image of $P$ by $h$ is a polygon (resp. standard polygon, rectangle, standard $2k$-gon) in the plane
\end{itemize}
\end{defi}
\begin{defi} \label{d.segment}
Let $F^s,F^u$ be the stable and unstable foliations of some pseudo-Anosov homeomorphism on a closed surface $\Sigma$. Any segment contained in a leaf of $F^s$ (resp. $F^u$) will be called a \emph{stable} (resp. \emph{unstable}) \emph{segment}. 

Moreover, if $P$ is a standard polygon in $\Sigma$, any connected component of $F^s\cap P$ (resp. $F^u\cap P$) will be called a \emph{stable} (resp. \emph{unstable}) \emph{leaf} of $P$. For every $z\in P$ we will denote by $P^s(z)$ (resp. $P^u(z)$) the stable (resp. unstable) leaf of $P$ containing $z$ and we will call the closure of each connected component of $P-(P^s(z)\cup P^u(z))$ a \emph{quadrant of $z$ in $P$}.
\begin{figure}[h]

  \begin{minipage}[ht]{0.4\textwidth}
    \centering 
     \vspace{0.8cm}
    \hspace{-0.5cm}
    \includegraphics[width=0.28\textwidth]{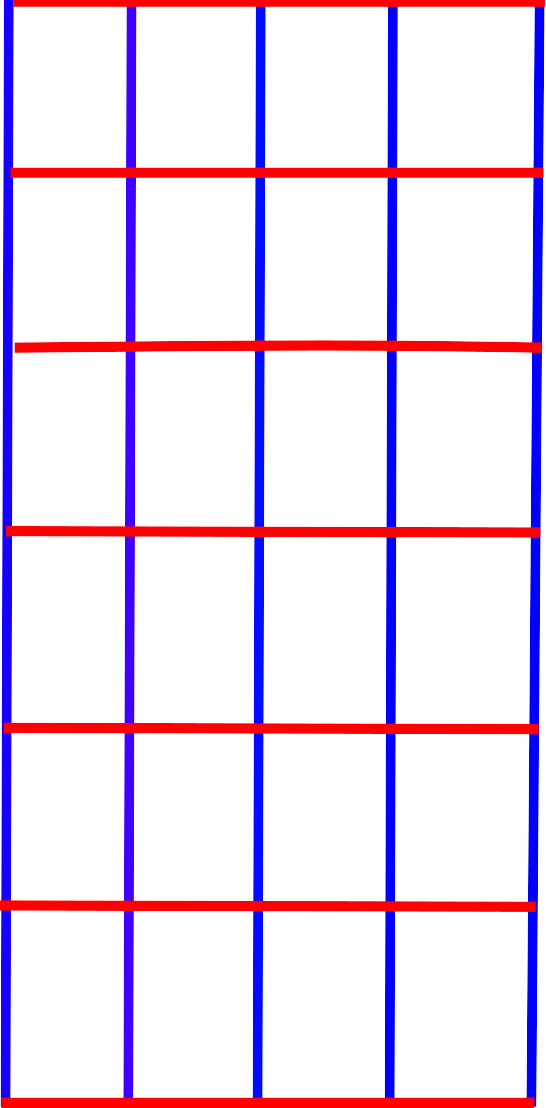}
  \hspace{-1cm}
    \caption*{\quad (a)}
    
  \end{minipage}
 \begin{minipage}[ht]{0.4\textwidth}
 \centering
 \vspace{0.9cm}
    \includegraphics[width=0.55\textwidth]{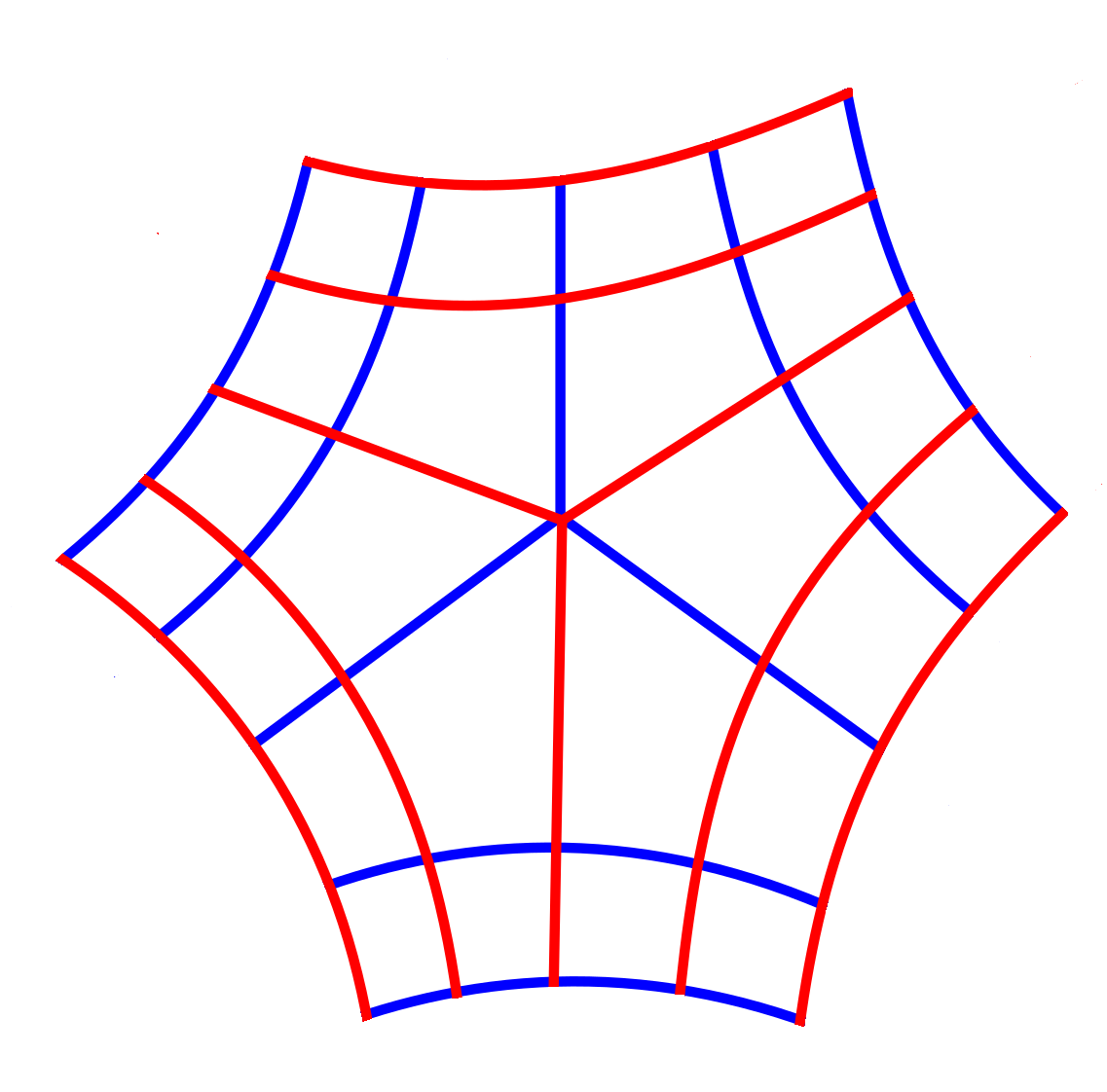}
    
    \caption*{(b)}
    
  \end{minipage}
\caption{Two standard polygons in $\Sigma$}
\label{f.polygons}
  \end{figure}
\end{defi}

\begin{proof}[Proof of Proposition \ref{p.gerberkatokdefinitionisequiv}]
    It suffices to prove that $f$ has the desired properties after one changes appropriately the differentiable structure of $\Sigma$. Denote by $F^s,F^u$ the stable and unstable foliations of $f$. Thanks to Definition \ref{d.transversefolisurfaces}, since $F^s,F^u$ are singular and transverse, every point in $\Sigma$ is contained in the interior of a standard polygon. Using the compactness of $\Sigma$, cover $\Sigma$ by a finite number of standard polygons $P_1,...,P_n$ such that 
\begin{itemize}
    \item every point in $\Sigma$ is contained in the interior of some $P_i$ and 
    \item each singularity of $F^{s,u}$ is contained in a unique $P_i$
    \item if $P_i$ and $P_j$ contain both a singularity of $F^{s,u}$, then $P_i\cap P_j=\emptyset$
    \item if $\inte{P_i}\cap\inte{P_j}\neq \emptyset$, then $\inte{P_i}\cap\inte{P_j}$ is connected
\end{itemize}
We will use the $P_1,...,P_n$ in order to define a new atlas on $\Sigma$. Let us construct a chart $\phi_i: P_i\rightarrow \mathbb{R}^2$ for every $P_i$

Assume first that $P_i$ is a rectangle (therefore it does not contain any singularity of $F^{s,u}$). Consider $s_1,s_2$ (resp. $u_1$,$u_2$) the stable (resp. unstable) segments in the boundary of $P_i$ and $O$ the corner of $P_i$ in $s_1\cap u_1$. Notice that every point $X\in P_i$ can be seen as the point of intersection of a unique pair of leaves $P^u(x(X)), P^s(y(X))$ in $P_i$,  where $x(X)\in s_1$ and $y(X)\in u_1$. We define $$\phi_i(X)=\left(\mu^u([O,x(X)]^s), \mu^s([O,y(X)]^u)\right)\in \mathbb{R}_{\geq0}^2$$ where $[O,x(X)]^s$ and $[O,y(X)]^u $ are the stable and unstable segments in $P$ going from $O$ to $x(X)$ and $y(X)$ respectively. Thanks to our assumptions on $(F^s,\mu^s),(F^u,\mu^u)$, it is easy to check that $\phi_i$ is a homeomorphism such that $(\phi_i)_{\star}(\mu^s)=|\text{Im}(dz)|$, $(\phi_i)_{\star}(\mu^u)=|\text{Re}(dz)|$. 

Assume now that $P_i$ is a $2k$-gon ($k\geq 3$) containing a $k$-prong singularity $O$ in its interior. By cutting $P_i$ along $P^s(O),P^u(O)$, we can view $P_i$ as a union of $2k$ rectangles, say $P_i^1,...,P_i^{2k}$ (see Figure \ref{f.polygons}). Similarly, consider $\mathcal{F}^{k}_h,\mathcal{F}^{k}_v$ the pair of transverse foliations defined in the beginning of this section with a $k$-prong singularity at $z=0$. By construction, the stable and unstable leaves of the singularity separate $\mathbb{R}^2$ into $2k$ sectors, namely $S_1,...,S_{2k}$, where $S_j:=\{r\cdot e^{i\theta}|r\in \mathbb{R}^{+}, \theta \in [\frac{2\pi(j-1)}{2k},\frac{2\pi j}{2k}]\}$. By an argument similar to the one given in the previous paragraph, we can construct a homeomorphism $\phi_i:P_i\rightarrow U\subset\mathbb{R}^2$, where $U$ a neighborhood of the origin, such that:
\begin{itemize}
    \item $\phi_i(O)=0$, $\phi_i(P_i^j)\subset S_j\cap U$
    \vspace{0.05cm}
    \item $\phi_i(F^s\cap P_i)=\mathcal{F}^{k}_h\cap U$, $\phi_i(F^u\cap P_i)=\mathcal{F}^{k}_v\cap U$
     \vspace{0.05cm}
    \item $ (\phi_i)_{*}(\mu^s)=|\text{Im}(dz^{k/2})|$, $(\phi_i)_{*}(\mu^u)=|\text{Re}(dz^{k/2})|$
\end{itemize}

Let us now show that $(\inte{P_1},\phi_1),...,(\inte{P_n},\phi_n)$ defines a $C^{\infty}$ atlas on $\Sigma$. Indeed, by our assumptions $\underset{i=1}{\overset{n}{\cup}}\inte{P_i}=\Sigma$. Next, consider $i\neq j$ such that $P_i\cap P_j\neq \emptyset$. 

Assume first that $P_i,P_j$ contain no singularities. In this case, $\phi_i\circ \phi_j^{-1}$ preserves the horizontal and vertical directions and the measures $|\text{Re}(dz)|, |\text{Im}(dz)|$. Hence,  $\phi_i\circ \phi_j^{-1}(x,y)=(\epsilon x, \epsilon' y) + v$, where $\epsilon, \epsilon'\in \{-,+\}$,  $v\in\mathbb{R}^2$ and thus $\phi_i\circ \phi_j^{-1}\in C^{\infty}$.

Assume now that $P_i$ contains a $k$-prong singularity. By our intial assumptions, $P_j$ does not contain any singularities of $F^{s,u}$ and $P_i\cap P_j$ contains only regular points of $F^{s,u}$. Denote by $P_i^1,...,P_i^{2k}$ the sub-rectangles of $P_i$ used during our construction of $\phi_i$. Since $P_j$ is also a rectangle and by hypothesis $\inte{P_j}\cap \inte{P_i}$ is connected, $P_j$ intersects the interior of at most two of the rectangles $P_i^1,...,P_i^{2k}$. We will consider the case where $\inte{P_j}\cap \inte{P_i}$ is contained in a unique sub-rectangle of $P_i$, say $P_i^1$ (the case where $P_j$ intersects the interior of two consecutive sub-rectangles of $P_i$ follows from a similar argument). Assume without any loss of generality that $\phi_i(P_i^1)\subset S_1=\{r\cdot e^{i\theta}|r\in \mathbb{R}^{+}, \theta \in [0,\frac{2\pi}{2k}]\}$. Consider $\Phi_{k}:S_1\setminus\{0\}\rightarrow \{z\in \mathbb{C}^*|\text{Re}(z)\geq 0, \text{Im}(z)\geq 0 \}$ the $C^{\infty}$-diffeomorphism given by $\Phi_{k}(z)=z^{k/2}$. Notice that $\Phi_k^{\star}|\text{Re}(dz)|=|\text{Re}(dz^{p/2})|$ and $\Phi_k^{\star}|\text{Im}(dz)|=|\text{Im}(dz^{p/2})|$. Thus, $\Phi_k\circ \phi_i\circ \phi_j^{-1}$ preserves the horizontal and vertical directions and the measures $|\text{Re}(dz)|, |\text{Im}(dz)|$. This implies that $\Phi_k\circ \phi_i\circ \phi_j^{-1}(x,y)=(\delta x, \delta' y) + v'$, where $\delta, \delta'\in \{-,+\}$, $v'\in\mathbb{R}^2$ and thus $\phi_i\circ \phi_j^{-1}\in C^{\infty}$, which finishes the proof.
\end{proof}

It is not difficult to check using the above atlas of $\Sigma$ that: 
\begin{coro}
    Let $f$ be a pseudo-Anosov homeomorphism on a closed smooth surface $\Sigma$ and $F^s,F^u$ be its stable and unstable foliations. We have that $f$ is conjugated to a pseudo-Anosov homeomorphism $g$: 
    \begin{enumerate}
       \item that preserves a pair of singular foliations $G^s,G^u$ that are $C^{\infty}$ on $\Sigma-\text{Sing}(G^{s,u})$ and such that the transverse measures $\nu^s$ and $\nu^u$ of $G^s$, $G^u$ equivalent to Lebesgue measures
        \item that is $C^{\infty}$ on $\Sigma-\text{Sing}(G^{s,u})$
    \end{enumerate}
\end{coro}

Even more, in \cite{GerberKatok} it is proved that: 
\begin{theorem}\label{t.gerberkatok}
    Any pseudo-Anosov homeomorphism on a smooth orientable closed surface $\Sigma$ is conjugated to a $C^{\infty}$ pseudo-Anosov homeomorphism through a homeomorphism of $\Sigma$ isotopic to the identity. 
\end{theorem}

The authors of \cite{GerberKatok} claim that the previous result can be generalized for pseudo-Anosov homeomorphisms with one prongs or pseudo-Anosov homeomorphisms on non-orientable surfaces.  

\subsection{Expansive vs pseudo-Anosov homeomorphisms: an equivalence}
\begin{defi}\label{d.expansivehomeo}
 Let $\Sigma$ be a smooth, closed surface and $h:\Sigma\rightarrow \Sigma$  a homeomorphism. Consider $d$ a distance given by some Riemannian metric on $\Sigma$. We will say that $h$ is \emph{expansive} if there exist $\epsilon>0$ such that for any $x,y \in \Sigma$ we have $$\forall n\in\mathbb{Z}\quad d(f^n(x), f^n(y))<\epsilon\implies x=y $$ 
\end{defi}

Notice that, as $\Sigma$ is closed, the expansive character of a homeomorphism does not depend on our choice of metric, as long as the topology given by the metric coincides with the underlying topology of the surface. This is why one could alter the above definition by asking that $d$ not necessarily originates from a Riemannian metric, but rather defines a topology on $\Sigma$ compatible with the topology arising from its manifold structure. 

The notion of expansivity is intimately related with any pseudo-Anosov homeomorphism. Indeed, according to a classical result in the theory of pseudo-Anosov homeomorphisms: 
\begin{prop}\label{p.pseudoanosovimpliesexpansive}
    Any pseudo-Anosov homeomorphism $f:\Sigma\rightarrow \Sigma$ is expansive. 
\end{prop}
 Conversely, J. Lewowicz and K.Hiraide simultaneously proved in \cite{Lewowicz} and \cite{Hiraide} respectively that : 
\begin{theorem}[Lewowicz\footnote{Contrary to Hiraide's proof, Lewowicz's proof of Theorem \ref{t.expansiveimpliespseudoanosov} supposes that the underlying surface is orientable}, Hiraide]\label{t.expansiveimpliespseudoanosov}
    Every expansive homeomorphism on a closed surface is pseudo-Anosov. 
\end{theorem}
Therefore, thanks to the two previous results, we get that a homeomorphism on a closed surface is pseudo-Anosov if and only if it is expansive. In order to provide some intuition for this close relationship between pseudo-Anosov homeomorphisms and expansiveness we will finish this section with a proof of Proposition \ref{p.pseudoanosovimpliesexpansive}. During the proof of Proposition \ref{p.pseudoanosovimpliesexpansive}, we will admit the following lemma, whose proof we leave to the reader:

\begin{lemm}\label{l.measuredfoli} Consider $f$ a pseudo-Anosov homeomorphism on a closed surface $\Sigma$, $F^s, F^u$ its invariant singular foliations and $\mu^s,\mu^u$ their transverse measures. Endow $\Sigma$ with a distance given by some Riemannian metric. We have that:
    \begin{enumerate}
        \item for any $\epsilon>0$ sufficiently small, there exists $c(\epsilon)>0$ with $ c(\epsilon)\underset{\epsilon\rightarrow 0}{\longrightarrow}0$ such that any stable (resp. unstable) segment $K$ contained inside an $\epsilon$-ball in $\Sigma$ verifies $\mu^u(K)\leq c(\epsilon)$ (resp. $\mu^s(K)\leq c(\epsilon)$)  
        \vspace{0.2cm}
    \item for every $\eta_0>0$ there exists $l(\eta_0)>0$ such that  every stable (resp. unstable) segment $K$ of measure  $\mu^u(K)<\eta_0$ (resp. $\mu^s(K)<\eta_0$) intersects at most once any unstable (resp. stable) segment $L$ of measure $\mu^s(L)<l(\eta_0)$ (resp. $\mu^u(L)<l(\eta_0)$)
    \end{enumerate}
\end{lemm}

\begin{proof}[Proof of Proposition \ref{p.pseudoanosovimpliesexpansive}]
     Denote by $F^s,F^u$ the stable and unstable foliations of $f$, by $\mu^s,\mu^u$ their associated transverse measures and by $B(x,\delta)$ the $\delta$-ball around any point $x$ in $\Sigma$. Endow $\Sigma$ with a distance $d$ given by some Riemannian metric and suppose that there exists a positive sequence $\delta_m \longrightarrow 0$ and $x_m, y_m\in \Sigma$ such that $x_m\neq y_m $ and  $d(f^k(x_m), f^k(y_m))<\delta_m$ for every $m\in \mathbb{N}$ and $k\in \mathbb{Z}$.

Since $\Sigma$ is compact and $F^s,F^u$ are transverse singular foliations, there exists a finite family of standard polygons $(P_i)_{i\in I}$ of $\Sigma$ such that the interiors of the $P_i$ cover $\Sigma$ and such that there exists $\eta>0$ satisfying the following: for every $x\in \Sigma$ we can find $i\in I$ for which $B(x,\eta) \subset \inte{P_i}$. It follows that for $m$ sufficiently big and for every $x\in \Sigma$ there exists $i\in I$ such that $B(x,\delta_m)\subset \inte{P_i}$. We will thus assume without any loss of generality that for every $m\in\mathbb{N}$ and $k\in \mathbb{Z}$ there exists $i\in I$ such that $f^k(x_m),f^k(y_m)\in \inte{P_i}$.

\textbf{Case 1:} There exists $m\in\mathbb{N}$ such that $x_m\in \inte{P_i}$ is a singular periodic point. 

By eventually considering a power of $f$, we can assume that $x_m$ is a singular fixed point and that $f$ preserves the quadrants of $x_m$ inside $P_i$. Since, with respect to the transverse measures $\mu^s$ and $\mu^u$, $f$ stretches $F^u$ and contracts $F^s$ by a factor $\lambda$ and $1/\lambda$ respectively, the dynamics around a fixed singular point behave as in Figure \ref{f.localperiodicsing}. It follows that for $m$ sufficiently big, no orbit of $f$ can stay $\delta_m$-close to a singular point in both the future and the past; hence Case 1 is impossible for $m$ sufficiently big. 

\begin{figure}
    \centering
    \includegraphics[scale=1.2]{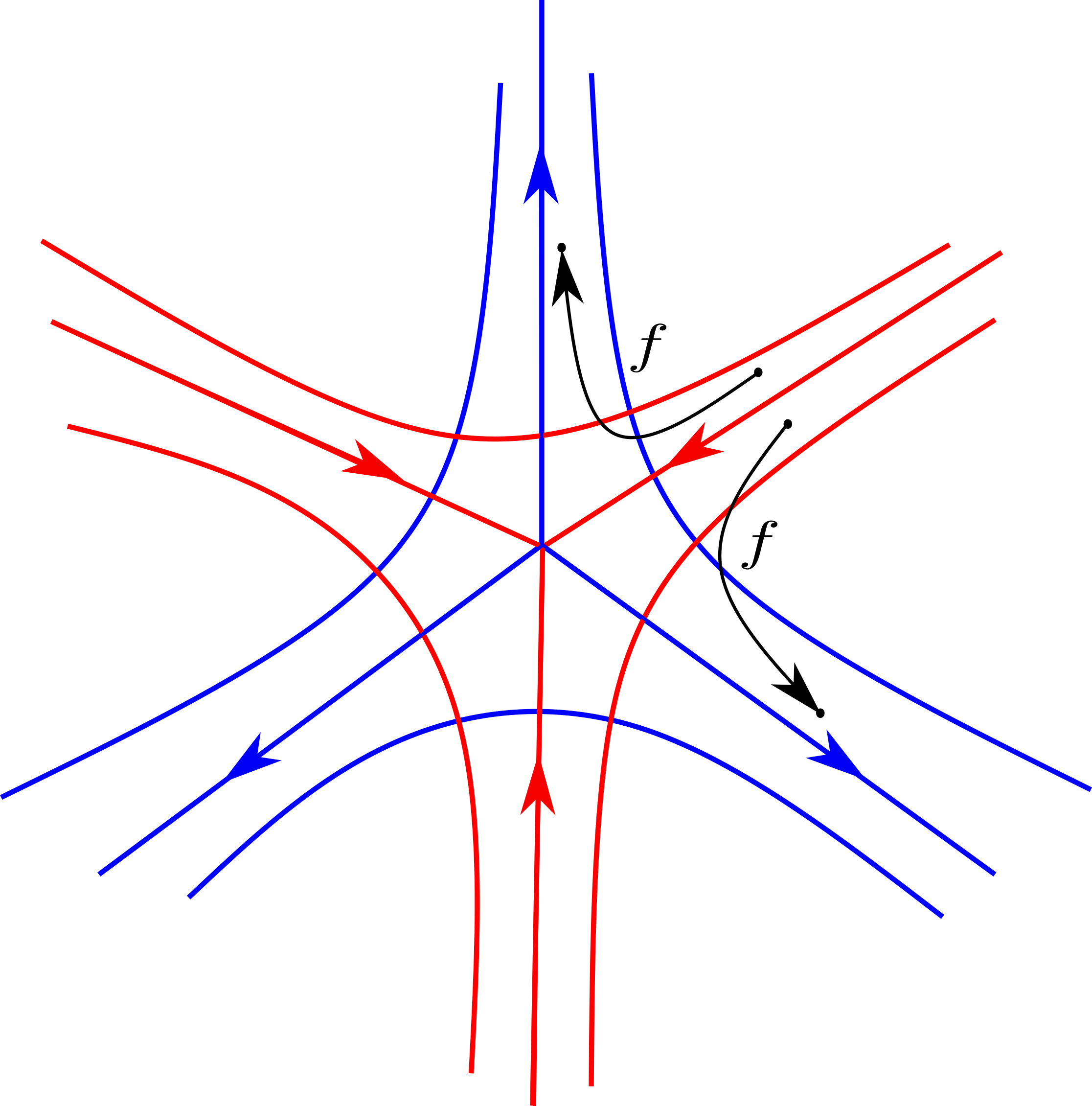}
    \caption{}
    \label{f.localperiodicsing}
\end{figure}
By the exact same argument one can show that for $m$ sufficiently big the following case is also impossible: 

\textbf{Case 2:} There exists $m\in\mathbb{N}$ and $k\in \mathbb{Z}$ such that $f^k(x_m),f^k(y_m)\in \inte{P_i}$, $P_i$ contains in its interior a singular periodic orbit $\gamma$ and there exist $Q,Q'$ two distinct quadrants of $\gamma$ in $P_i$ such that $f^k(x_m)\in Q$, $f^k(y_m)\in Q'$ and $\{f^k(x_m),f^k(y_m)\}\not\subset Q\cap Q'$

Assume $m$ sufficiently big so that Cases 1 and 2 are impossible. Cut every non rectangular polygon belonging in $(P_i)_{i\in I}$ along the quadrants of the unique singular point that it contains. By doing so and by keeping all the rectangles in $(P_i)_{i\in I}$, we have obtained a new family of rectangles $(S_j)_{j\in J}$ such that: 

\textbf{Final case:} For every $m\in\mathbb{N}$ and $k\in \mathbb{Z}$ there exists $j(k,m)\in J$ such that $f^k(x_m),f^k(y_m)\in S_{j(k,m)}$ (see Figure \ref{f.proofexpansive})

For every $z,z'\in S_j$, denote by $[z,z']$ the unique point in $S^s_j(z)\cap S^u_j(z')$ (see Definition \ref{d.segment} for the previous notations). Denote also by  $z_{m,k}$ the point $[f^k(x_m),f^k(y_m)]\in S_{j(k,m)}$ and by $I^s_{m,k}$ (resp. $I^u_{m,k})$ the stable (resp. unstable) segment of $S_{j(k,m)}$ going from $f^k(x_m)$ (resp. $f^k(y_m)$) to $z_{m,k}$. We will prove that for $m$ sufficiently big $f(z_{m,k})=z_{m,(k+1)}$, which will lead us to a contradiction.

\begin{figure}
    \centering
    \includegraphics[scale=2]{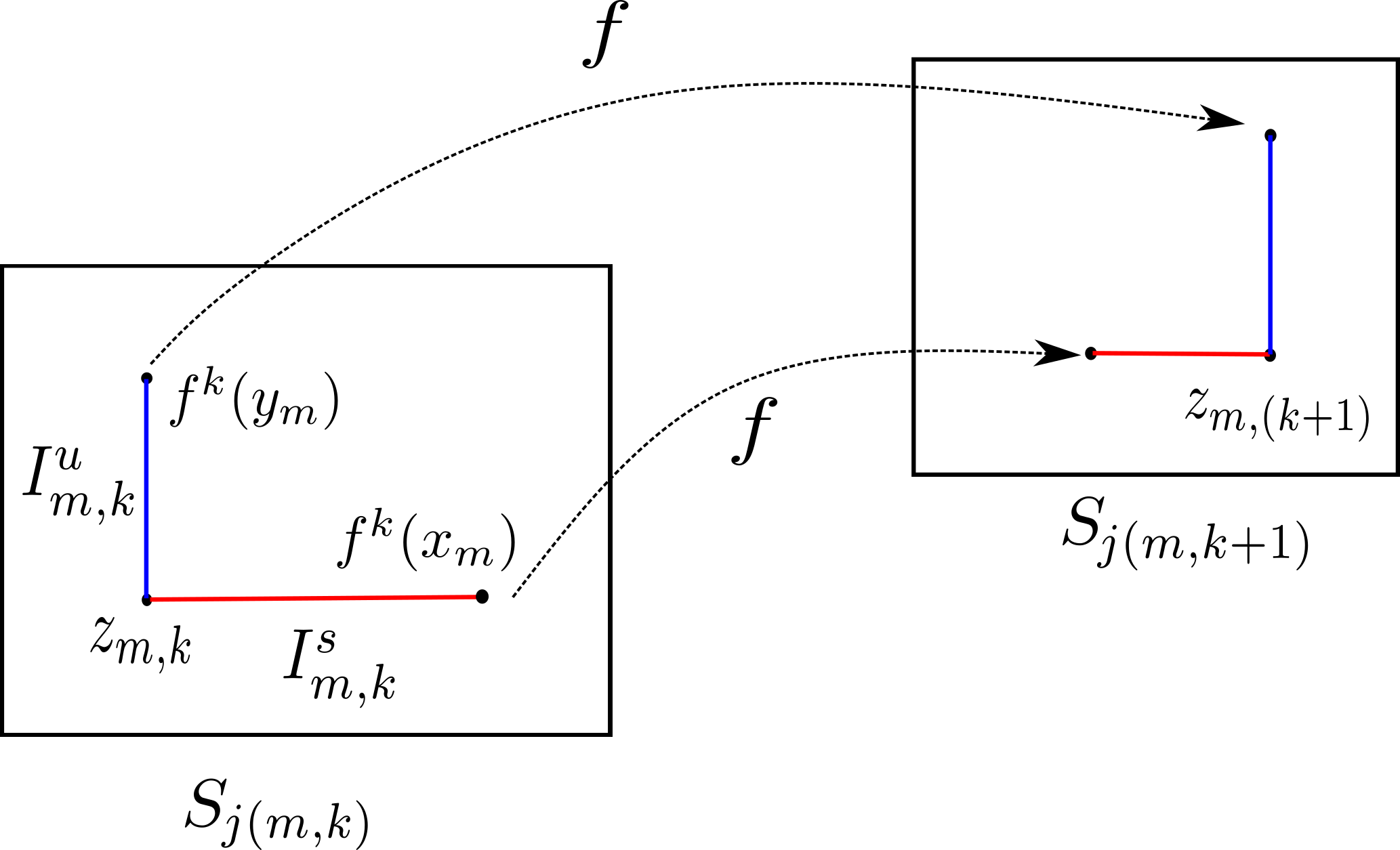}
    \caption{}
    \label{f.proofexpansive}
\end{figure}
Indeed, the function $[\cdot,\cdot]$ being continuous, there exists $\epsilon_m\underset{m\rightarrow +\infty}{\longrightarrow}0$ such that for every $x\in I^s_{m,k}$ and $ y\in I^u_{m,k}$ we have $ d(f^k(x_m), [x,y])<\epsilon_m$  and $d(f^k(y_m), [x,y])<\epsilon_m$. Also, by (1) of Lemma \ref{l.measuredfoli} we get that there exists $c(\epsilon_m)\underset{m\rightarrow +\infty}{\longrightarrow}0$ such that $\mu^u(I^s_{m,k})<c(\epsilon_m)$ and $\mu^s(I^u_{m,k})<c(\epsilon_m)$ for every $m,k$. If $f(z_{m,k})\neq z_{m,(k+1)}$, we get that $I^s_{m,(k+1)}\cup f(I^s_{m,k})$ intersects at least twice $I^u_{m,(k+1)}\cup f(I^u_{m,k})$. The first segment is of measure at most $(1/\lambda +1) \cdot c(\epsilon_m)$ with respect to $\mu^u$ and the second at most $(\lambda+1)\cdot c(\epsilon_m)$ 
 with respect to $\mu^s$, where $\lambda$ is the stretch factor of $f$. For $m$ sufficiently big this contradicts (2) of Lemma \ref{l.measuredfoli}. 
 
 Therefore, for $m$ sufficiently big $f(z_{m,k})=z_{m,(k+1)}$ and thus for every $k\in \mathbb{Z}$ we have $f^k(I^s_{m,0})=I_{m,k}\subset B(f^k(x_m),\epsilon_m)$. But, $\mu^u(f^k(I^s_{m,0}))=\lambda^{-k}\mu^u(I^s_{m,0})$ for every $k\in \mathbb{Z}$, which contradicts (1) of Lemma \ref{l.measuredfoli}.
\end{proof}

\begin{figure}[h!]
    \centering
    \includegraphics[scale=0.2]{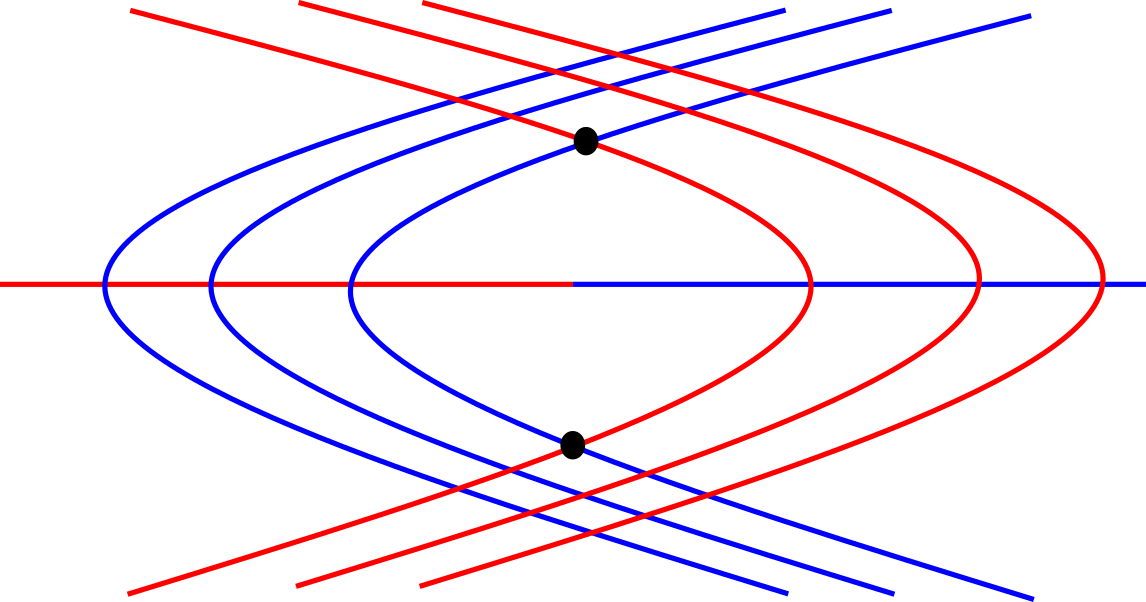}
    \caption{The two points above stay close in both the future and the past}
    \label{f.nonexpansivehomeo}
\end{figure}

Let us remark here that Proposition \ref{p.pseudoanosovimpliesexpansive} does not hold for pseudo-Anosov homeomorphisms with $1$-prongs (see Figure \ref{f.nonexpansivehomeo}). The existence of $1$-prongs changes significantly the dynamical behavior of the system and this is the main reason why we decided to define separately pseudo-Anosov homeomorphisms with and without $1$-prongs in Definition \ref{d.pseudohomeo}. 

\subsection{The definition of a pseudo-Anosov flow}\label{ss.defipseudoanosov}
Recall that in the previous section we introduced $\mathcal{F}^2_h,\mathcal{F}^2_v$ the horizontal and vertical foliations in $\mathbb{R}^2$, $\mathcal{D}_2$ the open euclidean square  $\{z\in\mathbb{C}||\text{Re}(z)|<1 \text{, } |{Im}(z)|<1\}$, $\pi_p(z)=z^p$, $\mathcal{D}_1= \pi_2(\mathcal{D}_2)$ and $\mathcal{D}_i= \pi_i^{-1}(\mathcal{D}_1)$ for any $i\geq 3$. We also defined $\mathcal{F}^1_{h,v}=\pi_2(\mathcal{F}^2_{h,v})$ and $\mathcal{F}^p_{h,v}=\pi_p^{-1}(\mathcal{F}^1_{h,v})$ ($p\geq 3$).

\begin{defi}\label{d.folisingular}

Let $M$ be a $3$-manifold with no boundary. We will say that $\mathcal{F}$ is \emph{a singular codimension one foliation of $M$} if it is a partition  of $M$ into sets, called \emph{leaves}, satisfying the following property: 

     For every $x\in M$ there exists $U_x$ a neighborhood of $x$ in $M$, $k\geq 1$ and $h:U_x\rightarrow \mathcal{D}_k\times [0,1]$ a homeomorphism verifying $h(x)=(0,\frac{1}{2})\in \mathbb{C}\times [0,1]$ and $$h(\mathcal{F}\cap U_x)=(\mathcal{F}_h^k\cap \mathcal{D}_k)\times [0,1]$$

 \noindent In other words, we ask that $h$ sends the restriction of the leaves of $\mathcal{F}$ on $U_x$ to the restriction of the leaves of $\mathcal{F}_h^k\times [0,1]$ on $\mathcal{D}_k\times [0,1]$. 

       Furthermore, following the above notations, if $k\neq 2$, then we will call $x$ a \emph{$k-$prong singularity of $\mathcal{F}$} or more simply \emph{a singularity of $\mathcal{F}$}. Any other point of $M$ will be called a \emph{regular point of $\mathcal{F}$}. We will denote by $\text{Sing}(\mathcal{F})$ the set of connected components of singularities of $\mathcal{F}$. Finally, any leaf of $\mathcal{F}$ containing a singular point of $\mathcal{F}$ will be called \emph{singular} and any other leaf of $\mathcal{F}$ will be called \emph{regular}. 
\end{defi}

As in the case of Definition \ref{d.singularfolisurfaces}, the above definition constitutes a very restrictive definition for a singular codimension one foliation on a 3-manifold. Despite its restrictive character, by convention, unless explicitely mentionned otherwise, whenever we will refer in this paper to a singular codimension one foliation, we will mean a foliation satisfying Definition \ref{d.folisingular}. 

Furthermore, it is not difficult to see using the previous definition that: 
\begin{remadef}\label{r.singularfolidim3}
    If $\mathcal{F}$ is a singular codimension one  foliation in the three manifold $M$ with $\partial M=\emptyset$: 

    \begin{itemize}
        \item each connected component of singularities in $\text{Sing}(\mathcal{F})$ is a properly embedded $1$-dimensional submanifold of $M$ contained in a single leaf of $\mathcal{F}$. It follows that $\text{Sing}(\mathcal{F})$ consists of circles and lines and that if $M$ is closed then $\text{Sing}(\mathcal{F})$ consists only of circles
        \item if $C$ is a connected component of singularities in $\text{Sing}(\mathcal{F})$ contained in the leaf $L$ of $\mathcal{F}$, then $L-C$ consists of finitely many connected components, the closure of each one of which will be called a \emph{separatrix} of $C$
        \item two different connected components of singularties in $\text{Sing}(\mathcal{F})$ can be separated by open sets. It follows that $\text{Sing}(\mathcal{F})$ consists of at most countably many circles and lines and that the set of singular leaves of $\mathcal{F}$ is at most countable. 
        \item if a $k$-prong singularity and a $k'$-prong singularity of $\mathcal{F}$ lie in the same connected component of $\text{Sing}(\mathcal{F})$ then $k=k'$. We will name those connected components \emph{circle or line $k$-prong singularities of $\mathcal{F}$}
        \item a singular leaf of $\mathcal{F}$ can a priori contain more than one connected components in $\text{Sing}(\mathcal{F})$
        \item outside the set of singularities, the singular foliation $\mathcal{F}$ defines a $C^0$ (regular) foliation on $M-\underset{C\in \text{Sing}(\mathcal{F})}{\cup}C$. We will say that $\mathcal{F}$ is \emph{orientable} (resp. \emph{transversely orientable}) if its restriction on $\Sigma-\text{Sing}(\mathcal{F})$ is orientable (resp. transversely orientable)
    \end{itemize}
\end{remadef}
\begin{defi}\label{d.transversedim3}
    Take $\mathcal{F},\mathcal{G}$ two singular codimension one foliations on $M^3$ with $\partial M=\emptyset$. We will say that $\mathcal{F}$ and $\mathcal{G}$ are \emph{transverse} if for every point $x\in M$ there exists $U_x$ a neighborhood of $x$ in $M$, $k\geq 1$ and $h:U_x\rightarrow \mathcal{D}_k\times [0,1]$ a homeomorphism such that 
        \begin{enumerate}
            \item $h(x)=(0, \frac{1}{2})\in \mathbb{C}\times [0,1]$
            \item $h(\mathcal{F}\cap U_x)=(\mathcal{F}_h^k\cap\mathcal{D}_k) \times [0,1]$ and $h(\mathcal{G}\cap U_x)=(\mathcal{F}_v^k\cap \mathcal{D}_k)\times [0,1]$
        \end{enumerate} 
\end{defi}

\begin{defi}[Topological pseudo-Anosov flow]\label{d.pseudoanosovflow}
  Let $M$ be a closed, smooth $3$-manifold, $d$ a distance on $M$ given by some Riemannian metric and $(X^t)_{t\in \mathbb{R}}$ a $C^0$-flow on $M$ (i.e. a one parameter family of homeomorphisms verifying $X^{t+s}=X^t\circ X^s$ and $X^0=id_M$). We will say that $(X^t)_{t\in \mathbb{R}}$ is a \emph{topological pseudo-Anosov} flow if: 
  \begin{enumerate}
      \item $(X^t)_{t\in \mathbb{R}}$ is non-singular (i.e. for every point $x$ in $M$ there exists $t\in \mathbb{R}$ such that $X^t$ does not fix $x$)
      \item $(X^t)_{t\in \mathbb{R}}$ preserves every leaf of a couple of transverse singular codimension one foliations $F^s$, $F^u$ with no circle $1$-prong singularities, called respectively the \emph{stable} and \emph{unstable foliations} of $(X^t)_{t\in \mathbb{R}}$
      \item for every $x\in M,y\in F^s(x)$ (resp. $y\in F^u(x)$) there exists an increasing homeomorphism $h:\mathbb{R}\rightarrow \mathbb{R}$ such that $$d(X^t(x),X^{h(t)}(y))\underset{t\rightarrow +\infty}{\longrightarrow} 0~~ \big(\text{resp. }d(X^t(x),X^{h(t)}(y))\underset{t\rightarrow -\infty}{\longrightarrow}0\big) $$

      \item there exist $\eta, \epsilon>0$ such that for any $x\in M$, $y\in F^s(x)$, $z\in F^u(x)$ and any increasing  homeomorphism $h:\mathbb{R}\rightarrow \mathbb{R}$ with $h(0)=0$ $$ \forall t\leq 0~ d(X^t(x),X^{h(t)}(y))<\eta \implies \exists |t_0|<\epsilon ~~ y=X^{t_0}(x)$$
      $$ \forall t\geq 0~ d(X^t(x),X^{h(t)}(z))<\eta \implies \exists |t_0|<\epsilon ~~z=X^{t_0}(x)$$
  \end{enumerate}
  If furthermore $\text{Sing}(F^s)=\text{Sing}(F^u)=\emptyset$, we will call $(X^t)_{t\in \mathbb{R}}$ a \emph{topological Anosov flow}. Finally, consider the following modification of  property (2): 
  \begin{enumerate} \renewcommand{\labelenumi}{(\arabic{enumi}$'$)}
  \setcounter{enumi}{1}
      \item $(X^t)_{t\in \mathbb{R}}$ preserves a pair of transverse singular codimension one foliations $F^s$, $F^u$ containing at least one circle $1$-prong singularity 
  \end{enumerate}
  A flow satisfying the properties (1), (2$'$) and (3) will be called a \emph{pseudo-Anosov flow with $1$-prong singularities}. 
\end{defi}
Let us point out that, as $M$ is closed and $d$ is given by a Riemannian metric on $M$, the previous definition does not depend on the choice of the distance $d$. Also, any pseudo-Anosov flow with $1$-prong singularities can not satisfy property (4); this is why we excluded property (4) from our definition of pseudo-Anosov flow with $1$-prong singularities.

As in the case of pseudo-Anosov homeomorphisms,  several definitions of the notion of pseudo-Anosov flow exist in the literature. One of the most important ones is the one of a \emph{smooth pseudo-Anosov flow} that we will now begin describing.  

Fix $\lambda_h,\lambda_v>1$. Denote by $\phi_2^{\lambda_h,\lambda_v}, \phi_2^{-\lambda_h,\lambda_v}, \phi_2^{\lambda_h,-\lambda_v}:\mathbb{R}^2\rightarrow \mathbb{R}^2$ the maps defined by $$\phi_2^{\lambda_h,\lambda_v}(x,y)= (\lambda_h^{-1}x, \lambda_v y)$$ 
$$\phi_2^{-\lambda_h,\lambda_v}(x,y)= (-\lambda_h^{-1}x, \lambda_v y)$$
$$\phi_2^{\lambda_h,-\lambda_v}(x,y)= (\lambda_h^{-1}x, -\lambda_v y)$$

Notice that the three previous maps commute with $-id$. We thus have that $\phi_2^{\lambda_h,\lambda_v},  \phi_2^{-\lambda_h,\lambda_v}, \phi_2^{\lambda_h,-\lambda_v}$ project to three uniquely defined  homeomorphisms  $\phi_1^{\lambda_h,\lambda_v},   \phi_1^{-\lambda_h,\lambda_v}, \phi_1^{\lambda_h,-\lambda_v}:\mathbb{R}^2\rightarrow \mathbb{R}^2$ respectively satisfying $$\pi_2\circ \phi_2^{\lambda_h,\lambda_v}=\phi_1^{\lambda_h,\lambda_v}\circ\pi_2$$
$$\pi_2\circ \phi_2^{-\lambda_h,\lambda_v}=\phi_1^{-\lambda_h,\lambda_v}\circ\pi_2$$
$$\pi_2\circ \phi_2^{\lambda_h,-\lambda_v}=\phi_1^{\lambda_h,-\lambda_v}\circ\pi_2$$

Notice that since $ \phi_2^{-\lambda_h,\lambda_v}=- \phi_2^{\lambda_h,-\lambda_v}$, we have that $\phi_1^{-\lambda_h,\lambda_v}= \phi_1^{\lambda_h,-\lambda_v}$.

For any $p\geq 3$, denote by $\phi_p^{\lambda_h,\lambda_v}:\mathbb{R}^2\rightarrow \mathbb{R}^2$ the unique homeomorphism preserving the pair of foliations $\mathcal{F}^p_{h,v}$, fixing the $\mathcal{F}^p_{h,v}$-prongs of the origin and satisfying $$\pi_p\circ \phi_p^{\lambda_h,\lambda_v}=\phi_1^{\lambda_h,\lambda_v}\circ\pi_p$$ Similarly, for any $p\geq 3$, denote by $\phi_p^{-\lambda_h,\lambda_v},  \phi_p^{-\lambda_h,-\lambda_v}:\mathbb{R}^2\rightarrow \mathbb{R}^2$ (resp. $:\mathbb{R}^2\rightarrow \mathbb{R}^2$) any homeomorphism preserving the pair of foliations $\mathcal{F}^p_{h,v}$, fixing an $\mathcal{F}^p_{v}$-prong (resp. $\mathcal{F}^p_{h}$-prong) of the origin and satisfying $$\pi_p\circ \phi_p^{-\lambda_h,\lambda_v}=\phi_1^{-\lambda_h,\lambda_v}\circ\pi_p 
\text{ (resp. }\pi_p\circ \phi_p^{\lambda_h,-\lambda_v}=\phi_1^{\lambda_h,-\lambda_v}\circ\pi_p)$$

It is not difficult to see that:
\begin{itemize}
     \item the homeomorphisms  $\phi^{\lambda_h,\lambda_v}_p,\phi^{-\lambda_h,\lambda_v}_p, \phi^{\lambda_h,-\lambda_v}_p$ contract (resp. stretch) each leaf of $\mathcal{F}^p_h$ (resp. $\mathcal{F}^p_v$) endowed with the measure $|\text{Re}(dz^{p/2})|$ (resp. $ |\text{Im}(dz^{p/2})|$ ) by a factor $\lambda_h^{-1}$ (resp. $\lambda_v$)
    \item contrary to $\phi_p^{\lambda_h,\lambda_v}$,  the homeomorphisms $\phi_p^{-\lambda_h,\lambda_v}, \phi_p^{\lambda_h,-\lambda_v}$ are not uniquely defined. However, despite their non uniqueness the previous homeomorphisms are completely determined by the $\mathcal{F}^p_{h}$-prong or $\mathcal{F}^p_{v}$-prong of the origin that they fix. It follows that any two choices of homeomorphisms satisfying our definition of $\phi_p^{-\lambda_h,\lambda_v}$ are conjugated by a rotation. Same for $\phi_p^{\lambda_h,-\lambda_v}$
   \item if $R_{\theta}(z)=e^{i2\pi\theta}z$, then for every $k\in \llbracket 0, p-1\rrbracket $ the map $R_{k/p}$ commutes with $\phi_p^{\lambda_h,\lambda_v}$
\end{itemize}
Let $k\in \llbracket 0, p-1\rrbracket$ and $l\in \{1,2\}$. We denote   

$$\phi_{p,k,+}^{\lambda_h,\lambda_v}:= \phi_p^{\lambda_h,\lambda_v}\circ R_{k/p}$$

$$\phi_{p,l,-}^{\lambda_h,\lambda_v}:=\left\{
		\begin{array}{ll}
			 \phi_p^{-\lambda_h,\lambda_v} & \mbox{when }  l=1 \\
             \phi_p^{\lambda_h,-\lambda_v}  & \mbox{when } l=2
			
		\end{array}
	\right.$$ Following the terminology of Mosher\footnote{Mosher defines the notion of local model and of smooth pseudo-Anosov flow in the case of orientable 3-manifolds. We generalize here Mosher's definition for non-orientable 3-manifolds.} (see \cite{Mosher3}), we will call $\phi_{p,k,+}^{\lambda_h,\lambda_v}$ the \emph{local model for a
pseudo-hyperbolic fixed point with stretching $\lambda^v$, compression $\lambda^h$, $p$ prongs,  rotation $k$ and positive orientation}. Similarly, we will call $\phi_{p,l,-}^{\lambda_h,\lambda_v}$ the \emph{local model of type $l$ for a
pseudo-hyperbolic fixed point, with stretching $\lambda^v$, compression $\lambda^h$, $p$ prongs and negative orientation}.

The previous maps are called models, as they model  the local behavior of any pseudo-Anosov homeomorphism around a fixed point. By suspending the previous local models, we will define a local model for a $p$-prong circle singularity. Take $o\in \{+,-\}$, $k\in \llbracket 0, p-1\rrbracket$ if $o=+$ and $k\in \{1,2\}$ if not. Consider the mapping torus
$$N_{p,k,o}:= \frac{\mathbb{R}^2\times \mathbb{R}}{((x,y),t+1) \sim (\phi_{p,k,o}^{\lambda_h,\lambda_v}(x,y),t) }$$
endowed with the constant speed vertical flow $(\Phi^t_{p,k,o})_{t\in\mathbb{R}}$ given by the vector field $\frac{\partial}{\partial t}$. When $o=+$ (resp. $o=-$), we will call $(N_{p,k,o},\Phi^t_{p,k,o})$ the \emph{local model for a pseudo-hyperbolic periodic orbit with stretching $\lambda^v$, compression $\lambda^h$, $p$ prongs, rotation (resp. type) $k$ and positive (resp. negative) orientation }. We will denote by $\gamma_{p,k,o}$ the periodic orbit of $\Phi^t_{p,k,o}$ defined by the suspension of the origin in $\mathbb{R}^2$ and by $F^s_{p,k,o},F^u_{p,k,o}$ the projection of $\mathcal{F}^p_h\times \mathbb{R}, \mathcal{F}^p_v\times \mathbb{R}$ on $N_{p,k,o}$. We will call $F^s_{p,k,o},F^u_{p,k,o}$ the \emph{local weak stable and unstable foliations of $(\Phi^t_{p,k,o})_{t\in\mathbb{R}}$}. 

Let us now define a metric on $N_{p,k,o}$. Since $z\rightarrow -z$ defines an isometry for the euclidiean metric $|\text{Re}(dz)|^2+|\text{Im}(dz)|^2$ and ${\pi_2}_{|\mathbb{R}^2-\{0\}}$ defines a cover of $\mathbb{R}^2-\{0\}$, $|\text{Re}(dz^{1/2})|^2+|\text{Im}(dz^{1/2})|^2$ defines a non-complete Riemannian metric on $\mathbb{R}^2-\{0\}$. Similarly, by lifting the previous metric, we have that $|\text{Re}(dz^{p/2})|^2+|\text{Im}(dz^{p/2})|^2$ defines also a non-complete Riemannian metric on $\mathbb{R}^2-\{0\}$ satisfying: 
$$ (\phi_{p,k,o}^{\lambda_h,\lambda_v})^{*}(|\text{Re}(dz^{p/2})|^2+|\text{Im}(dz^{p/2})|^2)= \lambda_h^{-2}\cdot|\text{Re}(dz^{p/2})|^2+\lambda_v^{2}\cdot|\text{Im}(dz^{p/2})|^2 $$

Notice that as the distance of any two points close to the origin is small for the above metric, $|\text{Re}(dz^{p/2})|^2+|\text{Im}(dz^{p/2})|^2$ defines a complete distance on $\mathbb{R}^2$ that we will denote by $d_p$. Consider now the Riemannian metric $\lambda_h^{-2t}\cdot|\text{Re}(dz^{p/2})|^2+\lambda_v^{2t}\cdot|\text{Im}(dz^{p/2})|^2+ |dt|^2$ on $\mathbb{R}^2\times \mathbb{R} - \{0\times \mathbb{R}\}$. Notice that the previous metric projects to a non-complete Riemannian metric on $N_{p,k,o}-\gamma_{p,k,o}$, which defines a complete distance $d_{p,k,o}$ on $N_{p,k,o}$. 
\begin{defi}[Mosher, \cite{Mosher3}]\label{d.smoothpseudoanosovflow}
Let $(X^t)_{t\in \mathbb{R}}$ be a non-singular flow on a closed, smooth $3$-manifold $M$.
We say that $(X^t)_{t\in \mathbb{R}}$ is a \emph{smooth pseudo-Anosov} flow if there exists a distance $d_M$ on $M$ such that the following are satisfied:
\begin{enumerate}
    \item There is a finite set $\Gamma$ of periodic orbits of $(X^t)_{t\in \mathbb{R}}$, called \emph{singular (periodic) orbits}, such that when restricted to $M-\Gamma$, the flow $(X^t)_{t\in \mathbb{R}}$ is smooth and $d_M$ corresponds to a distance given by a smooth Riemannian metric on $M-\Gamma$.
    \item Each orbit $\gamma$ in $ \Gamma$ is \emph{pseudo-hyperbolic}, i.e. there exist $\lambda_h,\lambda_v>1$, $p\geq 3$, $o\in \{+,-\}$, $k\in \llbracket 0, p-1 \rrbracket$ if $o=+$ or $k\in \{1,2\}$ if $o=-$, $U$ a neighborhood of $\gamma$ and $f_{\gamma}:U\rightarrow N_{p,k,o}$ an embedding such that: 
    \begin{itemize}
        \item $f_{\gamma}$ sends orbits of $(X^t)_{t\in \mathbb{R}}$ to orbits of $\Phi^t_{p,k,o}$ 
        \item $f_{\gamma}$ is smooth on $U-\gamma$
        \item $f_{\gamma}$ is bi-lipschitz with respect to the metrics $d_M$ and $d_{p,k,o}$
    \end{itemize}
    \item There exist two $dX^t$-invariant line bundles $E^s,E^u$ in $T(M-\Gamma)$ such that if $X$ is the vector field on $T(M-\Gamma)$ associated to $(X^t)_{t\in \mathbb{R}}$ and $\|\cdot\|$ the Riemmanian metric associated to $d_M$ on $M-\Gamma$ we have that: 
    \begin{itemize}
        \item $T(M-\Gamma)=E^s\oplus E^u\oplus \mathbb{R}X$
        \item there exists $C>0, \lambda\in (0,1)$ such that for every $u\in E^s$, $v\in E^u$ and every $t\geq 0$ we have $\|dX^t(u)\|\leq C\lambda^t\|u\|$ and $\|dX^{-t}(v)\|\leq C\lambda^t\|v\|$
    \end{itemize}
    \item For every $\gamma\in \Gamma$ the pre-image by $f_{\gamma}$ of the local weak stable and unstable foliations near $f_{\gamma}(\gamma)$ are tangent to $E^s$ and $E^u$ respectively
\end{enumerate}
Moreover, we will say that $(X^t)_{t\in \mathbb{R}}$ is a \emph{smooth Anosov flow} if $\Gamma=\emptyset$. 
\end{defi}
The relations between smooth and topological transitive pseudo-Anosov flows are nowadays well understood. The fact that any smooth Anosov flow is also topological, follows from several classical results of the stable manifold theory (see for instance \cite{stablemanifold}). The previous results can be adapted in order to also prove that smooth pseudo-Anosov flows are topological as well. Conversely, M.Shannon in his thesis \cite{Mariothese} proves that:
 \begin{defi}
     Let $M,N$ be any two manifolds. Two flows $(M,\Phi)$ and $(N,\Psi)$ are called \emph{orbitally equivalent} if there exists a homeomophism $h:M\rightarrow N$ sending positive orbits of $\Phi$ to positive orbits of $\Psi$. 
 \end{defi}
 \begin{theorem}\label{t.smoothenanosov}
     Any topological transitive Anosov flow is orbitally equivalent to a smooth Anosov flow.
 \end{theorem} 
The previous theorem was recently generalized by I.Agol and C.C. Tsang in \cite{Agol} 
  \begin{theorem}\label{t.smoothenpseudoanosov}
     Any topological transitive pseudo-Anosov flow is orbitally equivalent to a smooth pseudo-Anosov flow.
 \end{theorem} 
 
 Even though, thanks to Theorems \ref{t.smoothenanosov} and \ref{t.smoothenpseudoanosov}, ``topological implies smooth" is true in the transitive case, the smoothability of non-transitive Anosov or pseudo-Anosov flows still remains an open problem, the answer to which is unclear for the author. 

As the context of this paper is purely topological, from now on we will restrict ourselves to the study of topological pseudo-Anosov flows. 
\begin{conv*}
    Unless explicitly mentioned otherwise, from now on we will use the term pseudo-Anosov flow when referring to a topological pseudo-Anosov flow. 
\end{conv*}
\subsection{Expansive vs pseudo-Anosov flows: an equivalence}
\begin{defi}\label{d.expansiveflow}
    Let $M$ be a smooth, closed $3$-manifold endowed with a distance $d$ given by some Riemannian metric. Consider $(X^t)_{t\in \mathbb{R}}$ a non-singular $C^0$-flow on $M$. We will say that $(X^t)_{t\in \mathbb{R}}$ is \emph{expansive} if there exist $\epsilon, \eta>0$ such that for any $x,y\in M$ and any increasing  homeomorphism $h:\mathbb{R}\rightarrow \mathbb{R}$ with $h(0)=0$ $$ \forall t\in \mathbb{R}~ d(X^t(x),X^{h(t)}(y))<\eta \implies \exists |t_0|<\epsilon ~~ y=X^{t_0}(x)$$

\end{defi}

Exactly as in the case of pseudo-Anosov homeomorphisms, the notion of expansivity is intimately related with any pseudo-Anosov flow on a closed $3$-manifold. Indeed, thanks to a classical result in the theory of pseudo-Anosov flows: 
\begin{prop}\label{p.pseudoanosovflowimpliesexpansive}
    Any pseudo-Anosov flow $(X^t)_{t\in \mathbb{R}}$ on a closed, smooth manifold $M^3$ is expansive. 
\end{prop}
Conversely, we also have that  

\begin{theorem}[Inaba, Matsumoto, Oka]
    Every expansive flow on a closed 3-manifold is pseudo-Anosov. 
\end{theorem}

The above theorem can be deduced from Definition \ref{d.expansiveflow} and the main theorem of \cite{Inaba}, which relies essentially on the main theorem of \cite{Oka}. Thanks to the two above results, we get that a flow on a closed $3$-manifold is pseudo-Anosov if and only if it is expansive.

In order to better illustrate the close relationship between pseudo-Anosov flows and expansiveness, we will finish this section with the proof of Proposition \ref{p.pseudoanosovflowimpliesexpansive}. This proof relies on the existence of transverse standard polygons (see Definition \ref{d.standardpolygonflow}), on the understanding of the behavior of a pseudo-Anosov flow close to a periodic orbit (see Proposition \ref{p.aroundcircleprong}) -a result that will be used multiple times trhoughout this paper- and on two technical lemmas that will be provided here below.

\begin{defi}\label{d.standardpolygonflow}
 Consider $(X^t)_{t\in \mathbb{R}}$ a pseudo-Anosov flow on $M^3$ and $F^s,F^u$ its stable and unstable foliations. A \emph{transverse standard polygon} $P$ in $M$ is an embedded closed topological disk such that: 
 \begin{enumerate}
     \item there exists an embedded closed topological disk $U$, such that $P\subset \inte{U}$ and such that $U$ is \emph{transverse} to $(X^t)_{t\in \mathbb{R}}$, i.e there exists $\epsilon>0$ such that for every $x\in P$ the segment $\underset{t\in(-\epsilon, \epsilon)}{\cup}X^{t}(x)$ intersects $U$ once. 
     
    \noindent Let $F^s_U,F^u_U$ be the traces of the foliations $F^s,F^u$ on $U$.
     \item $P$ is a standard polygon for the surface $\inte{U}$ endowed with the pair of singular transverse foliations $(F^s_U,F^u_U) $. 
     
     \noindent If $P$ is a rectangle (resp. standard $2k$-gon) for the previous surface we will call $P$ a  \emph{transverse rectangle} (resp. \emph{transverse $2k$-gon}) of $(X^t)_{t\in \mathbb{R}}$ (see Figure \ref{f.polygons})

\end{enumerate} 
Furthermore, a connected component of $F^s\cap P$ (resp. $F^u\cap P$) will be called a \emph{stable} (resp. \emph{unstable}) \emph{leaf} of $P$. For any $x\in P$, we will denote by $P^s(x)$ (resp. $P^u(x)$) stable (resp. unstable) leaf of $P$ containing $x$ and we will call the closure of each connected component of $P-(P^s(x)\cup P^u(x))$ a \emph{quadrant of $x$ in $P$}.
\end{defi}

 Consider $(X^t)_{t\in \mathbb{R}}$ a pseudo-Anosov flow on $M^3$ and $F^s, F^u$ its stable and unstable foliations. Endow $M$ with a distance $d$ given by some Riemannian metric. By Theorem 17A of \cite{Whitney} and the fact that $(X^t)_{t\in \mathbb{R}}$ preserves a pair of transverse singular
foliations, it is easy to see that any $x\in M$ is contained in the interior of a transverse
standard polygon. Thanks to this fact and to the form of the foliations around a periodic orbit of $(X^t)_{t\in \mathbb{R}}$ it is not hard to prove that:
\begin{prop}\label{p.aroundcircleprong}
    Consider $(N_{p,r,o},\Phi_{p,r,o})$ the local model for a pseudo-hyperbolic periodic orbit $\gamma_{p,r,o}$ with stretching $\lambda^u$, compression $\lambda^h$, $p\geq 2$ prongs,  rotation $r \in \llbracket 0, p-1 \rrbracket$ and positive orientation (or type $r\in \{1,2\}$ and negative orientation). Denote by $F^s_{p,r,o}$ and $F^u_{p,r,o}$ the weak stable and unstable foliations of $\Phi_{p,r,o}$.
    
    For every periodic orbit $\gamma$ of $(X^t)_{t\in \mathbb{R}}$, there exist $p\geq 2$, $o\in \{+,-\}$, $r \in \llbracket 0, p-1 \rrbracket$ if $o=+$ (or $r\in \{1,2\}$ if $o=-$), $V_{p,r,o}$ a neighborhood of $\gamma_{p,r,o}$ in $N_{p,r,o}$, $U_\gamma$ a neighborhood of $\gamma$ in $M$ and a homeomorphism $H:U_\gamma\rightarrow V_{p,r,o}$ such that 
    \begin{itemize}
        \item $H$ defines an orbital equivalence between the restriction of  $(X^t)_{t\in \mathbb{R}}$ on $U_{\gamma}$ and the restriction of $\Phi_{p,r,o}$ on $V_{p,r,o}$ 
        \item for any $x\in U_{\gamma}$ the map $H$ takes the connected component of $F^s(x)\cap U_\gamma$ (resp. $F^u(x)\cap U_\gamma$) containing $x$ to the connected component of $F^s_{p,r,o}(H(x))\cap V_{p,r,o}$ (resp. $F^u_{p,r,o}(H(x))\cap V_{p,r,o}$) containing $H(x)$  
    \end{itemize}
\end{prop}
Fix a neighborhood $U_\gamma$ given by the above proposition for every circle prong singularity $\gamma$ of $(X^t)_{t\in \mathbb{R}}$. Using once again the existence of transverse standard polygons and the compactness of $M$, it is possible to
find for any sufficiently small $\beta > 0$, a finite number of transverse standard polygons $R_1, S_1, R_2, ..., R_n, S_n$ such that: 
\begin{itemize}
\item the $S_i$ are pairwise disjoint
\item $R_i\subset \inte{S_i}$, where $\inte{S_i}$ denotes the interior of $S_i$ 
\item $R_i$ and $S_i$ have the same number of stable (or unstable) boundary components for every $i\in I$
\item if $\inte{S_i}$ intersects a circle prong singularity $\gamma$, we have that $S_i\subset \inte{U_\gamma}$
\item every orbit of size $4\beta$ intersects at most once every $S_i$; hence the flow $(X^t)_{t\in \mathbb{R}}$ inside $\underset{t\in(-2\beta,2\beta)}{\cup}X^t(S_i)$ is conjugated to a constant speed vertical flow 
\item $\underset{i\in I}{\cup}\underset{t\in[-\beta,0]}{\cup}X^t(R_i)=\underset{i\in I}{\cup}\underset{t\in[0,\beta]}{\cup}X^t(R_i)=M$. In other words, the positive and negative orbit of any point in $M$ intersects the union $\underset{i\in I}{\cup}R_i$ in a time at most equal to $\beta$
\end{itemize}

Denote by $p^+:M \rightarrow \underset{i\in \llbracket 1,n \rrbracket}{\cup} R_i$ (resp. $p^-:M \rightarrow \underset{i\in \llbracket 1,n \rrbracket}{\cup} R_i$) the map associating to every point $x\in M$ the first point of intersection of the (strictly) positive orbit of $x$ (resp. negative orbit of $x$) with $\underset{i\in \llbracket 1,n \rrbracket}{\cup} R_i$. Denote also by $\pi_i:\underset{t\in(-2\beta,2\beta)}{\cup}X^t(S_i)\rightarrow S_i$ the projection map on $S_i$. By the above properties of the $R_i$ and thanks to the compactness of $M$, there exists $\delta>0$ such that for every $x\in M$ if $p^+(x)\in R_i$ or $p^-(x)\in R_i$, then $\pi_i$ is well defined on the ball $B(x,\delta):=\{y\in M|d(x,y)<\delta\}$. If $p^+(x)\in R_i$ (resp. $p^-(x)\in R_i$) denote by $p_x(y)$ (resp. $p^{-1}_x(y)$) the point $\pi_i(y)$ for every $y\in B(x,\delta)$. If for some $x\in M$ and $y\in B(x,\delta)$, the points $p^+(x)$ and $p_x(y)$ are $\delta$-close, we will denote by $p_x^2(y)$ the point $p_{p(x)}(p_x(y))$. We similarly define (when this is possible) $p_x^k(y)$ for every $k\in \mathbb{Z}$.
\begin{defi}\label{d.localstable}
    Consider $\epsilon>0$. For any $x\in M$ we define the \emph{$\epsilon$-local stable set of $x$} as $$F^s_\epsilon(x):=\{y\in M| \exists ~ h\in \text{Homeo}^+(\mathbb{R}) \text{ with } h(0)=0 \text{ such that } \forall t>0 ~ d(X^t(x),X^{h(t)}(y))<\epsilon\}$$
    We similarly define the \emph{$\epsilon$-local unstable set of $x$}, $F^u_{\epsilon}(x)$.
\end{defi}

Let us remark here that in the case of smooth pseudo-Anosov flows, the  $\epsilon$-local stable set of a point contains a neighborhood of this point inside its stable leaf (endowed with the leaf topology). A similar result holds for topological pseudo-Anosov flows, however this is a priori a non-trivial result. A proof of the previous fact, that uses the expansivity of a topological pseudo-Anosov flow, can be found in \cite{Inaba} (see Theorem 1.5) or \cite{Oka} (see Lemma 2.7).

Following the previous definition, two points $x,y\in M$ belong in the same local stable set if the positive orbit of $x$ shadows the positive orbit of $y$ for an appropriate reparametrization of the orbit of $y$. Since the previous reparametrization can be arbitrary, it is a priori difficult to decide whether two orbits belong in the same local stable set. According to the following two lemmas, this can be easily deduced by examining the successive intersections of the positive orbits of $x$ and $y$ with $R_1,S_1,...,R_n,S_n$. In order to simplify our notations, from now on $p^k$ will denote the $k$-th return map $(p^+)^k$ when $k> 0$ and $(p^-)^{|k|}$ when $k< 0$.   
\begin{lemm}\label{l.stablesetalternativedefinition}
     For every $\eta\in (0,\delta)$ there exists $\epsilon\in (0,\delta)$ such that for every $x\in M$, every $y\in F^s_{\epsilon}(x)$, and every $k\in \mathbb{N}^*$, $p_x^k(y)$ is well defined and $d(p_x^k(y),p^k(x))<\eta$. Furthermore, if $h\in \text{Homeo}^+(\mathbb{R})$ verifies $h(0)=0$ and $d(X^t(x),X^{h(t)}(y))<\epsilon$ for every $t\geq0$, then $p_x^k(y)=\pi_i(X^{h(s_k)}(y))$, where $p^k(x)=X^{s_k}(x)\in R_i$ and $i\in I$.
\end{lemm}
\begin{proof}
    Fix $\epsilon\in (0,\delta)$ sufficiently small so that: 
    \begin{enumerate}
        \item for every $i\in I$ and every $z,z'\in \underset{t\in(-2\beta,2\beta)}{\cup}X^t(S_i)$ with $d(z,z')<\epsilon$ we have $d(\pi_i(z),\pi_i(z'))<\eta$
        \item for every $i\in I$, the $\epsilon$-neighborhood of $\underset{t\in[-\beta,\beta]}{\cup}X^t(R_i)$ is contained in the interior of  $\underset{t\in[-2\beta,2\beta]}{\cup}X^t(S_i)$ 
    \end{enumerate}
    Consider $x\in M$, $y\in F^s_{\epsilon}(x)$ and $h\in \text{Homeo}^+(\mathbb{R})$ with $h(0)=0$ such that 
    \begin{equation}\label{eq.inequality}\tag{$\star$}
        d(X^t(x),X^{h(t)}(y))<\epsilon 
    \end{equation}

    Assume that $p(x)=X^{t_1}(x)\in R_i$, where $t_1\in (0,\beta]$ and $i\in I$. By (\ref{eq.inequality}), we have that the orbit segment $\underset{t\in[0,h(t_{1})]}{\cup}X^t(y)$ is contained in the interior of $\underset{t\in[-2\beta,2\beta]}{\cup}X^t(S_i)$; hence, $p_x(y)=\pi_i(y)=\pi_i(X^{h(t_1)}(y))$. It follows by our choice of $\epsilon$ that $d(p(x),p_x(y))=d(\pi_i(x),\pi_i(y))<\eta$. 

   Since $\eta<\delta$ we get that $p_x^2(y)$ is well defined. More precisely, assume that $p^2(x)=X^{t_2}(p(x))=X^{t_1+t_2}(x)\in R_j$, where $t_2\in (0,\beta]$ and $i\in I$. By (\ref{eq.inequality}), we have that the orbit segment $\underset{t\in[h(t_{1}), h(t_1+t_{2})]}{\cup}X^t(y)$ is contained in the interior of $\underset{t\in[-2\beta,2\beta]}{\cup}X^t(S_j)$; hence, $p^2_x(y)=p_{p(x)}(p_x(y))=p_{p(x)}(X^{h(t_1)}(y))=\pi_j(X^{h(t_1)}(y))$. It follows by our choice of $\epsilon$ that $d(p^2(x),p^2_x(y))=d\big(\pi_j(p(x)),\pi_j(X^{h(t_1)}(y))\big)<\eta$. We get the desired result by induction. 
\end{proof}
Conversely, using the continuity of $(X^t)_{t\in \mathbb{R}}$ and the fact that each orbit of $(X^t)_{t\in \mathbb{R}}$ intersects the $R_1,...,R_n$ in a time uniformly bounded by $\beta$ we have that:
\begin{lemm}\label{l.stablealternativedefinition2}
    For any $\epsilon>0$ sufficiently small, there exists $\eta\in (0,\delta)$ such that for every $x\in M$ and $y\in B(x,\eta)$ if $p_x^k(y)$ is well defined and $d(p_x^k(y),p^k(x))<\eta$ for every $k\in \mathbb{N}^*$, then $x\in F^s_\epsilon(y)$. 
\end{lemm}

We are now ready to begin the proof of Proposition \ref{p.pseudoanosovflowimpliesexpansive}. The following proof  is an adaptation of our proof of Proposition \ref{p.pseudoanosovimpliesexpansive}.
\begin{proof}[Proof of Proposition \ref{p.pseudoanosovflowimpliesexpansive}]
    Throughout this proof, we will follow the previously introduced notations: $R_1,S_1,...,R_n,S_n, \beta, \delta, p, p_x, U_\gamma$. Suppose that there exist two positive sequences $\eta_m \longrightarrow 0$, $\epsilon_m\longrightarrow +\infty$, a sequence of increasing homeomorphisms $h_m:\mathbb{R}\rightarrow \mathbb{R}$ with $h_m(0)=0$ and two sequences of points $x_m, y_m\in M$ such that for all $t\in \mathbb{R}$
    
    \begin{equation}\label{eq.shadowing1}
        d(X^t(x_m),X^{h_m(t)}(y_m))<\eta_m
    \end{equation}  
    
    \begin{equation}\label{eq.notinthesameorbit}
        y_m\notin \underset{s\in[-\epsilon_m,\epsilon_m]}{\cup}X^s(x_m)
    \end{equation}

\textbf{Case 1:} There exist infinitely many  $m\in\mathbb{N}$ such that $x_m$ is contained in a circle prong singularity $\gamma$. 

Consider $(N_{p,r,o},\Phi^t_{p,r,o})$ the local model for a pseudo-hyperbolic periodic orbit $\gamma_{p,r,o}$ with stretching $\lambda^v$, compression $\lambda^h$, $p$ prongs, rotation $r$ (resp. type $r$) and positive (resp. negative) orientation. By definition of $U_\gamma$, there exists $o\in \{+,-\}$, $p\geq 3$, $r\in \llbracket 0, p-1\rrbracket$ if $o=+$ or $r\in \{1,2\}$ if not, such that the flow $(X^t)_{t\in \mathbb{R}}$ inside $U_\gamma$ is orbitally equivalent to the flow $\Phi^t_{p,r,o}$ inside a compact neighborhood of $\gamma_{p,r,o}$. Denote by $H$ the previous orbital equivalence. Recall that in the previous section we endowed $N_{p,r,o}$ with a complete distance $d_{p,r,o}$. Using $H$ and our initial hypothesis, we get that there exists a positive sequence $\widetilde{\eta}_m\underset{m\rightarrow +\infty}{\longrightarrow} 0$, two sequences of points $\widetilde{x_m}\in \gamma_{p,r,o}$, $\widetilde{y_m}\in N_{p,r,o}$, a sequence of increasing homeomorphisms $\widetilde{h_m}:\mathbb{R}\rightarrow \mathbb{R}$ with $\widetilde{h_m}(0)=0$ such that $d_{p,r,o}(\Phi_{p,r,o}^t(\widetilde{x_m}),\Phi_{p,r,o}^{\widetilde{h_m}(t)}(\widetilde{y_m}))<\widetilde{\eta_m}$ for all $t\in \mathbb{R}$. One can easily verify that the previous fact implies that there exists $\widetilde{\epsilon}$ such that for $m$ sufficiently big $\widetilde{y_m}\in \underset{t\in[-\widetilde{\epsilon},\widetilde{\epsilon}]}{\cup}\Phi_{p,r,o}^t(\widetilde{x_m})$, which contradicts our initial hypothesis. 

\vspace{0.5cm}
It follows that Case 1 is impossible. From now on, by eventually considering a subsequence, we will assume  that $x_m$ is not contained in a circle prong singularity of $(X^t)_{t\in\mathbb{R}}$. By Lemma \ref{l.stablesetalternativedefinition} applied to $(X^t)_{t\geq0}$ and $(X^t)_{t\leq0}$, we get that  for $m$ sufficiently big, $p^k_{x_m}(y_m)$ is well defined for every $k\in \mathbb{Z}$ and that there exists a positive sequence $l_m\underset{m\rightarrow +\infty}{\longrightarrow}0$ such that for every $k\in \mathbb{Z}$
\begin{equation}\label{eq.closeiterates}
    d(p^k(x_m),p^k_{x_m}(y_m))<l_m
\end{equation} 

By eventually considering a subsequence assume that this is the case for every $m\in \mathbb{N}$. By an argument similar to the one we used in Case 1, one can show that the following case is also impossible: 

\textbf{Case 2:} There exist infinitely many $m\in\mathbb{N}$ for each one of which there exist  $k\in \mathbb{Z}^*$, $i\in I$, $\gamma$ a circle prong singularity of $(X^t)_{t\in\mathbb{R}}$ intersecting $R_i$ and $Q,Q'$ two distinct  quadrants of $\gamma\cap R_i$ in $R_i$ such that $p^k(x_m)\in Q$, $p^k_{x_m}(y_m)\in Q'$ and $\{p^k(x_m),p^k_{x_m}(y_m)\}\not\subset Q\cap Q'$

\vspace{0.5cm}
By eventually considering a subsequence, from now on we will assume that the situations described in Cases 1 and 2 do not arise for any $m\in \mathbb{N}$. We thus have that:

\textbf{Final case:} For every $m\in\mathbb{N}$ sufficiently big and every $k\in \mathbb{Z}^{*}$  either there exists a rectangle $R_{i(k,m)}$ with $i(k,m)\in I$ such that $p^k(x_m),p^k_{x_m}(y_m)\in R_{i(k,m)}$ or there exists a standard transverse $2p$-gon $R_{i(k,m)}$ with $i(k,m)\in I$ containing a singular point $s_{k,m}$ such that $p^k(x_m),p^k_{x_m}(y_m)$ belong in the same quadrant of $s_{k,m}$ in $R_{i(k,m)}$. 

Assume without any loss of generality that the above is true for every $m\in\mathbb{N}$. By eventually changing slightly our choice of $x_m$, we will assume that $x_m\notin \underset{i\in I}{\cup}R_i$ for every $m\in\mathbb{N}$.  

For any two points $z,z'$ in the same rectangle $R_i$ or in the same quadrant of a $2p$-gon $R_i$ denote by $[z,z']$ the unique point in $R^s_i(z)\cap R^u_i(z')$ (see Definition \ref{d.standardpolygonflow} for the previous notation). For every $k\in \mathbb{Z}^{*}$, denote by $z_{m,k}$ the point $[p^k(x_m),p^k_{x_m}(y_m)]\in R_{i(k,m)}$ and by $I^s_{m,k}$ (resp. $I^u_{m,k})$ the stable (resp. unstable) segment of $R_{i(k,m)}$ going from $p^k(x_m)$ (resp. $p^k_{x_m}(y_m)$) to $z_{m,k}$ (see Figure \ref{f.proofexpansiveflows}). We will prove that for $m$ sufficiently big we have that : 
\begin{itemize}
    \item for any $k\in -\mathbb{N}^*$, $p^{-1}_{p^k(x)}(z_{m,k})$ is well defined and \begin{equation}\label{eq.fact1}
       p^{-1}_{p^k(x_m)}(z_{m,k})=z_{m,(k-1)} 
    \end{equation}
    \item for any $k\in \mathbb{N}^*$, $p_{p^k(x)}(z_{m,k})$ is well defined and \begin{equation}\label{eq.fact2}
        p_{p^k(x_m)}(z_{m,k})=z_{m,(k+1)}
    \end{equation}
    \item $p^{-1}_{p(x)}(z_{m,1})$ is well defined and \begin{equation}\label{eq.fact3}
        p^{-1}_{p(x)}(z_{m,1})=z_{m,-1}
    \end{equation}
    \item $p_{p^{-1}(x)}(z_{m,-1})$ is well defined and \begin{equation}\label{eq.fact4}
        p_{p^{-1}(x)}(z_{m,-1})=z_{m,1}
    \end{equation}
\end{itemize}

\begin{figure}
    \centering
    \includegraphics[scale=1.9]{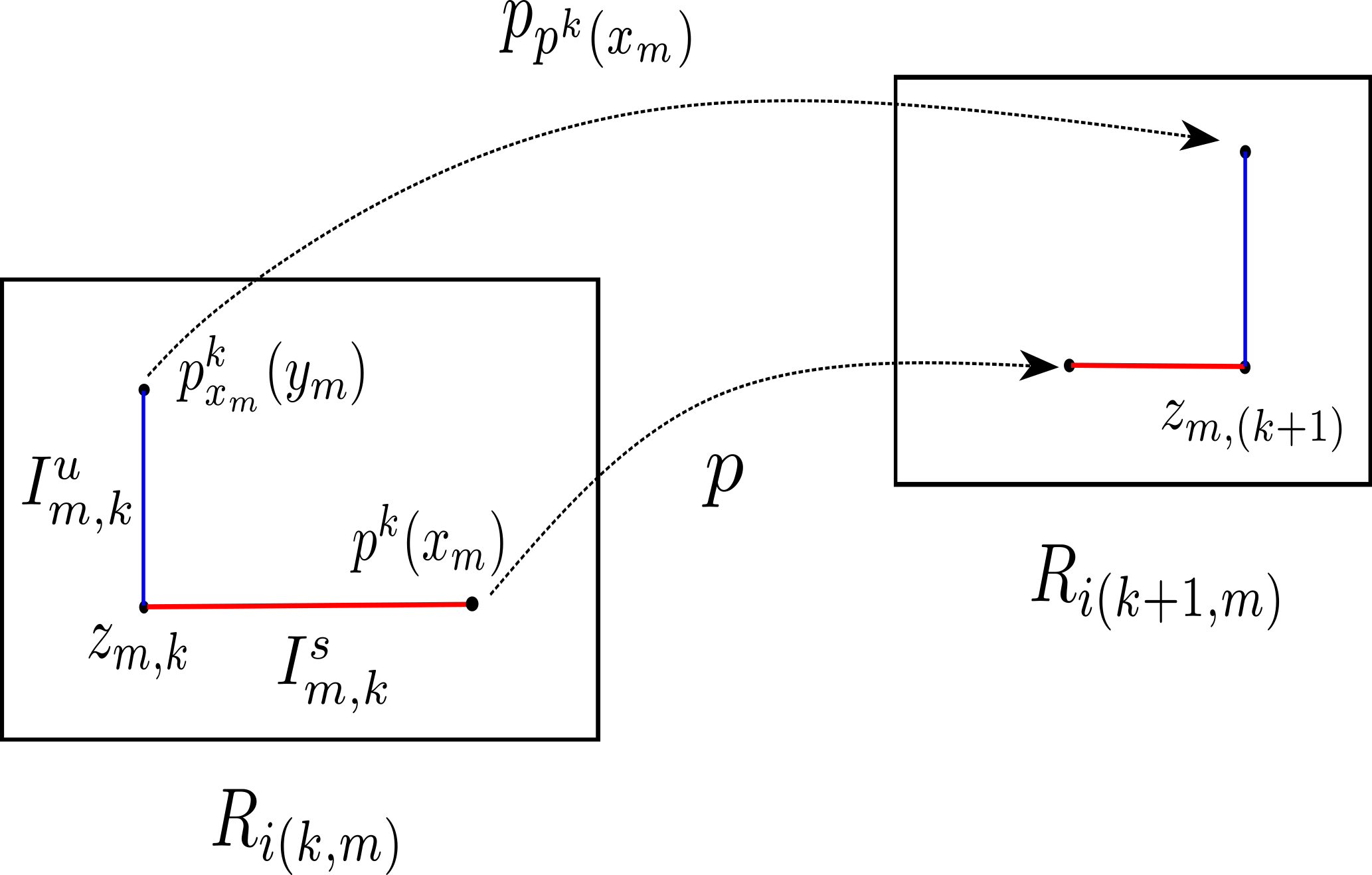}
    \caption{}
    \label{f.proofexpansiveflows}
\end{figure}
Let us prove \ref{eq.fact1} (a similar argument can be used to prove \ref{eq.fact2}, \ref{eq.fact3}, \ref{eq.fact4}). The function $[\cdot,\cdot]$ being continuous, thanks to Inequality \ref{eq.closeiterates}, there exists $\tau_m\underset{m\rightarrow +\infty}{\longrightarrow}0$ such that for every $k\in -\mathbb{N}^*$, $x\in I^s_{m,k}$ and $ y\in I^u_{m,k}$ we have 
\begin{equation}\label{eq.closeiterates2}
    d(p^k(x_m), [x,y])<\tau_m  \text{  and   }d(p^k_{x_m}(y_m), [x,y])<\tau_m
\end{equation} 

Take $m$ sufficiently big so that $\tau_m<\delta$. By our choice of $\delta$ (see our discussion prior to Definition \ref{d.localstable}), we have that $p^{-1}_{p^k(x)}$ is well defined on $I^s_{m,k}\cup I^u_{m,k}$ and since every stable and unstable segment in a transverse standard standard polygon intersect at most once, we also get that $p^{-1}_{p^k(x_m)}(I^{s,u}_{m,k})=I^{s,u}_{m,k-1}$ and $p^{-1}_{p^k(x_m)}(z_{m,k})=z_{m,(k-1)}$ for every $k\in -\mathbb{N}^*$. Our previous arguments imply that for $m$ sufficiently big and any $k\in -\mathbb{N}^*$, $p^k_{p^{-1}(x_m)}(z_{m,-1})=z_{m,k-1}$ and $$d(p^{k-1}(x_m),p^k_{p^{-1}(x_m)}(z_{m,-1}))=d(p^{k-1}(x_m),z_{m,k-1})<\tau_m$$ The previous inequality and Lemma \ref{l.stablealternativedefinition2} give that for $m$ sufficiently big there exists a positive sequence $\tau_m'\underset{m\rightarrow +\infty}{\longrightarrow}0$ such that $z_{m,-1}\in F^u_{\tau_m'}(p^{-1}(x_m))$. By proving and then using the facts \ref{eq.fact2}, \ref{eq.fact3}, \ref{eq.fact4} in the exact same way, we get that, for $m$ sufficiently big there exists a positive sequence $\tau_m''\underset{m\rightarrow +\infty}{\longrightarrow}0$ such that $z_{m,-1}\in F^s_{\tau_m''}(p^{-1}(x_m))$. It follows that for $m$ sufficiently big, there exists a sequence of increasing homeomorphisms $h'_m:\mathbb{R}\rightarrow \mathbb{R}$ with $h'_m(0)=0$ such that $$d(X^t(p^{-1}(x_m)),X^{h'_m(t)}(z_{m,-1}))<\max(\tau_m',\tau_m'')$$

Using the fact that the orbits of $y_m$ shadow the orbits of $x_m$, it is not hard to prove that for $m$ sufficiently big there exists a positive sequence $\tau_m'''\underset{m\rightarrow +\infty}{\longrightarrow}0$ and a sequence of increasing homeomorphisms $h''_m:\mathbb{R}\rightarrow \mathbb{R}$ with $h''_m(0)=0$ such that $$d(X^t(p^{-1}_{x_m}(y_m)),X^{h''_m(t)}(z_{m,-1}))<\tau_m'''$$

By our previous arguments, we have that for $m$ sufficiently big the orbit of $z_{m,-1}$ shadows the orbits of both $p^{-1}(x_m)$ and $p^{-1}_{x_m}(y_m)$ and also that by its construction $z_{m,-1}\in F^s(p^{-1}(x_m))\cap F^u(p^{-1}_{x_m}(y_m))$. Thanks to property (4) of Definition \ref{d.pseudoanosovflow} and to the previous facts, we get that there exists $\epsilon>0$ such that for any $m$ sufficiently big, $$z_{m,-1}\in \underset{t\in[-\epsilon,\epsilon]}{\cup}X^t(p^{-1}(x_m))=\underset{t\in[-\epsilon-\beta,\epsilon]}{\cup}X^t(x_m)$$ $$z_{m,-1}\in \underset{t\in[-\epsilon,\epsilon]}{\cup}X^t(p^{-1}_{x_m}(y_m))=\underset{t\in[-\epsilon-2\beta,\epsilon]}{\cup}X^t(y_m)$$ Finally, this implies that for $m$ sufficiently big $y_m\in \underset{t\in[-2\epsilon-\beta,2\epsilon+2\beta]}{\cup}X^t(x_m)$, which contradicts \ref{eq.notinthesameorbit} and gives us the desired result. 
\end{proof}

As in the case of pseudo-Anosov homeomorphisms, any pseudo-Anosov flow with circle $1$-prong singularities can not be expansive. This is one of the main reasons why 1-prongs are not allowed in our definition of pseudo-Anosov flows. Also, although property (4) of Definition \ref{d.pseudoanosovflow} is sometimes not included in the definition of a pseudo-Anosov flow, the reader can notice that it played a crucial role in the proof of Proposition \ref{p.pseudoanosovflowimpliesexpansive}. This naturally poses the following question : 

\begin{ques}
 Is a flow satisfying properties (1)-(3) of Definition \ref{d.pseudoanosovflow} necessarily expansive? 
\end{ques}

Very recently, the answer to this question was proven to be negative for non-transitive pseudo-Anosov flows (see Example 1.2.15 of \cite{planeapproach}). To this day, the author is unaware of an answer to the above question for transitive pseudo-Anosov flows. 

\subsection{Examples of topological pseudo-Anosov flows}

\subsubsection*{Smooth or topological Anosov flows} 
\quad
\vspace{0.5cm}

We refer the reader to \cite{Anosov}, \cite{Introanosov} or  \cite{thesisioannis} for examples of Anosov flows and an introduction to Anosov flow theory.

\subsubsection*{The suspension of a pseudo-Anosov homeomorphism}
\quad
\vspace{0.5cm}

Consider a pseudo-Anosov homeomorphism $f$ acting on a closed surface $\Sigma$ and the manifold: 
$$\mathbb{T}_f:= \frac{\Sigma \times [0,1]}{(f(y), 0) \sim (y, 1) }$$

 Let $(X^s)_{s\in\mathbb{R}}$ be the constant speed vertical flow given by the vector field $ \frac{\partial}{\partial t}$, where $t$ is the $[0,1]$ coordinate. 
 \begin{thdef}
     The flow $(X^s)_{s\in\mathbb{R}}$ is a pseudo-Anosov flow on $T_f$, called the \emph{suspension pseudo-Anosov flow of the homeomorphism $f$}. 
 \end{thdef}
 \begin{proof}
     Indeed, it is known that since $f$ is isotopic to a $C^{\infty}$ diffeomorphism, the manifold $T_f$ admits a smooth structure. Let $d$ be a distance on $T_f$ given by some Riemannian metric and $F^s_f,F^u_f$ be the stable and unstable foliations of $f$. 
     
    Clearly  $(X^s)_{s\in\mathbb{R}}$ satisfies property (1) of Definition \ref{d.pseudoanosovflow}. Next, since $F^s_f,F^u_f$ form a pair of $f$-invariant,  singular transverse foliations with no $1$-prong singularities, by pushing $F^s_f,F^u_f$ by  $(X^s)_{s\in\mathbb{R}}$, we obtain a pair of $X^s$-invariant singular transverse codimension one foliations $F^s,F^u$ with no circle $1$-prong singularities. Hence, $(X^s)_{s\in\mathbb{R}}$ satisfies property (2) of Definition \ref{d.pseudoanosovflow}. Finally, since $f$ (resp. $f^{-1}$) acts as a contraction on each leaf of $F^s_f$ (resp. $F^u_f$) and as an expansion on each leaf of $F^u_f$ (resp. $F^s_f$), we have that  $(X^s)_{s\in\mathbb{R}}$ satisfies properties (3) and (4) of Definition \ref{d.pseudoanosovflow}, which proves the desired result. 
 \end{proof}

\subsubsection*{Pseudo-Anosov flows in toroidal manifolds}
\quad
\vspace{0.3cm}

\begin{theorem}
    There exists a transitive pseudo-Anosov flow in dimension 3, whose every orbit except finitely many intersect a transverse torus.  
\end{theorem}

The above theorem is due to C.Bonatti and R.Langevin, who construct in \cite{BonattiLangevin} a transitive Anosov flow with the above properties. The previous example was next generalized by T.Barbot in \cite{Barbotbola} and gave rise to an infinite family of Anosov flows with similar properties, known as the family of \emph{generalized BL-flows}. Finally, by gluing simple dynamical blocks, T.Barbot and S.Fenley generalized the previous construction in \cite{BarbotFenley} for transitive and non-transitive pseudo-Anosov flows on graph manifolds giving rise to many new examples of Anosov and pseudo-Anosov flows in dimension 3. 

\subsubsection*{Pseudo-Anosov flows obtained by Dehn-Goodman-Fried surgery}
\quad
\vspace{0.5cm}

Starting from any pseudo-Anosov flow $\Phi$ in dimension 3, one can obtain infinitely many more pseudo-Anosov flows by performing Dehn-Goodman-Fried surgeries (see Section \ref{s.surgeries} for a defintion) on a finite set of periodic orbits of $\Phi$. A famous theorem by Fried (see \cite{Fried}) states that all transitive pseudo-Anosov flows can be obtained by performing a finite number of Dehn-Goodman-Fried surgeries on a suspension pseudo-Anosov flow. 

It is still not known whether there exists an analogue of Fried's theorem for non-transitive pseudo-Anosov flows. Another famous open conjecture in the fields of Anosov and pseudo-Anosov flows is the following one: 

\begin{conj}
 Given any two suspension pseudo-Anosov flows $\Phi_1$, $\Phi_2$, we can obtain $\Phi_2$ by performing a finite number of Dehn-Goodman-Fried surgeries on $\Phi_1$. 
\end{conj}

\subsection{The periodic orbits of a pseudo-Anosov flow}
As in the case of Anosov flows, every pseudo-Anosov flow has a rich dynamical behavior and many periodic orbits. For instance, as a direct consequence of Remark \ref{r.singularfolidim3} and our definition of pseudo-Anosov flow: 
\begin{prop}\label{p.firstperiodicorbits}
  Let $M$ be a smooth closed $3$-manifold, $(X^t)_{t\in\mathbb{R}}$ a pseudo-Anosov flow on $M$ and $F^s,F^u$ its stable and unstable foliations. Any singular leaf of $F^s,F^u$ contains a unique circle prong singularity forming a periodic orbit of $(X^t)_{t\in\mathbb{R}}$. 
\end{prop}

In fact, a much stronger result holds for a general pseudo-Anosov flow: 
\begin{prop}\label{p.periodicorbitsaredense}
    Let $(X^t)_{t\in\mathbb{R}}$ be a pseudo-Anosov flow on the closed manifold $M^3$ and $NW$ its non-wandering set. The set of periodic orbits of $(X^t)_{t\in\mathbb{R}}$ is dense inside $NW$. 
    
    In particular if $(X^t)_{t\in\mathbb{R}}$ is transitive (i.e. it admits a dense orbit), then the set of periodic orbits of $(X^t)_{t\in\mathbb{R}}$ is dense in $M$. 
\end{prop}
The above result was shown to be true for $C^1$ Anosov flows in \cite{Anosov} (see Lemma 13.1) and can be deduced for smooth pseudo-Anosov flows from the fact that a smooth pseudo-Anosov flow is semi-conjugated to a smooth structurally stable system (see \cite{Mosher1}). A proof of Proposition \ref{p.periodicorbitsaredense} in its general form can be found in \cite{markovpseudoanosov} (see also Theorem 1.8 and Remark 3.2 of \cite{nontransitiveanosovlike} for an alternative proof method).

The existence of periodic orbits together with the following theorem severely restrict the 3-manifolds that admit pseudo-Anosov flows. 
\begin{defi}
    Let $F$ be a singular codimension one foliation on $M$ with no circle $1$-prong singularities. Denote by $C_1,...,C_n$ the circle prong singularities of $F$. Take $x\in C_i$, $U$ a small open neighborhood of $x$ and $F_U(x)$ the connected component of $F(x)\cap U$ containing $x$. We will say that the closure of each connected component of $U-F_U(x)$ is a \emph{sector of $C_i$ in $U$}. 
    
    A loop $\gamma: \mathbb{S}^1\overset{C^0}{\rightarrow} M$ with no self-intersections will be called a \emph{closed quasi-transversal of $F$} if the following conditions are satisfied: 
    \begin{itemize}
         \item there exists a finite set $T\subset \mathbb{S}^1$ (possibly empty) such that for every $t\in T$, $\gamma(t)\in \underset{i\in \llbracket 1, n\rrbracket}{\cup}C_i$ 
        \item for every $t\in \mathbb{S}^1-T$ there exists a small neighborhood $J$ of $t$ such that $\gamma(J)\subset M-\underset{i\in \llbracket 1, n\rrbracket}{\cup}C_i$ is topologically transverse to $F$
       \item for every $t\in T$, if $\gamma(t)\in C_i$ and $V$ is a small neighborhood of $\gamma(t)$ in $M$, then for every small neighborhood $I$ of $t$ we have that $\gamma(I)$ is not contained in a unique sector of $C_i$ in $V$.  
    \end{itemize}
\end{defi}
\begin{theorem}\label{t.transversalhomotopicto0}
    Let $X^t$ be a pseudo-Anosov flow on $M^3$ and $F^s, F^u$ its stable and unstable foliations. Any closed quasi-transversal of $F^s$ or $F^u$ can not be freely homotopic to a trivial loop in $M$.
\end{theorem}

The above theorem constitutes a central result in the theory of smooth or topological pseudo-Anosov flows and admits several important corollaries. Historically, two different approaches have been used in order to prove Theorem \ref{t.transversalhomotopicto0} : the first approach relies on the one-sided holonomy theorem of A.Haefliger - a result that holds in arbitrary dimensions- and the second on the compact leaf theorem of S.P.Novikov - a result specific to the dimension 3-. 

More specifically, in the case of smooth Anosov flows with $C^2$ weak stable and unstable foliations, Theorem \ref{t.transversalhomotopicto0} can be deduced from Proposition 4.2 of \cite{Haefliger}, according to which if a $C^2$ foliation $\mathcal{F}$ of codimension one admits a homotopically trivial closed transversal, then there exists a closed loop inside a leaf of $\mathcal{F}$ that has trivial holonomy on one side and non-trivial on the other. More generally, in the case of a general smooth or topological Anosov flow, Theorem \ref{t.transversalhomotopicto0} can be deduced from the generalization of the one-sided holonomy theorem for $C^0$ foliations, proven by G.Hector and U.Hirsh in \cite{Hirsh} (see also Lemma 8 of \cite{Paternain}). Finally, Theorem \ref{t.transversalhomotopicto0} in its general form can be deduced from Lemmas 8 and 9 of \cite{Paternain}, which adapt the previous proof approach for $C^0$ singular foliations. 

Concerning the second family of proofs of Theorem \ref{t.transversalhomotopicto0}, S.P.Novikov proves in \cite{Novikov} that given a smooth foliation $\mathcal{F}$ on a closed $3$-manifold $V$ and $L$ a leaf of $\mathcal{F}$, if the inclusion map from $L$ to $V$ does not define a monomorphism on the level of the fundamental groups, then $\mathcal{F}$ admits a compact leaf. The previous theorem was generalized for $C^0$ foliations by V.V. Solodov (see \cite{Solodov}) and for essential laminations by D.Gabai and U.Oertel (see \cite{Gabai}). In the case of general smooth or topological Anosov flows, Theorem \ref{t.transversalhomotopicto0} can be deduced from the above result of V.V. Solodov (see also Theorem 1.5 of \cite{Inaba}). In the singular case,  S.Fenley and L.Mosher in \cite{Fenleymosher} utilize the above result of by D.Gabai and U.Oertel in order to prove Theorem \ref{t.transversalhomotopicto0} for smooth pseudo-Anosov flows. Although we are convinced that a similar approach applies to topological pseudo-Anosov flows, the author is not aware of a formal proof of this fact. 

\begin{coro}\label{c.periodicnontrivial}
     Consider $(X^t)_{t\in\mathbb{R}}$ a pseudo-Anosov flow on $M^3$ and $\gamma\in \pi_1(M)$ that is freely homotopic to a periodic orbit of $(X^t)_{t\in\mathbb{R}}$. For every $n\in \mathbb{Z}^*$, we have that $\gamma^n$ is non-trivial inside $\pi_1(M)$.  
 \end{coro}
 \begin{proof}
   By eventually considering the orientation cover of $M$, we will assume without any loss of generality that $M$ is orientable. Consider $\gamma$ a periodic orbit of $(X^t)_{t\in\mathbb{R}}$ and  $F^u$ the unstable foliation of $(X^t)_{t\in\mathbb{R}}$. By Theorem \ref{t.transversalhomotopicto0}, it suffices to prove that $\gamma^{K\cdot n}$ is freely homotopic to a quasi-transversal of $F^u$ for every $n\in \mathbb{N}^*$ and some $K\in \mathbb{N}$. 
   
   We restrict ourselves to a neighborhood of $\gamma$ for which we can apply Proposition \ref{p.aroundcircleprong}. Assume that inside this neighborhood the flow is orbitally equivalent to the flow defined by the local model for a pseudo-hyperbolic periodic orbit with $p$ prongs, rotation $r \in \llbracket 0,p - 1\rrbracket$ and positive orientation (the positivity of the orientation is implied by the orientability hypothesis on $M$). By the construction shown in Figure \ref{f.transversalperiodic} (for $p=3,r=1$) if $\text{gcd}(p,r)$ is the greatest common divisor of $p,r$ and $K=p/\text{gcd}(p,r)$, the loop $\gamma^{K\cdot n}$ is freely homotopic to a quasi-transversal of $F^u$ for every $n\in \mathbb{N}^*$. 
   \begin{figure}[h!]
       \centering
       \includegraphics[scale=0.4]{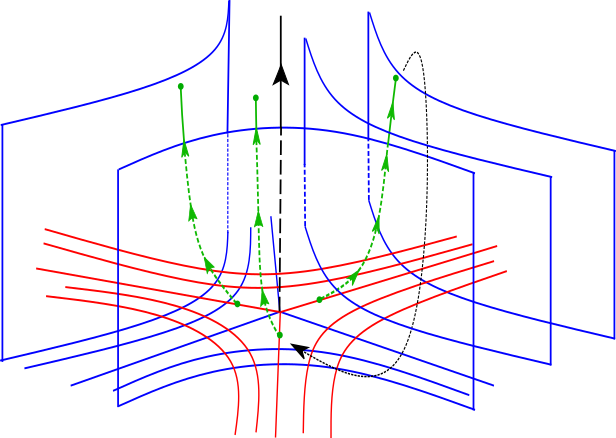}
       \caption{In the above figure the red and blue lines represent respectively weak stable and unstable leaves, the black line in the center corresponds to the periodic orbit $\gamma$ and the green line to the desired quasi-transversal of $F^u$}
       \label{f.transversalperiodic}
   \end{figure}
 \end{proof}
 As a result of the previous corollary and Proposition \ref{p.periodicorbitsaredense}, we get that: 
 \begin{coro}\label{c.infinitefund}
     For any $3$-manifold $M$ supporting a pseudo-Anosov flow, $\pi_1(M)$ is infinite.
 \end{coro}
 We should remark at this point that even if the fundamental group of a manifold sup-
porting a pseudo-Anosov flow needs to be “big”, the same result does not apply for the first homology group:
\begin{rema}
 There exist Anosov flows on 3-dimensional rational homology spheres. For instance, the geodesic flow of a genus 0 orbifold with 5 or more cone points is an  Anosov flow on a rational homology 3-sphere (see Lemma 2.1 in \cite{lefthand}).   
\end{rema}

\subsection{On the invariant foliations of a pseudo-Anosov flow}
Consider $\Phi=(\Phi^t)_{t\in \mathbb{R}}$ a pseudo-Anosov flow on a smooth $3$-manifold $M$ and $F^s,F^u$ the stable foliations of $\Phi$. 

Although it is not generally true that the set of periodic of $\Phi$ is dense in $M$, it turns out that the union of the stable (resp. unstable) manifolds of all the periodic orbits of $\Phi$ is dense in $M$. Even more, we have that 

\begin{prop}\label{p.densitystablemanifoldsperiodic}
    There exists a finite family $\Gamma$ of periodic orbits of $\Phi$ containing every singular periodic orbit of $\Phi$ and such that $\underset{\gamma \in \Gamma}{\cup}F^s(\gamma)$ is dense in $M$. 
\end{prop}
Naturally, the same proposition is true for $F^u$. The above result was shown to be true for smooth Anosov flows by D.V.Anosov (see Paragraph §14 of  \cite{Anosov}) and can be deduced for smooth pseudo-Anosov flows from the fact that a smooth pseudo-Anosov flow is semi-conjugated to a smooth structurally stable system (see \cite{Mosher1}). The fact that Proposition \ref{p.densitystablemanifoldsperiodic} holds also for topological pseudo-Anosov flows can be obtained by adapting the arguments used in Paragraph §14 of \cite{Anosov}. Another proof of this fact may be found in Proposition 7 of \cite{markovpseudoanosov}. 

In the case of transitive pseudo-Anosov flows (i.e. pseudo-Anosov flows admitting dense orbits) an even stronger result is true: 
\begin{prop}
    If $\Phi$ is transitive, then for every periodic orbit of $\Phi$, say $\gamma$, $F^s(\gamma)$ and $F^u(\gamma)$ are both dense in $M$. 
\end{prop}
The reader may find a proof of the above result in \cite{Plante} (see Theorem 1.3\footnote{Even though Theorem 1.3 in \cite{Plante} concerns $C^1$ Anosov flows, its proof can be applied in the exact same way in the cases of both smooth and topological pseudo-Anosov flows.}).

Recall that our definition of pseudo-Anosov flow restricts the dynamics of $\Phi$ inside any stable or unstable leaf. This imposes several restrictions on the topology of the leaves of $F^s$ and $F^u$:
\begin{defi}
    Let $\mathcal{O}$ be an orbit of $\Phi$ and $L$ its stable (resp. unstable) leaf. The union of $\mathcal{O}$ with a connected component of $L-\mathcal{O}$ (endowed with the leaf topology) will be called a \emph{stable} (resp. \emph{unstable})  {separatrix} of $\mathcal{O}$. 
\end{defi}
\begin{prop}\label{p.topologyleaves}
Fix $\gamma$ a periodic orbit of $\Phi$ and $S$ a stable (or unstable) separatrix of $\gamma$. 
   \begin{enumerate}
       \item  there exists a continuous immersion $\phi:\mathbb{S}^1\times [0,+\infty)\rightarrow M$ with the following properties: 
    \begin{itemize}
        \item $\phi_{|\mathbb{S}^1\times  (0,+\infty)}$ is a homeomorphism onto $S-\gamma$ (endowed with the leaf topology)
        \vspace{0.1cm}
        \item $\phi(\mathbb{S}^1\times \{0\})=\gamma$ and $\phi_{|\mathbb{S}^1\times \{0\}}$ defines a finite by one map 
    \end{itemize}
    \vspace{0.2cm}
    \item if $\gamma$ is non-singular, then its stable (or unstable) manifold endowed with the leaf topology is either homeomorphic to an open cylinder or an open M\"obius band. In both of the above cases, the fundamental group of $F^s(\gamma)$ (or $F^u(\gamma)$) is generated by the loop $\gamma$

   \end{enumerate}

\end{prop}

    \textit{Proof of (1).} Assume without any loss of generality that $S$ is a stable separatrix of $\gamma$. Recall that $S$ is invariant by $\Phi$ and that every $\Phi$-orbit in $S$ will accumulate to $\gamma$ in the future. Thanks to Proposition \ref{p.aroundcircleprong}, one can easily construct an immersion $\phi_0: \mathbb{S}^1\times [0,1]\rightarrow M$ such that 
    \begin{itemize}
        \item ${\phi_0}_{|\mathbb{S}^1\times  (0,1]}$ is a homeomorphism onto $N-\gamma$, where $N$ is a neighborhood of $\gamma$ in $S$ (with respect to the leaf topology)
        \item $\phi_0(\mathbb{S}^1\times \{0\})=\gamma$ and ${\phi_0}_{|\mathbb{S}^1\times \{0\}}$ defines a finite by one map 
        \item ${\phi_0}_{|\mathbb{S}^1\times \{1\}}$ is a circle transverse to the restriction of $\Phi$ on $S$
        \item for some $s>0$, the set ${\phi_0}_{|\mathbb{S}^1\times (s,1]}$ is a fundamental domain of the (forward) action by  $\Phi^1$, the time one map of $\Phi$, in ${\phi_0}(\mathbb{S}^1\times  [0,1])$
    \end{itemize}
 Moreover, by our previous construction, we have that for every $n, m\in\mathbb{Z}$ if $m\neq n$ then $$ \Phi^n({\phi_0}_{|\mathbb{S}^1\times (s,1]}) \cap \Phi^m({\phi_0}_{|\mathbb{S}^1\times (s,1]}) =
\emptyset$$ By therefore applying negatively $\Phi^1$, we can extend $\phi_0$ to an immersion $\phi$ with the desired properties. 

\textit{Proof of (2).} We will prove the desired result for $F^s(\gamma)$ (of course the same argument applies to $F^u(\gamma)$). Using Proposition \ref{p.aroundcircleprong}, we will consider the following four cases: 
\begin{enumerate}
    \item the dynamics of $\Phi$ around $\gamma$ is orbitally equivalent to the dynamics given by a local model for a pseudo-hyperbolic orbit with $2$ prongs, rotation $0$ and positive orientation
    \item the dynamics of $\Phi$ around $\gamma$ is orbitally equivalent to the dynamics given by a local model for a pseudo-hyperbolic orbit with $2$ prongs, rotation $1$ and positive orientation
     \item the dynamics of $\Phi$ around $\gamma$ is orbitally equivalent to the dynamics given by a local model for a pseudo-hyperbolic orbit with $2$ prongs, of type $1$ and negative orientation
    \item the dynamics of $\Phi$ around $\gamma$ is orbitally equivalent to the dynamics given by a local model for a pseudo-hyperbolic orbit with $2$ prongs, type $2$ and negative orientation
\end{enumerate}
We will prove Item (2) for the first two of the above cases. The last two cases follow from a similar argument. 

Following the notations of our proof of Item (1), in the first case, it is easy to see that ${\phi_0}_{|\mathbb{S}^1\times \{0\}}$ is injective, which implies, thanks to our previous construction, that $\phi$ is an injective immersion and even more a homeomorphism onto $S$ (endowed with the leaf topology). Furthermore, by construction, $\phi(\mathbb{S}^1\times \{0\})=\gamma$ and thus the fundamental group of $S$ is generated by $\gamma$. We obtain the desired result by remarking that $F^s(\gamma)$ is the union of two separatrices glued along the dynamically oriented loop $\gamma$.  

In the second case, it is easy to see that ${\phi_0}_{|\mathbb{S}^1\times \{0\}}$ is two to one, thus defining a double cover map from $\mathbb{S}^1\cong\mathbb{S}^1\times \{0\}$ to $\mathbb{S}^1\cong \gamma$. As $\phi_{|\mathbb{S}^1\times  (0,+\infty)}$ is a homeomorphism onto $S-\gamma$ (endowed with the leaf topology), it is not difficult to see that $\phi$ defines an immersion of $\mathbb{S}^1\times [0, +\infty)$ onto an open M\"obius band such that the fundamental group of $S$ is generated by $\gamma$. We obtain the desired result by remarking that in this case $F^s(\gamma)$ consists of a unique stable separatrix.

\subsection{Manifolds supporting pseudo-Anosov flows}

Besides the restriction imposed on the manifolds supporting pseudo-Anosov flows by Corollary \ref{c.infinitefund} and the existence of a pair of transverse singular foliations with the properties given by Proposition \ref{p.topologyleaves}, there exist more necessary conditions for a manifold to support pseudo-Anosov flows. 

By Theorem 1.7 of \cite{Inaba} we have that: 
\begin{theorem}\label{t.aspherical}
Any $3$-manifold $M$ supporting a pseudo-Anosov flow is aspherical.  
\end{theorem}
T.Inaba and S.Matsumoto prove the above theorem in the case where $M$ is orientable. However, using the fact that a finite lift of an expansive flow is expansive and the fact that a finite cover of a non-aspherical $3$-manifold is also non-aspherical, one easily obtains the above theorem. Using Theorem 5.1 of \cite{Scott} and Lemma 3.1 of \cite{Aspherical}, we also get that: 
\begin{coro}
    The universal cover of any $3$-manifold $M$ supporting a pseudo-Anosov flow is homeomorphic to $\mathbb{R}^3$. Hence, $M$ is irreducible. 
\end{coro}
\begin{coro}\label{c.torsionfree}
    The fundamental group of any $3$-manifold $M$ supporting a pseudo-Anosov flow is torsion-free.
\end{coro}
Next, as in the case of Anosov flows, not only the fundamental group of any $3$-manifold supporting a pseudo-Anosov flow is infinite, but it also has  \emph{exponential growth}: 

\begin{defi}
  Consider $G$ a finitely generated group and $S$ a finite symmetric set of generators of $G$ (i.e. if $s\in S$, then $s^{-1}\in S$). We will say that $G$ has \emph{exponential growth} if the cardinal of the following set:
$$\Gamma_S(n)=\{g\in G| g \text{ can be expressed as a word of length at most $n$ in $S$}\}$$
grows exponentially with $n$.

\end{defi}
It should be noted that being of exponential growth does not depend on the choice of the generating set $S$. 
\begin{theorem}[M. Paternain \cite{Paternain}]
    If $M^3$ supports a pseudo-Anosov flow, then $\pi_1(M)$ has exponential growth. 
\end{theorem}
Finally, an additional restriction arises from the fact that the orbit space (see Theorem-Definition \ref{thdef.bifoliatedplane}) of any pseudo-Anosov flow can be compactified to a closed disk: 
\begin{theorem} [S.Fenley \cite{Fe4}]
    If $M$ is a $3$-manifold supporting a pseudo-Anosov flow, then $\pi_1(M)$ admits a free action by homeomorphisms on the circle $\mathbb{S}^1$.
\end{theorem}

\subsection{On the foliations in the universal cover}
Consider $\Phi$ a pseudo-Anosov flow on a smooth $3$-manifold $M$, $F^s,F^u$ the stable foliations of $\Phi$, $\widetilde{M}\approx \mathbb{R}^3$ the universal cover of $M$ and $\widetilde{\Phi}=(\widetilde{\Phi}^t)_{t\in\mathbb{R}},\widetilde{F^s},\widetilde{F^u}$ the lifted flow and the lifted foliations on $\widetilde{M}$. Our goal in this section consists in stating several results concerning the topological properties of the orbits of $\widetilde{\Phi}$ and of the foliations $\widetilde{F^s},\widetilde{F^u}$.

First, thanks to Proposition \ref{p.firstperiodicorbits} and Corollary \ref{c.periodicnontrivial}, we have that 

\begin{prop}\label{p.noprongcircleinr3}
    The singular codimension one foliations $\widetilde{F^s},\widetilde{F^u}$ do not contain any circle prong singularities. 
\end{prop}

\begin{defi}
    An embedded closed topological disk $\widetilde{R}$ in  $\widetilde{M}$ will be called a \emph{transverse standard polygon} (resp. \emph{transverse rectangle}, \emph{transverse $2k$-gon}) if it projects in $M$ to a transverse standard polygon (resp. transverse rectangle, transverse $2k$-gon) of $\Phi$.

   Moreover, a connected component of $\widetilde{F^s}\cap \widetilde{R}$ (resp. $\widetilde{F^u}\cap \widetilde{R}$) will be called a \emph{stable} (resp. \emph{unstable}) \emph{leaf} of $\widetilde{R}$. Finally, for any $x\in \widetilde{R}$, we will denote by $\widetilde{R}^s(x)$ (resp. $\widetilde{R}^u(x)$) the stable (resp. unstable) leaf of $\widetilde{R}$ containing $x$ and we will call the closure of any connected component of $\widetilde{R}-(\widetilde{R}^s(x)\cup \widetilde{R}^u(x))$ \emph{a quadrant of $x$ in $\widetilde{R}$}.  
\end{defi}
\begin{defi}
    Consider $\widetilde{\mathcal{O}}$ an orbit of $\widetilde{\Phi}$. We define the \emph{stable  separatrix of $\widetilde{\mathcal{O}}$} as the union of $\widetilde{\mathcal{O}}$ with a connected component of $\widetilde{F^s}(\widetilde{\mathcal{O}})$ (with respect to the leaf topology). We similarly define the notion of \emph{unstable separatrix}. 
\end{defi}
In the next three propositions, we will prove that even though the orbits of $\Phi$ and the leaves of $F^s$ and $F^u$  have a rather complex topological behavior in $M$, when lifted on the universal cover their behavior becomes relatively trivial: 
\begin{prop}\label{p.closedfoluniversalcover}
    The foliations $\widetilde{F^s}$, $\widetilde{F^u}$ consist of closed leaves.
\end{prop}
\begin{proof}
    We will prove the above result for $\widetilde{F^u}$ (the exact same argument applies for $\widetilde{F^s}$). Consider $\widetilde{L}$ a leaf of $\widetilde{F^u}$ and suppose that $\text{Clos}(\widetilde{L})-\widetilde{L}\neq \emptyset$. Take $x\in \text{Clos}(\widetilde{L})-\widetilde{L}$. 
    
    There exists $\widetilde{R}$ a transverse standard polygon in $\widetilde{M}$ containing $x$ in its interior and such that  $\widetilde{L}\cap \widetilde{R}$ consists of infinitely many connected components. By eventually cutting $\widetilde{R}$ into a finite number of rectangles (when $\widetilde{R}$ is not a rectangle), we can assume that there exists $\widetilde{R}$ a transverse rectangle containing $x$ and such that  $\widetilde{L}\cap \widetilde{R}$ consists of infinitely many connected components. Denote by $\mathcal{U}$ the leaf space of all unstable leaves of $\widetilde{R}$ ( $\mathcal{U}\cong[0,1]$) and endow $\mathcal{U}$ with an orientation. 
    
    Assume first that $\widetilde{L}$ is a regular leaf of $\widetilde{F^u}$. Let $L_1,L_2$ two connected components of $\widetilde{R}\cap \widetilde{L}$. Since there exist infinitely many such connected components, we can suppose that there exists a continuous and injective path $\gamma:[0,1]\rightarrow \widetilde{L}$ with $\gamma(0)\in L_1,\gamma(1)\in L_2$ and for which the holonomy of $\widetilde{F^u}$ along $\gamma$ defines a map from a neighborhood of $L_1$ in $\mathcal{U}$ to a neighborhood of $L_2$ in $\mathcal{U}$ that respects the orientation. Using this, one can construct as in Figure \ref{f.closedleafproof}, a closed transversal to $\widetilde{F^u}$, which contradicts Theorem \ref{t.transversalhomotopicto0}. 

\begin{figure}
    \centering
    \includegraphics[scale=0.5]{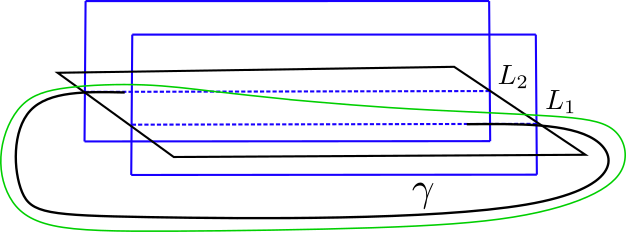}
    \caption{The green curve above represents a closed transversal to $\widetilde{F^u}$.}
    \label{f.closedleafproof}
\end{figure}
    Assume next that $\widetilde{L}$ is a singular leaf of $\widetilde{F^u}$. In this case, consider $\widetilde{L}'$ a separatrix of $\widetilde{L}$ intersecting $\widetilde{R}$ infinitely many times. By applying our previous argument for $\widetilde{L}'$, one can similarly construct a closed transversal to $\widetilde{F^u}$, which is impossible by Theorem \ref{t.transversalhomotopicto0}. 
\end{proof}

The arguments used during our proof of Proposition  \ref{p.closedfoluniversalcover}, also show that: 
\begin{rema}\label{r.finiteintersectionsleaves}
    Let $\widetilde{L}$ be a leaf of $\widetilde{F^s}$ or $\widetilde{F^u}$ and $\widetilde{R}$ a transverse standard rectangle (resp. $2k$-gon) in $\widetilde{M}$. We have that $\widetilde{R}\cap \widetilde{L}$ consists of at most two (resp. $4k$) connected components. 
\end{rema}

\begin{prop}\label{p.closedorbitsuniversalcover}
    The orbits of $\widetilde{\Phi}$ are all closed. 
\end{prop}
\begin{proof}
    Consider $\mathcal{O}$ a non-closed orbit of $\widetilde{\Phi}$ and $x\in \clos{\mathcal{O}}-\mathcal{O}$. There exists $\widetilde{R}$ a transverse standard polygon in $\widetilde{M}$ containing $x$ in its interior and intersecting $\mathcal{O}$ infinitely many times. This implies that the stable or unstable manifold of $\mathcal{O}$ intersects $\widetilde{R}$ along infinitely many stable or unstable leaves in $\widetilde{R}$, which is impossible by Remark \ref{r.finiteintersectionsleaves}. 
\end{proof}

\begin{prop}\label{p.regularleavesproperlyembedded}
     The regular leaves of the foliations $\widetilde{F^s}$, $\widetilde{F^u}$ are properly embedded surfaces in $\mathbb{R}^3$. 
\end{prop}

\begin{proof}
    Consider $\text{Sing}(\widetilde{F^s})=\text{Sing}(\widetilde{F^u})$ the set of all line prong singularities of  $\widetilde{F^s}$, $\widetilde{F^u}$ (see Proposition \ref{p.noprongcircleinr3}). Thanks to Remark \ref{r.singularfolidim3}, the singular foliations  $\widetilde{F^s}$, $\widetilde{F^u}$ define a pair of regular foliations in $\mathbb{R}^3 - \underset{C\in \text{Sing}(\widetilde{F^s})}{\cup}C$. Furthermore, thanks to Proposition \ref{p.closedfoluniversalcover},  any regular leaf of $\widetilde{F^s}$ or $\widetilde{F^u}$ is closed and in particular its closure cannot contain singular points of the foliations  $\widetilde{F^s}$ and $\widetilde{F^u}$. Therefore, in order to prove that the regular leaves of $\widetilde{F^s}$, $\widetilde{F^u}$ are properly embedded in $\mathbb{R}^3$, it suffices to show that the leaves of the restrictions of  $\widetilde{F^s}$, $\widetilde{F^u}$ on $\mathbb{R}^3 - \underset{C\in \text{Sing}(\widetilde{F^s})}{\cup}C$ are properly embedded. This follows from the fact that any closed leaf of a foliation in a second countable space is also properly embedded (see \cite{Chevallay}). 
\end{proof}
\begin{prop}\label{p.topologyliftperiodicleaves}
    Let $\widetilde{C}$ be an orbit of $\widetilde{\Phi}$ projecting to a periodic orbit of $\Phi$. If $\widetilde{C}$ is non-singular, then $\widetilde{F^s}(\widetilde{C})$ (resp. $\widetilde{F^u}(\widetilde{C})$)  is a properly embedded plane in $\mathbb{R}^3$ such that the restriction of $\widetilde{\Phi}$ on $\widetilde{F^s}(\widetilde{C})$ (resp. $\widetilde{F^u}(\widetilde{C})$) is conjugated to a trivial flow. 
\end{prop}
\begin{proof}
    We will prove the statement for $\widetilde{F^s}(\widetilde{C})$ (the same arguments apply for $\widetilde{F^u}(\widetilde{C})$). Thanks to Proposition \ref{p.topologyleaves}, if $C$ denotes the projection of $\widetilde{C}$ on $M$, then $F^s(C)$ is homeomorphic to an open cylinder or an open M\"obius band the fundamental group of which is generated by $C$. Recall that thanks to Coroallary \ref{c.periodicnontrivial}, $C$ is not homotopically trivial in $M$. It follows that by lifting to the universal cover, $F^s(C)$ is unwrapped to a plane. The fact that the previous plane is properly embedded follows from Proposition \ref{p.regularleavesproperlyembedded}. Finally, recall that the restriction of $\Phi$ on $F^s(C)$ defines a flow, for which every orbit is (topollogically) contracted by a unique periodic orbit generating the fundamental group of $F^s(C)$. Using our proof of Item (1) of Proposition \ref{p.topologyleaves}, one can easily construct a global transversal intersecting every orbit of the restriction of $\Phi$ on $F^s(C)$. The previous transversal lifts to a global transversal intersecting every orbit of the restriction of $\widetilde{\Phi}$ on $\widetilde{F^s}(\widetilde{C})$. Thanks to Poincar\'e-Bendixon theorem and to Corollary \ref{c.periodicnontrivial}, every orbit of the restriction of $\widetilde{\Phi}$  on $\widetilde{F^s}(\widetilde{C})$ intersects the previous transversal exactly once. It follows that the restriction of $\widetilde{\Phi}$ on $\widetilde{F^s}(\widetilde{C})$ is conjugated to a trivial flow, which finishes the proof of the proposition. 
\end{proof}
The following lemma constitutes an analogue of Proposition  \ref{p.topologyliftperiodicleaves} for  singular leaves of $\widetilde{F^s}$ and  $\widetilde{F^u}$ : 
 
\begin{prop}\label{p.complementsingleaves}
 Take $\widetilde{C}$ a line $p$-prong singularity of  $\widetilde{F^s}$. We have that $\widetilde{C}$ admits exactly $p$ distinct stable separatrices, each of which is homeomorphic to a closed half-plane. Furthermore, the union of any two stable separatrices of $\widetilde{C}$ forms a properly embedded plane in $\widetilde{M}$ and the restriction of $\widetilde{\Phi}$ on the previous plane is conjugated to a trivial flow. Finally, $\widetilde{M}-\widetilde{F^s}(\widetilde{C})$ consists of $p$ connected components homeomorphic to $\mathbb{R}^3$. 
\end{prop}
\begin{proof}
    Let  $\widetilde{S}$ be a stable separatrix of $\widetilde{C}$. Take $x\in \widetilde{C}$ and $\widetilde{R}$ a transverse $2p$-gon containing $x$ in its interior. Notice that $\widetilde{R}^s(x)-x$ consists of $p$ connected components. Since the (strictly) positive orbit of $x$ by $\widetilde{\Phi}$ can not revisit $x$ and $\widetilde{S}-\widetilde{C}$ has only one end terminating at $\widetilde{C}$ (see for instance Theorem 1.5 of \cite{Inaba}, this can also be deduced by  our proof of Item 1 of Proposition \ref{p.topologyleaves}), we get that $\widetilde{S}$ contains a unique connected component of $\widetilde{R}^s(x)-x$ and thus $\widetilde{C}$ admits exactly $p$ stable separatrices. 
    
Next, by an argument similar to the one we used in the proof of Proposition \ref{p.topologyliftperiodicleaves}, we have that 
\begin{itemize}
    \item $\widetilde{S}$ is a properly embedded copy of  $\mathbb{R}\times [0,+\infty)$ in $\mathbb{R}^3$
    \item the restriction of $\widetilde{\Phi}$ on $\widetilde{S}$ is conjugated to a horizontal  flow in $\mathbb{R}\times [0,+\infty)$
\end{itemize}
Hence, the union of two stable separatrices $\widetilde{S}_1, \widetilde{S}_2 $ of $\widetilde{C}$ defines a properly embedded plane in $\mathbb{R}^3$, which separates $\mathbb{R}^3$ in two connected components (see for instance \cite{Haefliger}) and such that the restriction of $\widetilde{\Phi}$ on the previous plane is conjugated to a trivial flow. 

Notice that thanks to our foliations charts,  $\widetilde{S}_1\cup \widetilde{S}_2 $ defines a locally flat submanifold of $\mathbb{R}^3$ (i.e. a submanifold for which for every $x\in \widetilde{S}_1\cup \widetilde{S}_2$, there exists $U_x$ a neighbordhood of $x$ in $\mathbb{R}^3$ and a homeomorphism $h:U_x\rightarrow \mathbb{R}^3$ such that $h\big(U_x\cap (\widetilde{S}_1\cup \widetilde{S}_2)\big)$ belongs in the affine plane $\mathbb{R}^2\times \{0\}$). By Theorem 3 of \cite{Brown}, there exists $V$ a neighborhood of $\widetilde{S}_1\cup \widetilde{S}_2$ and $g:\mathbb{R}^2\times(-1,1)\rightarrow V$ a homeomorphism such that $g(\mathbb{R}^2\times \{0\})=\widetilde{S}_1\cup \widetilde{S}_2$. Using this fact, one can easily deduce that the union of two stable separatrices $\widetilde{S}_1,\widetilde{S}_2 $ of $\widetilde{C}$ divide $\widetilde{M}$ into two  connected components homeomorphic to $\mathbb{R}^3$. By the same argument, by cutting one of the previous connected components along a third stable separatrix $\widetilde{S}_3$ of $\widetilde{C}$, we get that $\widetilde{M}- (\widetilde{S}_1\cup\widetilde{S}_2\cup \widetilde{S}_3)$ consists of three connected components homeomorphic to $\mathbb{R}^3$. We obtain the desired result by induction. 

\end{proof}
In the previous propositions, we proved that the topological behavior of the leaves of $\widetilde{F^s}$, $\widetilde{F^u}$ and the orbits of $\widetilde{\Phi}$ is rather easy to understand in the universal cover of $M$. In addition to the previous results, in the following pages we will show that, contrary to the case of their projections on $M$, a stable and an unstable leaf intersect along at most one orbit in the universal cover: 
\begin{prop}\label{p.oneintersection}
    Consider $\widetilde{C}$ an orbit of $\widetilde{\Phi}$ projecting in $M$ to a periodic orbit of $\Phi$. We have that $\widetilde{F^s}(\widetilde{C})\cap \widetilde{F^u}(\widetilde{C})=\widetilde{C}$.
\end{prop}
\begin{proof}
    Clearly, $\widetilde{C}\subset \widetilde{F^s}(\widetilde{C})\cap \widetilde{F^u}(\widetilde{C})$. Suppose now that $\widetilde{F^s}(\widetilde{C})\cap \widetilde{F^u}(\widetilde{C})$ contains a point that does not belong to $\widetilde{C}$. Since $\widetilde{F^s}(\widetilde{C})\cap \widetilde{F^u}(\widetilde{C})$ is invariant by the flow $\widetilde{\Phi}$, this would imply that there exists at least one $\widetilde{\Phi}$-orbit, say $\widetilde{C}'\neq\widetilde{C}$, contained in $\widetilde{F^s}(\widetilde{C})\cap \widetilde{F^u}(\widetilde{C})$. 

    Endow $M$ with a distance $d$ given by some Riemannian metric of $M$ and lift the previous Riemannian metric to a Riemannian metric on $\mathbb{R}^3$ giving rise to a distance $\widetilde{d}$. Take $x\in \widetilde{C}$ and $y\in \widetilde{C}'$.  Since $\widetilde{C}'$ and $\widetilde{C}$ are contained in the same stable manifold, by our definition of pseudo-Anosov flow, we can find an increasing homeomorphism $h:\mathbb{R}\rightarrow \mathbb{R}$ such that 
    \begin{equation}\label{eq.gettingclose}\widetilde{d}(\widetilde{\Phi}^t(x),\widetilde{\Phi}^{h(t)}(y))\underset{t\rightarrow +\infty}{\longrightarrow} 0
    \end{equation}

    It follows that there exists a transverse standard polygon in $\widetilde{M}$, say $\widetilde{R}$, and a stable leaf of $\widetilde{R}$, containing a point of both $\widetilde{C}$ and  $\widetilde{C}'$. Hence,  $\widetilde{F^u}(\widetilde{C})$ intersects $\widetilde{R}$ along at least two connected components. 

Consider now the action of $\pi_1(M)$ on $\widetilde{M}$ by deck transformations. By our definition of $d$ and $\widetilde{d}$, the previous action preserves the distance  $\widetilde{d}$. Recall that the projection of $\widetilde{C}$ on $M$, say $C$, is  homotopically non-trivial (see Corollary \ref{c.periodicnontrivial}). Fix $x_M$ the projection of $x$ on $M$ and endow $C$ with its dynamical orientation. The loop $C$ defines a  (non-trivial) element in $\pi_1(M,x_M)$, whose associated deck transformation preserving  $\widetilde{C}$ will be denoted by $g$. It is easy to see that $g(\widetilde{F^{s,u}}(\widetilde{C}))=\widetilde{F^{s,u}}(\widetilde{C})$ and that if we denote by $Per(C)$ the period of $C$, then $g^n(x)=\widetilde{\Phi}^{n\cdot Per(C)}(x)$ for every $n\in \mathbb{N}$. 

Recall that $\widetilde{R}$ intersected  $\widetilde{F^u}(\widetilde{C})$ along at least two connected components, say $L_1,L_2$. Let $$c_0=\underset{l_1\in L_1, l_2\in L_2}{\min}\widetilde{d}(l_1,l_2)$$

Since $g$ preserves $\widetilde{d}$ and  $\widetilde{F^u}(\widetilde{C})$,  we have that for every $N \in \mathbb{N}$ the standard transverse polygon $g^N(\widetilde{R})$ intersects $\widetilde{F^u}(\widetilde{C})$ along at least two connected components distanced by $c_0$. If $N$ is sufficiently big, thanks to \ref{eq.gettingclose}, the orbit $\widetilde{C}'$ will intersect $g^N(\widetilde{R})$ along a point that is closer to $g^N(x)=\widetilde{\Phi}^{N\cdot Per(C)}(x)$ than $c_0$. It follows that $g^N(\widetilde{R})$ intersects $\widetilde{F^u}(\widetilde{C})$ along at least three connected components or that $g^N(\widetilde{R})$ intersects $\widetilde{F^s}(\widetilde{C})$ along at least two connected components. By a finite repetition of this argument, we can construct a transverse standard polygon intersecting $\widetilde{F^u}(\widetilde{C})$ or $\widetilde{F^s}(\widetilde{C})$ along an arbitrarily big number of connected components, which contradicts Remark \ref{r.finiteintersectionsleaves}. 
\end{proof}

Using the above proposition and by repeating our proof of Proposition \ref{p.complementsingleaves}, one can prove the following : 

\begin{rema}\label{r.complementsingleaves}
    If $\widetilde{C}$ is a line $p$-prong singularity of $\widetilde{F^s}$ and $\widetilde{F^u}$, then we have that the union of a stable and unstable separatrix of $\widetilde{C}$ forms a properly embedded plane in $\widetilde{M}$ and that $\widetilde{M}-(\widetilde{F^s}(\widetilde{C})\cup\widetilde{F^u}(\widetilde{C}))$ consists of $2p$ connected components homeomorphic to $\mathbb{R}^3$. 
\end{rema}
The following theorem generalizes Proposition \ref{p.oneintersection} for any pair of stable and unstable leaves in $\widetilde{F^s}$ and $\widetilde{F^u}$:

\begin{theorem}\label{t.oneintersection}
    Let $L^s$ and $L^u$ be two leaves in $\widetilde{F^s}$ and $\widetilde{F^u}$ respectively. We have that either $L^s$ and $L^u$ do not intersect or they intersect along a unique orbit of $\widetilde{\Phi}$.
\end{theorem}

\begin{proof}
    Since both $L^s$ and $L^u$ are preserved by $\widetilde{\Phi}$, their intersections consist of orbits of $\widetilde{\Phi}$. Assume now that $L^s$ and $L^u$ intersect along more than one orbit of $\widetilde{\Phi}$. Take $\mathcal{O}$ and $\mathcal{O}'$ two $\widetilde{\Phi}$-orbits lying in $L^s\cap L^u$,  $x\in \mathcal{O}$ and $y\in \mathcal{O}'$. Endow $M$ with a distance $d$ given by some Riemannian metric of $M$ and lift the previous Riemannian metric to a Riemannian metric on $\mathbb{R}^3$ giving rise to a distance $\widetilde{d}$. Since $\mathcal{O}$ and $\mathcal{O}'$ belong in the same stable manifold, by our definition of pseudo-Anosov flow, we can find an increasing homeomorphism $h:\mathbb{R}\rightarrow \mathbb{R}$ such that 
    \begin{equation*}\widetilde{d}(\widetilde{\Phi}^t(x),\widetilde{\Phi}^{h(t)}(y))\underset{t\rightarrow +\infty}{\longrightarrow} 0
    \end{equation*}
    It follows that, by eventually changing our choice of $x,y$, there exists $\widetilde{R}$ a transverse standard polygon in $\widetilde{M}$ containing $x,y$ in its interior and such that $\widetilde{R}^s(x)= \widetilde{R}^s(y)$. By definition of $\widetilde{R}$, the stable leaf $ L^s$ intersects $\widetilde{R}$ along at least one stable leaf of $\widetilde{R}$ and $L^u$ along at least two unstable leaves of $\widetilde{R}$. Thanks to Proposition \ref{p.densitystablemanifoldsperiodic}, there exist $L_p^s\in \widetilde{F^s}$ and $L_p^u\in \widetilde{F^u}$ projecting on $M$ to the stable and unstable leaf of two  periodic orbits of $\Phi$ and such that $L_p^s$ (resp. $L_p^u$) intersects $\widetilde{R}$ along a stable leaf of $\widetilde{R}$ close to $\widetilde{R}^s(x)$ (resp.  two unstable leaves of $\widetilde{R}$ close to $\widetilde{R}^u(x)$ and $\widetilde{R}^u(y)$). By our definition of $L_p^s$ and $L_p^u$, we have that $L_p^s\cap L^u_p\cap \widetilde{R}$ contains two points $x_p$ and $y_p$ lying in the same stable leaf of $\widetilde{R}$ and that are respectively close to $x$ and $y$. If $x_p$ and $y_p$ lie in the same $\widetilde{\Phi}$-orbit, then the stable segment in $\widetilde{R}$ going from $x_p$ to $y_p$ intersects  twice the same $\widetilde{\Phi}$-orbit in $L^s_p$, while remaining transverse to the restriction of $\widetilde{\Phi}$ in $L_p^s$. This contradicts Propositions \ref{p.topologyliftperiodicleaves} and \ref{p.complementsingleaves}. By our previous argument, we get that $x_p$ and $y_p$ lie in two different orbits of $\widetilde{\Phi}$. We will therefore assume from now on that $L^s$ and $L^u$ project in $M$ to the stable and unstable manifolds of two periodic orbits of $\Phi$. 
    
    Notice that, because of Remark \ref{r.complementsingleaves}, two different stable 
    separatrices of a line prong singularity $\widetilde{C}$ 
    cannot intersect the same unstable leaf $\mathcal{L}$ in $\widetilde{F^u}$, except if $\mathcal{L}=\widetilde{F^u}(\widetilde{C})$. Indeed, if $\mathcal{L}\neq \widetilde{F^u}(\widetilde{C})$ and $\mathcal{L}$ intersects two distinct stable separatrices of $\widetilde{C}$, then thanks to Remark \ref{r.complementsingleaves}, $\mathcal{L}$ would be forced to intersect an unstable separatrix of $\widetilde{C}$, which is impossible. By our previous arguments, if $L^s$ is a singular stable leaf, then there would exist a stable separatrix of the line prong singularity in $L^s$ intersecting $L^u$ along two orbits. We can therefore distinguish two different cases: 
    \begin{enumerate}
        \item $L^s$ is a regular stable leaf of $\widetilde{F^s}$
        \item $L^s$ is the stable separatrix of some line prong singularity of $\widetilde{F^s}$
    \end{enumerate}
    We will prove that the first case is impossible. A similar argument can be used to prove that the second case is also impossible. Assume from now on that $L^s$ is a regular stable leaf.

    Take once again $x\in \mathcal{O}$ and $ y\in \mathcal{O}'$. Since $\mathcal{O}$ and $\mathcal{O}'$ belong in the same stable and unstable manifold, by our definition of pseudo-Anosov flow, we have that there exist two increasing homeomorphisms $h,h':\mathbb{R}\rightarrow \mathbb{R}$ such that $$\widetilde{d}(\widetilde{\Phi}^t(x),\widetilde{\Phi}^{h(t)}(y))\underset{t\rightarrow +\infty}{\longrightarrow} 0$$ $$\widetilde{d}(\widetilde{\Phi}^t(x),\widetilde{\Phi}^{h'(t)}(y))\underset{t\rightarrow -\infty}{\longrightarrow}0 $$ It follows from Proposition \ref{p.topologyliftperiodicleaves} that the orbits $\mathcal{O}$ and $\mathcal{O}'$ delimit a flow-band in $L^s$ the boundaries of which get progressively closer for $\widetilde{d}$ in both the future and the past. Using the fact that the restriction of $\widetilde{\Phi}$ in this band is conjugated to a trivial flow, we get that there exist two orbits in $L^s$ remaining arbitrarily close for $\widetilde{d}$ in both the future and the past. This contradicts the expansivity of $\Phi$ and finishes the proof of the theorem. 
\end{proof}

The previous theorem constitutes the most important step in the proof of the fact that $\widetilde{\Phi}$ admits an orbit space that is homeomorphic to a plane, a result that is crucial for the classification approach developed later in this paper. Before providing the proof of this theorem, let us finish this section with two important corollaries. The following corollary constitutes a stronger version of Remark \ref{r.finiteintersectionsleaves}:
\begin{coro}\label{c.complementsingleaves}
    Let $\widetilde{L}$ be a leaf in $\widetilde{F^s}$, $\widetilde{R}$ a transverse standard polygon of $\widetilde{\Phi}$ and $x\in \widetilde{L}\cap \widetilde{R}$. We have that $\widetilde{L}\cap \widetilde{R}=\widetilde{R}^s(x)$. 
\end{coro}
Naturally, the above lemma is also true for leaves in $\widetilde{F^u}$. 
\begin{proof}
Assume that $\widetilde{R}^s(x)\subsetneq \widetilde{L}\cap \widetilde{R}$ or equivalently that $\widetilde{L}$ intersects $\widetilde{R}$ along two or more stable leaves of $\widetilde{R}$. 
 
 \textbf{Case 1: $\widetilde{L}$ a regular leaf of $\widetilde{F^u}$ and $\widetilde{R}$ is a transverse rectangle}

 In that case, any unstable leaf $\widetilde{L'}$ in $\widetilde{F^u}$ that intersects $\widetilde{R}$, also intersects $\widetilde{L}\cap \widetilde{R}$ more than once. By eventually changing our choice of $\widetilde{L'}$, we may assume without any loss of generality that the projection of $\widetilde{L'}$ on $M$ is the unstable manifold of some periodic orbit of $\Phi$. The fact that  $\widetilde{L'}$ intersects $\widetilde{L}\cap \widetilde{R}$ more than once implies, thanks to Theorem \ref{t.oneintersection}, that there exists an orbit of $\widetilde{\Phi}$ in $\widetilde{L'}$ intersecting $\widetilde{R}$ along at least two points lying in the same unstable leaf of $\widetilde{R}$. Denote the previous two points by $x$ and $y$. The unstable segment in $\widetilde{R}$ going from $x$ to $y$ is a segment contained in $\widetilde{L}'$ intersecting twice the same $\widetilde{\Phi}$-orbit, while remaining transverse to the orbits of $\widetilde{\Phi}$. This contradicts Propositions  \ref{p.topologyliftperiodicleaves} and \ref{p.complementsingleaves}. 

\textbf{Case 2: $\widetilde{L}$ a singular leaf of $\widetilde{F^s}$ and $\widetilde{R}$ is a transverse rectangle}

By the exact same argument as before, one can prove that any  stable separatrix of the line prong singularity in $\widetilde{L}$ either does not intersect $\widetilde{R}$ or  intersects $\widetilde{R}$ along a unique stable leaf of $\widetilde{R}$.  

Assume now that there exist two distinct separatrices of the line prong singularity in $\widetilde{L}$ intersecting  $\widetilde{R}$. As $\widetilde{F^s}$ and $\widetilde{F^u}$ define on $\widetilde{R}$ a trivial bifoliation, this implies that there exists a regular unstable leaf in $\widetilde{F^u}$ which intersects two distinct separatrices of the line prong singularity in $\widetilde{L}$. This is impossible because of Remark \ref{r.complementsingleaves}, according to which any two stable separatrices of the line prong singularity in $\widetilde{L}$ are separated by an unstable separatrix of the same line prong singularity.

\textbf{Case 3: $\widetilde{L}$ is a regular leaf of $\widetilde{F^s}$ and $\widetilde{R}$ is a transverse $2p$-gon}
Let $C$ be the orbit by $\widetilde{\Phi}$ of the unique singular point in the interior of $\widetilde{R}$. By Lemma  \ref{p.complementsingleaves}, we have that $\mathbb{R}^3-\widetilde{F^s}(C)$ consists of $p$ connected components. Since $\widetilde{L}\cap\widetilde{F^s}(C)=\emptyset $, we get that $\widetilde{L}$ is contained inside a unique connected component of $\mathbb{R}^3-\widetilde{F^s}(C)$. 

Cut $\widetilde{R}$ along the quadrants of its singularity in order to obtain $2p$ transverse rectangles $R_1,...,R_{2p}$. By Case 2, $\widetilde{F^s}(C)$ intersects each of these rectangles along a unique stable segment (see Figure \ref{f.decomposeintorect}). This implies that the closure of any connected component of $\mathbb{R}^3-\widetilde{F^s}(C)$ contains exactly two rectangles among $R_1,...,R_{2p}$. Denote by $R_k$ and $R_l$ the two rectangles in  $\{R_1,...,R_{2p}\}$ that are contained in the closure of the same connected component of $\mathbb{R}^3-\widetilde{F^s}(C)$ as $\widetilde{L}$. Notice that, if a stable leaf of $\widetilde{R}$ intersects $R_k$ then it also intersects $R_l$ and vice versa. Using this fact, as $\widetilde{L}$ intersects $\widetilde{R}$ along two or more stable leaves, we get that $\widetilde{L}$ also intersects $R_k$ along two or more stable leaves, which is impossible by Case 1. 
\begin{figure}
    \centering
    \includegraphics[scale=0.45]{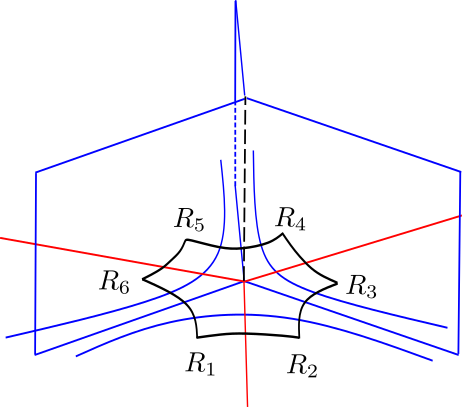}
    \caption{}
    \label{f.decomposeintorect}
\end{figure}

\textbf{Case 4: $\widetilde{L}$ a singular leaf of $\widetilde{F^u}$ and $\widetilde{R}$ is a transverse $2p$-gon}
Let $C$ be the orbit by $\widetilde{\Phi}$ of the unique singular point in the interior of $\widetilde{R}$. As in Case 3, cut $\widetilde{R}$ along the quadrants of its singularity in order to obtain $2p$ transverse rectangles $R_1,...,R_{2p}$. By Case 2, $\widetilde{F^s}(C)$ intersects each of these rectangles along a unique stable segment (see Figure \ref{f.decomposeintorect}). It follows that $\widetilde{L}\neq \widetilde{F^s}(C)$. Next, by Lemma  \ref{p.complementsingleaves}, we have that $\mathbb{R}^3-\widetilde{F^s}(C)$ consists of $p$ connected components. Since $\widetilde{L}\cap\widetilde{F^s}(C)=\emptyset $, we get that $\widetilde{L}$ is contained inside a unique connected component of $\mathbb{R}^3-\widetilde{F^u}(C)$. By the same argument that we used during the proof of Case 3, we get that $\widetilde{L}$ intersects a transverse rectangle $R_i$ along two or more unstable leaves of $R_i$, which is impossible by Case 2.

\end{proof}
\begin{coro}\label{c.orbitintersectpolygon}
    Consider $x\in \widetilde{M}$ and $\widetilde{R}$ a standard transverse polygon of $\widetilde{\Phi}$. There exists at most one $t\in \mathbb{R}$ such that $\widetilde{\Phi}^t(x)\in \widetilde{R}$. 
\end{coro}
    \begin{proof}
       Assume that there exist $t_1< t_2$ such that $X:=\widetilde{\Phi}^{t_1}(x)\in \widetilde{R}$ and $Y=\widetilde{\Phi}^{t_2}(x)\in \widetilde{R}$. Thanks to Corollary \ref{c.complementsingleaves} we get that the stable and unstable leaves of $\widetilde{R}$ containing $X$ and $Y$ coincide. Thus $X=Y$ and the orbit of $x$ by $\widetilde{\Phi}$ is periodic, which contradicts Corollary \ref{c.periodicnontrivial} and leads to an absurd.
    \end{proof}
\subsection{On the orbit space of a pseudo-Anosov flow}
Consider $\Phi$ a pseudo-Anosov flow on a smooth $3$-manifold $M$, $F^s,F^u$ the stable foliations of $\Phi$, $\widetilde{M}\approx \mathbb{R}^3$ the universal cover of $M$ and $\widetilde{\Phi},\widetilde{F^s},\widetilde{F^u}$ the lifted flow and the lifted stable and unstable foliations on $\mathbb{R}^3$. In this small section, we will define the \emph{orbit space} of a pseudo-Anosov flow and we will prove Barbot's theorem according to which the orbit space of $\Phi$ completely characterizes $\Phi$ up to orbital equivalence. 

\begin{lemm}\label{l.separatedorbits}
    Let $\mathcal{O}_1$ and $\mathcal{O}_2$ be two distinct orbits of $\widetilde{\Phi}$. There exist $U$ a neighborhood of $\mathcal{O}_1$ and $V$ a neighborhood of $\mathcal{O}_2$ such that $U,V$ are saturated by $\widetilde{\Phi}$ and $U\cap V= \emptyset$. 
\end{lemm}
\begin{proof}
    Assume that for every $\widetilde{\Phi}$-saturated neighborhood $U$ of $\mathcal{O}_1$ and every $\widetilde{\Phi}$-saturated neighborhood $V$ of $\mathcal{O}_2$, we have that $U\cap V\neq \emptyset$. By Theorem  \ref{t.oneintersection} the orbits $\mathcal{O}_1$ and $\mathcal{O}_2$ can not simultaneously belong in the same stable leaf in $\widetilde{F^s}$ and in the same unstable leaf in $\widetilde{F^u}$. Assume without any loss of generality that $\mathcal{O}_1$ and $\mathcal{O}_2$ \underline{belong in different stable leaves in $\widetilde{F^s}$}. 
    
Consider $o_1\in \mathcal{O}_1$, $o_2\in \mathcal{O}_2$ and $\widetilde{R}_1$, $\widetilde{R}_2$ two standard transverse polygons of $\widetilde{\Phi}$ containing respectively in their interior $o_1$ and $o_2$. Throughout this proof for any subset $A\subset \mathbb{R}^3$, we will denote by $\widetilde{\Phi}.A$ the set $\underset{x\in A}{\cup}\underset{t\in \mathbb{R}}{\cup}\widetilde{\Phi}^t(x)$. The sets $\widetilde{\Phi}.\widetilde{R}_1$ and $\widetilde{\Phi}.\widetilde{R}_2$
    are $\widetilde{\Phi}$-saturated neighborhoods of $\mathcal{O}_1$ and $\mathcal{O}_2$ respectively, thus $\widetilde{\Phi}.\widetilde{R}_1 \cap \widetilde{\Phi}.\widetilde{R}_2\neq \emptyset $. It follows that there exists $x_0\in \widetilde{R}_1$ whose orbit intersects $\widetilde{R}_2$ along a unique (see Corollary \ref{c.orbitintersectpolygon}) point, say $y_0\in \widetilde{R}_2$. 

    By taking a sequence of progressively smaller standard transverse polygons containing in their interiors $o_1$ and $o_2$, we can construct by the exact same argument a sequence of points $x_n\in \widetilde{R}_1$ such that: 
    \begin{itemize}
        \item $x_n\underset{n\rightarrow \infty}{\longrightarrow}o_1$
        \item the orbit of $x_n$ intersects $\widetilde{R}_2$ along a unique point, which we will denote by $y_n$
        \item $y_n\underset{n\rightarrow \infty}{\longrightarrow}o_2$
    \end{itemize}

    \textbf{Case 1: $\mathcal{O}_1$ and $\mathcal{O}_2$ belong in the same  unstable leaf in $\widetilde{F^u}$ }

In this case, we can find $o'_1\in \mathcal{O}_1$, $o'_2\in \mathcal{O}_2$ and $\widetilde{R}$ a standard transverse polygon of $\widetilde{\Phi}$ containing in its interior $o'_1$ and $o'_2$. For $n$ sufficiently big the $\widetilde{\Phi}$-orbit of $x_n$ (resp. $y_n$) intersects $\widetilde{R}$ along a point very close to $o'_1$ (resp. $o'_2$). Since the points $x_n,y_n$ lie in the same orbit, we get that there exists an orbit of $\widetilde{\Phi}$ that intersects twice $\widetilde{R}$, which contradicts Corollary \ref{c.orbitintersectpolygon}. 

\textbf{Case 2: $\mathcal{O}_1$ and $\mathcal{O}_2$ do not belong in the same stable or unstable leaf in $\widetilde{M}$ }

Following the previous notations, if $\widetilde{R}_1$ is a transverse $2k$-gon, by cutting $\widetilde{R}_1$ along the quadrants of its singularity, we obtain $2k$ transverse rectangles, one of which contains $o_1$ and a subsequence of  $(x_n)_{n\in \mathbb{N}}$. Therefore, by cutting if necessary $\widetilde{R}_1$ and $\widetilde{R}_2$, and by considering a subsequence of the $x_n$ we get the following: 
 \begin{itemize}
 \item $\widetilde{R}_1$ and $\widetilde{R}_2$ are transverse rectangles containing respectively $o_1$ and $o_2$ 
\item $x_n\in \widetilde{R}_1$ and $x_n\underset{n\rightarrow \infty}{\longrightarrow}o_1$
\item the orbit of $x_n$ intersects $\widetilde{R}_2$ along $y_n$
and $y_n\underset{n\rightarrow \infty}{\longrightarrow}o_2$
    \end{itemize}
    
Since $\widetilde{R}_1$ and $\widetilde{R}_2$ are rectangles, thanks to Theorem \ref{t.oneintersection}, each of the sets $\widetilde{F^u}(x_0)\cap \widetilde{F^s}(\mathcal{O}_1)$ and $\widetilde{F^u}(y_0)\cap \widetilde{F^s}(\mathcal{O}_2)$ consists of a unique orbit, which will be denoted by $\mathcal{O}_1'$ and $\mathcal{O}_2'$ respectively. By our initial hypothesis, $\mathcal{O}_1'\neq \mathcal{O}_2'$. Similarly, for every $n$ the set $\widetilde{F^u}(x_0)\cap \widetilde{F^s}(x_n)= \widetilde{F^u}(y_0)\cap \widetilde{F^s}(y_n)$ is non-empty and thus consists of a unique orbit $\mathcal{X}_n$ (see Theorem \ref{t.oneintersection}). The sequence of orbits $\mathcal{X}_n$ approaches simultaneously $\mathcal{O}_1'$ and $\mathcal{O}_2'$ as $n\rightarrow +\infty$. This implies that any $\widetilde{\Phi}$-saturated neighborhood $U'$ of $\mathcal{O}'_1$ intersects every $\widetilde{\Phi}$-saturated neighborhood $V'$ of $\mathcal{O}'_2$. We proved in Case 1 that this was impossible, which leads us to an absurd. 
\end{proof}
Thanks to the previous lemma, we are now ready to define the orbit space of a pseudo-Anosov flow: 
\begin{thdef}\label{thdef.bifoliatedplane}
    Let $\mathcal{P}$ be the space of orbits of $\widetilde{\Phi}$ endowed with the quotient topology. We have that 
    \begin{enumerate}
        \item $\mathcal{P}$ is homeomorphic to $\mathbb{R}^2$
        \item the projections on $\mathcal{P}$ of $\widetilde{F^s}$ and $\widetilde{F^u}$ define a pair of singular transverse foliations $\mathcal{F}^s,\mathcal{F}^u$ with no $1$-prong singularities 
        \item an action by deck transformations of $\pi_1(M)$ on $\mathbb{R}^3$ projects on $\mathcal{P}$ to an action by homeomorphisms preserving $\mathcal{F}^s,\mathcal{F}^u$
    \end{enumerate}
   We call $\mathcal{P}$ \emph{the bifoliated plane or the orbit space of $\Phi$} and $\mathcal{F}^s,\mathcal{F}^u$ the stable and unstable foliations of $\mathcal{P}$. 
\end{thdef}
    The above theorem constitutes a result of outmost importance in the theory of both smooth and topological pseudo-Anosov flows, as it reduces the study of pseudo-Anosov flows in dimension 3 to the study of a group action on a plane. Historically, Theorem-Definition \ref{thdef.bifoliatedplane} was simultaneously proven for $C^1$ Anosov flows by T.Barbot and S.Fenley in \cite{Barbotthese} and \cite{Fe1} respectively. In the case, of smooth pseudo-Anosov flows a proof of the above result may be found in \cite{Fenleymosher}. The existence of the bifoliated plane for topological pseudo-Anosov flows follows from an adaptation of one of the above arguments. For the sake of completeness, we provide here below a proof of the above theorem, by retracing the steps of the proof of Theorem 1.5.2 of \cite{Barbotthese}.
\begin{proof}
    Denote by $\pi^{\widetilde{\Phi}}:\mathbb{R}^3\rightarrow \mathcal{P}$ the projection from $\mathbb{R}^3$ to the space of orbits of $\widetilde{\Phi}$. Recall that, thanks to Corollary \ref{c.orbitintersectpolygon}, any orbit of $\widetilde{\Phi}$ intersects at most once a transverse standard polygon of $\widetilde{\Phi}$. Endow $\mathcal{P}$ with the quotient topology or equivalently the topology generated by the open sets of the form $\pi^{\widetilde{\Phi}}(\inte{P})$, where $\inte{P}$ is the interior of some transverse polygon $P$ of $\widetilde{\Phi}$. Thanks to Lemma \ref{l.separatedorbits}, the previous topology is Hausdorff.

    Consider $(P_i)_{i\in I}$ a collection of transverse standard polygons of $\widetilde{\Phi}$ such that any orbit of $\widetilde{\Phi}$ intersects the interior of some $P_i$. By eventually enlarging the previous family of transverse polygons, we will assume that $(P_i)_{i\in I}$ is invariant by the action of $\pi_1(M)$ by deck transformations on $\mathbb{R}^3$ and that for any two orbits  of $\widetilde{\Phi}$, say $\mathcal{O}_1$ and $\mathcal{O}_2$, there exist $i, j\in I$ such that $\mathcal{O}_1$ intersects the interior of $P_i$, $\mathcal{O}_2$ intersects the interior of $P_j$ and the $\widetilde{\Phi}$-saturations of $P_i$ and $P_j$ are disjoint. 
    
     Denote by $\mathcal{F}^2_h,\mathcal{F}^2_v$ the horizontal and vertical foliations in $\mathbb{R}^2$, $\mathcal{D}_2$ the closed euclidean square $\{z\in\mathbb{C}||\text{Re}(z)|\leq 1 \text{, } |{Im}(z)|\leq 1\}$, $\pi_p(z)=z^p$ for any $z\in \mathbb{C}$ and $p\in \mathbb{N}^*$, $\mathcal{D}_1= \pi_2(\mathcal{D}_2)$ and $\mathcal{D}_p= \pi_p^{-1}(\mathcal{D}_1)$ for any $p\geq 3$. The images by $\pi_2$ of the horizontal and vertical leaves of  $\mathcal{F}^2_{h}$ and $\mathcal{F}^2_{h}$ define respectively two singular foliations, that we will denote by $\mathcal{F}^1_{h}$ and $\mathcal{F}^1_{v}$. Similarly, for any $p\geq 3$ the pre-images by $\pi_p$ of the leaves of $\mathcal{F}^1_{h}$ and $\mathcal{F}^1_{v}$  define respectively two singular foliations, that we will denote by $\mathcal{F}^p_{h}$ and $\mathcal{F}^p_{v}$. 
    For every $i\in I$, if $P_i$ is a transverse rectangle (resp. $2p$-gon) denote by $\phi_i$ a homeomorphism from $P_i$ to $\mathcal{D}_2$ (resp. $\mathcal{D}_p$) taking the the stable and unstable leaves of $P_i$ to the leaves given by the restrictions of $\mathcal{F}^2_{h}$ and $\mathcal{F}^2_{v}$ on $\mathcal{D}_2$ (resp. $\mathcal{F}^p_{h}$ and $\mathcal{F}^p_{v}$ on $\mathcal{D}_p$).

    Denote by $\inte{P_i}$ the interior of the transverse standard polygon $P_i$. It is not difficult to prove that, thanks to our definition of $(P_i)_{i\in I}$, to Corollary \ref{c.orbitintersectpolygon} and to Lemma \ref{l.separatedorbits}, the set $$\mathcal{A}=\{(\pi^{\widetilde{\Phi}}(\inte{P_i}), \phi_i\circ(\pi^{\widetilde{\Phi}}_{|\inte{P_i}})^{-1})|i\in I\}$$ defines an atlas on $\mathcal{P}$, endowing $\mathcal{P}$ with a structure of a Hausdorff $C^0$ manifold. Notice now that since $\widetilde{\Phi}$ preserves $\widetilde{F^s}$, we can partition $\mathcal{P}$ in sets of the form $\pi^{\widetilde{\Phi}}(L)$, where $L\in \widetilde{F^s}$. The atlas $\mathcal{A}$ endows the previous partition with a structure of $C^0$ singular foliation. We will
    denote this foliation by $\mathcal{F}^s$. We can similarly define $\mathcal{F}^u$, the projection of $\widetilde{F^u}$ on $\mathcal{P}$. Once again thanks to the atlas $\mathcal{A}$, it is easy to see that the foliations $\mathcal{F}^s$ and $\mathcal{F}^u$ are transverse. 

    Next, since $\mathbb{R}^3$ is connected, so is $\mathcal{P}$. Moreover, as all orbits of $\widetilde{\Phi}$ are homeomorphic to $\mathbb{R}$ and intersect at most once every $P_i$ (see Corollary \ref{c.orbitintersectpolygon}) we get that $\pi^{\widetilde{\Phi}}$ defines a locally trivial fibration. We thus have the following short  exact sequence : 
    $$\pi_1(\mathbb{R}^3)\rightarrow \pi_1(\mathcal{P})\rightarrow \pi_0(\mathbb{R})$$ which implies that $\mathcal{P}$ is simply-connected. Notice that $\mathcal{P}$ can not be a sphere, since it either supports a regular foliation or a singular foliation without $1$-prong singularities. We conclude that $\mathcal{P}$ is homeomorphic to $\mathbb{R}^2$. 

    Finally, an action of $\pi_1(M)$ by deck transformations on $\mathbb{R}^3$ preserves the orbits of the lifted flow $\widetilde{\Phi}$, the lifted foliations $\widetilde{F^s}$ and $\widetilde{F^u}$ and sends transverse standard polygons of $\widetilde{\Phi}$ to transverse standard polygons of $\widetilde{\Phi}$. Thanks to our definition of the atlas $\mathcal{A}$, it is easy to see that the previous facts imply that any such action of  $\pi_1(M)$ on $\mathbb{R}^3$ descends to an action by homeomorphisms on $\mathcal{P}$, which finishes the proof of the theorem. 
\end{proof}

Let us make here an important comment on Item (3) of Theorem-Definition \ref{thdef.bifoliatedplane}. By classical covering space theory, $\pi_1(M)$ admits infinitely many actions on $\mathbb{R}^3$ by deck transformations. Furthermore, all the previous actions on $\mathbb{R}^3$ coincide up to performing an inner automorphism on $\pi_1(M)$. It follows that Item (3) of the previous theorem-definition defines infinitely many actions of $\pi_1(M)$ on $\mathcal{P}$ by homeomorphisms that are all the same up to performing an inner automorphism on $\pi_1(M)$. For the sake of simplicity and by an abuse of language, throughout this paper, any of the above actions of $\pi_1(M)$ on $\mathcal{P}$, will be called \emph{\textbf{the} action by $\pi_1(M)$ on $\mathcal{P}$ associated to $\Phi$}.

Next, the fact that the lifts on $\mathbb{R}^3$ of the stable and unstable foliations of $\Phi$ project to two singular foliations (in the sense of Definition \ref{d.singularfolisurfaces}) in $\mathcal{P}$ provides us with additional information on the topology of the leaves in $\widetilde{F}^s$ or $\widetilde{F}^u$. Indeed, thanks to Proposition \ref{p.regularleavesproperlyembedded} and to the fact that no singular foliation on the plane can admit compact leaves we have that 
\begin{coro}
    The regular leaves of $\widetilde{F}^s$ or $\widetilde{F}^u$ are properly embedded planes in $\mathbb{R}^3$
\end{coro}

We will finish this section with the proof of the following classical result of T.Barbot: 
\begin{theorem}\label{t.barbottheor}
  Let $M,N$ be two closed 3-manifolds and $\Phi=(\Phi^t)_{t\in \mathbb{R}}$, $\Psi=(\Psi^t)_{t\in \mathbb{R}}$ two pseudo-Anosov flows on $M$ and $N$ respectively. Denote by $\mathcal{P}_\Phi$ (resp. $\mathcal{P}_\Psi$) the orbit space of $\Phi$ (resp. $\Psi$), $\mathcal{F}^{s,u}_\Phi$ (resp. $\mathcal{F}^{s,u}_\Psi$) the stable and unstable foliations on $\mathcal{P}_\Phi$ (resp. $\mathcal{P}_\Psi$) and  $\rho_\Phi:\pi_1(M)\rightarrow \text{Homeo}(\mathcal{P}_\Phi)$ (resp. $\rho_\Psi:\pi_1(N)\rightarrow \text{Homeo}(\mathcal{P}_\Psi)$) the action of the fundamental group of $M$ (resp. $N$) on $\mathcal{P}_\Phi$ (resp. $\mathcal{P}_\Psi$) associated to $\Phi^t$ (resp. $\Psi^t$) by Theorem-Definition \ref{thdef.bifoliatedplane}. 
  
  The flows $\Phi$ and $\Psi$ are orbitally equivalent if and only if there exists a homeomorphism $h:\mathcal{P}_\Phi \rightarrow \mathcal{P}_\Psi$ and an isomorphism $\alpha:\pi_1(M)\rightarrow \pi_1(N)$ such that: 
  \begin{enumerate}
      \item $h(\mathcal{F}^{s}_\Phi)=\mathcal{F}^{s}_\Psi$ and $h(\mathcal{F}^{u}_\Phi)=\mathcal{F}^{u}_\Psi$
      \item $h$ is equivariant with respect to $\rho_\Phi$ and $\rho_\Psi$. More precisely, for every $g\in \pi_1(M)$ and every $x\in \mathcal{P}_\Phi$ we have that $$ h(\rho_\Phi(g)(x))=\rho_\Psi(\alpha(g))h(x)$$
  \end{enumerate}
\end{theorem}
The above theorem is of great importance in the theory of pseudo-Anosov flows as it opens the way for a systematic study of the class of pseudo-Anosov flows in dimension 3 (up to orbital equivalence) from the point of view of group actions on the plane. The author considers the second part of this paper to be one possible approach of the study of pseudo-Anosov from the plane perspective (an alternative approach was recently published in \cite{planeapproach}). We will now begin the proof of Theorem \ref{t.barbottheor}, which is an adaptation of the proof of Theorem 1.5.4 of \cite{Barbotthese} and Theorem 3.4 of \cite{Ba1}.
\begin{proof}
    Suppose first that there exists $H:M\rightarrow N$ an orbit equivalence between $\Phi^t$ and $\Psi^t$. Thanks to our definition of the stable and unstable foliations of a pseudo-Anosov flow, $H$ sends the stable (resp. unstable) foliation of $\Phi$ to the stable (resp. unstable) foliation of $\Psi$. Moreover, the map $H$ defines an isomorphism $H_*: \pi_1(M)\rightarrow \pi_1(N)$ and lifts to a homeomorphism $\widetilde{H}: \widetilde{M}\rightarrow \widetilde{N}$ with the following properties: 
    \begin{itemize}
\item $\widetilde{H}$ is an orbit equivalence between the lifted flows $\widetilde{\Phi}$ and $\widetilde{\Psi}$
\item $\widetilde{H}$ sends the lift of the stable (resp. unstable) foliation of $\Phi$ to the lift of the stable (resp. unstable) foliation of $\Psi$
\item  $\widetilde{H}$ is equivariant with respect to the actions of the fundamental groups by deck transformations
    \end{itemize}
    It is easy to see that $\widetilde{H}$ projects to a homeomorphism between $\mathcal{P}_\Phi$ and $ \mathcal{P}_\Psi$ with the desired properties.

\vspace{0.2cm}
Suppose now that there exists a homeomorphism $h$ and an isomorphism $\alpha$ as in the statement of the theorem. Let us prove that $\Phi$ and $\Psi$ are orbitally equivalent. Recall that Theorem 17A of \cite{Whitney} and the fact that $\Phi$ preserves a pair of transverse singular
foliations imply that any $x\in M$ is contained in the interior of a transverse
standard polygon in $M$. Using the previous fact and the compactness of $M$,  for any sufficiently small $\beta > 0$ there exists a finite number of transverse standard polygons $R_1, S_1, R_2, ..., R_n, S_n$ in $M$ with the following properties: 
\begin{itemize}
\item the $S_i$ are pairwise disjoint
\item $R_i\subset \inte{S_i}$, where $\inte{S_i}$ denotes the interior of $S_i$
\item $R_i$ and $S_i$ have the same number of stable (or unstable) boundary components for every $i\in \llbracket 1, n\rrbracket$
\item every orbit of size $4\beta$ intersects at most once every $S_i$; thus the flow $\Phi$ inside $\underset{t\in(-2\beta,2\beta)}{\cup}\Phi^t(S_i)$ is conjugated to a constant speed vertical flow 
\item $\underset{\llbracket 1, n\rrbracket}{\cup}\underset{t\in(-\beta,0)}{\cup}\Phi^t(R_i)=\underset{\llbracket 1, n\rrbracket}{\cup}\underset{t\in(0,\beta)}{\cup}\Phi^t(R_i)=M$ 
\end{itemize}

Fix a sufficiently small $\beta$ and such a family of polygons $R_1, S_1, R_2, ..., R_n, S_n$ in $M$. By lifting the previous family of polygons on the universal cover, we obtain an infinite family of transverse standard polygons $(\widetilde{R_i})_{i\in I}$, $(\widetilde{S_i})_{i\in I}$ such that $\widetilde{R_i}\subset \widetilde{S_i}$, the $\widetilde{S_i}$ are pairwise disjoint and $\underset{i\in I}{\cup}\underset{t\in[-\beta,\beta]}{\cup}\Phi^t(\widetilde{R_i})=\widetilde{M}$. 

Thanks to Corollary \ref{c.orbitintersectpolygon} the projection of a transverse standard polygon of $\widetilde{\Phi}$ on $\mathcal{P}_\Phi$ is a standard polygon with respect to the foliations $\mathcal{F}^{s}_\Phi$ and $\mathcal{F}^{u}_\Phi$ (see Definition \ref{d.standardpolygon}). Since $h(\mathcal{F}^{s}_\Phi)=\mathcal{F}^{s}_\Psi$ and $h(\mathcal{F}^{u}_\Phi)=\mathcal{F}^{u}_\Psi$, the image of a standard polygon in $\mathcal{P}_\Phi$ is a standard polygon in $\mathcal{P}_\Psi$. Now, as $\widetilde{N}$ defines a trivial line bundle over $\mathcal{P}_\Psi$, any standard polygon in $\mathcal{P}_\Psi$ can be lifted to a transverse standard polygon of $\widetilde{\Psi}$. We can define in this way a procedure that associates any transverse standard polygon of $\widetilde{\Phi}$ to some transverse standard polygon of $\widetilde{\Psi}$. Throughout this proof, we will say that a transverse standard polygon $r$ of $\widetilde{\Phi}$ and a transverse standard polygon $r'$ of $\widetilde{\Psi}$ are \emph{related} if one can be obtained from the other by the above procedure or equivalently if the image by $h$ of the projection of $r$ on $\mathcal{P}_\Phi$ corresponds to the projection of $r'$ on $\mathcal{P}_\Psi$. 
\addtocontents{toc}{\protect\setcounter{tocdepth}{2}}
\subsubsection*{Constructing of a family of polygons for $\widetilde{\Psi}$}

Let us associate to each polygon in $(\widetilde{R_i})_{i\in I}$ a standard transverse polygon of $\widetilde{\Psi}$ with which they are related. We will denote by $\widetilde{\rho_\Phi}$ and $\widetilde{\rho_\Psi}$ the actions by deck transformations of $\pi_1(M)$ and $\pi_1(N)$ on $\widetilde{M}$ and $\widetilde{N}$. By our construction, each of the families $(\widetilde{R_i})_{i\in I}$, $(\widetilde{S_i})_{i\in I}$ consists of $n$ polygons up to the action of $\widetilde{\rho_\Phi}$. Pick $\widetilde{R_1}\subset\widetilde{S_1},...,\widetilde{R_n}\subset\widetilde{S_n}$ a set of representatives of each orbit of polygons. Thanks to the  
procedure described in the previous paragraph, we can associate to $\widetilde{R_1}$ and $\widetilde{S_1}$ two transverse standard polygons $\widetilde{R'_1}$ and $\widetilde{S'_1}$ of $\widetilde{\Psi}$ such that  $\widetilde{R'_1}\subset\widetilde{S'_1}$. We similarly define $\widetilde{R'_2},\widetilde{S'_2},...,\widetilde{R'_n},\widetilde{S'_n}$. For every $g\in \pi_1(M)$ and $j\in \llbracket 1,n\rrbracket$ we now associate to the polygons $\widetilde{\rho_\Phi}(g)(\widetilde{R_j}),\widetilde{\rho_\Phi}(g)(\widetilde{S_j})$ the polygons $\widetilde{\rho_\Psi}(\alpha(g))(\widetilde{R'_j}), \widetilde{\rho_\Psi}(\alpha(g))(\widetilde{S'_j})$ (notice that $\widetilde{\rho_\Psi}(\alpha(g))(\widetilde{R'_j}) \subset \widetilde{\rho_\Psi}(\alpha(g))(\widetilde{S'_j})$). We thus constructed two families  $(\widetilde{R'_i})_{i\in I}$ and $(\widetilde{S'_i})_{i\in I}$ of transverse standard polygons of $\widetilde{\Psi}$ that are invariant by $\widetilde{\rho_\Psi}$, for which for every $i\in I$ the polygons $\widetilde{R_i}$ (resp. $\widetilde{S_i}$) and $\widetilde{R_i'}$ (resp. $\widetilde{S'_i}$) are related -this results from the fact that $h$ is equivariant with respect to the actions $\rho_\Phi$ and $\rho_\Psi$- and such that the association $\widetilde{R_i}\rightarrow \widetilde{R_i'}$ (resp. $\widetilde{S_i}\rightarrow \widetilde{S_i'}$) is equivariant with respect to the group actions $\widetilde{\rho_\Phi}$ and $\widetilde{\rho_\Psi}$.

\subsubsection*{Constructing a map $F:\widetilde{M}\rightarrow\widetilde{N}$ sending orbits of $\widetilde{\Phi}$ to orbits $\widetilde{\Psi}$}

Denote by $\widetilde{U_i}$ the neighborhood $\underset{t\in(-\beta,\beta)}{\cup}\Phi^t(\widetilde{S_i})$. As by our construction the family $(\widetilde{S_i})_{i\in I}$ projects to a family of $n$ polygons in $M$, the neighborhoods $(\widetilde{U_i})_{i\in I}$ project to $n$ neighborhoods in $M$, say $U_1$,...,$U_n$, whose interiors cover $M$ thanks to our hypotheses on $R_1,S_1,...,R_n,S_n$. Let $(v_i)_{i\in \llbracket 1, n\rrbracket}$ be a $C^0$ partition of unity subordinate to the cover $\inte{U_1},...,\inte{U_n}$. The partition of unity $(v_i)_{i\in \llbracket 1, n\rrbracket}$ lifts to a partition of unity subordinate to $(\inte{\widetilde{U_i}})_{i\in I}$, that we will denote by $(\widetilde{v_i})_{i\in I}$. 

Thanks to the fact that every orbit of $\Phi$ of size $4\beta$ intersects at most once every $S_i$, we get that any point $y\in \widetilde{M}$ can belong to at most $n$ neighborhoods of the family $(\widetilde{U_i})_{i\in I}$. Denote by 
$\widetilde{U_{i_1}},...,\widetilde{U_{i_s}}$ the previous neighborhoods. By definition of the $\widetilde{U_i}$, we get that there exist $x_{i_1}\in \widetilde{S_{i_1}},...,x_{i_s}\in \widetilde{S_{i_s}}$ and $t_{i_1},...,t_{i_s}\in (-\beta,\beta)$ such that for every $j\in \{i_1,...,i_s\}$ we have that $y=\widetilde{\Phi}^{t_j}(x_j)$. Denote by $\mathcal{O}$ the $\widetilde{\Phi}$-orbit of $y$ and by $\mathcal{O}'$ the orbit of $\widetilde{\Psi}$ associated to $\mathcal{O}$ by $h$. Fix $o'\in \mathcal{O}'$. Since the polygons $\widetilde{S_{i_1}}$ and $\widetilde{S'_{i_1}}$ are related, the orbit $\mathcal{O}'$ intersects $\widetilde{S'_{i_1}}$ at a unique point $\widetilde{\Psi^{T_1}}(o')$ (see Corollay \ref{c.orbitintersectpolygon}). We similarly define $\widetilde{\Psi^{T_2}}(o')\in \widetilde{S'_{i_2}},...,\widetilde{\Psi^{T_s}}(o')\in \widetilde{S'_{i_s}}$. Let  $F:\widetilde{M}\rightarrow\widetilde{N}$ be the map defined by: $$F(y):=\widetilde{\Psi}^{(T_1+t_{i_1})\widetilde{v_{i_1}}(y)+ ...+ (T_s+t_{i_s})\widetilde{v_{i_s}}(y)}(o')$$

Notice that $F(y)$ does not depend on our choice of $o'$. Let us further explain our choice of definition for $F$: if $y$ belongs in a neighborhood $\widetilde{U_{i_1}}$, then it would be natural to define ``$F(y)=\widetilde{\Psi}^{T_1+t_{i_1}}(o')$". Since in general $y$ belongs to multiple neighborhoods of $(\widetilde{U_i})_{i\in I}$, if we were to define $F$ in this way (by using only one randomly selected neighborhood containing $y$), then $F$ would certainly fail to be continuous. Our definition of $F$ defines the image of $y$ as the mean value (with respect to the partition of unity $(\widetilde{v_i})_{i\in I}$) of all its previous ``natural images". 

\subsubsection*{On the properties of the map $F$} It is not very difficult to prove that $F$ is continuous. Moreover, notice that by our construction $F$ sends any orbit $\mathcal{U}$ of $\widetilde{\Phi}$ to a subset of the orbit of $\widetilde{\Psi}$ associated to $\mathcal{U}$ by $h$. Hence, on the level of the bifoliated planes $F$ projects to the homeomorphism $h$. It follows that $F$ sends a stable/unstable leaf $L$ of $\widetilde{\Phi}$ to a subset of the stable/unstable leaf of $\widetilde{\Psi}$ associated to $L$ by $h$. 

Next, following our previous notations, if $g\in \pi_1(M)$, then changing $y$ by $\rho_\Phi(g)(y)$ :
\begin{enumerate}
    \item changes the neighborhoods $\widetilde{U_{i_1}},...,\widetilde{U_{i_s}}$ to $\rho_\Phi(g)(\widetilde{U_{i_1}}),...,\rho_\Phi(g)(\widetilde{U_{i_s}})$
    \item changes the polygons $\widetilde{S_{i_1}},...,\widetilde{S_{i_s}}$ to $\rho_\Phi(g)(\widetilde{S_{i_1}}),...,\rho_\Phi(g)(\widetilde{S_{i_s}})$. Hence, it also changes the polygons $\widetilde{S'_{i_1}},...,\widetilde{S'_{i_s}}$ to $\rho_\Psi(\alpha(g))(\widetilde{S'_{i_1}}),...,\rho_\Psi(\alpha(g))(\widetilde{S'_{i_s}})$
    \item changes the points $x_{i_1},...,x_{i_s}$ to $\rho_\Phi(g)(x_{i_1}),...,\rho_\Phi(g)(x_{i_s})$ 
    \item leaves invariant the times $t_{i_1},...,t_{i_s}$ 
    \item changes $o'$ to $\rho_\Psi(\alpha(g))(o')$
    \item changes the points $\widetilde{\Psi^{T_1}}(o'),...\widetilde{\Psi^{T_s}}(o')$ to $\widetilde{\Psi^{T_1}}(\rho_\Psi(\alpha(g))(o')),...,\widetilde{\Psi^{T_s}}(\rho_\Psi(\alpha(g))(o'))$
\end{enumerate}
From the above it follows that $F$ is equivariant with respect to the actions of $\rho_\Phi$ and $\rho_\Psi$. Thus, $F$ defines a continuous map $f:M\rightarrow N$ sending each orbit of $\Phi$ to a subset of an orbit of $\Psi$ and each stable/unstable leaf of $\Phi$ to a subset of a stable/unstable leaf of $\Psi$. 

The map $f$ induces an isomorphism between $\pi_1(M)$ and $\pi_1(N)$ and since $M,N$ are aspherical (see Theorem \ref{t.aspherical}) it also induces an isomoprhism between the $k$-th homotopy groups of $M$ and $N$ for every $k\geq 0$. By Whitehead's theorem, this implies that $f$ is a homotopy equivalence. This implies in particular that $f$ is surjective. 

By the construction of $F$ and $f$, two distinct orbits of $\Phi$ can not be sent by $f$ to subsets of the same orbit of $\Psi$. However, $f$ might fail to be injective along an orbit of $\Phi$. We will thus say that $f$ is \emph{orbitally injective}. Let us now modify $f$ in order to make $f$ injective. 

Let $u:M\times \mathbb{R}\rightarrow \mathbb{R}$ be the continuous map defined as follows: $$\forall t\in \mathbb{R} ~ \forall x\in M~\quad f(\Phi^t(x))=\Psi^{u(t,x)}(f(x))$$ 
It is not difficult to see that $u$ verifies the following properties: 
\begin{equation}\label{eq.behavioru1}
   \forall s,t\in \mathbb{R} ~ \forall x\in M~\quad u(t+s,x)=u(t,\Phi^s(x))+u(s,x) 
\end{equation} 
\begin{equation}\label{eq.behavioru2}
    \forall x\in M~\quad u(0,x)=0
\end{equation}

\subsubsection*{A condition producing injectivity}
Suppose that there exists $T>0$ sufficiently big such that $u(T,x)>0$ for every $x\in M$. By altering slightly our definition of $f$, we will produce an orbital equivalence between $\Phi$ and $\Psi$. 

Indeed, consider $u_T(t,x)=\frac{1}{T} \int_{t}^{t+T}u(s,x)\,ds$ and let $f_T:M\rightarrow N$ be the map defined by $$f_T(x)=\Psi^{u_T(0,x)}(f(x))$$ Clearly, $f_T$ is homotopic to $f$, thus $f_T$ remains surjective, $f_T$ sends each orbit of $\Phi$ to a unique orbit of $\Psi$, $f_T$ is orbitally injective and also satisfies: $$ \forall t\in \mathbb{R} ~ \forall x\in M~\quad f_T(\Phi^t(x))=\Psi^{u_T(t,x)}(f(x))$$ 

Notice that the derivative of $u_T(t,x)$ with respect to $t$ is positive for every $t\in \mathbb{R}$ and every $x\in M$: $$\frac{\partial}{\partial t} u_T(t,x)=\frac{1}{T}[u(t+T,x)-u(t,x)]=\frac{1}{T}u(T,\Phi^t(x))>0$$
Thanks to this fact, we get that $f_T$ is not only orbitally injective, but also injective along the orbits of $\Phi$. This shows that $f_T$ is injective and as $f_T$ sends positively oriented orbits of $\Phi$ to positively oriented orbits of $\Psi$, we get that $f_T$ is an orbital equivalence between $\Phi$ and $\Psi$. 

Therefore, in order to finish the proof of Theorem \ref{t.barbottheor}, it suffices to prove that there exists $T>0$ sufficiently big so that $u(T,x)>0$ for every $x\in M$. We will name the previous condition $(\star)$.
\subsubsection*{Proof of $(\star)$ condition}
Suppose  that there exists a sequence of points in $M$, say $(x_k)_{k\in \mathbb{N}}$, and a sequence of positive times $t_k\underset{k\rightarrow +\infty}{\longrightarrow} +\infty$ such that 
\begin{equation}\label{eq.negativehypothesis}
  u(t_k,x_k)\leq 0  
\end{equation} 

Consider $d_M$ (resp. $d_N$) the distance associated to a Riemannian metric on $M$ (resp. $N$) and $\widetilde{d_M}$ (resp. $\widetilde{d_N}$) the distance associated to the lift of the previous Riemannian metric on $\widetilde{M}$ (resp. $\widetilde{N}$). Recall that, in the beginning of this proof, we defined in $M$ two families of transverse standard polygons $R_1,...,R_n$, $S_1,...,S_n$ (with $R_i\subset S_i$) such that any orbit segment of $\Phi$ of size $\beta$ intersects $\underset{i\in \llbracket 1, n\rrbracket}{\cup}R_i$. Recall also that we denoted by $(\widetilde{R_i})_{i\in I}, (\widetilde{S_i})_{i\in I}$ (with $\widetilde{R_i}\subset \widetilde{S_i}$) the lifts of the previous two families of polygons on the universal cover. 

\textbf{Claim.} Consider $A,\delta>0$, where $\delta$ is sufficiently small. There exist $c_1(\delta),c_2(\delta),c_3(\delta,A)>0$ with $ c_1(\delta),c_2(\delta), c_3(\delta,A)\underset{\delta\rightarrow 0}{\longrightarrow}0 $ such that for every $i\in I$ and every $\widetilde{u},\widetilde{v}\in \widetilde{S_i}$ with $\delta/2 \leq \widetilde{d_M}(\widetilde{u},\widetilde{v})\leq \delta$ we have that 
\begin{equation}\label{eq.distance}
    c_1(\delta)\leq \widetilde{d_N}(F(\widetilde{u}),F(\widetilde{v}))\leq c_2(\delta)
\end{equation}
\begin{equation}\label{eq.distancemin}
    \min_{t_1, t_2\in [-A,A]}\widetilde{d_N}(\widetilde{\Psi}^{t_1}(F(\widetilde{u})),\widetilde{\Psi}^{t_2}(F(\widetilde{v})))\geq c_3(\delta,A)
\end{equation}

\begin{proof}[Proof of the claim.] Let us prove that Inequality \ref{eq.distance} holds (the proof of Inequality \ref{eq.distancemin} can be obtained by a similar argument). The existence of $c_2(\delta)$ is an immediate consequence of the fact that $F$ is a lift of a continuous map between two compact manifolds. Next, recall that $F$ is orbitally injective and that every $\widetilde{\Phi}$-orbit intersects at most once any transverse standard polygon in $\widetilde{M}$ (see Corollary \ref{c.orbitintersectpolygon}). Hence,  $F(\widetilde{u}), F(\widetilde{v})$ belong to different $\widetilde{\Psi}$-orbits. Suppose now that there exist two sequences $(\widetilde{u_j})_{j\in \mathbb{N}}, (\widetilde{v_j})_{j\in \mathbb{N}}\in \widetilde{M}^{\mathbb{N}}$ such that $\delta/2 \leq \widetilde{d_M}(\widetilde{u_j},\widetilde{v_j})\leq \delta$, the points $\widetilde{u_j},\widetilde{v_j}$ belong in $\widetilde{S_j}\in (\widetilde{S_i})_{i\in I}$ for every $j\in \mathbb{N}$ and $$\widetilde{d_N}(F(\widetilde{u_j}),F(\widetilde{v_j}))\underset{j \rightarrow +\infty}{\longrightarrow}0 $$ 

Since the family $(\widetilde{S_i})_{i\in I}$ is finite up to the action of $\rho_\Phi$, $F$ is equivariant with respect to the actions $\rho_\Phi, \rho_\Psi$  and $\rho_\Phi$ (resp. $\rho_\Psi$) acts on $(\widetilde{M}, \widetilde{d_M})$ (resp. $(\widetilde{N}, \widetilde{d_N})$) by isometries, by eventually considering subsequences and by changing if necessary our choice of $\widetilde{u_j},\widetilde{v_j}$, we may assume without any loss of generality that 
\begin{enumerate}
    \item for every $j\in \mathbb{N}$ the points $\widetilde{u_j},\widetilde{v_j}$ belong in a unique transverse standard polygon $\widetilde{S}\in (\widetilde{S_i})_{i\in I}$
    \item $\widetilde{d_N}(F(\widetilde{u_j}),F(\widetilde{v_j}))\underset{j \rightarrow +\infty}{\longrightarrow}0$
    \item $ \widetilde{u_j}\underset{j \rightarrow +\infty}{\longrightarrow}\widetilde{U}\in \widetilde{S}$ and $\widetilde{v_j}\underset{j \rightarrow +\infty}{\longrightarrow}\widetilde{V}\in \widetilde{S}$
\end{enumerate} Notice that because of the inequality $\delta/2 \leq \widetilde{d_M}(\widetilde{u_j},\widetilde{v_j})\leq \delta$, the points $\widetilde{U}$ and $\widetilde{V}$ can not be equal. Also, thanks to the continuity of $F$, $F(\widetilde{U})=F(\widetilde{V})$ and therefore, since $F$ is orbitally injective, the points $\widetilde{U},\widetilde{V}$ lie on the same $\widetilde{\Phi}$-orbit, which intersects twice the transverse polygon $\widetilde{S}$. This is impossible by Corollary \ref{c.orbitintersectpolygon} and finishes the proof of the claim. 
\end{proof}

Let us now resume the proof of property $(\star)$. Let $p:\widetilde{M} \rightarrow \underset{i\in I}{\cup} \widetilde{R_i}$ be the map associating to every point $\widetilde{x}\in \widetilde{M}$ the first point of intersection of the (strictly) positive $\widetilde{\Phi}-$orbit of $\widetilde{x}$ with $\underset{i\in I}{\cup} \widetilde{R_i}$. For every $i\in I$, denote by $\pi_i:\underset{t\in(-2\beta,2\beta)}{\cup}\widetilde{\Phi}^t(\widetilde{S_i})\rightarrow \widetilde{S_i}$ the projection on $\widetilde{S_i}$ along the flowlines of $\widetilde{\Phi}$. Since the projection of the family $ (\widetilde{R_i})_{i\in I}$ on $M$ consists of a finite number of polygons and $M$ is compact, there exists $\epsilon>0$ such that for every $i\in I$ and every $\widetilde{x}\in \widetilde{M}$ such that $p(\widetilde{x})\in \widetilde{R_i}$ we have that $\pi_i$ is well defined on the $\epsilon$-ball around $\widetilde{x}$. We will denote by $p_{\widetilde{x}}$ the continuous map  associating the points of the $\epsilon$-ball of $\widetilde{x}$ to their projection on $\widetilde{S_i}$. If for two $\epsilon$-close points $\widetilde{x}, \widetilde{y}\in \widetilde{M}$, the points $p(\widetilde{x})$ and $p_{\widetilde{x}}(\widetilde{y})$ are also  $\epsilon$-close, we will denote by $p_{\widetilde{x}}^2(\widetilde{y})$ the point $p_{p(\widetilde{x})}(p_{\widetilde{x}}(\widetilde{y}))$. We similarly define (when this is possible) $p_{\widetilde{x}}^k(\widetilde{y})$ for every $k\in \mathbb{N}$.

Recall that we previously defined $(x_k)_{k\in \mathbb{N}}, (t_k)_{k\in \mathbb{N}}$ such that  $t_k\underset{k\rightarrow +\infty}{\longrightarrow} +\infty$ and $u(t_k,x_k)\leq 0$. Fix $\delta>0$ sufficiently small and choose for every $k\in \mathbb{N}$ a lift of $x_k$ on $\widetilde{M}$, say $\widetilde{x_k}$, and a point  $\widetilde{y_k}\in \widetilde{F^s}(\widetilde{x_k})$ such that the two following conditions are verified $$\widetilde{d_M}(\widetilde{x_k},\widetilde{y_k})\leq \epsilon$$ $$ \delta/2 \leq \widetilde{d_M}(p(\widetilde{x_k}),p_{\widetilde{x_k}}(\widetilde{y_k}))\leq \delta$$ 

Notice that since the points $p(\widetilde{x_k}),p_{\widetilde{x_k}}(\widetilde{y_k})$ belong in the same standard transverse polygon in $(\widetilde{S_i})_{i\in I}$, we have that the $\widetilde{\Phi}$-orbits of $p(\widetilde{x_k}),p_{\widetilde{x_k}}(\widetilde{y_k})$ (and thus also the $\widetilde{\Phi}$-orbits of $\widetilde{x_k}, \widetilde{y_k}$) are different. We also remark that, by Theorem 1.5 of \cite{Inaba} and by Lemma \ref{l.stablesetalternativedefinition} (applied for the families $R_1,...,R_n,S_1,...S_n$ on $M$), if $\delta$ is taken sufficiently small, then for every $k\in \mathbb{N}$ and $m\in \mathbb{N}$, the point $p^m_{\widetilde{x_k}}({\widetilde{y_k}})$ is well defined. Even more, for $\delta$ sufficiently small, thanks to Lemma 2.7 of \cite{Oka} (applied once again for $R_1,...,R_n,S_1,...S_n$) the distance $\widetilde{d_M}(p^m(\widetilde{x_k}), p^m_{\widetilde{x_k}}({\widetilde{y_k}}))$ goes to zero uniformly with $k$ when $m\rightarrow +\infty$. Thanks to the previous fact and the fact that any positive $\widetilde{\Phi}$-orbit intersects $\underset{i\in I}{\cup}\widetilde{R_i}$ in a uniformly bounded time, it is not hard to prove that there exists for every $k\in \mathbb{N}$ an increasing homeomorphism $h_k\in \text{Homeo}^+(\mathbb{R})$ such that for every $t\geq 0$  
$$\widetilde{d_M}(\widetilde{\Phi}^t(\widetilde{x_k}), \widetilde{\Phi}^{h_k(t)}({\widetilde{y_k}}))\leq C(t)$$
where $C:\mathbb{R}^+\rightarrow \mathbb{R}^+$ is a continuous function (that does not depend on $k$) going to zero when  $t\rightarrow +\infty$.

By applying $F$ and by using the fact that $F$ is a lift of a continuous map between two compact manifolds, we get that for every $t\geq 0$

$$\widetilde{d_N}\big(F(\widetilde{\Phi}^t(\widetilde{x_k})), F(\widetilde{\Phi}^{h_k(t)}({\widetilde{y_k}}))\big)=\widetilde{d_N}\big(\widetilde{\Psi}^{\widetilde{u}(t,\widetilde{x_k})}(F(\widetilde{x_k})), \widetilde{\Psi}^{\widetilde{u}(h_k(t),\widetilde{y_k})}(F({\widetilde{y_k}}))\big)\leq D(t)$$
where $\widetilde{u}$ is the lift of $u$ on $\widetilde{M}\times \mathbb{R}$ and $D:\mathbb{R}^+\rightarrow \mathbb{R}^+$ is a continuous function (that does not depend on $k$) going to zero when  $t\rightarrow +\infty$.

In particular, 
\begin{equation}\label{eq.absurd}
    \widetilde{d_N}\big(\widetilde{\Psi}^{\widetilde{u}(t_k,\widetilde{x_k})}(F(\widetilde{x_k})), \widetilde{\Psi}^{\widetilde{u}(h_k(t_k),\widetilde{y_k})}(F({\widetilde{y_k}}))\big)\leq D(t_k)\underset{k\rightarrow +\infty}{\longrightarrow}0
\end{equation}

We will now prove that for $k$ sufficiently big the $\widetilde{\Psi}$-orbit of $F(\widetilde{y_k})$ intersects twice a transverse standard polygon of  $\widetilde{\Psi}$, which will lead us to an absurd. Take $P_1,...,P_m,T_1,...,T_m $ a family of standard transverse polygons in $N$ such that 
\begin{itemize}
\item the $T_i$ are pairwise disjoint
\item $P_i\subset \inte{T_i}$, where $\inte{T_i}$ denotes the interior of $T_i$
\item $P_i$ and $T_i$ have the same number of stable (or unstable) boundary components for every $i\in \llbracket 1, m\rrbracket$
\item every orbit of size $4\alpha$ intersects at most once every $T_i$; thus the flow $\Psi$ inside $\underset{t\in(-2\alpha,2\alpha)}{\cup}\Psi^t(T_i)$ is conjugated to a constant speed vertical flow 
\item $\underset{\llbracket 1, m\rrbracket}{\cup}\underset{t\in(-\alpha,0)}{\cup}\Psi^t(P_i)=\underset{\llbracket 1, m\rrbracket}{\cup}\underset{t\in(0,\alpha)}{\cup}\Psi^t(P_i)=N$ 
\end{itemize}

Consider $(\widetilde{P_j})_{j\in J}$ and $(\widetilde{T_j})_{j\in J}$ the lifts on $\widetilde{N}$ of the previous families of standard transverse polygons. We will assume without any loss of generality that for every $j\in J$ we have that $\widetilde{P_j}\subset \widetilde{T_j}$. Thanks to Inequality \ref{eq.distance} and our definition of $\widetilde{x_k}, \widetilde{y_k}$, if $\delta$ is sufficiently small, then  $$c_1(\delta) \leq \widetilde{d_N}(o_k,o'_k)\leq c_2(\delta)$$ where $o_k:=F(p(\widetilde{x_k}))$ and $o'_k:=F(p_{\widetilde{x_k}}(\widetilde{y_k}))$. Denote by $\widetilde{X_k}\in \widetilde{P_k}$ the first point of intersection between the (strictly) positive $\widetilde{\Psi}$-orbit of $o_k$ and $\underset{j\in J}{\cup}\widetilde{P_j}$. For every $j\in J$ denote also by $\Pi_j:\underset{t\in(-2\alpha,2\alpha)}{\cup}\widetilde{\Psi}^t(\widetilde{T_j})\rightarrow \widetilde{T_j}$ the projection on $\widetilde{T_j}$ along the flowlines of $\widetilde{\Psi}$. Thanks to our definition of the families of polygons $(\widetilde{P_j})_{j\in J}$,  $(\widetilde{T_j})_{j\in J}$ and to the compactness of $N$, if $\delta$ (and thus also $c_2(\delta)$) is sufficiently small, then $\widetilde{Y_k}:=\Pi_j(o'_k)\in \widetilde{T_k}$ is well defined for any $k\in \mathbb{N}$. Finally, thanks to Inequality \ref{eq.distancemin} and our definition of $\Pi_k$, $$\widetilde{d_N}(\widetilde{X_k},\widetilde{Y_k})\geq c_3(\delta, 2\alpha)$$ 

We have thus constructed a first point of intersection between the $\widetilde{\Psi}$-orbit of $F(\widetilde{y_k})$ and $\widetilde{T_k}$. Let us now move on to the construction of a second point of intersection between the $\widetilde{\Psi}$-orbit of $F(\widetilde{y_k})$ and $\widetilde{T_k}$. Similarly to our construction of the homeomorphisms $h_k$, thanks to Theorem 1.5 of \cite{Inaba}, Lemma \ref{l.stablesetalternativedefinition} and Lemma 2.7 of \cite{Oka}, for any $\widetilde{x}\in \widetilde{N}$ and any $\widetilde{y}\in \widetilde{F^s}(\widetilde{x})$ such that $\widetilde{d_N}(\widetilde{x},\widetilde{y})$ is sufficiently small there exists an increasing homeomorphism $h\in \text{Homeo}^+(\mathbb{R})$ such that  
$\widetilde{d_N}(\widetilde{\Psi}^t(\widetilde{x}), \widetilde{\Psi}^{h(t)}(\widetilde{y}))$ goes uniformly to zero when  $t\rightarrow +\infty$ independently of our choice of $\widetilde{x},\widetilde{y}$ (provided that $\widetilde{d_N}(\widetilde{x},\widetilde{y})$ is uniformly small). Recall that $\widetilde{y_k}\in \widetilde{F^s}(\widetilde{x_k})$ and that $F$ sends the stable/unstable leaves of $\widetilde{\Phi}$ to the stable/unstable leaves of $\widetilde{\Psi}$. Using the results of this paragraph, the fact that $\widetilde{u}(t_k,\widetilde{x_k})=u(t_k,x_k)\leq 0$ (see Equation \ref{eq.negativehypothesis}) and Equation  \ref{eq.absurd}, we get that for $k$ sufficiently big there exists a point in the $\widetilde{\Psi}$-orbit of $F({\widetilde{y_k}})$, whose distance from $F(\widetilde{x_k})$ is arbitrarily small compared to $c_3(\delta, 2\alpha)$.

Next, by our definition of $p$ and $\beta$, we have that for every $k\in \mathbb{N}$ there exists $\beta>s_k>0$ such that $\widetilde{\Phi}^{s_k}(\widetilde{x_k})=p(\widetilde{x_k})$. It follows that $$o_k=F(p(\widetilde{x_k}))=\widetilde{\Psi}^{\widetilde{u}(s_k,p(\widetilde{x_k}))}(F(\widetilde{x_k}))$$ where $|\widetilde{u}(s_k,p(\widetilde{x_k}))|$ is uniformly bounded with respect to $k$, since $\widetilde{u}$ lifts a continuous function from  $\mathbb{R}\times M$ to $\mathbb{R}$. 
By our previous arguments, we get that for $k$ sufficiently big there exists a point in the $\widetilde{\Psi}$-orbit of $F({\widetilde{y_k}})$, whose distance from $o_k=F(p(\widetilde{x_k}))$ is arbitrarily small compared to $c_3(\delta, 2\alpha)$. The projection by $\Pi_k$ of the previous point on $\widetilde{T_k}$ defines a point $\widetilde{Z_k}$ lying in the intersection of the orbit of $F({\widetilde{y_k}})$ with $\widetilde{T_k}$. For $k$ sufficiently big and $\delta$ sufficiently small, $$\widetilde{d_N}(\widetilde{Z_k},\widetilde{X_k})< c_3(\delta, 2\alpha)< \widetilde{d_N}(\widetilde{Y_k},\widetilde{X_k})$$ thus $\widetilde{Z_k}\neq \widetilde{Y_k}$. We have thus proven that for $k$ sufficiently big and $\delta$ sufficiently small, the $\widetilde{\Psi}$-orbit of $F({\widetilde{y_k}})$ intersects the transverse standard polygon $\widetilde{T_k}$ along the two distinct points $\widetilde{Y_k}$ and $\widetilde{Z_k}$. This is impossible because of Corollary \ref{c.orbitintersectpolygon} and finishes the proof of the theorem. 

\end{proof}
\addtocontents{toc}{\protect\setcounter{tocdepth}{3}}
\subsection{On the stable and unstable foliations in the bifoliated plane}

Our goal in this section consists in proving the analogues of Theorem \ref{t.transversalhomotopicto0}, Lemma \ref{p.complementsingleaves}, Corollary \ref{c.complementsingleaves} and Theorem \ref{t.oneintersection} for the stable and unstable foliations in the bifoliated plane of a pseudo-Anosov flow. The proofs of this section do not really rely on the nature of the pseudo-Anosov dynamics, but rather on the topological restrictions arising when considering pairs of singular foliations on the plane. Throughout this section, we will thus consider pairs of singular foliations on the plane not necessarily associated with some pseudo-Anosov flow in dimension $3$.

\begin{defi}
    Let $F$ be a singular foliation on $\mathbb{R}^2$ (in the sense of Definition \ref{d.singularfolisurfaces}) and $S$ its set of singular points. Take $y\in \mathbb{R}^2$. The union of $y$ with a connected component of $F(y)-\{y\}$ (with respect to the leaf topology) will be called \emph{an $F$-separatrix of $y$}.

    Next, take $x\in S$, $U$ a small open neighborhood of $x$ and $F_U(x)$ the arc connected component of $F(x)\cap U$ containing $x$. We will say that the closure of each connected component of $U-F_U(x)$ is a \emph{sector of $x$ in $U$}.

    A simple closed continuous curve $\gamma: \mathbb{S}^1\rightarrow \mathbb{R}^2$ will be called a \emph{simple closed quasi-transversal of $F$} if the following conditions are satisfied: 
    \begin{itemize}
        \item there exists a finite set $T\subset \mathbb{S}^1$ such that for every $t\in T$, $\gamma(t)\in S$ 
        \item for every $t\in \mathbb{S}^1-T$ and any small neighborhood $J$ of $t$, $\gamma(J)\subset \mathbb{R}^2-S$ is topologically transverse to $F$
        
        \item for every $t\in T$, if $V$ is a small neighborhood of $\gamma(t)$ in $\mathbb{R}^2$, then for every small neighborhood $I$ of $t$ we have that $\gamma(I)$ is not contained in a unique sector of $\gamma(t)$ in $V$.  
    \end{itemize}
\end{defi}
Consider $\mathcal{F}^+$ and $\mathcal{F}^-$ two transverse singular foliations on $\mathbb{R}^2$ (in the sense of Definition \ref{d.singularfolisurfaces}). Denote by $\Sigma$ the set of singular points of $\mathcal{F}^+$ and $\mathcal{F}^-$ and assume that \underline{no two points in $\Sigma$ lie on the same leaf of $\mathcal{F}^+$ or $\mathcal{F}^-$}.

\begin{prop}\label{p.noquasitransvontheplane}
    The foliations $\mathcal{F}^+$ and $\mathcal{F}^-$ do not admit simple closed quasi-transversals. 
\end{prop}
\begin{proof}
    This is an application of the Euler-Poincaré formula. Indeed, assume that there exists $\gamma$ a simple closed quasi-transversal for $\mathcal{F}^+$. Denote by $D$ the closed disk obtained from the union of $\gamma$ and its  interior (i.e. bounded complementary region of $\gamma$ in $\mathbb{R}^2$). Notice that $\mathcal{F}^+$ defines a  foliation $F_D$ in $D$ with finitely many singularities such that: 
    
    \begin{itemize}
        \item the singularities of $F_D$ in the interior of $D$ are all $p$-prong singularities with $p\geq 3$
        \item the singularities of $F_D$ on the boundary of $D$ are all $p$-prong singularities (see Figure \ref{f.singularityboundary}) with $p\geq 1$
        \item $F_D$ is transverse to $\partial D$ except on the singularities
    \end{itemize} 

    \begin{figure}[h]
        \centering
        \includegraphics[scale=0.7]{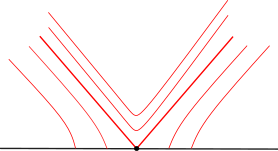}
        \caption{A 2-prong singularity in the boundary of $D$}
        \label{f.singularityboundary}
    \end{figure}

    The fact that $p\geq 1$ in the second case, follows from the third condition in our  definition of a quasi-transversal curve. Consider a second copy of $(D, F_D)$, say $(D', F_{D'})$, and glue $D$ and $D'$ along their boundaries by the identity map. We thus obtain a sphere endowed with a singular foliation $F_{\mathbb{S}^2}$ (in the sense of Definition \ref{d.singularfolisurfaces}). Notice that the singularities of $F_{\mathbb{S}^2}$ consist of the singularities of $F_D$ in the interior of $D$, the singularities of $F_{D'}$ in the interior of $D'$ and the common singularities of $F_D$ and $F_{D'}$ along the common boundary of $D$ and $D'$. The fact that the foliation $F_D$ is transverse to $\partial D$ except at a finite set of singularities, implies that $F_{\mathbb{S}^2}$ contains no other singularity except from the previous ones. Denote by $S$ the set of singularities of $F_{\mathbb{S}^2}$ and for every $s\in S$ denote the total number of prongs of $s$ by $p_s$. By the Euler-Poincaré formula (see for instance Proposition 11.4 in \cite{Primer}) we have that 
    
    \begin{equation}\label{eq.caracteristic}
        \sum_{s\in S}(2-p_s)=2\chi(\mathbb{S}^2)=4
    \end{equation}

    For any $s\in S$ that does not intersect the boundary of $D$ or $D'$, by definition of $F_{\mathbb{S}^2}$ we have that $p_s\geq 3$. For any $s\in S$ in the boundary of $D$, by construction of $F_{\mathbb{S}^2}$, we have that $s$ has a total number of prongs that is equal to two times the number of prongs of $s$ in $F_D$. This implies that $p_s\geq 2$. It follows that the left term of Equation \ref{eq.caracteristic} is negative and we thus obtained an absurd. 
\end{proof}
Thanks to the previous proposition, we will now prove several important results concerning the topological properties of the leaves of $\mathcal{F}^+$ and $\mathcal{F}^-$. 
\begin{prop}\label{p.propertiesoffoliformarkovianactions}
    \begin{enumerate}
        \item The foliations $\mathcal{F}^+$ and $\mathcal{F}^-$ do not admit compact leaves. Moreover, if $\sigma\in \Sigma$, then $\sigma$ does not admit compact $\mathcal{F}^+$-separatrices or compact $\mathcal{F}^-$-separatrices
        \item Each leaf of $\mathcal{F}^+$ or $\mathcal{F}^-$ is closed and properly embedded in $\mathbb{R}^2$. More specifically, a regular leaf of $\mathcal{F}^+$ or $\mathcal{F}^-$ is a properly embedded line in $\mathbb{R}^2$, whereas if $\sigma\in \Sigma$ is a $p$-prong singularity, both $\mathcal{F}^+(\sigma)$ and $\mathcal{F}^-(\sigma)$ consist of $p$ separatrices, each one of which is a properly embedded closed half-line in $\mathbb{R}^2$ 
        \item Any regular leaf of $\mathcal{F}^+$ or $\mathcal{F}^-$ separates $\mathbb{R}^2$ into 2 simply connected components. Similarly, if $\sigma\in \Sigma$ is a $p$-prong singularity of $\mathcal{F}^+$ and $\mathcal{F}^-$, then $$\mathbb{R}^2-\mathcal{F}^+(\sigma) ~~\big(\text{resp. }\mathbb{R}^2-\mathcal{F}^-(\sigma),~ \mathbb{R}^2-(\mathcal{F}^+(\sigma)\cup \mathcal{F}^-(\sigma))\big)$$ consists of $p$ (resp. $p$, $2p$) simply connected components, the closure of each one of which is homemorphic to a closed half-plane. 
        \item Every leaf of  $\mathcal{F}^+$ intersects at most once any leaf of $\mathcal{F}^-$
        \item Any leaf of  $\mathcal{F}^+$ or $\mathcal{F}^-$ intersects a standard polygon in $(\mathbb{R}^2, \mathcal{F}^+, \mathcal{F}^-)$ at most along one connected component

    \end{enumerate}
\end{prop}

\begin{proof}[Proof of (1)]
    We will prove the statement for $\mathcal{F}^+$. The same argument applies to $\mathcal{F}^-$. If any regular leaf of $\mathcal{F}^+$ is compact, then it would form a simple closed quasi-transversal of  $\mathcal{F}^-$, which would contradict Proposition \ref{p.noquasitransvontheplane}. Assume now that there exists $\sigma\in \Sigma$ such that $\mathcal{F}^+(\sigma)$ is compact. In this case, any $\mathcal{F}^+$-separatrix of $\sigma$ forms a simple closed curve that returns back to $\sigma$ and that is quasi-transverse to $\mathcal{F}^-$. This contradicts Proposition \ref{p.noquasitransvontheplane}.
\end{proof}

\begin{proof}[Proof of (2)]
    We will begin by proving that the leaves of $\mathcal{F}^-$ are closed and properly embedded. The same argument can naturally be applied for $\mathcal{F}^+$. Consider $L$ a leaf of $\mathcal{F}^-$ and suppose that $\text{Clos}(L)-L\neq \emptyset$. Take $x\in \text{Clos}(L)-L$. 
    
    By definition of a singular foliation, there exists $R$ a standard polygon in $\mathbb{R}^2$ containing $x$ in its interior and such that  $L\cap R$ consists of infinitely many arc connected components. If $R$ is not a rectangle, by eventually cutting $R$ into a finite number of rectangles, we can assume that there exists $R$ a rectangle containing $x$ and such that  $L\cap R$ consists of infinitely many connected components. 
    
    Denote by $\mathcal{U}$ the leaf space of the restriction of $\mathcal{F}^-$  on $R$ - $\mathcal{U}$ is homeomorphic to a segment- and endow $\mathcal{U}$ with an orientation. 
    
    \textbf{Case 1: $\mathcal{F}^-$ is orientable and $L$ is a regular leaf of $\mathcal{F}^-$}
    
    Since $\mathcal{F}^-$ is orientable, we have that there exist $L_1,L_2$ two connected components of $R\cap L$ and a continuous path $\gamma:[0,1]\rightarrow L$ such that 
    \begin{itemize}
        \item $\gamma(0)\in L_1,\gamma(1)\in L_2$ and $\gamma$ is injective
        \item the $\mathcal{F}^-$-holonomy along $\gamma$ defines an orientation preserving  map from a neighborhood of $L_1$ in $\mathcal{U}$ to a neighborhood of $L_2$ in $\mathcal{U}$
        \item $\gamma\cap R \subset L_1\cup L_2$
    \end{itemize} Using this, one can construct as in Figure \ref{f.transversecurve}, a simple closed transversal to $\mathcal{F}^-$, which contradicts Proposition \ref{p.noquasitransvontheplane}. 

\begin{figure}[h!]
    \centering
    \includegraphics[scale=0.5]{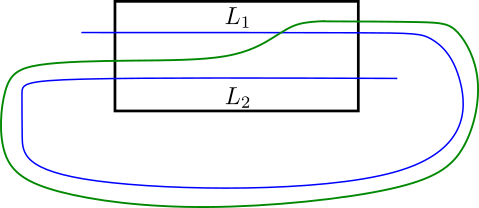}
    \caption{The green curve above represents a closed transversal to $\widetilde{F^u}$.}
    \label{f.transversecurve}
\end{figure}
    \textbf{Case 2: $\mathcal{F}^-$ is orientable and $L$ is a singular leaf of $\mathcal{F}^-$}
    
    In this case, take $\sigma\in L\cap \Sigma$ and consider $L'$ an  $\mathcal{F}^-$-separatrix of $\sigma$ intersecting $R$ infinitely many times. By applying our previous argument for $L'$, one can similarly construct a closed transversal to $\mathcal{F}^-$, which is impossible by Proposition \ref{p.noquasitransvontheplane}.

   \textbf{Case 3: $\mathcal{F}^-$ is non-orientable} 

   By lifting $\mathcal{F}^-$ on the double (ramified) orientation cover, we obtain a singular orientable foliation $\widetilde{\mathcal{F}^-}$ (in the sense of Definition \ref{d.singularfolisurfaces}) for which Proposition \ref{p.noquasitransvontheplane} still holds. Take $\widetilde{R}_1, \widetilde{R}_2$ the two lifts of the  rectangle $R$ and $\widetilde{L}$ a leaf of $\widetilde{\mathcal{F}^-}$ lifting $L$ on the orientation cover. Since $L$ intersects $R$ infinitely many times,  $\widetilde{L}$ will also intersect at least one of the rectangles  $\widetilde{R}_1, \widetilde{R}_2$ infinitely many times. By the exact same argument as in Cases 1 and 2, we conclude that this case is also impossible, which proves that the leaves of $\mathcal{F}^-$ are closed. 
   
   Next, let us show that the leaves of $\mathcal{F}^-$ are properly embedded in $\mathbb{R}^2$. Notice first that by our previous arguments, 
   
   \begin{itemize}
       \item for any regular leaf $L$ of $\mathcal{F}^-$ there cannot be a sequence of points in $L$ converging to a point in $\Sigma$ 
       \item for any $\mathcal{F}^-$-separatrix $L$ of a singular point $\sigma\in \Sigma$, if $I$ is an arbitrarily small segment in $L$ containing $\sigma$, then there cannot be a sequence of points in $L-I$ converging to a point in $\Sigma$  
   \end{itemize}
   
   By the above facts, in order to prove that the leaves of $\mathcal{F}^-$ are properly embedded in $\mathbb{R}^2$, it suffices to show that the leaves of the restriction of $\mathcal{F}^-$ on  $\mathbb{R}^2-\Sigma$ are properly embedded. This follows from the facts that $\mathcal{F}^-$ defines a non-singular foliation on $\mathbb{R}^2-\Sigma$ and that any closed leaf of a foliation in a second countable space is properly embedded (see \cite{Chevallay}). Finally, the second part of the statement of Item (2) follows immediately from our previous arguments and from Item (1). 
\end{proof}

\begin{proof}[Proof of (3)]
Any regular leaf of $\mathcal{F}^+$ or $\mathcal{F}^-$ is a properly embedded line in $\mathbb{R}^2$; hence by the Jordan curve theorem, it separates $\mathbb{R}^2$ into two connected components the closure of each one of which is a closed half-plane. 

Consider now $\sigma\in \Sigma$ a $p$-prong singularity of $\mathcal{F}^+$ and $\mathcal{F}^-$. Once again, thanks to Item (2) and the Jordan curve theorem, the union of two $\mathcal{F}^+$-separatrices or two $\mathcal{F}^-$-separatrices of $\sigma$ forms a properly embedded line in $\mathbb{R}^2$ separating $\mathbb{R}^2$ into two connected components the closure of each one of which is a closed half-plane. It follows that both $\mathbb{R}^2-\mathcal{F}^+(\sigma)$ and $\mathbb{R}^2-\mathcal{F}^-(\sigma)$ consist of $p$ connected components the closure of each one of which is a closed half-plane. 

Let us now study the connected components of $\mathbb{R}^2-(\mathcal{F}^+(\sigma)\cup\mathcal{F}^-(\sigma))$. We will begin by showing that $$\mathcal{F}^+(\sigma)\cap \mathcal{F}^-(\sigma)=\{\sigma\}$$

Indeed, consider $L^+$ an $\mathcal{F}^+$-separatrix of $\sigma$ intersecting $\mathcal{F}^-(\sigma)$ at a point different than $\sigma$. Starting from $\sigma$ and following $L^+$, consider $\sigma'$ the first intersection point between $L^+$ and $\mathcal{F}^-(\sigma)$. Denote by $[\sigma,\sigma']^+$ the unique segment in $L^+$ going from $\sigma$ to $\sigma'$ and by $[\sigma',\sigma]^-$ the unique segment in $\mathcal{F}^-(\sigma)$ going from $\sigma'$ to $\sigma$. The concatenation of the two previous segments defines a bigon $B_0$ in the plane, whose boundary consists of one segment in $L^+\subset \mathcal{F}^+(\sigma)$ containing $\sigma$ and one segment in $\mathcal{F}^+(\sigma)$ also containing $\sigma$. Notice that thanks to Item (2), by eventually considering a smaller bigon, we can assume that no leaf in $\mathcal{F}^+$ or $\mathcal{F}^-$ exits the interior of $B_0$ by intersecting $\sigma$. Also, recall that $\sigma$ is the unique singular point in $\mathcal{F}^+(\sigma)\cup\mathcal{F}^-(\sigma)$; hence $\sigma$ is the only singular point in the boundary of $B_0$. Thanks to the previous fact and to the fact that the singularities of $\mathcal{F}^+$ and $\mathcal{F}^-$ form a discrete set in $\mathbb{R}^2$, we have that $B_0$ contains a finite number of singular points. 

If the number of singular points in the interior of $B_0$ is equal to $0$, then this would contradict the fact that $\mathcal{F}^+$ and $\mathcal{F}^-$ are transverse singular foliations. It follows that there exists $\sigma_1\neq \sigma$ a singular point in the interior of $B_0$. Using Item (2), we have that $\mathcal{F}^+(\sigma_1)$ and $\mathcal{F}^-(\sigma_1)$ must intersect $\partial B_0$. In particular, two consecutive $\mathcal{F}^+$-separatrices of $\sigma_1$ intersect $\partial B_0$. By approximating the previous pair of separatrices by a regular leaf in $\mathcal{F}^+$, we construct a bigon $B_1\subset B_0$ containing strictly less singular points in its interior than $B_0$ and such that $B_1$ is bounded by one segment in a leaf of $\mathcal{F}^+$ and one segment in a leaf of $\mathcal{F}^-$. By a repeated application of this argument, we construct a bigon $B$ containing no singular points in its interior and such that $B$ is bounded by one segment in a leaf of $\mathcal{F}^+$ and one segment in a leaf of $\mathcal{F}^-$. Once again, this contradicts the fact that $\mathcal{F}^+$ and $\mathcal{F}^-$ are transverse singular foliations and thus leads to a contradiction. We deduce that $\mathcal{F}^+(\sigma)\cap \mathcal{F}^-(\sigma)=\{\sigma\}$. 

Finally, we deduce the desired result on the connected components of  $\mathbb{R}^2-(\mathcal{F}^+(\sigma)\cup\mathcal{F}^-(\sigma))$ by an argument similar to the one we used in the beginning of the proof of Item (3). 
\end{proof}
\begin{proof}[Proof of (4)]
Take $L^+\in \mathcal{F}^+$ and $L^-\in \mathcal{F}^-$ such that $L^+\cap L^-$ consists of at least two points, say $x$ and $y$. 

\textbf{Case 1: $L^+$ and $L^-$ are regular leaves of $\mathcal{F}^+$ and $\mathcal{F}^-$}

In this case, thanks to Item (2), $L^+$ and $L^-$ are properly embedded lines in $\mathbb{R}^2$. Consider $[x,y]^+$ (resp.  $[x,y]^-$ ) the unique segment in $L^+$ (resp. $L^-$) going from $x$ to $y$. The concatenation of the two previous segments defines a bigon $B_0$ in $\mathbb{R}^2$ bounded by one segment in $L^+$ and one segment in $L^-$. The bigon $B_0$ contains finitely many singular points of $\mathcal{F}^+$ and $\mathcal{F}^-$ in its interior, but does not contain any singular point in its boundary. By Item (2), if $\sigma\in \Sigma$ is a singular point in the interior of $B_0$, all its $\mathcal{F}^+$-separatrices intersect $[x,y]^-$. In particular, any two consecutive $\mathcal{F}^+$-separatrices of $\sigma$ intersect $[x,y]^-$. By approximating the previous pair of separatrices by a regular leaf in $\mathcal{F}^+$, we construct a new bigon $B_1\subset B_0$ containing strictly less singular points in its interior than $B_0$ and such that $B_1$ is bounded by one segment in a regular leaf of $\mathcal{F}^+$ and one segment in a regular leaf of $\mathcal{F}^-$. By applying this argument finitely many times, we obtain a bigon $B$ containing no singular points in its interior and such that $B$ is bounded by one segment in a regular leaf of $\mathcal{F}^+$ and one segment in a regular leaf of $\mathcal{F}^-$. This contradicts the fact that $\mathcal{F}^+$ and $\mathcal{F}^-$ are transverse singular foliations, which makes this case impossible.  

\textbf{Case 2: $L^+$ is a regular leaf of $\mathcal{F}^+$ and $L^-$ is a singular leaf of $\mathcal{F}^-$}

Denote by $\sigma$ the singular point in $L^-$ and notice that $L^+$ cannot 
intersect two distinct $\mathcal{F}^-$-separatrices of $\sigma$, because  Item (3) implies that any two $\mathcal{F}^-$-separatrices of $\sigma$ are separated by an $\mathcal{F}^+$-separatrix of $\sigma$. It follows that $L^+$ intersects at least twice the same $\mathcal{F}^-$-separatrix of $\sigma$. By approximating the previous separatrix by a regular leaf of $\mathcal{F}^-$, we go back to Case 1; hence, this case is also impossible. 

Naturally, for the same exact reason, the case where $L^+$ is a singular leaf of $\mathcal{F}^+$ and $L^-$ is a regular leaf of $\mathcal{F}^-$ is also impossible. It therefore remains to consider the following case: 

\textbf{Case 3: $L^+$ and $L^-$ are singular leaves of $\mathcal{F}^+$ and $\mathcal{F}^-$}
Let $\sigma^+,\sigma^-\in \Sigma$ be the unique singular points contained in $L^+$ and $L^-$ respectively. Notice that we have already proven during the proof of Item (3) that the points $\sigma^+$ and $\sigma^-$ can not coincide. Notice also that because of Item (3), two distinct $\mathcal{F}^+$-separatrices of $\sigma^+$ can not both intersect $L^-$. Similarly, two distinct $\mathcal{F}^-$-separatrices of $\sigma^-$ can not intersect both $L^+$. It follows that there exist an $\mathcal{F}^+$-separatrix of $\sigma^+$ and an $\mathcal{F}^-$-separatrix of $\sigma^-$ intersecting twice in $\mathbb{R}^2$. By approximating the previous separatrices by regular leaves in $\mathcal{F}^+$ and $\mathcal{F}^-$, we go back to Case 1; hence, this case is also impossible, which finishes the proof of Item (4).

\end{proof}
\begin{proof}[Proof of (5)]
 Consider $R$ a standard polygon in $(\mathbb{R}^2, \mathcal{F}^+, \mathcal{F}^-)$ and $L^+$ a leaf in $\mathcal{F}^+$. If $R$ is a rectangle and $L$ intersects $R$ along at least two connected components, then $L^+$ would intersect at least twice any leaf in $\mathcal{F}^-$ intersecting $R$, which is impossible because of Item (4). 
 
 Suppose now that $R$ is a standard $2p$-gon and let $\sigma$ be the $p$-prong singular point of $\mathcal{F}^+$ in $R$. Cut $R$ along the quadrants of $\sigma$ (see Definition \ref{d.segment}), in order to obtain $2p$ rectangles $R_1,...,R_{2p}$. By our previous arguments, $\mathcal{F}^+(\sigma)$ intersects each of these rectangles along a unique segment. Thanks to Item (3), this implies that the closure of any connected component of $\mathbb{R}^2-\mathcal{F}^+(\sigma)$ contains exactly two rectangles among $R_1,...,R_{2p}$. Denote by $R_k$ and $R_l$ the two rectangles in  $\{R_1,...,R_{2p}\}$ that are contained in the closure of the same connected component of $\mathbb{R}^2-\mathcal{F}^+(\sigma)$ as $L^+$. Notice that, if a leaf in $\mathcal{F}^+$ intersects $R_k$ then it also intersects $R_l$ and vice versa. Using this fact, as $L^+$ intersects $R$ along two or more connected components, we get that $L^+$ also intersects $R_k$ along two or more connected components, which is impossible by our previous arguments.   
\end{proof}

\section{Preliminaries}\label{s.preliminaries}
\subsection{Markov partitions and geometric types} \label{s.markovpartitionsgeomtypesdefi}
Let $M$ be a closed 3-manifold carrying a pseudo-Anosov flow $\Phi$. Let $F^s$ and $F^u$ be the stable and unstable foliations of $\Phi$. 
\begin{defi}
Consider $R$ a transverse rectangle of $\Phi$ (see Definition \ref{d.standardpolygonflow}). The rectangle $R$ is trivially bifoliated by stable and unstable segments. We will denote by $\partial^{s}R$ (resp. $\partial^{u}R$) the union of the two stable (resp. unstable) segments in $\partial R$, which we will call the \emph{stable boundary} (resp. \emph{unstable boundary}) of $R$. 

Furthermore, any subrectangle of $R$, say $R'$, for which $\partial^s R'\subset \partial^s R$ will be called a \emph{vertical subrectangle} of $R$. Similarly, any subrectangle $R'$ of $R$ for which $\partial^u R'\subset \partial^u R$ will be called a \emph{horizontal  subrectangle} of $R$.
\end{defi}

\begin{defi}\label{d.markovpartition}
A \emph{Markov partition} of $\Phi$ is a finite family of transverse rectangles of $\Phi$, say $R_1,...,R_n$, such that: 
\begin{enumerate}
\item The rectangles are pairwise disjoint
\item The (striclty) positive orbit by $\Phi$ of any point $x\in  \underset{i=1}{\overset{n}{\cup}}R_i$ intersects $\underset{i=1}{\overset{n}{\cup}}R_i$ along a point that we will denote by $f(x)$. Furthermore, there exists $T>0$ such that every orbit segment of $\Phi$ of length $T$ intersects $\underset{i=1}{\overset{n}{\cup}}R_i$. 
\item For any two $i,j$ the closure of each  connected component of  $f(\inte{R_i})\cap \inte{R_j}$  (the previous set can be empty) is a vertical subrectangle of $R_j$. 
\item For any two $i,j$ the closure of each connected component of  $f^{-1}(\inte{R_i})\cap \inte{R_j}$ (the previous set can be empty) is a horizontal subrectangle of $R_j$. 
\end{enumerate}
Furthermore, we will call the family $R_1,...,R_n$ a \emph{reduced} Markov partition if for every $i\neq j$, there does not exist a continuous function $\tau: R_i \rightarrow \mathbb{R}$ such that $\Phi^{\tau}(R_i)\subseteq R_j$ or $\Phi^{\tau}(R_i)\supseteq R_j$. 
\end{defi}
\begin{rema}\label{r.defmarkovpartition}
\begin{itemize}
    \item The compactness of $M$ and point $(2)$ of the above definition imply that the number of connected components of $f(\inte{R_i})\cap \inte{R_j}$  or $f^{-1}(\inte{R_i})\cap \inte{R_j}$ is finite for any $i,j$.
    \item Since $f$ preserves the stable and unstable segments of  $\underset{i=1}{\overset{n}{\cup}}\inte{R_i}$, for any Markov partition one can generalize by induction the points (3) and (4) of the above definition to the following statement: for every $i,j\in \llbracket 1,n\rrbracket$,  $N\geq 0$ (resp. $N\leq 0$) and $V$ vertical (resp. horizontal) subrectangle of $R_i$ the set $f^N(\inte{V})\cap \inte{R_j}$ consists of a finite number of connected components the closure of each one of which is a vertical (resp. horizontal) subrectangle of $R_j$.
\end{itemize} 
\end{rema}
 The reason why we use the interiors of the rectangles $R_1,...,R_n$ in points (3) and (4) of Definition \ref{d.markovpartition}, instead of the rectangles themselves, is that $f$ often fails to be continuous around a point in $\underset{i=1}{\overset{n}{\cup}}\partial R_i$, whereas  for any $x\in f^{-1}(\inte{R_i})\cap R_j$ there exists $U$ a small neighborhood of $x$ in $\inte{R_j}$ such that $f(U)\subset \inte{R_i}$ and $f_{|U}$ is continuous. One can easily deduce from this fact that 
 \begin{rema}\label{r.defmarkovpartition1}
\begin{itemize}
\item the image by $f$ of every connected component of $f^{-1}(\inte{R_i})\cap \inte{R_j}$ is a connected component of $f(\inte{R_j})\cap \inte{R_i}$
\item $ f(\underset{i=1}{\overset{m}{\cup}}\partial^sR_i)\cap \underset{i=1}{\overset{m}{\cup}}\inte{R_i}=\emptyset \text{ and } f^{-1}(\underset{i=1}{\overset{m}{\cup}}\partial^uR_i)\cap \underset{i=1}{\overset{m}{\cup}}\inte{R_i}=\emptyset$
\end{itemize}
\end{rema}
The following two results together with Proposition \ref{p.densitystablemanifoldsperiodic}  ensure the existence of Markov partitions and reduced Markov partitions for all pseudo-Anosov flows.  
\begin{theorem}[\cite{markovpseudoanosov}] \label{t.existenceofmarkovpartitions}
Let $\Phi$ be any pseudo-Anosov flow on a closed 3-manifold $M$ and $F^s,F^u$ its stable and unstable foliations. For any finite set of periodic orbits $\Gamma$ of $\Phi$ such that $\Gamma$ contains every singular periodic orbit of $\Phi$ and such that $\underset{\gamma\in\Gamma}{\cup}F^s(\gamma)$ and  $\underset{\gamma\in\Gamma}{\cup}F^u(\gamma)$ are dense in $M$, there exists a Markov partition of $\Phi$ formed by rectangles, whose stable and unstable boundaries are contained respectively in $\underset{\gamma\in\Gamma}{\cup}F^s(\gamma)$ and  $\underset{\gamma\in\Gamma}{\cup}F^u(\gamma)$.
\end{theorem}
The existence of infinitely many sets of periodic orbits satisfying the above hypothesis was shown in \cite{markovpseudoanosov}. 

\begin{theorem}\label{t.reducedmarkovpartitionexists}
Let $\Phi$ be any pseudo-Anosov flow on a closed 3-manifold $M$ and $F^s,F^u$ its stable and unstable foliations. For any finite set of periodic orbits $\Gamma$ of $\Phi$ such that $\Gamma$ contains every singular periodic orbit of $\Phi$ and such that $\underset{\gamma\in\Gamma}{\cup}F^s(\gamma)$ and  $\underset{\gamma\in\Gamma}{\cup}F^u(\gamma)$ are dense in $M$, there exists a reduced Markov partition of $\Phi$ formed by rectangles, whose stable and unstable boundaries are contained respectively in $\underset{\gamma\in\Gamma}{\cup}F^s(\gamma)$ and  $\underset{\gamma\in\Gamma}{\cup}F^u(\gamma)$.
\end{theorem}
\begin{proof}
Same argument as in Theorem 1.3.16 of \cite{thesisioannis}.

\end{proof}

To every Markov partition $R_1,...R_n$ of $\Phi$ we can associate (this association will be explained in a more general and detailed way in Definition \ref{d.geometrictypemarkovfamily}) a family of finite combinatorial objects, called \emph{geometric types}, that describe the way in which the first return map acts on $R_1,...,R_n$. Geometric types constitute one of the most important objects in this paper, thanks to which we will later classify pseudo-Anosov flows up to Dehn-Goodman-Fried surgeries. 
\begin{defi}\label{d.geometrictype}
Take $n\in \mathbb{N}^*$ and  $(h_i)_{i\in \llbracket 1, n\rrbracket},(v_i)_{i\in \llbracket 1, n\rrbracket} \in (\mathbb{N}^*)^n$ such that $$\sum_i h_i =\sum_i v_i$$ 

Consider now for every $i\in \llbracket 1,n \rrbracket$ two (abstract) finite sets of the form $\{H_i^j, j\in \llbracket 1, h_i\rrbracket\}$ and $\{V_i^j, j\in \llbracket 1, v_i\rrbracket\} $, a bijection $\phi$ between  $$\mathcal{H}:=\lbrace H^j_i| i\in\llbracket 1,n \rrbracket, j \in \llbracket 1,h_i \rrbracket  \rbrace  \text{ and } \mathcal{V}:=\lbrace V^j_i| i\in\llbracket 1,n \rrbracket, j \in \llbracket 1,v_i \rrbracket  \rbrace$$ and a function $u$ from $\mathcal{H}$ to $\lbrace -1,+1\rbrace$.

The data $(n,(h_i)_{i \in \llbracket 1,n \rrbracket}, (v_i)_{i\in \llbracket 1,n \rrbracket}, \mathcal{H}, \mathcal{V},\phi, u)$ will be called a \emph{geometric type}. 
\end{defi}

\textbf{A geometric interpretation of the geometric type} 

Let $G=(n,(h_i)_{i \in \llbracket 1,n \rrbracket}, (v_i)_{i\in \llbracket 1,n \rrbracket}, \mathcal{H}, \mathcal{V},\phi, u)$ be a geometric type with $\mathcal{H}=\lbrace H^j_i| i\in\llbracket 1,n \rrbracket, j \in \llbracket 1,h_i \rrbracket  \rbrace  \text{ and } \mathcal{V}=\lbrace V^j_i| i\in\llbracket 1,n \rrbracket, j \in \llbracket 1,v_i \rrbracket  \rbrace$. 

Consider $n$ copies of $[0,1]^2$, say $R_1,...,R_n$, trivially bifoliated by horizontal and vertical segments (oriented from bottom to top and left to right respectively). Consider also inside every $R_i$, $h_i$ (resp. $v_i$) 
mutually disjoint horizontal (resp. vertical) subrectangles, that we are going to denote by $\tilde{H}_i^1,...,\tilde{H}_i^{h_i}$ (resp. $\widetilde{V}_i^1,...,\widetilde{V}_i^{v_i}$). 

\begin{conv}\label{c.drawinggeomtype}
    Throughout this paper, whenever we use the geometric version of a geometric type, we will assume that the subrectangles $\widetilde{H}_i^1,...,\widetilde{H}_i^{h_i}$ (resp. $\widetilde{V}_i^1,...,\widetilde{V}_i^{v_i}$) are ordered inside $R_i$ from bottom to top (from left to right). 
\end{conv}
Let $\widetilde{\phi}: \underset{i\in \llbracket 1,n\rrbracket}{\cup}\underset{j\in \llbracket 1,h_i\rrbracket}{\cup}\widetilde{H}^j_i \rightarrow \underset{i\in \llbracket 1,n\rrbracket}{\cup}\underset{j\in \llbracket 1,v_i\rrbracket}{\cup}\widetilde{V}^j_i$ be a homeomorphism :
\begin{itemize}
    \item sending in an affine way each horizontal subrectangle $\widetilde{H}_i^j$ to the vertical subrectangle $\widetilde{V}^{j'}_{i'}$, where $i',j'$ are such that $\phi(H_i^j)=V^{j'}_{i'}$
    \item respecting (resp. reversing) the orientation of the horizontal and vertical foliation when restricted to  $\widetilde{H}_i^j$ if $u(H_i^j)=1$ (resp. $u(H_i^j)=-1$)
\end{itemize}    
We will call $(R_1,...,R_n, (\widetilde{H}_i^j)_{i\in\llbracket 1,n \rrbracket, j \in \llbracket 1,h_i \rrbracket}, (\widetilde{V}_i^j)_{i\in\llbracket 1,n \rrbracket, j \in \llbracket 1,v_i \rrbracket}, \tilde{\phi})$ the \emph{geometrisation of $G$}. 

\begin{conv}
    Throughout this paper, we will not distinguish between geometric types and their geometrisations. A geometric type will be thought at the same time as a combinatorial object, but also as a set of rectangles endowed with dynamics. 
\end{conv}

\textbf{Equality and equivalence of geometric types} 

\begin{figure}[h!]
\includegraphics[scale=0.75]{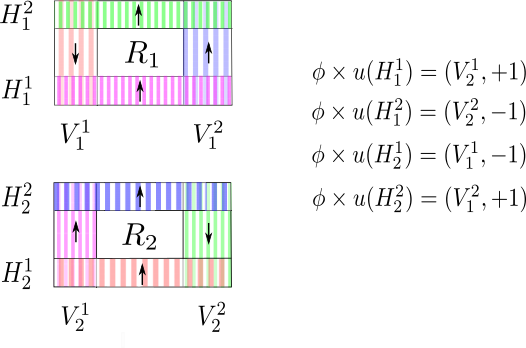}
\caption{In the above example the ``first return map" sends a horizontal rectangle to a vertical rectangle of the same color. The arrows indicate how the ``first return map" acts on the orientation of the vertical foliation.}
\label{f.examplegeometrictype}
\end{figure}

 A geometric type $(n,(h_i)_{i \in \llbracket 1,n \rrbracket}, (v_i)_{i\in \llbracket 1,n \rrbracket}, \mathcal{H}, \mathcal{V},\phi, u)$ is a purely combinatorial  version of a Markov partition, where the function $\phi$  plays the role of the first return map and the function $u$ indicates whether the first return map preserves or reverses the orientation of the horizontal and vertical foliations (see Figure \ref{f.examplegeometrictype}). 

\begin{defi}\label{d.equivalentgeomtypes}
Take $(n,(h_i)_{i \in \llbracket 1,n \rrbracket}, (v_i)_{i\in \llbracket 1,n \rrbracket}, \mathcal{H}, \mathcal{V},\phi,u)$ and $(n',(h'_i)_{i \in \llbracket 1,n' \rrbracket}, (v'_i)_{i\in \llbracket 1,n' \rrbracket}, \mathcal{H}', \mathcal{V}',\phi', u')$ two geometric types. Suppose that
\begin{align*}
    \mathcal{H}=\{H_{i}^j, i\in \llbracket 1,n \rrbracket, j\in \llbracket 1, h_i \rrbracket\} \quad & \text{and} \quad  \mathcal{V}=\{V_{i}^{j}, i\in \llbracket 1,n \rrbracket, j\in \llbracket 1, v_i \rrbracket\} \\  \mathcal{H}'=\{{H'_{i}}^{j}, i\in \llbracket 1,n' \rrbracket, j\in \llbracket 1, h'_i \rrbracket\} \quad & \text{and} \quad  \mathcal{V}'=\{{V'_{i}}^{j}, i\in \llbracket 1,n' \rrbracket, j\in \llbracket 1, v'_i \rrbracket\}
\end{align*} We will say that the two previous geometric types are \emph{equivalent} if 

    \begin{enumerate}
    \item $n=n'$
    \item up to reindexing the pairs $((h_i,v_i))_{i \in \llbracket 1,n \rrbracket}$ we have that $(h_i,v_i)=(h_i', v_i')$ for every $i\in \llbracket 1,n \rrbracket$
    \item there exists a bijection $H: \mathcal{H}\cup \mathcal{V} \rightarrow \mathcal{H}'\cup \mathcal{V}'$ such that
    \begin{itemize}
        \item $H(\mathcal{H})=\mathcal{H}'$ and $H$ respects the order of the elements in $\mathcal{H}$ and $ \mathcal{H}'$. In other words, for every $i\in \llbracket 1, n\rrbracket$, $H$ defines a bijection between $\{H_{i}^j, j\in \llbracket 1,h_i \rrbracket \}$ and $\{{H'_{i}}^{j}, j\in \llbracket 1,h'_i \rrbracket\}$ that is monotonous with respect to $j$. We define $\epsilon_i=+1$ if the previous map is increasing and $\epsilon_i=-1$ if not.
        \item $H(\mathcal{V})=\mathcal{V}'$ and $H$ respects the order of the elements in $\mathcal{V}$ and $ \mathcal{V}'$. In other words, for every $i\in \llbracket 1, n\rrbracket$, $H$ defines a bijection between $\{V_i^j, j\in \llbracket 1,v_i \rrbracket \}$ and $\{{V'_i}^{j}, j\in \llbracket 1,v'_i \rrbracket\}$ that is monotonous with respect to $j$. We define $\epsilon'_i=+1$ if the previous map is increasing and $\epsilon'_i=-1$ if not.
        \item 
        either  $\epsilon_i\cdot \epsilon_i'=-1$ for all $i$ or  $\epsilon_i\cdot \epsilon_i'=+1$ for all $i$
        
        \item $H$ respects $\phi$ and $\phi'$. In other words, for every $h\in \mathcal{H}$ $\phi'(H(h))= H(\phi(h))$ 
        \item for every $h\in \mathcal{H}$ we have $u'(H(h))=\epsilon_i \cdot \epsilon_j \cdot u(h)=\epsilon'_i \cdot \epsilon'_j \cdot u(h)$, where $i,j$ are such that $h$ and $\phi(h)$ are respectively of the form $H_i^{ \bullet}$ and $V_j^{ \bullet}$
        \end{itemize}
    \end{enumerate}

If furthermore all the above $\epsilon_i$ and $\epsilon_i'$ are equal to $+1$, then we will say that the two geometric types are \emph{equal}.  
\end{defi}

\textbf{A geometric interpretation of the relation of equivalence.} Using the geometric interpretation of a geometric type, an equivalence $H$ between two geometric types can be thought as a homeomorphism between the rectangles of the two geometric types sending horizontal/vertical subrectangles to horizontal/vertical subrectangles and respecting the maps $\phi$ and $u$. Keep in mind, that $H$ can a priori change the orientation of the vertical/horizontal foliations of a rectangle.

For instance, changing the orientation of the stable and unstable foliations of the rectangle $R_1$ in the example of Figure \ref{f.examplegeometrictype} yields an equivalent geometric type. Indeed, recall that the horizontal and vertical subrectangles of $R_1$, namely $H_1^1,H_1^{2}$ and $V_1^1,V_1^{2}$, are respectively ordered from bottom to top and from left to right. Changing the orientation of the foliations of $R_1$  amounts to reindexing the $H_1^j$ and the $V_1^j$ and hence leads to a different geometric type, where $\phi$ and $u$ are no longer the same. It is easy to check that the ``identity map" sending trivially $R_1$ and $R_2$ to themselves satisfies the axioms of Definition \ref{d.equivalentgeomtypes}, thus defining an equivalence between the two (non equal) geometric types.

More generally, it is not difficult to see that given a geometric type $$G=(n,(h_i)_{i \in \llbracket 1,n \rrbracket}, (v_i)_{i\in \llbracket 1,n \rrbracket}, \mathcal{H}, \mathcal{V},\phi, u)$$  reindexing the pairs $(h_i,v_i)$ or changing the orientation of the stable or unstable foliation in every rectangle of $G$ or changing the orientation of both  foliations of any rectangle of $G$, produces a geometric type equivalent to $G$ and even more by a finite repetition of such operations we can produce every geometric type (up to equality) that is equivalent to $G$. 

Finally, if we consider geometric types up to equality, we have that 
\begin{rema}\label{r.finitenumbergeometrictypes}
The relation of equivalence between geometric types is an equivalence relation for which every equivalence class of geometric types is finite.
\end{rema}

\textbf{A technical remark.} Once again, our definition of equivalence is a purely combinatorial one, but throughout this paper we will be often thinking of an equivalence between two geometric types in a geometric way. Our previous geometric interpretation explains the reason why we defined the notion of equivalence in this way. There is though one axiom in Definition \ref{d.equivalentgeomtypes} that we have not yet thoroughly explained. The reason why we added the condition ``either  $\epsilon_i\cdot \epsilon_i'=-1$ for all $i$ or  $\epsilon_i\cdot \epsilon_i'=+1$ for all $i$" is a technical one and is closely related to the fact that in this paper we will use geometric types in order to classify pseudo-Anosov flows on \emph{orientable} 3-manifolds. Should one wish to generalize our approach for non-orientable manifolds, our definition of a geometric type should be modified and also one should drop the previous axiom from our definition of equivalence. 

\subsection{Markovian families of a pseudo-Anosov flow}
In this section, we will define the notion of Markovian family in the bifoliated plane $\mathcal{P}$ of a pseudo-Anosov flow $\Phi$ and we will prove that general Markovian families in $\mathcal{P}$ share several properties with reduced Markov partitions of $\Phi$. The existence of Markovian families (see Proposition \ref{p.projectionmarkovpartition}) for any pseudo-Anosov flow in dimension $3$ constitutes the starting point of the classification approach developped in this paper. 

Let $M$ be a closed, connected 3-manifold, $\Phi$ a pseudo-Anosov flow on $M$ and $F^s,F^u$ the stable and unstable manifolds of $\Phi$. Consider $\mathcal{P}$ the bifoliated plane of $\Phi$ endowed with its stable and unstable foliations $\mathcal{F}^s, \mathcal{F}^u$ and with its natural action by $\pi_1(M)$ preserving the foliations $\mathcal{F}^s$ and $\mathcal{F}^u$ (see Theorem-Definition \ref{thdef.bifoliatedplane}). 

\begin{conv}
    In order to differentiate between the action of $\pi_1(M)$ on the universal cover of $M$, say $\widetilde{M}$, by deck transformations and its induced action on the orbit space $\mathcal{P}$, throughout this section we will use the following notations: 
    \newline{}
\begin{minipage}[b]{0.5\linewidth}
\begin{align*}
&(\pi_1(M),\mathcal{P})\rightarrow \mathcal{P}\\
& \quad (g,x)\rightarrow g(x)
\end{align*}
\end{minipage}
\begin{minipage}[b]{0.5\linewidth}
\begin{align*}
&(\pi_1(M),\widetilde{M})\rightarrow \widetilde{M}\\
& \quad (g,x)\rightarrow g.x
\end{align*}
\end{minipage}
Recall that the projection from $\widetilde{M}$ to $\mathcal{P}$ is equivariant with respect to the previous two actions. 
\end{conv}
\begin{defi}\label{d.subrectangleplane}
Let $R$ be a rectangle in $\mathcal{P}$ (see Definition \ref{d.standardpolygon}). We will call the union of the two stable segments in the boundary of $R$ the \emph{stable boundary} of $R$ and we will denote it by $\partial^sR$. We similarly define $\partial^uR$, the \emph{unstable boundary} of $R$. 

Moreover, any subrectangle $R'$ of $R$ such that $\partial^sR'\subset \partial^sR $ will be called a \emph{vertical subrectangle} of $R$. If furthermore $R' \neq R$, then we will call $R'$ a \emph{non-trivial} vertical subrectangle of $R$. We define similarly horizontal subrectangles and non-trivial horizontal subrectangles.

\end{defi}

\begin{defi}\label{d.markovfamily}
A \emph{Markovian family} in $\mathcal{P}$ is a $\pi_1(M)$-invariant set of  rectangles $(R_i)_{i \in I}$ covering $\mathcal{P}$ (i.e. $\cup_{i \in I} R_i = \mathcal{P}$) such that 
\begin{enumerate}
\item (Finiteness axiom) $(R_i)_{i \in I}$ is the union of a finite number of orbits of rectangles of the action by $\pi_1(M)$ (i.e. there exist $i_1,...,i_n\in I$ such that for every $R\in (R_i)_{i \in I}$ there exist $g\in \pi_1(M)$ and $p\in \llbracket 1, n\rrbracket$ such that $g(R_{i_p})=R$) 
\item (Markovian intersection axiom) For every two distinct rectangles $R_i,R_j$ in $(R_i)_{i \in I}$, if $\overset{\circ}{R_i } \cap \overset{\circ}{R_j} \neq \emptyset$, then, up to changing the roles of $R_i$ and $R_j$, $R_i \cap R_j$ is a non-trivial horizontal subrectangle of $R_i$ and a non-trivial vertical subrectangle of $R_j$ 
\item (Finite return time axiom) Take any point $x\in \mathcal{P}$ and $U$ any compact neighborhood of $x$ in $\mathcal{P}$. The neighborhood $U$ can be covered a finite number of rectangles of the family $(R_i)_{i \in I}$ 
\end{enumerate} 
\end{defi}
Thanks to Theorem \ref{t.reducedmarkovpartitionexists} and to the following proposition, Markovian families exist in the bifoliated plane of any pseudo-Anosov flow: 
\begin{prop}\label{p.projectionmarkovpartition}
Take $R_1,R_2,...,R_n$ a reduced Markov partition of $\Phi$. Consider $(\widetilde{R_i})_{i\in I}$ the family of lifts of the $R_i$ on $\widetilde{M}=\mathbb{R}^3$ and $(\overline{R_i})_{i\in I}$ the projection of $(\widetilde{R_i})_{i\in I}$ on $\mathcal{P}$. The set $(\overline{R_i})_{i\in I}$ is a Markovian family. 
\end{prop}
\begin{proof}
We will check that $(\overline{R_i})_{i\in I}$ satisfies the axioms of Definition \ref{d.markovfamily}. Denote by $\widetilde{\Phi}$ the lift of $\Phi$ on $\widetilde{M}$. First, since the family of rectangles $(\widetilde{R_i})_{i\in I}$ intersects every orbit of the lifted flow $\widetilde{\Phi}$ and is $\pi_1(M)$-invariant (with respect to the action of $\pi_1(M)$ by deck transformations), we get that the family $(\overline{R_i})_{i\in I}$ covers $\mathcal{P}$ and is also $\pi_1(M)$-invariant (with respect to the natural action of $\pi_1(M)$ on $\mathcal{P}$). Moreover, since the rectangles in  $(\widetilde{R_i})_{i\in I}$ are transverse to $\widetilde{\Phi}$,   Corollary \ref{c.orbitintersectpolygon} implies that the $\overline{R_i}$ are rectangles in $\mathcal{P}$.

\vspace{0.2cm}
\textit{Finiteness axiom} 

The family of rectangles  $(\overline{R_i})_{i\in I}$ satisfies  the finiteness axiom, since by construction it consists of the projections on $\mathcal{P}$ of the lifts on $\widetilde{M}$ of a finite number of rectangles in $M$. 

\textit{Markovian intersection axiom} 

Take two distinct $i,j \in I$ such that $\overset{\circ}{\overline{R_i} } \cap \overset{\circ}{\overline{R_j}} \neq \emptyset$. Consider the lifts $\widetilde{R_i},\widetilde{R_j}\in (\widetilde{R_l})_{l\in I}$ of  $\overline{R_i},  \overline{R_j}$ on $\widetilde{M}$. Denote by $A$ (resp. $B$) the set of points of $\inte{\widetilde{R_i}}$, whose positive (resp. negative) orbits by $\widetilde{\Phi}$ intersect $\inte{\widetilde{R_j}}$. The fact that $\overset{\circ}{\overline{R_i} } \cap \overset{\circ}{\overline{R_j}} \neq \emptyset$ implies that $A\cup B\neq \emptyset$. Without any loss of generality, let us assume that $A\neq \emptyset$. Let us show that the closure of $A$ is a non-trivial horizontal subrectangle of $\widetilde{R_i}$ and that $B=\emptyset$.

Consider $A_k$ the set of points in $A\subset \inte{\widetilde{R_i}}$ whose positive orbits intersect $k\in \mathbb{N}$ rectangles in $(\widetilde{R_l})_{l\in I}$ before intersecting $\widetilde{R_j}$. Take $k\in \mathbb{N}$ such that $A_k\neq \emptyset$ and $\widetilde{x}\in A_k$. Denote by $R_i,R_j\in (R_l)_{l\in \llbracket 1, n\rrbracket}$ the projections of $\widetilde{R_i}$, $\widetilde{R_j}$ in $M$, by $x_M$ the projection of $\widetilde{x}$ on $R_i$ and by $f: \underset{l=1}{\overset{n}{\cup}}R_l \rightarrow \underset{l=1}{\overset{n}{\cup}}R_l$ the first return map on $\underset{l=1}{\overset{n}{\cup}}R_l$. By Remark \ref{r.defmarkovpartition},  $f^{-(k+1)}(\inte{R_j})\cap \inte{R_i}$ consists of a finite number of connected components the closure of each one of which is a horizontal subrectangle of $R_i$.

Consider $h_M$ the closure of the connected component of $f^{-(k+1)}(\inte{R_j})\cap \inte{R_i}$ containing $x_M$. Since, by Remark \ref{r.defmarkovpartition1}, $f$ (resp. $f^{-1}$) preserves $\underset{l=1}{\overset{n}{\cup}}\partial^s R_l$ (resp. $\underset{l=1}{\overset{n}{\cup}}\partial^u R_l$) and since $\inte{h_M}\subset \inte{R_i}$, $f^{(k+1)}(\inte{h_M})\subset \inte{R_j}$, we have that for every $m\in \llbracket 1,k+1\rrbracket$, $f^{m}(\inte{h_M})\cap \underset{l=1}{\overset{n}{\cup}}\partial R_l= \emptyset$. Hence, for every $m\in \llbracket 1,k+1\rrbracket$, 
\begin{enumerate}
    \item by our discussion prior to Remark \ref{r.defmarkovpartition1}, the $m$-th return map $f^m$ is continuous on $\inte{h_M}$, as a composition of $m$ continuous functions
    \item $f^{m}(\inte{h_M})$ belongs in a unique rectangle of our Markov partition $R_1,...,R_n$
    \item $f^{(k+1)}(\inte{h_M})$ is a connected component of $f^{(k+1)}(\inte{R_i})\cap \inte{R_j}$ and thus by Remark \ref{r.defmarkovpartition} its closure is a vertical subrectangle of $R_j$
\end{enumerate}  

 Consequently, by lifting $h_M$ on the universal cover, we get that $\widetilde{x}$ belongs in a horizontal subrectangle $h$ of $\widetilde{R_i}$ whose interior is a subset of $A_k$ and whose positive orbit intersects $\widetilde{R_j}$ along a vertical subrectangle. By repeating this argument a finite number of times, we get that the closure of $A_k$ can be written as the union of a finite number of horizontal subrectangles of $\widetilde{R_i}$ with pairwise disjoint interiors. This number is necessarily equal to one since any stable leaf in $\widetilde{F^s}$ intersecting  $\widetilde{R_j}$ can not intersect $\widetilde{R_i}$ along more than one connected component (see Corollary \ref{c.complementsingleaves}). For the same reason, $k$ is the unique positive integer for which $A_k$ is non-empty. We thus have proven that the closure of $A$ is a non-trivial horizontal subrectangle of $\widetilde{R_i}$. 

 If $B\neq \emptyset$, by the same argument, $B$ would be a non-trivial vertical subrectangle of $\widetilde{R_i}$. Therefore, $A\cap B\neq \emptyset$, which would imply the existence of an orbit of $\widetilde{\Phi}$ intersecting more than once $\overline{R_j}$, which is impossible by Theorem \ref{t.oneintersection}. Finally, by projecting $A$ on the bifoliated plane we get that $\overline{R_i} \cap \overline{R_j}$ is a non-trivial horizontal subrectangle of $\overline{R_i}$ and a non-trivial vertical subrectangle of $\overline{R_j}$, which gives us the desired result.

\vspace{0.2cm}
\textit{Finite return time axiom} 

Take $x\in \mathcal{P}$ and a compact neighborhood $U$ of $x$ in $\mathcal{P}$. Assume without any loss of generality that $U$ is homeomorphic to a closed disk and consider $\widetilde{U}$ a lift of $U$ on $\mathbb{R}^3$ that is also homeomorphic to a closed disk. Since the positive orbit by $\Phi$ of any point in $M$ intersects a rectangle in $R_1,...,R_n$ in a uniformly bounded time, the positive orbit by $\widetilde{\Phi}$ of any point in $\widetilde{U}$ also intersects a rectangle in $(\widetilde{R_i})_{i\in I}$ in a uniformly bounded time, say $T$. Moreover, since the action of $\pi_1(M)$ on $\mathbb{R}^3$ is properly discontinuous, for any compact set $K\subset \mathbb{R}^3$ only a finite number of rectangles in $(\widetilde{R_i})_{i\in I}$ can intersect $K$. By taking $K=\underset{t\in [0,T]}{\cup}\widetilde{\Phi}(\widetilde{U})$, we get that there exists a finite family of rectangles in $(\widetilde{R_i})_{i\in I}$ such that the positive orbit of any point in $\widetilde{U}$ intersects a rectangle in this family in a uniformly bounded time. This proves the desired result and finishes the proof of the proposition. 

\end{proof}

Similarly to the case of Anosov flows (see \cite{preprint}), as we will show later in this paper, a general Markovian family $\mathcal{R}$ shares many properties with the Markovian families constructed in Proposition \ref{p.projectionmarkovpartition}: the boundaries of all the rectangles in $\mathcal{R}$ belong to stable/unstable leaves of $\mathcal{F}^s, \mathcal{F}^u$ with non-trivial stabilizers in $\pi_1(M)$, every point in $\mathcal{P}$ belongs to infinitely many rectangles of the family, etc.  Because of the abundance of such similarities the reader may choose to think of a Markovian family as the analogue of a Markov partition inside the bifoliated plane of a pseudo-Anosov flow and also we conjecture that:  
\begin{conj}\label{conj.markov}
Every Markovian family of $\Phi$ corresponds to the projection on $\mathcal{P}$ of the family of lifts on $\mathbb{R}^3$ of the rectangles of a reduced Markov partition of $\Phi$.
\end{conj}

Before finishing this section, let us prove two important lemmas concerning the Markovian families of $\Phi$ and providing evidence that support Conjuecture \ref{conj.markov}. 

\begin{defi}
Consider $\mathcal{R}$ a Markovian family of $\Phi$. For any $x\in \mathcal{P}$, we will call the closure of any connected component of $\mathcal{F}^s(x)-\{x\}$ a \emph{stable separatrix of $x$}. We define similarly unstable separatrices.    

Take $x\in \mathcal{P}$. By Proposition \ref{p.propertiesoffoliformarkovianactions}, $\mathcal{P}-(\mathcal{F}^s(x)\cup \mathcal{F}^u(x))$ consists of a finite number of connected components, the closure of each one of which is homeomorphic to a closed half-plane. Each of the previous closed half-planes will be called a \emph{quadrant of $x$ in $\mathcal{P}$}. A neighborhood of $x$ inside a quadrant of $x$ in $\mathcal{P}$ will be called a \emph{germ of quadrant of $x$}.

\end{defi}

\begin{lemm}\label{l.vertsubrectangleexists}
Take $R\in\mathcal{R}$, $x\in R$ and  $G$ a sufficiently small  germ of quadrant of $x$ contained in $R$. There exists $R_v\in \mathcal{R}$ such that $G\subset R_v$ and $R_v\cap R$ is a non-trivial vertical subrectangle of $R$.

\end{lemm}
\begin{proof}
Take $x,R, G$ as in the statement. Denote by $S,U$ the stable and unstable separatrices of $x$ delimiting the quadrant of $x$ in $\mathcal{P}$ containing $G$. Let $\widetilde{\Phi}$ be the lift of $\Phi$ on the universal cover $\widetilde{M}$ and $\pi$ the natural projection from $\widetilde{M}$ to the orbit space $\mathcal{P}$. 

By the finiteness axiom, $\mathcal{R}$ is the union of a finite number of orbits of rectangles by the action of $\pi_1(M)$ (with respect to the action of $\pi_1(M)$ on $\mathcal{P}$ defined in Theorem-Definition \ref{thdef.bifoliatedplane}). Take $R_1,...,R_n$ a representative of each orbit of rectangles in $\mathcal{R}$. Lift each of those rectangles of $\mathcal{P}$ to a rectangle $\widetilde{R_i}$ in $\mathbb{R}^3$ transverse to the lifted flow $\widetilde{\Phi}$. Using the equivariance for the projection $\pi$ of the action of $\pi_1(M)$ on $\mathcal{P}$ and the action of $\pi_1(M)$ on $\mathbb{R}^3$ by deck transformations, we have the following: for every $i\in \llbracket 1, n\rrbracket$ and every $g\in \pi_1(M)$ the rectangle $g.\widetilde{R_i}$ is a rectangle transverse to $\widetilde{\Phi}$, whose image by $\pi$ is equal to $g(R_i)$. We define $$\widetilde{\mathcal{R}}:= \{g.\widetilde{R_i}|g\in \pi_1(M), i\in \llbracket 1, n\rrbracket\}$$ Notice that $\widetilde{\mathcal{R}}$ forms a set of rectangles transverse to $\widetilde{\Phi}$ that lifts $\mathcal{R}$ on $\mathbb{R}^3$.

Take $\tilde{R}\in \widetilde{\mathcal{R}}$ and $\tilde{x}\in \tilde{R}$ such that $\pi(\tilde{R})=R$ and $\pi(\tilde{x})=x$. Consider also the stable and unstable separatrices  $\widetilde{S}:=\pi^{-1}(S)$ and $\widetilde{U}:=\pi^{-1}(U)$. It suffices to show that there exist infinitely many rectangles $\tilde{r}$ in $\widetilde{\mathcal{R}}$ such that



\begin{enumerate}
    \item the negative orbit of $\tilde{x}$ by $\widetilde{\Phi}$ intersects $\tilde{r}$. We will denote this unique intersection point by 
    $\tilde{x}_{\tilde{r}}$ (see Corollary \ref{c.orbitintersectpolygon}). 
    \item the separatrices $\widetilde{S}$ and $\widetilde{U}$ intersect $\tilde{r}$ along non-trivial segments containing $\tilde{x}_{\tilde{r}}$
\end{enumerate}   We will name the existence of infinitely many such rectangles, property  $(\star)$.

\vspace{0.2cm}
\textit{Proof that $(\star)$ implies the lemma} 

Assume that $(\star)$ holds. Since the action of $\pi_1(M)$ is properly discontinuous on $\mathbb{R}^3$, property $(\star)$ implies that  for every $T>0$ there exists $t<-T$ such that $\widetilde{\Phi}^t(\tilde{x})$ belongs in a rectangle $\widetilde{r_T}$ in $\widetilde{\mathcal{R}}$ that has a non-trivial intersection with $\widetilde{S}$ and $\widetilde{U}$. We will show that for $T$ sufficiently big the projection of $\widetilde{r_T}$ on $\mathcal{P}$ contains $G$ and intersects $R$ along a non-trivial vertical subrectangle.

Take $K$ a compact neighborhood of $\widetilde{R}$ in $\mathbb{R}^3$ and let  $\widetilde{A_T}$ be the set of points of $\widetilde{R}$, whose orbits by $\widetilde{\Phi}$ intersect $\widetilde{r_T}$. By taking $T$ sufficiently big and by using the fact that the action of $\pi_1(M)$ on $\mathbb{R}^3$ is properly discontinuous, we can assume that $\widetilde{r_T}\cap K =\emptyset$. 

Furthermore, $\widetilde{A_T}\not\subset \partial\widetilde{R}$, since the rectangles $\widetilde{R}$ and $\widetilde{r_T}$ intersect $ \widetilde{S}$ and $ \widetilde{U}$ along non-trivial segments containing $\tilde{x}$ and $\tilde{x}_{\widetilde{r_T}}$ respectively. Thanks to the previous fact and the Markovian intersection axiom applied for the projections on $\mathcal{P}$ of  $\widetilde{R}$ and $\widetilde{r_T}$, we get that the set $\widetilde{A_T}$ is a non-trivial horizontal or vertical subrectangle of $\widetilde{R}$. The connectedness of $\widetilde{A_T}$, Corollary \ref{c.orbitintersectpolygon} and the fact that $\widetilde{r_T}\cap\widetilde{A_T}=\emptyset$ imply that the negative orbit of every point of $\widetilde{A_T}$ intersects $\widetilde{r_T}$ and the positive orbit of $\widetilde{A_T}$ does not intersect $\widetilde{r_T}$. Even more, there exist a continuous map  $\phi_T: \widetilde{A_T}\rightarrow \mathbb{R}^{-}$  such that for every $z\in \widetilde{A_T}$ we have $\widetilde{\Phi}^{\phi_T(z)}(z)\in \widetilde{r_T}$ (see Figure \ref{f.pushbyflow}).

\begin{figure}[h!]
    \centering
    \includegraphics[scale=0.5]{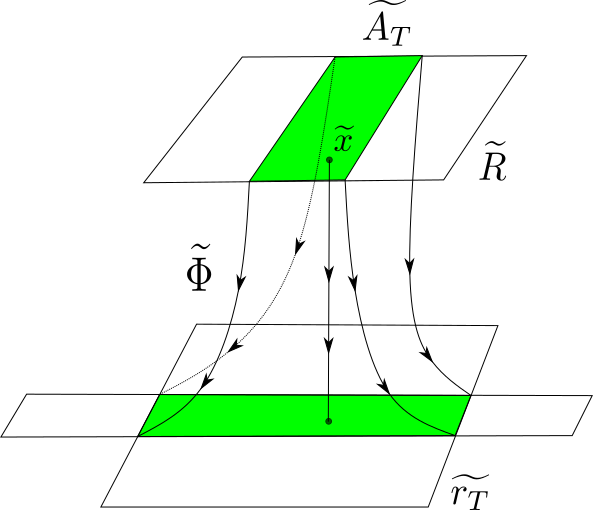}
    \caption{}
    \label{f.pushbyflow}
\end{figure}

Let $M_T:=\max_{z\in \widetilde{A_T}}\phi_T(z)<0$. By taking progressively bigger $T$ (and thus by considering rectangles $\widetilde{r_T}$ that avoid progressively bigger neighborhoods $K$ of $\widetilde{R}$) we get that $M_T\rightarrow -\infty$ as $T\rightarrow +\infty$. Also, by the same argument as in the previous paragraph, $\widetilde{\Phi}^{\phi_T}(\widetilde{A_T})$ is a non-trivial vertical or horizontal subrectangle of $\widetilde{r_T}$. Let us show that for $T$ sufficiently big $\widetilde{\Phi}^{\phi_T}(\widetilde{A_T})$ cannot be a vertical subrectangle of $\widetilde{r_T}$. 

Suppose the contrary. By the Markovian intersection axiom applied for the projections on $\mathcal{P}$ of  $\widetilde{R}$ and $\widetilde{r_T}$, this would imply that $\widetilde{A_T}$ is a horizontal subrectangle of $\widetilde{R}$. Let $u$ be an unstable leaf of $\widetilde{A_T}$ and $u_1,u_2$ be its two extremities. Consider a Riemannian metric on $M$, its lift on $\widetilde{M}$ and $d_M,\widetilde{d_M}$ their associated distances on $M$ and  $\widetilde{M}$ respectively. The fundamental group acts on $(\mathbb{R}^3,\widetilde{d_M}) $ by isometries. Since $\widetilde{\mathcal{R}}$ is the union of a finite number of orbits of rectangles by the action of $\pi_1(M)$, the distance between the two stable boundaries of any rectangle in $\widetilde{\mathcal{R}}$ is uniformly bounded from below by a constant $c_{min}$. For $T$ sufficiently big and thus also for $|M_T|$ sufficiently big, there exists a point in the negative orbit of $u_2$ that is very close to $\widetilde{\Phi}^{\phi_T(u_1)}(u_1)\in \widetilde{r_T}$. By using Corollary \ref{c.orbitintersectpolygon}, we get that for $T$ sufficiently big, there exists a point in the negative orbit of $u_2$ that belongs in $\widetilde{r_T}$ and is very close to $\widetilde{\Phi}^{\phi_T(u_1)}(u_1)\in \widetilde{r_T}$. Hence, for $T$ sufficiently big, the negative orbit of $u$ intersects $\widetilde{r_T}$ along a segment (see Corollary \ref{c.complementsingleaves}) of diameter smaller than $c_{min}$, which leads to a contradiction. 

Take $r_T$ and $A_T$ to be respectively the projections of $\widetilde{r_T}$ and $\widetilde{A_T}$ on $\mathcal{P}$. By definition, both $R$ and $r_T$ contain $x$ and $G$, provided that $G$ is taken sufficiently small. We also  proved that $R\cap r_T=A_T$ is a non-trivial horizontal subrectangle of $r_T$ and thus a non-trivial vertical subrectangle of $R$, which gives us the desired result.

\vspace{0.2cm}
\textit{Proof of property $(\star)$} 

If $\tilde{x}$ belongs in the lift on $\mathbb{R}^3$ of a periodic orbit of $\Phi$, then one can easily prove that property $(\star)$ is true. Suppose from now on that $\tilde{x}$ does not belong in the lift on $\mathbb{R}^3$ of a periodic orbit of $\Phi$. In this case, we will  prove $(\star)$ by contradiction.

Suppose that the negative orbit of $\tilde{x}$ by $\widetilde{\Phi}$ intersects finitely many rectangles in $\widetilde{\mathcal{R}}$ that intersect non-trivially the separatrices $ \widetilde{S}$ and $ \widetilde{U}$.  This implies that there exists $T>0$ such that $(\widetilde{\Phi}^t(\tilde{x}))_{t\leq -T}$ does not intersect any rectangle $\tilde{r}$ in $\widetilde{\mathcal{R}}$ intersecting $ \widetilde{S}$ and $ \widetilde{U}$ non-trivially. 

Take $x_M$ to be the projection of $\tilde{x}$ in $M$, $y_M\in M$ a point in the $\alpha$-limit of the $\Phi$-orbit of $x_M$ and $W$ a small rectangle in $M$ transverse to $\Phi$  containing $y_M$ in its interior. There exists $(t_n)_{n\in \mathbb{N}}$, an increasing sequence in $\mathbb{R}^+$ going to infinity such that $\Phi^{-t_n}(x_M)\in W$ and $\Phi^{-t_n}(x_M) \underset{n\rightarrow +\infty}{\rightarrow}y_M$. Since the orbit of $x_M$ is by hypothesis non-closed, we can assume that the points $\Phi^{-t_n}(x_M)$ are pairwise distinct. 

Let us lift everything on $\mathbb{R}^3$. Take $\tilde{y}$ to be a lift of $y_M$ on $\mathbb{R}^3$ and $\widetilde{W}$ the lift of $W$ containing $\tilde{y}$. By Corollary \ref{c.orbitintersectpolygon}, the unique lifts of the points $\Phi^{-t_n}(x_M)$ that are contained in $\widetilde{W}$ belong to different orbits of $\widetilde{\Phi}$. Therefore, there exists a sequence of $g_n\in \pi_1(M)$ such that $\widetilde{\Phi}^{-t_n}(g_n.\tilde{x})\in \widetilde{W}$ and $\widetilde{\Phi}^{-t_n}(g_n.\tilde{x}) \underset{n\rightarrow +\infty}{\rightarrow}\tilde{y}$. 

By projecting everything on $\mathcal{P}$, we have that there exists a sequence of $g_n\in \pi_1(M)$, such that $g_n(x)\underset{n\rightarrow +\infty}{\rightarrow}y$, where $y$ is the projection of $\tilde{y}$ on $\mathcal{P}$. By eventually considering a subsequence, we can assume that all the $g_n(x)$ are contained in the interior of a germ of quadrant of $y$ or inside a stable or unstable separatrix of $y$. One can check that in every case, by the finite return time axiom and by eventually considering another subsequence, there exists $r\in \mathcal{R}$ containing $y$, containing also all the $g_n(x)$ and intersecting non-trivially all $g_n(S)$ and $g_n(U)$. 

This implies that there exists $\tilde{r}\in \widetilde{\mathcal{R}}$ such that for every $n \in \mathbb{N}$ the orbits by $\widetilde{\Phi}$ of the points $\tilde{y}\in \widetilde{W}$ and $\widetilde{\Phi}^{-t_n}(g_n.\tilde{x})\in \widetilde{W}$ intersect $\tilde{r}$ and also   $g_n.\widetilde{S}$ and $g_n.\widetilde{U}$ intersect $\tilde{r}$ non-trivially. Moreover, thanks to Corollary \ref{c.orbitintersectpolygon}, using the facts that the orbit of $\tilde{y}$ intersects $\tilde{r}$ and $\widetilde{\Phi}^{-t_n}(g_n.\tilde{x}) \underset{n\rightarrow +\infty}{\rightarrow}\tilde{y}$, we get that the orbits of the $\widetilde{\Phi}^{-t_n}(g_n.\tilde{x})\in \widetilde{W}$ intersect $\tilde{r}$ in a uniformly bounded time with respect to $n$. By our previous arguments, we have that for every $n$, in a uniformly bounded time the point $\widetilde{\Phi}^{-t_n}(\tilde{x})$ will cross a rectangle in $\widetilde{\mathcal{R}}$ intersecting $\widetilde{S}$ and $\widetilde{U}$ non-trivially. By hypothesis $t_n\underset{n\rightarrow +\infty}{\rightarrow} +\infty$, which contradicts the fact that $(\widetilde{\Phi}^t(\tilde{x}))_{t\leq -T}$ does not intersect any rectangle in $\widetilde{\mathcal{R}}$ intersecting $ \widetilde{S}$ and $ \widetilde{U}$ non-trivially.  
\end{proof}
Naturally, by symmetry we also have that:  
\begin{lemm}\label{l.horizsubrectangleexists}
Consider $\mathcal{R}$ a Markovian family of $\Phi$. Take $R\in\mathcal{R}$, $x\in R$ and  $G$ a sufficiently small germ of quadrant of $x$ contained in $R$. There exists $R_h\in \mathcal{R}$ such that $G\subset R_h$ and $R_h\cap R$ is a non-trivial horizontal subrectangle of $R$.
\end{lemm}
Using the notations of our previous proof, during the proof of property $(\star)$, we showed that thanks to the third axiom in the definition of a Markovian family (see Definition \ref{d.markovfamily}), the negative orbit by $\widetilde{\Phi}$ of any point $\tilde{x}\in \mathbb{R}^3$ will intersect in finite time a rectangle in $\widetilde{\mathcal{R}}$. This is  the reason why we called this axiom the finite return time axiom. 

Also, during the proof of the sufficiency of property $(\star)$ we proved that: 
\begin{rema}\label{r.negativeorbitrectangles}
    Consider $\mathcal{R}$ is a Markovian family of $\Phi$ and $\widetilde{\mathcal{R}}$ a $\pi_1(M)$-invariant family of rectangles in $\mathbb{R}^3$ transverse to $\widetilde{\Phi}$ such that 
    
    \begin{enumerate}
        \item if $\pi$ is the projection from $\mathbb{R}^3$ to $\mathcal{P}$, then for any $\tilde{R}\in \widetilde{\mathcal{R}}$ we have that $\pi(\widetilde{R})\in \mathcal{R}$
        \item for any $r \in \mathcal{R}$ there exists a unique $\widetilde{r}\in \widetilde{\mathcal{R}}$ such that  $\pi(\widetilde{r})=r$
    \end{enumerate}
     If the negative orbit by $\widetilde{\Phi}$ of $\tilde{R}\in \widetilde{\mathcal{R}}$ intersects along non-trivial subrectangles the rectangles $\widetilde{R_1},...,\widetilde{R_n},...$ then for $n$ sufficiently big we have that $ \pi(\widetilde{R})\cap \pi(\widetilde{R_n})$ is a horizontal subrectangle of $\pi(\widetilde{R_n})$ and a vertical subrectangle of $\pi(\widetilde{R})$.
\end{rema}
\begin{lemm}\label{l.infiniteintersectionverticalrectangles}
Consider $\mathcal{R}$ a Markovian family of $\Phi$ and a sequence of rectangles $(r_n)_{n\in \mathbb{N}}$ in $\mathcal{R}$ such that for every $k\in \mathbb{N}$, $r_{k+1}\cap r_k$ is a non-trivial vertical subrectangle of $r_k$. We have that $\overset{+\infty}{\underset{k=0}{\cap}} r_{k}$ is an unstable segment of $r_0$.
\end{lemm}
\begin{proof}
Notice that thanks to the Markovian intersection property, for any strictly increasing function $k:\mathbb{N}\rightarrow \mathbb{N}$, the intersection $r_{k(n+1)}\cap r_{k(n)}$ is a vertical subrectangle of $r_{k(n)}$. It therefore suffices to show the lemma for any subsequence of $(r_n)_{n\in \mathbb{N}}$. 

By the finiteness axiom and by eventually considering a subsequence, we can assume that all the rectangles in our sequence belong to the same orbit of rectangles in $\mathcal{P}$. In other words, for any $n \in \mathbb{N}$ there exists $g_{n} \in \pi_1(M)$ such that $g_{n}(r_0)=r_n$. 

Consider now $\widetilde{r_0}$ a rectangle in $\mathbb{R}^3$, transverse to the lifted flow $\widetilde{\Phi}$, whose projection on $\mathcal{P}$ is $r_0$. Denote by $\widetilde{r_n}$ the rectangle $g_n.\widetilde{r_0}$. 

Since the action of $\pi_1(M)$ on $\mathbb{R}^3$ is properly discontinuous, by eventually considering a subsequence $(r_n)_{n\in \mathbb{N}}$, we can assume that $\widetilde{r_n} \cap \widetilde{r_m} = \emptyset$ for all $n\neq m$. Moreover, by hypothesis, for every $n \in \mathbb{N}$ there are orbits of $\widetilde{\Phi}$ intersecting both $\widetilde{r_0}$ and $g_{n}.\widetilde{r_0}$. Let us denote by $A_n$ the intersection of all these orbits with $\widetilde{r_n}$. By Corollary \ref{c.orbitintersectpolygon}, the set $A_n$ projects homeomorphically to $r_0 \cap r_n \subset \mathcal{P}$; hence, by the Markovian intersection axiom $A_n$ is a horizontal subrectangle of $\widetilde{r_n}$. 

Even more, thanks once again to Corollary \ref{c.orbitintersectpolygon}, there exists $t_n: A_n \rightarrow \mathbb{R}$ a continuous function such that $\widetilde{\Phi^{t_n(x)}}(x) \in \widetilde{r_0}$ for every $x \in A_n$ and $\widetilde{\Phi^{t_n}}(A_n)$ is a vertical subrectangle of $\widetilde{r_0}$. Since $\widetilde{r_n} \cap \widetilde{r_0} = \emptyset$, we have either $t_n(x)>0$ for all $x \in A_n$ or $t_n(x)<0$ for all $x \in A_n$.  By Remark \ref{r.negativeorbitrectangles}, by eventually removing a finite number of $r_n$, we can assume that $t_n>0$ for every $n\in \mathbb{N}$. Furthermore, we have that 
\begin{equation}\label{eq.goestoinfinity}
\min_{x\in A_n}t_n(x)\underset{n\rightarrow+\infty}{\longrightarrow} +\infty
\end{equation} 
Indeed, if we take any $c>0$ and $K_c:=\underset{t\in[-c,c]}{\cup}\widetilde{\Phi^t}(\widetilde{r_0})$, then, since the action by $\pi_1(M)$ on $\mathbb{R}^3$ is properly discontinuous, we can find $N$ sufficiently big, such that for every $n\geq N$ we have $g_{n}.\widetilde{r_0} \cap K_c = \emptyset$. By taking $c$ progressively bigger, one can easily prove (\ref{eq.goestoinfinity}).

Let us now consider $\overset{+\infty}{\underset{k=0}{\cap}} r_{k}$. As a decreasing intersection of vertical subrectangles, the set $\overset{+\infty}{\underset{k=0}{\cap}} r_{k}$ is either an unstable segment of $r_0$ or a vertical subrectangle of $r_0$. Take $x\in \overset{+\infty}{\underset{k=0}{\cap}} r_{k}$. Denote by $\widetilde{x_n}$ the lift of $x$ inside $\widetilde{r_n}$, by $\widetilde{s_n}$ the stable leaf of $\widetilde{r_n}$ containing $\widetilde{x_n}$ and by $\widetilde{s_n^1},\widetilde{s_n^2}$ the two extremities of $\widetilde{s_n}$. 

Let $\widetilde{d_M}$ be the distance associated to the lift of a Riemannian metric of $M$ on  $\mathbb{R}^3$. Recall that $\pi_1(M)$ acts by isometries on $(\mathbb{R}^3,\widetilde{d_M})$; hence, the diameter of any stable leaf of $\widetilde{r_n}$ -in particular the diameter of $\widetilde{s_n}$- is uniformly (with respect to $n$) bounded from above. For $n$ sufficiently big, there exists a point in the positive orbit of $\widetilde{s_n^2}$ that is very close to $\widetilde{\Phi}^{t_n(\widetilde{s_n^1})}(\widetilde{s_n^1})\in \widetilde{r_0}$. Even more, by Corollary \ref{c.orbitintersectpolygon}, we get that for $n$ sufficiently big, there exists a point in the positive orbit of $\widetilde{s_n^2}$ that belongs in $\widetilde{r_0}$ and is very close to $\widetilde{\Phi}^{t_n(\widetilde{s_n^1})}(\widetilde{s_n^1})\in \widetilde{r_0}$. Hence, for $n$ sufficiently big, the positive orbit of $\widetilde{s_n}$ intersects $\widetilde{r_0}$ along a segment (see Corollary  \ref{c.complementsingleaves}) of diameter that goes to $0$ as $n\rightarrow +\infty$. This proves that the intersection of all the $\widetilde{\Phi^{t_n}}(A_n)$ is an unstable segment of $\widetilde{r_0}$. By projecting everything on $\mathcal{P}$, we get the desired result. 
\end{proof}

Naturally, by symmetry the following is also true: 
\begin{lemm}\label{l.infiniteintersectionhorizontalrectangles}
Consider $\mathcal{R}$ a Markovian family of $\Phi$ and a sequence of rectangles $(r_n)_{n\in \mathbb{N}}$ in $\mathcal{R}$ such that for every $k\in \mathbb{N}$, $r_{k+1}\cap r_k$ is a non-trivial horizontal subrectangle of $r_k$. We have that $\overset{+\infty}{\underset{k=0}{\cap}} r_{k}$ is a stable segment of $r_0$.
\end{lemm}

\begin{rema}\label{r.comparison}
In the previous lemmas, we established a ``bridge"  between general Markovian families and Markov partitions. Indeed, statements analogous to the previous ones were already known to be true for Markov partitions. Fix $\mathcal{M}$ a Markov partition of a pseudo-Anosov flow $\Phi$ in $M^3$: 

\vspace{0.2cm}
\hspace{-0.8cm}
\begin{tabular}{ | p{4.6cm} | p{12cm}| } 
  \hline
 \hspace{0.5cm} Markovian families & \hspace{4cm}   Markov partitions \vspace{0.1cm} \\
 
  \hline
  Lemmas \ref{l.vertsubrectangleexists} and \ref{l.horizsubrectangleexists}& Every positive and negative orbit (by $\Phi$) in $M$ intersects in bounded time a rectangle of $\mathcal{M}$ \\ 

  \hline
 Lemmas \ref{l.infiniteintersectionverticalrectangles} and \ref{l.infiniteintersectionhorizontalrectangles}& Let $\mathcal{M}:=
 \{R_1,...,R_m\}$ and $f_n$ the $n-$th return map on $\underset{i\in\llbracket 1, m \rrbracket}{\cup}R_i$ with $n\in \mathbb{Z}$. For any $i, j\in \llbracket 1, m \rrbracket$, we have that the closures of the connected components of $f_n(\inte{R_i})\cap \inte{R_j}$ are vertical (resp. horizontal) subrectangles of $R_j$, whose width (resp. height) goes uniformly to $0$ as $n\longrightarrow +\infty$ (resp. $n\longrightarrow -\infty$)\\ 
  \hline
\end{tabular}
\end{rema}
The previous analogies between general Markovian families and Markov partitions provide some evidence supporting Conjecture \ref{conj.markov}. Several more analogies between the two previous notions will be proven later on in a more general setting (see for instance Lemmas \ref{l.periodicboundary}, \ref{l.existenceofpredecessors}, \ref{l.crossingrectanglesnoperiodicpoints}, \ref{l.crossingrectangleswithperiodicpoints}). 

\subsection{Action of the fundamental group on the bifoliated plane}
Consider $\Phi$ a pseudo-Anosov flow on a smooth, closed and connected $3$-manifold $M$, and  $F^s,F^u$ its stable and unstable foliations. Denote by $\widetilde{M}$ the universal cover of $M$ and by $\mathcal{P}$ the bifoliated plane of $\Phi$ together with its stable and unstable foliations $\mathcal{F}^s,\mathcal{F}^u$. Our goal in this section consists in defining the notion of \emph{Markovian action}, a group action on a bi-foliated plane having several properties in common with the group actions on the plane arising from pseudo-Anosov flows. 

\begin{defi}\label{d.periodicpseudoanosovaction}
We will say that the point $x\in \mathcal{P}$ \emph{corresponds to the orbit $\gamma$} of $\Phi$ or that \emph{the orbit $\gamma$ corresponds to the point $x$} if there exists a lift of $\gamma$ on $\widetilde{M}$, whose 
projection on $\mathcal{P}$ is $x$. 

A point $x$ in $\mathcal{P}$ will be called \emph{periodic} if there exists   $g\in \pi_1(M)-\{id\}$ 
such that $g(x)=x$. We similarly define \emph{periodic stable} or \emph{unstable} leaves in $\mathcal{P}$. 
\end{defi}
\begin{defi}
    Let $g\in \text{Hom}(\mathbb{R})$. We will say that $g$ is a \emph{topological contraction} (resp. \emph{expansion}) if there exists $x_0\in \mathbb{R}$ such that for every $x\in \mathbb{R}$ the sequence $g^n(x)$ converges to $x_0$ as $n\rightarrow +\infty$ (resp. as $n\rightarrow -\infty$). 
\end{defi}
Even though the action by $\pi_1(M)$ on $\mathcal{P}$ is not well understood globally for most pseudo-Anosov flows, there several known results concerning the local dynamics  around periodic points and periodic stable/unstable leaves.  
\begin{prop}\label{p.periodicinbifoliated}
Consider $\Phi$ a pseudo-Anosov flow on $M^3$, $\mathcal{P}$ its orbit space, $x\in \mathcal{P}$ a periodic point and $\mathcal{L}$ a periodic stable or unstable leaf in $\mathcal{P}$. 
We have the following:
\begin{enumerate}
\item The point $x$ corresponds to a unique periodic orbit $\gamma$ of $\Phi$. Conversely, any periodic orbit of $\Phi$ corresponds to infinitely  many periodic points in $\mathcal{P}$
\item Consider $x_0\in \gamma$ and $\gamma_g$ the element in $\pi_1(M,x_0)$ represented by the dynamically oriented orbit $\gamma$. The stabilizer of $x$ in $\pi_1(M)$ is isomorphic to $<\gamma_g> \cong \mathbb{Z}$
\item The leaf $\mathcal{L}$ contains a unique periodic point
\item If $h\in \pi_1(M)-\{ id\}$ fixes $\mathcal{L}$, then the action of $h$ on $\mathcal{L}$ is a topological contraction or expansion 

\vspace{0.2cm}
\item Consider $\phi^{\lambda_h,\lambda_v}_{p,r,o}:\mathbb{R}^2\rightarrow \mathbb{R}^2$ the local model for a pseudo-hyperbolic fixed point $\gamma_{p,r,o}$ with stretching $\lambda^u$, compression $\lambda^h$, $p\geq 2$ prongs,  rotation $r \in \llbracket 0, p-1 \rrbracket$ (resp. type $r
\in \{1,2\}$) and positive (resp. negative) orientation. If $g\in \pi_1(M)-\{ id\}$ fixes the point $x\in \mathcal{P}$, then there exist $p\geq 2$, $o\in \{+,-\}$, $r \in \llbracket 0, p-1 \rrbracket$ if $o=+$ (resp. $r\in \{1,2\}$ if $o=-$) and $U_x$ a neighborhood of $x$ in $\mathcal{P}$ such that the action of $g$ in $U_x$ is conjugated to $\phi^{\lambda_h,\lambda_v}_{p,r,o}$ restricted on a compact neighborhood of $\gamma_{p,r,o}$ inside $\mathbb{R}^2$

\end{enumerate}
\end{prop}
\begin{proof}[\quad{}Proof of property $(1)$]
Denote by $\widetilde{\Phi}$ the lift of $\Phi$ on the universal cover $\widetilde{M}$ and consider $g\in \pi_1(M)-\{id\}$ fixing $x\in \mathcal{P}$. The point $x$ corresponds to an orbit $\mathcal{O}$ of $\widetilde{\Phi}$ that is preserved by the action of $g$. As the action of $g$ on $\widetilde{M}$ does not admit fixed points ($g\neq id$ is a deck transformation), $\mathcal{O}$ projects to a periodic orbit of $\Phi$. Conversely, given a periodic orbit $\gamma$ of $\Phi$, there are infinitely many orbits of $\widetilde{\Phi}$ that lift $\gamma$, each one of which, thanks to Corollary \ref{c.periodicnontrivial}, is preserved by the action of a non-trivial element of $\pi_1(M)$. We get the desired result by projecting everything on $\mathcal{P}$. 

\vspace{0.2cm}
\textit{Proof of property $(2)$}
Take $x,x_0,\gamma,\gamma_g$ as in the statement. Consider also $g\in \text{Stab}(x)-\{id\}$ and $\widetilde{x_0}$ a lift of $x_0$ on $\mathbb{R}^3$, whose projection on $\mathcal{P}$ is $x$. The action of $g$ on $\mathbb{R}^3$ preserves the $\widetilde{\Phi}$-orbit of $\widetilde{x_0}$. The orbit segment of $\widetilde{\Phi}$ going from $\widetilde{x_0}$ to $g.\widetilde{x_0}$ projects to a loop in $M$, whose homotopy class in $\pi_1(M,x_0)$ is equal to $\gamma_g^k$ for some $k\in \mathbb{Z}$. It follows that $\text{Stab}(x)\subset <\gamma>$. Conversely, it is easy to see that $<\gamma_g>\subset \text{Stab}(x)$. Finally, thanks to Corollaries \ref{c.torsionfree} and \ref{c.periodicnontrivial}, we get that $\text{Stab}(x)=<\gamma_g>\cong\mathbb{Z}$. 

\vspace{0.2cm}
\textit{Proof of property $(5)$}
Using our previous notations, property $(5)$ is an immediate consequence of the fact that $\text{Stab}(x)=<\gamma_g>$ and of Proposition \ref{p.aroundcircleprong}. 

\vspace{0.2cm}
\textit{Proof of property $(3)$}
As a result of Definition \ref{d.pseudoanosovflow}, a stable or unstable leaf of a pseudo-Anosov flow contains at most one periodic orbit. By Item (1) and our construction of $\mathcal{F}^s, \mathcal{F}^u$ (see Theorem-Definition \ref{thdef.bifoliatedplane}) the same applies for the leaves in $\mathcal{F}^s, \mathcal{F}^u$. Therefore, it suffices to prove that $\mathcal{L}$ contains a periodic point. If $\mathcal{L}$ is a singular leaf, then the result is trivial (see Proposition \ref{p.firstperiodicorbits}); we can therefore assume without any loss of generality that $\mathcal{L}$ is a regular stable leaf in $\mathcal{P}$. 

Consider $\mathcal{R}$ a Markovian family of $\Phi$ (the existence of such a family follows from Theorem \ref{t.reducedmarkovpartitionexists} and Proposition \ref{p.projectionmarkovpartition}). Let $y\in \mathcal{L}$, $h\in \text{Stab}(\mathcal{L})-\{id\}$ and $\mathcal{L}^+(y)$ be the stable separatrix of $y$ containing $h(y)$. Thanks to Item (2) of Proposition \ref{p.propertiesoffoliformarkovianactions}, $\mathcal{L}^+(y)$ is a properly embedded closed half-line in $\mathcal{P}$. 

By the finite return axiom, there exists $R_0\in \mathcal{R}$ such that $y\in R_0$ and such that $R_0$ intersects non-trivially $\mathcal{L}^+(y)$ along a segment containing $y$ (this intersection is indeed a segment because of Item (5) of Proposition \ref{p.propertiesoffoliformarkovianactions}). Thanks to Lemma \ref{l.horizsubrectangleexists}, there exists $R_1\in \mathcal{R}$ such that $R_0\cap R_1$ is a horizontal subrectangle of $R_0$ and such that $\mathcal{L}^+\cap R_0\subset R_1$. By applying infinitely many times this argument, we construct $R_0,R_1,...,R_i,...$ such that $R_i\cap R_{i+1}$ is a horizontal subrectangle of $R_i$ and $\mathcal{L}^+(y)\cap R_i\subset R_{i+1}$ for every $i\in \mathbb{N}$. By the finiteness axiom, by eventually considering a subsequence, we may assume that all the $R_i$ belong in the same orbit of rectangles by the action of $\pi_1(M)$. Take $g_i\in \pi_1(M)$ such that $g_i(R_i)=R_{i+1}$. 

Denote by $J_i$ the segment $R_i\cap \mathcal{L}^+$. By construction, $J_i$ is of the form $[y,y_i]^s$, where $y_i\in \mathcal{L}^+$. Since $R_i\cap R_{i+1}$ is a horizontal subrectangle of $R_i$, we have that $J_i\subseteq J_{i+1}$ for every $i\in \mathbb{N}$. 

Assume first that there exists $I$ such that $J_i=J_{i+1}$ for every $i\geq I$.  Denote by $u_i$ the unstable boundary component of $R_i$ containing $y_i$. By eventually considering a subsequence of the $(R_i)_{i\in\mathbb{N}}$, we can assume without any loss of generality that $g_i(u_i)=u_{i+1}$ for every $i\in \mathbb{N}$. For every $i\geq I$, by our hypothesis, the segments $u_{i}, u_{i+1}$ belong in the same unstable leaf $U$ and by the Markovian intersection axiom $u_{i+1}\subset u_{i}$. It follows that $g_{i}$ has a fixed point inside $u_{i+1}$. Since $U$ can contain at most one periodic point, it follows that there exists a periodic point inside $\underset{i\geq I}{\cap}R_j\cap U$. By Lemma \ref{l.infiniteintersectionhorizontalrectangles}, the previous set contains a unique point, which in our case is none other than $y_I\in \mathcal{L}^+$. It follows that $y_I$ is periodic, which gives us the desired result in this case. 

By our previous argument, by eventually considering a subsequence, we may assume that for every $i\in \mathbb{N}$ we have $J_i\subsetneq J_{i+1}$. Assume now that $\underset{i\in \mathbb{N}}{\cup}J_i=[y,y_{lim})$, where $y_{lim}\in \mathcal{L}^+$ and $[y,y_{lim})$ denotes the semi-closed segment in $\mathcal{L}^+$ going from $y$ to $y_{lim}$. In this case, take $\mathcal{L}^-(y_{lim})$ the stable separatrix of $y_{lim}$ containing $y$. By the finite return axiom, there exists $R_{lim}\in \mathcal{R}$ such that $R_{lim}$ contains $y_{lim}$ and intersects $\mathcal{L}^-(y_{lim})$ along a non-trivial segment containing $y_{lim}$. For $i$ sufficiently big $\inte{R_i}\cap \inte{R_{lim}}\neq \emptyset$ and also, thanks to Lemma \ref{l.infiniteintersectionhorizontalrectangles},  $R_i$ is a very thin rectangle along the unstable direction. Therefore, by the Markovian intersection property, when $i$ is sufficiently big, $R_i\cap R_{lim}$ is a horizontal subrectangle of $R_{lim}$. This is impossible, since $y_{lim}\notin J_i$ and thus $\mathcal{L}\cap R_{lim}\not\subset \mathcal{L}\cap R_i$. It follows that this case is impossible and that $\underset{i\in \mathbb{N}}{\cup}J_i=\mathcal{L}^+(y)$. 

By our previous arguments, there exists $i\in \mathbb{N}$ such that $R_i$ contains both $y$ and $h(y)$, and such that $h(y)\notin \partial^u R_i$. By considering if necessary $h^2$, we may assume that $h$ preserves the connected components of $\mathcal{P}-\mathcal{L}$ (see Item (2) of Proposition \ref{p.propertiesoffoliformarkovianactions}). By the Markovian intersection axiom, this implies that $h(R_i)\cap R_i$ is a horizontal or vertical subrectangle of $R_i$. The index of $h$ along $\partial R_i$ being non-zero, we get that $h$ has a fixed point $Y$ in $R_i$. If $Y\notin \mathcal{L}$, since $\mathcal{L}$ is preserved by $h$, we can find another fixed point by $h$ on the unstable leaf of $Y$, which is impossible since every unstable leaf in $\mathcal{F}^u$ contains at most one periodic point. This proves that $Y\in \mathcal{L}$ and gives us the desired result.

\vspace{0.2cm}
\textit{Proof of property $(4)$}
 Take $Y$ the unique periodic point in $\mathcal{L}$. By eventually considering a power of $h$, we can assume that $h$ preserves every separatrix of $Y$. Take $S:=[Y,\infty)$ a separatrix of $Y$. Since $\mathcal{L}$ does not contain any periodic point other than $Y$ and thanks to Item (5) of this proposition, by eventually considering $h^{-1}$, we get that for every $z\in S$ close to $Y$ the sequence $h^{n}(z)$ goes to $\infty$ when $n\rightarrow +\infty $ and goes to $Y$ when $n\rightarrow -\infty $. It follows that for every $z\in S$, $h^n(z)\underset{n\rightarrow -\infty}{\longrightarrow}Y$, which finishes the proof of the proposition. 
\end{proof}

As we have previously stated, a pseudo-Anosov flow contains infinitely many periodic orbits. According to the following proposition, this remains true for its associated action on the bifoliated plane. 
\begin{prop}\label{p.densityperiodic}
    Let $\Phi$ be a pseudo-Anosov flow on $M^3$ and  $\mathcal{P}$ its orbit space. We have that: 
    \begin{enumerate}
        \item There exists a finite set of periodic stable or unstable leaves in $\mathcal{P}$, whose orbit by $\pi_1(M)$ is dense in $\mathcal{P}$ 
    \end{enumerate}
    If furthermore $\Phi$ is transitive (i.e. $\Phi$ has a dense orbit in $M$), 
    \begin{enumerate}
    \setcounter{enumi}{1}
        \item The set of periodic points (or periodic stable/unstable leaves) in $\mathcal{P}$ is dense
        \item The orbit by $\pi_1(M)$ of every stable or unstable leaf in $\mathcal{P}$ is dense 
        \item There exists a point in $\mathcal{P}$ whose orbit by $\pi_1(M)$ is dense
    \end{enumerate}
\end{prop}
Property (1) is a consequence of Proposition 7 of \cite{markovpseudoanosov}. Property (2) follows from Proposition \ref{p.periodicorbitsaredense}. Property (3) is a result of the fact that for any transitive pseudo-Anosov flow the stable and unstable foliations are minimal (this follows from the same proof as Proposition 7 of \cite{markovpseudoanosov}). Finally, property (4) results immediately from the definition of transitivity.

Given a plane $\mathcal{P}$, a closed and smooth $3$-manifold $M$ and a continuous action $\rho:\pi_1(M)\rightarrow \text{Hom}(\mathcal{P})$ preserving a pair of transverse singular foliations, it very natural to ask under which conditions this action is associated to a pseudo-Anosov flow on $M$:
\begin{ques}\label{q.actionplane}
Give a necessary and sufficient condition for $\rho$ to be realised (up to conjugation or semi-conjugation) as the action of $\pi_1(M)$ on the orbit space of a pseudo-Anosov flow on $M$. 
\end{ques}

Significant progress has been made recently in the direction of the previous problem with the definition of the notion of Anosov-like actions on the plane by T.Barthelmé, K.Mann and S.Frankel in \cite{circleatinfinity}. In this paper, we will introduce and make use of a slightly different family of actions on the plane that conjecturally provide an answer to the Question \ref{q.actionplane}: 
\begin{defi}\label{d.markovianaction}
    A $C^0$ action $\rho$ of a countable group $G$ on a plane $\mathcal{P}$ is called a \emph{Markovian action} if it satisfies the following properties:
    \begin{enumerate}
        \item The action of $G$ preserves a pair of transverse singular foliations $(\mathcal{F}^s,\mathcal{F}^u)$, which we will call the \emph{stable} and \emph{unstable} foliations of $\rho$
        \item  If $l$ and $l'$ are two stable (resp. unstable) leaves of   that are not separated in the leaf space of $\mathcal{F}^s$ (resp. $\mathcal{F}^u$), then there exists a non-trivial element in $G$ fixing simultaneously $l$ and $l'$
        \item If a non-trivial element $g \in G$ preserves a leaf $l\in \mathcal{F}^s$ or $l\in \mathcal{F}^u$, then it has a unique fixed point in $l$. Moreover, this fixed point is the same for every non-trivial element in $G$ that preserves $l$
        \item The stabilizer in $G$ of any point in $\mathcal{P}$ is either trivial or isomorphic to $\mathbb{Z}$
        \item Every singularity of $\mathcal{F}^s$ and $\mathcal{F}^u$ admits a non-trivial stabilizer in $G$ 
        \item For any non-trivial $g\in G$ and any $x\in \mathcal{P}$ fixed by $g$, up to changing $g$ to $g^{-1}$, we have that $g$ is 
        topologically expanding on $\mathcal{F}^u(x)$ and topologically contracting on $\mathcal{F}^s(x)$
        \item There exists a finite set of stable  (resp. unstable) leaves with non-trivial stabilizers in $G$, whose orbit by $G$ is dense in $\mathcal{P}$ 
        \item  $G$ preserves a Markovian family $\mathcal{R}$ in $(\mathcal{P},\mathcal{F}^s,\mathcal{F}^u)$  
        \end{enumerate} 

        If furthermore the Markovian family $\mathcal{R}$ verifies the results of the Lemmas \ref{l.vertsubrectangleexists}, \ref{l.horizsubrectangleexists}, \ref{l.infiniteintersectionverticalrectangles} and \ref{l.infiniteintersectionhorizontalrectangles}, then we will say that $\mathcal{R}$ is a \emph{strong Markovian family} and that the action of $G$ is a \emph{strong Markovian action}.

\end{defi}

It is clear that, thanks to Lemmas \ref{l.vertsubrectangleexists}, \ref{l.horizsubrectangleexists}, \ref{l.infiniteintersectionverticalrectangles} and \ref{l.infiniteintersectionhorizontalrectangles} and Propositions  \ref{p.projectionmarkovpartition}, \ref{p.periodicinbifoliated} and \ref{p.densityperiodic}: 

\begin{rema} \label{r.pseudoanosovactionsarestrongmarkovian}
 Any Markovian family associated to a pseudo-Anosov flow is a strong Markovian family. It follows that the action of the fundamental group on the bifoliated plane of any pseudo-Anosov flow is a strong Markovian action. 
\end{rema} Conversely, we conjecture that: 
\begin{conj}\label{conj.markovianactions}
    If $\rho$ is a strong Markovian action, then there exists a pseudo-Anosov flow $(M,\Phi)$, such that the action of $\pi_1(M)$ on the bifoliated plane of $\Phi$ is conjugated to $\rho$. 
\end{conj}

\subsection{Dehn-Goodman-Fried surgery}\label{s.surgeries}
Let $M$ be a closed, smooth and \emph{orientable} 3-manifold carrying a pseudo-Anosov flow $\Phi$. Let $F^s$ and $F^u$ be the stable and unstable foliations of $\Phi$. In this section, we will define the operation of Dehn-Goodman-Fried surgery, a powerful tool allowing to construct new pseudo-Anosov flows starting from any pseudo-Anosov flow. Contrary to its original definition for Anosov flows (see \cite{Fried}), one of the main difficulties in the definition of a Dehn-Goodman-Fried surgery in the pseudo-Anosov flow setting consists in the lack of differentiability of a pseudo-Anosov flow.

Following the original definition of Fried (see \cite{Fried}), any Dehn-Goodman-Fried surgery on $\Phi$ consists of two distinct geometric/dynamical operations: 
\begin{enumerate}
    \item the \emph{blow-up operation}, which consists in exploding a periodic orbit $\gamma$ to a torus. In more geometric terms, 
    if the curve $\gamma$ was smooth, then the blow-up operation would  consist in replacing $\gamma$ by the unit normal bundle of $\gamma$ 
    \item the \emph{blow-down operation}, which consists in crushing the torus obtained after the previous explosion to a new circle
\end{enumerate}
In the ensuing pages, we will define in detail each of the previous operations first for hyperbolic models and then for general pseudo-Anosov flows. 

Notice that, contrary to the previous sections, we assumed here $M$ to be orientable. The reason behind this is that a Dehn-Goodman-Fried surgery can only be defined for periodic orbits with orientable tubular neighborhoods (see Remark \ref{r.surgery}). By assuming that $M$ is orientable, we can therefore guarantee that this condition is verified for every periodic orbit of $\Phi$.

\subsubsection{Blow-up operation for hyperbolic models}\label{s.blowupmodels}
\quad

Let us first recall the definition of a hyperbolic model given in  Section \ref{ss.defipseudoanosov} and introduce at the same time some useful notations. Denote by $\phi_2:\mathbb{R}^2\rightarrow \mathbb{R}^2$ the map $\phi_2(x,y)= (\frac{1}{2}x, 2 y)$, whose stable and unstable foliations, say $\mathcal{F}^s_2,\mathcal{F}_2^u$, coincide  respectively with the horizontal and vertical foliations in $\mathbb{R}^2$. Denote also by $\pi_p:\mathbb{C}\rightarrow \mathbb{C}$ the map $\pi_p(z)=z^p$ (with $p\in \mathbb{N}_{\geq 3}$). 

The images of the foliations $\mathcal{F}^s_2,\mathcal{F}_2^u$ by $\pi_2$ define a pair of  transverse singular foliations $\mathcal{F}^s_1,\mathcal{F}_1^u$  with a $1-$prong singularity at $0$. Similarly, for every $p\geq 3$ the pre-images of $\mathcal{F}^s_1,\mathcal{F}_1^u$ by $\pi_p$ define a pair of transverse singular  foliations $\mathcal{F}^s_p,\mathcal{F}_p^u$ with a $p$-prong singularity at $0$. 

Furthermore, the map $\phi_2$ projects by $\pi_2$ to a unique homeomorphism $\phi_1:\mathbb{R}^2\rightarrow \mathbb{R}^2$ preserving $\mathcal{F}_1^u$ and $\mathcal{F}_1^s$ and satisfying $\pi_2\circ \phi_2=\phi_1\circ\pi_2$. Similarly, for any $p\geq 3$, denote by $\phi_p:\mathbb{R}^2\rightarrow \mathbb{R}^2$ the unique homeomorphism preserving the foliations $\mathcal{F}^s_p,\mathcal{F}_p^u$, preserving each stable and unstable prong of the origin and satisfying $\pi_p\circ \phi_p=\phi_1\circ\pi_p$.

We define for every $p\geq 1$ and $k\in \llbracket 0, p-1\rrbracket $ $$\phi_{p,k}:=\phi_p\circ R_{k/p}=R_{k/p} \circ \phi_p$$ where $R_{\theta}:\mathbb{R}^2\rightarrow \mathbb{R}^2$ is the rotation around the origin of angle $2\pi\theta$. Notice that for every $p\geq 1$ and $k\in \llbracket 0, p-1\rrbracket $,  the map $\phi_{p,k}$ preserves the foliations $\mathcal{F}^s_p$ and $\mathcal{F}_p^u$, acts as a contraction along $\mathcal{F}^s_p$ and as an expansion along $\mathcal{F}^u_p$.  The map $\phi_{p,k}$ was called in Section \ref{ss.defipseudoanosov} the local model for a pseudo-hyperbolic fixed point with stretching $2$, compression $2$, $p$ prongs, rotation $k$ and positive orientation.

We similarly defined in Section \ref{ss.defipseudoanosov} a local model for a $p$-prong circle singularity. Consider the mapping torus
$$N_{p,k}:= \frac{\mathbb{R}^2\times \mathbb{R}}{((x,y),t+1) \sim (\phi_{p,k}(x,y),t) }$$
endowed with the constant speed vertical flow $\Phi_{p,k}=(\Phi^t_{p,k})_{t\in\mathbb{R}}$ given by the vector field $\frac{\partial}{\partial t}$. We called  $(N_{p,k},\Phi^t_{p,k})$ the local model for a pseudo-hyperbolic periodic orbit with stretching $2$, compression $2$, $p$ prongs, rotation $k$ and positive orientation. We will denote by $\gamma_{p,k}$ the periodic orbit of $\Phi_{p,k}$ defined by the suspension of the origin and by $F^s_{p,k},F^u_{p,k}$ the projection of $\mathcal{F}^s_p\times \mathbb{R}, \mathcal{F}^u_p\times \mathbb{R}$ on $N_{p,k}$. We called $F^s_{p,k},F^u_{p,k}$ the local weak stable and unstable foliations of $\Phi_{p,k}$. 

We are now ready to define the blow-up operation for a hyperbolic model. First, notice that for every $p\in \mathbb{N}\setminus\{0,2\}$ the map $\phi_p$ defines a homeomorphism that is $C^{\infty}$ at every point of $\mathbb{R}^2$ except at the origin ($\phi_2$ is everywhere smooth). Despite its lack of differentiability at the origin (except when $p=2$),  $\phi_p$ acts on the set of half-lines based at the origin for every $p\in \mathbb{N}^*$. It follows that after blowing up the origin in $\mathbb{R}^2$, we obtain a surface diffeomorphic to $\mathbb{S}^1\times [0,1)$ on which $\phi_p$ defines a $C^{\infty}$ diffeomorphism $\phi_p^*$. 

Trivially, the diffeomorphism $\phi_p^*$ restricted on $\mathbb{S}^1\times (0,1)$ is smoothly conjugated to $\phi_p$ restricted on $\mathbb{R}^2\setminus\{0\}$. More specifically, there exists a $C^{\infty}$ map $\Pi: \mathbb{S}^1\times [0,1)\rightarrow \mathbb{R}^2$ such that $\Pi(\mathbb{S}^1\times \{0\})=0$, $\Pi_{|\mathbb{S}^1\times (0,1)}: \mathbb{S}^1\times (0,1)\rightarrow \mathbb{R}^2-\{0\}$ is a diffeomorphism and $$\Pi\circ\phi^*_p = \phi_p\circ \Pi$$ 

\begin{figure}
    \centering
    \includegraphics[scale=0.12]{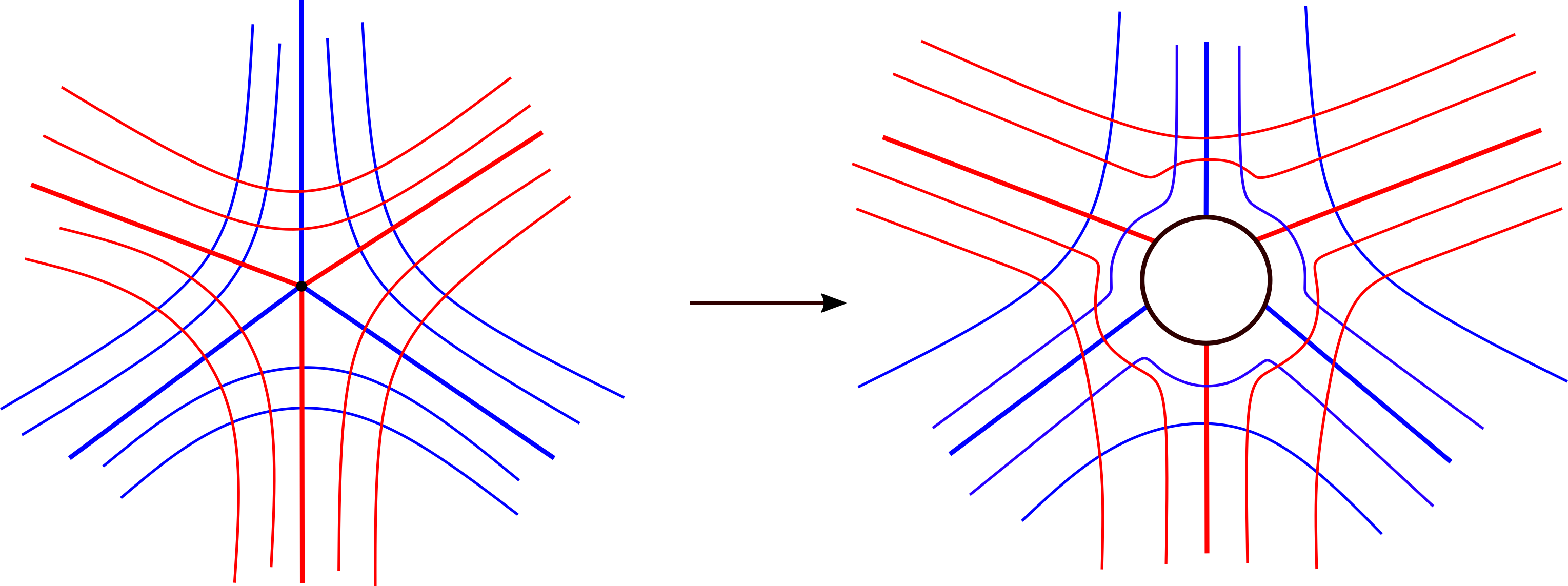}
    \caption{Foliations before and after blow-up}
    \label{f.folibeforeafter}
\end{figure}

 Let $$(\mathcal{F}^{s,u}_{p})^*:=\overline{\text{Conn}}(\Pi^{-1}(\mathcal{F}^{s,u}_{p}(0)))\cup \{\Pi^{-1}(L)|L\in\mathcal{F}^{s,u}_{p} \text{ and } 0\notin L \}$$ be the \emph{stable and unstable foliations of $\phi_p^*$}, where $\overline{\text{Conn}}(A)$ denotes the set of connected components of the closure of $A$. It is not hard to prove the following :
\begin{rema}\label{r.phistar}\quad
\newline{\quad}
\begin{itemize}
    \item $(\mathcal{F}^s_{p})^*,(\mathcal{F}^u_{p})^*$ are invariant by $\phi^*_p$ and form a pair of transverse foliations on $\mathbb{S}^1\times (0,1)$ (see Figure \ref{f.folibeforeafter})
    \item $\phi_p^*$ is a Morse-Smale diffeomorphism on $\mathbb{S}^1\times \{0\}$ with $2p$ hyperbolic fixed points ($p$ attracting and $p$ repelling). Each of the previous attracting (resp. repelling) fixed points is contained in a unique leaf of $(\mathcal{F}^u_{p})^*$ (resp. $(\mathcal{F}^s_{p})^*$), whose image by $\Pi$ is an unstable (resp. stable) prong of the origin in $(\mathbb{R}^2, \mathcal{F}^s_{p},\mathcal{F}^u_{p})$ 
\end{itemize}
\end{rema}

Let $R_{\theta}^*$ be the rotation of angle $2\pi\theta$ on $\mathbb{S}^1\times [0,1)$, or equivalently the unique diffeomorphism on $\mathbb{S}^1\times [0,1)$ such that $\Pi\circ R_{\theta}^*=R_{\theta}\circ \Pi$. We define for every $p\geq 1$ and $k\in \llbracket 0, p-1\rrbracket$ $$\phi^*_{p,k}:=\phi^*_p\circ R_{k/p}^*= R_{k/p}^*\circ\phi^*_p$$ Consider now the smooth manifold defined as follows: 
$$N^*_{p,k}:= \frac{\mathbb{S}^1\times [0,1)\times \mathbb{R}}{((x,y),t+1) \sim (\phi_{p,k}^*(x,y),t) }$$
endowed with the constant speed vertical flow $\Phi_{p,k}^*=\big((\Phi^*_{p,k})^t\big)_{t\in \mathbb{R}}$ given by the vector field $\frac{\partial}{\partial t}$. We will call $(N^*_{p,k},\Phi_{p,k}^*) $ the \emph{flow obtained from $\Phi_{p,k}$ after blowing up $\gamma_{p,k}$}. 

Just as in the case of $\phi_p$ and $\phi_p^*$, the flow $\Phi_{p,k}^*$ restricted on the interior of $N_{p,k}^*$ is $C^{\infty}$ conjugated to $\Phi_{p,k}$ restricted on $N_{p,k}-\gamma_{p,k}$. More specifically, thanks to $\Pi$, we can define a $C^{\infty}$ map $\Pi^*: N^*_{p,k}\rightarrow N_{p,k}$ such that $\Pi^*_{|\inte{N_{p,k}^*}}: \inte{N_{p,k}^*}\rightarrow N_{p,k}-\gamma_{p,k}$ is a $C^{\infty}$ diffeomorphism and for every $t\in \mathbb{R}$ 
\begin{equation} \label{eq.conjugationflowsblowup}
    \Pi^*\circ(\Phi^*_{p,k})^t= \Phi_{p,k}^t\circ \Pi^*
\end{equation}
We will call $\Pi^*$ the \emph{blow-down map associated to the blow-up of $\gamma_{p,k}$}. Moreover, similarly to $\Phi_{p,k}$, the flow $\Phi^*_{p,k}$ preserves the projections of $(\mathcal{F}^s_p)^*\times \mathbb{R}, (\mathcal{F}^u_p)^*\times \mathbb{R}$ on $N^*_{p,k}$, which we will denote respectively by $(F^s_{p,k})^*,(F^u_{p,k})^*$ and we will call the \emph{stable and unstable foliations of $\Phi_{p,k}^*$}. Finally, by our construction and thanks to Remark \ref{r.phistar}, it is easy to see that: 
\begin{rema}\label{r.conjugation1}
\quad \quad
\begin{itemize}
    \item the flow $\Phi^*_{p,k}$ restricted on the boundary of $N^*_{p,k}$ is a Morse-Smale flow with $gcd(k,p)$ (gcd stands for the greatest common divisor) attracting and $gcd(k,p)$ repelling hyperbolic periodic orbits
    \item $(F^s_{p,k})^*,(F^u_{p,k})^*$ define a pair of transverse foliations in the interior of $N^*_{p,k}$
   \item $(F^{s,u}_{p,k})^*=\text{Conn}((\Pi^*)^{-1}(F^{s,u}_{p,k}(\gamma_{p,k}))\cup \text{Per}_{\partial N^*_{p,k}})\cup \{(\Pi^*)^{-1}(L)|L\in F^{s,u}_{p,k} \text{ and } \gamma_{p,k}\not \subset L \}$, where $\text{Conn}(A)$ denotes the set of connected components $A$ and $\text{Per}_{\partial N^*_{p,k}}$ the union of all periodic orbits of $\Phi^*_{p,k}$ on $\partial N^*_{p,k}$
\end{itemize}
\end{rema}
\subsubsection{Blow-up operation for pseudo-Anosov flows}\label{s.blowupflows}

\quad

Consider $M$ an orientable, smooth, closed 3-manifold, $\Phi=(\Phi^t)_{t\in \mathbb{R}}$ a pseudo-Anosov flow on $M$, $F^s,F^u$ its stable and unstable foliations and $\gamma$ a periodic orbit of $\Phi$. Recall that, thanks to the orientability hypothesis, the tubular neighborhoods of $\gamma$ are solid tori.  

Consider $(N_{p,k},\Phi_{p,k})$ the local model for a pseudo-hyperbolic periodic orbit $\gamma_{p,k}$ with $p\geq 2$ prongs,  rotation $k \in \llbracket 0, p-1 \rrbracket$ and positive orientation. By Proposition \ref{p.aroundcircleprong}, there exist $p\geq 2$, $k \in \llbracket 0, p-1 \rrbracket$, $U_{\gamma}$ a neighborhood of $\gamma$ homeomorphic to a solid torus, $V_{p,k}$ a neighborhood of $\gamma_{p,k}$ and $H: U_{\gamma}\rightarrow V_{p,k}$ a homeomorphism defining an orbital equivalence between $(\Phi,U_\gamma)$ and $(\Phi_{p,k},V_{p,k})$ and sending stable/unstable leaves in  $U_\gamma$ to stable/usntable leaves in $V_{p,k}$. 

Let $(N^*_{p,k},\Phi^*_{p,k})$ be the flow obtained from $\Phi_{p,k}$ by blowing up $\gamma_{p,k}$, $\Pi^*:N^*_{p,k}\rightarrow N_{p,k}$ be the blow-down map associated to the blow-up of $\gamma_{p,k}$ and $V^*_{p,k}:=(\Pi^*)^{-1}(V_{p,k})$. Thanks to Equation \ref{eq.conjugationflowsblowup}, the map $H^{-1}\circ\Pi^*$ defines an orbital equivalence between the flows $(\Phi_{p,k}^*,V^*_{p,k}-\partial N^*_{p,k})$ and $(\Phi,U_\gamma-\gamma)$. Furthermore, thanks to Remark 
\ref{r.conjugation1}, for any $x\in V^*_{p,k}$ the map $H^{-1}\circ\Pi^*$ takes the connected component of $(F^{s,u}_{p,k})^*(x)\cap (V^*_{p,k}-\partial N^*_{p,k})$ containing $x$ to the connected component of $F^{s,u}(H^{-1}\circ\Pi^*(x))\cap (U_{\gamma}-\gamma)$ containing $H^{-1}\circ\Pi^*(x)$.  

Consider $M^*$ the $C^0$ manifold obtained by glueing $V^*_{p,k}$ and $M-\inte{U_\gamma}$ along $(\Pi^*)^{-1}(\partial V_{p,k})$ and $\partial U_\gamma$ via the map $H^{-1}\circ\Pi^*$. We will call $M^*$ the \emph{manifold obtained after blowing up $\gamma$}. By using the previous properties of $H^{-1}\circ\Pi^*$, we get that the flows $(\Phi_{p,k}^*,V^*_{p,k})$ and $(\Phi,M-U_\gamma)$ define on $M^*$ a $C^0$ orientable foliation of dimension $1$ that can be parametrized by a flow $\Phi^*=((\Phi^*)^t)_{t\in \mathbb{R}}$, thanks to Theorem 27A of \cite{Whitney}. We will call $\Phi^*$ a \emph{flow obtained after blowing up $\gamma$}. Similarly, by glueing the leaves of the foliations  $(F^{s,u}_{p,k})^*$ on $ V^*_{p,k}$ with the leaves of the foliations $F^{s,u}$ on $M-U_\gamma$ we define on $M^*$ two new sets of leaves, that will be denoted by $(F^s)^*$ and  $(F^u)^*$ and that will be called respectively the \emph{stable and unstable foliations of $(M^*,\Phi^*)$}. It is not difficult to see that, thanks to Remark \ref{r.conjugation1}: 
\begin{rema}\label{r.conjugation2}
\quad \quad 
    \begin{enumerate}
        \item $(F^{s,u})^*$ define a pair of transverse foliations in the interior of $M^*$, which are invariant by $\Phi^*$
        \item $gcd(k,p)$ leaves of $(F^s)^*$ and $gcd(k,p)$ leaves of $(F^u)^*$ intersect $\partial M^*$. Each of the previous intersections corresponds to a unique periodic orbit of $\Phi^*$ on $\partial M^*$
        \item contrary to $\Phi^*_{p,k}$, even if we fix our choices of $U_\gamma$ and $H$, a flow obtained after blowing up $\gamma$ is uniquely defined only up to reparametrization
    \end{enumerate}
\end{rema}

As in the case of hyperbolic models, we can define a blow-down map for $M^*$. Consider $\Pi^*_M: M^*\rightarrow M$ the map defined as the identity on $M-U_\gamma\subset M^*$ and as $H^{-1}\circ \Pi^*$ on $ V^*_{p,k}\subset M^*$. This map is trivially continuous and will be called the \emph{blow-down map associated to $(M^*,\Phi^*)$}. By construction, 
\begin{rema}\label{r.conjugationpstar}
\quad

\begin{enumerate}
\item $(F^{s,u})^*=\text{Conn}((\Pi^*_M)^{-1}(F^{s,u}(\gamma))\cup \text{Per}_{\partial M^*})\cup \{(\Pi^*_M)^{-1}(L)|L\in F^{s,u} \text{ and } \gamma\not \subset L \}$, where $\text{Conn}(A)$ denotes the set of connected components $A$ and $\text{Per}_{\partial M^*}$ the union of all the periodic orbits of $\Phi^*$ on $\partial M^*$
 \item $\Pi^*_M$ defines an orbital equivalence between $(\inte{M^*},\Phi^*)$ and  $(M-\gamma,\Phi)$
 \end{enumerate}
\end{rema}

\subsubsection{Dehn-Goodman-Fried surgery for pseudo-Anosov flows}\label{s.dehn-frieddefi}
\quad

Consider $M$ an orientable, smooth, closed 3-manifold, $\Phi=(\Phi^t)_{t\in \mathbb{R}}$ a pseudo-Anosov flow on $M$, $F^s,F^u$ its stable and unstable foliations and $\gamma$ a periodic orbit of $\Phi$. Let $(M^*,\Phi^*)$ be a flow obtained from $\Phi$ after blowing up $\gamma$. Denote by $\Phi^*_{|\partial M^*}$ the restriction of $\Phi^*$ on $\partial M^*$. 

We will say that a foliation by circles $\mathcal{F}$ on $\partial M^*\cong \mathbb{T}^2$ is \emph{collapsible} if it is invariant by $\Phi^*_{|\partial M^*}$ and every leaf of $\mathcal{F}$ is a \emph{global section} of $\Phi^*_{|\partial M^*}$ (i.e. every leaf $l$ of $\mathcal{F}$ is topologically transverse to $\Phi^*_{|\partial M^*}$ and intersects every orbit of $\Phi^*_{|\partial M^*}$) on which the orbits of $\Phi^*_{|\partial M^*}$ return in a uniformly bounded time. If such a foliation exists for $\Phi^*$, then by crushing its leaves to points, we can define a new blow-down operation giving rise to a new flow on a closed manifold. We will call this operation a Dehn-Goodman-Fried surgery. 

The next lemma that was proven in \cite{surgeriesnontransitive} states  that by eventually reparametrizing $\Phi^*$ close to $\partial M^*$, $\Phi^*$  always admits collapsible foliations. 

\begin{lemm}[\cite{surgeriesnontransitive}]\label{l.goodfoli}
   Consider $P$ the homotopy (or homology) class of any periodic orbit of $\Phi^*$ on $\partial M^*$. For every indivisible element $\sigma\in \pi_1(\partial M^*)$ such that $\sigma \neq \pm P$, there exists $\widetilde{\Phi}^*$, a flow obtained by reparametrizing $\Phi^*$ close to $\partial M^*$, such that $\widetilde{\Phi}^*$ admits a collapsible foliation on $\partial M^*$, whose every (parametrized) leaf is a circle of homotopy (or homology) class $\pm\sigma$. 
\end{lemm}

Choose now any indivisible class $\sigma\neq \pm P\in \pi_1(\partial M^*)$. By the above lemma, there exists $\widetilde{\Phi^*}$ a reparametrization of $\Phi^*$ close to $\partial M^*$, such that $\widetilde{\Phi^*}$ admits a  collapsible foliation on $\partial M^*$, say $\mathcal{F}_\sigma$, whose (parametrized) leaves are circles of homotopy class $\pm\sigma$. Let $$M_{\mathcal{F}_\sigma}:=\quotient{M^*}{\mathcal{F}_\sigma}$$ be the manifold obtained from $M^*$ by identifying every leaf of $\mathcal{F}_\sigma$ to a point. Consider $\Pi_{\mathcal{F}_\sigma}$ the natural projection from $M^*$ to $M_{\mathcal{F}_\sigma}$. The map $\Pi_{\mathcal{F}_\sigma}$ is clearly continuous and will be called the \emph{blow-down map associated to $\mathcal{F}_\sigma$}. Since $\mathcal{F}_\sigma$ is $\widetilde{\Phi^*}$ invariant, the flow $\widetilde{\Phi^*}$ projects to a $C^0$ flow $\Phi_{\mathcal{F}_\sigma}=(\Phi_{\mathcal{F}_\sigma}^t)_{t\in \mathbb{R}}$ on $M_{\mathcal{F}_\sigma}$ for which $$\Pi_{\mathcal{F}_\sigma}\circ(\widetilde{\Phi^*})^t= \Phi_{\mathcal{F}_\sigma}^t\circ \Pi_{\mathcal{F}_\sigma}$$

Notice that $\Phi_{\mathcal{F}_\sigma}$ was constructed thanks to a choice of a flow $\Phi^*$ obtained after blowing-up of $\gamma$, a choice of a reparametrization $ \widetilde{\Phi^*} $ of $\Phi^*$ and a choice of foliation $\mathcal{F}_\sigma$ preserved by $\widetilde{\Phi^*}$. Despite the previous facts, it is not difficult to prove that
\begin{prop}
    The orbital equivalence class of $(M_{\mathcal{F}_\sigma},\Phi_{\mathcal{F}_\sigma})$ depends only on $\pm\sigma$ and not our choice of $\Phi^*$ or $\widetilde{\Phi^*}$ or $\mathcal{F}_\sigma$.
\end{prop}
This justifies the following definition: 
\begin{defi}\label{d.defisurgery}
    We will say that the flow $(M_{\mathcal{F}_\sigma},\Phi_{\mathcal{F}_\sigma})$ is obtained from $(M,\Phi)$ by the \emph{Dehn-Goodman-Fried surgery along $\gamma$ associated to the homotopy class $\sigma$ (or $-\sigma$)}.
\end{defi}

\begin{rema}\label{r.surgery}
\begin{itemize}

    \item Let us remark immediately that the flows obtained by Dehn-Goodman-Fried surgeries along $\gamma$ associated to the holonomy classes $\sigma$ and $-\sigma$ are orbitally equivalent. 
    \item Thanks to Proposition \ref{p.aroundcircleprong}, any periodic orbit with non-orientable tubular neighborhoods admits tubular neighborhoods with a Klein bottle boundary. There exists a unique (up to homotopy) foliation on the Klein bottle with no one sided leaf. It follows that we can only perform a trivial Dehn-Goodman-Fried surgery on a periodic orbit with non-orientable tubular neighborhoods. 
\end{itemize}
\end{rema}

In most of the cases, the flows obtained from Dehn-Goodman-Fried surgeries on pseudo-Anosov flows remain pseudo-Anosov :

\begin{theorem}[\cite{surgeriesnontransitive}]\label{t.surgerygivespseudoanosov}
    Let $\Phi$ be a pseudo-Anosov flow on a smooth, orientable, 3-manifold $M$. Consider $\gamma$ a periodic orbit of $\Phi$, $(M^*, \Phi^*)$ a flow obtained after blowing-up $\Phi$ along $\gamma$ together with its stable and unstable foliations $(F^s)^*, (F^u)^*$.

    Let $P\in \pi_1(\partial M^*)$ be the homotopy (or homology) class of a periodic orbit of $\Phi^*$ on $\partial M^*$ and $\sigma$ any indivisible class inside $\pi_1(\partial M^*)$ such that $\sigma\neq \pm P$. Denote by $(M_\sigma,\Phi_\sigma)$ the flow obtained from $\Phi$ after a Dehn-Goodman-Fried surgery along $\gamma$ associated to the homotopy class $\sigma$ and by $\Pi_\sigma: M^* \rightarrow M_\sigma$ the blow-down map associated to the previous surgery. We have that: 
    \begin{enumerate}
        \item the projections by $\Pi_\sigma$ of $(F^s)^*, (F^u)^*$ define a pair of singular transverse foliations, denoted by $F_\sigma^s,F_\sigma^u$
        \item if $\gamma_{\sigma}:=\Pi_\sigma(\partial M^*)$, then the flows $(M-\gamma,\Phi)$ and $(M_\sigma-\gamma_\sigma,\Phi_\sigma)$ are orbitally equivalent
    \end{enumerate}
Assume that $\gamma$ is associated by Proposition \ref{p.aroundcircleprong} to a local model for a pseudo-hyperbolic periodic orbit with $p\geq 2$ prongs and rotation $k\in \llbracket 0,p-1 \rrbracket$. Let also $K$ be the absolute value of the algebraic intersection number of the two homotopy classes $P$ and $\sigma$.
    \begin{enumerate}
        \setcounter{enumi}{2}
        \item $\gamma_{\sigma}=\Pi_\sigma(\partial M^*)$ is a circle $K\times gcd(k,p)$-prong singularity of $F^{s,u}_{\sigma}$, where $gcd(k,p)$ is the greatest common divisor of $k$ and $p$
        \item $\Phi_\sigma$ is a pseudo-Anosov flow whose stable and unstable foliations are respectively $F_\sigma^s,F_\sigma^u$ if and only if $\gamma_{\sigma}$ is not a $1-$prong singularity.
        \item there exists a flow obtained from $\Phi_\sigma$ by a Dehn-Goodman-Fried surgery along  $\gamma_\sigma$ that is orbitally equivalent to $\Phi$ 
    \end{enumerate}

\end{theorem}
By construction, a Dehn-Goodman-Fried surgery modifies not only the behavior of the flow around a periodic orbit, but also the topology of the underlying manifold. Even though the orbits of the flows before and after surgery are in bijective correspondance (see the following definition) and the flows outside the orbit on which the surgery was performed  are orbitally equivalent, the flows before and after surgery can have completely different behaviors. 
\begin{defi}\label{d.associatedorbits}
    Let $\Pi_M^{*}: M^*\rightarrow M$ and $\Pi_{\sigma}: M^*\rightarrow M_\sigma$ be the blow-down maps associated to $(M^*,\Phi^*)$ and to $\mathcal{F}_\sigma$ respectively. For any orbit $\mathcal{O}$ of $\Phi$, we will say that the orbit of $\Phi_{\sigma}$ given by $\Pi_{\sigma}\circ (\Pi_M^{*})^{-1}(\mathcal{O})$ is the orbit of $\Phi_{\sigma}$ \emph{corresponding to $\mathcal{O}$} after surgery. 
\end{defi}

\subsubsection{Describing surgeries by pairs of integers }\label{s.surgeryaspairsof integers}
\quad

Just as in the case of Dehn surgeries on $3-$manifolds, it will be practical for us to describe a Dehn-Goodman-Fried surgery by a pair of integers instead of an element in the fundamental group of the torus obtained after blowing-up a periodic orbit. 

Consider $M$ an orientable, smooth, closed 3-manifold, $\Phi=(\Phi^t)_{t\in \mathbb{R}}$ a pseudo-Anosov flow on $M$, $F^s,F^u$ its stable and unstable foliations and $\gamma$ a periodic orbit of $\Phi$. Let $(M^*,\Phi^*)$ be a flow obtained from $\Phi$ after blowing up $\gamma$. By eventually changing the manifold atlas of $M^*$, we may assume without any loss of generality that $M^*$ is smooth. Also, using the fact that $M$ is orientable, it is easy to see that $M^*$ is also orientable. Fix for the rest of this section an arbitrary orientation on $M^*$. Using this orientation, we will define a canonical base in $\pi_1(\partial M^*) \approx \mathbb{Z}^2$, thanks to which Dehn-Goodman-Fried surgeries on $\gamma$ will be associated to pairs of integers. 

\vspace{0.2cm}
\textit{A first canonical element in $\pi_1(\partial M^*)$: the meridian.}

Take $X\in \partial M^*$, $\nu\in T_XM^*$ a vector pointing towards the interior of $M^*$ and $(u,v)$ a pair of linearly independent vectors in $T_X(\partial M^*)\subset T_XM^*$. We choose an orientation on $\partial M^*$ such that $(u,v)$ is a positively oriented frame in $\partial M^*$ if and only if the triplet $(\nu, u, v)$ defines an orientation on $M^*$ compatible with our original choice of orientation. Notice that the previous orientation on $\partial M^*$ does not depend on our choice of $X, u, v, \nu$, but only on our original choice of orientation of $M^*$.  

Denote by $i_g(\gamma_1,\gamma_2)$ (resp. $i_a(\gamma_1,\gamma_2)$) the geometric (resp. algebraic) number of intersections of two curves (resp. curves or homotopy classes) in $\partial M^*$. We remark here that $i_a$ is well defined thanks to our choice of orientation of $\partial M^*$.

Consider now any periodic orbit $\gamma_P$ of $\Phi^*$ on $\partial M^*$ endowed with its dynamical orientation and denote by $P$ the homotopy (or homology) class of $\gamma_P$ (this class does not depend on our choice of periodic orbit on $\partial M^*$). Next, let $\Pi^*_M: M^* \rightarrow M$ be the blow-down map associated to $(M^*,\Phi^*)$. The restriction of $\Pi^*_M$ on $\partial M^*$ defines a fibration by circles topologically transverse to $\Phi^*$ over the circle $\gamma$. Take $\gamma_m$ one of the previous circles and choose an orientation on $\gamma_m$ so that $i_a(\gamma_P,\gamma_m)> 0$. Denote by $m$ the homotopy (or homology) class of $\gamma_m$ endowed with the previous orientation. We are going to call $m$ the \emph{meridian} of $\partial M^*$ (the class $m$ does not depend on our choice of $\gamma_P$ or $\gamma_m$, but only on our initial choice of orientation of $M^*$).

\begin{figure}
    \centering
    \includegraphics[scale=0.2]{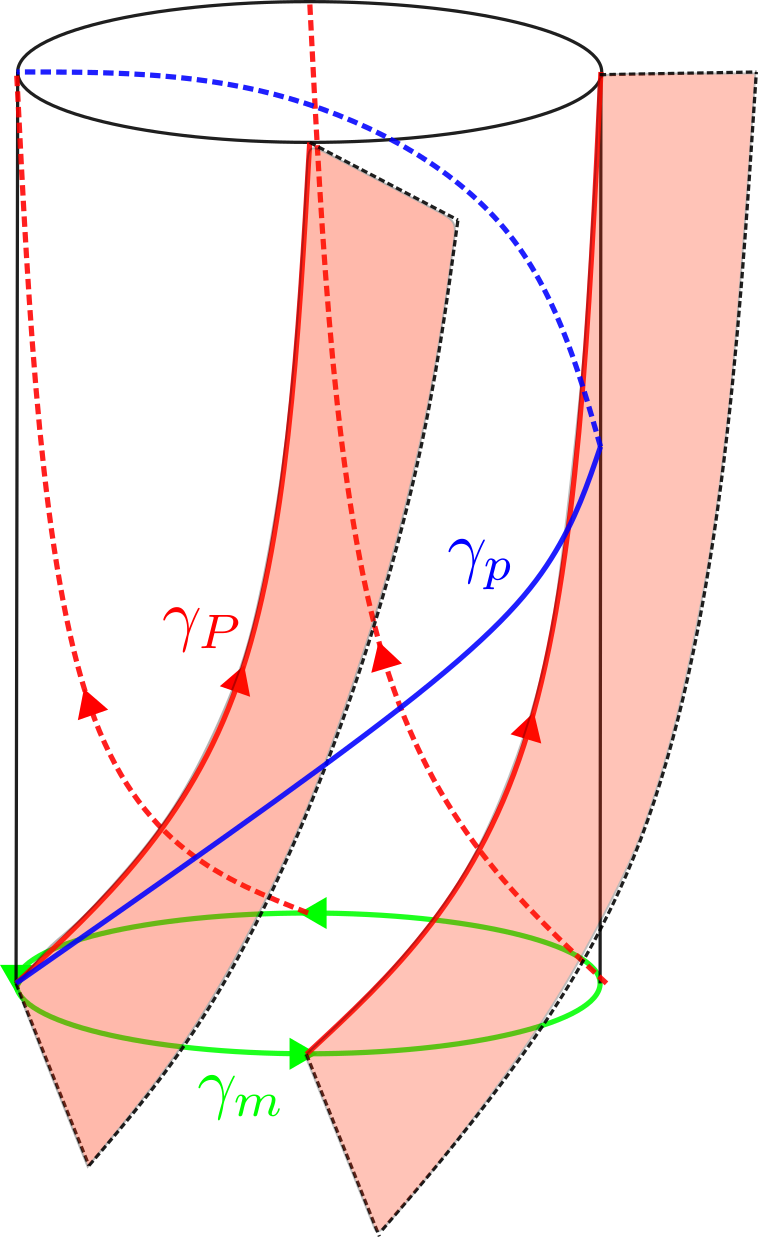}
    \caption{Blowing up a local model for a pseudo-hyperbolic orbit with $4$ prongs and rotation $1$}
    \label{f.parallelmeridian}
\end{figure}

Notice that $\gamma_m$ and $\gamma_P$ are simple closed curves, that intersect transversely  and that verify $$i_a(\gamma_P,\gamma_m)=i_g(\gamma_P,\gamma_m)>0$$ In other words, $\gamma_m$ and $\gamma_P$ intersect at least once and performing a homotopy on either $\gamma_m$ or  $\gamma_P$ can only increase the number of intersections between the two curves. Even though $m$ and $P$ are linearly independent in $\pi_1(\partial M^*) \approx \mathbb{Z}^2$, they do not always form a basis of $\pi_1(\partial M^*)$, as $\gamma_m$ and $\gamma_P$ can intersect multiple times (see Figure \ref{f.parallelmeridian}). More precisely, if $\gamma$ is associated by Proposition \ref{p.aroundcircleprong} to a local model for a pseudo-hyperbolic orbit with $n\geq 2$ prongs and rotation $r \in \llbracket 0,p-1\rrbracket$, then by using Remark \ref{r.conjugation2}, we get that the total number of intersections between $\gamma_P$ and $\gamma_m$ is $\frac{n}{gcd(n,r)}$. For example, when $r=0$ we get that $P,m$ form a basis of $\pi_1(\partial M^*)$ and when $r=1$ we get that $P,m$ do not form a basis of $\pi_1(\partial M^*)$, as $\gamma_m$ and $\gamma_P$ intersect $n$ times.

\vspace{0.2cm}
\textit{A second canonical element in $\pi_1(\partial M^*)$: the parallel.}

By construction, $\gamma_m$ intersects transversely along $\frac{n}{gcd(n,r)}$ points each of the $gcd(n,r)$ periodic orbits of $\Phi^*$ on $\partial M^*$. Consider $x_0,...,x_{n-1}$ the points of intersection of the periodic orbits of $\Phi^*$ on $\partial M^*$ and $\gamma_m$. By eventually reindexing the $x_i$, we can assume that $x:=x_0$ is a point of intersection between $\gamma_P$ and $\gamma_m$ and that when following positively $\gamma_m$ starting from $x_i$ the first point in $\{x_1,...,x_{n-1}\}$ that we encounter is $x_{i+1}$ (the indexes are considered modulo $n$). Following positively $\gamma_P$ starting from $x$, let $x_k$ be the first point of intersection between $\gamma_P$ and $\gamma_m$ after $x$, where $k\in \llbracket 0, n-1\rrbracket$. Let us remark here that 
\begin{itemize}
    \item $k=r$ or $k=n-r$ depending on our choice of orientation of $M^*$
    \item following positively $\gamma_P$ starting from $x_i$, the point  $x_{i+k}$ (the indexes are considered modulo $n$) is the first point of intersection between $\gamma_P$ and $\gamma_m$ after $x_i$
\end{itemize}

Consider $y:=x_{n-k}$,  $[y,x]_P$ the segment contained in $\gamma_P$ going from $y$ to $x$ (following $\gamma_P$ positively) and  $[x,y]_m$ the segment  contained in $\gamma_m$ going from $x$ to $y$ (following $\gamma_m$ positively). By convention, if $x=y$ (and thus $k=0$), we take $[x,y]_m=\gamma_m$ and $[x,y]_P=\gamma_P$. The juxtaposition of $[x,y]_m$ followed by $[y,x]_P$ defines a closed curve, whose homotopy (or homology) class will be denoted by $p$. Notice that by performing a homotopy on the previous curve we can obtain a curve $\gamma_p$ intersecting $\gamma_m$ once (see Figure \ref{f.parallelmeridian}). It follows that $m$ and $p$ form a basis of $\pi_1(\partial M^*)$. We will call $p$ the \emph{parallel} of $\partial M^*$. 

\begin{lemm}\label{l.parallelmeridian}
    We have that $$ P+\frac{(n-k)\cdot m}{gcd(n,k)}=\frac{n\cdot p}{gcd(n,k)} $$ where $gcd(n,k)$ denotes the greatest common divisor of $n$ and $k$ or equivalently the number of closed orbits of $\Phi^*$ on $\partial M^*$. 
\end{lemm}
\begin{proof}
The second part of the lemma is an immediate consequence of our definition of $k$. Let us now prove the above relation between $m,p$ and $P$. Following the previous notations, the curve $[x,y]_m$ performs a rotation around the torus of ``angle" $\frac{2\pi (n-k)}{n}$ in the direction of the meridian. Furthermore, since the total number of intersections between $\gamma_P$ and $\gamma_m$ is $\frac{n}{gcd(n,k)}$, the curve $[y,x]_P$ allows us to traverse only a fraction of the total ``length" of $\gamma_P$, more precisely $\frac{gcd(n,k)}{n}$ of its total ``length". Using the previous facts, it is not difficult to see there exists a loop representative inside the homotopy class $\frac{n}{gcd(n,k)}\cdot p$ that starts from $x=x_0$ goes to $x_{k}$ following positively $\gamma_P$, then goes from $x_{k}$ to $x_{2k}$ (the indexes are considered modulo $n$) following positively $\gamma_P$, until it reaches $x_0$, thus following positively once the curve $\gamma_P$, and at the end performs $\frac{n}{gcd(n,k)}\times \frac{(n-k)}{n}=\frac{(n-k)}{gcd(n,k)}\in \mathbb{N}$ turns around the torus following positively $\gamma_m$. This finishes the proof of the lemma. 
\end{proof}

Using the fact that $m$ and $p$ form a basis of $\partial M^*$, we can now describe all Dehn-Goodman-Fried surgeries on $\gamma$ by pairs of integers that are relatively prime. More specifically, let $\sigma$ be an indivisible class inside $\pi_1(\partial M^*)$ such that $\sigma \neq \pm P$. There exists a unique pair of relatively prime integers $(s,t)\in \mathbb{Z}^2$ such that $\sigma=s\cdot m + t\cdot p$. 

\begin{defi}
    Let $\Phi'$ be the flow obtained from $\Phi$ by a Dehn-Goodman-Fried surgery along $\gamma$ associated to the homotopy class $\sigma$. In addition to its previous name, we will call $\Phi'$ the flow obtained from $\Phi$ by an \emph{$(s,t)$-Dehn-Goodman-Fried surgery on $\gamma$} or an \emph{$(s,t)-$surgery on $\gamma$}. 
\end{defi}

\begin{rema}
\begin{itemize}
    \item Thanks to Remark \ref{r.surgery}, the flows obtained from $\Phi$ by an $(s,t)$-surgery or an $(-s,-t)$-surgery on $\gamma$ are orbitally equivalent
    \item Recall that by construction, if a periodic orbit $\gamma$ of $\Phi$ is associated by Proposition \ref{p.aroundcircleprong} to a pseudo-hyperbolic model with $n\geq 2$ prongs and rotation $k\in \llbracket 0, n-1\rrbracket$, then thanks to Lemma \ref{l.parallelmeridian}, the $(\pm\frac{(n-k)}{gcd(n,k)}, \mp \frac{n}{gcd(n,k)})-$surgery on $\gamma$ is not well-defined 
    \item Even if a surgery is well defined, it is possible because of Theorem \ref{t.surgerygivespseudoanosov}, that it leads to the construction of a circle $1$-prong singularity and thus to the construction of a flow that is not pseudo-Anosov 
\end{itemize}
\end{rema}
\section{Markovian families of rectangles}\label{s.markovianfamilies}
From this moment and on, we will begin the exposition of our approach for classifying pseudo-Anosov flows in dimension 3. This section is dedicated to the proof of Theorem A. In order to make our results as general as possible, but also as a means to provide some evidence supporting Conjecture \ref{conj.markovianactions}, we will prove the following generalization of Theorem-Definition A:

 \textbf{Theorem-Definition A'.}\textit{
Let $\rho$ be any  orientation preserving strong Markovian action on the plane (see Definition \ref{d.markovianaction}) leaving invariant a strong Markovian family $\mathcal{R}$ (see Definition \ref{d.markovianaction}). We can canonically associate to $\mathcal{R}$ a finite set of pairwise equivalent geometric types, called the geometric types of $\mathcal{R}$ or the geometric types associated to $\mathcal{R}$.}     
\vspace{0.20cm}
 
For a more precise statement, see Definition \ref{d.geometrictypemarkovfamily} and Theorem \ref{t.associatemarkovfamiliestogeometrictype}.  

We remind the reader that thanks to Remark \ref{r.pseudoanosovactionsarestrongmarkovian}, any Markovian family associated to a pseudo-Anosov flow is strong and that the action by the fundamental group on the bifoliated plane of any pseudo-Anosov flow is a strong Markovian action; hence, Theorem-Definition A' indeed generalizes Theorem-Definition A.

\subsection{Some general results for strong Markovian actions and families}
Fix $\mathcal{P}$ a plane, $\rho: G\rightarrow \text{Homeo}(\mathcal{P})$ an orientation preserving strong Markovian action, leaving invariant a strong Markovian family $\mathcal{R}$ and preserving the pair of singular foliations $\mathcal{F}^s$ and $\mathcal{F}^u$, that we will call respectively the stable and unstable foliations of $\rho$. 
\begin{lemm}\label{l.periodicboundary}
The boundary of any $R\in \mathcal{R}$ consists of stable/unstable segments belonging to stable/unstable leaves with non-trivial stabilizers in $G$.  
\end{lemm}
\begin{proof}
Indeed, take $r_0\in \mathcal{R}$ and consider one of its stable boundary segments, say $s$. Let us denote by $S$ the stable leaf in $\mathcal{F}^s$ containing $s$.

Take $x\in s$. Since $\mathcal{R}$ satisfies the results of Lemma \ref{l.horizsubrectangleexists}, there exists $r_1\in \mathcal{R}$ containing $x$ such that $r_1 \cap r_0$ is a non-trivial horizontal subrectangle of $r_0$. Therefore, $r_1$ contains $s$. By repeatedly applying  this argument, we can construct $r_0,r_1,...,r_n,...$ a sequence of rectangles in $\mathcal{R}$ containing $s$ in their boundaries and such that $r_n \cap r_{n+1}$ is a non-trivial horizontal subrectangle of $r_n$. 

By the finiteness axiom, there exist $i,j$ two distinct integers and $g\in \pi_1(M)-\{id\}$ such that $g(r_i)=r_j$ and $g(S)=S$, which gives  the desired result.  
\end{proof}

The following lemma shows that as for any Markov partition, we can define a notion of first return map for any strong Markovian family : 
\begin{lemm}\label{l.existenceofpredecessors}
For any rectangle $R \in \mathcal{R}$ there exists a unique finite collection of rectangles $R_1,...,R_n \in \mathcal{R}$ intersecting $R$ along non-trivial vertical subrectangles and such that:
\begin{enumerate}
\item $R_1,...,R_n$ are maximal for the previous property: any $R' \in \mathcal{R}$ intersecting $R$ along a non-trivial vertical subrectangle satisfies $R' \cap R \subseteq R_i \cap R$ for some $i \in \llbracket 1, n \rrbracket$ 
\item $R_1,...,R_n$ have disjoint interiors 
\item The $R_1,...,R_n$ cover $R$: $\overset{n}{\underset{i=1}{\cup}} R_i \cap R = R$ 
\end{enumerate}

\end{lemm}
The analogue of the previous lemma for horizontal subrectangles is also true. 
\begin{proof} 
Fix $R \in \mathcal{R}$. Let's call the property of intersecting $R$ along a non-trivial vertical subrectangle, property $(\star)$. Since $\mathcal{R}$ verifies the results of Lemma \ref{l.vertsubrectangleexists}, by definition there exist rectangles in $\mathcal{R}$ satisfying $(\star)$. Let us begin by showing that for any point $x\in R$ there exists at least one rectangle maximal for $(\star)$ containing $x$. 

Indeed, let us fix $x\in R$. By definition of a strong Markovian family, there exists $r_0$ satisfying $(\star)$ containing $x$. Suppose there is no maximal rectangle for $(\star)$ containing $x$. Therefore, for any rectangle $S\in \mathcal{R}$ containing $x$ and  intersecting $R$ along a vertical subrectangle,  there exists $S'\in \mathcal{R}$ satisfying $(\star)$ such that $R \cap S \subsetneq R\cap S'$. We can thus construct by induction an infinite sequence $r_0,r_1,r_2,....,r_n,...$ of rectangles in $\mathcal{R}$ satisfying $(\star)$ and such that $r_k\cap R \subsetneq r_{k+1}\cap R$ for all $k\in \mathbb{N}$. By the Markovian intersection axiom, $r_k \cap r_{k+1}$ is a horizontal subrectangle of $r_k$ and therefore, since $\mathcal{R}$ verifies the results of Lemma \ref{l.infiniteintersectionhorizontalrectangles}, we get that $\overset{+\infty}{\underset{k=0}{\cap}} r_{k}$ is a stable segment of $r_0$. But $r_k\cap R \subsetneq r_{k+1}\cap R$ for every $k$, so $r_0\cap R \subset \overset{+\infty}{\underset{k=0}{\cap}} r_{k}\cap R$, which is impossible, since $\overset{+\infty}{\underset{k=0}{\cap}} r_{k}$ is a segment. We deduce the existence of a maximal rectangle for $(\star)$ containing $x$.

Next, let us prove that maximal rectangles for $(\star)$ are either identical or they have disjoint interiors. This will imply in particular that for any $S\in \mathcal{R}$ that satisfies $(\star)$,   the maximal rectangle for $(\star)$ containing $S\cap R$ is unique. Indeed, by the Markovian intersection axiom, two distinct rectangles $R_i,R_j$ in $\mathcal{R}$ satisfying $(\star)$ have either disjoint interiors or they satisfy one of the following: $R \cap R_i \subseteq R\cap R_j$ or $R \cap R_i \subseteq R\cap R_j$. One can easily check that the Markovian intersection axiom, according to which, up to interchanging $R_i$ and  $R_j$, the set $R_i \cap R_j$ is a   \emph{non-trivial} vertical subrectangle of $R_i$ and a \emph{non-trivial} horizontal subrectangle of $R_j$, implies that the case $R_i\cap R=R_j\cap R$ is impossible. We thus get the desired result.

Finally, let us prove that the set of maximal rectangles for $(\star)$ is finite. Suppose the contrary.  Under this hypothesis, there exists a sequence $(r_i)_{i\in \mathbb{N}}$ of maximal rectangles for $(\star)$ such that the $\overset{\circ}{r_i}$ are pairwise disjoint. By compactness, we can assume that the rectangles $r_i \cap R$ converge to an unstable segment of $R$, say $s$. Without any loss of generality, assume that the rectangles $r_i \cap R$ accumulate to $s$ from the right. Take $x\in s$ and suppose that $\mathcal{F}^s_-(x)$ is the stable separatrix on the right of $x$. Since Lemma \ref{l.vertsubrectangleexists} is true for $\mathcal{R}$, there exists $R' \neq R$ containing $x$ and a small neighborhood of $x$ inside $\mathcal{F}^s_-(x)$ such that $R'\cap R$ is a non-trivial vertical subrectangle of $R$ (hence $R'$ contains $s$). Since $R'$ intersects $\mathcal{F}^s_-(x)$ and contains $s$, it also contains all the points of $R$ on the right of $s$ and that are sufficiently close to $s$. Therefore, for $i$ sufficiently big, $r_i \cap R \subset R'\cap R$, which contradicts the fact that the $r_i$ are maximal. Therefore, the set of maximal rectangles for $(\star)$ is finite.

\end{proof}

\begin{defi}\label{d.successor}
For any $R\in \mathcal{R}$, we will say that $R'$ is a \emph{predecessor} (resp. \emph{successor}) of $R$ if $R'\cap R$ is a non-trivial vertical (resp. horizontal) subrectangle of $R$ and $R'$ is maximal for this property in the sense of the previous lemma. 

We will say that $R'$ is a predecessor of \emph{$2$-nd generation} of $R$, if $R'$ is a predecessor of a predecessor of $R$. We define similarly a predecessor (resp. successor) of \emph{$n$-th generation} for any $n \in \mathbb{N}^*$. By convention, a $1$-st generation predecessor (resp. successor) of $R$, is the same thing as a predecessor (resp. successor) of $R$. 

\end{defi}

\begin{figure}[h!]
\includegraphics[scale=0.4]{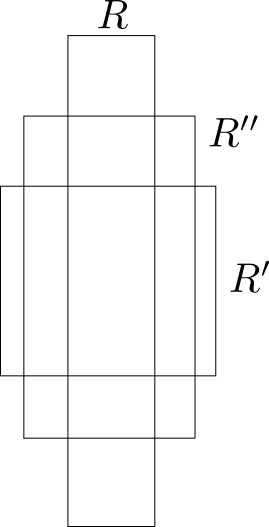}
\caption{}
\label{f.precedentsuivant}
\end{figure}

\begin{rema}\label{r.precedentsuivant}
If $R\in \mathcal{R}$ is a predecessor of $R'\in \mathcal{R}$ and $g\in G$, then 

\begin{itemize}
    \item $\rho(g)(R)$ is a predecessor of $\rho(g)(R')$ and 
    \item $R'$ is a successor of $R$ 
\end{itemize}
\end{rema}
The first statement is an easy consequence of Definition \ref{d.successor} and the fact that $\mathcal{R}$ is preserved by the action of $G$. Concerning the second statement, assume that $R''\neq R'$ is the successor of $R$ containing $R\cap R'$. By the Markovian intersection property, since $R$ is a predecessor of $R'$ and $R''$ is a successor of $R$, $R\cap R'$ is a non-trivial horizontal subrectangle of $R$ and $R\cap R''$ is a non-trivial vertical subrectangle of $R''$ (see Figure \ref{f.precedentsuivant}). Also, since $R''$ contains $R\cap R'$, we have that $$R\cap R' \subseteq R''\cap R'$$ and thus $R''$ intersects $R'$ along a non-trivial vertical subrectangle. As in our proof of Lemma \ref{l.existenceofpredecessors}, using the Markovian intersection axiom, the fact that   $R''\neq R$ and the fact that both $R$ and $R''$ intersect $R'$ along a non-trivial vertical subrectangle, one can easily show that $R\cap R' \subsetneq R''\cap R'$. It follows that  $R$ can not be a predecessor of $R'$, which leads to an absurd and gives the desired result.

\begin{lemm}\label{l.npredecessor}
Take any two $R, R' \in \mathcal{R}$ such that $R'\cap R$ is a non-trivial vertical subrectangle of $R$. We have that $R'$ is a predecessor of $n$-th generation of $R$ for some $n\in \mathbb{N}^*$. Furtheremore, there exists a unique sequence in $\mathcal{R}$, say $R_0=R, R_1,...,R_n=R'$, starting from $R$, ending at $R'$ and such that $R_{i+1}$ is a predecessor of $R_i$ for every $i\in \llbracket 1, n-1\rrbracket$.  

\end{lemm}
\begin{proof}
By maximality, if $R'$ is not one of the predecessors of $R$, then there exists $R_1$ a predecessor of $R$ containing $R'\cap R$. Notice that since the predecessors of $R$ have disjoint interiors, $R_1$ is the unique predecessor of $R$ containing $R'\cap R$. Furthermore, by the Markovian intersection property, $R'$ intersects $R_1$ along a non-trivial vertical subrectangle. Again, if $R'$ is not one of the predecessors of $R_1$, then $R'\cap R_1$ is contained in a unique predecessor of $R_1$, say $R_2$. We construct in this way a sequence $R_0=R,R_1,...,R_i,...$ such that for every $i$, the rectangle  $R_{i+1}$ is the predecessor of $R_i$ containing $R'\cap R_i$. If there exists no $i$ such that $R_i\cap R=R'\cap R$, the previous sequence is infinite and since Lemma \ref{l.infiniteintersectionverticalrectangles} is true for $\mathcal{R}$, $\overset{+\infty}{\underset{k=0}{\cap}} R_{k}$
is an unstable segment of $R$ containing $R'\cap R$, which is impossible. We conclude that there exists $n$ such that $R_n\cap R=R'\cap R$. As we have previously stated, the Markovian intersection property implies in this case that $R_n=R'$, which gives us the desired result.
\end{proof}
\begin{lemm}\label{l.intersectionsuccessorspredecessors}
    Take $R\in \mathcal{R}$ and $R', R''$ two  predecessors (resp. successors) of $R$ of  generations $n$ and $m$ with $n,m \in \mathbb{N}^*$. We have that $R'\cap R''\neq \emptyset$ if and only if $R'\cap R''\cap R \neq \emptyset$. Similarly, $\inte{R'}\cap \inte{R''}\neq \emptyset$ if and only if $\inte{R'}\cap \inte{R''}\cap R \neq \emptyset$.
\end{lemm}
\begin{proof}
    We will prove the first part of the lemma. The second follows from a similar argument. The fact that $R'\cap R''\cap R \neq \emptyset$ implies $R'\cap R''\neq \emptyset$ is obvious. Assume now that $R'\cap R''\neq \emptyset$.  This implies that there exists a stable and an unstable leaf in $\mathcal{F}^s$ and $\mathcal{F}^u$ respectively, crossing both $R'$ and $ R''$. If $R'\cap R''\cap R = \emptyset$, this would imply that one of two previous leaves intersects $R$ along two connected components, which is impossible by Item (5) of Proposition  \ref{p.propertiesoffoliformarkovianactions}. 
\end{proof}
\subsection{Geometric types of a strong Markovian family (Theorem-Definition A')} \label{s.prooftheoremA}
In this section, we will define the set of \emph{geometric types of a strong Markovian family $\mathcal{R}$} (see Definition \ref{d.geometrictypemarkovfamily}), a finite set of geometric types each characterizing the pattern of intersection of the rectangles in $\mathcal{R}$, and we will show that the geometric types of any strong Markovian family are pairwise equivalent (see Theorem \ref{t.associatemarkovfamiliestogeometrictype}). 

Fix $\mathcal{P}$ a plane, $\rho: G\rightarrow \text{Homeo}(\mathcal{P})$ an orientation preserving strong Markovian action, leaving invariant a strong Markovian family $\mathcal{R}$ and preserving the pair of singular foliations $\mathcal{F}^s$ and $\mathcal{F}^u$, that we will call respectively the stable and unstable foliations of $\rho$. 

\vspace{0.2cm}
\textit{Constructing a geometric type from $\mathcal{R}$}

Fix an orientation on $\mathcal{P}$, $\{r_1,...,r_n\} \subset \mathcal{R}$ a set of representatives of every rectangle orbit in $\mathcal{R}$ and choose an orientation of the stable and unstable foliations inside every $r_i$ so that the positive stable direction followed by the positive unstable direction inside every $r_i$ defines an orientation compatible with our choice of orientation of $\mathcal{P}$. 
 
Let $h_i$ (resp. $v_i$) be the number of successors (resp. predecessors) of the rectangle $r_i$. Associate to each $r_i$ a copy of $[0,1]^2$ trivially bifoliated by horizontal and vertical segments, that we will denote by $R_i$. Inside every $R_i$ consider $h_i$ (resp. $v_i$) pairwise disjoint horizontal (resp. vertical) subrectangles $H_i^1,...,H_i^{h_i}$ (resp. $V_i^1,...,V_i^{v_i}$) ordered from bottom to top (resp. from left to right). Thanks to our choice of orientations of the foliations inside the $r_i$, we can identify each successor (resp. predecessor) of $r_i$ with a rectangle in $\{H_i^1,...,H_i^{h_i}\}$ (resp. $\{V_i^1,...,V_i^{v_i}\}$).

\begin{figure}[h!]
    \centering
    \includegraphics[scale=0.4]{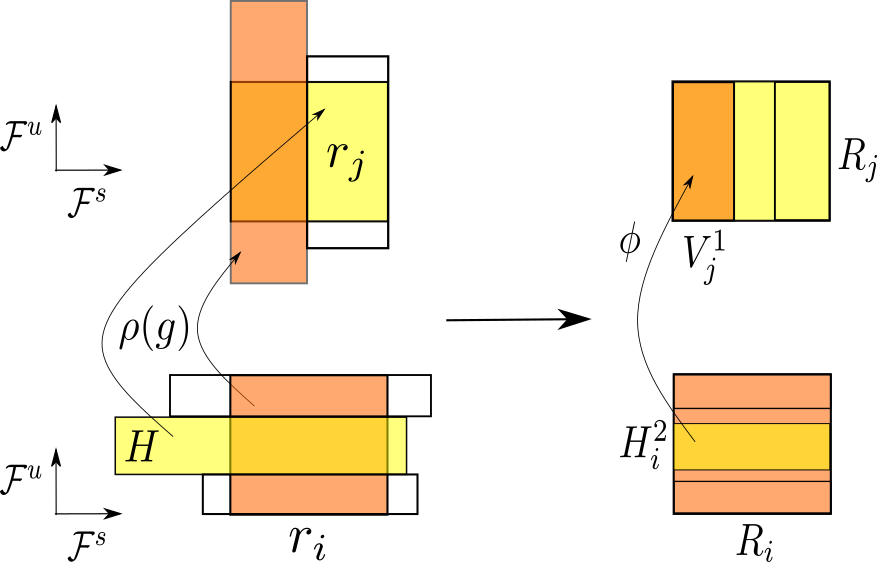}
    \caption{Associating to a Markovian family a geometric type}
    \label{f.canonicalgeom}
\end{figure}

Consider $H$ a successor of $r_i$ in the $G$-orbit of $r_j$ (see Figure \ref{f.canonicalgeom}). There exists $g\in G$ such that $\rho(g)(H)=r_j$. By Remark \ref{r.precedentsuivant}, $V:=\rho(g)(r_i)$ is a predecessor of $r_j=\rho(g)(H)$. Using our previous choices of orientations, if $H$ is the $k$-th successor of $r_i$ 
(for the bottom to top order) and $\rho(g)(r_i)$ is the $l$-th predecessor of $r_j$ (for the left to right
order), we define $\phi(H^k_i)=V_j^l$. Furthermore, since $H$ intersects $r_i$ along a horizontal subrectangle, the foliations inside $H$ inherit a natural orientation from our choice of orientations inside $r_i$. We define $u(H^k_i)=1$  (resp. $u(H^k_i):=-1$) if $\rho(g)$ preserves (resp. reverses) the orientation of both foliations when it sends $H$ to $r_j$ (we use here the fact that our action is orientation preserving). $\phi$ defines a map from $\mathcal{H}:=\lbrace H^j_i| i\in\llbracket 1,n \rrbracket, j \in \llbracket 1,h_i \rrbracket  \rbrace$ to $\mathcal{V}:=\lbrace V^j_i| i\in\llbracket 1,n \rrbracket, j \in \llbracket 1,v_i \rrbracket  \rbrace$ and $u$ defines a function from $\mathcal{H}$ to $\{-1,1\}$. 
\begin{defi}[Geometric type of a Markovian family] \label{d.geometrictypemarkovfamily}
   We will call $$\mathcal{G}=(n,(h_i)_{i \in \llbracket 1,n \rrbracket}, (v_i)_{i\in \llbracket 1,n \rrbracket}, \mathcal{H}, \mathcal{V},\phi, u)$$ a \emph{geometric type associated to  $\mathcal{R}$} or more simply a \emph{geometric type of $\mathcal{R}$}. Similarly, a geometric type associated to a strong Markovian family of $\rho$ will be called a \emph{geometric type associated to $\rho$} or more simply a \emph{geometric type of $\rho$}.
\end{defi}

Our previous construction of $\mathcal{G}$ depends on our choice of representatives $r_1,...,r_n$, our choice of orientations of $\mathcal{P}$ and of the foliations inside every $r_i$ and finally on our choice of elements $g\in G$ in the definition of $\phi$ and $u$. Following our previous notations, 
\begin{prop}\label{p.independencegeomtype}
    The geometric type associated to $\mathcal{R}$ $$\mathcal{G}=(n,(h_i)_{i \in \llbracket 1,n \rrbracket}, (v_i)_{i\in \llbracket 1,n \rrbracket}, \mathcal{H}, \mathcal{V},\phi, u)$$ is a geometric type in the sense of Definition \ref{d.geometrictype}. Furthermore, 
    \begin{enumerate}
        \item if $k_i\in G$, then changing $r_i$ to $\rho(k_i)(r_i)$ does not change the $h_i,v_i$
        \item the elements $g\in G$ used in the definition of $\phi$ and $u$ are uniquely defined, thus $\mathcal{G}$ depends only on our choice of representatives $r_1,...,r_n$ and of orientations of $\mathcal{P}$ and of the foliations inside every $r_i$
    \end{enumerate}
\end{prop}
\begin{proof}
First, the invariance of the $h_i$ and $v_i$ when replacing $r_i$ by another rectangle in its orbit by $G$, is an immediate consequence of Remark \ref{r.precedentsuivant}.  

  Next, notice that in order to show that $\mathcal{G}$ is a geometric type it suffices to prove that $\phi$ defines a bijection (this also proves that $|\mathcal{H}|={\sum_{i=1}^{n}} h_i =|\mathcal{V}|= {\sum_{i=1}^n} v_i$). Suppose that there exist $H_p^q, H_i^k\in \mathcal{H}$ such that $H_p^q \neq H_i^k$ and $\phi(H_p^q)=\phi(H_i^k)=V_j^l\in \mathcal{V}$. Let $H$ (resp. $H'$) be the successor of $r_i$ (resp. $r_p$) associated to $H_i^k$ (resp. $H_p^q$) thanks to our choice of unstable orientation inside $r_i$ (resp. $r_p$). By our construction of $\phi$, there exist  $g,g'\in G$ such that $\rho(g)(H)=\rho(g')(H')=r_j$ and such that both $\rho(g)(r_i)$ and $\rho(g')(r_p)$ correspond to the predecessor of $r_j$ associated to $V_j^l$ for our choice of stable orientation inside $r_j$; hence $\rho(g)(r_i)=\rho(g')(r_p)$. This implies that $i=p$. Moreover, since $H_p^q \neq H_i^k$, we also have that $g\neq g'$. It follows that there exists $h\in G-\{id\}$ such that $\rho(h)(r_i)=r_i$. The previous property implies that $\rho(h)$ must have a fixed point in $r_i$ and by definition of a strong Markovian action, $\rho(h)$ must act as a topological expansion or contraction on the stable manifold of this fixed point. This contradicts the fact that $\rho(h)(r_i)=r_i$; hence $\phi$ is injective. 
  
  Let us remark here, before proving the surjectivity of $\phi$, that by using the same argument as in the previous paragraph, one obtains that the elements $g\in G$ used in the definitions of $\phi$ and $u$ are uniquely defined. 
  
  Finally, let us now prove that $\phi$ is surjective. Indeed, take $V_j^l\in \mathcal{V}$ and denote by $V$ the predecessor of $r_j$ associated to $V_j^l$ (thanks to our choice of stable orientation in $r_j$). There exists a unique $r_i\in \{r_1,...,r_n\}$ and a unique $g\in G$ such that $\rho(g)(V)=r_i$. Thanks to Remark \ref{r.precedentsuivant}, $H:=\rho(g)(r_j)$ is a successor of $r_i=\rho(g)(V)$. Denote by $H_i^k$ the element in $\mathcal{H}$ corresponding to $H$, thanks to our choice of unstable orientation in $r_i$. It is easy to see that $\phi(H_i^k)=V_j^l$, which finishes the proof of the surjectivity of $\phi$ and of the proposition.

\end{proof}

    \begin{rema}\label{r.canonicalassociationgeometrictype}
    By Definition \ref{d.geometrictypemarkovfamily} and Proposition \ref{p.independencegeomtype}, given $\mathcal{R}$ a strong Markovian family preserved by $\rho$, the following triplet of choices of representatives and orientations defines a unique geometric type of $\mathcal{R}$: 
    \begin{enumerate}
        \item a choice of representatives $r_1,...,r_n$ of every rectangle orbit in $\mathcal{R}$
        \item a choice of orientation of $\mathcal{P}$
        \item a choice of orientation of the stable and unstable foliations inside every $r_i$ such that the positive stable direction followed by the positive unstable direction inside every $r_i$ define an orientation compatible with our choice of orientation of $\mathcal{P}$
    \end{enumerate}
    \end{rema}

By the previous remark, there exist a priori more than one geometric types associated to $\mathcal{R}$. However, all the previous choices of orientations or representatives lead  to the construction of a unique geometric type up to equivalence: 

\begin{theorem}[Theorem-Definition A'] 
 Fix $\mathcal{P}$ a plane, $\rho: G\rightarrow \text{Homeo}(\mathcal{P})$ an orientation preserving strong Markovian action, preserving the pair of singular foliations $\mathcal{F}^s$ and $\mathcal{F}^u$. For any strong Markovian family $\mathcal{R}$ of $\rho$, the set of geometric types associated to $\mathcal{R}$ is included in a unique equivalence class of geometric types. 
\label{t.associatemarkovfamiliestogeometrictype}
\end{theorem}

\begin{proof}
Let $Rep:=\{r_1,...,r_n\}$ be a set of representatives of every rectangle orbit in $\mathcal{R}$ and $$\mathcal{G}=(R_1,...,R_n, (h_i)_{i \in \llbracket 1,n \rrbracket}, (v_i)_{i\in \llbracket 1,n \rrbracket}, \mathcal{H}, \mathcal{V},\phi, u)$$ be the geometric type associated to $\mathcal{R}$ for this choice of representatives, some choice of orientation of $\mathcal{P}$ and a choice of orientation of the stable and unstable  foliations inside every $r_i$ such that the positive stable direction followed by the positive unstable direction in every $r_i$ define an orientation compatible with our original choice of orientation on $\mathcal{P}$. Notice that we are using here the geometric interpretation of a geometric type, where 
\begin{itemize}
    \item $R_1,...,R_n$ are abstract copies of $[0,1]^2$ endowed with the horizontal and vertical foliations 
    \item the $(R_i)_{i\in\llbracket 1, n\rrbracket}$ are bijectively associated to the elements in $Rep$. If $R_i$ is associated to $r_j\in Rep$, then $h_i$ (resp. $v_i$)  corresponds to the number of successors (resp. predecessors) of $r_j$
    \item every $R_i$ contains $h_i$ (resp. $v_i$) horizontal (resp. vertical) subrectangles, say $H_{i}^{1},...,H_{i}^{h_i}$ (resp. $V_{i}^{1},...,V_{i}^{v_i}$),  that are pairwise  disjoint and ordered from bottom to top (resp. left to right)
    \item $\mathcal{H}=\{H_i^j|i\in \llbracket 1,n\rrbracket ~j\in \llbracket 1,h_i\rrbracket\}$ and $ \mathcal{V}=\{V_i^j|i\in \llbracket 1,n\rrbracket ~j\in \llbracket 1,v_i\rrbracket\}$
\end{itemize} 

\vspace{0.2cm}
Consider now $$\mathcal{G}'=(R'_1,...,R'_n, (h'_i)_{i \in \llbracket 1,n \rrbracket}, (v'_i)_{i\in \llbracket 1,n \rrbracket}, \mathcal{H'},\mathcal{V'},\phi', u')$$ the geometric type associated to $\mathcal{R}$ for the same choice of representatives, the same choice of orientation of $\mathcal{P}$, the same choice of orientation of the foliations inside every $r_1,...,r_{n-1}$ and finally the opposite choice of orientation for both the stable and unstable foliations inside $r_n$. We will now show that $\mathcal{G}$ and $\mathcal{G}'$ are equivalent. 

Indeed, by eventually reindexing the $(R_i)_{i\in\llbracket 1, n\rrbracket}, (R'_i)_{i\in\llbracket 1, n\rrbracket}$, we can assume that for every $i\in \llbracket 1,n \rrbracket $ the rectangles $R_i$ and $R'_i$ are associated to the representative $r_i$ (see Definition \ref{d.geometrictypemarkovfamily}) and thus that $h_i=h'_i$ and $v_i=v'_i$.

We would now like to define a bijection from $\mathcal{H}\cup \mathcal{V}$ to $ \mathcal{H'}\cup \mathcal{V'}$. There exists a natural way to associate to any rectangle in $\mathcal{H}$ a rectangle in $\mathcal{H}'$. Indeed, recall that $\mathcal{H'}:=\{{H'_i}^j|i\in \llbracket 1,n\rrbracket ~j\in \llbracket 1, h'_i\rrbracket\}$ (resp. $\mathcal{V'}:=\{{V'_i}^j|i\in \llbracket 1,n\rrbracket ~j\in \llbracket 1,v'_i\rrbracket\}$) is a set of mutually disjoint horizontal (resp. vertical) subrectangles of $R_1',...,R_n'$. Recall also that by Definition \ref{d.geometrictypemarkovfamily}, for every $i \in\llbracket 1,n \rrbracket$ and every $j \in \llbracket 1,h_i \rrbracket$, the rectangle $H_i^j$ (resp. ${H'_i}^j$) corresponds to the $j-$th successor of $r_i$ for the bottom to top order defined by the unstable orientation in $r_i$ used in the definition of $\mathcal{G}$ (resp. $\mathcal{G}'$). Denote by  $R_{i,j}\in \mathcal{R}$ (resp. $R'_{i,j}\in \mathcal{R}$) the successor of $r_i$ associated to $H_i^j$ (resp. ${H'_i}^j$). Notice that by our choices of orientations, if $i\neq n$, then  $R_{i,j}=R'_{i,j}$ and if $i=n$, then $R_{n,j}=R'_{n,(h_n-j)}$. We can therefore associate in the first case ${H'_i}^j$ to $H_i^j$ and in the second ${H'_n}^{(h_n-j)}$ to $H_n^j$. It is clear that the previous map defines a bijection $H_h:\mathcal{H}\rightarrow \mathcal{H'}$. Moreover, notice that for every $i \in\llbracket 1,n \rrbracket$ we have that $H_h(\{H_i^j, j\in \llbracket 1, h_i\rrbracket \})= \{{H'_i}^{j'}, j'\in \llbracket 1, h_i\rrbracket\}$ and that $H_h$ restricted on $\{H_i^j, j\in \llbracket 1, h_i\rrbracket \} $ is increasing with respect to $j$ when $i\neq n$ and decreasing when $i=n$. In the exact same way, we can construct a bijection $H_v$ from $\mathcal{V}$ to $\mathcal{V'}$, whose restriction on $\{V_i^j, j\in \llbracket 1, v_i\rrbracket \} $ is increasing with respect to $j$ if and only if $i\neq n$. By combining $H_v$ and $H_h$, we define a bijection $H_{eq}:\mathcal{H}\cup \mathcal{V}\rightarrow  \mathcal{H'}\cup \mathcal{V'}$. Let us now show that $H_{eq}$ is an equivalence between $\mathcal{G}$ and $\mathcal{G}'$. 
 


Indeed, following the notations of Definition \ref{d.equivalentgeomtypes}, we have that $\epsilon(R_i)=\epsilon'(R_i)=+1$ if $i\neq n$ and $\epsilon(R_n)=\epsilon'(R_n)=-1$; hence for every $i\in \llbracket 1,n \rrbracket $ we have that $\epsilon(R_i)\cdot \epsilon'(R_i)=+1$. Thanks to the previous fact, it suffices to show that for any rectangle $H_i^j\in \mathcal{H}$ whose image by $\phi$ is a vertical subrectangle of $R_k$, we have  $$H_{eq}\circ \phi(H_i^j)= \phi'\circ H_{eq}(H_i^j) \text{ and }u'(H_{eq}(H_i^j))=\epsilon(R_i)\epsilon(R_k)u(H_i^j)$$ Let us start by proving the first equality. By Definition \ref{d.equivalentgeomtypes}, the subrectangles in $\mathcal{V}'$ are bijectively identified with the predecessors of the rectangles in $Rep=\{r_k|k\in \llbracket 1,n\rrbracket\}$. Therefore, in order to show that both $H_{eq}\circ \phi(H_i^j) \in \mathcal{V}'$ and $\phi'\circ H_{eq}(H_i^j) \in \mathcal{V}'$ correspond to the same vertical subrectangle of $R'_k$, it suffices to show that they correspond to the same predecessor of $r_k$.

Recall that by Definition \ref{d.geometrictypemarkovfamily}, $H_i^j\in \mathcal{H}$ corresponds to the $j-$th successor of $r_i$, say $R$, for the bottom to top order defined  by our choice of unstable orientation used in the definition $\mathcal{G}$. Next, by Proposition \ref{p.independencegeomtype} and Definition \ref{d.geometrictypemarkovfamily}, there exists a unique $g\in G$ such that $\rho(g)(R)=r_k \in Rep$ and such that the rectangle $\phi(H_i^j)$ corresponds to the predecessor $\rho(g)(r_i)$ of $r_k$. Assume that $\phi(H_i^j)=V_k^l\in \mathcal{V}$ and thus that $\rho(g)(r_i)$ is the $l$-th predecessor of $r_k$ in the left to right order defined  by our choice of stable orientation used in the definition of $\mathcal{G}$. 

If $k\neq n$, then by our definition of $H_{eq}$, the rectangle $H_{eq}(\phi(H_i^j))={V'_k}^l$ is also the $l$-th predecessor of $r_k$ in the left to right order defined by our choice of stable orientation in $r_k$ used in the definition of $\mathcal{G}'$. Since $k\neq n$, our choices of stable orientations in $r_k$ coincide for $\mathcal{G}$ and $\mathcal{G}'$, thus $H_{eq}(\phi(H_i^j))\in \mathcal{V}'$ corresponds to the predecessor $\rho(g)(r_i)$ of $r_k$. 

If $k=n$, then the rectangle $H_{eq}(\phi(H_i^j))={V'_k}^{(v_k-l)}$ corresponds to the $(v_k-l)$-th predecessor of $r_k$ in the left to right order defined by our choice of stable orientation in $r_k$ used in the definition of $\mathcal{G}'$. Since $k=n$, our choices of stable orientations in $r_k$ are opposite for $\mathcal{G}$ and $\mathcal{G}'$, thus $H_{eq}(\phi(H_i^j))$ corresponds to the predecessor $\rho(g)(r_i)$ of $r_k$.
 
Similarly, by the exact same argument, we get that the rectangle $H_{eq}(H_i^j)$ corresponds to the successor $R$ of $r_i$. Moreover, by Proposition \ref{p.independencegeomtype}, $g$ is the unique element in $G$ for which $\rho(g)(R)=r_k \in Rep$ and thus, by Definition \ref{d.geometrictypemarkovfamily}, $\phi'(H_{eq}(H_i^j))$ corresponds to the predecessor $\rho(g)(r_i)$ of $r_k$. This proves that $H_{eq}\circ \phi(H_i^j)= \phi'\circ H_{eq}(H_i^j)$. 

Recall now that by Definition \ref{d.geometrictypemarkovfamily}, $u'(H_{eq}(H_i^j))=+1$ (resp. $u(H_i^j)=+1$) if and only if $\rho(g)$ sends positively oriented stable/unstable segments in $r_i\cap R$ -with respect to the stable and unstable orientations used in the definition of $\mathcal{G}$ (resp. $\mathcal{G}'$)- to positively oriented segments in $r_k$ -with respect to the stable and unstable orientations used in the definition of $\mathcal{G}$ (resp. $\mathcal{G}'$)- By our choices of orientations, we have that $u(H_i^j)=u'(H_i^j)$ if $i,k\neq n$ or $i=k=n$ and $u(H_i^j)=-u'(H_i^j)$ if $i\neq k=n$ or $k\neq i=n$. The previous statement is equivalent to  $u'(H_{eq}(H_i^j))=\epsilon(R_i)\epsilon(R_k)u(H_i^j)$.

We have thus shown that, starting from any geometric type $\mathcal{G}$ of $\mathcal{R}$, given by a choice of representatives $Rep$ and orientations as in Remark \ref{r.canonicalassociationgeometrictype}, changing the orientation of the stable and unstable foliations inside a unique rectangle leads to the construction of an equivalent geometric type. By applying this a finite number of times and by using the fact that the equivalence of geometric types is an equivalence relation, we get that changing the orientation of the stable and unstable foliations inside any number of rectangles in $Rep$  leads to the construction of a geometric type $\mathcal{G}''$ that is equivalent to $\mathcal{G}$. By a similar argument, we can also show that changing our choice of orientation of $\mathcal{P}$ and the stable orientation inside every representative, leads also to the construction of an equivalent geometric type. It follows that any changes in our choices of orientations only lead to geometric types in the same equivalence class. Finally, it can be easily checked that changing a representative $r_i$ by $\rho(w)(r_i)$ where $w\in G$, leads to the construction of a geometric type equal to $\mathcal{G}$ when the stable and unstable orientations inside $\rho(w)(r_i)$ are the push-forward by  $\rho(w)$ of the stable and unstable orientations inside $r_i$. This concludes the proof of the theorem. 

\end{proof}

\section{Boundary periodic points and arc points}\label{s.boundarypoints}

A strong Markovian family does not cover all the points of a plane in the same way: for every strong Markovian family there exist points that do not belong in the interior of any rectangle of the family. Among those points, we can distinguish non-periodic points that we will call \textit{boundary arc points} and periodic points that we will call  \textit{boundary periodic points}. 

The sets of boundary arc and boundary periodic points will be of great importance to us throughout this paper, since they are the only points of the plane that a Markovian family cannot ``see". This is the main reason why, as we will later prove, a class of geometric types associated to a pseudo-Anosov flow contains all the information of the original flow, except from the behavior of the flow around the boundary periodic orbits (see Theorem B). In this section, our goal is to show that for any strong Markovian family of a strong Markovian action $\rho$ the set of boundary arc points and also the set of boundary periodic points are non-empty and are the union of a finite number of orbits by $\rho$.

As in the previous section, fix $\mathcal{P}$ a plane, $\rho: G\rightarrow \text{Homeo}(\mathcal{P})$ an orientation preserving strong Markovian action, preserving the pair of singular foliations $\mathcal{F}^s$ and $\mathcal{F}^u$ and leaving invariant a strong Markovian family $\mathcal{R}$.

\begin{defi}\label{d.boundaryperiodic}
We will say that $x\in \mathcal{P}$ is a \emph{periodic point for $\rho$} if and only if there exists $g\in G-\{id\}$ such that $\rho(g)(x)=x$. We similarly define a \emph{periodic leaf for $\rho$} in  $\mathcal{F}^s$ or $\mathcal{F}^u$.

Consider $x\in\mathcal{P}$ that is not contained in the interior of any rectangle of $\mathcal{R}$. If $x$ is a periodic point of $\rho$, we will call $x$ a \emph{boundary periodic point}. If not, we will call $x$ a \emph{boundary arc point}. 

\end{defi}

\begin{prop}\label{p.boundaryperiodicpoints}
The set of boundary periodic points of $\mathcal{R}$ is non-empty, discrete and is the union of a finite number of orbits by the action of $\rho$. 
\end{prop}
\begin{proof}
By Lemma \ref{l.periodicboundary}, the stable and unstable boundary segments of any rectangle in $\mathcal{R}$ belong to periodic stable and unstable leaves. By the finiteness axiom, we have that the set $A$ of periodic stable or unstable leaves containing a stable or an unstable boundary component of some rectangle in $\mathcal{R}$ corresponds to the union of the $\rho$-orbits of a finite number of leaves. Therefore, by the definition of a Markovian action, the $\rho$-invariant set $B$ of periodic points belonging to a leaf in $A$ is also finite up to the action of $\rho$. We conclude that the set of boundary periodic points is a subset of $B$ and since it clearly forms a $\rho$-invariant set, we obtain a proof of the fact that the set of boundary periodic points of $\mathcal{R}$ is finite up to the action of $\rho$. 

In order to show that there exists at least one boundary periodic point, it suffices to show that every point of $B$ is a boundary periodic orbit. Take $p\in B$. Suppose $p$ is contained in the interior of some rectangle $R\in \mathcal{R}$. Since $p\in B$, we can assume without any loss of generality that $\mathcal{F}^s(p)$ contains a stable boundary component, say $s$, of some rectangle $R' \in \mathcal{R}$. Let $g \in \text{Stab}(p)$ such that $\rho(g)(s)\subset \inte{R}$ (such a $g$ exists thanks to the definition of a Markovian action). The intersection of $\rho(g)(R')$ and $R$ does not satisfy the Markovian intersection axiom. Therefore, $p$ is not contained in the interior of any rectangle in $\mathcal{R}$ and is a boundary periodic point. 

Finally, the fact that the set of boundary periodic points of $\mathcal{R}$ is discrete follows from Definition \ref{d.boundaryperiodic} and the finite return axiom. 
\end{proof}

The following remark provides an alternative proof of the existence of boundary periodic points when $\mathcal{F}^{s,u}$ contain a singularity: 
\begin{rema}\label{r.singularitiesareboundaries}
    Any prong singularity of $\mathcal{F}^{s,u}$ is a boundary periodic point of $\mathcal{R}$. 
\end{rema}
Let us now prove the analogue of Proposition \ref{p.boundaryperiodicpoints} for boundary arc points. 

\begin{prop}\label{p.arcpointsexist}
The set of boundary arc points of $\mathcal{R}$ is non-empty and is the union of a finite number of orbits by the action of $\rho$. 
\end{prop}

In order to prove the above proposition, we will use the following lemma that resembles closely Lemma \ref{l.existenceofpredecessors}. 
\begin{lemm}\label{l.crossingrectanglesnoperiodicpoints}
Take a rectangle $R \in \mathcal{R}$ and $s$ one of its stable boundary components. Suppose that $s$ does not contain any periodic point of $\rho$. There exists a unique family of rectangles $R_1,...,R_n \in \mathcal{R}$ such that:
\begin{enumerate}
\item  for every $ i\in \llbracket 1,n \rrbracket$ $\partial^s R_i \cap s =\emptyset$ and $R\cap R_i$ is a non-trivial vertical subrectangle of $R$
\item $R_1,...,R_n$ are maximal for the previous property. In other words, any $R' \in \mathcal{R}$ intersecting $R$ along a non-trivial vertical subrectangle and such that $\partial^s R' \cap s =\emptyset$ verifies $R' \cap R \subseteq R_i \cap R$ for some $i \in \llbracket 1, n \rrbracket$ 
\item $R_1,...,R_n$ have disjoint interiors 
\item $R_1,...,R_n$ cover $R$: $\overset{n}{\underset{i=1}{\cup}} R_i \cap R = R$ 
\end{enumerate}
\end{lemm}
\begin{proof}
The proof of this lemma is an adaptation of the proof of Lemma \ref{l.existenceofpredecessors}. Let us call the property of ``intersecting $R$ along a non-trivial vertical subrectangle and having stable boundaries disjoint from $s$" property $(\star)$. Fix $x\in s$ and $\mathcal{F}_{+}^s(x)$ a stable separatrix of $x$ intersecting  $s$ along a non-trivial segment. We will show that there exists $r \in \mathcal{R}$ satisfying  $(\star)$ and intersecting $\mathcal{F}_{+}^s(x)$ along a non-trivial segment containing $x$. 

Assume the opposite. By Lemma \ref{l.existenceofpredecessors}, $x$ belongs to exactly one predecessor of $R$ intersecting non-trivially $\mathcal{F}^s_+(x)$, say $R_1$. By our assumption, $R_1$ does not verify $(\star)$, hence $x$ belongs to the stable boundary of $R_1$. By the same argument, $x$ belongs to exactly one predecessor of $R_1$ intersecting non-trivially $\mathcal{F}^s_{+}(x)$, say $R_2$. Again, by our assumption $R_2$ does not verify $(\star)$. By induction, we can construct $R=R_0,R_1,...,R_n...$ such that $R_{n+1}$ is a predecessor of $R_n$, every $R_n$ contains $x$ in its stable boundary and intersects non-trivially $\mathcal{F}^s_+(x)$.  

By the finiteness axiom, there exists $(R_{i(n)})_{n\in \mathbb{N}}$ a subsequence of $(R_n)_{n\in \mathbb{N}}$ containing rectangles in the same orbit of $\rho$. By eventually considering another subsequence, we can assume that there exists $g_0\in G$ preserving  $\mathcal{F}^s(x)$ and such that $\rho(g_0)(R_{i(0)})=R_{i(1)}$. It follows that if $s'\subset s$ is the stable boundary of $R_{i(0)}$ containing $x$ (see Figure \ref{f.proof34}), then $\rho(g_0)(s')\subset s'\subset s$. Hence, $s$ contains a periodic point, which contradicts our original hypothesis. 
\begin{figure}[h!]
\centering
\includegraphics[scale=0.55]{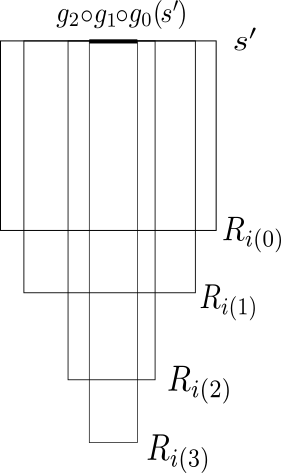}
\caption{}
\label{f.proof34}
\end{figure}

Finally, by repeating our proof of Lemma  \ref{l.existenceofpredecessors}, one can easily prove that: 
\begin{itemize}
    \item for every $x\in s$ and any stable separatrix of $x$, say $\mathcal{F}^s_+(x)$, that intersects $s$ along a non-trivial segment, there  exists a maximal rectangle for $(\star)$ containing  $x\in s$ and intersecting $\mathcal{F}^s_+(x)$ along a non-trivial segment. This proves in particular that the maximal rectangles for $(\star)$ cover $R$.

    \item two maximal rectangles for $(\star)$ are either equal or have disjoint interiors. 
    
    \item the set of maximal rectangles for $(\star)$ is finite.

\end{itemize}
   We remind the reader that in our proof of Lemma  \ref{l.existenceofpredecessors}, thanks to the Markovian intersection axiom, the non-existence of a maximal rectangle for $(\star)$ satisfying the properties of the first point implied the existence of an infinite sequence of rectangles in $\mathcal{R}$ that contradicted Lemma \ref{l.infiniteintersectionhorizontalrectangles}. Next, the second point was a consequence of the Markovian intersection axiom. Finally, our proof of the third point relied on the fact that if there existed infinitely many rectangles maximal for $(\star)$, then we could find arbitrarily thin rectangles along the stable direction that would be maximal for $(\star)$; this would contradict the first of the above points. 
\end{proof}

Of course the analogue of the previous lemma for horizontal subrectangles is also true. By a similar argument and by restricting to a fundamental domain, we can also prove the following: 
\begin{lemm}\label{l.crossingrectangleswithperiodicpoints}
Take a rectangle $R \in \mathcal{R}$ and $s$ one of its stable boundary components. Suppose that $s$ contains a periodic point $p$. There exists a family of rectangles $R_1,...,R_n \in \mathcal{R}$ such that: 
\begin{enumerate}
\item for every $i\in \llbracket 1,n \rrbracket$ $\partial^sR_i \cap s=\emptyset$ and $R_i\cap R$ is a non-trivial vertical subrectangle of $R$
\item the $R_1,...,R_n$ are ``maximal" for the previous property: any $R' \in \mathcal{R}$ intersecting $R$ along a non-trivial vertical subrectangle and whose stable boundary is disjoint from $s$ verifies $R' \cap R \subseteq \rho(g)(R_i) \cap R$ for some $i \in \llbracket 1, n \rrbracket$ and some $g\in \text{Stab}(p)$
\item For every $i\neq j$ and $g\in \text{Stab}(p)$ the rectangles $\rho(g)(R_i)$ and $R_j$ have disjoint interiors 
\item If $\mathcal{F}^s_+(p)$ is a stable separatrix of $p$ intersecting $s$ non-trivially, then the union of $R_1\cap \mathcal{F}^s_+(p),...,R_n \cap \mathcal{F}^s_+(p)$ does not contain $p$ and forms a fundamental domain of $\mathcal{F}_+^s(p)-p$ for the action of $\text{Stab}(\mathcal{F}_+^s(p))\cong \mathbb{Z}$ 
\item For any other family $R_1',...,R'_{n'}\in \mathcal{R}$ with the previous properties we have that $n'=n$ and also, up to reindexing, we have that there exist $g_1,...,g_n\in \text{Stab}(p)$ such that $R_1'=\rho(g_1)(R_1)$,...,$R_n'=\rho(g_n)(R_n)$ 
\end{enumerate}
\end{lemm}

\begin{figure}[h!]
    \centering
    \includegraphics[scale=0.8]{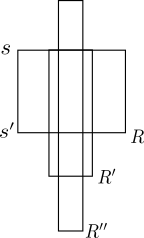}
    \caption{The rectangle $R''$ is an $s$-crossing predecessor of $R$ and the rectangle $R'$ is a crossing successor of $R''$. }
    \label{f.crossingexample}
\end{figure}

\begin{defi}\label{d.crossingsuccessor}
Take a rectangle $R \in \mathcal{R}$ and $s$ one of its stable boundary components. We will call $R'$ an \emph{$s$-crossing predecessor} of $R$ if it is maximal for the following property (in the sense of the previous lemmas): $R'$ intersects $R$ along a non-trivial vertical subrectangle and $s\cap \partial^s R' = \emptyset$. We define similarly a \emph{$u$-crossing successor} where $u$ in an unstable boundary component of $R$. 
\end{defi}
\begin{rema}\label{r.crossingpredsucc}
Consider $R\in \mathcal{R}$, $s,s'$ the two stable boundary components of $R$ and $R'$ a predecessor of $R$.
    \begin{itemize}
        \item  The rectangle $R'$ is not necessarily an $s$-crossing predecessor of $R$. However, since $R\cap R'$ is a \emph{non-trivial} horizontal subrectangle of $R'$, we have that if $R'$ is not an $s$-crossing predecessor of $R$, then it is necessarily an $s'$-crossing predecessor of $R$ (see Figure \ref{f.crossingexample}).
        \item For simplicity purposes, we will often call an $s$-crossing or $s'$-crossing predecessor of $R$, a crossing predecessor of $R$. Same for successors.
        \item Contrary to the case of predecessors/successors, if $R''$ is a crossing predecessor of $R$, then $R$ is not necessarily a crossing successor of $R''$ (see Figure \ref{f.crossingexample}). Moreover, a crossing successor/predecessor of $R$ can be a successor/predecessor of $R$ of a big generation.
    \end{itemize}
\end{rema}
\begin{proof}[Proof of Proposition \ref{p.arcpointsexist}]
Consider $R\in \mathcal{R}$ and $s$ a stable boundary component of $R$. Notice that, by definition, the set of boundary arc points in $s$ is contained in the intersection of $s$ with the unstable boundaries of the $s$-crossing predecessors of $R$. Hence, by Lemma \ref{l.crossingrectanglesnoperiodicpoints}, if $s$ contains no periodic point, the number of boundary arc points in $s$ is finite. If $s$ contains a periodic point, by Lemma \ref{l.crossingrectangleswithperiodicpoints} the number of boundary arc points in $s$ is finite up to the action of the stabilizer of the periodic point. Naturally, the previous facts remain true if $s$ were an unstable boundary component of $R$. Therefore, since the number of rectangles of any Markovian family is finite up to the action of $\rho$, the number of boundary arc points up to the action of $\rho$ is also finite.

We will now show that the set of boundary arc points is not empty. It suffices to show that for every $R\in \mathcal{R}$, if $s$ is a stable (resp. unstable) boundary component of $R$, then the intersection of $s$ with the unstable (resp. stable)  boundary of one  $s$-crossing predecessor (resp.  successor) of $R$ consists of two boundary arc points. Indeed if this is the case, by Lemmas \ref{l.crossingrectanglesnoperiodicpoints} and \ref{l.crossingrectangleswithperiodicpoints}, we get the desired result. 
 

Take a rectangle $R\in \mathcal{R} $, $s$ a stable boundary component of $R$, $R'$ a $s$-crossing predecessor of $R$ and $x \in \partial^u R' \cap s$. Let us show that $x$ is a boundary arc point. Notice first that thanks to Lemmas \ref{l.crossingrectanglesnoperiodicpoints} and \ref{l.crossingrectangleswithperiodicpoints}, $x$ cannot be a periodic point. Hence, it suffices to show that $x$ is not in the interior of some rectangle in $\mathcal{R}$. Suppose that $x$ is contained in the interior of $R'' \in \mathcal{R}$ (see Figure \ref{f.Proofprop33}). In this case, it is not difficult to see that  the Markovian intersection axiom implies that $R'\cap R''$ is a vertical subrectangle of $R''$ and $R\cap R''$ is also a vertical subrectangle of $R$. This contradicts the maximality of the $s$-crossing predecessor $R'$ and finishes the proof of our initial claim. 

\begin{figure}[h]

    \includegraphics[scale=0.4]{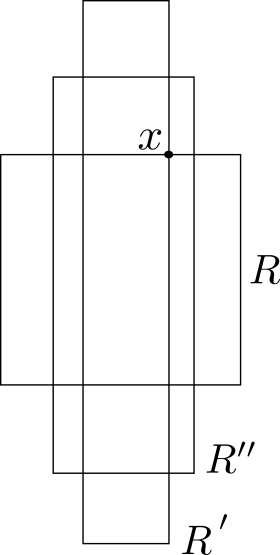}
   
  \caption{This is the only possible configuration for which $x\in \inte{R''}$}
  \label{f.Proofprop33}
\end{figure}
  
\end{proof}




\begin{rema}\label{r.caracterisationarcpoints}In the proofs of Propositions \ref{p.boundaryperiodicpoints} and \ref{p.arcpointsexist} we established that:
\begin{itemize}
\item A periodic point that is contained in the boundary of one rectangle in $\mathcal{R}$ is a boundary periodic point. 
\item More generally, a periodic point in the plane is a boundary periodic point if and only if its stable or unstable leaf contain a stable or unstable boundary component of some rectangle in $\mathcal{R}$. 
\item Consider $R\in \mathcal{R}$ and $s$ a stable (resp. unstable) boundary component of $R$. A point in the intersection of $s$ with the unstable (resp. stable)  boundary of one $s$-crossing predecessor (resp.  successor) of $R$ is a boundary arc point. Conversely, a predecessor (resp. successor) of $R$ that crosses $s$ and whose unstable (resp. stable) boundary intersects $s$ along two boundary arc points is an $s$-crossing predecessor of $R$. 
\end{itemize}
\end{rema}

\section{The bifoliated plane up to surgeries}\label{s.constructionPbar}
Let $\mathcal{P}_1$,  $\mathcal{P}_2$ be two planes endowed with two orientation preserving strong Markovian actions $\rho_1: G_1\rightarrow \text{Homeo}{(\mathcal{P}_1)}$, $\rho_2: G_2\rightarrow \text{Homeo}{(\mathcal{P}_2)}$ preserving respectively the singular foliations $\mathcal{F}^{s,u}_1$, $\mathcal{F}^{s,u}_2$ and the strong Markovian families $\mathcal{R}_1$, $\mathcal{R}_2$. Consider $\Gamma_1$,  $\Gamma_2$ the boundary periodic orbits of $\mathcal{R}_1$ and $\mathcal{R}_2$ (those sets are non-empty thanks to Lemma \ref{l.periodicboundary}). Assume now that the equivalence classes of geometric types associated to $\mathcal{R}_1$ and $\mathcal{R}_2$ are the same (see Definition \ref{d.geometrictypemarkovfamily}). 

In view of Theorem B, we would like to prove that, when the two previous actions correspond to the actions on the bifoliated planes of two pseudo-Anosov flows $\Phi_1,\Phi_2$, then $
\Phi_1$ can be obtained one from $\Phi_2$ by a finite number of Dehn-Goodman-Fried surgeries on specific periodic orbits. In this section, we will introduce the \emph{bifoliated plane up to surgeries}, the main tool that will allow us to compare pseudo-Anosov flows up to specific surgeries. Once again, in order to describe our results in their most general form, we will define the bifoliated plane up to surgeries for any strong Markovian action. 

More specifically, even though $\mathcal{R}_1$ and $\mathcal{R}_2$ are associated to the same equivalence class of geometric types, the actions $\rho_1$ and $\rho_2$ can significantly differ inside $\mathcal{P}_1$ and $\mathcal{P}_2$. However, as we will later prove, this is no longer the case for the lifts of the previous actions on the universal covers of $\mathcal{P}_1-\Gamma_1$ and $\mathcal{P}_2-\Gamma_2$, denoted by $\widetilde{\mathcal{P}_1}$ and $\widetilde{\mathcal{P}_2}$. Intuitively, the reason why this is happening is because we erased from the bifoliated plane all the points that the Markovian families $\mathcal{R}_1,\mathcal{R}_2$ can not ``see". In the case where  $\rho_1,\rho_2$ are associated to pseudo-Anosov flows, the previous result implies that the lifts of $\rho_1$ and $\rho_2$ on $\widetilde{\mathcal{P}_1}$ and $\widetilde{\mathcal{P}_2}$ are the same (up to conjugation). This motivates our interest in the spaces $\widetilde{\mathcal{P}_1}$ and $\widetilde{\mathcal{P}_2}$. 

From the practical standpoint the study of $\widetilde{\mathcal{P}_1}$ and $\widetilde{\mathcal{P}_2}$ poses a serious problem: the lifts of the rectangles in $\mathcal{R}_1$, $\mathcal{R}_2$ on  $\widetilde{\mathcal{P}_{1}}$, $\widetilde{\mathcal{P}_{2}}$ do not always correspond to  rectangles in the sense of Definition \ref{d.standardpolygon}; some rectangles lift to ``rectangles" with corners at infinity and thus cease being compact. In order to avoid this difficulty, by adding some points at infinity, we will complete  $\widetilde{\mathcal{P}_{1}}$, $\widetilde{\mathcal{P}_{2}}$ to ``branched" covering spaces of $\mathcal{P}_{1}$, $\mathcal{P}_{2}$ that we will denote by $\clos{\mathcal{P}_{1}}$, $\clos{\mathcal{P}_{2}}$. The ``ramification points" of these covers will correspond to $\Gamma_{1}$, $\Gamma_{2}$ and their ``ramification indexes" will be infinite. The spaces $\clos{\mathcal{P}_{1}}$, $\clos{\mathcal{P}_{2}}$ will be called the \emph{bifoliated planes of $\rho_1,\rho_2$ up to surgeries along $\Gamma_1,\Gamma_2$}.

\subsection{The bifoliated plane up to surgeries: the construction} \label{s.constructionPbar12} 
Fix $\mathcal{P}$ a plane endowed with an orientation, $\rho: G\rightarrow \text{Homeo}(\mathcal{P})$ an orientation preserving strong Markovian action, preserving the pair of singular foliations $\mathcal{F}^s$ and $\mathcal{F}^u$ (see Definition \ref{d.folisingular}) and leaving invariant a strong Markovian family $\mathcal{R}$. Denote by $\Gamma$ the boundary periodic points of $\mathcal{R}$.  

The construction of  $\clos{\mathcal{P}}$ is very similar to the construction of the universal cover $\widetilde{\mathcal{P}}$ of $\mathcal{P}-\Gamma$. Let us first recall the latter one. 
Fix $x_0\in \mathcal{P}-\Gamma$ and $\text{Curv}_1= \lbrace \gamma:[0,1]\overset{C^0}{\rightarrow} \mathcal{P} | ~ \gamma(0)=x_0, \gamma[0,1]\cap \Gamma= \emptyset \rbrace$.  For any  $\gamma_1,\gamma_2\in \text{Curv}_1$, we will say that $\gamma_1$ and $\gamma_2$ are \emph{homotopic for $\sim_1$} or $\gamma_1 \sim_1 \gamma_2$ if there exists $H:[0,1]^2\overset{C^0}{\rightarrow} \mathcal{P}$ such that $H(0,\cdot)=\gamma_1$, $H(1,\cdot)=\gamma_2$, $H(\cdot,1)$ is constant and for every $t\in [0,1]$ we have $H(t,\cdot)\in  \text{Curv}_1$. By a classical result of covering space theory, we have that $$ \widetilde{\mathcal{P}} = \text{Curv}_1 \big{/}\sim_1 $$

Any point $\tilde{x}\in \widetilde{\mathcal{P}}$ corresponds to a class of arcs of $\text{Curv}_1$ starting from $x_0$ and ending at the same point of $\mathcal{P}-\Gamma$, say $x$. We define the projection $\tilde{\pi}$ from $\widetilde{\mathcal{P}}$ on $\mathcal{P}-\Gamma$ as the function associating to $\tilde{x}$ the point $x$. The previous function is clearly a surjection. Conversely,  for any $\delta \in \text{Curv}_1$, we will denote by $\widetilde{\delta}$ the equivalence class for $\sim_1$ containing the curve $\delta$.

Recall now that $\widetilde{\mathcal{P}}$ is naturally endowed with a topology. Consider a point $\tilde{x}\in \widetilde{\mathcal{P}}$, $\gamma: [0,1]\rightarrow \mathcal{P}$ an arc inside the equivalence class defined by $\tilde{x}$, $V$ an open disk in $\mathcal{P}-\Gamma$ containing $x$ and $s\in (0,1)$ such that $\gamma([s,1])\subset V$. We define the set $\widetilde{V}:=\lbrace\widetilde{\delta}| \delta\in \text{Curv}_1, \forall t\in [0,s] \quad \delta(t)=\gamma(t) , \delta([s,1])\subset V \rbrace$ to be an open neighborhood of $\tilde{x}$ in $\widetilde{\mathcal{P}}$. It is easy to see that $\widetilde{V}$ does not depend on our choice of $\gamma$ or $s$. We endow $\widetilde{\mathcal{P}}$ with the topology generated by the previous open sets. It is known that for our previous choice of topology the map $\tilde{\pi}$ is continuous and that for every $\widetilde{x}\in \widetilde{\mathcal{P}}$ the set $$\{\widetilde{W}| W \text{ open disk in $\mathcal{P}-\Gamma$ containing $\tilde{\pi}(x)$} \}$$ forms a neighborhood basis of $\tilde{x}$ in $\widetilde{\mathcal{P}}$. 

Now let us define the space $\clos{\mathcal{P}}$. Take $\text{Curv}_2= \lbrace \gamma:[0,1]\overset{C^0}{\rightarrow} \mathcal{P} | ~ \gamma(0)=x_0, \gamma[0,1)\cap \Gamma= \emptyset \rbrace $. For any $\gamma_1,\gamma_2\in \text{Curv}_2$, we will say that $\gamma_1$ and $\gamma_2$ are \emph{homotopic for $\sim_2$} or $\gamma_1 \sim_2 \gamma_2$  if there exists $H:[0,1]^2\overset{C^0}{\rightarrow} \mathcal{P}$ such that $H(0,\cdot)=\gamma_1$, $H(1,\cdot)=\gamma_2$, $H(\cdot,1)$ is constant and for every $t\in [0,1]$ we have $H(t,\cdot)\in  \text{Curv}_2$. We define
$$ \clos{\mathcal{P}} := \text{Curv}_2 \big{/}\sim_2 $$

We are going to call $\clos{\mathcal{P}}$ the \emph{bifoliated plane of $\rho$ up to surgeries on $\Gamma$}. Any point $\clos{x}\in \clos{\mathcal{P}}$ corresponds to an equivalence class of arcs of $\text{Curv}_2$ starting from $x_0$ and ending at the same point of $\mathcal{P}$, say $x$. We define the projection $\clos{\pi}$ from $\clos{\mathcal{P}}$ on $\mathcal{P}$ as the function associating to $\clos{x}$ the point $x$. The function $\clos{\pi}$ is clearly a surjection. 

Let us remark at this moment that by definition, $\widetilde{\mathcal{P}}$ is naturally included in $\clos{\mathcal{P}}$. Indeed, following our previous notations, notice that $\text{Curv}_1\subset \text{Curv}_2$. More precisely, $\text{Curv}_2- \text{Curv}_1$ consists exactly of the curves in $\text{Curv}_2$ ending at a point of $\Gamma$. Moreover, if $\gamma_1,\gamma_2 \in \text{Curv}_1\subset \text{Curv}_2$, then $\gamma_1\sim_1\gamma_2$ if and only if $\gamma_1\sim_2\gamma_2$. This implies the existence of an inclusion map $\phi: \widetilde{\mathcal{P}} \rightarrow \clos{\mathcal{P}}$ such that $\widetilde{\pi}=\clos{\pi}\circ \phi$ and $\clos{\mathcal{P}}-\phi(\widetilde{\mathcal{P}})=\clos{\Gamma}$, where $\overline{\Gamma}:=\clos{\pi}^{-1}(\Gamma)$. One can thus think of the points in $\overline{\Gamma}$ as points in the boundary at infinity of $\widetilde{\mathcal{P}}$. 

Let us now define a topology on $\clos{\mathcal{P}}$. In order to do so, we will provide an alternative definition of $\clos{\mathcal{P}}$ using a smooth structure on $\mathcal{P}$.

\vspace{0.3cm}
\subsubsection{A geometric construction of $\clos{\mathcal{P}}$ and its topology} \label{s.geometricconstr}
~\vspace{0.3cm}

 Endow $\mathcal{P}$ with a smooth manifold structure and blow-up every point in $\Gamma$ to a circle. By doing so, since $\Gamma$ is discrete (see Proposition \ref{p.boundaryperiodicpoints}), we obtain a surface $\mathcal{P}_{exp}$\footnote{exp stands for explosion} homeomorphic to a plane minus a discrete infinite set of open disks. Notice that there exists a smooth projection map $p_{exp}:\mathcal{P}_{exp}\rightarrow \mathcal{P}$, whose restriction in the interior of  $\mathcal{P}_{exp}$ is a diffeomorphism. 

Denote by $X_0\in \mathcal{P}_{exp}$ the unique preimage by $p_{exp}$ of $x_0$ (the base point of $\mathcal{P}$ used in the definition of $\text{Curv}_1$, $\text{Curv}_2$) and by $\text{Curv}_3$ the set of continuous arcs in $\mathcal{P}_{exp}$ starting from $X_0$. By classical covering space theory, the space $\widetilde{\mathcal{P}_{exp}}$ obtained when we identify any two arcs in $\text{Curv}_3$ up to homotopy relatively to their boundaries is the universal cover of $\mathcal{P}_{exp}$. Notice that $\widetilde{\mathcal{P}_{exp}}$ is a surface with infinite many line boundaries and that every point in $\widetilde{\mathcal{P}_{exp}}$, or otherwise said every class of arcs in  $\text{Curv}_3$ up to homotopy relatively to their boundaries, is projected by $p_{exp}$ to unique class of homotopic arcs for $\sim_2$ in $\text{Curv}_2$ or otherwise said a unique element in $\clos{\mathcal{P}}$. 

The previous map from $\widetilde{\mathcal{P}_{exp}}$ to $\clos{\mathcal{P}}$ is surjective, but not injective. Indeed, two different curves in $\text{Curv}_3$ that are very close and that end at different points of the same boundary circle in $\mathcal{P}_{exp}$ always fail to be homotopic relatively to their boundaries in $\mathcal{P}_{exp}$; however their projections in $\mathcal{P}$ by $p_{exp}$ are close and end at the same point in $\Gamma$, thus they are  homotopic for $\sim_2$. 

Let us now define a homotopy $\sim_3$ between arcs in $\text{Curv}_3$ in order to repair this injectivity default. Consider $\gamma_1,\gamma_2 \in \text{Curv}_3$. We will say that $\gamma_1$ and $\gamma_2$ are \emph{homotopic for $\sim_3$} or $\gamma_1\sim_3\gamma_2$ if one of the two following properties is satisfied: 
\begin{itemize}
    \item $\gamma_1(1)$ belongs in the interior of $\mathcal{P}_{exp}$ and  
    there exists $H:[0,1]^2\overset{C^0}{\rightarrow} \mathcal{P}_{exp}$ such that $$H(0,\cdot)=\gamma_1, H(1,\cdot)=\gamma_2, H(\cdot,0)=X_0, \text{ and } H(\cdot,1)=\gamma_1(1)$$
    \item $\gamma_1(1)$ belongs in a circle boundary component $C$ of $\mathcal{P}_{exp}$ and there exists $H:[0,1]^2\overset{C^0}{\rightarrow} \mathcal{P}$ such that $$H(0,\cdot)=\gamma_1, H(1,\cdot)=\gamma_2, H(\cdot,0)=X_0, \text{ and } H(\cdot,1)\in C$$
\end{itemize}

One can easily check that : 
\begin{prop}\label{p.sim2sim3}
    Let $\gamma_1,\gamma_2\in \text{Curv}_3$. We have that $\gamma_1\sim_3\gamma_2$ if and only if $p_{exp}(\gamma_1)\sim_2 p_{exp}(\gamma_2)$
\end{prop}
Thanks to the above proposition we get that the map $p_{exp}$ is well defined on the quotient by $\sim_3$; hence there exists a bijection 
\begin{equation}\label{eq.identification}
    \phi_{exp}: \text{Curv}_3 \big{/}\sim_3 \longrightarrow \clos{\mathcal{P}}=\text{Curv}_2 \big{/}\sim_2
\end{equation} Moreover, it is easy to see that the space $\text{Curv}_3$ quotiented by $\sim_3$ is the surface obtained by collapsing all boundary lines of $\widetilde{\mathcal{P}_{exp}}$ to points. Denote by $proj$ the projection map from $\widetilde{\mathcal{P}_{exp}}$ to $\text{Curv}_3 \big{/}\sim_3$. We endow $\text{Curv}_3 \big{/}\sim_3$ with the quotient topology for $proj$ and $\clos{\mathcal{P}}$ with the topology rendering $ \phi_{exp}$ a homeomorphism. 

The previous construction of $\clos{\mathcal{P}}$ and its topology relies on the usage of a smooth manifold structure on $\mathcal{P}$. Keeping in mind that Markovian actions are only continuous, the previous definition of the topology on $\clos{\mathcal{P}}$ is not always easy to utilize in the $C^0$ setting. In the following lines, we would like to provide an alternative and also useful construction of the previous topology on $\clos{\mathcal{P}}$, which resembles closely the construction of the topology of $\widetilde{\mathcal{P}}$ that we described in Section \ref{s.constructionPbar12}.

\vspace{0.3cm}
\subsubsection{A second construction of the topology of $\clos{\mathcal{P}}$} \label{s.secondconstruction}
~\vspace{0.3cm}

\begin{figure}[h!]
    \centering
    \includegraphics[scale=0.35]{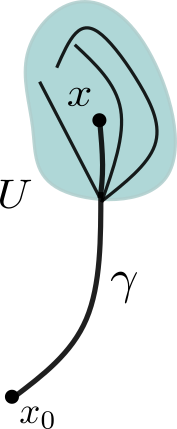}
    \caption{A neighborhood of a point in $\clos{\mathcal{P}}-\clos{\Gamma}$}
    \label{f.neibhregularpointclosP}
\end{figure}

Let  $\clos{x}\in \clos{\mathcal{P}} - \clos{\Gamma}$, $x=\clos{\pi}(\clos{x})$ and $U$ an open disk in $\mathcal{P}-\Gamma$ containing $x$. As in the case of $\widetilde{\mathcal{P}}$, for any $\delta \in \text{Curv}_2$, we will denote by $\clos{\delta}$ the equivalence class for $\sim_2$ containing the curve $\delta$. Consider $\gamma\in \text{Curv}_2$ a curve contained in the equivalence class defined by $\clos{x}$ and $s\in [0,1]$ such that $\gamma([s,1])\subset U$ (see Figure \ref{f.neibhregularpointclosP}). We define the set $$\clos{V}(U)=\lbrace\clos{\delta}| \delta\in \text{Curv}_2, \forall t\in [0,s] \quad \delta(t)=\gamma(t) , \delta([s,1])\subset U \rbrace$$ to be a neighbourhood of $\clos{x}$ in $\clos{\mathcal{P}}$. It is easy to see that $\clos{V}(U)$ depends only on our choice of $U$ and not on our choice of $\gamma$ or $s$. 

Special care is needed when defining the neighborhoods of the points in $\clos{\Gamma}$. Fix $\clos{y}\in \clos{\Gamma}$, $y=\clos{\pi}(\clos{y})$, $S$ any standard polygon containing $y$ in its interior, $(U_n)_{n\in \mathbb{Z}}$ a sequence of (embedded) closed disks in $S$ containing $y$ such that 
\begin{enumerate}
    \item for every $n\in \mathbb{Z}$ there exists $Q_n$ a quadrant of $y$ in $S$ such that $U_n\subseteq Q_n$ and $U_n$ is a neighborhood of $y$ in $Q_n$ 
    \item for every $n\in \mathbb{Z}$ we have that $Q_{n+1}$ is the quadrant right after $Q_n$ in the counter-clockwise order given by our orientation of $\mathcal{P}$
\end{enumerate} Any sequence of disks $(U_n)_{n\in \mathbb{Z}}$ in $\mathcal{P}$ satisfying the previous conditions will be called \emph{a sequence of quadrant neighborhoods of $y\in \Gamma$}. 

Fix $(U_n)_{n\in \mathbb{Z}}$ of $\mathcal{P}$ a sequence of quadrant neighborhoods of $y$. Consider $\gamma_0$ a curve in $\text{Curv}_2$ going from $x_0$ (the base point in $\text{Curv}_2$) to a point in the interior of $U_0$ such that for any embedded arc $e$ in $U_0$ going from the endpoint of $\gamma_0$ to $y$, the juxtaposition of $\gamma_0$ followed by $e$ gives a curve in $\text{Curv}_2$ inside the equivalence class defined by $\clos{y}$. Any curve $\gamma_0$ in  $\text{Curv}_2$ satisfying the previous conditions will be called \emph{a reference curve} for the sequence of quadrant neighborhoods $(U_n)_{n\in \mathbb{Z}}$ and the point $\clos{y}$.
\begin{figure}[h!]
    \centering
    \includegraphics[scale=0.35]{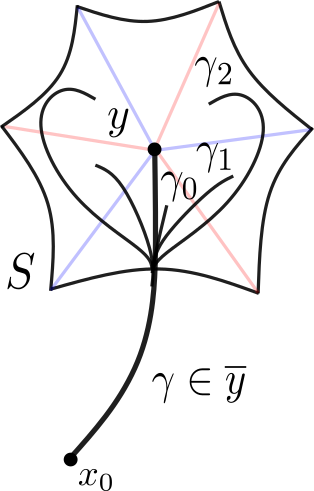}
    \caption{A neighborhood of a point in $\clos{\Gamma}$}
    \label{f.neighbconicpointclosP}
\end{figure}  

We define by induction the curves $\gamma_n$ for any $n\in \mathbb{Z}-\{0\}$ as follows (see Figure \ref{f.neighbconicpointclosP}): 

\begin{itemize}
    \item if $n> 0$ then $\gamma_n: [0,1]\rightarrow \mathcal{P}$ is a curve obtained by the juxtaposition of $\gamma_{n-1}$ followed by an arc contained in the interior of $U_{n-1}\cup U_n$ and going from the endpoint of $\gamma_{n-1}$ to any point in the interior of $U_n$ 
    \item if $n<0$ then $\gamma_n: [0,1]\rightarrow \mathcal{P}$  is a curve obtained by the juxtaposition of $\gamma_{n+1}$ followed by an arc contained in the interior of $U_{n+1}\cup U_n$ and going from the endpoint of $\gamma_{n+1}$ to a point in the interior of $U_n$ 
\end{itemize}

Let $s_n\in (0,1)$ be such that $\gamma_n([s_n,1])\subset U_n$. We define the following set to be a neighborhood of $\clos{y}$ : 

$$\clos{V}_{\clos{y}}(\gamma_0,(U_n)_{n\in \mathbb{Z}}) = \underset{n\in \mathbb{Z}}{\bigcup}\{\clos{\delta}|\delta\in \text{Curv}_2, \forall t\in [0,s_n] \quad \delta(t)=\gamma_n(t) , \delta([s_n,1])\subset U_n \}$$ 

\begin{figure}[h]

  \begin{minipage}[ht]{0.4\textwidth}
    \centering 
     \vspace{0.8cm}
    \includegraphics[width=1\textwidth]{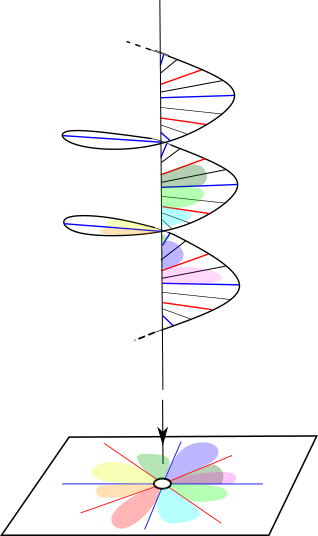}
  \hspace{-1cm}

  \end{minipage}
 \begin{minipage}[ht]{0.4\textwidth}
 \centering
    \includegraphics[width=0.9\textwidth]{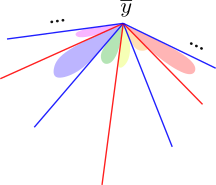}
    \hspace{-1cm}

  \end{minipage}
  
  \caption{In the figure in the left, we see a representation of a sequence of quadrant neighborhoods of $y$ in $\mathcal{P}_{exp}$ and $\widetilde{\mathcal{P}_{exp}}$. In the figure in the right, we collapse the line boundary of the left picture to a point in order to obtain a representation of $\clos{V}_{\clos{y}}(\gamma_0,(U_n)_{n\in \mathbb{Z}})$ as a neighborhood of $\clos{y}$.}
  \label{f.conicpointneighborhood}
\end{figure}

 A representation of $\clos{V}_{\clos{y}}(\gamma_0,(U_n)_{n\in \mathbb{Z}})$ using our geometric construction of $\clos{\mathcal{P}}$ can be found in Figure \ref{f.conicpointneighborhood}. One can check that $\clos{V}_{\clos{y}}(\gamma_0,(U_n)_{n\in \mathbb{Z}})$ only depends on our choice of $S$, $(U_n)_{n\in \mathbb{Z}}$ and of $\gamma_0$ and does not depend on our choice of $s_n$ (for $n\in \mathbb{Z}$) or $\gamma_n$ (for $n\neq 0$).

We endow $\clos{\mathcal{P}}$ with the topology generated by all the previous neighborhoods. 

 The reason why we defined the neighborhoods of the points in $\clos{\Gamma}$ so much differently than the neighborhoods of the points in $\clos{\mathcal{P}}-\clos{\Gamma}$ can be explained by using our geometric construction of $\clos{\mathcal{P}}$: a point $\clos{y}\in \overline{\Gamma}$ corresponds to a boundary line $\widetilde{C}$ in $\widetilde{\mathcal{P}_{exp}}$ that was collapsed to a point. If $\widetilde{C}$ projects to the boundary circle $C$ in $\mathcal{P}_{exp}$, then the lifts on $\widetilde{\mathcal{P}_{exp}}$ of the neighborhoods of $C$ inside $\mathcal{P}_{exp}$ do not form a basis of neighborhoods of $\widetilde{C}$ inside  $\widetilde{\mathcal{P}_{exp}}$. In other words, if we would like that the topologies defined in Sections \ref{s.geometricconstr} and \ref{s.secondconstruction} coincide, then contrary to the case of the points inside $\clos{\mathcal{P}} - \clos{\Gamma}$, one can not construct a basis of neighborhoods of $\clos{y}\in \overline{\Gamma}$ by simply lifting on $\clos{\mathcal{P}}$ different neighborhoods of its projection $\clos{\pi}(\clos{y}) \in \mathcal{P}$. 

\subsection{$\clos{\mathcal{P}}$ as a topological space} Denote by $\mathcal{T}$ the topology on $\clos{\mathcal{P}}$ generated by the neighborhoods defined in Section \ref{s.secondconstruction} and by $\mathcal{T}_{geom}$ the topology on $\clos{\mathcal{P}}$ defined in Section \ref{s.geometricconstr} thanks to our geometric construction.  Our main goals in this section consist in constructing a base for $\mathcal{T}$ (see Proposition \ref{p.baseoftopology}) and in proving that the topologies $\mathcal{T}$ and $\mathcal{T}_{geom}$ are equivalent (see Theorem \ref{t.topologicalpropertiesbifplaneuptosurgery}). In order to do so, we will first need to prove two useful propositions.

\begin{prop}\label{p.neighborhoodconicpoints} Consider $\clos{y}\in \clos{\Gamma}, \clos{\pi}{(\clos{y})}=y$, $S$ a standard polygon containing $y$ in its interior, $x\in S-y$, $(U_n)_{n\in \mathbb{Z}}$ a sequence of quadrant neighborhoods of $y$ in $S$ and $\gamma_0$ a reference curve for $(U_n)_{n\in \mathbb{Z}}$ and $\clos{y}$. Denote by  $\gamma,\gamma'$ any two curves in $\text{Curv}_2$ going from $x_0$ (the base point in $\text{Curv}_2$) to $x$ and such that for any embedded arc $e\subset S$ going from $x$ to $y$, the juxtaposition of $\gamma$ followed by $e$ and the juxtaposition of $\gamma'$ followed by $e$ define two curves in $\text{Curv}_2$ that belong in the equivalence class defined by $\clos{y}$. 

Denote by $\delta$ a curve in $S-y$ starting and ending at $x$ and going once around $y$ in the counter-clockwise direction. We have that there exists $k\in \mathbb{N}$ such that one of the two following situations arises: 
\begin{itemize}
    \item the juxtaposition of the curve $\gamma$ followed $k$ times by the curve $\delta$ defines a curve homotopic to $\gamma'$ for $\sim_2$ 
    \item  the juxtaposition of the curve $\gamma'$ followed $k$ times by the curve $\delta$ defines a curve homotopic to $\gamma$ for $\sim_2$ 
\end{itemize}.
\end{prop}
\begin{prop}\label{p.changingreferencecurve}
      Consider $\clos{y}\in \clos{\Gamma}, \clos{\pi}{(\clos{y})}=y$, $S$ a standard polygon containing $y$ in its interior, $2p$ be the total number of quadrants of $y$ in $S$, $(U_n)_{n\in \mathbb{Z}}$ a sequence of quadrant neighborhoods of $y$ in $S$ and $\gamma_0$ a reference curve for $(U_n)_{n\in \mathbb{Z}}$ and $\clos{y}$. Fix $k\in \mathbb{N}$ (resp. $k\in -\mathbb{N}$) and $l\in \llbracket 0, 2p-1 \rrbracket$ (resp. $l\in -\llbracket 0, 2p-1 \rrbracket$). Denote by $\gamma_0'$ a curve in $\mathcal{P}$ obtained by the juxtaposition of $\gamma_0$ followed by a curve in $S-y$ turning counter-clockwise (resp. clockwise) around $y$ and stopping at a point in the interior of $U_{2pk+l}$ after making $k$ (resp. $-k$) full turns around $y$. We have that 
         $$\clos{V}_{\clos{y}}(\gamma_0',(U_{n+2pk+l})_{n\in \mathbb{Z}})= \clos{V}_{\clos{y}}(\gamma_0,(U_n)_{n\in \mathbb{Z}})$$
  \end{prop}   
Proposition \ref{p.changingreferencecurve} is a relatively easy consequence of the construction of  $\clos{V}_{\clos{y}}(\gamma_0,(U_n)_{n\in \mathbb{Z}})$; we thus leave its proof to the reader. Let us now prove Proposition \ref{p.neighborhoodconicpoints}. 
\begin{proof}[Proof of Proposition \ref{p.neighborhoodconicpoints}]
    In order to prove the desired result, we will make use of our geometric construction of $\clos{\mathcal{P}}$ (see Section \ref{s.geometricconstr}). Endow $\mathcal{P}$ with a smooth structure and blow-up every point in $\Gamma$ in order to obtain $\mathcal{P}_{exp}$, a surface homeomorphic to a plane minus a countable discrete set of open disks. Denote by $p_{exp}$ the projection map from $\mathcal{P}_{exp}$ to $\mathcal{P}$. Recall that the restriction of $p_{exp}$ in the interior of $\mathcal{P}_{exp}$ is a diffeomorphism. 
    
   Take $e:[0,1]\rightarrow S$ an arc with no self-intersections going from $x$ to $y$ and denote by $\gamma_{exp}, \gamma_{exp}':[0,1]\rightarrow \mathcal{P}_{exp}$ the curves $p_{exp}^{-1}\circ \gamma$ and $p_{exp}^{-1}\circ \gamma'$. By eventually performing a homotopy in $\sim_2$, we  can assume that $p_{exp}^{-1}\circ e_{|[0,1)}$ extends continuously to a curve $e_{exp}:[0,1]\rightarrow \mathcal{P}_{exp}$ ending at some point in the boundary component $C$ of $\mathcal{P}_{exp}$. Let $(\gamma e) _{exp}$  (resp. $\gamma e$,  $(\gamma' e) _{exp}$, $\gamma'e$) be the curve obtained by the juxtaposition of $\gamma_{exp}$ (resp. $\gamma $,  $\gamma'_{exp}$, $\gamma'$) followed by $e_{exp}$ (resp. $e$, $e_{exp}$, $e$). The curves $(\gamma e) _{exp}$ and $(\gamma' e) _{exp}$ belong in $\text{Curv}_3$, begin at $X_0=p_{exp}^{-1}(x_0)$ and end at the same point in $C$. By Proposition \ref{p.sim2sim3}, since $\gamma e\sim_2 \gamma' e$, the curves $(\gamma e) _{exp}$ and $(\gamma' e) _{exp}$ are homotopic for $\sim_3$. Therefore, $(\gamma e) _{exp}$ can be obtained from $(\gamma' e) _{exp}$ by performing a homotopy forcing $(\gamma' e) _{exp}$ to go  around $C$ a certain number of times. In other words, if $x$ is the common endpoint of $\gamma$ and $\gamma'$, $X=p_{exp}^{-1}(x)$ and $\delta_{exp}= p_{exp}^{-1}\circ \delta$, then up to exchanging the roles of $\gamma_{exp}$ and $\gamma'_{exp}$, there exists $k\in \mathbb{N}$ such that the juxtaposition of $\gamma_{exp}$ followed $k$ times by $\delta_{exp}$ and finally by $e_{exp}$ is homotopic relatively to its boundary to $(\gamma'e) _{exp}$. This gives us the desired result. 
\end{proof}

We are now ready to begin proving the two main results of this section. 

\begin{prop}\label{p.baseoftopology}
\quad 

\begin{enumerate}
     \item Fix $\clos{x}\in \clos{\mathcal{P}}-\clos{\Gamma}$ and $x=\clos{\pi}(x)$. For any open disk $W$ inside $\mathcal{P}-\Gamma$ containing $x$, denote by $\clos{V}_{\clos{x}}(W)$ the neighborhood of $\clos{x}$ associated to $W$ that was constructed in Section \ref{s.secondconstruction}. The following set forms a neighborhood basis of $\clos{x}$   $$\{\clos{V}_{\clos{x}}(W)| W \text{ open disk in $\mathcal{P}-\Gamma$ containing $x$}\}$$  
    \item Fix $\clos{y}\in \clos{\Gamma}$, $y= \clos{\pi}(y)$ and $S$ a standard polygon containing $y$. For any sequence of quadrant neighborhoods $(U_n)_{n\in \mathbb{Z}}$ of $y$ in $S$ and any reference curve $\gamma_0$ denote by $\clos{V}_{\clos{y}}(\gamma_0,(U_n)_{n\in \mathbb{Z}})$ the neighborhood of $\clos{y}$ associated to $\gamma_0$ and $(U_n)_{n\in \mathbb{Z}}$ that was constructed in Section \ref{s.secondconstruction}. The following set forms a neighborhood basis of $\clos{y}$
    \begin{align*}
        \Bigl\{ \clos{V}_{\clos{y}}(\gamma_0, (U_n)_{n\in \mathbb{Z}})|  &(U_n)_{n\in \mathbb{Z}} \text{ is a sequence of quadrant neighborhoods of $y$ in $S$ } \\ &\gamma_0 \text{ is a reference curve for $(U_n)_{n\in \mathbb{Z}}$ and $\clos{y}$} \Bigl\}
    \end{align*} 
\end{enumerate}
\end{prop}
\begin{proof}
    Take $\clos{x}\in \clos{\mathcal{P}}-\clos{\Gamma}$ and $\mathcal{V}_{\clos{x}}$ the set of subsets $\mathcal{N}$ of $\clos{\mathcal{P}}$ for each one of which there exists $W$ an open disk in $\mathcal{P}-\Gamma$ containing $x$  such that $\clos{V}_{\clos{x}}(W)\subset \mathcal{N}$. The set $\mathcal{V}_{\clos{x}}$ satisfies the following properties:
    \begin{enumerate}
    \item $\clos{\mathcal{P}}\in \mathcal{V}_{\clos{x}}$
    \item every element of $\mathcal{V}_{\clos{x}}$ contains $\clos{x}$ 
    \item if $\mathcal{N}\in \mathcal{V}_{\clos{x}}$ and $\mathcal{N}\subset \mathcal{N}'$ then $\mathcal{N}'\in \mathcal{V}_{\clos{x}}$
    \item for every element $\mathcal{N} $ of $\mathcal{V}_{\clos{x}}$ there exists $\mathcal{N}'\in \mathcal{V}_{\clos{x}}$ a subset of $\mathcal{N}$ such that for all $\clos{u}\in \mathcal{N}'$ we have that $\mathcal{N}\in\mathcal{V}_{\clos{u}}$ 
     \item $\mathcal{V}_{\clos{x}}$ is stable by finite intersections
    \end{enumerate}
    Indeed, the proofs of (1)-(3) are direct applications of the definition of $\mathcal{V}_{\clos{x}}$. Let us now prove (4). Take $\mathcal{N}\in \mathcal{V}_{\clos{x}}$ containing $\clos{V}_{\clos{x}}(W)$, where $W$ is an open disk in $\mathcal{P}-\Gamma$ containing $x$. It is easy to see that $\clos{V}_{\clos{x}}(W)\in\mathcal{V}_{\clos{u}}$  for every $\clos{u}\in\clos{V}_{\clos{x}}(W)$. Hence, thanks to (3), $\mathcal{N}\in\mathcal{V}_{\clos{u}}$ for every $\clos{u}\in\clos{V}_{\clos{x}}(W)$. Finally, the proof of (5) follows from the fact that for any $W, W'$ open disks inside $\mathcal{P}-\Gamma$ containing $x$, we have $\clos{V}_{\clos{x}}(D)\subset \clos{V}_{\clos{x}}(W)\cap\clos{V}_{\clos{x}}(W')$, where $D\subset W\cap W'$ is an open disk containing $x$. 
    
    Thanks to the above facts and by definition of the topology generated by a family of neighborhoods, the set $\mathcal{V}_{\clos{x}}$ is exactly the set of all neighborhoods of $\clos{x}$ in $\clos{\mathcal{P}}$ with respect to the topology $\mathcal{T}$, which finishes the proof of the first part of the proposition. 
    
    Take now $\clos{y}\in \clos{\Gamma}$, $S$ a standard polygon of $y$ and $\mathcal{V}_{\clos{y}}$ the set of subsets $\mathcal{N}$ of $\clos{\mathcal{P}}$ for each one of which there exist $(U_n)_{n\in \mathbb{Z}}$ a sequence of quadrant neighborhoods of $y$ in $S$ and $\gamma_0$ a reference curve for $(U_n)_{n\in \mathbb{Z}}$ and $\clos{y}$ such that $\clos{V}_{\clos{y}}(\gamma_0,(U_n)_{n\in \mathbb{Z}})\subset \mathcal{N}$. 

    Similarly to before, let us prove that the set $\mathcal{V}_{\clos{y}}$ satisfies the following properties: 
    
    \begin{enumerate}
    \item $\clos{\mathcal{P}}\in \mathcal{V}_{\clos{y}}$
    \item every element of $\mathcal{V}_{\clos{y}}$ contains $\clos{y}$ 
    \item if $\mathcal{N}\in \mathcal{V}_{\clos{y}}$ and $\mathcal{N}\subset\mathcal{N}'$, then $\mathcal{N}'\in \mathcal{V}_{\clos{y}}$
    \item for every element $\mathcal{N} $ of $\mathcal{V}_{\clos{y}}$ there exists $\mathcal{N}'\in \mathcal{V}_{\clos{y}}$ a subset of $\mathcal{N}$ such that for all $\clos{u}\in \mathcal{N}'$ we have that $\mathcal{N}\in\mathcal{V}_{\clos{u}}$ 
     \item $\mathcal{V}_{\clos{y}}$ is stable by finite intersections
     \item $\mathcal{V}_{\clos{y}}$ does not depend on our choice of standard polygon $S$
    \end{enumerate}
    As before, the proofs of (1)-(3) follow directly from the definition of $\mathcal{V}_{\clos{y}}$. Let us now explain why (4) is true. Consider $\mathcal{N}$ an element of $\mathcal{V}_{\clos{y}}$ containing $\clos{V}_{\clos{y}}(\gamma_0,(U_n)_{n\in \mathbb{Z}})$, where $(U_n)_{n\in \mathbb{Z}}$ is a sequence of quadrant neighborhoods of $y$ in $S$ and $\gamma_0$ is a reference curve for $(U_n)_{n\in \mathbb{Z}}$ and $\clos{y}$. Denote by $Q_n$ the quadrant of $y$ in $S$ containing $U_n$. For every $n\in \mathbb{Z}$ consider $U'_n\subset U_n$ a very small closed disk forming a neighborhood of $y$ in $Q_n$. The sequence $(U'_n)_{n\in \mathbb{Z}}$ defines a sequence of quadrant neighborhoods of $y$ such that $\clos{V}_{\clos{y}}(\gamma_0,(U'_n)_{n\in \mathbb{Z}})\subset \clos{V}_{\clos{y}}(\gamma_0,(U_n)_{n\in \mathbb{Z}})$. Moreover, by eventually taking the disks in $(U'_n)_{n\in \mathbb{Z}}$ smaller, it is easy to see that $\clos{V}_{\clos{y}}(\gamma_0,(U_n)_{n\in \mathbb{Z}})\in \mathcal{V}_{\clos{u}}$ for every $\clos{u}\in \clos{V}_{\clos{y}}(\gamma_0,(U'_n)_{n\in \mathbb{Z}})$, which implies thanks to (3) that $\mathcal{N} \in \mathcal{V}_{\clos{u}}$ for every $\clos{u}\in \clos{V}_{\clos{y}}(\gamma_0,(U'_n)_{n\in \mathbb{Z}})$. This finishes the proof of (4). 
    
   Next, let us prove (5). Consider $\mathcal{N},\mathcal{N}'$ two sets in $\mathcal{V}_{\clos{y}}$ containing respectively $\clos{V}_{\clos{y}}(\gamma_0,(U_n)_{n\in \mathbb{Z}})$ and $\clos{V}_{\clos{y}}(\gamma'_0,(U'_n)_{n\in \mathbb{Z}})$, where $(U_n)_{n\in \mathbb{Z}}, (U'_n)_{n\in \mathbb{Z}}$ are two sequences of quadrant neighborhoods of $y$ in $S$, $\gamma_0$ is a reference curve for $(U_n)_{n\in \mathbb{Z}}$ and $\clos{y}$, and $\gamma'_0$ is a reference curve for $(U'_n)_{n\in \mathbb{Z}}$ and $\clos{y}$. Thanks to Proposition \ref{p.changingreferencecurve}, by an eventual reindexation of $(U'_n)_{n\in \mathbb{Z}}$ and a modification of our choice of reference curve $\gamma_0'$, we can assume that $\inte{U_0}\cap \inte{U'_0}\neq \emptyset$. Furthermore, thanks to the same proposition (applied for $k=l=0$) by eventually changing once again our choices of $\gamma_0$ and $\gamma_0'$, we can assume without any loss of generality that $\gamma_0$ and $\gamma_0'$ end at the same point. Finally, 
   by Proposition \ref{p.changingreferencecurve} and Proposition  \ref{p.neighborhoodconicpoints}, by eventually reindexing $(U'_n)_{n\in \mathbb{Z}}$ we can also assume that $\gamma_0= \gamma'_0$. It follows that there exist $(U_n)_{n\in \mathbb{Z}}, (U'_n)_{n\in \mathbb{Z}}$  two sequences of quadrant neighborhoods of $y$ and $\gamma_0$ a reference curve for both $(U_n)_{n\in \mathbb{Z}}$ and $ (U'_n)_{n\in \mathbb{Z}}$ such that $$\clos{V}_{\clos{y}}(\gamma_0,(U_n)_{n\in \mathbb{Z}})\subset \mathcal{N} \text{ and } \clos{V}_{\clos{y}}(\gamma_0,(U'_n)_{n\in \mathbb{Z}})\subset \mathcal{N}'$$ 
   
   Denote by $Q_n$ the quadrant of $y$ in $S$ containing $U_n$ and $U'_n$. For every $n\in \mathbb{Z}$ consider $U''_n\subset U_n\cap U'_n$ a closed disk forming a neighborhood of $y$ inside $Q_n$. It is easy to see that $$\clos{V}_{\clos{y}}(\gamma_0,(U''_n)_{n\in \mathbb{Z}})\subset \clos{V}_{\clos{y}}(\gamma_0,(U_n)_{n\in \mathbb{Z}})\cap \clos{V}_{\clos{y}}(\gamma_0,(U'_n)_{n\in \mathbb{Z}})$$ which implies that $\mathcal{N}\cap \mathcal{N}'\in \mathcal{V}_{\clos{y}}$. This finishes the proof of (5). 

   Finally, consider $\mathcal{V}_{\clos{y},S}:=\mathcal{V}_{\clos{y}}$ and $\mathcal{V}_{\clos{y},S'}$ the set defined in the same way as  $\mathcal{V}_{\clos{y},S}$ after replacing $S$ by another standard polygon $S'$ containing $y$ in its interior. Similarly to our proof of (5),  for any two elements $\mathcal{N}\in \mathcal{V}_{\clos{y},S}$ and $\mathcal{N}'\in \mathcal{V}_{\clos{y},S'}$, we can construct a small neighborhood of $\clos{y}$ of the form $\clos{V}_{\clos{y}}(\gamma_0,(U''_n)_{n\in \mathbb{Z}})$, such that  $\clos{V}_{\clos{y}}(\gamma_0,(U''_n)_{n\in \mathbb{Z}})\in \mathcal{V}_{\clos{y},S}\cap \mathcal{V}_{\clos{y},S'}$ and $\clos{V}_{\clos{y}}(\gamma_0,(U''_n)_{n\in \mathbb{Z}})\subset \mathcal{N}\cap \mathcal{N}'$. Thanks to (3), this implies that $\mathcal{V}_{\clos{y},S}=\mathcal{V}_{\clos{y},S'}$. 

    From all the above facts and by definition of the topology generated by a family of neighborhoods, the set $\mathcal{V}_{\clos{y}}$ is exactly the set of all neighborhoods of $\clos{y}$ in $\clos{\mathcal{P}}$ with respect to the topology  $\mathcal{T}$, which proves the second part of the proposition. 
    
\end{proof}

\begin{theorem}\label{t.topologicalpropertiesbifplaneuptosurgery}
    Denote by $\widetilde{\mathcal{P}}$ the universal cover of $\mathcal{P}-\Gamma$, $\clos{\mathcal{P}}$ the bifoliated plane of $\rho$ up to surgeries on $\Gamma$ endowed with the topology $\mathcal{T}$, $\phi: \widetilde{\mathcal{P}}\rightarrow \clos{\mathcal{P}}$ and $\clos{\pi}:\clos{\mathcal{P}}\rightarrow \mathcal{P}$ the inclusion map and the projection map defined in Section \ref{s.constructionPbar12}.
    \begin{enumerate}  
    \item The map $\phi$ is a continuous embedding
    \item The set $\clos{\Gamma}$ is discrete 
    \item The map $\clos{\pi}$ is continuous on  $\clos{\mathcal{P}}$ and locally injective on  $\clos{\mathcal{P}}-\clos{\Gamma}$
    \item The topologies $\mathcal{T}$ and $\mathcal{T}_{geom}$ are equivalent
   
\end{enumerate}
 
\end{theorem}

\begin{proof}[Proof of (1)]
    Recall that by our discussion in Section \ref{s.constructionPbar12}, we have that $\text{Curv}_1\subset \text{Curv}_2$, that if $\gamma_1,\gamma_2\in \text{Curv}_1$ then $\gamma_1\sim_1\gamma_2$ if and only if $\gamma_1\sim_2\gamma_2$ (thus $\phi$ is injective) and finally that $\phi(\widetilde{\mathcal{P}})=\clos{\mathcal{P}}-\clos{\Gamma}$. Next, recall that we defined open sets for $\clos{\mathcal{P}}-\clos{\Gamma}$ and for $\widetilde{\mathcal{P}}$ in the exact same way. Using this fact together with Proposition \ref{p.baseoftopology}, we get that $\phi$ is a homeomorphism onto its image, which proves the desired result.
\end{proof}
\begin{proof} [Proof of (2)]
    Consider for every point $y\in \Gamma$ a standard polygon $S_y$ containing $y$ in its interior and such that no two distinct standard polygons $S_y, S_{y'}$ intersect (the existence of such a collection of polygons follows from the discreteness of $\Gamma$ in $\mathcal{P}$, see Proposition \ref{p.boundaryperiodicpoints}). Take for every $y\in \Gamma$ a sequence of quadrant neighborhoods of $y$ in $S$, say $(U^y_n)_{n\in \mathbb{Z}}$, and for every $\clos{y}\in \clos{\Gamma}$ denote by $\gamma_{\clos{y}}$ a reference curve of $\clos{y}$ for the sequence of quadrant neighborhoods $(U^{\clos{\pi}(\clos{y})}_n)_{n\in \mathbb{Z}}$. It suffices to prove that the neighborhoods $(\clos{V}_{\clos{y}}(\gamma_{\clos{y}},(U^{\clos{\pi}(\clos{y})}_n)_{n\in \mathbb{Z}}))_{\clos{y}\in \clos{\Gamma}}$ are two by two disjoint. 
    
    Assume that there exist $\clos{y},\clos{y'}\in \clos{\Gamma}$ such that $\clos{V}_{\clos{y}}(\gamma_{\clos{y}},(U^{\clos{\pi}(\clos{y})}_n)_{n\in \mathbb{Z}})\cap \clos{V}_{\clos{y'}}(\gamma_{\clos{y'}},(U^{\clos{\pi}(\clos{y'})}_n)_{n\in \mathbb{Z}})\neq \emptyset$. By our construction, since the $(S_y)_{y\in \Gamma}$ are two by two disjoint, we have that $\clos{\pi}(\clos{y})=\clos{\pi}(\clos{y'})$. This implies that there exists a  curve $\delta$ whose equivalence class for $\sim_2$ belongs in $\clos{V}_{\clos{y}}(\gamma_{\clos{y}},(U^{\clos{\pi}(\clos{y})}_n)_{n\in \mathbb{Z}})$ and two embedded arcs $e,e'$ in $S_y$ going from the endpoint of $\delta$ to $y$, such that the equivalence classes for $\sim_2$ of the juxtaposition of $\delta$ followed by $e$ and the juxtaposition of $\delta$ followed by $e'$ are not the same. This is impossible. 
\end{proof}
\begin{proof} [Proof of (3)]
    Indeed, if $\widetilde{\pi}$ is the projection map from $\widetilde{\mathcal{P}}$ to $\mathcal{P}$, by classical covering space theory $\widetilde{\pi}$ is continuous and locally injective on $\widetilde{\mathcal{P}}$. Thanks to Item (1) and the fact that  $\widetilde{\pi}=\clos{\pi}\circ \phi$ (see Section \ref{s.constructionPbar12}), we also get that $\clos{\pi}$ is continuous and locally injective on $\clos{\mathcal{P}}-\clos{\Gamma}$. Take now $y\in \Gamma$ and $S_y$ a small standard polygon containing $y$ in its interior. The quadrants of $y$ in $S_y$ define a periodic sequence of quadrant neighborhoods of $y$, say $(U_n)_{n\in \mathbb{Z}}$. For every $\clos{y}\in \clos{\pi}^{-1}(y)$ denote by $\gamma_{\clos{y}}$ a reference curve of $\clos{y}$ for $(U_n)_{n\in \mathbb{Z}}$. It is easy to see that $$\clos{\pi}^{-1}(S_y)=\underset{\clos{y}\in \clos{\pi}^{-1}(y)}{\cup}\clos{V}_{\clos{y}}(\gamma_{\clos{y}},(U_n)_{n\in \mathbb{Z}})$$ which gives the desired result. 
\end{proof}
\begin{proof}[Sketch of proof of (4)]
    We remind some of our notations used in Section \ref{s.geometricconstr} (see also Figure \ref{f.commutativediagram}). Let $\mathcal{P}_{exp}$ be the manifold obtained from $\mathcal{P}$ by blowing-up every point of $\Gamma$ to a circle, $\widetilde{\mathcal{P}_{exp}}$ its universal cover, $p_{exp}$ the projection from  $\mathcal{P}_{exp}$ to $\mathcal{P}$, $X_0\in \mathcal{P}_{exp}$ the unique preimage by $p_{exp}$ of the point $x_0$ (the base point of $\text{Curv}_2$) and $\text{Curv}_3$ the set of continuous arcs in $\mathcal{P}_{exp}$ starting from $X_0$. 

    \begin{figure}
        \centering
        \includegraphics[scale=0.17]{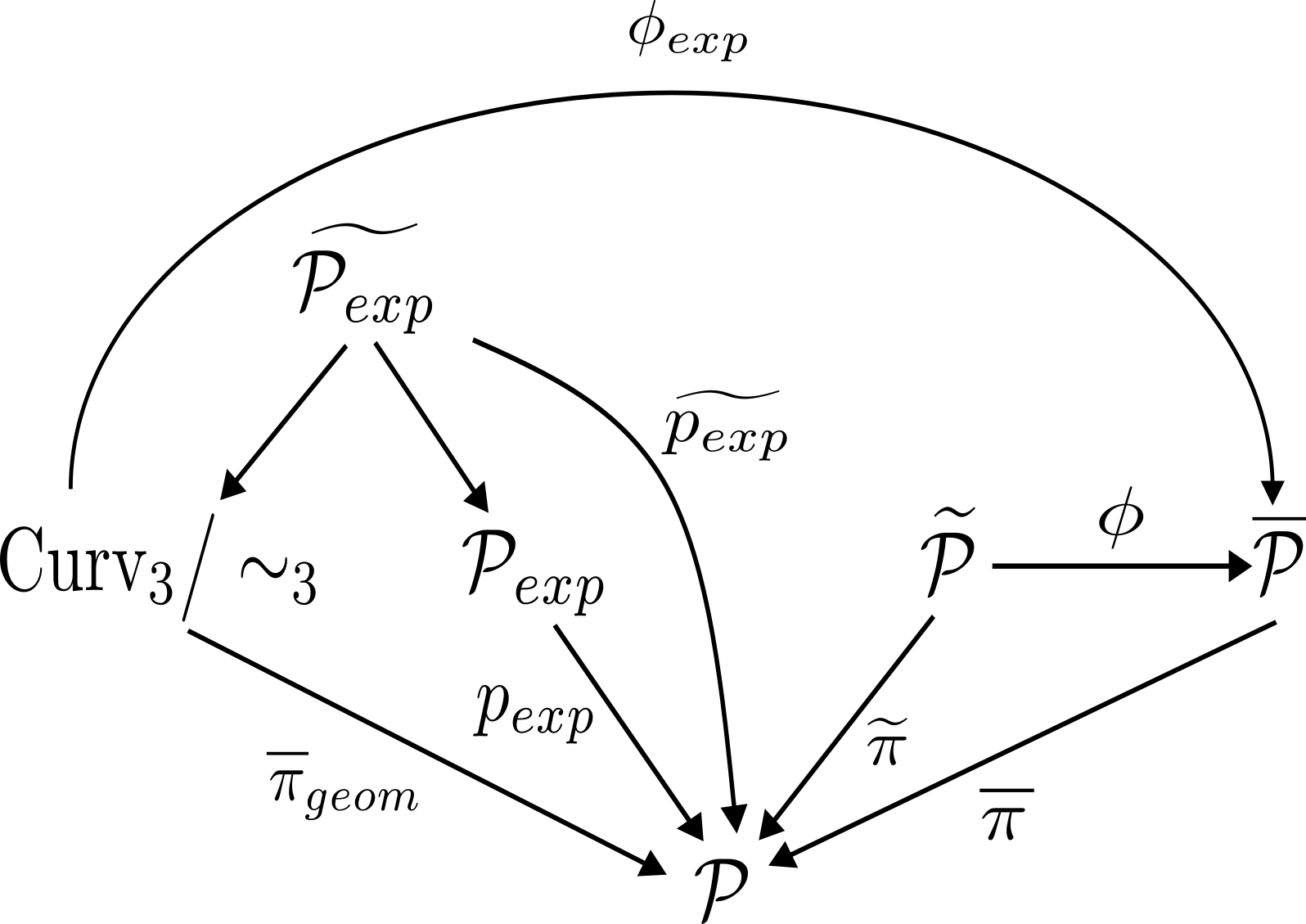}
        \caption{The above diagram is commutative}
        \label{f.commutativediagram}
    \end{figure}

    Thanks to Proposition \ref{p.sim2sim3} and the definition of $\sim_3$, we defined in Section \ref{s.geometricconstr} a bijection $\phi_{exp}$ from  $\text{Curv}_3 \big{/}\sim_3$ to $\clos{\mathcal{P}}$. Recall that $\text{Curv}_3 \big{/}\sim_3$ is homeomorphic to the space obtained by collapsing all boundary lines of $\widetilde{\mathcal{P}_{exp}}$ to points.  Also, by definition of $\mathcal{T}_{geom}$, if we endow $\clos{\mathcal{P}}$ with the topology $\mathcal{T}_{geom}$, then $\phi_{exp}$ becomes a homeomorphism. It therefore suffices to prove that $\phi_{exp}$ is also a homeomorphism when we endow $\clos{\mathcal{P}}$ with $\mathcal{T}$. 

    From now on we endow $\clos{\mathcal{P}}$ with the topology $\mathcal{T}$. Denote by $\Gamma_{geom}$ the set of points of $\text{Curv}_3 \big{/}\sim_3$ whose lifts on $\widetilde{\mathcal{P}_{exp}}$ are boundary lines and by $\widetilde{p_{exp}}$ the continuous projection from $\widetilde{\mathcal{P}_{exp}}$ to $\mathcal{P}$. Every boundary line of $\widetilde{\mathcal{P}_{exp}}$ projects by $\widetilde{p_{exp}}$ to a unique point in $\Gamma$. It follows that $\widetilde{p_{exp}}$ defines a projection map $\clos{\pi}_{geom}$ from $\text{Curv}_3 \big{/}\sim_3$ to $\mathcal{P}$. One can easily check that $\clos{\pi}_{geom}$ is continuous on $\text{Curv}_3 \big{/}\sim_3$, locally injective on $\text{Curv}_3 \big{/}\sim_3- ~ \Gamma_{geom}$ and finally that $\clos{\pi}_{geom}=\clos{\pi}\circ\phi_{exp}$. Thanks to Item (3) and the previous facts,  we get that the restriction of $\phi_{exp}$ on $\text{Curv}_3 \big{/}\sim_3- ~ \Gamma_{geom}$ defines a homeomorphism onto $\clos{\mathcal{P}}-\clos{\Gamma}$.

    It remains to prove that $\phi_{exp}$ is continuous on $\Gamma_{geom}$ and that its inverse is continuous on $\clos{\Gamma}$. Consider $y_{geom}\in \Gamma_{geom}$, $\clos{y}=\phi_{exp}(y_{geom})\in \clos{\Gamma}$, $y=\clos{\pi}_{geom}(y_{geom})=\clos{\pi}(\clos{y})\in \Gamma$,  $S$ a standard polygon in $\mathcal{P}$ containing $y$ in its interior and  $U_{geom}$ a neighborhood of $y_{geom}$ such that $\clos{\pi}_{geom}(U_{geom})\subset \inte{S}$. By eventually taking $U_{geom}$ smaller, we can assume that, if $C$ is the line boundary of $\widetilde{\mathcal{P}_{exp}}$ that was collapsed to $y_{geom}$, then the lift of $U_{geom}$ on  $\widetilde{\mathcal{P}_{exp}}$ is homeomorphic to a band $[0,1]\times \mathbb{R}$, one boundary component of which is $C$. The map $\widetilde{p_{exp}}$ sends $C$ to $y$ and the other boundary component of the previous band to a curve $\gamma$ in $S-y$ going infinitely many times around $y$. Endow $\gamma$ with some parametrization. By eventually taking $U_{geom}$ even smaller, we can assume that if $(t_n)_{n\in \mathbb{Z}}$ are the successive times in which $\gamma$ intersects the stable or unstable leaves of $y$ in $S$, then $\gamma_{|[t_n,t_{n+1}]}$ bounds a neighborhood of $y$ in a quadrant of $S$ that is homeomorphic to a closed disk. Using those neighborhoods one can construct a sequence of quadrant neighborhoods of $y$ in $S$ and prove that $\phi_{exp}(U_{geom})$ is a neighborhood of $\clos{y}$ for the topology $\mathcal{T}$. This proves the continuity of $\phi_{exp}^{-1}$ on $\clos{\Gamma}$. The continuity of  $\phi_{exp}$ on $\Gamma_{geom}$ follows from a similar argument. 
\end{proof}

\subsection{$\clos{\mathcal{P}}$ as a bifoliated plane}\label{s.foliationsinpbarconstruction}
Fix $\mathcal{P}$ a plane endowed with an orientation, $\rho: G\rightarrow \text{Homeo}(\mathcal{P})$ an orientation preserving strong Markovian action, preserving the pair of singular foliations $\mathcal{F}^s$ and $\mathcal{F}^u$ (see Definition \ref{d.folisingular}) and leaving invariant a strong Markovian family $\mathcal{R}$. Denote by $\Gamma$ the boundary periodic points of $\mathcal{R}$. Similarly to $\mathcal{P}$, the bifoliated plane of $\rho$ up to surgeries on $\Gamma$, that we will denote by $\clos{\mathcal{P}}$, can be naturally endowed with a pair of transverse ``singular" \footnote{Contrary to our definition of singular foliations (see Definition \ref{d.singularfolisurfaces}), those foliations will always admit singular points with infinitely many prongs.} foliations.

Recall that by definition of a pair of singular transverse foliations (see Definitions \ref{d.singularfolisurfaces} and \ref{d.transversefolisurfaces}), for every $x\in \mathcal{P}$ there exist a small neighborhood $U_x$ of $x$ and a homeomorphism $H_x: U_x\rightarrow \mathbb{R}^2$ conjugating the restriction on $U_x$ of the pair $(\mathcal{F}^s,\mathcal{F}^u)$ to a smooth (potentially singular) model of singular transverse foliations. Using the fact that $\Gamma$ is discrete, we can endow $\mathcal{P}$ with a smooth structure such that 
for every $\gamma \in \Gamma$ the map $H_\gamma$ is a $C^{\infty}$ diffeomorphism. 

\begin{figure}[h!]
    \centering
    \includegraphics[width=0.5\linewidth]{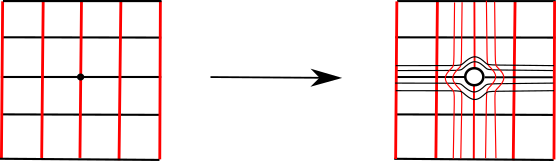}
    \caption{The foliations  $\mathcal{F}_{exp}^s$ and  $\mathcal{F}_{exp}^u$ after blowing up a non-singular point in $\Gamma$}
    \label{f.explodefoli}
\end{figure}

Using the previous smooth structure on $\mathcal{P}$, blow-up each point of $\Gamma$ to a circle, in order to obtain a plane minus a countable set of open disks, say $\mathcal{P}_{exp}$, and a pair of foliations in the interior of $\mathcal{P}_{exp}$, say $\mathcal{F}_{exp}^s$ and $\mathcal{F}_{exp}^u$, such that:  
\begin{itemize}
    \item Thanks to Remark \ref{r.singularitiesareboundaries} and Proposition \ref{p.propertiesoffoliformarkovianactions}, any leaf in $\mathcal{F}_{exp}^s$ or $\mathcal{F}_{exp}^u$ is a properly embedded line in the interior of $\mathcal{P}_{exp}$ (see Figure \ref{f.explodefoli})
    \item Let $p_{exp}$ is the projection map from $\mathcal{P}_{exp}$ to $\mathcal{P}$. The preimage by $p_{exp}$ of any leaf in $\mathcal{F}^s$ (resp. $\mathcal{F}^u$) that does not intersect $\Gamma$ is a leaf in $\mathcal{F}_{exp}^s$ (resp. $\mathcal{F}_{exp}^u$)
    \item If $\gamma\in \Gamma$ and $\mathcal{F}_+^s(\gamma)$ (resp. $\mathcal{F}_+^u(\gamma)$) is a stable (resp. unstable) separatrix of $\gamma$, then the preimage by $p_{exp}$ of $\mathcal{F}_+^s(\gamma)-\{\gamma\}$ (resp. $\mathcal{F}_+^u(\gamma)-\{\gamma\}$) is a leaf in $\mathcal{F}_{exp}^s$ (resp. $\mathcal{F}_{exp}^u$). Furthermore, thanks to our choice of smooth structure on $\mathcal{P}$, any such leaf $l$ of $\mathcal{F}_{exp}^s$ or $\mathcal{F}_{exp}^u$, can be continuously extended to the boundary of $\mathcal{P}_{exp}$ via the addition of one point (see Figure \ref{f.explodefoli})
    \item $\mathcal{F}_{exp}^s$ and $\mathcal{F}_{exp}^u$ define a pair of transverse foliations (with no singularities) in the interior of $\mathcal{P}_{exp}$ 
\end{itemize}

\begin{figure}[h]

  \begin{minipage}[ht]{0.4\textwidth}
    \centering 
     \vspace{0.8cm}
    \includegraphics[width=1\textwidth]{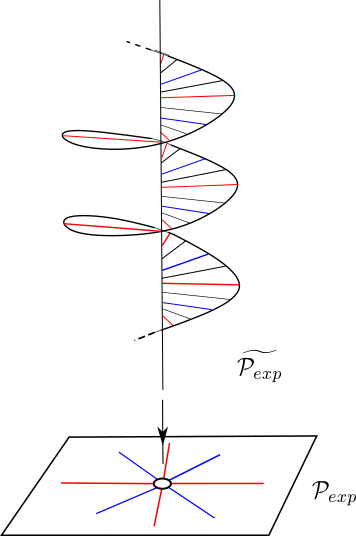}
  \hspace{-1cm}

  \end{minipage}
 \begin{minipage}[ht]{0.4\textwidth}
 \centering
    \includegraphics[width=0.8\textwidth]{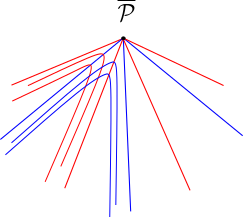}
    \hspace{-1cm}

  \end{minipage}
  
  \caption{In the figure in the left, we see $\widetilde{\mathcal{P}_{exp}}$ endowed with the lifts of the foliations $\mathcal{F}_{exp}^s$ and $\mathcal{F}_{exp}^u$. In the figure in the right, we collapse the line boundary of the left picture to a point in order to obtain $\clos{\mathcal{P}}$ together with $\clos{\mathcal{F}^s}$ and $\clos{\mathcal{F}^u}$.}
  \label{f.foliationspbar}
\end{figure}

Consider $\widetilde{\mathcal{P}_{exp}}$ (see Figure \ref{f.foliationspbar}), the universal cover of $\mathcal{P}_{exp}$, endowed with the lifts of the foliations $\mathcal{F}_{exp}^s$ and $\mathcal{F}_{exp}^u$. After collapsing every boundary line of $\widetilde{\mathcal{P}_{exp}}$ to a point, thanks to our geometric construction of $\clos{\mathcal{P}}$ and Theorem  \ref{t.topologicalpropertiesbifplaneuptosurgery}, we get a space homeomorphic to $\clos{\mathcal{P}}$ and a pair of transverse foliations, say $\clos{\mathcal{F}_{inc}^s}$\footnote{inc stands for incomplete} and $\clos{\mathcal{F}_{inc}^u}$, such that: 
\begin{itemize}
    \item Denote by $\clos{\pi}$ the projection map from $\clos{\mathcal{P}}$ to $\mathcal{P}$. By identifying with $\clos{\mathcal{P}}$ the space obtained by collapsing every boundary line of $\widetilde{\mathcal{P}_{exp}}$ to a point, we have that every leaf in $\clos{\mathcal{F}_{inc}^s}$ or $\clos{\mathcal{F}_{inc}^u}$ is a properly embedded line in $\clos{\mathcal{P}}-\clos{\Gamma}$, where $\clos{\Gamma}:=\clos{\pi}^{-1}(\Gamma)$ 
    \item Each connected component of the preimage by $\clos{\pi}$ of any leaf in $\mathcal{F}^s$ (resp. $\mathcal{F}^u$) that does not intersect $\Gamma$ is a leaf in $\clos{\mathcal{F}_{inc}^s}$ (resp. $\clos{\mathcal{F}_{inc}^u}$)
    \item If $\gamma\in \Gamma$  and $\mathcal{F}_+^s(\gamma)$ (resp. $\mathcal{F}_+^u(\gamma)$) is a stable (resp. unstable) separatrix of $\gamma$, then each connected component of the preimage by $\clos{\pi}$ of $\mathcal{F}_+^s(\gamma)-\{\gamma\}$ (resp. $\mathcal{F}_+^u(\gamma)-\{\gamma\}$) is a leaf in  $\clos{\mathcal{F}_{inc}^s}$ (resp. $\clos{\mathcal{F}_{inc}^u}$). Furthermore,  any such leaf of $\clos{\mathcal{F}_{inc}^s}$ or $\clos{\mathcal{F}_{inc}^u}$ can be continuously extended to $\clos{\mathcal{P}}$ via the addition of one point
    \item $\clos{\mathcal{F}_{inc}^s}$ and $\clos{\mathcal{F}_{inc}^u}$ define a pair of transverse foliations (with no singularities) in $\clos{\mathcal{P}}-\clos{\Gamma}$
\end{itemize}

 By continuously extending in $\clos{\mathcal{P}}$ every leaf of $\clos{\mathcal{F}_{inc}^s}$ and $\clos{\mathcal{F}_{inc}^u}$ that projects via $\clos{\pi}$ to the interior of a stable or unstable separatrix of a point in $\Gamma$, we produce in $\clos{\mathcal{P}}$ a pair of ``transverse singular foliations" $\clos{\mathcal{F}^s}$ and $\clos{\mathcal{F}^u}$. More precisely,  
\begin{rema}\label{r.classificationofleavesoffolipbar}
Every leaf $\clos{L}$ of $\clos{\mathcal{F}^s}$ (resp. $\clos{\mathcal{F}^u}$) satisfies exactly one of the following :  

\begin{itemize}
    \item $\clos{L}$ is homeomorphic to $\mathbb{R}$, $\clos{\pi}$ is injective on $\clos{L}$ and $\clos{\pi}(\clos{L})$ is a leaf of $\mathcal{F}^s$ (resp. $\mathcal{F}^u$) that does not intersect $\Gamma$ 
    \item $\clos{L}$ is homeomorphic to $\mathbb{R}_{\geq 0}$, $\clos{\pi}$ is injective on $\clos{L}$ and $\clos{\pi}(\clos{L})$ is a stable (resp. unstable) separatrix of a point in $\Gamma$ 
\end{itemize}
\end{rema}
We will call the foliations $\clos{\mathcal{F}^s}$ and $\clos{\mathcal{F}^u}$ the \emph{stable} and \emph{unstable  foliations} of $\clos{\mathcal{P}}$. By construction, $\clos{\mathcal{F}^s}$ and $\clos{\mathcal{F}^u}$ define a pair of transverse foliations inside $\clos{\mathcal{P}}-\clos{\Gamma}$ and also every point in $\clos{\Gamma}$ is singular, as it admits infinitely many prongs. Despite the fact that $\clos{\mathcal{F}^s}$ and $\clos{\mathcal{F}^u}$ do not satisfy the criteria of Definitions \ref{d.singularfolisurfaces} and \ref{d.pseudohomeo}, by an abuse of language, we will say that $\clos{\mathcal{F}^s}$ and $\clos{\mathcal{F}^u}$ form a pair of transverse singular foliations. 

Moreover, even though $\clos{\mathcal{F}^s}$ and $\clos{\mathcal{F}^u}$ may be much more complex than the foliations $\mathcal{F}^s$ or $\mathcal{F}^u$, we have that: 
\begin{rema}\label{r.propertiesfolipbar}
    \begin{enumerate}
        \item Every singular point of $\clos{\mathcal{F}^s}$ and $\clos{\mathcal{F}^u}$ belongs  in $\clos{\Gamma}$ and conversely every point in $\clos{\Gamma}$ is singular and admits infinitely many prongs
        \item The restrictions on $\clos{\mathcal{P}}-\clos{\Gamma}$ of the foliations $\clos{\mathcal{F}^s}$ and $\clos{\mathcal{F}^u}$ are orientable and transversely orientable. 
        
        \item The leaves of the restrictions on $\clos{\mathcal{P}}-\clos{\Gamma}$ of the foliations  $\clos{\mathcal{F}^s}$ and $\clos{\mathcal{F}^u}$ are closed, non-compact and properly embedded into $\clos{\mathcal{P}}-\clos{\Gamma}$. 
        \item Any stable leaf in $\clos{\mathcal{F}^s}$ intersects at most once any unstable leaf in $\clos{\mathcal{F}^u}$. 
        \item The complement of any stable or unstable leaf in $\clos{\mathcal{P}}$ consists of exactly two simply connected components.
    \end{enumerate}
\end{rema}

Indeed, Item (1) is a direct consequence of our construction of $\clos{\mathcal{F}^s}$ and $\clos{\mathcal{F}^u}$. Item (2) results from the fact that $\clos{\mathcal{P}}-\clos{\Gamma}$ is homeomorphic to a plane and that the restrictions on $\clos{\mathcal{P}}-\clos{\Gamma}$ of the foliations $\clos{\mathcal{F}^s}$ and $\clos{\mathcal{F}^u}$ define two line foliations. Item (3) is a consequence of Proposition  \ref{p.propertiesoffoliformarkovianactions} and Item (4) is a consequence of the same proposition and  Remark \ref{r.classificationofleavesoffolipbar}. Finally, one can prove Item (5) by using the fact that $\clos{\mathcal{P}}$ is a topological plane with some points at infinity and the fact that $\clos{\mathcal{F}^s}$ and $\clos{\mathcal{F}^u}$ were obtained by a continuous extension of two foliations on the plane $\clos{\mathcal{P}}-\clos{\Gamma}$ consisting of properly embedded leaves.

Finally, before finishing this section let us describe the form of the foliations $\clos{\mathcal{F}^s}$ and $\clos{\mathcal{F}^u}$ near a singular point in $\clos{\Gamma}$:

\begin{prop}\label{p.singularitiesoffoliations}
For every $\clos{\gamma}\in \clos{\Gamma}$, the stable leaves in $\clos{\mathcal{F}^{s}}$ that intersect $\clos{\gamma}$  form a countable set of leaves that can be indexed by $\mathbb{Z}$ --denote by $...,s_{-2},s_{-1},s_0,s_1,s_2...$ the previous indexation-- so that:
\begin{itemize}
\item if $\gamma:=\clos{\pi}(\clos{\gamma})$ and $\gamma$ is a $p$-prong singularity in $\mathcal{F}^s$ (with $p\geq 2$), then for all $k, l\in \mathbb{Z}$  such that $(k-l)$ is a multiple of $p$ we have that $\clos{\pi}(s_k)=\clos{\pi}(s_l)$
\item for all $k\in \mathbb{Z}$ $s_k\cap \clos{\Gamma}= \lbrace \clos{\gamma} \rbrace$
\item for all $k,m\in \mathbb{Z}$ $s_k$ is not separated from $s_{m}$ in $\clos{\mathcal{P}}$ (i.e. $s_k-\lbrace \clos{\gamma} \rbrace$ and  $s_{m}-\lbrace \clos{\gamma} \rbrace$ are non-separated stable leaves in $\clos{\mathcal{P}}-\clos{\Gamma}$) if and only if $|k-m|=1$
\end{itemize}
\end{prop} 
\begin{proof}
Indeed, take $\clos{\gamma}\in \clos{\Gamma}$ and $\gamma=\clos{\pi}(\clos{\gamma})$. Recall that when we collapsed every boundary line of $\widetilde{\mathcal{P}_{exp}}$ to a point, we obtained a space homeomorphic to $\clos{\mathcal{P}}$. Using this fact, denote by  $\widetilde{C}$ the boundary line in $\widetilde{\mathcal{P}_{exp}}$ that was collapsed to $\clos{\gamma}$ and by $C$ the projection of $\widetilde{C}$ in $\mathcal{P}_{exp}$ (see Figure \ref{f.commutativediagram}). Assume that $\gamma$ is a $p$-prong singularity in $\mathcal{F}^s$. By our construction of $\mathcal{F}^s_{exp}$ and thanks to the fact that the leaves of $\mathcal{F}^s$ are closed, we have that only leaves in $\mathcal{F}^s_{exp}$ that get arbitrarily close to the boundary circle $C$ are the $p$ leaves in $\mathcal{F}^s_{exp}$ that project via $p_{exp}$ to the stable separatrices of $\gamma$ (see Figure \ref{f.explodefoli}). Extend continuously the previous leaves by a point so that they intersect $C$ and denote by $L_1,...,L_p$ the previous collection of extended leaves. Assume without any loss of generality that $L_1,...,L_p$ are cyclically ordered. Lift everything on $\widetilde{\mathcal{P}_{exp}}$, the universal cover of $\mathcal{P}_{exp}$, and consider the subset of all the lifts of $L_1,...,L_p$ on $\widetilde{\mathcal{P}_{exp}}$ that intersect $\widetilde{C}$ (see Figure \ref{f.foliationspbar}). We can index the previous set of lifts by $\mathbb{Z}$ --denote by  $(\widetilde{L_n})_{n\in\mathbb{Z}}$ the previous indexation-- so that :

\begin{itemize}
    \item for any two distinct $i,j\in \mathbb{Z}$ we have that $\widetilde{L_i}\cap \widetilde{L_j}=\emptyset$
    \item if $n\equiv i \text{ mod }p$, where $i\in \llbracket 1, p \rrbracket$, then $\widetilde{L_n}$ is a lift of $L_i$ intersecting $\widetilde{C}$
    \item for any $n\in \mathbb{Z}$, if $\widetilde{x_i}$ the unique point in $\widetilde{L_i}\cap \widetilde{C}$, then there does not exist $m\in \mathbb{Z}$ such that $\widetilde{x_m}$ lies in the interior of the interval in $\widetilde{C}$ defined by $\widetilde{x_n}$ and $\widetilde{x_{n+1}}$. In other words, the points in $(\widetilde{x_n})_{n\in\mathbb{Z}}$ form the successive intersections of the line $\widetilde{C}$ with the lifts of $L_1,...,L_p$. 
\end{itemize}

The first point of the above proposition is obtained by simply collapsing $\widetilde{C}$ into a point and by using the fact that the diagram in Figure \ref{f.commutativediagram} is commutative. Next, if $s_k$ intersects two points in $\clos{\Gamma}$, then so does its projection by $\clos{\pi}$. This is impossible, because by our definition of a Markovian action, a stable leaf in $\mathcal{F}^s$ contains at most one periodic point. This proves the second point of the proposition. Finally, the fact that $s_k$ is not separated from $s_{k+1}$ can be seen in Figure \ref{f.foliationspbar}. Conversely, if $|l-k|\neq 1$, then $s_k$ and $s_l$ are separated because for any $n\in \mathbb{Z}$ the line $s_n-\{\clos{\gamma}\}$ is properly embedded in  $\clos{\mathcal{P}}-\clos{\Gamma}$ and thus separates  $\clos{\mathcal{P}}-\clos{\Gamma}$ into two connected components. 
\end{proof}
\subsection{A group action on $\clos{\mathcal{P}}$}\label{s.groupactiondefi}
Fix $\mathcal{P}$ a plane endowed with an orientation, $\rho: G\rightarrow \text{Homeo}(\mathcal{P})$ an orientation preserving strong Markovian action, preserving the pair of singular foliations $\mathcal{F}^s$ and $\mathcal{F}^u$ and leaving invariant a strong Markovian family $\mathcal{R}$. Denote by $\Gamma$ the boundary periodic points of $\mathcal{R}$ and by $\clos{\mathcal{P}}$ the bifoliated plane of $\rho$ up to surgeries on $\Gamma$.

Our goal in this section consists in showing that, in addition to a pair of transverse singular foliations $\clos{\mathcal{F}^s},\clos{\mathcal{F}^u}$,  the bifoliated plane of $\rho$ up to surgeries on $\Gamma$ admits a natural group action by homeomorphisms. Intuitively, recall that our construction of $\clos{\mathcal{P}}$ (see Section \ref{s.constructionPbar12}) resembled very closely the classical construction of $\widetilde{\mathcal{P}-\Gamma}$, the universal cover of $\mathcal{P}-\Gamma$. By lifting on $\widetilde{\mathcal{P}-\Gamma}$ all the homeomorphisms in $\rho(G)$, we obtain a group action (this group is generally not isomorphic to $G$) on $\widetilde{\mathcal{P}-\Gamma}$ by homeomorphisms. By treating $\clos{\mathcal{P}}$ as the universal cover of $\mathcal{P}$ with respect to the homotopy $\sim_2$, we will similarly define a group action on $\clos{\mathcal{P}}$ by a group that is larger than $G$.

More precisely, recall that by our construction in Section \ref{s.constructionPbar12}, $\clos{\mathcal{P}}=\text{Curv}_2 \big{/}\sim_2$, where $$\text{Curv}_2= \lbrace \gamma:[0,1]\overset{C^0}{\rightarrow} \mathcal{P} | ~ \gamma(0)=x_0, \gamma[0,1)\cap \Gamma= \emptyset \rbrace $$ and $x_0\in \mathcal{P}-\Gamma$. We define 
\begin{align*}
   \text{Curv}_2(x_0)&:=\lbrace \gamma:[0,1]\overset{C^0}{\rightarrow} \mathcal{P} | ~ \gamma(0)=x_0, \gamma(1)\in \rho(G)(x_0), \gamma[0,1)\cap \Gamma= \emptyset \rbrace \\ &:=\{\gamma\in \text{Curv}_2|\gamma(1)\in \rho(G)(x_0)\} 
\end{align*} Consider $\gamma\in \text{Curv}_2$, $\gamma_0\in \text{Curv}_2(x_0)$
and $g\in G$ such that $\rho(g)(\gamma_0(0))=\rho(g)(x_0)=\gamma_0(1)$. We define $\clos{\rho}(g,\gamma_0)(\gamma)$ to be the curve in $\text{Curv}_2$ obtained by the juxtaposition of $\gamma_0$ followed by $\rho(g)(\gamma)$.

It is easy to see that 
\begin{prop} Consider $\gamma_0,\gamma_0'\in \text{Curv}_2(x_0)$, $\gamma,\gamma'\in \text{Curv}_2$ and $g\in G$ such that  $\rho(g)(\gamma_0(0))=\rho(g)(\gamma'_0(0))=\rho(g)(x_0)=\gamma_0(1)=\gamma'_0(1)$
    \begin{itemize}
        \item $\gamma_0\sim_2\gamma_0'$ if and only if $\clos{\rho}(g,\gamma_0)(\gamma)\sim_2 \clos{\rho}(g,\gamma'_0)(\gamma)$
        \item $\gamma\sim_2\gamma'$ if and only if $\clos{\rho}(g,\gamma_0)(\gamma)\sim_2 \clos{\rho}(g,\gamma_0)(\gamma')$
    \end{itemize}
\end{prop} For any $\delta\in \text{Curv}_2$, denote by $\clos{\delta}$ the class in $\text{Curv}_2 \big{/}\sim_2=\clos{\mathcal{P}}$ containing $\delta$. Thanks to the above proposition, $\clos{\rho}(g,\clos{\gamma_0})$ defines an injective map from $\clos{\mathcal{P}}$ to $\clos{\mathcal{P}}$. The previous map is also surjective. Indeed, consider $\delta\in \text{Curv}_2$ and $\delta'$ a curve in $\mathcal{P}$ starting at the endpoint of $\gamma_0$ and such that the juxtaposition of $\gamma_0$ followed by $\delta'$ defines a curve in $\text{Curv}_2$ homotopic to $\delta$ for $\sim_2$. By construction, we have that $\clos{\rho}(g,\clos{\gamma_0})(\clos{\rho(g^{-1})(\delta')})=\clos{\delta}$. 

Furthermore, we have that: 
\begin{prop}\label{p.groupactiondefi} Let $\gamma_0\in \text{Curv}_2(x_0)$ and $g\in G$ be   such that $\rho(g)(\gamma_0(0))=\rho(g)(x_0)=\gamma_0(1)$. We have that 
    \begin{enumerate}
        \item $\clos{\pi}\circ \clos{\rho}(g,\clos{\gamma_0})=\rho(g)\circ \clos{\pi}$
        \item for any $x\in \mathcal{P}$, $\clos{x},\clos{x'}\in \clos{\pi}^{-1}(x)$ there exists $\delta_0$ a closed curve in $\text{Curv}_2(x_0)$ such that $\clos{x'}=\clos{\rho}(e,\clos{\delta_0})(\clos{x})$, where $e$ the trivial element in $G$
        \item if $\delta_0$ is a closed in $\text{Curv}_2(x_0)$, then $\clos{\rho}(e,\clos{\delta_0})=id$ if and only if $\delta_{0}$ is homotopically trivial for $\sim_2$
        \item if $g'\in G$ and $\gamma_1\in \text{Curv}_2(x_0)$ are such that $\rho(g')(\gamma_1(0))=\rho(g')(x_0)=\gamma_1(1)$, then $$\clos{\rho}(g,\clos{\gamma_0})\circ \clos{\rho}(g',\clos{\gamma_1})=\clos{\rho}(gg',\clos{\rho(g)(\gamma_1)\gamma_0})$$ where $\rho(g)(\gamma_1)\gamma_0$ denotes the curve in $\text{Curv}_2(x_0)$ obtained by the juxtaposition of $\gamma_0$ followed by $\rho(g)(\gamma_1)$. In particular, $$\clos{\rho}(g,\clos{\gamma_0})\circ \clos{\rho}(g^{-1},\clos{\rho(g^{-1})(\gamma_0^-)})= \clos{\rho}(g^{-1},\clos{\rho(g^{-1})(\gamma_0^-)})\circ \clos{\rho}(g,\clos{\gamma_0})=id$$ where $\gamma_0^-:[0,1] \rightarrow \mathcal{P}$  denotes the curve defined by $\gamma_0^-(t)=\gamma_0(1-t)$ for every $t\in [0,1]$.
        \item $\clos{\rho}(g,\clos{\gamma_0}):\clos{\mathcal{P}}\rightarrow \clos{\mathcal{P}}$ is a homeomorphism 
    \end{enumerate}
\end{prop}
\begin{proof}[Proof of (1)]
    Consider $\gamma\in \text{Curv}_2$. By construction, $\clos{\rho}(g,\clos{\gamma_0})(\clos{\gamma})\in \clos{\mathcal{P}}=\text{Curv}_2 \big{/}\sim_2$ is a class of curves homotopic for $\sim_2$ containing the curve obtained by the juxtaposition of $\gamma_0$ followed by $\rho(g)(\gamma)$. By definition of $\clos{\pi}$, the point $\clos{\pi}( \clos{\rho}(g,\clos{\gamma_0})(\clos{\gamma}))$ is equal to the endpoint of the previous curve, which is equal to the image by $\rho(g)$ of $\gamma(1)$. This gives us the desired result. 
\end{proof}
\begin{proof}[Proof of (2)]
Consider $\gamma\in \text{Curv}_2$ and $\gamma'\in \text{Curv}_2$ belonging in the equivalence classes defined by $\clos{x}$ and $\clos{x'}$ respectively. Notice that $\gamma(0)=\gamma'(0)=x_0$ as $\gamma,\gamma'\in \text{Curv}_2$ and $\gamma(1)=\gamma'(1)=x$ as $\clos{\pi}(\clos{x})=\clos{\pi}(\clos{x'})=x$. It is easy now to see that there exists $\delta_0$ a closed curve in $ \text{Curv}(x_0)$ such that the curve obtained by the juxtaposition of $\delta_0$ followed by $\gamma$ is homotopic to $\gamma'$ for $\sim_2$. This gives us the desired result. 
\end{proof}
\begin{proof}[Proof of (3)] 
     Consider $\gamma\in \text{Curv}_2$ and $\delta_0$ a closed curve in $ \text{Curv}(x_0)$. By construction, $\clos{\rho}(e,\clos{\delta_0})(\clos{\gamma})\in \clos{\mathcal{P}}=\text{Curv}_2 \big{/}\sim_2$ corresponds to a class of curves homotopic for $\sim_2$ containing the curve obtained by the juxtaposition of $\delta_0$ followed by $\gamma$. If  $\delta_0$ is homotopically trivial for $\sim_2$, the previous curve is homotopic to $\gamma$ for $\sim_2$. Conversely, if for any curve $\gamma$ in $\text{Curv}_2$, the curve obtained as the juxtaposition of $\delta_0$ followed by $\gamma$ is homotopic to $\gamma$ for $\sim_2$, then $\delta_0$ is homotopically trivial for $\sim_2$. 
\end{proof}
\begin{proof}[Proof of (4)]
    Consider $\gamma\in \text{Curv}_2$. By construction, $\clos{\rho}(g',\clos{\gamma_1})(\clos{\gamma})\in \clos{\mathcal{P}}=\text{Curv}_2 \big{/}\sim_2$ is a class of curves homotopic for $\sim_2$ containing the curve obtained by the juxtaposition of $\gamma_1$ followed by $\rho(g')(\gamma)$. Similarly, $\clos{\rho}(g,\clos{\gamma_0})\circ \clos{\rho}(g',\clos{\gamma_1})(\clos{\gamma})$ is a class of curves homotopic for $\sim_2$ containing the curve obtained by the juxtaposition of $\gamma_0$ followed by $\rho(g)(\gamma_1)$ followed by $\rho(g)(\rho(g')(\gamma))$, which gives us the desired result. 
\end{proof}
\begin{proof}[Proof of (5)]
   We have already established that $\clos{\rho}(g,\clos{\gamma_0}):\clos{\mathcal{P}}\rightarrow \clos{\mathcal{P}}$ is a bijection. Let us now prove that it is a homeomorphism by using Proposition \ref{p.baseoftopology}. 
   
   Consider $\clos{x}\in \clos{\mathcal{P}}-\clos{\Gamma}$, $\clos{y}\in \clos{\Gamma}$, $x=\clos{\pi}(\clos{x})$, $y=\clos{\pi}(\clos{y})$, $W$ a small open disk in $\mathcal{P}-\Gamma$ around $x$, $S$ a standard polygon containing $y$, $(U_n)_{n\in \mathbb{Z}}$ a sequence of quadrant neighborhoods of $y$ in $S$ and finally $\delta_0$ a reference curve of $\clos{y}$ for the sequence of neighborhoods $(U_n)_{n\in \mathbb{Z}}$. Denote respectively by $\clos{V}_{\clos{x}}(W)$ and by $\clos{V}_{\clos{y}}(\delta_0,(U_n)_{n\in \mathbb{Z}})$ the neighborhoods of $\clos{x}$ and $\clos{y}$ that were constructed in Section \ref{s.secondconstruction}. Recall that if $\gamma_1\in\text{Curv}_2$ belongs in the equivalence class defined by $\clos{x}\in \clos{\mathcal{P}}=\text{Curv}_2 \big{/}\sim_2$ and $s\in (0,1)$ is such that $\gamma_1([s,1])\subset W$, then $$\clos{V}_{\clos{x}}(W)=\lbrace\clos{\gamma}| \gamma\in \text{Curv}_2, \forall t\in [0,s] \quad \gamma(t)=\gamma_1(t) , \gamma([s,1])\subset W \rbrace$$ It follows that if $\rho(g)(\gamma_1)\gamma_0: [0,1] \rightarrow \mathcal{P}$ is the curve obtained by the juxtaposition of $\gamma_0$ followed by $\rho(g)(\gamma_1)$ and $s'\in (0,1)$ is such that $(\rho(g)(\gamma_1)\gamma_0)([s',1])\subset \rho(g)(W)$ then 
   \begin{align*}
     \clos{\rho}(g,\clos{\gamma_0})(\clos{V}_{\clos{x}}(W))&=\lbrace\clos{\gamma}| \gamma\in \text{Curv}_2, \forall t\in [0,s'] \quad \gamma(t)=(\rho(g)(\gamma_1)\gamma_0)(t) , \gamma([s',1])\subset \rho(g)(W) \rbrace\\ &=\clos{V}_{\clos{\rho}(g,\clos{\gamma_0})(\clos{x})}(\rho(g)(W))  
   \end{align*}
   
   Similarly, using the definition of $\clos{\rho}(g,\clos{\gamma_0})$ and $\clos{V}_{\clos{y}}(\delta_0,(U_n)_{n\in \mathbb{Z}})$, one can easily check that
    $\clos{\rho}(g,\clos{\gamma_0})(\clos{V}_{\clos{y}}(\delta_0,(U_n)_{n\in \mathbb{Z}}))$ is equal to the following neighborhood of $\clos{\rho}(g,\clos{\gamma_0})(\clos{y})$ $$\clos{V}_{\clos{\rho}(g,\clos{\gamma_0})(\clos{y})}(\rho(g)(\delta_0)\gamma_0,(\rho(g)(U_n))_{n\in \mathbb{Z}})$$ where $\rho(g)(\delta_0)\gamma_0$ denotes the curve obtained by the juxtaposition of $\gamma_0$ followed by $\rho(g)(\delta_0)$. Thanks to Proposition \ref{p.baseoftopology}, our previous arguments prove the continuity of the inverse of $\clos{\rho}(g,\clos{\gamma_0})$. The fact that $\clos{\rho}(g,\clos{\gamma_0})$ defines a homeomorphism follows from Item (4). 
\end{proof}

Consider $$\clos{G}=\{(g,\clos{\gamma})|g\in G \text{, } \clos{\gamma}\in \text{Curv}_2 (x_0)\big{/}\sim_2\text{ and } \rho(g)(x_0)=\clos{\pi}(\clos{\gamma})\}$$ We endow $\clos{G}$ with the group structure given by $$(g,\clos{\gamma_0})\cdot(g',\clos{\gamma_1})=(gg',\clos{\rho(g)(\gamma_1)\gamma_0})$$ where $\clos{\rho(g)(\gamma_1)\gamma_0}$ is the equivalence class in $\text{Curv}_2 (x_0)\big{/}\sim_2$ obtained by the juxtaposition of a representative of $\clos{\gamma_0}$ followed by the image by $\rho(g)$ of a representative of $\clos{\gamma_1}$. It is easy to check, thanks to Proposition \ref{p.groupactiondefi}, that the previous operation defines indeed group structure on $\clos{G}$ such that: 
\begin{itemize}
    \item if $e$ is the trivial element in $G$ and $\gamma_{triv}\in \text{Curv}_2(x_0)$ a homotopically trivial loop for $\sim_2$ then $(e,\clos{\gamma_{triv}})$ is the identity element in $\clos{G}$ 
    \item $(g,\clos{\gamma_0})^{-1}= (g^{-1},\clos{\rho(g^{-1})(\gamma_0^-)})$, where $\gamma_0^{-}:[0,1]\rightarrow \mathcal{P}$ is defined as follows : if $\gamma_0:[0,1]\rightarrow \mathcal{P}$ is a representative in $\text{Curv}_2$ if the class $\clos{\gamma_0}\in \text{Curv}_2 (x_0)\big{/}\sim_2$, then $\gamma_0^{-}(t)=\gamma_0(1-t)$ for every $t\in [0,1]$ 
\end{itemize}
We will call the previous group, the \emph{extension of $G$ associated to $\Gamma$}. Notice that since $G$ is countable (see Definition \ref{d.markovianaction}) and also since the homotopy classes in $\text{Curv}_2(x_0)$ for $\sim_2$ are countable, we get that 

\begin{rema}\label{r.extensioncountable}
    The group $\clos{G}$ is a countable group. 
\end{rema}

Denote by $\mathfrak{p}$ the map associating every $(g,\clos{\gamma})\in \clos{G}$ to $g\in G$. The map defines a homomorphism from $\clos{G}$ to $G$, whose kernel consists of the elements of $\clos{G}$ of the form $(e,\clos{\gamma})$, where $\clos{\gamma}$ is a class of closed curves in $\text{Curv}_2(x_0)$ that are homotopic for $\sim_2$. Since $x_0\in \mathcal{P}-\Gamma$ and two closed curves in $\text{Curv}_2(x_0)$ are homotopic in $\mathcal{P}-\Gamma$ if and only if they are homotopic for $\sim_2$, we have that the kernel of $\mathfrak{p}$ is identified with $\pi_1(\mathcal{P}-\Gamma,x_0)$. Recall now that by our geometric construction of $\clos{\mathcal{P}}$, the space $\clos{\mathcal{P}}-\clos{\Gamma}$ can be identified with the universal cover of $\mathcal{P}-\Gamma$. Furthermore, notice that  when $\gamma$ is a closed curve in $\text{Curv}_2(x_0)$, our definition of $\clos{\rho}(e,\clos{\gamma})\in \text{Homeo}(\clos{\mathcal{P}}-\clos{\Gamma})$ coincides with the definition of the deck transformation on the universal cover of $\mathcal{P}-\Gamma$ associated to the homotopy class of $\gamma$ in $\pi_1(\mathcal{P}-\Gamma,x_0)$; hence, the action of $\pi_1(\mathcal{P}-\Gamma,x_0)\leq \clos{G}$ on $\clos{\mathcal{P}}-\clos{\Gamma}$ can be identified with the action of $\pi_1(\mathcal{P}-\Gamma,x_0)$ by deck transformations on $\clos{\mathcal{P}}-\clos{\Gamma}$. We have thus obtained that: 
\begin{prop}\label{p.kernelprojectiongroup}
   Let $\mathfrak{p}$ be the projection from $\clos{G}$ to $G$. We have that $$\text{ker}(\mathfrak{p})=\pi_1(\mathcal{P}-\Gamma,x_0)$$ In other words, $\clos{G}\big{/}\pi_1(\mathcal{P}-\Gamma,x_0) \cong G$.
    
    Moreover, the action of $\text{ker}(\mathfrak{p})$ on $\clos{\mathcal{P}}-\clos{\Gamma}$ can be identified with the action of $\pi_1(\mathcal{P}-\Gamma,x_0)$ by deck transformations on the universal cover of $\mathcal{P}-\Gamma$. In particular, for any $\clos{g}\in \text{ker}(\mathfrak{p})$, we have that $\clos{\rho}(\clos{g})=id$ or $\clos{\rho}(\clos{g})$ does not fix any point in $\clos{\mathcal{P}}-\clos{\Gamma}$.
\end{prop}

Finally, thanks to Proposition \ref{p.groupactiondefi}, 

\begin{prop}\label{p.actionequivariance}
    There exists a faithful action $\clos{\rho}:\clos{G}\rightarrow \text{Homeo}(\clos{\mathcal{P}})$ such that $$\clos{\pi}\circ \clos{\rho}(\clos{g})=\rho(\mathfrak{p}(\clos{g}))\circ \clos{\pi}$$ 
\end{prop} 

In conclusion, by combining the results of the two previous sections we get that, similarly to $\mathcal{P}$, the bifoliated plane of $\rho$ up to surgeries on $\Gamma$ can be naturally endowed with a pair of transverse singular foliations and a group action. As we will prove in Section \ref{s.theoremC}, contrary to $\mathcal{P}$, when $\rho$ is associated to a pseudo-Anosov flow, the previous singular foliations and group action do not describe the original flow up to orbital equivalence, but rather up to specific surgeries (see Theorem C).  In other words, once we consider the universal cover of $\mathcal{P}-\Gamma$, we lose track of how the original flow behaves around a finite number of periodic orbits, namely the periodic orbits associated to $\Gamma$. 

\subsection{On the properties of the action $\clos{\rho}$}

Following the notations of the previous section, consider $(\clos{\mathcal{P}},\clos{\mathcal{F}^s}, \clos{\mathcal{F}^u}, \clos{\rho})$ the bifoliated plane of $\rho$ up to surgeries on $\Gamma$ endowed with its stable and unstable foliations and group action. By construction $\clos{\mathcal{P}}$ is homeomorphic to the universal cover of $\mathcal{P}-\Gamma$ together with some points at infinity. It is therefore natural for the action $\clos{\rho}$ to inherit several properties from $\rho$ related to its strong Markovian action nature:

\begin{prop}\label{p.propertiespbaraction}
    \begin{enumerate}
        \item $\clos{\rho}$ preserves the pair of singular foliations $(\clos{\mathcal{F}^s},\clos{\mathcal{F}^u})$ 
        \item Two points $\clos{x},\clos{y}$ in $\clos{\mathcal{P}}$ are in the same $\clos{G}$-orbit if and only if the points $\clos{\pi}(\clos{x}), \clos{\pi}(\clos{y})$ are in the same $G$-orbit. 
        \item Let $\clos{g}$ be a non-trivial element in $ \clos{G}$. If $\clos{\rho}(\clos{g})$ preserves a leaf $\clos{L}\in \clos{\mathcal{F}^s}$ or $\clos{L}\in \clos{\mathcal{F}^u}$, then it has a unique fixed point in $\clos{L}$. Moreover, this fixed point is the same for every non-trivial element in $\clos{G}$ that preserves $\clos{L}$
        \item The stabilizer in $\clos{G}$ of any point $\clos{x}$ in $\clos{\mathcal{P}}-\clos{\Gamma}$ is either trivial or isomorphic to $\mathbb{Z}$. Moreover, $$\mathfrak{p}(\text{Stab}_{\clos{\rho}}(\clos{x}))=\text{Stab}_{\rho}(\clos{\pi}(\clos{x}))$$
        \item The stabilizer in $\clos{G}$ of any leaf $\clos{L}$ in $\clos{\mathcal{F}^s}$ or $\clos{\mathcal{F}^u}$ is either trivial or isomorphic to $\mathbb{Z}$. Moreover, $$\mathfrak{p}(\text{Stab}_{\clos{\rho}}(\clos{L}))=\text{Stab}_{\rho}(\clos{\pi}(\clos{L}))$$
        \item For any non-trivial $\clos{g}\in \clos{G}$ and any $\clos{x}\in \clos{\mathcal{P}}-\clos{\Gamma}$ (resp. $\clos{x}\in \clos{\Gamma}$) fixed by $\clos{\rho}(\clos{g})$, up to changing $\clos{g}$ to $\clos{g}^{-1}$, we have that $\clos{\rho}(\clos{g})$ is 
        topologically expanding on $\clos{\mathcal{F}^u}(\clos{x})$ (resp. any of the unstable leaves of $\clos{x}$) and topologically contracting on $\clos{\mathcal{F}^s}(\clos{x})$ (resp. any of the stable leaves of $\clos{x}$)
        \item There exists a finite set of stable (resp. unstable) leaves in $\clos{\mathcal{F}^s}$ (resp. $\clos{\mathcal{F}^u}$) with non-trivial stabilizers in $\clos{G}$, whose orbit by $\clos{G}$ forms a dense subset of $\clos{\mathcal{P}}$ 
        \end{enumerate} 

\end{prop}

\begin{proof}[Proof of (1)]
   Consider $\clos{L}$ a leaf of $\clos{\mathcal{F}^s}$ and $\clos{g}\in \clos{G}$. Assume for the sake of simplicity that $\clos{L}$ does not intersect $\clos{\Gamma}$ (the general case follows from a similar argument). First of all notice that, since $\clos{\rho}(\clos{g})$ is a homeomorphism of $\clos{\mathcal{P}}$, $\clos{\rho}(\clos{g})(\clos{L})$ is connected. Next, using Remark \ref{r.classificationofleavesoffolipbar}, Proposition \ref{p.actionequivariance} and our original hypothesis, we get that both $\clos{\pi}(\clos{L})$ and $\clos{\pi}(\clos{\rho}(\clos{g})(\clos{L}))$ are stable leaves in $\mathcal{F}^s$ that do not intersect $\Gamma$. Thanks to our discussion prior to Remark \ref{r.classificationofleavesoffolipbar}, every connected component of the preimage by $\clos{\pi}$ of $\clos{\pi}(\clos{\rho}(\clos{g})(\clos{L}))\in \mathcal{F}^s$ defines a leaf in $\clos{\mathcal{F}^s}$, which proves that $\clos{\rho}(\clos{g})(\clos{L})\in \clos{\mathcal{F}^s}$. 
\end{proof}
\begin{proof}[Proof of (2)]
    Indeed, if $\clos{x},\clos{y}$ are in the same $\clos{G}$-orbit, then their projections are in the same $G$-orbit thanks to Proposition \ref{p.actionequivariance}. Conversely, recall that by Proposition \ref{p.kernelprojectiongroup} for every $\gamma\in \mathcal{P}-\Gamma$ the group $\text{ker}(\mathfrak{p})$ acts transitively on $\clos{\pi}^{-1}(\gamma)$. It therefore suffices to show that the same result is true if $\gamma\in \Gamma$. Indeed, take $\gamma\in \Gamma$, $\clos{\gamma},\clos{\gamma'}\in \clos{\Gamma}$ any two lifts of $\gamma$ on $\clos{\mathcal{P}}$, $\clos{L},\clos{L'}$ two leaves of $\clos{\mathcal{F}^s}$ containing $\clos{\gamma},\clos{\gamma'}$ respectively and such that $\clos{\pi}(\clos{L})=\clos{\pi}(\clos{L'})=L$, where $L$ is a stable separatrix of $\gamma$. Thanks to Proposition \ref{p.kernelprojectiongroup}, there exists $\clos{g}\in \text{ker}(\mathfrak{p})$ such that $\clos{\rho}(\clos{g})(\clos{L})=\clos{L'}$. By Propositions \ref{p.actionequivariance} and \ref{p.singularitiesoffoliations}, we deduce that $\clos{\rho}(\clos{g})(\clos{\gamma})=\clos{\gamma'}$. 
\end{proof}
\begin{proof}[Proof of (3)]
    Consider $\clos{g}$ a non-trivial element in $ \clos{G}$ such that $\clos{\rho}(\clos{g})$ preserves a leaf $\clos{L}\in \clos{\mathcal{F}^s}$. Let us first show that $\mathfrak{p}(\clos{g})\neq e \in G$, where $e$ is the trivial element in $G$. Indeed, if $\mathfrak{p}(\clos{g})=e$, then thanks to Proposition \ref{p.actionequivariance} and to the injectivity of $\clos{\pi}$ on $\clos{L}$ (see Remark \ref{r.classificationofleavesoffolipbar}), we get that  $\clos{\rho}(\clos{g})$ fixes every point in $\clos{L}$. Using Proposition \ref{p.kernelprojectiongroup} and the faithfulness of the action $\clos{\rho}$, the previous fact implies that $\clos{\rho}(\clos{g})=id$ and thus that $\clos{g}$ is trivial in $\clos{G}$, which is absurd. 
    
    We can therefore suppose from now on that $\mathfrak{p}(\clos{g})\neq e$. In that case, since $\clos{\rho}(\clos{g})$ preserves $\clos{L}$, thanks to Proposition \ref{p.actionequivariance}, we get that $\rho(\mathfrak{p}(\clos{g}))\neq id$ preserves $\clos{\pi}(\clos{L})$. Recall now that by Remark \ref{r.classificationofleavesoffolipbar}, $\clos{\pi}(\clos{L})$ could either be a stable leaf of $\mathcal{F}^s$ that does not intersect $\Gamma$ or a stable separatrix of a point in $\Gamma$. In both cases, since $\rho(\mathfrak{p}(\clos{g}))$ preserves $\clos{\pi}(\clos{L})$ and since $\rho$ is a Markovian action, $\rho(\mathfrak{p}(\clos{g}))$ fixes a point in $\clos{\pi}(\clos{L})$. Using Proposition \ref{p.actionequivariance} and the injectivity of $\clos{\pi}$ on $\clos{L}$ (see Remark \ref{r.classificationofleavesoffolipbar}), we get that $\clos{\rho}(\clos{g})$ also fixes a point in $\clos{L}$. 

    By a similar argument, we can also prove  that any two elements in $\clos{G}$ that preserve $\clos{L}$ fix the same point in $\clos{L}$, which finishes the proof of Item (3). 

\end{proof}
\begin{proof}[Proof of (4)]
    Take $\clos{x}\in \clos{\mathcal{P}}-\clos{\Gamma}$. By Proposition \ref{p.kernelprojectiongroup}, for any $\clos{g}\in \text{ker}(\mathfrak{p})$ we have that $\clos{\rho}(\clos{g})=id$ or $\clos{\rho}(\clos{g})$ does not fix any point in $\clos{\mathcal{P}}-\clos{\Gamma}$. The previous fact together with Proposition \ref{p.actionequivariance} imply that $\mathfrak{p}$ defines a monomorphism from $\text{Stab}_{\clos{\rho}}(\clos{x})$ into $\text{Stab}_{\rho}(\clos{\pi}(\clos{x}))$. Therefore, $\text{Stab}_{\clos{\rho}}(\clos{x})$ can be seen as a subgroup of $\text{Stab}_{\rho}(\clos{\pi}(\clos{x}))$, which is either trivial or isomorphic to $\mathbb{Z}$, by the definition of a Markovian action. This proves that $\text{Stab}_{\clos{\rho}}(\clos{x})$ is either trivial or isomorphic to $\mathbb{Z}$. 
    
    It remains to prove that $\mathfrak{p}$ defines an epimorphism from $\text{Stab}_{\clos{\rho}}(\clos{x})$ to $\text{Stab}_{\rho}(\clos{\pi}(\clos{x}))$. Take $g\in \text{Stab}_{\rho}(\clos{\pi}(\clos{x}))$ and $\clos{g}$ such that $\mathfrak{p}(\clos{g})=g$. Thanks to Proposition \ref{p.actionequivariance}, we have that $\clos{\rho}(\clos{g})(\clos{x})\in \clos{\pi}^{-1}(\clos{\pi}(\clos{x}))$ and since $\text{ker}(\mathfrak{p})$ acts transitively on $\clos{\pi}^{-1}(\clos{\pi}(\clos{x}))$ (see Proposition \ref{p.kernelprojectiongroup}) we have that there exists $\clos{h}\in \text{ker}(\mathfrak{p})$ such that $\clos{\rho}(\clos{g}\clos{h})$ fixes $\clos{x}$. Since $\mathfrak{p}(\clos{g}\clos{h})=g$, this proves the desired result. 
\end{proof}

\begin{proof}[Proofs of (5),(6)]
    The proof of (4) follows from the same argument that was used in the proof of (3). The proof of (5) is a direct application of (3) and of Proposition \ref{p.actionequivariance}. 
\end{proof}
\begin{proof}[Proof of (7)]    
    Take $L_1,...,L_n$ a set of stable leaves in $\mathcal{F}^s$, whose orbits by $G$ form altogether a dense subset of $\mathcal{P}$. Assume for the sake of simplicity that $L_1,...,L_n$ do not intersect any point in $\Gamma$ (the general case follows from a similar argument). In that case, we claim that if $\clos{L_1},...,\clos{L_n}\in \clos{\mathcal{F}^s}$ are lifts of $L_1,...,L_n$ (see our discussion prior to Remark \ref{r.classificationofleavesoffolipbar}), then their orbits by $\clos{G}$ form a dense subset of $\clos{\mathcal{P}}$. 
    
    If the orbits by $\clos{G}$ of  $\clos{L_1},...,\clos{L_n}$ do not form a dense set in $\clos{\mathcal{P}}$, then there would exist $U$ an open set in  $\clos{\mathcal{P}}-\clos{\Gamma}$ that would not intersect the set $\clos{\rho}(\clos{G})(\underset{i\in\llbracket 1, n \rrbracket}{\cup}\clos{L_i})$. By Proposition \ref{p.actionequivariance} and thanks to Item (2), this would imply that the neighborhood $\clos{\pi}(U)\subset \mathcal{P}$ (this is indeed a neighborhood thanks for instance to our geometric construction of $\clos{\mathcal{P}}$) does not intersect $\rho(G)(\underset{i\in\llbracket 1, n \rrbracket}{\cup}L_i)$, which leads to a contradiction. 
\end{proof}

Thanks to the previous proposition, we can see that $\clos{\rho}$ shares several properties with $\rho$, as it satisfies most of the axioms defining a Markovian action. The main property that differentiates $\clos{\rho}$ from $\rho$ and any other Markovian action consists in the fact that the stabilizers in $\clos{G}$ of any point in $\clos{\Gamma}$ is isomorphic to a group that is bigger than $\mathbb{Z}$ :

\begin{prop}\label{p.stabilizersz2}
The stabilizer in $\clos{G}$ of any point in $\clos{\Gamma}$ is isomorphic to $\mathbb{Z}^2$. 
\end{prop}
\begin{proof}
Take $\clos{\gamma}\in \clos{\Gamma}$. By Proposition \ref{p.singularitiesoffoliations}, the set of stable leaves in $\clos{\mathcal{F}^s}$ intersecting  $\clos{\gamma}$ is countable and ordered along $\mathbb{Z}$. We will denote the previous set of leaves by  $...s_{-2},s_{-1},s_0,s_1,s_2...$ Since the elements of $\clos{G}$ preserve $\clos{\mathcal{F}^{s,u}}$, any element in $\text{Stab}(\clos{\gamma})$ permutes the set of leaves $\lbrace s_k|k\in \mathbb{Z}\rbrace$ and even more associates non-separated leaves to non-separated leaves. Therefore, by Proposition \ref{p.singularitiesoffoliations} any element $\clos{g}\in \text{Stab}(\clos{\gamma})$ acts on the set $\lbrace s_k|k\in \mathbb{Z}\rbrace \equiv \mathbb{Z}$ as the composition of some symmetry on $\mathbb{Z}$ and some translation. But since the action $\rho$ is orientation preserving, by Proposition \ref{p.actionequivariance}, $\clos{G}$ acts on $\clos{\mathcal{P}}$ by orientation preserving homeomorphisms too. Hence, $\clos{\rho}(\clos{g})$ can  only act as a translation on $\lbrace s_k|k\in \mathbb{Z}\rbrace$. 

Consider first the set of elements in $\text{Stab}(\clos{\gamma})$ acting on $\lbrace s_k|k\in \mathbb{Z}\rbrace$ as zero translations and thus the set of elements in $\clos{G}$ that fix every stable and unstable leaf of $\clos{\gamma}$. Thanks to Proposition \ref{p.propertiespbaraction}, the stabilizer in $\clos{G}$ of a leaf $s_i$ is isomorphic to $\mathbb{Z}$. It follows the elements of $\text{Stab}(\clos{\gamma})$ stabilizing every $s_i$ form a subgroup of $\clos{G}$ isomorphic to $\mathbb{Z}$. Take $z$ the generator of $\text{Stab}(s_0)$ that acts on $s_0$ as a topological expansion (see Proposition \ref{p.propertiespbaraction}).

Next, by Proposition \ref{p.singularitiesoffoliations}, there exist stable leaves of $\clos{\gamma}$ that project to the same separatrix in $\mathcal{P}$. It follows from Proposition \ref{p.propertiespbaraction}  that there exist elements in $\clos{G}$ acting as non-zero translations on $\lbrace s_k|k\in \mathbb{Z}\rbrace$. Take $t$ an element in $\clos{G}$ acting with the smallest possible strictly positive translation on $\lbrace s_k|k\in \mathbb{Z}\rbrace$.  Recall now that since $\clos{\mathcal{P}}-\clos{\Gamma}$ is the universal cover of $\mathcal{P}-\Gamma$, any element of $\clos{G}$ that fixes one point in $\clos{\mathcal{P}}-\clos{\Gamma}$ is the identity (see Proposition \ref{p.kernelprojectiongroup}). Hence, we have that $\text{Stab}(\clos{\gamma})=<z,t>$. It therefore suffices to show that $z$ and $t$ commute. 

Indeed, by the definition of $z$ and $t$ we have $tzt^{-1}\in \text{Stab}(s_0)$. Since $tzt^{-1}$ is a conjugate of $z$, it also acts as an expansion on $s_0$. Furthermore, by Proposition \ref{p.actionequivariance}, $\mathfrak{p}(t)$ and $\mathfrak{p}(z)$ fix $\clos{\pi}(\clos{\gamma})\in \mathcal{P}$ and therefore $\mathfrak{p}(tzt^{-1})=\mathfrak{p}(t)\mathfrak{p}(z)\mathfrak{p}(t)^{-1}=\mathfrak{p}(z)$, since $\text{Stab}(\clos{\pi}(\clos{\gamma}))=\mathbb{Z}$.

By our previous arguments, both $z$ and $tzt^{-1}$ preserve $s_0$ and $\mathfrak{p}(tzt^{-1}z^{-1})=e $, where $e$ is the trivial element in $G$. By Proposition \ref{p.propertiespbaraction}, it follows that $tzt^{-1}z^{-1}$ is equal to the trivial element in $\clos{G}$, which proves that $z$ and $t$ commute. Finally, it is clear that neither $z$ or  $t$ are of finite order. It follows that $<z,t>=\mathbb{Z}^2$.
\end{proof}

\section{Proof of Theorem C}\label{s.theoremC}

Our goal in this section consists in proving Theorem C. Contrary to Theorems A, D and E that admit natural generalizations for strong Markovian actions, Theorems B and C consider pseudo-Anosov flows up to surgery. Although I believe that it is technically possible to define a ``surgery" operation for a strong Markovian action and using this to generalize Theorems B and C for Markovian actions, this goes well beyond the scope of this paper. We will therefore restrict ourselves to proving Theorem C for strong Markovian actions arising from pseudo-Anosov flows.

\subsection{$\clos{\mathcal{P}}-\clos{\Gamma}$ as the bifoliated plane of a flow}\label{s.ptildebifoliated}
In the previous section, we defined, for any general strong Markovian action $\rho: G\rightarrow \text{Homeo}(\mathcal{P})$, the bifoliated plane $(\clos{\mathcal{P}}, \clos{\mathcal{F}^{s,u}},\clos{\rho})$ of $\rho$ up to surgeries on $\Gamma$ together with its stable/unstable foliations and natural group action by $\clos{G}$, the extension of $G$ associated to $\Gamma$. Our previous definition does not require $\rho$ to be associated with some pseudo-Anosov flow. However, when $\rho$ arises from some pseudo-Anosov flow $(M,\Phi)$, our definition does not provide a direct link between $(\clos{\mathcal{P}}, \clos{\mathcal{F}^{s,u}},\clos{\rho})$ and the dynamics of $\Phi$ inside $M$. Our goal in this section is to clarify the relation between $(M,\Phi)$ and $(\clos{\mathcal{P}}, \clos{\mathcal{F}^{s,u}},\clos{\rho})$ and to provide an alternative definition for $\clos{G}$ using the fundamental group of $M$.

Consider $M$ to be a closed, orientable, connected 3-manifold, $\Phi$ a pseudo-Anosov flow on $M$, $\mathcal{P}$ its bifoliated plane endowed with an orientation together with its pair of singular foliations $\mathcal{F}^s$ and $\mathcal{F}^u$ and its associated orientation-preserving action $\rho: \pi_1(M)\rightarrow \text{Homeo}(\mathcal{P})$. Consider $\mathcal{R}$ a Markovian family preserved by $\rho$ (recall that any Markovian family preserved by $\rho$ is strong thanks to Lemmas \ref{l.vertsubrectangleexists}, \ref{l.horizsubrectangleexists}, \ref{l.infiniteintersectionverticalrectangles} and \ref{l.infiniteintersectionhorizontalrectangles}). Denote by $\Gamma$ the boundary periodic points of $\mathcal{R}$, $\Gamma_M$ the finite set of periodic orbits of $\Phi$ associated to $\Gamma$, $\clos{\mathcal{P}}$ the bifoliated plane of $\rho$ up to surgeries on $\Gamma$ together with its pair of singular foliations $\clos{\mathcal{F}^s}$ and $\clos{\mathcal{F}^u}$ and its associated orientation-preserving action $\clos{\rho}: \clos{\pi_1(M)}\rightarrow \text{Homeo}(\mathcal{P})$, where $\clos{\pi_1(M)}$ is the extension of $\pi_1(M)$ associated to $\Gamma$.

\begin{prop}\label{p.actionpi1minusorbits}
Consider $(M-\Gamma_M,\Phi_{\Gamma_M})$ the restriction of $\Phi$ on $M-\Gamma_M$. We have that the universal cover of $M-\Gamma_M$ is homeomorphic to $\mathbb{R}^3$ and that the orbit space of the lift on $\mathbb{R}^3$ of the flow $(M-\Gamma_M,\Phi_{\Gamma_M})$ can be identified with the universal cover of $\mathcal{P}-\Gamma$ or equivalently with the space $\clos{\mathcal{P}}-\clos{\Gamma}$. The projections on $\clos{\mathcal{P}}-\clos{\Gamma}$ of the lifts on $\mathbb{R}^3$ of the stable and unstable leaves of $\Phi_{\Gamma_M}$ define on $\clos{\mathcal{P}}-\clos{\Gamma}$ a pair of transverse regular foliations. Finally, the action of $\pi_1(M-\Gamma_M)$ on $\mathbb{R}^3$ by deck transformations projects to a faithful action by homeomorphisms on $\clos{\mathcal{P}}-\clos{\Gamma}$. 
\end{prop}
\begin{proof}
Denote by $\widetilde{\Phi}$ the lift of $\Phi$ on $\widetilde{M}=\mathbb{R}^3$. Let us  first show that the universal cover of $M-\Gamma_M$ is homeomorphic to $\mathbb{R}^3$. Let $\widetilde{\Gamma_M}$ be the lift of $\Gamma_M$ on $\widetilde{M}$. Since $\widetilde{M}-\widetilde{\Gamma_M} $ is a covering space of $M-\Gamma_M$, the universal cover of $M-\Gamma_M$, denoted by $\widetilde{M-\Gamma_M}$, coincides with  the universal cover of $\widetilde{M}-\widetilde{\Gamma_M} $. The fact that $\widetilde{M-\Gamma_M}\cong \mathbb{R}^3$ follows from the fact that the universal cover of $\mathbb{R}^3$ minus a countable and transversely discrete set of properly embedded lines is homeomorphic to $\mathbb{R}^3$. 

Similarly, the lift $\widetilde{\Phi_{\Gamma_M}}$ of $\Phi_{\Gamma_M}$ on $\widetilde{M-\Gamma_M}$ can be identified with the lift of $\widetilde{\Phi}$ minus the orbits $\widetilde{\Gamma_M}$ on the universal cover of  $\widetilde{M}-\widetilde{\Gamma_M}$. This proves that the orbit space of $\widetilde{\Phi_{\Gamma_M}}$ can be identified with the universal cover of $\mathcal{P}-\Gamma$. Furthermore, since $\Gamma_M$ contains all the circle prong singularities of $\Phi$ (see Remark \ref{r.singularitiesareboundaries}), the restrictions of the stable and unstable foliations of $\Phi$ on $M-\Gamma_M$ define a pair of transverse regular foliations. Hence, their lifts on $\widetilde{M-\Gamma_M}$ project on the orbit space of $\widetilde{\Phi_{\Gamma_M}}$, thus endowing $\clos{\mathcal{P}}-\clos{\Gamma}$ with a pair of transverse regular foliations. 

Finally, the action of $\pi_1(M-\Gamma_M)$ on $\widetilde{M-\Gamma_M}$ by deck transfomrations preserves the orbits of $\widetilde{\Phi_{\Gamma_M}}$ and thus defines on its orbit space an action by homeomorphisms. Since $\Phi$ admits non-closed orbits in $M-\Gamma_M$, any element of $\pi_1(M-\Gamma_M)$ preserving every orbit of $\widetilde{\Phi_{\Gamma_M}}$ is necessarily trivial. This proves that the action of $\pi_1(M-\Gamma_M)$ on $\clos{\mathcal{P}}-\clos{\Gamma}$ is faithful and finishes the proof of the proposition.

\end{proof}

Similarly to our discussion after Theorem-Definition \ref{thdef.bifoliatedplane}, by classical cover space theory, $\pi_1(M-\Gamma_M)$ admits infinitely many actions on $\widetilde{M-\Gamma_M}$ by deck transformations. Furthermore, all the previous actions on $\widetilde{M-\Gamma_M}$ coincide up to performing an inner automorphism on $\pi_1(M-\Gamma_M)$. It follows that Proposition \ref{p.actionpi1minusorbits} defines infinitely many actions of $\pi_1(M-\Gamma_M)$ on $\clos{\mathcal{P}}-\clos{\Gamma}$ that are all the same up to performing an inner automorphism on $\pi_1(M-\Gamma_M)$. For the sake of simplicity and by an abuse of language, each of the previous actions of $\pi_1(M-\Gamma_M)$ on $\clos{\mathcal{P}}-\clos{\Gamma}$, will be called \emph{\textbf{the} action by $\pi_1(M-\Gamma_M)$ on $\clos{\mathcal{P}}-\clos{\Gamma}$ associated to $\Phi$}.

In addition to the action by $\pi_1(M-\Gamma_M)$ on $\clos{\mathcal{P}}-\clos{\Gamma}$ associated to $\Phi$, thanks to Propositions  \ref{p.kernelprojectiongroup} and \ref{p.actionequivariance}, we know that $\clos{\rho}:\clos{\pi_1(M)}\rightarrow \text{Homeo}(\clos{\mathcal{P}})$ defines another faithful action on $\clos{\mathcal{P}}-\clos{\Gamma}$ preserving a pair of transverse foliations. In the next proposition, we will show that the previous actions coincide: 
 
\begin{prop}\label{p.relationbetweenactions}
  The group $\pi_1(M-\Gamma_M)$ is isomorphic to $\clos{\pi_1(M)}$. Furthermore, the foliations in $\clos{\mathcal{P}}-\clos{\Gamma}$ that were defined in Proposition \ref{p.actionpi1minusorbits} coincide with the restrictions of $\clos{\mathcal{F}^s}$ and $\clos{\mathcal{F}^u}$ on $\clos{\mathcal{P}}-\clos{\Gamma}$. Finally, if $\clos{\rho}_{\pi_1(M-\Gamma_M)}$ denotes the action of $\pi_1(M-\Gamma_M)$ on $\clos{\mathcal{P}}-\clos{\Gamma}$ associated to $\Phi$, there exists $H$ an isomorphism from $\pi_1(M-\Gamma_M)$ to $\clos{\pi_1(M)}$ such that for every $g\in \pi_1(M-\Gamma_M)$  $$\clos{\rho}_{\pi_1(M-\Gamma_M)}(g)=\clos{\rho}(H(g))$$
  \end{prop}
\begin{proof}
Let $\widetilde{M}$ be the universal cover of $M$ and $\widetilde{\Gamma_M}$ the lift of $\Gamma_M$ on $\widetilde{M}$. Let us first show that 

\textbf{Claim.} $\pi_1(\widetilde{M}-\widetilde{\Gamma_M}) = \text{ker}(\mathfrak{i})$, where $\mathfrak{i}$ is the natural morphism from $\pi_1(M-\Gamma_M)$ to $\pi_1(M)$ given by the inclusion $M-\Gamma_M\hookrightarrow M$. 

\textit{Proof of the claim.} Indeed, consider $O$ a point in $\widetilde{M}-\widetilde{\Gamma_M}$ and its projection $O_M$ in $M-\Gamma_M$. Take $\gamma$ a non-homotopically trivial loop in $\widetilde{M}-\widetilde{\Gamma_M}$ starting and ending at $O$. Its projection $\gamma_M$  in $M-\Gamma_M$ is also non-homotopically trivial. However, by adding back the lines $\widetilde{\Gamma_M}$ in $\widetilde{M}-\widetilde{\Gamma_M}$ the loop $\gamma$ becomes homotopically trivial (as $\widetilde{M}$ is simply connected) and therefore by adding back the circles $\Gamma_M$ in $M-\Gamma_M$, the loop $\gamma_M$ becomes homotopically trivial. Hence, the map $\gamma \rightarrow \gamma_M$ defines a morphism $i: \pi_1(\widetilde{M}-\widetilde{\Gamma_M}, O)\rightarrow \text{ker}\big(\pi_1(M-\Gamma_M,O_M)\rightarrow \pi_1(M,O_M)\big)=ker(\mathfrak{i})$. 

Conversely, consider a loop $\gamma_M\in ker(\mathfrak{i})=\text{ker}\big(\pi_1(M-\Gamma_M,O_M)\rightarrow \pi_1(M,O_M)\big)$. Take the lift $\gamma$ of $\gamma_M$ that starts at $O$. By definition, if we add $\Gamma_M$ in $M-\Gamma_M$, then $\gamma_M$ becomes homotopically trivial. Therefore, if we add $\widetilde{\Gamma_M}$ in $\widetilde{M}-\widetilde{\Gamma_M}$ the curve $\gamma$ must also become homotopically trivial. We deduce that $\gamma$ is a loop that starts and ends at $O$. The map $\gamma_M \rightarrow \gamma$ defines a morphism $j:  ker(\mathfrak{i})\rightarrow \pi_1(\widetilde{M}-\widetilde{\Gamma_M}, O)$. It is easy to check that $i$ and $j$ are inverses of each other and therefore $\pi_1(\widetilde{M}-\widetilde{\Gamma_M}) = ker(\mathfrak{i})$. This finishes the proof of the claim.

\vspace{0.2cm}
From the previous claim we get that the cover $\widetilde{M}-\widetilde{\Gamma_M}$ of $M-\Gamma_M$ is a normal covering of $M-\Gamma_M$ and also that $$\text{Aut}(\widetilde{M}-\widetilde{\Gamma_M}\rightarrow M-\Gamma_M) = \quotient{\pi_1(M-\Gamma_M)}{ker(\mathfrak{i})} \cong \pi_1(M)$$

Notice that the restriction on $\widetilde{M}-\widetilde{\Gamma_M}$ of the action of $\text{Aut}(\widetilde{M}\rightarrow M)\cong \pi_1(M)$ coincides (up to inner automorphism) with the action of $\text{Aut}(\widetilde{M}-\widetilde{\Gamma_M}\rightarrow M-\Gamma_M)\cong \pi_1(M)$ on $\widetilde{M}-\widetilde{\Gamma_M}$. 

\vspace{0.5cm}
\noindent\textit{Construction of the group isomorphism $H$}.
\vspace{0.2cm}

First of all, let us recall the definition of $\clos{\pi_1(M)}$, the extension of $\pi_1(M)$ associated to $\Gamma$ (see Section \ref{s.groupactiondefi} for more details). Consider $x_0\in \mathcal{P}-\Gamma$ and $$\text{Curv}_2:= \lbrace \gamma:[0,1]\overset{C^0}{\rightarrow} \mathcal{P} | ~ \gamma(0)=x_0, \gamma[0,1)\cap \Gamma= \emptyset \rbrace $$ 
$$ \text{Curv}_2(x_0):=\{\gamma\in \text{Curv}_2|\gamma(1)\in \rho(\pi_1(M))(x_0)\}$$  
 Recall that every element of $\clos{\pi_1(M)}$ is a couple consisting of: 
\begin{enumerate}
    \item an element $h\in \pi_1(M)$
    \item the homotopy class relatively to $\sim_2$ (or equivalently relatively to its boundary) of a curve in $\text{Curv}_2(x_0)$  starting from $x_0$ and ending at $\rho(h)(x_0)$
\end{enumerate} 

Before beginning the construction of the desired isomorphism, let us introduce some notations. Denote by $\pi_{\widetilde{M}-\widetilde{\Gamma_M}}$ the projection map from $\widetilde{M-\Gamma_M}$ to $ \widetilde{M}-\widetilde{\Gamma_M}$ and by $\pi_{\mathcal{P}}$ the projection map from $\widetilde{M}$ to the orbit space $\mathcal{P}$ of $\widetilde{\Phi}$. 

By Theorem-Definition \ref{thdef.bifoliatedplane}, there exists $$\psi_{\widetilde{M}-\widetilde{\Gamma_M}}: \pi_1(M)\rightarrow \text{Homeo}(\widetilde{M}-\widetilde{\Gamma_M})$$ an action of $\text{Aut}(\widetilde{M}-\widetilde{\Gamma_M}\rightarrow M-\Gamma_M)\cong\pi_1(M)$ on $\widetilde{M}$ by deck transformations verifying for every $h\in \pi_1(M)$
\begin{equation}\label{eq.equivariancepsirho}
    \pi_{\mathcal{P}}\circ \psi_{\widetilde{M}-\widetilde{\Gamma_M}}(h)=\rho(h)\circ \pi_{\mathcal{P}}
\end{equation}

Similarly, by lifting $\psi_{\widetilde{M}-\widetilde{\Gamma_M}}$ on $\widetilde{M-\Gamma_M}$, we get that there exists $$\psi_{\widetilde{M-\Gamma_M}}: \pi_1(M-\Gamma_M)\rightarrow \text{Homeo}(\widetilde{M-\Gamma_M})$$ an action by deck transformations of $\text{Aut}(\widetilde{M-\Gamma_M}\rightarrow M-\Gamma_M)\cong \pi_1(M-\Gamma_M)$ on $\widetilde{M-\Gamma_M}$ verifying  for every $h\in \pi_1(M-\Gamma_M)$ 
\begin{equation}\label{eq.equivariancepsipsi}
    \pi_{\widetilde{M}-\widetilde{\Gamma_M}}\circ \psi_{\widetilde{M-\Gamma_M}}(h)=\psi_{\widetilde{M}-\widetilde{\Gamma_M}}(\mathfrak{i}(h))\circ \pi_{\widetilde{M}-\widetilde{\Gamma_M}}
\end{equation} where $\mathfrak{i}$ is the natural morphism from $\pi_1(M-\Gamma_M)$ to $\pi_1(M)$ given by the inclusion $M-\Gamma_M\hookrightarrow M$.
 \vspace{0.5cm}

 Fix $g\in\pi_1(M-\Gamma_M)$. We will now associate to $g$ an element of $\clos{\pi_1(M)}$. Take $X_0\in\widetilde{M}-\widetilde{\Gamma_M}$ such that $\pi_{\mathcal{P}}(X_0)=x_0$, $\widetilde{X_0}$ a lift of $X_0$ on $\widetilde{M-\Gamma_M}$, $\widetilde{\gamma_g}$ a continuous curve in $\widetilde{M-\Gamma_M}$ going from $\widetilde{X_0}$ to $\psi_{\widetilde{M-\Gamma_M}}(g)(\widetilde{X_0})$ and $\gamma_g:=\pi_{\widetilde{M}-\widetilde{\Gamma_M}}(\widetilde{\gamma_g})$. Notice that the homotopy class of $\gamma_g$ relatively to its boundaries depends only on $g$ and not on our choice of $\widetilde{\gamma_g}$. Also, by Equation \ref{eq.equivariancepsipsi}, $\gamma_g$ is a curve in $\widetilde{M}-\widetilde{\Gamma_M}$ starting from $X_0$ and ending at $\psi_{\widetilde{M}-\widetilde{\Gamma_M}}(\mathfrak{i}(g))(X_0)$. 

 We define $$H(g):=(\mathfrak{i}(g), \overline{\pi_{\mathcal{P}}(\gamma_g)})\in \clos{\pi_1(M)}$$ where $\overline{\pi_{\mathcal{P}}(\gamma_g})$ corresponds to the homotopy class of  $\pi_{\mathcal{P}}(\gamma_g)$ relatively to its boundaries in $\mathcal{P}-\Gamma$. The fact that $(\mathfrak{i}(g), \overline{\pi_{\mathcal{P}}(\gamma_g)})\in \clos{\pi_1(M)}$ is an immediate result of Equation \ref{eq.equivariancepsirho}. Finally, the proof of the fact that $H$ defines an isomorphism is left as an exercise to the reader. 

\vspace{0.5cm}
 \noindent{\textit{On the foliations defined in Proposition \ref{p.actionpi1minusorbits}}}
\vspace{0.2cm}

Denote by $\Phi_{\Gamma_M}$ the restriction of $\Phi$ on $M-\Gamma_M$ and by $\widetilde{\Phi_{\Gamma_M}}$ the lift of the previous flow on $\widetilde{M-\Gamma_M}$. Denote also by $\widetilde{\mathcal{F}^{s,u}}$ the lifts on $\widetilde{M}$ of the stable and unstable foliations of $\Phi$. Notice that the lifts on $\widetilde{M-\Gamma_M}$ of the stable and unstable foliations of $\Phi_{\Gamma_M}$ can be identified with the lifts on $\widetilde{M-\Gamma_M}$ of the restrictions of $\widetilde{\mathcal{F}^{s,u}}$ on $\widetilde{M}-\widetilde{\Gamma_M}$. Using the previous fact together with the fact that $\clos{\mathcal{P}}-\clos{\Gamma}$ is the universal cover of $\mathcal{P}-\Gamma$, we get that the foliations defined in Proposition \ref{p.actionpi1minusorbits}  correspond to the lifts on $\clos{\mathcal{P}}-\clos{\Gamma}$ of the stable and unstable foliations in $\mathcal{P}-\Gamma$ or equivalently, thanks to our construction in Section \ref{s.foliationsinpbarconstruction}, to the restrictions of $\clos{\mathcal{F}^{s}}$ and $\clos{\mathcal{F}^{u}}$ on $\clos{\mathcal{P}}-\clos{\Gamma}$. 

\vspace{0.9cm}
 \noindent{The action \textit{$\clos{\rho}$ coincides with the action of $\pi_1(M-\Gamma_M)$ on $\clos{\mathcal{P}}-\clos{\Gamma}$ associated to $\Phi$}}
\vspace{0.2cm}

 Let us first recall the definition of $\clos{\rho}$ given in Section  \ref{s.groupactiondefi}. Consider $\text{Curv}_1$ the set of continuous curves in $\mathcal{P}-\Gamma$ starting from $x_0$ and $\text{Curv}_1\big{/}\sim$ the set of homotopy classes relatively to their boundaries of the curves in  $\text{Curv}_1$. For every $\gamma\in \text{Curv}_1$ denote by $[\gamma]_1$ the class of $\gamma$ in $\text{Curv}_1\big{/}\sim$. In Section \ref{s.constructionPbar12}, we established that $\text{Curv}_1\big{/}\sim$ can be identified with $\clos{\mathcal{P}}-\clos{\Gamma}$. Furthermore, thanks to our discussion in Section \ref{s.groupactiondefi}, for every $g\in \pi_1(M-\Gamma_M)$ the action of $\clos{\rho}\big((\mathfrak{i}(g), \overline{\pi_{\mathcal{P}}(\gamma_g}))\big)= \clos{\rho}(H(g))$ can be seen as the operation in $\text{Curv}_{1}\big{/}\sim$ associating to every $[\gamma]_1\in \text{Curv}_{1}\big{/}\sim$, the homotopy class in  $\text{Curv}_{1}\big{/}\sim$ containing the curve obtained by the juxtaposition of $\pi_\mathcal{P}(\gamma_g)$ followed by $\rho(\mathfrak{i}(g))(\gamma)$. 

Similarly, let $\text{Curv}_{\widetilde{M}-\widetilde{\Gamma_M}}$ be the set of of continuous curves in $\widetilde{M}-\widetilde{\Gamma_M}$ starting at $X_0$ and $\text{Curv}_{\widetilde{M}-\widetilde{\Gamma_M}}\big{/}\sim$ be the set of homotopy classes relatively to their boundaries of the curves in  $\text{Curv}_{\widetilde{M}-\widetilde{\Gamma_M}}$. For every curve $\gamma$ in $\text{Curv}_{\widetilde{M}-\widetilde{\Gamma_M}}$  denote by $[\gamma]$ its homotopy class in $\text{Curv}_{\widetilde{M}-\widetilde{\Gamma_M}}\big{/}\sim$. Recall that, by classical covering space theory, $\widetilde{M-\Gamma_M}$ can be identified with $\text{Curv}_{\widetilde{M}-\widetilde{\Gamma_M}}\big{/}\sim$. Furthermore, by eventually performing an inner automorphism on $\pi_1(M-\Gamma_M)$, for every $g\in\pi_1(M-\Gamma_M)$ the deck transformation  $\psi_{\widetilde{M-\Gamma_M}}(g)$ can be seen as the operation in $\text{Curv}_{\widetilde{M}-\widetilde{\Gamma_M}}\big{/}\sim$ associating to every $[\gamma]\in \text{Curv}_{\widetilde{M}-\widetilde{\Gamma_M}}\big{/}\sim$, the homotopy class in $\text{Curv}_{\widetilde{M}-\widetilde{\Gamma_M}}\big{/}\sim$ containing the curve obtained by the juxtaposition of $\gamma_g$ followed by $\psi_{\widetilde{M}-\widetilde{\Gamma_M}}(\mathfrak{i}(g))(\gamma)$.

The desired result is an immediate consequence of our previous arguments and of Equation \ref{eq.equivariancepsirho}.

\end{proof}

Take $\clos{\rho}_{\pi_1(M-\Gamma_M)}$ the action of $\pi_1(M-\Gamma_M)$ on $\clos{\mathcal{P}}-\clos{\Gamma}$ associated to $\Phi$. Thanks to the previous proposition, $\clos{\rho}_{\pi_1(M-\Gamma_M)}$ coincides with the restriction of $\clos{\rho}$ on $\clos{\mathcal{P}}-\clos{\Gamma}$ and thus can be continuously extended on $\clos{\mathcal{P}}$. We will name the previous action of $\pi_1(M-\Gamma_M)$ on $\clos{\mathcal{P}}$, which coincides with $\clos{\rho}$, \emph{the action of $\pi_1(M-\Gamma_M)$ on $\clos{\mathcal{P}}$ associated to $\Phi$}. 

\subsection{On the elements in $\pi_1(M-\Gamma_M)$ fixing a point in $\clos{\Gamma}$}\label{s.elementsfixingpointsingammabar}
\quad{~}

Consider $M$ to be a closed, orientable, connected 3-manifold, $\Phi$ a pseudo-Anosov flow on $M$, $\mathcal{P}$ its bifoliated plane endowed with an orientation together with its pair of singular foliations $\mathcal{F}^s$ and $\mathcal{F}^u$ and its associated orientation-preserving action $\rho: \pi_1(M)\rightarrow \text{Homeo}(\mathcal{P})$. Consider next $\mathcal{R}$ a Markovian family preserved by $\rho$. Denote by $\Gamma$ the boundary periodic points of $\mathcal{R}$, $\Gamma_M$ the finite set of periodic orbits of $\Phi$ associated to $\Gamma$, $\clos{\mathcal{P}}$ the bifoliated plane of $\rho$ up to surgeries on $\Gamma$ together with its pair of singular foliations $\clos{\mathcal{F}^s}$ and $\clos{\mathcal{F}^u}$ and its associated orientation-preserving action $\clos{\rho}: \pi_1(M-\Gamma_M)\rightarrow \text{Homeo}(\mathcal{P})$. Denote by $\clos{\pi}$ the projection map from $\clos{\mathcal{P}}$ to $\mathcal{P}$ and by $\clos{\Gamma}$ the preimage of $\Gamma$ by $\clos{\pi}$.

Thanks to Proposition \ref{p.relationbetweenactions}, the action  $\clos{\rho}$ satisfies the results of Propositions \ref{p.kernelprojectiongroup}, \ref{p.actionequivariance}, \ref{p.propertiespbaraction} and \ref{p.stabilizersz2}. In particular, by Propositions  \ref{p.propertiespbaraction} and \ref{p.stabilizersz2} we have that the stabilizer in $\pi_1(M-\Gamma)$ of any point in $\clos{\mathcal{P}}$ is either trivial, isomorphic to $\mathbb{Z}$ or isomorphic to $\mathbb{Z}^2$. Next, thanks to  Propositions  \ref{p.periodicinbifoliated} and \ref{p.actionpi1minusorbits}, exactly as for the bifoliated plane of a pseudo-Anosov flow,  
\begin{itemize}
    \item the orbits of $\Phi$ that are not periodic correspond to points in $\clos{\mathcal{P}}$, whose stabilizers in $\pi_1(M-\Gamma)$ are trivial 
    \item the periodic orbits of $\Phi$ that do not belong in $\Gamma_M$ correspond to points in $\clos{\mathcal{P}}$, whose stabilizers in $\pi_1(M-\Gamma)$ are isomorphic to $\mathbb{Z}$
    \item the periodic orbits of $\Gamma_M$ correspond to points in $\clos{\Gamma}$, whose stabilizers in $\pi_1(M-\Gamma)$ are isomorphic to $\mathbb{Z}^2$
\end{itemize} 
In the second case, we already know that the stabilizer (up to conjugacy) of the point is generated by the homotopy class of its associated periodic orbit in $M$. Our goal in this section, consists in describing the homotopy classes in $\pi_1(M-\Gamma)$ forming the $\mathbb{Z}^2$ stabillizer of a point in $\clos{\mathcal{P}}$. 

\vspace{0.5cm}
Denote by $\widetilde{\Phi}$ the lift of $\Phi$ on the universal cover of $M$. Consider an orientation on $M$. By lifting the previous orientation on $\widetilde{M}$ and using the direction of the flow $\widetilde{\Phi}$, we can canonically orient any local transverse section of $\widetilde{\Phi}$. It follows that our choice of orientation on $M$ defines an orientation on the orbit space $\mathcal{P}$. Even more, since $\clos{\mathcal{P}}-\clos{\Gamma}$ is the universal cover of $\mathcal{P}-\Gamma$, the previous orientation on $\mathcal{P}$ also defines an orientation on $\clos{\mathcal{P}}-\clos{\Gamma}$. 

Fix for the rest of this section $\clos{\gamma}\in\clos{\Gamma}$. Let $...,s_{-1},s_0,s_1,s_2,...$ be the stable leaves in $\clos{\mathcal{F}^s}$ containing $\clos{\gamma}$. By eventually reindexing the $s_i$, we will assume without any loss of generality that for every $i\in \mathbb{Z}$ the leaf $s_{i+1}$ is not separated from $s_i$ (see Proposition \ref{p.singularitiesoffoliations}) and that for any point $x\in s_i-\{\clos{\gamma}\}$, the direction at $x$ pointing inside the connected component of $\clos{\mathcal{P}}-s_i$ containing $s_{i+1}$ (see Remark \ref{r.propertiesfolipbar}) followed by the direction at $x$ pointing inside $s_i$ and away from $\clos{\gamma}$ defines an orientation compatible with our choice of orientation of $\clos{\mathcal{P}}-\clos{\Gamma}$.\footnote{Formally speaking, as the $s_i$ are all $C^0$, this sentence makes little sense. However, one can make sense of this by using a foliation chart.}   

In Proposition \ref{p.stabilizersz2}, using the fact that $\clos{\rho}$ is orientation preserving, we showed that the stabilizer of $\clos{\gamma}$ is equal to $<z,t>\cong \mathbb{Z}^2$, where 

\begin{itemize}
    \item $z$ in the unique element in $\pi_1(M-\Gamma)$ preserving every $s_i$, generating $\text{Stab}(s_i)$ and acting as an expansion on $s_i$ for every $i\in \mathbb{Z}$
    \item $t$ is any element in $\pi_1(M-\Gamma)$ taking for every $i\in \mathbb{Z}$ the leaf $s_i$ to $s_{l+i}$, where $l$ is the smallest positive integer such that $s_0 $ and $s_{l}$ belong in the same $\clos{\rho}$-orbit
\end{itemize}

Notice that $t$ is uniquely defined up to a multiplication by a power of $z$. We would now like to define a different basis of $\text{Stab}(\clos{\gamma})$, whose elements will be uniquely defined by their actions on $\clos{\mathcal{P}}$.

\vspace{0.5cm}
\textit{A different basis of $\text{Stab}(\clos{\gamma})$.} Consider  $\gamma:=\clos{\pi}(\clos{\gamma})$, $n$ the smallest positive integer such that $\clos{\pi}(s_n)=\clos{\pi}(s_0)$ ($n$ corresponds to the number of stable prongs of $\gamma$) and $\mathfrak{p}$ the natural morphism from $\pi_1(M-\Gamma)$ to $\pi_1(M)$. Recall that $\mathfrak{p}$ is equivariant with respect to $\clos{\rho}$ and $\rho$ (see Proposition \ref{p.actionequivariance}). We now define  

\begin{itemize}
    \item $m$ as the unique element in $\text{ker}(\mathfrak{p})\cap \leq \pi_1(M-\Gamma)$ taking $s_0$ to $s_n$ 
    \item $p$ as the unique element in $\pi_1(M-\Gamma)$ taking $s_0$ to a stable leaf in $\{s_1,...,s_{n-1}, s_n\}$ and such that $\mathfrak{p}(p)=g$, where $g$ is the generator of $\text{Stab}(\gamma)$ in $\pi_1(M)$ acting as an expansion on $\mathcal{F}^s(\gamma)$
\end{itemize}
Both $m$ and $p$ are well defined, belong in $\text{Stab}(\clos{\gamma})$  and are unique, thanks to Propositions  \ref{p.singularitiesoffoliations}, \ref{p.kernelprojectiongroup},   \ref{p.actionequivariance}, \ref{p.propertiespbaraction}. Furthermore, 
\begin{prop}\label{p.mpformbasis}
    We have that $<m,p>=<z,t>=\text{Stab}(\clos{\gamma})$
\end{prop}
\begin{proof}
    
Indeed, by our definition of $p$ and $m$, we have that $<p,m>\leq \text{Stab}(\clos{\gamma})$. Conversely, using Proposition \ref{p.periodicinbifoliated} in addition to the fact that the action of $\pi_1(M)$ on $\mathcal{P}$ is orientation preserving and that $n$ corresponds to the number of prongs of $\gamma$, there exists $r\in \llbracket 0, n-1 \rrbracket$ such that the action of $g$ close to $\gamma$ is conjugated to $\phi_{n,r}$, a local model for a pseudo-hyperbolic fixed point with $n$ prongs, rotation $r$ and positive orientation. Using the previous fact, we get that $$ \text{Stab}_{\pi_1(M)}(\clos{\pi}(s_0))=<g^{{\frac{n}{\text{gcd(n,r)}}}}>$$ where $\text{gcd(n,r)}$ is the greatest common divisor of $n$ and $r$ (by convention $\text{gcd}(n,0)=n$). Since both $g^{{\frac{n}{\text{gcd(n,r)}}}}\in \pi_1(M)$ and $\mathfrak{p}(z)\in \pi_1(M)$ generate the group $\text{Stab}_{\pi_1(M)}(\clos{\pi}(s_0))\cong \mathbb{Z}$ and act as expansions on $\mathcal{F}^s(\gamma)$, we get that $$\mathfrak{p}(z)= g^{{\frac{n}{\text{gcd(n,r)}}}}$$. It follows that $z\cdot p^{-\frac{n}{\text{gcd(n,r)}}}$ is equal to a power of $m$ and thus $z\in <m,p>$. Next, since $g$ generates $\text{Stab}(\gamma)$, it is not difficult to see that the group $<p,m>$ acts transitively on the set of stable leaves $s_i$ that belong in the $\clos{\rho}$-orbit of $s_0$. Using the previous fact, the fact that $t$ is only defined up to a power of $z$ and that $z\in <m,p>$, we get that $t\in <m,p>$, which gives us the desired result. 
\end{proof}

Recall that our goal in this section consists in describing the homotopy classes in $M-\Gamma_M$ belonging to $\text{Stab}(\clos{\gamma})$. In order to do so, it suffices to describe the homotopy classes in $M-\Gamma_M$ corresponding to $m$ and $p$. 

\vspace{0.5cm}
\textit{On the description of the homotopy classes $m$ and $p$.} Consider $(M^*, \Phi^*)$ a flow obtained after blowing-up the orbits of $\Gamma_M$ (see our construction in Section \ref{s.blowupflows}). Notice that $\pi_1(M^*)\cong\pi_1(M-\Gamma_M)$ and that our choice of orientation of $M$ defines an orientation on $M^*$. 

The point $\clos{\gamma}\in \clos{\mathcal{P}}$ corresponds to a periodic orbit $\gamma_M$ in $\Gamma_M$ that was blown-up to a torus $T$ in $M^*$. In Section \ref{s.surgeryaspairsof integers}, thanks to our choice of orientation of $M^*$, we defined a meridian and a parallel in $T$, that we will denote by $m_\gamma$ and $p_\gamma$, forming a basis of $\pi_1(T)\leq \pi_1(M^*)\cong\pi_1(M-\Gamma_M)$. After recalling the definition of $m_\gamma$ and $p_\gamma$, we will show that (up to conjugation) $m_\gamma$ and $p_\gamma$ coincide respectively with $m$ and $p$. 

\vspace{0.5cm}
\textit{Definition of $m_\gamma$.}  Using our choice of orientation on $M^*$, endow $T$ with its inward orientation. Next, denote by $i_g(\gamma_1,\gamma_2)$ (resp. $i_a(\gamma_1,\gamma_2)$) the geometric (resp. algebraic) number of intersections of two curves (resp. curves or homotopy classes) in $T$. We remark here that $i_a$ is well defined thanks to our choice of orientation of $T$.

Consider now any periodic orbit $\gamma_P$ of $\Phi^*$ on $T$ endowed with its dynamical orientation and denote by $P_\gamma\in \pi_1(T)$ the homotopy class of $\gamma_P$ in $T$ (this class does not depend on our choice of periodic orbit on $T$). Let $\Pi^*_M: M^* \rightarrow M$ be the blow-down map associated to $(M^*,\Phi^*)$. The restriction of $\Pi^*_M$ on $T$ defines a fibration by circles topologically transverse to $\Phi^*$ over the circle $\gamma_M$. Take $\gamma_m$ one of the previous circles and choose an orientation on $\gamma_m$ so that $i_a(\gamma_P,\gamma_m)> 0$. Denote by $m_\gamma$ the homotopy class in $T$ of $\gamma_m$ endowed with the previous orientation. Recall that we called $m_\gamma\in \pi_1(T)$ the meridian of $T$ and that $m_\gamma$ depends only on our initial choice of orientation of $M^*$.

\begin{figure}[h]
    \centering
    \includegraphics[scale=0.2]{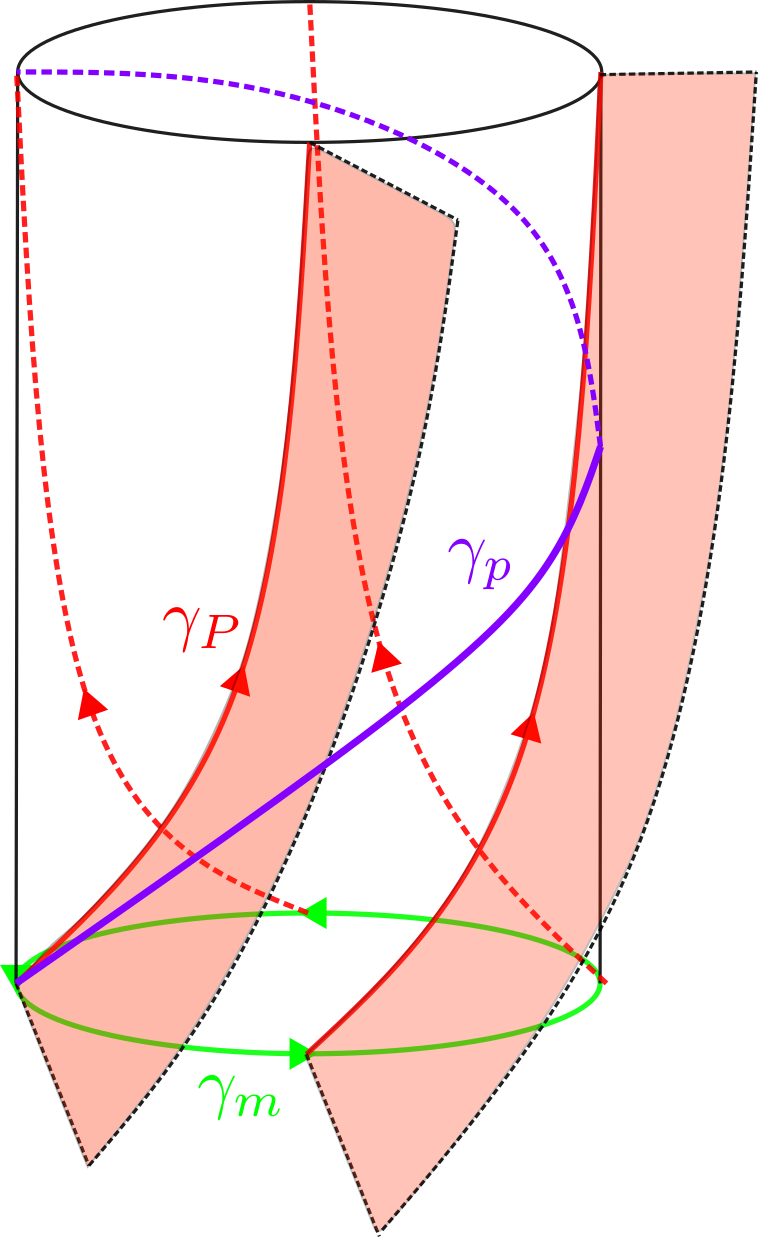}
    \caption{An example of blow-up of a periodic orbit of $\Phi$  with $4$ prongs and rotation $1$}
    \label{f.figurecurves}
\end{figure}

\vspace{0.5cm}
\textit{Definition of $p_\gamma$.} 
Recall that $\gamma_m$ and $\gamma_P$ are simple closed curves, that intersect transversely and that verify  $$i_a(\gamma_P,\gamma_m)=i_g(\gamma_P,\gamma_m)>0$$
In other words, $\gamma_m$ and $\gamma_P$ intersect at least once and performing a homotopy on either $\gamma_m$ or  $\gamma_P$ can only increase the number of intersections between the two curves.
More precisely, if $\gamma_M$ is associated by Proposition \ref{p.aroundcircleprong} to a local model for a pseudo-hyperbolic orbit with $n\geq 2$ prongs and rotation $r \in \llbracket 0,p-1\rrbracket$, then thanks to Remark \ref{r.conjugation2}, we get that the total number of intersections between $\gamma_P$ and $\gamma_m$ is $\frac{n}{gcd(n,r)}$. 

 Consider now $x_0,...,x_{n-1}$ the points of intersection of $\gamma_m$ with all the periodic orbits of $\Phi^*$ on $T$. By eventually reindexing the $x_i$, we can assume that $x:=x_0$ is a point of intersection between $\gamma_P$ and $\gamma_m$ and that when following positively $\gamma_m$ starting from $x$, we encounter first $x_1$, then $x_2$,..., and so on till we return to $x$. Following positively $\gamma_P$ starting from $x$, let $x_k$ be the first point of intersection between $\gamma_P$ and $\gamma_m$ after $x$, where $k\in \llbracket 0, n-1\rrbracket$. Let us remark here that 
\begin{itemize}
    \item $k=r$ or $k=n-r$ depending on our choice of orientation of $M^*$
    \item following positively $\gamma_P$ starting from $x_i$, the point $x_{i+k}$ (the indexes are considered modulo $n$) is the first point of intersection between $\gamma_P$ and $\gamma_m$ after $x_i$
\end{itemize}

Consider $y:=x_{n-k}$, $[y,x]_P$ the segment contained in $\gamma_P$ going from $y$ to $x$ (following $\gamma_P$ positively) and  $[x,y]_m$ the segment  contained in $\gamma_m$ going from $x$ to $y$ (following $\gamma_m$ positively). By convention, if $x=y$ (and thus $k=0$), we take $[x,y]_m=\gamma_m$ and $[y,x]_P=\gamma_P$. The juxtaposition of $[x,y]_m$ followed by $[y,x]_P$ defines a closed curve, whose homotopy class in $T$ will be denoted by $p_\gamma\in \pi_1(T)$. Notice that by performing a homotopy on the previous curve inside $T$, we can obtain a curve $\gamma_p$ intersecting $\gamma_m$ only once (see Figure \ref{f.parallelmeridian}). It follows that $m_\gamma$ and $p_\gamma$ form a basis of $\pi_1(T)$. Recall that we called $p_\gamma$ the parallel of $T$ and that $p_\gamma$ depends only on our choice of orientation on $M^*$.  

\vspace{0.7cm}

We are now ready to state the main result of this section. Recall that in the beginning of this section, we defined $\clos{\gamma}\in \clos{\Gamma}$, $(s_i)_{i\in \mathbb{Z}}$ the set of ordered stable manifolds of $\clos{\gamma}$ in $\clos{\mathcal{P}}$,  $\gamma=\clos{\pi}(\clos{\gamma})$, $n$ the number of stable prongs of $\gamma$ and $\gamma_M$ the periodic orbit of $\Phi$ corresponding to $\clos{\gamma}$. 

\begin{figure}
    \centering
    \includegraphics[scale=0.17]{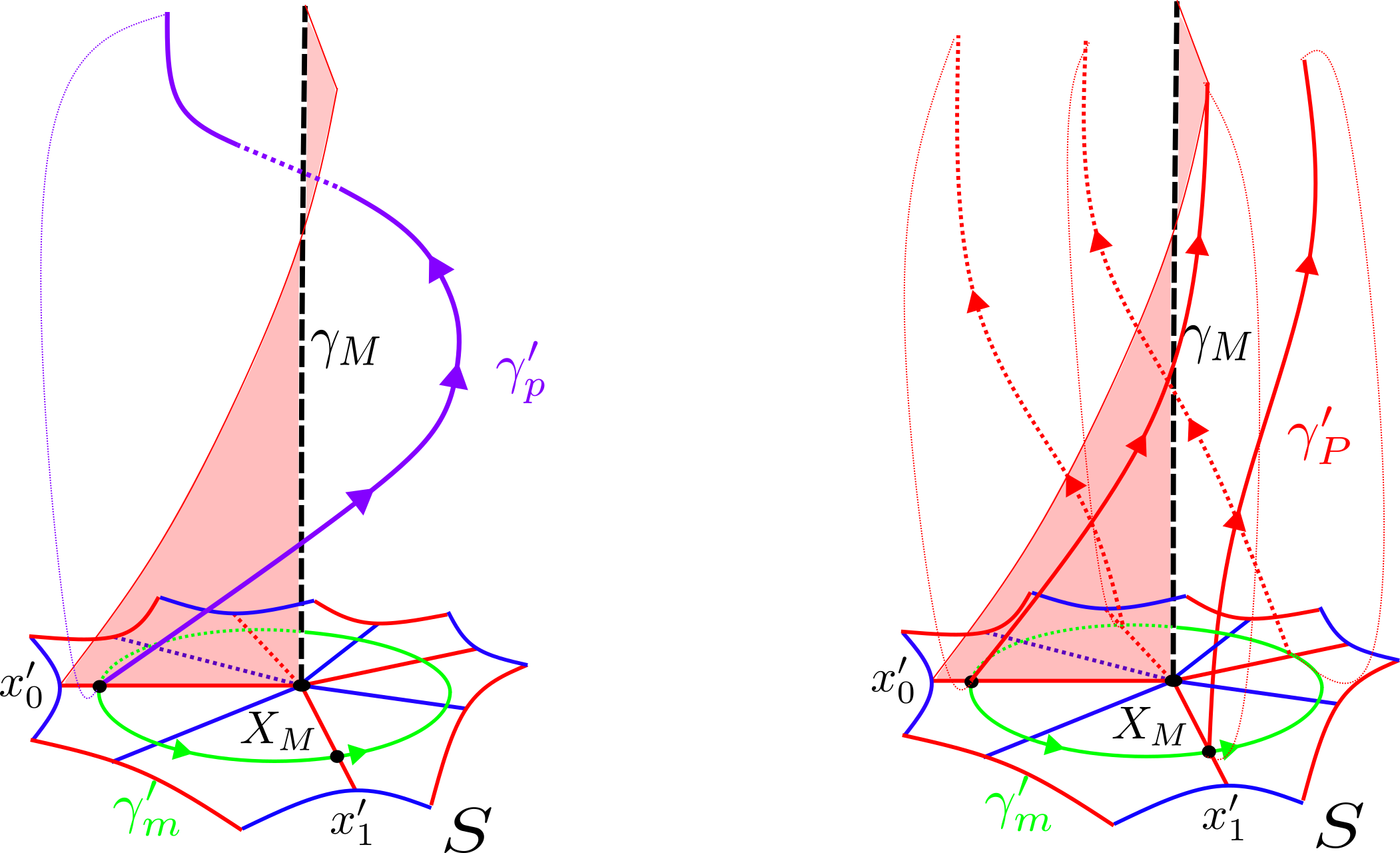}
    \caption{}
    \label{f.curvesafterblowup}
\end{figure}

Identify from now on the interior of $M^*$ with $M-\Gamma_M$. Push by a homotopy $\gamma_P$, $\gamma_m$ and $\gamma_p$ away from $T\subset \partial M^*$ in order to obtain three new curves  $\gamma'_P$, $\gamma'_m$ and $\gamma'_p$ in $M-\Gamma_M$ (see Figure \ref{f.curvesafterblowup}). By eventually changing our choice of  $\gamma'_P$, $\gamma'_m$, $\gamma'_p$ assume without any loss of generality that 
\begin{enumerate}

    \item there exists a transverse standard polygon $S$ for $\Phi$ intersecting $\gamma_M$ at $X_M$ and such that $\gamma'_m\subset S$
    \item $\gamma'_m$ intersects transversely and at a unique point any stable or unstable separatrix of $X_M$ in $S$. Denote by $x'_0,...,x_{n-1}'$ the successive intersections of $\gamma_m'$ with the stable leaf of $X_M$ in $S$.
    \item $\gamma'_P$, $\gamma'_m$, $\gamma'_p$ intersect at $x_0'\in M-\Gamma_M$
\end{enumerate}
By a small abuse of language, denote respectively by $P_\gamma$, $m_\gamma$ and $p_\gamma$ the homotopy classes of $\gamma'_P$, $\gamma'_m$ and  $\gamma'_p$ in $\pi_1(M-\Gamma_M,x_0')$. 

\begin{prop}\label{p.secondproofofcommutation}
We have that there exists $h\in \pi_1(M-\Gamma_M,x_0')$ such that $m_\gamma=h\cdot m \cdot h^{-1}$, $p_\gamma=h \cdot p  \cdot h^{-1}$ and $P_\gamma=h \cdot z  \cdot h^{-1}$. In particular, for every $i\in \mathbb{Z}$   

\begin{itemize}
    \item $\clos{\rho}(m_{\gamma})(s_i)=s_{i+n}$ 
    \item  $\clos{\rho}(p_{\gamma})(s_i)=s_{i+n-k}$ 
    \item  $\clos{\rho}(P_{\gamma})(s_i)=s_{i}$ 
\end{itemize}

\end{prop}
\begin{proof}
Thanks to the fact that $p,m, z\in \pi_1(M-\Gamma_M)$ are uniquely defined by their action on $\clos{\mathcal{P}}$, in order to prove the above proposition, it suffices to understand how the elements $p_\gamma,m_\gamma, P_\gamma$ act on $\clos{\mathcal{P}}$ via $\clos{\rho}$. First of all, notice that the three curves $\gamma'_P$, $\gamma'_m$ and $\gamma'_p$ can be moved inside $\gamma_M$ by a homotopy in $M$. It follows, thanks to Proposition \ref{p.relationbetweenactions}, that $\clos{\rho}(m_\gamma), \clos{\rho}(p_\gamma)$ and $\clos{\rho}(P_\gamma)$ fix a point in $\clos{\Gamma}$ that is in the same $\clos{\rho}$-orbit as $\clos{\gamma}$. By eventually conjugating $P_\gamma$, $m_\gamma$ and $p_\gamma$ by an element in $\pi_1(M-\Gamma_M,x_0')$ assume without any loss of generality that  $\clos{\rho}(m_\gamma), \clos{\rho}(p_\gamma)$ and $\clos{\rho}(P_\gamma)$ fix $\clos{\gamma}$. Using our previous hypothesis, we will now show that $m_\gamma=m$, $p_\gamma=p$ and $P_\gamma=z$.

First, notice that if we add $\Gamma_M$ back to $M-\Gamma_M$ the curve $\gamma'_m$ becomes homotopically trivial. Hence, $m_\gamma\in \text{ker}(\mathfrak{p})$, where $\mathfrak{p}$ is the natural morphism from $\pi_1(M-\Gamma_M,x'_0)$ to $\pi_1(M,x'_0)$. Also, take $L$ to be any stable separatrix of $X_M$ in $S$ and $Y$ the unique point of intersection of $\gamma'_m$ and $L$. Thanks to our choice of orientation of $\gamma_m$, we have that the direction of the flow $\Phi$ at $Y$ followed by the direction of $\gamma_m'$ at $Y$,  followed finally by the direction in $L$ pointing away from $X_M$ defines an orientation on $M$ compatible with our original choice of orientation\footnote{Once again, this sentence makes little sence as we are working in a $C^0$ setting. However, one can make sense of this by using a well chosen chart}. Thanks to the previous fact and the fact that $\gamma'_m$ goes once around $\gamma_M$, we get that  $\clos{\rho}(m_\gamma)(s_i)=s_{n+i}$ for every $i\in \mathbb{Z}$; hence, $m_\gamma=m$.

Next, if we add $\Gamma_M$ back to $M-\Gamma_M$ the curve $\gamma'_p$ becomes homotopic to $\gamma_M$, which implies thanks to Proposition \ref{p.periodicinbifoliated} that $\mathfrak{p}(p_\gamma)$ is a generator of $\text{Stab}_{\pi_1(M)}(\gamma)$, where $\gamma=\clos{\pi}(\clos{\gamma})$. More specifically, it is the generator of $\text{Stab}(\gamma)$ acting on $\mathcal{F}^s(\gamma)$, the stable leaf of $\gamma$ in $\mathcal{P}$, as an expansion. Recall that $\gamma'_p$ is homotopic to a curve consisting first of a segment in $\gamma'_m$ going from $x'_0$ to $x'_{n-k}$ followed by a segment inside a stable separatrix of $\gamma_M$. By looking at the action of $p_\gamma$ as a deck transformation on the universal cover of $M-\Gamma_M$, we get that $\clos{\rho}(p_\gamma)(s_i)=s_{i+n-k}$ for every $i\in \mathbb{Z}$; hence, $p_\gamma=p$.

Finally, if we add $\Gamma_M$ back to $M-\Gamma_M$ the curve $\gamma_P$ becomes homotopic to $\gamma_M$ followed positively $\frac{n}{\text{gcd}(n,k)}$ times. Thanks to the previous fact, the fact that $\gamma_M$ is associated by Proposition \ref{p.aroundcircleprong} to a local model for a pseudo-hyperbolic orbit with $n\geq 2$ prongs and rotation $r \in \llbracket 0,p-1\rrbracket$ and since $k\in \{r, n-r\}$, we get that $\mathfrak{p}(P_\gamma)=\mathfrak{p}(p_\gamma)^{\frac{n}{\text{gcd}(n,k)}}$ is the generator of $\text{Stab}_{\pi_1(M)}(\clos{\pi}(s_0))$ acting on $\clos{\pi}(s_0)$ as an expansion. Recall that $\gamma'_P$ is homotopic to a loop inside a stable separatrix of $\gamma_M$. By looking at the action of $P_\gamma$ as a deck transformation on the universal cover of $M-\Gamma_M$, we get that $\clos{\rho}(P_\gamma)(s_i)=s_i$  for every $i\in \mathbb{Z}$; hence, $P_\gamma=z$.

\end{proof}

\subsection{Proof of Theorem C}\label{s.barbotgeneral}
Let $(\Phi_1, M_1)$, $(\Phi_2, M_2)$ be  two pseudo-Anosov flows on two closed, orientable 3-manifolds, $(\mathcal{P}_1,\mathcal{F}_1^{s,u},\rho_1), (\mathcal{P}_2,\mathcal{F}_2^{s,u},\rho_2)$ their bifoliated planes endowed with their natural foliations and group actions, $\mathcal{R}_1$, $\mathcal{R}_2$ two Markovian families preserved by $\rho_1$ and $\rho_2$ and $\Gamma_1,\Gamma_2$ their boundary periodic points. Denote by $\Gamma^{M_1}_1$, $\Gamma^{M_2}_2$ the two finite sets of periodic orbits of $\Phi_1$,  $\Phi_2$ corresponding to $\Gamma_1$ and $\Gamma_2$. Denote also by $(\clos{\mathcal{P}_1}, \clos{\mathcal{F}_1^{s,u}}, \clos{\rho_1})$, $(\clos{\mathcal{P}_2}, \clos{\mathcal{F}_2^{s,u}}, \clos{\rho_2})$ the bifoliated planes of $\Phi_1$, $\Phi_2$ up to surgeries on $\Gamma_1$, $\Gamma_2$ endowed with their natural foliations and group actions by $\pi_1(M_1-\Gamma^{M_1}_1)$ and $\pi_1(M_2-\Gamma^{M_2}_2)$ respectively.

\begin{theorem}[Theorem C]\label{t.generalbarbot}
 We can perform Dehn-Goodman-Fried surgeries on the orbits $\Gamma^{M_1}_1$ of $(\Phi_1, M_1)$ so that the flow obtained after surgery is orbitally equivalent to $(\Phi_2, M_2)$ and so that $\Gamma^{M_2}_2$ are the orbits of $\Phi_2$ corresponding to $\Gamma^{M_1}_1$ after surgery (see Definition \ref{d.associatedorbits}) if and only if there exists a homeomorphism $h: \clos{\mathcal{P}_1}\rightarrow \clos{\mathcal{P}_2}$ such that: 

\begin{enumerate}
    \item the image by $h$ of any stable (resp. unstable) leaf in $\clos{\mathcal{F}_1^{s}}$ (resp. $\clos{\mathcal{F}_1^{u}}$) is a stable (resp. unstable) leaf in $\clos{\mathcal{F}_2^{s}}$ (resp. $\clos{\mathcal{F}_2^{u}}$)
    \item there exists an isomorphism $\alpha: \pi_1(M_1-\Gamma^{M_1}_1) \rightarrow \pi_1(M_2-\Gamma^{M_2}_2)$ such that for every $g\in \pi_1(M_1-\Gamma^{M_1}_1)$ and every $x\in \clos{\mathcal{P}_1}$ we have $$h(\clos{\rho_1}(g)(x))= \clos{\rho_2}(\alpha(g))(h(x))$$ 
\end{enumerate}

\end{theorem}

\begin{proof}
Suppose that after performing surgeries on $(\Phi_1, M_1)$ along the orbits $\Gamma_1^{M_1}$, we obtain a flow that is orbitally equivalent to $(\Phi_2, M_2)$ and that $\Gamma_2^{M_2}$ is the set of periodic orbits of $\Phi_2$ corresponding to $\Gamma_1^{M_1}$ after surgery. Denote by  $\Phi_1-\Gamma^{M_1}_1$ and $\Phi_2-\Gamma^{M_2}_2$ the restrictions of the flows $\Phi_1$ and $\Phi_2$ on $M_1-\Gamma^{M_1}_1$ and $M_2-\Gamma^{M_2}_2$ and  by $\widetilde{\Phi_1-\Gamma^{M_1}_1}$ and $\widetilde{\Phi_2-\Gamma^{M_2}_2}$ the lifts of $\Phi_1-\Gamma^{M_1}_1$ and $\Phi_2-\Gamma^{M_2}_2$ on $\widetilde{M_1-\Gamma^{M_1}_1}$ and $\widetilde{M_2-\Gamma^{M_2}_2}$, the universal covers of $M_1-\Gamma^{M_1}_1$ and $M_2-\Gamma^{M_2}_2$. Thanks to Theorem \ref{t.surgerygivespseudoanosov}, the flows $(\Phi_1-\Gamma^{M_1}_1, M_1-\Gamma^{M_1}_1)$ and $(\Phi_2-\Gamma^{M_2}_2, M_2-\Gamma^{M_2}_2)$ are orbitally equivalent. The orbital equivalence, say $H$, between $\Phi_1-\Gamma^{M_1}_1$ and $\Phi_2-\Gamma^{M_2}_2$ defines an isomorphism $\alpha: \pi_1(M_1-\Gamma^{M_1}_1)\rightarrow \pi_1(M_2-\Gamma^{M_2}_2)$ and lifts to an orbital equivalence $\widetilde{H}$ between $\widetilde{\Phi_1-\Gamma^{M_1}_1}$ and $\widetilde{\Phi_2-\Gamma^{M_2}_2}$, such that: 
\begin{enumerate}
    \item $\widetilde{H}$ sends the lift on $\widetilde{M_1-\Gamma^{M_1}_1}$ of a stable (resp. unstable) leaf of $\Phi_1-\Gamma_1$ to the lift on $\widetilde{M_2-\Gamma^{M_2}_2}$ of a stable (resp. unstable) leaf  of $\Phi_2-\Gamma_2$ 
    \item if $\psi_1$ and $\psi_2$ denote the actions by deck transformations of $\pi_1(M_1-\Gamma_1^{M_1})$ and $\pi_1(M_2-\Gamma_2^{M_2})$ on $\widetilde{M_1-\Gamma^{M_1}_1}$ and $\widetilde{M_2-\Gamma^{M_2}_2}$, then for every $g\in \pi_1(M_1-\Gamma_1^{M_1})$ and every $x\in\widetilde{M_1-\Gamma^{M_1}_1}$ we have $$\widetilde{H}(\psi_1(g)(x))= \psi_2(\alpha(g))(\widetilde{H}(x))$$
\end{enumerate}
 Since $\widetilde{H}$ sends orbits of $\widetilde{\Phi_1-\Gamma^{M_1}_1}$ to orbits of $\widetilde{\Phi_2-\Gamma^{M_2}_2}$, it induces a homeomorphism $h:\clos{\mathcal{P}_1}-\clos{\Gamma_1} \rightarrow \clos{\mathcal{P}_2}-\clos{\Gamma_2}$  between the orbit spaces of the lifted flows (see Proposition \ref{p.actionpi1minusorbits}) such that 
\begin{enumerate}
    \item the image by $h$ of any stable/unstable leaf in $\clos{\mathcal{F}_1^{s,u}}$ is a stable/unstable leaf in $\clos{\mathcal{F}_2^{s,u}}$ 
    \item for every $g\in \pi_1(M_1-\Gamma^{M_1}_1)$ and every $x\in \clos{\mathcal{P}_1}-\clos{\Gamma_1}$ we have $$h(\clos{\rho_1}(g)(x))= \clos{\rho_2}(\alpha(g))(h(x))$$
\end{enumerate}
Since $H$ sends a neighborhood of an orbit in $\Gamma^{M_1}_1$ to the neighborhood of a unique orbit in $\Gamma_2$, the homeomorphism $h$ sends a punctured neighborhood of any point in $\clos{\Gamma_1}$ to a punctured neighborhood of a unique point in $\clos{\Gamma_2}$. Therefore, $h$ can be extended to a homeomorphism from $\clos{\mathcal{P}_1}$ to  $\clos{\mathcal{P}_2}$, which gives us the desired result. 

Let us now show the converse. Assume that there exists a homeomorphism $h: \clos{\mathcal{P}_1}\rightarrow \clos{\mathcal{P}_2}$ with the above properties. Fix an orientation on  $\clos{\mathcal{P}_1}$ and $ \clos{\mathcal{P}_2}$, so that $h$ is orientation preserving. First, notice that 

\vspace{0.2cm}
\noindent\textbf{Claim.} $h(\clos{\Gamma_1})=\clos{\Gamma_2}$; hence the number of orbits in $\Gamma^{M_1}_1$ is equal to the number of orbits in $\Gamma^{M_2}_2$

\vspace{0.2cm}
\noindent \textit{Proof of Claim.} Indeed, using the Proposition \ref{p.stabilizersz2}, the number of orbits in $\Gamma^{M_1}_1$ (resp. $\Gamma^{M_2}_2$) is equal to the number of $\clos{\rho_1}$-orbits (resp. $\clos{\rho_2}$-orbits)  in $\clos{\Gamma_1}$ (resp. $\clos{\Gamma_2}$) or equivalently to the number of  $\clos{\rho_1}$-orbits (resp. $\clos{\rho_2}$-orbits) of points $\clos{\mathcal{P}_1}$ (resp. $\clos{\mathcal{P}_2}$), whose stabilizers are isomorphic to $\mathbb{Z}^2$. By hypothesis, $h$ is equivariant with respect to the action of the fundamental groups and sends points with $\mathbb{Z}^2$ stabilizers to points with $\mathbb{Z}^2$ stabilizers. We conclude that $\Gamma^{M_1}_1$ and $\Gamma^{M_2}_2$ contain the same number of orbits. \qed

\vspace{0.2cm}

Recall that $\clos{\mathcal{P}_1}-\clos{\Gamma_1}$ is the bifoliated plane of $\Phi_1-\Gamma_1$. By taking a pair of directions on $\clos{\mathcal{P}_1}-\clos{\Gamma_1}$ that define an orientation on $\clos{\mathcal{P}_1}$ compatible with our original choice of orientation and by adding at the end the direction of the flow in $\widetilde{M_1-\Gamma_1^{M_1}}$, we define an orientation on $\widetilde{M_1-\Gamma^{M_1}_1}$ and therefore an orientation on $M_1$. Similarly, we can define an orientation on $M_2$. 

Denote by $...,s_{-1, \clos{\gamma}},s_{0, \clos{\gamma}},s_{1, \clos{\gamma}},s_{2, \clos{\gamma}},...$ the stable leaves in $\clos{\mathcal{F}^s}$ containing $\clos{\gamma}\in\clos{\Gamma_1}$ (resp.  $\clos{\gamma}\in\clos{\Gamma_2}$) ordered clockwisely for our choice of orientation of $\clos{\mathcal{P}_1}-\clos{\Gamma_1}$ (resp. $\clos{\mathcal{P}_2}-\clos{\Gamma_2}$). In other words, we assume that the  $...,s_{-1, \clos{\gamma}},s_{0, \clos{\gamma}},s_{1, \clos{\gamma}},s_{2, \clos{\gamma}},...$ are ordered in such a way that for every $i\in \mathbb{Z}$ the leaf $s_{i, \clos{\gamma}}$ is not separated from $s_{i+1, \clos{\gamma}}$ in the stable leaf space of $\clos{\mathcal{F}^s_1}$ (resp. $\clos{\mathcal{F}^s_2}$) and that for any $x\in s_{i, \clos{\gamma}}$ the direction at $x$ pointing inside the connected component of $\clos{\mathcal{P}_1}-s_{i, \clos{\gamma}}$ (resp.  $\clos{\mathcal{P}_2}-s_{i, \clos{\gamma}}$) containing $s_{i+1, \clos{\gamma}}$ followed by the direction at $x$ pointing inside $s_{i, \clos{\gamma}}$ and away from $\clos{\gamma}$ defines an orientation compatible with our original orientation of $\clos{\mathcal{P}_1}-\clos{\Gamma}_1$ (resp.  $\clos{\mathcal{P}_2}-\clos{\Gamma_2}$). Recall that thanks to the previous orientations on $M_1$ and $M_2$, one can construct as in Section \ref{s.elementsfixingpointsingammabar}, for every element $\clos{\gamma}$ in $\clos{\Gamma_1}$ (resp. $\clos{\Gamma_2}$) a basis $p^{M_1}_{\clos{\gamma}},m^{M_1}_{\clos{\gamma}}$ (resp. $p^{M_2}_{\clos{\gamma}},m^{M_2}_{\clos{\gamma}}$) of $\text{Stab}(\clos{\gamma})$ satisfying Proposition \ref{p.secondproofofcommutation}. 

Using our previous Claim and the fact that $h$ is equivariant with respect to the group actions on $\clos{\mathcal{P}_1}$ and $\clos{\mathcal{P}_2}$, we get that for every $\clos{\gamma}\in \clos{\Gamma_1}$, the isomorphism $\alpha$ sends a basis of $\text{Stab}(\clos{\gamma})$ to a basis of $\text{Stab}(h(\clos{\gamma}))$, where $h(\clos{\gamma})\in \clos{\Gamma_2}$. Fix $\clos{\gamma}\in \clos{\Gamma_1}$ and $\gamma\in \Gamma^{M_1}_1$ the periodic orbit in $M_1$ that is associated to $\clos{\gamma}$ (see Proposition \ref{p.periodicinbifoliated}). Thanks to the previous fact, there exists a matrix $A_\gamma\in \text{GL}_2(\mathbb{Z})$ such that $A_\gamma=\begin{pmatrix}
a_\gamma & c_\gamma \\
b_\gamma & d_\gamma
\end{pmatrix}$ and 
\begin{equation}\label{eq.imageofm}
    \alpha(m^{M_1}_{\clos{\gamma}})=(m^{M_2}_{h(\clos{\gamma})})^{a_\gamma}\cdot (p^{M_2}_{h(\clos{\gamma})})^{b_\gamma}
\end{equation}
\begin{equation}\label{eq.imageofp}
    \alpha(p^{M_1}_{\clos{\gamma}})=(m^{M_2}_{h(\clos{\gamma})})^{c_\gamma}\cdot (p^{M_2}_{h(\clos{\gamma})})^{d_\gamma}
\end{equation}

As indicated by our choice of notation, 

\vspace{0.2cm}
\noindent\textbf{Claim.} $A_\gamma$ only depends on $\gamma$ and not on the choice of lift of $\gamma$ in $\clos{\mathcal{P}_1}$

\vspace{0.2cm}
\noindent \textit{Proof of Claim.} Indeed, take $g\in \pi_1(M_1-\Gamma^{M_1}_1)$. By our definition of $m^{M_1}_{\clos{\gamma}}, p^{M_1}_{\clos{\gamma}}$ (see Section \ref{s.elementsfixingpointsingammabar}) and since the actions $\clos{\rho_1}$ and $\clos{\rho_2}$ are orientation preserving, we have that 
\begin{itemize}
    \item $m^{M_1}_{\clos{\rho_1}(g)(\clos{\gamma})}= gm^{M_1}_{\clos{\gamma}}g^{-1}$,  $p^{M_1}_{\clos{\rho_1}(g)(\clos{\gamma})}=gp^{M_1}_{\clos{\gamma}}g^{-1}$ 
    \vspace{0.7cm}
    \item $m^{M_2}_{h(\clos{\rho_1}(g)(\clos{\gamma}))}= \alpha(g)m^{M_2}_{h(\clos{\gamma})}\alpha(g)^{-1}$, $p^{M_2}_{h(\clos{\rho_1}(g)(\clos{\gamma}))}= \alpha(g)p^{M_2}_{h(\clos{\gamma})}\alpha(g)^{-1}$  \end{itemize}
Therefore,

$$\alpha(m^{M_1}_{\clos{\rho_1}(g)(\clos{\gamma})})= (m^{M_2}_{h(\clos{\rho_1}(g)(\clos{\gamma}))})^{a_\gamma} \cdot (p^{M_2}_{h(\clos{\rho_1}(g)(\clos{\gamma}))})^{b_\gamma}$$
$$\alpha(p^{M_1}_{\clos{\rho_1}(g)(\clos{\gamma})})= (m^{M_2}_{h(\clos{\rho_1}(g)(\clos{\gamma}))})^{c_\gamma} \cdot (p^{M_2}_{h(\clos{\rho_1}(g)(\clos{\gamma}))})^{d_\gamma}$$

\qed 

\vspace{0.2cm}

By the above, to every $\gamma\in \Gamma^{M_1}_1$ we can associate a matrix $$A_\gamma=\begin{pmatrix}
a_\gamma & c_\gamma \\
b_\gamma & d_\gamma
\end{pmatrix}\in \text{GL}_2(\mathbb{Z})$$

In addition to belonging to $\text{GL}_2(\mathbb{Z})$ the matrices $A_\gamma$ satisfy several additional properties that will crucial for the construction of our proof. 

\vspace{0.5cm}
\textit{On the properties of the coefficients of $A_\gamma$.} Once again fix $\clos{\gamma}\in \clos{\Gamma_1}$ and $\gamma\in \Gamma^{M_1}_1$ the periodic orbit in $M_1$ that is associated to $\clos{\gamma}$. First let us show that, 

\vspace{0.2cm}

\begin{empheq}[box=\widefbox]{equation}\label{eq.generatortogenerator}
\alpha(P_{\clos{\gamma}}^{M_1})=P_{h(\clos{\gamma})}^{M_2}
\end{empheq}

\vspace{0.2cm}

\noindent \textit{Proof of Claim.} Since the action $\clos{\rho_1}$ preserves the orientation of $\clos{\mathcal{P}_1}$, we have that $\text{Stab}(s_{i,\clos{\gamma}})=\text{Stab}(s_{j,\clos{\gamma}})$ for every $i,j\in \mathbb{Z}$. The same result of course applies for the action of $\clos{\rho_2}$ on $\clos{\mathcal{P}_2}$. As $h$ is equivariant with respect to  $\clos{\rho_1},  \clos{\rho_2}$, we have that $\alpha(\text{Stab}(s_{0,\clos{\gamma}}))=\text{Stab}(h(s_{0,\clos{\gamma}}))$. By Proposition \ref{p.propertiespbaraction}, the previous stabilizers are isomorphic to $\mathbb{Z}$ and thanks to Proposition \ref{p.secondproofofcommutation}, they are respectively generated by $P_{\clos{\gamma}}^{M_1}$ and $P_{h(\clos{\gamma})}^{M_2}$. Finally, using the equivariance of $h$ and the fact that by definition $P_{\clos{\gamma}}^{M_1}$ and $P_{h(\clos{\gamma})}^{M_2}$ are the unique generators $\text{Stab}(s_{0,\clos{\gamma}})$ and $\text{Stab}(h(s_{0,\clos{\gamma}}))$ acting on $s_{0,\clos{\gamma}}$ and $h(s_{0,\clos{\gamma}})$ as topological expansions, we get the desired result. \qed 

\vspace{0.2cm}

Equation \ref{eq.generatortogenerator} translates to a condition on the coefficients of $A_\gamma$ (see Equation \ref{eq.relationxyn1n2}). Indeed, recall that thanks to Proposition \ref{p.secondproofofcommutation} and Lemma \ref{l.parallelmeridian}, there exist $n_1, n_2\in \mathbb{N}^*$, $k_1\in \llbracket 1, n_1\rrbracket$ and  $k_2\in \llbracket 1, n_2\rrbracket$ such that 
\begin{equation}\label{eq.valueofP1}
    P_{\clos{\gamma}}^{M_1}=(m_{\clos{\gamma}}^{M_1})^{-\frac{(n_1-k_1)}{gcd(n_1,k_1)}} \cdot (p_{\clos{\gamma}}^{M_1})^{\frac{n_1}{gcd(n_1,k_1)}}
\end{equation}
\begin{equation}\label{eq.valueofP2}
    P_{h(\clos{\gamma})}^{M_2}=(m_{h(\clos{\gamma})}^{M_2})^{-\frac{(n_2-k_2)}{gcd(n_2,k_2)}} \cdot (p_{h(\clos{\gamma})}^{M_2})^{\frac{n_2}{gcd(n_2,k_2)}}
\end{equation}
where $gcd(n_1,k_1)$ and $gcd(n_2,k_2)$ denote respectively the greatest common divisors of $n_1, k_1$ and $n_2, k_2$. We remind the reader that prior to Lemma \ref{l.parallelmeridian}, we provided an explicit definition of $n_1,k_1$ thanks to which we have that if $\gamma$ is associated by Proposition \ref{p.aroundcircleprong} to a local model for a pseudo-hyperbolic orbit with $n(\gamma)\geq 2$ prongs and rotation $r(\gamma) \in \llbracket 0,p-1\rrbracket$, then $n_1=n(\gamma)$ and $r(\gamma)\in \{k_1,n_1-k_1\}$. Furthermore, thanks to Proposition \ref{p.secondproofofcommutation}, we have that $\clos{\rho_1}(m_{\clos{\gamma}}^{M_1})$ and $\clos{\rho_1}(p_{\clos{\gamma}}^{M_1})$ act as translations by $n_1$ and $n_1-k_1$ on the set of $\{s_{i,\clos{\gamma}}|i\in \mathbb{Z}\}$. Similar results hold for $n_2,k_2$. It follows that $\{s_{i,\clos{\gamma}}|i\in \mathbb{Z}\}$ (resp. $\{s_{i,h(\clos{\gamma})}|i\in \mathbb{Z}\}$) consists of $gcd(n_1,k_1)$ (resp. $gcd(n_2,k_2)$) stable leaves up to the action of $\text{Stab}(\clos{\gamma})$ (resp. $\text{Stab}(h(\clos{\gamma}))$). 

By our previous arguments, we can conclude, thanks to the equivariance of $h$, that  

\begin{equation}\label{eq.equalgcd12}
    gcd(n_1,k_1)=gcd(n_2,k_2)
\end{equation} 

Moreover, thanks to Equations \ref{eq.imageofm}, \ref{eq.imageofp} and  \ref{eq.valueofP1} we have that 

\begin{align*}
    \alpha(P^{M_1}_{\gamma})&=\alpha((m_{\clos{\gamma}}^{M_1})^{-\frac{(n_1-k_1)}{gcd(n_1, k_1)}} \cdot (p_{\clos{\gamma}}^{M_1})^{\frac{n_1}{gcd(n_1,k_1)}}) \\ 
    &= (m_{h(\clos{\gamma})}^{M_2})^{-\frac{a_\gamma\cdot (n_1-k_1)}{gcd(n_1,k_1)}}\cdot (p_{h(\clos{\gamma})}^{M_2})^{-\frac{b_\gamma\cdot (n_1-k_1)}{ gcd(n_1,k_1)}} (m_{h(\clos{\gamma})}^{M_2})^{\frac{c_\gamma\cdot n_1}{ gcd(n_1,k_1)}}\cdot (p_{h(\clos{\gamma})}^{M_2})^{\frac{d_\gamma\cdot n_1}{gcd(n_1,k_1)}}
\end{align*}

By combining the above equation with Equations \ref{eq.generatortogenerator}, \ref{eq.valueofP2}, we get that

\begin{equation*}
  \left\{
    \begin{aligned}
      & {-\frac{(n_2-k_2)}{gcd(n_2,k_2)}}= -\frac{a_\gamma\cdot (n_1-k_1)}{gcd(n_1,k_1)}+ \frac{c_\gamma\cdot n_1}{ gcd(n_1,k_1)}\\
      & \frac{n_2}{gcd(n_2,k_2)}= -\frac{b_\gamma\cdot (n_1-k_1)}{gcd(n_1,k_1)}+ \frac{d_\gamma\cdot n_1}{gcd(n_1,k_1)}
    \end{aligned}
  \right.
\end{equation*}

Finally, by using Equation \ref{eq.equalgcd12} and by multiplying everywhere by $gcd(n_1,k_1)$, we get that 
\begin{equation}\label{eq.relationxyn1n2}
    \left\{
    \begin{aligned}
      & -(n_2-k_2)= -a_\gamma\cdot (n_1-k_1)+ c_\gamma\cdot n_1\\
      & n_2= -b_\gamma\cdot (n_1-k_1)+ d_\gamma\cdot n_1
    \end{aligned}
  \right.
\end{equation}

The previous system of equations forms an important relation between the coefficients of the matrix $A_\gamma$ and $n_1,k_1,n_2,k_2$. In particular, the following inequality (which is an immediate consequence of the fact that $\Phi_2$ is pseudo-Anosov and therefore does not admit circle $1$-prong singularities) will be used multiple times in this proof. 

\begin{equation}\label{eq.conditiononeprong}
-b_\gamma \cdot (n_1-k_1)+ d_\gamma\cdot n_1 \geq 2 
\end{equation}

Finally, before moving on with the proof of Theorem C, let us show that, by using the previous conditions on the coefficients of $A_\gamma$ and also our choices of orientations of $\clos{\mathcal{P}_1}-\clos{\Gamma_1}$ and $\clos{\mathcal{P}_2}-\clos{\Gamma_2}$, we get that

\vspace{0.2cm}

\begin{empheq}[box=\widefbox]{equation}\label{eq.positivedet}
\text{det}(A_\gamma)=1
\end{empheq}

\vspace{0.2cm}

\textit{Proof of Claim.} Indeed, recall that we oriented $\clos{\mathcal{P}_1}-\clos{\Gamma_1}$ and $\clos{\mathcal{P}_2}-\clos{\Gamma_2}$ so that $h$ be orientation preserving. This guarantees that if we order clockwisely the stable leaves of $\clos{\gamma}$ for our choice of orientation of $\clos{\mathcal{P}_1}-\clos{\Gamma_1}$, then their images by $h$ are also ordered clockwisely for our choice of orientation of $\clos{\mathcal{P}_2}-\clos{\Gamma_2}$. By Proposition \ref{p.secondproofofcommutation}, $\clos{\rho_2}(m_{h(\clos{\gamma})}^{M_2})$ acts as a positive translation on the set $\{s_{i, h(\clos{\gamma})}|i\in \mathbb{Z}\}$. By the equivariance of $h$, this implies that $\clos{\rho_1}(\alpha^{-1}(m_{h(\clos{\gamma})}^{M_2}))$ also acts as a positive translation on the set of $\{s_{i, \clos{\gamma}}|i\in \mathbb{Z}\}$. It is easy to check that $$\alpha^{-1}(m^{M_2}_{h(\clos{\gamma})})=(m^{M_1}_{\clos{\gamma}})^{\frac{d_\gamma}{\text{det}(A_\gamma)}}\cdot (p^{M_1}_{\clos{\gamma}})^{-\frac{b_\gamma}{\text{det}(A_\gamma)}}$$

Using Proposition \ref{p.secondproofofcommutation}, we deduce that $$ \frac{d_\gamma\cdot n_1}{\text{det}(A_\gamma)}-\frac{b_\gamma\cdot (n_1-k_1)}{\text{det}(A_\gamma)}>0$$

Thanks to the second equation of \ref{eq.relationxyn1n2}, we get the desired result. \qed 

\vspace{0.2cm}

Let us now move on with the proof of Theorem C. Consider $(M_3,\Phi_3)$ the flow obtained from $(M_1,\Phi_1)$ by performing a $(d_\gamma,-b_\gamma)$-surgery on $\gamma$ for every $\gamma\in \Gamma^{M_1}_1$. In order to finish the proof of the theorem, it suffices to show that $(M_3,\Phi_3)$ is a well defined pseudo-Anosov flow that is orbitally equivalent to $(M_2,\Phi_2)$.  

Let us begin by showing that 

\vspace{0.2cm}
\noindent\textbf{Claim.} $(M_3,\Phi_3)$ is a pseudo-Anosov flow 

\vspace{0.2cm}

\noindent \textit{Proof of Claim.} Consider $(M^*, \Phi^*)$ a flow obtained after blowing-up the orbits of $\Gamma_1^{M_1}$ (see our construction in Section \ref{s.blowupflows}). The points in $\clos{\Gamma_1}$ correspond to periodic orbits in $\Gamma_1^{M_1}$ that were blown-up to tori in $M^*$. Let $\Pi^*_M: M^* \rightarrow M$ be the blow-down map associated to the previous blow-up. Fix $\clos{\gamma}\in \clos{\Gamma_1}$, $\gamma$ its associated periodic orbit in $\Gamma_1^{M_1}$ and $T_\gamma$ its associated torus in $M^*$. Let $\gamma_P,\gamma_p,\gamma_m$ be the curves in $T_\gamma$ constructed in Section \ref{s.elementsfixingpointsingammabar}, whose homotopy classes in $\pi_1(M^*)\cong\pi_1(M_1-\Gamma_1)$ are respectively $P^{M_1}_{\clos{\gamma}},p^{M_1}_{\clos{\gamma}}$ and $m^{M_1}_{\clos{\gamma}}$. Performing a $(d_\gamma,-b_\gamma)$-surgery on $\gamma$ consists in crushing  $T_\gamma$ into a circle by collapsing to points the leaves of a collapsible foliation $\mathcal{F}_{\sigma_{\clos{\gamma}}}$ (see our discussion prior to Lemma \ref{l.goodfoli} for a definition of a collapsible foliation), whose every leaf is freely homotopic to $\sigma_{\clos{\gamma}}:=(m_{\clos{\gamma}}^{M_1})^{d_\gamma} \cdot (p_{\clos{\gamma}}^{M_1})^{-b_\gamma}$. 

Thanks to Equations \ref{eq.valueofP1} and \ref{eq.conditiononeprong}, $\sigma_{\clos{\gamma}}$ is not freely homotopic to $P_{\clos{\gamma}}$ and every leaf of $\mathcal{F}_{\sigma_{\clos{\gamma}}}$ will intersect $\gamma_P$ at least twice. It follows from Theorem \ref{t.surgerygivespseudoanosov} that $(M_3,\Phi_3)$ is indeed a pseudo-Anosov flow. \qed

\vspace{0.2cm}

It suffices now to show that $(M_3,\Phi_3)$ is orbitally equivalent to $(M_2,\Phi_2)$. Denote by $\Gamma^{M_3}_3\subset M_3$ the periodic orbits corresponding to $\Gamma^{M_1}_1$ after surgery, by $\clos{\mathcal{P}_3}$ the bifoliated plane up of $\Phi_3$ up to surgeries on $\Gamma^{M_3}_3$, $\clos{\mathcal{F}_3^{s,u}}$ the stable and unstable foliations of $\clos{\mathcal{P}_3}$, $\clos{\Gamma_3}$ the lifts on $\clos{\mathcal{P}_3}$ of the orbits in $\Gamma^{M_3}_3$ and by $\clos{\rho_3}$ the action of $\pi_1(M_3-\Gamma_3^{M_3})$ on $\clos{\mathcal{P}_3}$. Since $M_1-\Gamma_1^{M_1}$ and $M_3-\Gamma_3^{M_3}$ are homeomorphic, $M_3$ is naturally endowed with an orientation. Using this orientation, define as in Section \ref{s.elementsfixingpointsingammabar} for every $\clos{\gamma}\in \clos{\Gamma_3}$ the elements $m^{M_3}_{\clos{\gamma}}, p^{M_3}_{\clos{\gamma}} , P^{M_3}_{\clos{\gamma}}\in \text{Stab}(\clos{\gamma})$. 

By the first part of this proof, there exists a homeomorphism  $H:\clos{\mathcal{P}_1}\rightarrow \clos{\mathcal{P}_3}$ satisfying the following: 
\begin{enumerate}
    \item the image by $H$ of any stable (resp. unstable) leaf in $\clos{\mathcal{F}_1^{s}}$ (resp. $\clos{\mathcal{F}_1^{u}}$) is a stable (resp. unstable) leaf in $\clos{\mathcal{F}_3^{s}}$ (resp. $\clos{\mathcal{F}_3^{u}}$ )
    \item there exists $\delta:\pi_1(M_1-\Gamma^{M_1}_1)\rightarrow \pi_1(M_3-\Gamma^{M_3}_3)$ a group isomorphism such that  for every $g\in \pi_1(M_1-\Gamma_1)$ and every $x\in \clos{\mathcal{P}_1}$, we have $$H(\clos{\rho_1}(g)(x))= \clos{\rho_3}(\delta(g))(H(x))$$ 
\end{enumerate}

Furthermore, let us show that for any $\clos{\gamma}\in \clos{\Gamma_1}$ we have that  
\begin{equation}\label{eq.bar3}
    \delta(P^{M_1}_{\clos{\gamma}})= P^{M_3}_{H(\clos{\gamma})}
\end{equation}
\begin{equation}\label{eq.bar4}
    \delta^{-1}(m^{M_3}_{H(\clos{\gamma})})= (m_{\clos{\gamma}}^{M_1})^{d_\gamma} \cdot (p_{\clos{\gamma}}^{M_1})^{-b_\gamma}
\end{equation}
\begin{equation}\label{eq.bar5}
    \delta^{-1}(p^{M_3}_{H(\clos{\gamma})})= (m_{\clos{\gamma}}^{M_1})^{-c_\gamma} \cdot (p_{\clos{\gamma}}^{M_1})^{a_\gamma}
\end{equation}
     Indeed, take $\clos{\gamma}\in \clos{\Gamma_1}$, $\gamma$ its associated orbit in $\Gamma_1^{M_1}$ and $\gamma'\in\Gamma_3^{M_3}$ the periodic orbit of $\Phi_3$ corresponding to $\gamma$ after surgery. One can prove Equation \ref{eq.bar3} either by repeating the argument that lead us to the proof of Equation \ref{eq.generatortogenerator} or by noticing that $\Phi^*$ is a flow obtained from $\Phi_1$ (resp. $\Phi_3$) after blowing-up every orbit in $\Gamma_1^{M_1}$ (resp. $\Gamma_3^{M_3}$) and that $P^{M_1}_{\clos{\gamma}}$ (resp. $P^{M_3}_{H(\clos{\gamma})}$) corresponds to the homotopy class of a periodic orbit of $\Phi^*$ on $T_\gamma$ (resp. $T_{\gamma'}=T_\gamma$).
    
\vspace{0.2cm}
\textit{Proof of \ref{eq.bar4}.} Notice that since after the surgery every leaf of the collapsible foliation $\mathcal{F}_{\sigma_{\clos{\gamma}}}$ became homotopically trivial, the natural meridian on $T_{\gamma'}=T_\gamma$ when blowing up $\gamma'$ is homotopic to a leaf of $\mathcal{F}_{\sigma_{\clos{\gamma}}}$, whereas the natural meridian on $T_\gamma$ when blowing up $\gamma$ is freely homotopic to $\gamma_m$. This implies that one of the following is true : 

\begin{itemize}
    \item $\delta^{-1}(m^{M_3}_{H(\clos{\gamma})})= \sigma_{\clos{\gamma}}=(m_{\clos{\gamma}}^{M_1})^{d_\gamma} \cdot (p_{\clos{\gamma}}^{M_1})^{-b_\gamma}$ or 
    \item $\delta^{-1}(m^{M_3}_{H(\clos{\gamma})})= \sigma_{\clos{\gamma}}^{-1}=(m_{\clos{\gamma}}^{M_1})^{-d_\gamma} \cdot (p_{\clos{\gamma}}^{M_1})^{b_\gamma}$
\end{itemize}

Fix $\clos{\gamma}\in \clos{\Gamma_1}$. Recall that, thanks to our choice of orientations on $M_1$ and $M_3$, the orbital equivalence between $\Phi_1$ minus $\Gamma_1^{M_1}$ and $\Phi_3$ minus $\Gamma_3^{M_3}$ is orientation preserving. By lifting the previous orientations, one can canonically orient the universal covers of $M_1-\Gamma_1^{M_1}$ and $M_3-\Gamma_3^{M_3}$. Moreover, using the direction of the lifted flows the previous orientations define respectively an orientation on $\clos{\mathcal{P}_1}-\clos{\Gamma_1}$ and on $\clos{\mathcal{P}_3}-\clos{\Gamma_3}$. It is easy to see that $H$ is orientation preserving for our previous choice of orientations; hence $H$ takes the set of stable leaves of $\clos{\gamma}$ ordered clockwisely to the set of stable leaves of $H(\clos{\gamma})$ ordered clockwisely. Since by Proposition \ref{p.secondproofofcommutation}, $\clos{\rho_3}(m^{M_3}_{H(\clos{\gamma})})$ acts as a positive translation on the clockwise ordered set of stable leaves of $H(\clos{\gamma})$, the equivariance of $H$ forces $\clos{\rho_1}(\delta^{-1}(m^{M_3}_{H(\clos{\gamma})}))$ to also acts as a positive translation on the clockwise ordered set of stable leaves of $\clos{\gamma}$. Thanks to Proposition \ref{p.secondproofofcommutation},  $\clos{\rho_1}(\sigma_{\clos{\gamma}})= \clos{\rho_1}((m_{\clos{\gamma}}^{M_1})^{d_\gamma} \cdot (p_{\clos{\gamma}}^{M_1})^{-b_\gamma})$ acts as a translation by $$d_\gamma\cdot n_1 -b_\gamma\cdot (n_1-k_1)$$  on the clockwise ordered set of stable leaves of $\clos{\gamma}$. This number is positive thanks to Equation \ref{eq.relationxyn1n2}. It follows that $\delta^{-1}(m^{M_3}_{H(\clos{\gamma})})= (m_{\clos{\gamma}}^{M_1})^{d_\gamma} \cdot (p_{\clos{\gamma}}^{M_1})^{-b_\gamma}$. \qed 

\vspace{0.2cm}

\vspace{0.2cm}
\textit{Proof of \ref{eq.bar5}.} First, notice that since $H$ is equivariant with respect to the group actions, we have that $\delta^{-1}(\text{Stab}(H(\clos{\gamma})))=\text{Stab}(\clos{\gamma})$ and thus the isomorphism $\delta^{-1}$ takes a base of $\text{Stab}(H(\clos{\gamma}))$ to a base of $\text{Stab}(\clos{\gamma})$. It follows from Equation \ref{eq.bar4} that there exist $x,y \in \mathbb{Z}$ such that 
\begin{equation}\label{eq.preimagedelta}
    \delta^{-1}(p^{M_3}_{H(\clos{\gamma})})=(m_{\clos{\gamma}}^{M_1})^x \cdot (p_{\clos{\gamma}}^{M_1})^y
\end{equation}

\begin{equation}\label{eq.matrixgl}
\Delta_\gamma:=  \begin{pmatrix}
d_\gamma & x \\-b_\gamma & y 
\end{pmatrix}\in \text{GL}_2(\mathbb{Z})
\end{equation}

The fact that $\Delta_\gamma\in\text{GL}_2(\mathbb{Z})$, implies that $d_\gamma\cdot y + b_\gamma \cdot x =\text{det}(\Delta_\gamma) = \pm1$. Furthermore, using the fact that $$  \text{det}\begin{pmatrix}
d_\gamma & -c_\gamma\\-b_\gamma & a_\gamma
\end{pmatrix}=1$$
we get that there exists $l\in \mathbb{Z}$ such that 

\begin{equation}\label{eq.valuesxy}
  \left\{
    \begin{aligned}
      & x= -\frac{c_\gamma}{\text{det}(\Delta_\gamma)}+ l \cdot d_\gamma\\
      & y= \frac{a_\gamma}{ \text{det}(\Delta_\gamma)} - l \cdot b_\gamma
    \end{aligned}
  \right.
\end{equation}
It therefore suffices to prove that $\text{det}(\Delta_\gamma)=1$ and that $l=0$. As before, thanks to Proposition \ref{p.secondproofofcommutation} and Lemma \ref{l.parallelmeridian}, we have that there exist $n_3\in \mathbb{N}^*$ and $k_3\in \llbracket 1, n_1\rrbracket$ such that 
\begin{equation}\label{eq.stabilizer3}
    P_{H(\clos{\gamma})}^{M_3}=(m_{H(\clos{\gamma})}^{M_3})^{-\frac{(n_3-k_3)}{gcd(n_3,k_3)}} \cdot (p_{H(\clos{\gamma})}^{M_3})^{\frac{n_3}{gcd(n_3,k_3)}}
\end{equation} 
where $gcd(n_3,k_3)$ denotes the greatest common divisors of $n_3, k_3$. Once again, if $\gamma'\in \Gamma_3^{M_3}$ is the periodic orbit of $\Phi_3$ associated to $H(\clos{\gamma})$, then if $\gamma'$ is associated by Proposition \ref{p.aroundcircleprong} to a local model for a pseudo-hyperbolic orbit with $n(\gamma')$ prongs and rotation $r(\gamma')$, then $n_3=n(\gamma')$ and $r(\gamma')\in \{k_3,n_3-k_3\}$. Also, as before, by Proposition \ref{p.secondproofofcommutation}, $\clos{\rho_3}(m_{H(\clos{\gamma})}^{M_3})$ and $\clos{\rho_3}(p_{H(\clos{\gamma})}^{M_3})$ act as translations by $n_3$ and $n_3-k_3$ on the set of clockwise ordered stable leaves of $H(\clos{\gamma})$ and $gcd(n_3,k_3)$ corresponds to the total number of stable leaves of $H(\clos{\gamma})$ up to the action of $\text{Stab}(H(\clos{\gamma}))$. 

The equivariance of $H$ implies that 

\begin{equation}\label{eq.equalgcd}
  gcd(n_3,k_3)=gcd(n_1,k_1)  
\end{equation}

Moreover, thanks to Equations \ref{eq.bar3} and \ref{eq.preimagedelta}, we have that 
\begin{align*}
    \delta^{-1}(P^{M_3}_{H(\clos{\gamma})})&=\delta^{-1}((m_{\clos{H(\gamma)}}^{M_3})^{-\frac{(n_3-k_3)}{gcd(n_3,k_3)}} \cdot (p_{\clos{H(\gamma)}}^{M_3})^{\frac{n_3}{gcd(n_3,k_3)}}) \\ 
    &= (m_{\clos{\gamma}}^{M_1})^{-\frac{d_\gamma\cdot (n_3-k_3)}{gcd(n_3,k_3)}}\cdot (p_{\clos{\gamma}}^{M_1})^{\frac{b_\gamma\cdot (n_3-k_3)}{gcd(n_3,k_3)}} (m_{\clos{\gamma}}^{M_1})^{\frac{x\cdot n_3}{ gcd(n_3,k_3)}}\cdot (p_{\clos{\gamma}}^{M_1})^{\frac{y\cdot n_3}{gcd(n_3,k_3)}}
\end{align*}

By combining the above equation with Equations \ref{eq.valueofP1}, \ref{eq.bar3}, \ref{eq.equalgcd},  we get that

\begin{equation*}
  \left\{
    \begin{aligned}
      & {-\frac{(n_1-k_1)}{gcd(n_1,k_1)}}= -\frac{d_\gamma\cdot (n_3-k_3)}{gcd(n_1,k_1)}+ \frac{x\cdot n_3}{ gcd(n_1,k_1)}\\
      & \frac{n_1}{gcd(n_1,k_1)}= \frac{b_\gamma\cdot (n_3-k_3)}{gcd(n_1,k_1)}+ \frac{y\cdot n_3}{gcd(n_1,k_1)}
    \end{aligned}
  \right.
\end{equation*}

Multiplying everywhere by $gcd(n_1,k_1)$, we get that 
\begin{equation}\label{eq.relationxyn1n3}
    \left\{
    \begin{aligned}
      & -(n_1-k_1)= -d_\gamma\cdot (n_3-k_3)+ x\cdot n_3\\
      & n_1= b_\gamma\cdot (n_3-k_3)+ y\cdot n_3
    \end{aligned}
  \right.
\end{equation}

Thanks to the above, we can show that $\text{det}(\Delta_\gamma)=1$. Indeed, by multiplying the first of the above equations by $b_\gamma$, the second by $d_\gamma$ and by adding the resulting equations we get that:

\begin{equation*}
    (x\cdot b_\gamma+y\cdot d_\gamma)n_3 =d_\gamma\cdot n_1 - b_\gamma\cdot (n_1-k_1)
\end{equation*}

The previous equation, together with Equation \ref{eq.relationxyn1n2} and the definition of $\Delta_\gamma$ imply that 

\begin{equation}\label{eq.valuen3}
   \text{det}(\Delta_\gamma)\cdot n_3=n_2
\end{equation}

Since both $n_2$ and $n_3$ are positive, we get that $\text{det}(\Delta_\gamma)=1$ and that $n_3=n_2$. Let us now show that $l=0$.



Notice now that since $\Delta_\gamma$ is invertible $b_\gamma$ and $d_\gamma$ cannot both be equal to zero. Assume without any loss of generality that $b_\gamma\neq 0$. By solving the second equation of \ref{eq.relationxyn1n3} with respect to $k_3$ and then by using Equation \ref{eq.valuesxy} together with the fact that $\text{det}(\Delta_\gamma)=1$, we get that: 

\begin{align*}
    k_3 &= n_3-\frac{n_1}{b_\gamma}+\frac{ y \cdot n_3}{b_\gamma} \\ &=n_3-\frac{n_1}{b_\gamma}+ \frac{a_\gamma \cdot n_3}{b_\gamma}- l \cdot n_3 \\ &= (1+\frac{a_\gamma}{b_\gamma}-l)\cdot n_3 -\frac{n_1}{b_\gamma} 
\end{align*}

Using Equation \ref{eq.relationxyn1n2} it is easy to check that $$ n_1 = b_\gamma \cdot (n_2 - k_2)+ a_\gamma\cdot  n_2$$ By solving the above equation with respect to $k_2$ we get that $$  k_2= (1+\frac{a_\gamma}{b_\gamma})\cdot n_2 -\frac{n_1}{b_\gamma} $$

By definition, $k_2\in \llbracket 1, n_2\rrbracket $. Therefore, since $n_2=n_3$ (see Equation \ref{eq.valuen3}), $k_3\in \llbracket 1, n_3\rrbracket$ if and only if $k_3=k_2$ or equivalently $l=0$. We thus get the desired result. \qed

\vspace{0.3cm}
We are now ready to prove that $(M_3,\Phi_3)$ is orbitally equivalent to $(M_2,\Phi_2)$. Consider the homeomorphism $K=h\circ H^{-1}: \clos{\mathcal{P}_3}\rightarrow \clos{\mathcal{P}_2}$. Using our definition of $H$ and $h$ together with Equations  \ref{eq.imageofm}, \ref{eq.imageofp}, \ref{eq.generatortogenerator}, \ref{eq.bar4} and \ref{eq.bar5} we have that $K$ satisfies the following: 
\begin{enumerate}
    \item the image by $K$ of any stable (resp. unstable) leaf in $\clos{\mathcal{F}_3^{s}}$ (resp. $\clos{\mathcal{F}_3^{u}}$) is a stable (resp. unstable) leaf in $\clos{\mathcal{F}_2^{s}}$ (resp. $\clos{\mathcal{F}_2^{u}}$ )
    \item there exists an isomorphism $\beta:=\alpha \circ \delta^{-1}: \pi_1(M_3-\Gamma_3) \rightarrow \pi_1(M_2-\Gamma_2)$ such that for every $g\in \pi_1(M_3-\Gamma_3)$, $\clos{\gamma}\in \clos{\Gamma_3}$ and $x\in \clos{\mathcal{P}_3}$ we have $$K(\clos{\rho_3}(g)(x))= \clos{\rho_2}(\beta(g))(K(x))$$ 
    \begin{align}\label{eq.beta1}
       \nonumber \beta(m^{M_3}_{\clos{\gamma}})&= \alpha((m_{H^{-1}(\clos{\gamma})}^{M_1})^{d_\gamma} \cdot (p_{H^{-1}(\clos{\gamma})}^{M_1})^{-b_\gamma})\\   &= (m_{K(\clos{\gamma})}^{M_2})^{a_\gamma\cdot d_\gamma} \cdot (p_{K(\clos{\gamma})}^{M_2})^{d\gamma\cdot b_\gamma} \cdot (m_{K(\clos{\gamma})}^{M_2})^{-b_\gamma\cdot c_\gamma} \cdot (p_{K(\clos{\gamma})}^{M_2})^{-d\gamma\cdot b_\gamma}\\\nonumber&=(m_{K(\clos{\gamma})}^{M_2})^{\text{det}(A_\gamma)}=m_{K(\clos{\gamma})}^{M_2}
    \end{align}
     \begin{align}\label{eq.beta2}
       \nonumber \beta(p^{M_3}_{\clos{\gamma}})&= \alpha((m_{H^{-1}(\clos{\gamma})}^{M_1})^{-c_\gamma} \cdot (p_{H^{-1}(\clos{\gamma})}^{M_1})^{a_\gamma})\\   &= (m_{K(\clos{\gamma})}^{M_2})^{-a_\gamma\cdot c_\gamma} \cdot (p_{K(\clos{\gamma})}^{M_2})^{-c_\gamma\cdot b_\gamma} \cdot (m_{K(\clos{\gamma})}^{M_2})^{a_\gamma\cdot c_\gamma} \cdot (p_{K(\clos{\gamma})}^{M_2})^{d\gamma\cdot a_\gamma}\\\nonumber&=(p_{K(\clos{\gamma})}^{M_2})^{\text{det}(A_\gamma)}=p_{K(\clos{\gamma})}^{M_2}
    \end{align}
\end{enumerate}

Consider now $\clos{\gamma_1},...,\clos{\gamma_n}$ a representative of each $\clos{\rho_3}$-orbit in $\clos{\Gamma_3}$. By the equivariance of $K$ and the fact that  $K(\clos{\Gamma_3})=\clos{\Gamma_2}$, $K(\clos{\gamma_1}),...,K(\clos{\gamma_n})$ are representatives of each $\clos{\rho_2}$-orbit in $\clos{\Gamma_2}$. Let $\gamma_i$ be the periodic orbit in $\Gamma_3^{M_3}$ associated to $\clos{\gamma_i}$. Recall that $m_{K(\clos{\gamma_i})}^{M_3}$ corresponds to the meridian of the torus obtained after blowing up $\gamma_i$. Therefore, by the Seifert-Van Kampen theorem, 
$$\text{ker}(\pi_1(M_3-\Gamma_3)\rightarrow \pi_1(M_3))=<m^{M_3}_{\clos{\gamma_1}},...,m^{M_3}_{\clos{\gamma_n}}>^{\pi_1(M_3-\Gamma_3)}$$
where $<A>^{\pi_1(M_3-\Gamma_3)}$ stands for the normal subgroup of $\pi_1(M_3-\Gamma_3)$ generated by $A$.

Similarly, 
$$\text{ker}(\pi_1(M_2-\Gamma_2)\rightarrow \pi_1(M_2))=<m^{M_2}_{K(\clos{\gamma_1})},...,m^{M_2}_{K(\clos{\gamma_n})}>^{\pi_1(M_2-\Gamma_2)}$$ 
    Using the two previous equalities  together with (\ref{eq.beta1}), we get  that 
\begin{equation}\label{eq.imageofkernel}
\beta(\text{ker}(\pi_1(M_3-\Gamma_3)\rightarrow \pi_1(M_3)))= \text{ker}(\pi_1(M_2-\Gamma_2)\rightarrow \pi_1(M_2))
\end{equation}
Also, by Propositions \ref{p.kernelprojectiongroup} and \ref{p.actionequivariance}, if $(\mathcal{P}_2, \mathcal{F}_2^{s,u},\rho_2)$ and $(\mathcal{P}_3, \mathcal{F}_3^{s,u},\rho_3)$ are the bifoliated planes of $\Phi_2$ and $\Phi_3$ together with their stable/unstable foliations and group actions \begin{equation}\label{eq.p2bar}
    \quotient{\clos{\mathcal{P}_2}}{\text{ker}(\pi_1(M_2-\Gamma_2)\rightarrow \pi_1(M_2))}=\mathcal{P}_2
\end{equation}
\begin{equation}\label{eq.p3bar}
    \quotient{\clos{\mathcal{P}_3}}{\text{ker}(\pi_1(M_3-\Gamma_3)\rightarrow \pi_1(M_3))}=\mathcal{P}_3
\end{equation}
Finally, thanks to Proposition \ref{p.actionequivariance}, to Equations  (\ref{eq.imageofkernel}), (\ref{eq.p2bar}), (\ref{eq.p3bar}) and to the equivariance of $K$ with respect to $\clos{\rho_2}$ and $\clos{\rho_3}$, we get that $K$ projects to a homeomorphism $k:\mathcal{P}_3 \rightarrow \mathcal{P}_2 $ that satisfies the following: 
\begin{enumerate}
    \item the image by $k$ of any stable (resp. unstable) leaf in $\mathcal{F}_3^{s}$ (resp. $\mathcal{F}_3^{u}$) is a stable (resp. unstable) leaf in $\mathcal{F}_2^{s}$ (resp. $\mathcal{F}_2^{u}$)
    \item there exists an isomorphism $\mu: \pi_1(M_3) \rightarrow \pi_1(M_2)$ such that for every $g\in \pi_1(M_3)$ and $x\in \mathcal{P}_3$ we have $$k(\rho_3(g)(x))= \rho_2(\mu(g))(k(x))$$ 
\end{enumerate}
We conclude that $\Phi_2$ is orbitally equivalent to $\Phi_3$ thanks to Theorem \ref{t.barbottheor}. 
\end{proof}

\section{Markovian families as coordinate systems in $\clos{\mathcal{P}}$}\label{s.homotopiesofpaths}
Having proven Theorem C, from now on we will turn our attention to the proof of Theorem B. Thanks to Theorem C, in order to prove Theorem B, it suffices to show that the geometric type of a strong Markovian family $\mathcal{R}$ contains enough information to completely describe not only the bifoliated plane up to surgeries along the boundary periodic points of $\mathcal{R}$, but also its associated foliations and its associated group action. 

The proof of the above result is rather long and technical and thus will be split into three sections. In this section we will define several tools, which are going to be important for us throughout the proof of Theorem B. More specifically, in Section \ref{s.liftmarkovianfamilies}, we will define the notion of strong Markovian family on the bifoliated plane up to surgeries and  its set of associated geometric types. We will also show that the previous set coincides with the set of geometric types associated to the original Markovian family on the bifoliated plane (see Proposition \ref{p.liftedmarkovfamiliessamegeomtype}). Next, in Section \ref{ss.rectanglepaths}, we will define the notion of rectangle path, which can be thought as a discretized curve in the bifoliated plane up to surgeries that can be described via a geometric type. In Section \ref{s.rectanglepathscoordinates}, we will explain how these rectangle paths can be used as ``discrete coordinates" in the bifoliated plane up to surgeries and how they can be used in order to compare two strong Markovian actions preserving two strong Markovian families that are associated to the same equivalence class of geometric types. 

\subsection{Strong Markovian families on the bifoliated plane up to surgeries}\label{s.liftmarkovianfamilies}
Once again, in this section we will work in most general setting in which our results hold, namely strong Markovian actions. Fix $\mathcal{P}$ a plane endowed with an orientation, $\rho: G\rightarrow \text{Homeo}(\mathcal{P})$ an orientation preserving strong Markovian action, preserving the pair of singular foliations $\mathcal{F}^s$ and $\mathcal{F}^u$ (see Definition \ref{d.folisingular}) and leaving invariant a strong Markovian family $\mathcal{R}$. Denote by $\Gamma$ the boundary periodic points of $\mathcal{R}$, by $\clos{\mathcal{P}}$ the bifoliated plane of $\rho$ up to surgeries on $\Gamma$, by $\clos{\mathcal{F}^s}$ and $\clos{\mathcal{F}^u}$ the stable and unstable foliations in $\clos{\mathcal{P}}$, by $\clos{\rho}:\clos{G}\rightarrow \clos{\mathcal{P}}$ the lift of $\rho$ on $\clos{\mathcal{P}}$, by $\clos{\Gamma}$ the lift of $\Gamma$ on $\clos{\mathcal{P}}$ and by $\clos{\pi}$ the projection from $\clos{\mathcal{P}}$ to $\mathcal{P}$. Using as a starting point our definition of strong Markovian family in $\mathcal{P}$, we will now  extend the notion of strong Markovian family in $\clos{\mathcal{P}}$. 

We define rectangles in $\clos{\mathcal{P}}$ in the exact same way as for $\mathcal{P}$ (see Definition \ref{d.standardpolygon}). Notice that any rectangle in $\mathcal{P}$ whose interior does not contain any point of $\Gamma$ lifts to a rectangle in $\clos{\mathcal{P}}$. This is why, we will  define a \emph{strong Markovian family}  in $\clos{\mathcal{P}}$ as the lift on $\clos{\mathcal{P}}$ of a strong Markovian family of $\mathcal{P}$, whose boundary periodic points are exactly $\Gamma$. Fix for the rest of this section $\clos{\mathcal{R}}$ to be the strong Markovian family in $\clos{\mathcal{P}}$ obtained by lifting $\mathcal{R}$.

\begin{defi}\label{d.successorinpbar} For any two rectangles $R_1,R_2 \in \overline{\mathcal{R}}$, we will say that $R_1$ is a \emph{predecessor} (resp. \emph{successor}) of $R_2$ if $\overset{\circ}{R_1}\cap \overset{\circ}{R_2} \neq \emptyset$ and $\clos{\pi}(R_1)\subset \mathcal{P}$ is a predecessor (resp. successor) of $\clos{\pi}(R_2)$. 

We similarly define \emph{predecessors/successors of generation $k\geq 2$} and \emph{$s$-crossing predecessors/successors} on $\clos{\mathcal{P}}$. 
\end{defi}
It is not hard to show that 

\begin{rema}\label{r.equivdefi}
\begin{itemize}
    \item Defining a strong Markovian family (resp. a predecessor/successor or a predecessor/succesor of generation $k$ or a crossing predecessor/successor) as in Definition \ref{d.markovfamily} (resp. \ref{d.successor} or \ref{d.crossingsuccessor}) would lead to a definition equivalent to the one given here above. 
    \item Predecessor/successors, predecessor/successors of generation $k$, crossing predecessor/successors in $\clos{\mathcal{R}}$ still satisfy the results of Lemmas \ref{l.infiniteintersectionverticalrectangles}, \ref{l.infiniteintersectionhorizontalrectangles},  \ref{l.existenceofpredecessors},  \ref{l.npredecessor}, \ref{l.intersectionsuccessorspredecessors},  \ref{l.crossingrectanglesnoperiodicpoints} and \ref{l.crossingrectangleswithperiodicpoints}.  
\end{itemize}

\end{rema}

\begin{defi}\label{d.boundarypbar}
A point in $\overline{\mathcal{P}}$ will be called \emph{periodic} (resp. \emph{boundary arc}, \emph{boundary periodic}) if its projection on $\mathcal{P}$ is a periodic (resp. boundary arc, boundary periodic) point. A stable/unstable leaf $f$ of $\overline{\mathcal{P}}$ will be called \emph{periodic} if there exists a non-trivial $g\in \clos{G}$ such that $\clos{\rho}(g)(f)=f$

\end{defi}

Once again, it is not hard to show that 

\begin{rema}\label{r.equivdefiboundarypoints}
\begin{itemize}
    \item Defining a periodic point (resp. boundary periodic point or a boundary arc point) as in Definition \ref{d.boundaryperiodic} would lead to a definition equivalent to Definition \ref{d.successorinpbar}.
    \item Periodic stable/unstable leaves, boundary arc and boundary periodic points still satisfy the results of Lemma \ref{l.periodicboundary}, Propositions \ref{p.boundaryperiodicpoints}, \ref{p.arcpointsexist}, Remark \ref{r.caracterisationarcpoints}. 
    \item The boundary periodic points of $\clos{\mathcal{R}}$ coincide with $\clos{\Gamma}$ and every boundary arc point of $\clos{\mathcal{R}}$ is contained in the stable or unstable leaf of some point in $\clos{\Gamma}$. 
\end{itemize}

\end{rema}

Next, by repeating the same construction as in Definition \ref{d.geometrictypemarkovfamily}, we can associate to the strong Markovian family $\clos{\mathcal{R}}$ a finite number of geometric types, that we will call once again \emph{the geometric types of $\clos{\mathcal{R}}$} or \emph{the geometric types associated to $\clos{\mathcal{R}}$}. Furthermore, by the same arguments as in the proof of Theorem \ref{t.associatemarkovfamiliestogeometrictype}, we can show that the geometric types of $\clos{\mathcal{R}}$ are pairwise equivalent and that Remark \ref{r.canonicalassociationgeometrictype} remains true for Markovian families in $\clos{\mathcal{P}}$. In other words, 

\begin{rema}\label{r.canonicalassociationgeometrictypeinpbar}
    Given $\clos{\mathcal{R}}$ a strong Markovian family preserved by $\clos{\rho}$, the following triplet of choices of representatives and orientations defines a unique geometric type of $\clos{\mathcal{R}}$: 
    \begin{enumerate}
        \item a choice of representatives $r_1,...,r_n$ of every rectangle orbit in $\clos{\mathcal{R}}$
        \item a choice of orientation of $\clos{\mathcal{P}}-\clos{\Gamma}$
        \item a choice of orientation of the foliations inside every $r_i$ such that the positive stable direction followed by the positive unstable direction inside every $r_i$ define an orientation compatible with our choice of orientation of $\clos{\mathcal{P}}-\clos{\Gamma}$
    \end{enumerate}
    \end{rema}

Recall that thanks to Remark \ref{r.propertiesfolipbar}, the foliations $\clos{\mathcal{F}^s}$ and $\clos{\mathcal{F}^u}$ are always orientable and transversely orientable (even if $\mathcal{F}^s$ and $\mathcal{F}^u$ are not). Thanks to Remark \ref{r.canonicalassociationgeometrictypeinpbar}, the following couple  of choices of representatives and orientations defines a unique geometric type of $\clos{\mathcal{R}}$:

\begin{itemize}
    \item a choice of representatives of every rectangle orbit in $\clos{\mathcal{R}}$
    \item a choice of orientations of $\clos{\mathcal{F}^s}$ and $\clos{\mathcal{F}^u}$ 
\end{itemize}

The set of geometric types of $\clos{\mathcal{R}}$ obtained by all the above pairs of choices will be called the set of \emph{special geometric types of $\clos{\mathcal{R}}$}. Before finishing this section, let us prove two results providing respectively a better description of the set of geometric types and of the set of special geometric types of $\clos{\mathcal{R}}$.  

\begin{prop}\label{p.liftedmarkovfamiliessamegeomtype}
Let $\clos{\mathcal{R}}$ be a strong Markovian family in $\clos{\mathcal{P}}$ preserved by $\clos{\rho}$ and $\mathcal{R}$ its projection on $\mathcal{P}$. The set of geometric types associated to the Markovian families $\clos{\mathcal{R}}$ and $\mathcal{R}$ coincide.
\end{prop}
\begin{lemm}[Realisability lemma] \label{l.geomtypeinclass}
Consider the unique equivalence class of geometric types containing all the geometric types of $\clos{\mathcal{R}}$ and $$\mathcal{G}=(n,(h_i)_{i \in \llbracket 1,n \rrbracket}, (v_i)_{i\in \llbracket 1,n \rrbracket}, \mathcal{H}, \mathcal{V},\phi, u)$$ a geometric type in this class. If $\clos{\rho}$ does not preserve the orientations of the foliations $\clos{\mathcal{F}^s}$ and $\clos{\mathcal{F}^u}$, then $\mathcal{G}$ is a special geometric type associated to $\clos{\mathcal{R}}$. If we also assume that $u= +1$, the same result is true when $\clos{\rho}$ preserves the orientations of the foliations $\clos{\mathcal{F}^s}$ and $\clos{\mathcal{F}^u}$. 
\end{lemm}

\begin{proof}[Proof of the proposition]
  Consider a choice of representatives $r_1,...,r_n$ of every rectangle orbit in $\mathcal{R}$ and a lift $\clos{r_i}\in \clos{\mathcal{R}}$ of every $r_i$. Thanks to Propositions \ref{p.actionequivariance} and  \ref{p.propertiespbaraction}, the rectangles $\clos{r_i}$ are representatives of every rectangle orbit in $\clos{\mathcal{R}}$. Consider next a choice of orientation of $\mathcal{P}$ and of the stable and unstable foliations of every $r_i$ so that the positive stable direction followed by the positive unstable direction inside every $r_i$ define an orientation compatible with our choice of orientation of $\mathcal{P}$. By pulling back the previous choices of orientations by $\clos{\pi}$, we can define an orientation on $\clos{\mathcal{P}}$ and an orientation of the stable and unstable foliations inside every $\clos{r_i}$. It is easy to see that the positive stable direction followed by the positive unstable direction inside every $\clos{r_i}$ define an orientation compatible with our choice of orientation of $\clos{\mathcal{P}}$.

  According to Remarks \ref{r.canonicalassociationgeometrictype} and \ref{r.canonicalassociationgeometrictypeinpbar}, together with these choices of representatives and orientations, the Markovian family  $\clos{\mathcal{R}}$ (resp. $\mathcal{R}$) defines a unique geometric type  $\clos{\mathcal{G}}=(\clos{n},(\clos{h_i})_{i \in \llbracket 1,n \rrbracket}, (\clos{v_i})_{i\in \llbracket 1,n \rrbracket}, \clos{\mathcal{H}}, \clos{\mathcal{V}},\clos{\phi},\clos{u})$ (resp.  $\mathcal{G}=(n,(h_i)_{i \in \llbracket 1,n \rrbracket}, (v_i)_{i\in \llbracket 1,n \rrbracket}, \mathcal{H}, \mathcal{V},\phi, u)$).  

It suffices to show that $\mathcal{G}$ is equal to $\clos{\mathcal{G}}$. In order to do so, let us go back to the construction of $\mathcal{G}$ and $\clos{\mathcal{G}}$ (see Definition \ref{d.geometrictypemarkovfamily}). The number $n$ (resp. $\clos{n}$) in $\mathcal{G}$ (resp. $\clos{\mathcal{G}}$) is equal to the number of distinct orbits of rectangles in $\mathcal{R}$ (resp. $\clos{\mathcal{R}}$) by the action of $\rho$ (resp. $\clos{\rho}$). Once again, thanks to Propositions \ref{p.actionequivariance} and \ref{p.propertiespbaraction}, $n=\clos{n}$. 

Next, the number $h_i$ (resp. $v_i$)  corresponds to the number of successors (resp. predecessors) of $r_i$. Similarly, the number $\clos{h_i}$ (resp. $\clos{v_i}$)  corresponds to the number of successors (resp. predecessors) of $\clos{r_i}$. By Definition \ref{d.successorinpbar}, 
$\clos{r_i}$ and $r_i$ have the same number of predecessors and successors. In other words, $h_i=\clos{h}_i$ and $v_i=\clos{v}_i$.

Now, let us order the successors (resp. predecessors) of every $r_i$ from bottom to top (resp. from left to right) using our choice of orientation of the unstable  (resp. stable) foliation inside $r_i$. We will denote by $H_i^k$ (resp. $V_i^k$) the $k$-th successor (resp. predecessor) of $r_i$ for our previous choice of order. Recall that $\mathcal{H}=\{H^k_i, i\in \llbracket 1, n\rrbracket, k\in \llbracket 1, h_i\rrbracket\}$ and $\mathcal{V}=\{V^k_i, i\in \llbracket 1, n\rrbracket, k\in \llbracket 1, v_i\rrbracket\}$. We similarly define the families $(\clos{V_i^k})_{k\in \llbracket 1, v_i\rrbracket}$, $(\clos{H^k_i})_{k\in \llbracket 1, h_i\rrbracket}$ for every rectangle $\clos{r_i}$ and also the sets  $\clos{\mathcal{H}}, \clos{\mathcal{V}}$. By our choice of orientations of the stable and unstable foliations inside every $\clos{r_i}$, we have that $\clos{\pi}(\clos{H^k_i})=H_i^k$ and $\clos{\pi}(\clos{V_i^k})=V_i^k$. Therefore, we can canonically identify  $\mathcal{H}$ with $\clos{\mathcal{H}}$ and $\mathcal{V}$ with  $\clos{\mathcal{V}}$. 

It remains to prove that $\phi$ is equal to $\clos{\phi}$ and that $u$ is equal to $\clos{u}$ (modulo the identification of $\mathcal{H}$ with $\clos{\mathcal{H}}$ and of $\mathcal{V}$ with  $\clos{\mathcal{V}}$). First, recall 

\begin{enumerate}
    \item that $\phi(H^k_i)=V^l_j$ if and only if there exists $g\in G$ such that $\rho(g)(H^k_i)=r_j$ and $\rho(g)(r_i)$ is the $l$-th predecessor (from left to right) of $r_j$
    \item that the previous $g$ is uniquely defined by the property $\rho(g)(H^k_i)=r_j$ (see Proposition \ref{p.independencegeomtype})
    \item that as $H^k_i$ intersects $r_i$, it inherits a natural orientation for both its stable and unstable foliations and   that $u(H^k_i)=+1$ if and only if  $\rho(g)$ sends the previous oriented stable and unstable foliations in $H^k_i$ to the oriented stable and unstable foliations in $r_j$, according to our initial choices of orientations
\end{enumerate} The functions $\clos{\phi}$ and $\clos{u}$ are similarly defined. Thanks to Propositions \ref{p.actionequivariance} and \ref{p.propertiespbaraction}, we deduce that $\phi(H^k_i)=V^l_j$ if and only if $\clos{\phi}(\clos{H^k_i})=\clos{V^l_j}$. Similarly, it is not difficult to check using Proposition \ref{p.actionequivariance} that $u(H^k_i)=1 $ if and only if $\clos{u}(\clos{H^k_i})=1$. Finally, using the identifications between $\mathcal{H}$ and $\clos{\mathcal{H}}$, $\mathcal{V}$ and $\clos{\mathcal{V}}$, we have that $\phi$ and $\clos{\phi}$ (resp. $u$ and $\clos{u}$) define the same functions from $\mathcal{H}$ to $\mathcal{V}$ (resp. from $\mathcal{H}$ to $\{\pm1\}$), which proves that $\mathcal{G}$ is equal to $\clos{\mathcal{G}}$ and thus finishes the proof of the proposition. 
\end{proof}

\begin{proof}[Proof of the realisability lemma]
Throughout this proof, we will assume that $\clos{\rho}$ does not preserve the orientations of 
$\clos{\mathcal{F}^{s,u}}$. The case where $\clos{\rho}$ preserves the orientations of the foliations follows from similar arguments. 

Take $\mathcal{G}=(n,(h_i)_{i \in \llbracket 1,n \rrbracket}, (v_i)_{i\in \llbracket 1,n \rrbracket}, \mathcal{H}, \mathcal{V},\phi, u)$ a geometric type associated to $\clos{\mathcal{R}}$. We will show that $\mathcal{G}$ is a special geometric type of $\clos{\mathcal{R}}$. In order to do so, it suffices to show that there exists a set of representatives of every $\clos{\rho}$-orbit of rectangles in $\clos{\mathcal{R}}$ and a choice of orientations of $\clos{\mathcal{F}^{s,u}}$ such that $\mathcal{G}$ is the special geometric type associated to $\clos{\mathcal{R}}$ for this choice of representatives and orientations. 

Choose arbitrarily a set of representatives $r_1,...,r_n$ of every rectangle orbit in $\clos{\mathcal{R}}$ and choose also arbitrary orientations for $\clos{\mathcal{F}^s}$ and $\clos{\mathcal{F}^u}$. By Theorem \ref{t.associatemarkovfamiliestogeometrictype} applied for $\clos{\mathcal{R}}$, the family $\clos{\mathcal{R}}$ together with our previous choices of representatives and orientations defines a unique special geometric type $$\mathcal{G}'=(n,(h'_i)_{i \in \llbracket 1,n' \rrbracket}, (v'_i)_{i\in \llbracket 1,n' \rrbracket}, \mathcal{H'}, \mathcal{V'},\phi', u')$$ that is equivalent to $\mathcal{G}$. Our goal from now on consists in adapting our choices of representatives and orientations in order to produce, thanks to Remark \ref{r.canonicalassociationgeometrictypeinpbar}, a special geometric type equal to $\mathcal{G}$.

First, notice that, thanks to Proposition \ref{p.independencegeomtype}, up to reindexing the $(h'_i,v'_i)$, we have that $h'_i=h_i$ and $v'_i=v_i$. Hence, $\mathcal{G}'$ is of the form $(n,(h_i)_{i \in \llbracket 1,n \rrbracket}, (v_i)_{i\in \llbracket 1,n \rrbracket}, \mathcal{H'}, \mathcal{V'},\phi', u')$. 

Using the geometric interpretation of a geometric type (see Section \ref{s.markovpartitionsgeomtypesdefi}), from now on we will think of a geometric type as a set of rectangles and subrectangles. More specifically, given any geometric type $\mathcal{K}=(n,(h_i)_{i \in \llbracket 1,n \rrbracket}, (v_i)_{i\in \llbracket 1,n \rrbracket}, \mathcal{H}_\mathcal{K}, \mathcal{V}_\mathcal{K},\phi_\mathcal{K}, u_\mathcal{K})$, we will denote by $R_1(\mathcal{K}),$..., $R_n(\mathcal{K})$ the rectangles associated to $\mathcal{K}$, by  $\{H_i^j(\mathcal{K})|j\in \llbracket 1, h_i\rrbracket\}\subset \mathcal{H}_\mathcal{K}$ the horizontal subrectangles of $R_i(\mathcal{K})$ ordered from bottom to top and by  $\{V_i^j(\mathcal{K})|j\in \llbracket 1, v_i\rrbracket\}\subset \mathcal{V}_\mathcal{K}$ the vertical subrectangles of $R_i(\mathcal{K})$ ordered from left to right. Recall that thanks to the choice of orientations of $\clos{\mathcal{F}^{s,u}}$ that gave rise to $\mathcal{G}'$, each of the elements of $\{H_i^j(\mathcal{K})|j\in \llbracket 1, h_i\llbracket\}$ (resp. $\{V_i^j(\mathcal{K})|j\in \llbracket 1, v_i\llbracket\}\subset \mathcal{V}_\mathcal{K}$) corresponds to a unique successor (resp. predecessors) of $r_i$.

Denote by $H$ the equivalence between $\mathcal{G}$ and $\mathcal{G}'$. Following the notations of Definition \ref{d.equivalentgeomtypes}, consider for every $i\in \llbracket 1 , n \rrbracket$ the integers $\epsilon_i(\mathcal{G}'),\epsilon_i'(\mathcal{G}')\in \{-1,+1\}$. Recall that  

\begin{enumerate}
    \item for every $i\in \llbracket 1 , n \rrbracket$, the equivalence $H$ defines a bijection between the sets of horizontal/vertical subrectangles of $R_i(\mathcal{G})$ and $R_i(\mathcal{G'})$. More specifically, $H$ defines a monotonous bijection between the sets $\{V_i^j(\mathcal{G})|j\in \llbracket 1, v_i\rrbracket\}$ (resp. $\{H_i^j(\mathcal{G})|j\in \llbracket 1, h_i\rrbracket\}$) and  $\{V_i^j(\mathcal{G'})|j\in \llbracket 1, v_i\rrbracket\}$ (resp. $\{H_i^j(\mathcal{G'})|j\in \llbracket 1, h_i\rrbracket\}$)
    \item if the previous bijection respects the order of the vertical (resp. horizontal) subrectangles, then $\epsilon_i(\mathcal{G}')= +1$ (resp. $\epsilon_i'(\mathcal{G}')= +1$) and if it reverses it, then  $\epsilon_i(\mathcal{G}')= -1$ (resp. $\epsilon_i'(\mathcal{G}')= -1$)
    \item thanks to Definition \ref{d.equivalentgeomtypes} we have either that $\epsilon_i(\mathcal{G}') \cdot \epsilon'_i(\mathcal{G}')=-1$ for all $i\in \llbracket 1 , n \rrbracket$ or that $\epsilon_i(\mathcal{G}') \cdot \epsilon_i'(\mathcal{G}')=+1$ for all $i\in \llbracket 1 , n \rrbracket$
\end{enumerate} 

We will first show that by eventually changing one of our initial choices of orientations, we can assume that $\epsilon_i(\mathcal{G}') \cdot \epsilon_i'(\mathcal{G}')=+1$ for all $i\in \llbracket 1 , n \rrbracket$. Indeed, changing the orientation of  $\clos{\mathcal{F}^s}$ without modifying our choice of representatives or of orientation of $\clos{\mathcal{F}^u}$, produces a new special geometric type $$\mathcal{G}''=(n,(h_i)_{i \in \llbracket 1,n \rrbracket}, (v_i)_{i\in \llbracket 1,n \rrbracket}, \mathcal{H''}, \mathcal{V''},\phi'', u'')$$ associated to $\clos{\mathcal{R}}$ that can be obtained from $\mathcal{G}'$ by reindexing all the vertical subrectangles in $\mathcal{G}'$ as follows: $$V_i^j(\mathcal{G}')\rightarrow V_i^{v_i-j}(\mathcal{G}'')$$

In more geometric terms, going from left to right, for any $i\in \llbracket 1, n\rrbracket$ the first vertical subrectangle of $R_i(\mathcal{G}')$, becomes the last vertical subrectangle of the rectangle $R_i(\mathcal{G}'')$  and vice versa.

By construction of $\mathcal{G}''$ and by Theorem \ref{t.associatemarkovfamiliestogeometrictype}, there exists $H^*$ an equivalence between $\mathcal{G}$ and $\mathcal{G}''$. Once again following the notations of Definition \ref{d.equivalentgeomtypes}, consider for every $i\in \llbracket 1 , n \rrbracket$ the integers $\epsilon_i(\mathcal{G}''),\epsilon_i'(\mathcal{G}'')\in \{-1,+1\}$. Since the order of the vertical subrectangles has changed with respect to $\mathcal{G}'$, but the order of the horizontal subrectangles remains the same, we have that: 
$$\epsilon_i(\mathcal{G}'')=\epsilon_i(\mathcal{G}') \text{ and } \epsilon'_i(\mathcal{G}'')=-\epsilon'_i(\mathcal{G}') $$

Therefore, by eventually changing our initial choice of orientation on $\clos{\mathcal{F}^s}$, we can assume that for every $i \in \llbracket 1 , n \rrbracket$ we have $\epsilon_i(\mathcal{G}') \cdot \epsilon_i'(\mathcal{G}')=+1$. 

Suppose that there exists $i_0\in \llbracket 1, n\rrbracket$ such that $\epsilon_{i_0}(\mathcal{G}')= \epsilon'_{i_0}(\mathcal{G}')=-1$. Take $g\in \clos{G}$ such that $\clos{\rho}$ reverses the orientation of the foliations $\clos{\mathcal{F}^{s,u}}$ and replace in our initial choice of representatives $r_{i_0}$ by $\clos{\rho}(g)(r_{i_0})$. Changing our representatives in this way without modifying the orientations of the foliations, produces a new special geometric type $$\mathcal{G}'''=(n,(h_i)_{i \in \llbracket 1,n \rrbracket}, (v_i)_{i\in \llbracket 1,n \rrbracket}, \mathcal{H'''}, \mathcal{V'''},\phi''', u''')$$ associated to $\clos{\mathcal{R}}$ that can be obtained from $\mathcal{G}'$ by reindexing the horizontal and vertical subrectangles of the rectangle $R_{i_0}(\mathcal{G}')$ in $\mathcal{G}'$ as follows:$$V_{i_0}^j(\mathcal{G}')\rightarrow V_{i_0}^{v_{i_0}-j}(\mathcal{G}''') \text{ and } H_{i_0}^j(\mathcal{G}')\rightarrow H_{i_0}^{h_{i_0}-j}(\mathcal{G}''')$$

In more geometric terms, going from left to right (resp. bottom to top), the first vertical (resp. horizontal) subrectangle of $R_{i_0}(\mathcal{G}')$ becomes the last vertical (resp. horizontal) subrectangle of the rectangle $R_{i_0}(\mathcal{G}''')$ and vice versa.

By construction of $\mathcal{G}'''$ and by Theorem \ref{t.associatemarkovfamiliestogeometrictype}, there exists $H^{**}$ an equivalence between $\mathcal{G}$ and $\mathcal{G}'''$. Once again following the notations of Definition \ref{d.equivalentgeomtypes}, consider for every $i$ the integers $\epsilon_i(\mathcal{G}'''),\epsilon_i'(\mathcal{G}''')\in \{-1,+1\}$. Since the orders of the vertical and horizontal  subrectangles in $R_{i_0}(\mathcal{G}')$ have both changed and the order of the horizontal (or vertical) subrectangles of any other rectangle of $\mathcal{G}'$ remained the same, we have that : 
$$\epsilon_i(\mathcal{G}''')=\epsilon_i(\mathcal{G}') \text{ and } \epsilon'_i(\mathcal{G}''')=\epsilon'_i(\mathcal{G}') \text{ for $i\neq i_0$} $$
$$\epsilon_{i_0}(G''')=-\epsilon_{i_0}(G') \text{ and } \epsilon'_{i_0}(\mathcal{G}''')=-\epsilon'_{i_0}(\mathcal{G}') $$

By eventually changing our initial choice of representatives, we can thus assume that for every $i\in \llbracket 1 , n \rrbracket$ the numbers $\epsilon_i(\mathcal{G}')$ and $ \epsilon_i'(\mathcal{G}')$ are all equal to $+1$ and thus that $\mathcal{G}'=\mathcal{G}$, which gives us the desired result. 

\end{proof}
We end this section with an important corollary of the two previous lemmas. 

\begin{coro}\label{c.specialgeomtypes}
   Consider $\rho_1,\rho_2$ two strong Markovian actions on the planes $\mathcal{P}_1,\mathcal{P}_2$ preserving the strong Markovian families $\mathcal{R}_1,\mathcal{R}_2$, whose boundary periodic points are respectively $\Gamma_1, \Gamma_2$. Denote by $\clos{\mathcal{P}_1},\clos{\mathcal{P}_2}$ the bifoliated planes of $\rho_1,\rho_2$ up to surgeries on $\Gamma_1, \Gamma_2$ and by $\clos{\mathcal{R}_1}, \clos{\mathcal{R}_2}$ the lifts of $\mathcal{R}_1,\mathcal{R}_2$ on $\clos{\mathcal{P}_1},\clos{\mathcal{P}_2}$. 
   
   If $\mathcal{R}_1,\mathcal{R}_2$ are associated to the same equivalence class of geometric types, then $\clos{\mathcal{R}_1}, \clos{\mathcal{R}_2}$ are also associated to the same equivalence class of geometric types. Furthermore, the special geometric types associated to $\clos{\mathcal{R}_1}, \clos{\mathcal{R}_2}$ coincide. 
\end{coro}

\subsection{Rectangle paths in the bifoliated plane up to surgeries}\label{ss.rectanglepaths}
Following the notations of the previous section,

\begin{defi}\label{d.rectanglepath}
A finite sequence of rectangles in $\clos{\mathcal{R}}$ of the form $R_0$,$R_1$,...,$R_n$ will be called a \emph{rectangle path} (of \emph{length} $n+1$) going from $R_0$ to $R_n$ if for every $i \in \llbracket 0, n-1 \rrbracket$ the rectangle $R_{i+1}$ is either a successor or a predecessor of $R_i$. We will say that the rectangle path $R_0$,$R_1$,...,$R_n$ is 
\begin{itemize}
    \item \emph{trivial} if it is of length one 
    \item \emph{closed} if it also verifies $R_n=R_0$
    \item \emph{increasing} (resp. \emph{decreasing}) if it verifies that $R_{i+1}$ is a successor (resp. predecessor) of $R_i$ for every $i$ 
    \item \emph{monotonous} if it is increasing or decreasing
\end{itemize}
\end{defi}

The reader may think of a rectangle path as a discretized curve inside the bifoliated plane up to surgeries. In order to illustrate the close connection between curves and rectangle paths, our goal in this section consists in defining the notion of rectangle path associated to a (good polygonal) curve in $\clos{\mathcal{P}}$ and the notion of (good polygonal) curve associated to some rectangle path.

\begin{defi}\label{d.polygonalcurve}
A continuous curve defined by a topological immersion $\gamma: [0,1]\rightarrow \clos{\mathcal{P}}$  will be called a \emph{polygonal curve} if it is a finite juxtaposition of alternating stable and unstable segments in $\clos{\mathcal{F}^{s,u}}$. Furthermore, if the curve $\gamma$ is closed, we will call it a \emph{closed} polygonal curve.
\end{defi}
 
\begin{defi}\label{d.goodcurve}
By definition, for any polygonal curve $\gamma$ there exists a unique finite sequence $d_0=0<d_1<...<d_n=1$ such that the $\gamma([d_i,d_{i+1}])$ are alternating stable and unstable segments. Any such stable or unstable segment will be called a \emph{side} of $\gamma$. The number $n$ will be called the \emph{length} of the polygonal curve $\gamma$. Also, a polygonal curve $\gamma$ will be called \emph{good} if 
\begin{enumerate}
\item none of its sides belongs to the stable or unstable leaf of a point in $\clos{\Gamma}$ 
\item there do not exist $i\neq j\in \llbracket 0,n-1\rrbracket$ such that $\gamma([d_i,d_{i+1}])$ and $\gamma([d_j,d_{j+1}])$ belong to the same stable or unstable leaf in $\clos{\mathcal{F}^{s,u}}$. 
\end{enumerate}

\end{defi}
\begin{rema}\label{r.gamma0}
 Notice that by our above definition if $\gamma$ is a closed and good polygonal curve, then $\gamma(0)$ is necessarily a corner point of $\gamma$ (i.e. not in the interior of any stable or unstable segment in $\gamma([0,1])$). 
 \end{rema}
\begin{lemm}\label{l.createpolygonalcurve}
Any immersed curve $\gamma: [0,1]\rightarrow \clos{\mathcal{P}}$ that is disjoint from $\clos{\Gamma}$ and for which $\gamma(0)$ and $\gamma(1)$ do not belong to a stable or unstable leaf of a point in $\clos{\Gamma}$ is homotopic relatively to its boundary to a good polygonal curve in $\clos{\mathcal{P}}$

\end{lemm}
\begin{proof}
Indeed, we can locally deform $\gamma$ to a finite juxtaposition of alternating stable and unstable segments. By compactness, $\gamma$ is homotopic relatively to its boundary to a polygonal curve $\gamma'$ in $\clos{\mathcal{P}}-\clos{\Gamma}$. Without any loss of generality, we can assume that $\gamma'$ is the juxtaposition of a stable segment $s_1$, followed by an unstable segment $u_1$,..., ending with an unstable segment $u_n$. 

Suppose that one of the stable segments in $\gamma'$, say $s_k$,  belongs to a stable leaf of a point in $\clos{\Gamma}$. Notice that $s_k\neq s_1$, because of our initial hypothesis. The points in $\clos{\Gamma}$ forming a finite set of $\clos{\rho}$-orbits (see for instance Proposition \ref{p.propertiespbaraction}) and $\clos{G}$ being countable (see Remark \ref{r.extensioncountable}), we get that $\clos{\Gamma}$ is countable and thus their stable and unstable leaves too. Since $\gamma'$ does not intersect $\clos{\Gamma}$, by changing a little bit the ``lengths" of the segments $u_{k-1}$ and $u_k$, we can replace $s_k$ by another stable segment in $s_k$'s neighborhood that does not belong to the countable set of stable leaves of $\clos{\Gamma}$. By a repeated application of the previous argument, we get that $\gamma'$ is homotopic relatively to its boundary to a polygonal curve satisfying the first axiom of Definition \ref{d.goodcurve}.

By the exact same procedure, we can show that $\gamma'$ is homotopic relatively to its boundary to a polygonal curve satisfying also the second axiom of Definition \ref{d.goodcurve}, which gives us the desired result.  

\end{proof}

\begin{lemm}\label{l.finiterepetition}
    Consider $s$ a closed stable segment in $\clos{\mathcal{F}^s}$ that does not belong to the stable leaf of any point in $\clos{\Gamma}$. There can not exist a sequence $(R^{(i)})_{i\in \mathbb{N}}$ of rectangles in $\clos{\mathcal{R}}$ such that
    \begin{enumerate}
        \item an extremity of $s$ lies in the interior of $R^{(0)}$
        \item $R^{(i+1)}$ is a crossing successor of $R^{(i)}$ for every $i\in \mathbb{N}$
        \item $R^{(0)}\cap s\subsetneq... \subsetneq R^{(i)}\cap s \subsetneq s $ for every $i\in \mathbb{N}$ 
    \end{enumerate}
\end{lemm} 
\begin{proof}
    Assume that a sequence of rectangles with the above properties exists. Denote by $s_0$ the extremity of $s$ in $\inte{R^{(0)}}$ and by $s_1$ the other extremity of $s$ (see Figure \ref{f.Proof44}). Since the sequence of closed segments $R^{(i)} \cap s$ is strictly increasing, there exists $x$ a first point of $s$ that is not contained in any of the $R^{(i)}$. Suppose that the $R^{(i)} \cap s$ approach $x$ from the left and consider $[t,x]$ a small stable segment in $s$ on the left of $x$. Since by Remark \ref{r.equivdefi}, $\clos{\mathcal{R}}$ satisfies the conditions of Definition \ref{d.markovfamily}, by the finite return time axiom applied for $\clos{\mathcal{R}}$, there exists $R\in \clos{\mathcal{R}}$ containing $[t,x]$. Notice also that, since $s$ does not belong to the stable leaf of a point in $\clos{\Gamma}$, by Remark \ref{r.equivdefiboundarypoints}, the segment $s$ cannot contain a non-trivial segment of the stable boundary of $R$. It follows that $[t,x]-\lbrace x \rbrace \subset \inte{R}$.

On the other hand, for every $i$, $R^{(i+1)}\cap R^{(i)}$ is a non-trivial horizontal subrectangle of $R^{(i)}$; hence, by Remark \ref{r.equivdefiboundarypoints}, $R^{(n)}$ will be a very thin along the vertical direction and a very long along the horizontal direction rectangle approaching $x$ from the left. Hence, for $n$ big $\overset{\circ}{R^{(n)}}\cap \overset{\circ}{R} \neq \emptyset$, which contradicts the Markovian intersection property  for $\clos{\mathcal{R}}$ (see Remark \ref{r.equivdefi} and Definition \ref{d.markovfamily}) and leads to an absurd. 

\begin{figure}[h!]
\centering 
\includegraphics[scale=0.5]{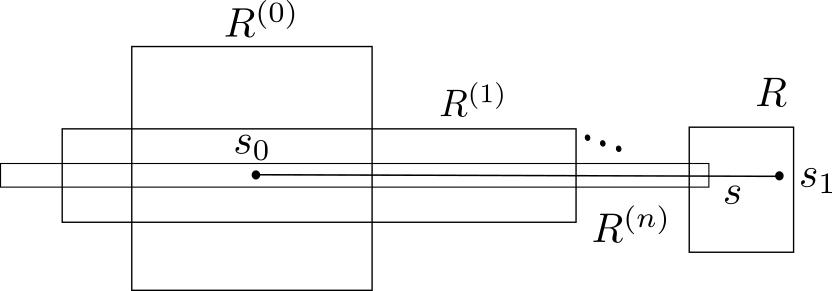}
\caption{}
\label{f.Proof44}
\end{figure}
\end{proof}
Naturally, a similar lemma holds for unstable segments in $\clos{\mathcal{F}^u}$ too.

Consider $\gamma:[0,1]\rightarrow \clos{\mathcal{P}}$ a good polygonal curve in $\clos{\mathcal{P}}$ and $r_0\in \clos{\mathcal{R}}$ such that $\gamma(0)\in \overset{\circ}{r_0}$. Assume without any loss of generality that $\gamma$ is the juxtaposition of the stable segment $s$, followed by the unstable segment $u$, ..., followed at the end by the segment $f$.    

Suppose now that $s$ crosses the unstable boundary component $U$ of $r_0$. Because of Definition \ref{d.goodcurve} and Remark \ref{r.equivdefiboundarypoints}, $s$ does not belong to the stable leaf of a boundary arc or boundary periodic point of $\clos{\mathcal{R}}$. Therefore, by Lemmas \ref{l.crossingrectanglesnoperiodicpoints}, \ref{l.crossingrectangleswithperiodicpoints} applied for $\clos{\mathcal{R}}$,  there exists a unique $U$-crossing successor of $r_0$, say $R^{(1)}$, containing $r_0\cap s$. If $s$ exits $R^{(1)}$, by applying  the same argument, $s$ will enter a unique crossing successor of $R^{(1)}$, say $R^{(2)}$. Thanks to Lemma \ref{l.finiterepetition}, by applying a finite number of times the previous procedure, we will eventually obtain a rectangle $R^{(n)}$, a successor of some generation of $r_0$, for which $s\subset R^{(n)}$. By  Lemma \ref{l.npredecessor} applied for $\clos{\mathcal{R}}$ (see Remark \ref{r.equivdefi}), there exists a unique increasing rectangle path $R_0=r_0,R_1,...,R_m=R^{(n)}$ going from $R_0$ to $R^{(n)}$.
\begin{defi}\label{d.curvetorectanglepath}
We will call the rectangle path $R_0=r_0,R_1,...,R_m=R^{(n)}$, \emph{the rectangle path associated to $s$ starting from $r_0$}. One defines similarly the \emph{rectangle path associated to an unstable segment starting from a rectangle in $\mathcal{R}$}. 

We now define recursively \emph{the rectangle path associated to $\gamma$ starting from $r_0$}, say $path_{r_0}(\gamma)$, as follows: 
\begin{itemize}
    \item when $\gamma$ has length one, $path_{r_0}({\gamma})$ is defined as the rectangle path associated to $s$ starting from $r_0$ 
    \item when $\gamma$ has length $n>1$, consider  $d_0=0<d_1<...<d_n=1$ such that the $\gamma([d_i,d_{i+1}])$ are alternating stable and unstable segments (the $d_i$ are uniquely defined by the previous properties). Denote by $\gamma':=\gamma_{|[d_1,1]}$ the polygonal curve of length $n-1$ obtained by removing from $\gamma$ its first segment, namely the segment $s$. We define $path_{r_0}(\gamma)$ as the concatenation of $path_{r_0}(s)$ followed by $path_{R^{(n)}}(\gamma')$.
    
\end{itemize} 
\end{defi}
 Following our previous notations, let us remark here that when the length of $\gamma$ is  greater than $1$, $path_{r_0}(\gamma)$ is well defined, since by Definition \ref{d.goodcurve} and Remark \ref{r.equivdefiboundarypoints} both $\gamma':[d_1,1]\rightarrow \clos{\mathcal{P}}$ and $s:[0,d_1]\rightarrow \clos{\mathcal{P}}$ are good polygonal curves and $\gamma(d_{1})\in \inte{R^{(n)}}$.

\begin{prop}\label{p.rectpathsstartingending}
For every two rectangles $R_0, R\in \clos{\mathcal{R}}$ there exists a rectangle path starting from $R_0$ and ending at $R$. 
\end{prop}
\begin{proof}
We consider an immersed curve $\gamma: [0,1] \rightarrow \clos{\mathcal{P}}$ such that $\gamma(0)\in \overset{\circ}{R_0}$, $\gamma(1)\in \overset{\circ}{R}$, $\gamma$ does not intersect $\clos{\Gamma}$ and  $\gamma(0),\gamma(1)$ don't belong to the stable or unstable leaf of a point in $\clos{\Gamma}$. By Lemma \ref{l.createpolygonalcurve}, there exists $\gamma'$ a good polygonal curve homotopic to $\gamma$ relatively to its boundary. Consider $R_0,R_1,...R_n$ the rectangle path associated to $\gamma'$ starting from $R_0$.  Let us remark that $R_n$ is not necessarily equal to $R$. However, thanks to Remark \ref{r.equivdefiboundarypoints}, since $\gamma(1)$ does not belong to the stable or unstable leaf of a point in $\clos{\Gamma}$, $\gamma(1)\in \overset{\circ}{R}\cap \inte{R_n}$ and by Lemma \ref{l.npredecessor} applied for $\clos{\mathcal{R}}$, the rectangle $R$ is a $k$-th successor or predecessor of $R_n$ for some $k \in \mathbb{N}$. Assume without any loss of generality that $R$ is a $k$-th successor of $R_n$. In this case, by Lemma \ref{l.npredecessor}, there exists a unique increasing rectangle path $R_n,...,R_{n+k}=R$ from $R_n$ to $R$. The rectangle path $R_0,...,R_{n+k}=R$ is the desired rectangle path.
\end{proof}

Of course the rectangle path of the previous corollary is not unique; it depends on the choice of the original curve $\gamma$ and on how we deformed it into a good polygonal curve $\gamma'$. In fact, for any two rectangles  $R_0, R\in \clos{\mathcal{R}}$ there exist infinitely many distinct rectangle paths going from $R_0$ to $R$: if $R_0,R_1,...,R$ is one such rectangle path, then $R_0,R_1,R_0,R_1,...,R$ is also a rectangle path from $R_0$ to $R$. 

As we have previously explained, it is possible to associate a set of rectangle paths to any good polygonal curve in $\clos{\mathcal{P}}$. The inverse association is also possible:  
\begin{prop}\label{p.rectanglepathtocurve}
Fix a rectangle path $R_0,...,R_n$. There exists a good polygonal curve $\gamma$ such that $R_0,...,R_n$ is the  rectangle path starting from $R_0$ associated to $\gamma$. Furthermore, if $R_0,...,R_n$ is closed, then we can choose $\gamma$ to also be closed. 
\end{prop}
\begin{proof}
Consider $R_0,...,R_n$ a rectangle path and $x\in \inte{R_0}$ that does not belong to the stable or unstable leaf of a point in $\clos{\Gamma}$. Consider now a  polygonal path starting from $x\in R_0$, exiting $R_0$ in order to enter $R_1$ (we use here the fact that $R_1\not \subset R_0$ see Remark \ref{r.crossingpredsucc}), then exiting $R_1$ in order to enter $R_2$ and so on until we reach $R_n$. We ask that our polygonal path ends to a point that does not belong to the stable or unstable leaf of a point in $\clos{\Gamma}$ and that the path does not intersect any point of $\clos{\Gamma}$. If $R_n=R_0$, we will also ask that our polygonal path ends at $x$. 

By the same  arguments as in the proof of Lemma \ref{l.createpolygonalcurve}, using the fact that $x$ does not belong in the stable or unstable leaf of a point in $\clos{\Gamma}$, it is possible to prove that the previous (resp. closed when $R_0=R_n$) polygonal curve can be deformed to a good (resp. and closed) polygonal curve $\gamma$ starting from $x$, exiting $R_0$ in order to enter $R_1$, then exiting $R_1$ in order to enter $R_2$ and so on until we reach $R_n$.

It is now easy to check, using  Definition \ref{d.curvetorectanglepath}, that the rectangle path starting from $R_0$ associated to $\gamma$ is exactly $R_0,...,R_n$. 
\end{proof}

\begin{rema}\label{r.polygonalcurverectassociation}
\begin{itemize}
    \item The construction of $\gamma$ in Proposition  \ref{p.rectanglepathtocurve} is not unique. 
    \item If $r_0,...,r_n$ is a rectangle path associated to a good polygonal curve $\gamma:[0,1]\rightarrow \clos{\mathcal{P}}$, then there exist $0=c_0<c_1<...<c_{n+1}=1$ such that $\gamma([c_i,c_{i+1}])\subset r_i$. The $c_i$ are not unique. Using our choice of $c_i$, we define $\textit{Rect}_{\gamma, r_0}:[0,1]\rightarrow \left\{ r_0,...,r_n \right\}$ as the function sending every interval of the form $[c_i,c_{i+1})$ to $r_i$ and such that  $\textit{Rect}_{\gamma, r_0}(1)=r_n$. 
 \end{itemize}
 \end{rema}
 
 As we have previously proven, good polygonal curves in $\clos{\mathcal{P}}$ are intimately related with rectangle paths in $\clos{\mathcal{R}}$. The function $\textit{Rect}_{\gamma, r_0}$ provides a way to draw the link between a continuous curve and a rectangle path, which the reader can think as the discretization of a continuous good polygonal curve in $\clos{\mathcal{P}}$.

 \begin{conv}\label{conv.rectanglepathimage}
     Even though $\textit{Rect}_{\gamma, r_0}$ is a function that associates to every point of $\gamma$ a rectangle that contains it, from now on and unless explicitly said otherwise we will think of the image by $\textit{Rect}_{\gamma, r_0}$ of any interval in $[0,1]$ as a rectangle path in $\clos{\mathcal{R}}$.
 \end{conv}
 Even though $\textit{Rect}_{\gamma, r_0}$, is not uniquely defined, in the following pages, it will be useful to us to define the function $\textit{Rect}_{\gamma, r_0}$ as in the following lemma: 
 \begin{lemm}\label{l.choiceci}
     Let $\gamma:[0,1]\rightarrow \clos{\mathcal{P}}$ be a good polygonal curve and  $r_0,...,r_n$ its associated  rectangle path starting from $r_0$. We can choose $0=c_0<c_1<...<c_{n+1}=1$ and define $\textit{Rect}_{\gamma, r_0}: [0,1]\rightarrow \{r_0,r_1,...,r_n\}$ so that: 
\begin{itemize}
    \item $\gamma([c_i,c_{i+1}])\subset r_i$ for every $i$
    \item for any stable (resp. unstable) segment $S$ of $\gamma$ the rectangle path $\textit{Rect}_{\gamma,r_0}(\gamma^{-1}(S))$ is either trivial or increasing (resp. decreasing)
    \item if $d_0=0<d_1<...<d_n=1$ are such that the segments $\gamma([d_i,d_{i+1}])$ are alternating stable and unstable segments in $\gamma$, then $\textit{Rect}_{\gamma,r_0}([0,d_i])$ coincides with the rectangle path starting from $R_0$ associated to the good polygonal curve $\gamma_{|[0,d_i]}$ 
\end{itemize}
 \end{lemm}
 \begin{proof}
  Recall that to any choice of $0=c_0<c_1<...<c_{n+1}=1$ such that $\gamma([c_i,c_{i+1}])\subset r_i$, we can associate a function $\textit{Rect}_{\gamma, r_0}$. Let us now construct a set of $c_i$, such that $\textit{Rect}_{\gamma, r_0}$ verifies the desired properties.

  Assume without any loss of generality, that $\gamma$ is the juxtaposition of the stable segment $s_1$, followed by the unstable segment $u_1$,..., ending with the unstable segment $u_m$. Let $c_0=0$, $c_{n+1}=1$, $t_{r_0}= \text{max}\{t\in [0,1]~|~\gamma([0,t])\subset r_0\}$, $t_{s_l}$, $t_{u_l}$ be the unique elements in $[0,1]$ such that $$\gamma([0,t_{s_l}])=s_1\cup u_1\cup s_2\cup... \cup s_l$$ $$\gamma([0,t_{u_l}])=s_1\cup u_1\cup s_2\cup... \cup s_l \cup u_l$$

  Fix $\epsilon>0$ a very small positive number. We begin by defining $c_1$. In order to do so, we will consider two cases. 

  \vspace{0.5cm}
  \textit{Case 1.} Suppose first that \underline{$s_1$ does not exit $r_0$}. If $u_1, s_2,...,u_m$ do not exit $r_0$, then by Definition \ref{d.curvetorectanglepath} the rectangle path associated to $\gamma$ is trivial; hence $n=0$ and  $c_1=c_{n+1}= 1$. If $u_1,s_2,...,s_l$ do not exit $r_0$ and $u_l$ exits $r_0$, then we define $c_1=t_{s_l}+ \epsilon$. Similarly, if $u_1,s_2,...,u_l$ do not exit $r_0$ and $s_{l+1}$ exits $r_0$, then we define $c_1=t_{u_l}+ \epsilon$. Notice that in both of the previous cases, if $\epsilon$ is sufficiently small,  thanks to the fact that $\gamma$ is good, we have that $\gamma([c_0,c_{1}])\subset r_0$. Indeed, by definition none of the segments $s_1,u_1,s_2,...,u_l$ can begin or end at the boundary of a rectangle in $\clos{\mathcal{R}}$. 

   \vspace{0.5cm} 
  \textit{Case 2.} Suppose now that \underline{$s_1$ exits $r_0$ in order to visit $r_{k_1}$, a crossing successor of $r_0$}, then exits $r_{k_1}$ in order to visit $r_{k_2}$, a crossing successor of $r_{k_1}$ and so on till (see Lemma \ref{l.finiterepetition}) it reaches $r_{k_l}$ that verifies $s_1\subset r_{k_l}$ (see Figure \ref{f.choiceoftheci}). Recall that by Definition \ref{d.curvetorectanglepath}, every time that $\gamma$ leaves a rectangle $R$ in order to visit $R'$, a crossing successor of $R$, we add to the rectangle path associated to $\gamma$ the unique increasing rectangle path from $R$ to $R'$. Therefore, $r_{k_1},...,r_{k_l}$ belong in the rectangle path starting from $r_0$ associated to $\gamma$. We will assume that for every $i\in \llbracket 1,l\rrbracket $ the rectangle $r_{k_i}$ corresponds to $(k_i+1)$-th rectangle of the rectangle path starting from $r_0$ associated to $\gamma$. It follows from our previous remark that for every $i\in \llbracket 0,k_l-1\rrbracket $ the rectangle $r_{i+1}$ is a successor of $r_i$ containing $\gamma([0,t_{r_0}])$. Choose now arbitrarily $c_1,...,c_{k_l-1}\in [0,1]$ such that $0<c_1<c_2<c_3<...<c_{k_l-1}<t_{r_0}$ and let $c_{k_l}:=t_{s_1}+\epsilon$. Notice that by our previous arguments, if $\epsilon$ is sufficiently small, then we have that $\gamma([c_{k_l-1},c_{k_l}]) \subset r_{k_l}$ and also that for every $i\in \llbracket 0,k_l-2\rrbracket $ we have $\gamma([c_i,c_{i+1}])\subset \gamma([0,t_{r_0}]) \subset r_i$.

\begin{figure}[h!]
    \centering
    \includegraphics[scale=0.06]{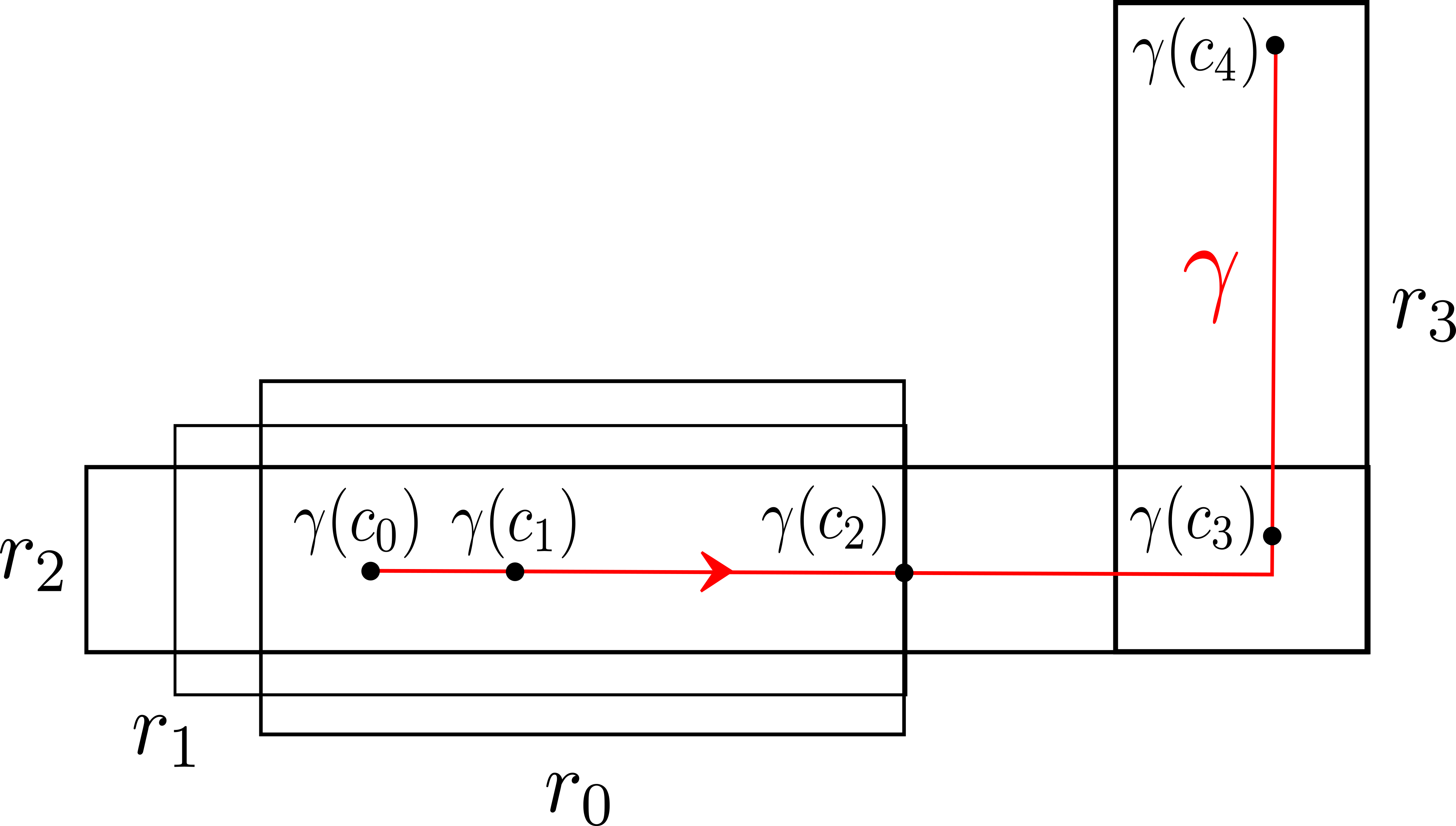}
    \caption{The red curve in the above figure represents a good polygonal curve, whose associated rectangle path starting from $r_0$ is $r_0,r_1,r_2,r_3$. }
    \label{f.choiceoftheci}
\end{figure}

  By our previous construction of the $c_i$ and by Remark \ref{r.polygonalcurverectassociation}, $\textit{Rect}_{\gamma,r_0}(\gamma^{-1}(s_1))$ corresponds to the rectangle path starting from $R_0$ associated to the good polygonal curve $s_1$ and for every stable segment $s\subset s_1$, the rectangle path $\textit{Rect}_{\gamma,r_0}(\gamma^{-1}(s))$ is either trivial or increasing. Therefore, our construction verifies the desired properties. Once the previous procedure ends for $s_1$, we will repeat it for $u_1$,$s_2$,..., till we finally reach $u_m$; we thus get the desired set of $c_i$. 
 \end{proof}

\subsection{Rectangle paths as coordinate systems} \label{s.rectanglepathscoordinates}

Let $\rho_1:G_1 \rightarrow \text{Homeo}(\mathcal{P}_1)$, $\rho_2: G_2 \rightarrow  \text{Homeo}(\mathcal{P}_2)$ be two strong Markovian actions acting on the planes $\mathcal{P}_1$ and $\mathcal{P}_2$ together with their stable and unstable foliations $\mathcal{F}_{1}^{s}$,$\mathcal{F}_{1}^{u}$ (resp. $\mathcal{F}_{2}^{s}$,$\mathcal{F}_{2}^{u}$). Assume that $\rho_1$ and $\rho_2$ preserve  two strong Markovian families $\mathcal{R}_1$ and $\mathcal{R}_2$, whose associated classes of geometric types coincide. 

For every $i\in \{1,2\}$, let us denote by $\Gamma_{i}$ the boundary periodic points associated to $\mathcal{R}_{i}$ in $\mathcal{P}_{i}$, $\clos{\mathcal{P}_{i}}$ the bifoliated plane of $\rho_i$ up to surgeries on $\Gamma_{i}$, $\clos{\rho_i}:\clos{G_i}\rightarrow  \text{Homeo}(\clos{\mathcal{P}_{i}})$ the lift of $\rho_i$ on $\clos{\mathcal{P}_{i}}$, $\clos{\mathcal{F}^{s}_{i}}, \clos{\mathcal{F}^{u}_{i}}$ the stable and unstable  foliations of $\clos{\rho_i}$, $\clos{\mathcal{R}_{i}}$ the lift of $\mathcal{R}_{i}$ on $\clos{\mathcal{P}_{i}}$ and $\clos{\Gamma_{i}}$ the lift of $\Gamma_{i}$ on $\clos{\mathcal{P}_{i}}$. 

Thanks to Corollary, \ref{c.specialgeomtypes}, by appropriately choosing representatives for every rectangle orbit in $\clos{\mathcal{R}_1}$ and $\clos{\mathcal{R}_2}$, and orientations for the foliations $\clos{\mathcal{F}_{1}^{s,u}}$,$\clos{\mathcal{F}_{2}^{s,u}}$, we may assume that $\clos{\mathcal{R}_1}$ and $\clos{\mathcal{R}_2}$ are associated to the same (special) geometric type $$\mathcal{G}=(R_1,...,R_n,(h_i)_{i \in \llbracket 1,n \rrbracket}, (v_i)_{i\in \llbracket 1,n \rrbracket}, \mathcal{H}, \mathcal{V},\phi, u)$$

Notice that we are using here the geometric interpretation of a geometric type and thus we are thinking of $\mathcal{G}$ as a set of rectangles and subrectangles. Using our previous choice of representatives and orientations, we would like to associate to any rectangle path in $\clos{\mathcal{P}_{1}}$ a rectangle path in $\clos{\mathcal{P}_{2}}$. This will allow us to simultaneously navigate into $\clos{\mathcal{P}_{1}}$ and $\clos{\mathcal{P}_{2}}$ by using rectangle paths and to compare (locally) the actions $\clos{\rho_1}$ and $\clos{\rho_2}$.

Recall that every $\clos{\rho_1}$-orbit of rectangles in $\clos{\mathcal{R}_{1}}$ corresponds to a unique rectangle of $\mathcal{G}$. Same for $\clos{\mathcal{R}_{2}}$.

\begin{defi}
We will say that $R^1 \in \clos{\mathcal{R}_1}$ (resp. $R^2 \in \clos{\mathcal{R}_2}$) is of \emph{type $i$} if the $\clos{G_1}$-orbit (resp. $\clos{G_2}$-orbit) of $R^1$ (resp. $R^2$) in $\clos{\mathcal{P}_1}$ (resp. $\clos{\mathcal{P}_2}$) corresponds to the $i$-th rectangle of $\mathcal{G}$. 
\end{defi}

Fix $r_0^1\in \clos{\mathcal{R}_1}$ and $r_0^2 \in \clos{\mathcal{R}_2}$ the unique rectangles of type $1$ among our choices of representatives in $ \clos{\mathcal{R}_1}$ and $ \clos{\mathcal{R}_2}$, thanks to which we previously defined $\mathcal{G}$. From now on, we will call $r_0^{1}$, $r_0^{2}$ the \textit{origin rectangles} in $\clos{\mathcal{P}_{1}}$, $\clos{\mathcal{P}_{2}}$ and we will also call any rectangle path in $\clos{\mathcal{P}_{1}}$, $\clos{\mathcal{P}_{2}}$ starting from $r_0^{1}$, $r_0^{2}$ a \textit{centered rectangle path}.

To any centered rectangle path $r^1_0,...,r^1_N$ in $\clos{\mathcal{P}_{1}}$ it is possible to canonically associate a unique centered rectangle path $r^2_0,...,r^2_N$ in $\clos{\mathcal{P}_{2}}$ such that for every $i\in \llbracket 0,N \rrbracket$ $r^1_i$ and $r^2_i$ are of the same type. 

Indeed, take $r^1_0,...,r^1_N$ a centered rectangle path in $\clos{\mathcal{P}_{1}}$. Let us start by constructing $r^2_1$. Since $r^1_0,...,r^1_N$ is a rectangle path $r^1_1$ is a predecessor or a successor of $r^1_0$. Assume without any loss of generality that it is a successor. Using our choice of orientation of $\clos{\mathcal{F}^{u}_{1}}$, assume that $r^1_1$ is the $k-$th successor of $r^1_0$ from bottom to top. Take $r^2_1$ to be the $k$-th successor of $r^2_0$ from bottom to top for our choice of orientation on $\clos{\mathcal{F}^{u}_{2}}$. By our definition of $r^1_0, r^2_0$ and also thanks to the fact that $\clos{\mathcal{R}_{1}}$ and $\clos{\mathcal{R}_{2}}$ correspond to the same geometric type (for our choices of representatives and orientations), we have that both $r^1_0$ and $r^2_0$ correspond in $\mathcal{G}$ to the rectangle $R_1$, that the orientations of  $\clos{\mathcal{F}^{u}_{1}}$ and $\clos{\mathcal{F}^{u}_{2}}$ inside $r^1_0$ and $r^2_0$ correspond both to the bottom to top orientation in $R_1$ and therefore $r^1_1$ and $r^2_1$ both correspond to $k$-th horizontal subrectangle of $R_1$ from bottom to top. Hence,  $r^1_1$ and $r^2_1$ are of the same type. 

Consider next $\clos{g_1}$ (resp. $\clos{g_2}$) the unique element in $\clos{G_1}$ (resp. $\clos{G_2}$) such that $\clos{\rho_1}(\clos{g_1})$ (resp. $\clos{\rho_2}(\clos{g_2})$) takes $r_1^1$ (resp. $r_1^2$) to the representative of its orbit that was initially chosen during the construction of $\mathcal{G}$. Once again, since $\clos{\mathcal{R}_{1}}$ and $\clos{\mathcal{R}_{2}}$ correspond to the same geometric type (for our choices of representatives and orientations), we have that $\clos{\rho_1}(\clos{g_1})$ reverses the orientation of the foliations if and only if $\clos{\rho_2}(\clos{g_2})$ does (see our definition of the function $u$ in Definition \ref{d.geometrictypemarkovfamily}). It follows that $r^1_1$ and $r^2_1$ correspond both to some rectangle $R_i$ in $\mathcal{G}$ and that the orientations of $\clos{\mathcal{F}^{u}_{1}}$ and $\clos{\mathcal{F}^{u}_{2}}$ inside $r^1_1$ and $r^2_1$ correspond both to the bottom to top orientation in $R_i$ or both to the top to bottom orientation in $R_i$. Using this fact, we can easily construct by the same exact argument the rectangles $r^2_2,...,r_N^2$ by induction.

\begin{defi}\label{d.associatingrectanglepaths}
  We will call the rectangle path $r^2_0,...,r^2_N$ \emph{the (centered) rectangle path in $\clos{\mathcal{P}_2}$ associated to the (centered) rectangle path} $r^1_0,...,r^1_N$.  
\end{defi}

The main idea behind the proof of Theorem B lies in the fact that we can use rectangle paths and associated rectangle paths as coordinate systems in order to compare $\clos{\mathcal{P}_1}$ and $\clos{\mathcal{P}_2}$.  The most important part of the proof of Theorem B consists in proving the following theorem, according to which the two  coordinate systems given by rectangle paths and associated rectangle paths are compatible:
\begin{theorem}\label{t.endingbysamerectangle}
Fix any two centered rectangle paths in $\clos{\mathcal{P}_1}$ ending by the same rectangle. Their associated centered rectangle paths in $\clos{\mathcal{P}_2}$ also end by the same rectangle.
\end{theorem} 
The previous theorem is equivalent to the following:
\begin{theorem}\label{t.closedrectanglespathscorrespondtoclosedpaths}
Fix any closed centered rectangle path in $\clos{\mathcal{P}_1}$. Its associated centered rectangle path in $\clos{\mathcal{P}_2}$ is also closed.
\end{theorem} 
\begin{proof}[Proof of the equivalence of Theorems \ref{t.endingbysamerectangle} and \ref{t.closedrectanglespathscorrespondtoclosedpaths}]$\quad$

(\ref{t.endingbysamerectangle} $\Rightarrow$ \ref{t.closedrectanglespathscorrespondtoclosedpaths}): Take $r^1_0$ the trivial rectangle path and $C_{path}$ any closed centered rectangle path in $\clos{\mathcal{P}_1}$. Since $C_{path}$ is closed, it also ends at $r^1_0$. By Theorem \ref{t.endingbysamerectangle}, the rectangle paths in $\clos{\mathcal{P}_2}$ associated to $C_{path}$ and $r_0^1$ end by the same rectangle. But the trivial rectangle path $r^1_0$ is associated (see our discussion prior to Definition \ref{d.associatingrectanglepaths}) to the trivial rectangle path $r^2_0$. Therefore, the centered rectangle path in  $\clos{\mathcal{P}_2}$ associated to $C_{path}$ is closed.

(\ref{t.closedrectanglespathscorrespondtoclosedpaths} $\Rightarrow$ \ref{t.endingbysamerectangle}): 
Take two centered rectangle paths $r^1_0,...r^1_N$ and $R^1_0,R^1_1,...,R^1_m$ in $\clos{\mathcal{P}_1}$ such that $r^1_0=R^1_0$ and $R^1_m=r^1_N$. The rectangle path $r^1_0,...r^1_N, R^1_{m-1}, R^1_{m-2},...,R^1_1,r^1_0$ is a closed centered rectangle path in $\clos{\mathcal{P}_1}$. Therefore, by Theorem \ref{t.closedrectanglespathscorrespondtoclosedpaths} its associated rectangle path in $\clos{\mathcal{P}_2}$ is also closed. Denote the previous  rectangle path by $r^2_0,...r^2_N, R^2_{m-1}, R^2_{m-2},...,R^2_1,r^2_0$. It is not difficult to see, that $r^2_0,...r^2_N$ is the centered rectangle path associated to $r^1_0,...r^1_N$ and that $r^2_0,R^2_1,...,R^2_{m-1},r^2_N$ is the centered rectangle path associated to $R^1_0,R^1_1,...,R^1_m$. We deduce that the centered rectangle paths in $\clos{\mathcal{P}_2}$ associated to $r^1_0,...r^1_N$ and $R^1_0,R^1_1,...,R^1_m$ end by the same rectangle. 

\end{proof}

\section{A homotopy theory for rectangle paths - Proof of Theorem \ref{t.closedrectanglespathscorrespondtoclosedpaths}} \label{s.rectanglepaths}

Our goal in this section will be to prove Theorem \ref{t.closedrectanglespathscorrespondtoclosedpaths}, the most difficult step in the proof of Theorem B. Our proof of Theorem \ref{t.closedrectanglespathscorrespondtoclosedpaths} will consist in three steps: 
\begin{enumerate}
    \item deforming a closed rectangle path in $\clos{\mathcal{P}_1}$ to a trivial one, by the use of three combinatorial operations that we call \emph{homotopies}
    \item showing that a rectangle path that can be deformed to the trivial one is necessarily closed 
    \item  proving that a rectangle path in $\clos{\mathcal{P}_1}$ can be deformed to a trivial one if and only if its associated rectangle path in $\clos{\mathcal{P}_2}$ can be deformed to a trivial one
\end{enumerate}

Concerning the organization of this section, in Section \ref{s.singularities} we will define the notion of \emph{tangency} for a good polygonal curve and we will state a combinatorial result involving the  number of different tangencies of a simple closed polygonal curve. The previous result will be crucial for us during the proof of Theorem \ref{t.closedrectanglespathscorrespondtoclosedpaths}. Next, in Section \ref{s.closedpathsarehomotopicallytr} we will define the notion of homotopy for rectangle paths in the bifoliated plane up to surgeries. Finally, the rest of this section will be dedicated to the proof of the above mentioned three steps.

\subsection{Singularities of polygonal curves}\label{s.singularities}
Fix $\rho: G\rightarrow \text{Homeo}(\mathcal{P})$ an orientation preserving strong Markovian action, preserving the pair of singular foliations $\mathcal{F}^s$ and $\mathcal{F}^u$ and leaving invariant a strong Markovian family $\mathcal{R}$. We denote by $\Gamma$ the set of boundary periodic points of $\mathcal{R}$, by $\clos{\mathcal{P}}$ the bifoliated plane of $\rho$ up to surgeries on $\Gamma$, by $\clos{\mathcal{R}}$ the lift of $\mathcal{R}$ on $\clos{\mathcal{P}}$, by $\clos{\Gamma}$ the lift of $\Gamma$ on $\clos{\mathcal{P}}$ and by $\clos{\rho}: \clos{G}\rightarrow \text{Homeo}(\clos{\mathcal{P}})$ the lift of the action $\rho$ on $\clos{\mathcal{P}}$. 

Let us immediately remark, that contrary to the case of a plane, Jordan's theorem  does not hold in $\clos{\mathcal{P}}$: 

\begin{rema}\label{r.jordan}
    There exists a simple closed curve $\gamma$ in $\overline{\mathcal{P}}$, such that $\overline{\mathcal{P}}-\gamma$ contains more than two connected components. 
\end{rema}
Indeed, take $s$ an arc in $\overline{\mathcal{P}}$ connecting two different points in $\clos{\Gamma}$. Notice that since $\overline{\mathcal{P}}$ is topologically a plane to which we have added a countable set of points at infinity, namely the points $\clos{\Gamma}$, the arc $s$ disconnects $\overline{\mathcal{P}}$. Consider now an arc $s'$ close to $s$, having the same extremities as $s$ and intersecting nowhere else $s$. By parametrizing $s\cup s'$, we can obtain a simple, closed curve $\gamma$ in $\overline{\mathcal{P}}$ such that $\overline{\mathcal{P}}-\gamma$ consists of two unbounded connected components and one bounded connected component. 

Despite the existence of such curves and since $\clos{\mathcal{P}}-\clos{\Gamma}$ is a topological plane, Jordan's theorem applies for simple, closed curves in $\overline{\mathcal{P}}$ that do not intersect $\clos{\Gamma}$. This justifies the following definition :

\begin{defi}\label{d.simple}
A closed polygonal curve $\gamma:[0,1]\rightarrow \clos{\mathcal{P}}$ will be called \emph{simple} if the function $\gamma$ restricted on $[0,1)$ is injective. 

A simple closed polygonal curve $\gamma$ in $\clos{\mathcal{P}}-\clos{\Gamma}$ defines two complementary regions: a bounded region that we will name the \emph{interior} of $\gamma$ and an unbounded region that we will name the \emph{exterior} of $\gamma$.

\end{defi}
Notice that, since $\clos{\mathcal{P}}$ is homeomorphic to a plane together with a countable set of points at infinity:

\begin{rema}\label{r.gammabarneverinterior}
    The interior of a simple closed polygonal curve in $\clos{\mathcal{P}}-\clos{\Gamma}$ can never contain a point of $\clos{\Gamma}$
\end{rema}

Take $\gamma$ a simple closed polygonal curve in $\clos{\mathcal{P}}-\clos{\Gamma}$, $D$ the interior of $\gamma$ and $U$ a neighborhood of $D\cup \gamma$. Endow $\clos{\mathcal{F}^s}$ and $\clos{\mathcal{F}^u}$ with an orientation (see Remark \ref{r.propertiesfolipbar}). Using the previous orientations, we will define 4 types of tangencies of $\gamma$ to the foliation $\clos{\mathcal{F}^s}$.

\begin{defi}
Consider $\gamma$ as a function from $\mathbb{S}^1$ (endowed with an orientation) to $\clos{\mathcal{P}}$. Fix $s$ a stable segment of $\gamma$. We will denote by $u^{-}$ the unstable segment in $\gamma$ coming right before $s$ (for our choice of orientation of $\mathbb{S}^1$) and by $u^+$ the unstable segment in $\gamma$ coming right after $s$. Using the orientation on $\mathbb{S}^1$, we can naturally orient $u^{-}$ and $u^+$. We will say that $s$ is a \emph{stable tangency} of $\gamma$ if one of the two previous segments is negatively oriented and the other is positively oriented with respect to the orientation of $\clos{\mathcal{F}^{u}}$. 

Furthermore, we will say that a stable tangency $s$ of $\gamma$ is a tangency of type (a) (resp. (d)) if for any small positively oriented unstable segment $u$ starting from a point in $s$ we have that $u\subset \clos{D}$ and if any stable leaf of $\clos{\mathcal{F}^{s}}$ intersecting $u$ intersects also (resp. does not intersect) $u^{-}$ and $u^+$. 

By changing the positively oriented segment $u$ to a negatively oriented segment, we can similarly define a tangency of type (b) or (c). We define in an analogous way \emph{unstable tangencies} of type (a), (b), (c) or (d) for $\gamma$.
    
    
    
    
\begin{figure}
    \centering
    \hspace{-2cm}
    \includegraphics[scale=0.5]{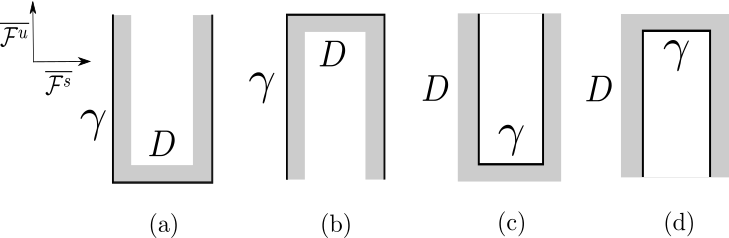}
    \caption{Types of singularities for $\gamma$}
    \label{f.singularities}
\end{figure}

\end{defi}

Endow the plane $\clos{\mathcal{P}}-\clos{\Gamma}$ with a smooth structure. Using the fact that the leaves of $\clos{\mathcal{F}^s}$ are closed inside  $\clos{\mathcal{P}}-\clos{\Gamma}$ (see Remark \ref{r.propertiesfolipbar}), we have that up to conjugating $\clos{\mathcal{F}^s}$ by a homeomorphism, we can assume that $\clos{\mathcal{F}^s}$ is a smooth foliation inside $\clos{\mathcal{P}}-\clos{\Gamma}$ (see more generally Exercise 8.4.9 of \cite{Candel}). Consider now $D'$ a small perturbation of the closed disk $\clos{D}$ inside $U$ such that $\partial D'$ is smooth and in general position with respect to $\clos{\mathcal{F}^s}$. The curve  $\partial D'$ has a finite number of tangencies with $\clos{\mathcal{F}^s}$.   

Take $x$ a point of tangency between $\partial D'$ and $\clos{\mathcal{F}^s}$. We will say that $x$ is an \textit{interior} (resp. \textit{exterior}) \textit{tangency} if any vector based at $x$ and tangent to $\clos{\mathcal{F}^s(x)}$ points towards the interior (resp. exterior) of $D'$.

By considering any non-singular vector field $X$ on $U$ tangent to $\clos{\mathcal{F}^s}$ -such a vector field exists, since $\clos{\mathcal{F}^s}$ is orientable (see Remark \ref{r.propertiesfolipbar})- and by applying  the Poincar\'e Index theorem for the vector field $X$ and  the disk $D'$, we get  that: $$(\text{number of exterior tangencies of }\partial D') - (\text{number of interior  tangencies of }\partial D')=2 $$ 
Therefore, there are always more exterior tangencies in $\partial D'$ than interior ones. 

By eventually changing our initial perturbation $D'$, we may assume that: 
\begin{enumerate}
    \item every stable tangency of type (a) or (b) (resp. (c) or (d)) of $\gamma$ becomes after perturbation an exterior (resp. interior) tangency of $\partial D'$
    \item every stable segment of $\gamma$ that is not a tangency becomes after perturbation transverse to  $\clos{\mathcal{F}^s}$
    \item the tangencies of $\partial D'$ are bijectively associated to stable tangencies of $\gamma$
\end{enumerate}
As a consequence of this, we obtain the following result:
\begin{lemm}\label{l.numbertangencies} Consider the set of stable tangencies of any simple closed polygonal curve $\gamma$ in $\clos{\mathcal{P}}-\clos{\Gamma}$. We have that:
$$(\text{number of tangencies of type (a) or (b)}) - (\text{number of tangencies of type (c) or (d)})=2 $$ 
\end{lemm}
\begin{rema}
Notice that the type (a) (resp. (b),(c) or (d)) of a stable tangency of $\gamma$, depends on our choice of orientation of $\clos{\mathcal{F}^{u}}$. By changing our choice of orientation, tangencies of type (a) (resp. (b), (c), (d)) become tangencies of type (b) (resp. (a), (d), (c)). This is why, the previous lemma remains true independently from our choices of orientation. 
\end{rema}
Naturally, the above lemma is also true for unstable tangencies. Let us finish this section with the following result: 
\begin{lemm}\label{l.canonicalneighbourhoodtangency}
Let $\gamma$ be a simple closed polygonal curve in $\clos{\mathcal{P}}-\clos{\Gamma}$, $s$ a stable tangency of $\gamma$ of type (a) or (b) and  $D$ the interior of $\gamma$. There exists a unique rectangle $R$ contained in $D\cup \gamma$ such that
\begin{itemize}
   \item $s$ is a stable boundary component of $R$
   \item the unstable boundaries of $R$ are contained in $\gamma$ 
   \item $R$ is maximal for the previous properties
\end{itemize} 
\end{lemm}
\begin{proof}
Consider $\gamma$ as a function from $\mathbb{S}^1$ (endowed with an orientation) to $\clos{\mathcal{P}}-\clos{\Gamma}$ and $u^{-}$ (resp. $u^{+}$) the unstable segment in $\gamma$ coming right before (resp. after) $s$. Denote by $I^-$ (resp. $I^+ $) the set of points $x$ in $u^-$ (resp. $u^+$) for which there exists a stable segment inside $D \cup \gamma$ going from $x$ to $u^+$ (resp. $u^-$). Denote also by $x_0^-$ (resp. $x_0^+$) the unique point in $s\cap u^-$ (resp. $s\cap u^+$). Notice that $I^-$ (resp. $I^+$) contains a neighbordhood of $x_0^-$ (resp. $x_0^+$) in $u^-$ (resp. $u^+$), since $s$ intersects both $u^-$ and $u^+$. Let us also remark that, thanks to Remark \ref{r.gammabarneverinterior}, if $x\in I^-$, then all points of $u^-$ between $x$ and $x^-_0$ belong in $I^-$, hence $I^-$ is a segment. The same applies for $I^+$. 

Next, by the definitions of $I^-$ and $I^+$, for every $x\in I^-$ the stable segment in $ D\cup \gamma$ going from $x$ to $u^+$ intersects $u^+$ along a point of $I^+$ and vice versa. Let us now show that $I^-$ is a closed segment in $u^-$. Assume that $I^-$ is of the form $[x_0^-, X^-)\subset u^-$. By our previous arguments, this implies that $I^+$ is also of the form $[x_0^+,X^+)$. 

If the stable leaf of $X^-$ intersects $u^+$, since it can be accumulated by stable segments in $D\cup\gamma$ going from $u^-$ to $u^+$, then there exists a stable segment in $D\cup \gamma$ going from $X^-$ to $u^+$ and thus $X^-\in I^-$. We deduce that the stable leaf of $X^-$ (resp. $X^+$) does not intersect $u^+$ (resp. $u^-$). Moreover, the stable leaves of $X^-$ and $X^+$ can be accumulated from one side by segments going from $u^-$ to $u^+$; thus the previous two leaves are non-separated stable leaves. 

Denote by $\clos{\pi}$ the projection from $\clos{\mathcal{P}}$ to $\mathcal{P}$ and by $s_x$ the stable segment in $D\cup \gamma$ going from $x \in I^-$ to $u^+$. By applying Corollary 2.19 of \cite{nontransitiveanosovlike}\footnote{We remark here that Corollary 2.19 of \cite{nontransitiveanosovlike} is proven for non-transitive Anosov-like actions. The definition of such an action in \cite{nontransitiveanosovlike} consists of 5 axioms, denoted by A1-A5. The axiom A5 is not used in the proof of Corollary 2.19, the axioms A1, A2, A4 are included in our definition of strong Markovian action and finally A3 follows from Remark \ref{r.singularitiesareboundaries}} for the non-separated stable leaves in $\mathcal{P}$ containing $\clos{\pi}(\clos{\mathcal{F}^s}(X^+))$ and $\clos{\pi}(\clos{\mathcal{F}^s}(X^-))$, we get that there exists an unstable leaf $L\in \clos{\mathcal{F}^u}$ such that the set $L\cap \underset{x\in I^-}{\cup}s_x\subset D\cup \gamma$ contains the entire  unstable separatrix of some point in $L$. Since  $D\cup \gamma $ is compact and does not intersect $\clos{\Gamma}$, the previous statement contradicts the fact that the leaves of $\clos{\mathcal{F}^u}$ are properly embedded in $\clos{\mathcal{P}}-\clos{\Gamma}$ (see Remark \ref{r.propertiesfolipbar}). We deduce that $I^-$ and $I^+$ are closed. The segments in $D\cup \gamma$ going from $I^-$ to $I^+$ form the  rectangle with the desired properties.

\end{proof} 
The rectangle constructed in the previous lemma forms a canonical neighborhood of the tangency $s$ inside $D\cup \gamma$. We will call this neighborhood the \textit{domain} of $s$. Furthermore, if the domain of $s$ covers completely either $u^-$ or $u^+$, we will call it a \textit{complete domain}. In any other case, we will call it an \emph{incomplete domain}. By our proof of Lemma \ref{l.canonicalneighbourhoodtangency} we can deduce that: 
\begin{rema}\label{r.incompletedomain}
 If the domain $\mathcal{D}$ of $s$ is  incomplete, there exists a  stable segment of $\gamma$ intersecting the interior of the stable boundary component of $\mathcal{D}$ that is not $s$.  
\end{rema}

\subsection{Three types of homotopies of rectangle paths and good polygonal curves}\label{s.closedpathsarehomotopicallytr}
Fix $\rho: G\rightarrow \text{Homeo}(\mathcal{P})$ an orientation preserving strong Markovian action, preserving the pair of singular foliations $\mathcal{F}^s$ and $\mathcal{F}^u$ and leaving invariant a strong Markovian family $\mathcal{R}$. We denote by $\Gamma$ the boundary periodic point of $\mathcal{R}$, by $\clos{\mathcal{P}}$ the bifoliated plane of $\rho$ up to surgeries on $\Gamma$, by $\clos{\mathcal{R}}$ the lift of $\mathcal{R}$ on $\clos{\mathcal{P}}$, by $\clos{\Gamma}$ the lift of $\Gamma$ on $\clos{\mathcal{P}}$ and by $\clos{\rho}: \clos{G}\rightarrow \text{Homeo}(\clos{\mathcal{P}})$ the lift of the action $\rho$ on $\clos{\mathcal{P}}$. Consider $\clos{\mathcal{F}^s}$ and $\clos{\mathcal{F}^u}$ the stable and unstable foliations of $\clos{\rho}$ endowed with an orientation.

Fix for the rest of this section $R_0,...,R_n$ a rectangle path in $\clos{\mathcal{R}}$.
\begin{defi}
We will say that two rectangle paths are \emph{homotopic by a homotopy of type $\mathcal{A}$} if they are equal. 

Two good polygonal curves $\gamma,\delta:[0,1]\rightarrow \clos{\mathcal{P}}$ whose associated rectangle paths starting from $R_0$ are equal to  $R_0,...,R_n$ will be called \emph{homotopic by a homotopy of type $\mathcal{A}$} (\emph{relatively to $R_0,...,R_n$}).

\end{defi}

\begin{conv}
 By an abuse of language, whenever the context makes it clear, we will omit the rectangle paths relatively to which we will perform homotopies of good polygonal curves.  
\end{conv}

\begin{figure}[h!]
\centering
\includegraphics[scale=0.75]{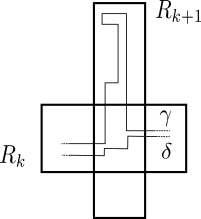}
\caption{$\gamma$ and $\delta$ are homotopic by a homotopy of type $\mathcal{B}$}
\label{f.homotopyoftypeB}
\end{figure}
\begin{defi}\label{d.homotopyB}
Consider $R_0,...,R_n$ a rectangle path such that  there exists $k\in \llbracket 0,n-2\rrbracket$ for which $R_k=R_{k+2}$. We are going to say that the rectangle paths $R_0,..,R_k,R_{k+3},R_{k+4},...$ and $R_0,...,R_n$ are \emph{homotopic by a homotopy of type $\mathcal{B}$}. 

Consider $\gamma$, $\delta$ two good polygonal curves whose associated rectangle paths starting from $R_0$ are equal to $R_0,...,R_n$ and $R_0,..,R_k,R_{k+3},R_{k+4},...,R_n$ respectively (see Figure \ref{f.homotopyoftypeB}). We will say that $\gamma$ and $\delta$ are \emph{homotopic by a homotopy of type $\mathcal{B}$} (\emph{relatively to $R_0,...,R_n$ and $R_0,..,R_k,R_{k+3},$ $ R_{k+4},...,R_n$}).  
\end{defi}
Homotopies of type $\mathcal{A}$ or type $\mathcal{B}$ correspond to ``movements" of good polygonal curves in the interior of the rectangles of $\clos{\mathcal{R}}$. As it turns out, the previous two elementary movements are not enough to deform any closed rectangle path to a trivial one. This is why we will now define a third type of homotopy allowing good polygonal curves to cross the boundary of some rectangle in $\clos{\mathcal{R}}$ and more specifically to cross a boundary arc point in $\overline{\mathcal{P}}$. We will name this homotopy, \emph{homotopy of type $\mathcal{C}$}. In order to introduce this homotopy, we will first need to define the notion of \textit{cycle} around a boundary arc point.

\begin{figure}[h!]
\centering
\includegraphics[scale=1]{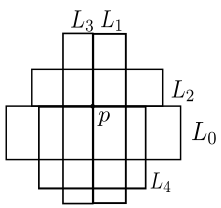}
\caption{}
\label{f.cycle}
\end{figure}

Take $p\in \clos{\mathcal{P}}$ a boundary arc point of $\clos{\mathcal{R}}$ and $L_0\in \clos{\mathcal{R}}$ such that $L_0$ contains $p$ and two germs of quadrants of $p$ (i.e. two small neighborhoods of $p$ inside two distinct quadrants of $p$). Assume without any loss of generality that $p\in \partial^s L_0$ (see Figure \ref{f.cycle}) and denote by $r$ the stable boundary component of $L_0$ containing $p$. By 
Lemmas \ref{l.crossingrectanglesnoperiodicpoints} and \ref{l.crossingrectangleswithperiodicpoints} applied for $\clos{\mathcal{R}}$ (see Remark \ref{r.equivdefi}), there exists $L_1$ an $r$-crossing predecessor of $L$ such that $p\in \partial{L_0}\cap \partial{L_1}$. 

Consider $r_1$ the unstable boundary component of $L_1$ containing $p$. Again by 
Lemmas \ref{l.crossingrectanglesnoperiodicpoints} and \ref{l.crossingrectangleswithperiodicpoints}, there exists $L_2$ an $r_1$-crossing successor of $L_1$ such that $p\in \partial{L_2}\cap \partial{L_1}$ and $\overset{\circ}{L_2}\cap \overset{\circ}{L_0}=\emptyset$. By the same exact procedure, we construct $L_3$ and $L_4$ such that $\overset{\circ}{L_3}\cap \overset{\circ}{L_1}=\emptyset$ and $\overset{\circ}{L_4}\cap \overset{\circ}{L_2}=\emptyset$.

\begin{defi}\label{d.cyclearcpoint}
 We are going to call $(L_0,L_1,L_2,L_3,L_4)$ a \emph{cycle} around $p$ starting from $L_0$.
\end{defi}

\begin{rema}\label{r.existenceofcycle}
For any boundary arc point $p\in \clos{\mathcal{P}}$, there exists a cycle around $p$.\end{rema}

Indeed, by our discussion prior to Definition \ref{d.cyclearcpoint}, it suffices to show that there exists $L \in\clos{\mathcal{R}}$ containing $p$ and the germs of two quadrants of $p$. Take $L\in\clos{\mathcal{R}}$ containing $p$. Since $p$ is a boundary arc point, $ p\in \partial L$. If $p$ is not a corner point of $L$, then $L$ contains the germs of two quadrants of $p$. If $p$ is a corner point of $L$, let us name $s$ the stable boundary of $L$ in which it is contained and $L'$ the $s$-crossing predecessor of $L$ that contains $p$. The rectangle $L'$ contains the germs of two quadrants of $p$.

The set of boundary arc points of $\clos{\mathcal{R}}$, by definition, is disjoint from $\clos{\Gamma}$; hence boundary arc points are  regular points of the foliations $\clos{\mathcal{F}^{s}}$ and $\clos{\mathcal{F}^{u}}$. Take $p$ a boundary arc point of $\clos{\mathcal{R}}$ and $(L_0,L_1,L_2,L_3,L_4)$ a cycle around of $p$. Using the orientations of $\clos{\mathcal{F}^{s}}$ and $\clos{\mathcal{F}^{u}}$ fixed in the beginning of this section, denote by $(+,+), (+,-), (-,+), (-,-)$ the 4 quadrants of $p$, where $(+,+)$ (resp. $(+,-), (-,+), (-,-)$) denotes the quadrant of $p$ bounded by the positive (resp. positive, negative, negative) stable separatrix of $p$ and the positive (resp. negative, positive, negative) unstable separatrix of $p$.  Denote also by  $G_{(\epsilon,\epsilon')}$ a germ of the $(\epsilon,\epsilon')$ quadrant of $p$, where $\epsilon,\epsilon'\in \lbrace -,+ \rbrace$. By the definition of a cycle, up to cyclic permutation of the indices of the $L_i$, we have one of the two following statements: 

\vspace{0.2cm}

\begin{minipage}{0.5\linewidth}
\center{\hspace{-0.7cm}\underline{Positive cycle}}
\begin{itemize}
    \item $\overset{\circ}{L_0}\cap G_{(-,-)}\neq \emptyset$, $\overset{\circ}{L_0}\cap G_{(+,-)}\neq \emptyset$ 
    \item $\overset{\circ}{L_1}\cap G_{(+,-)}\neq \emptyset$, $\overset{\circ}{L_1}\cap G_{(+,+)}\neq \emptyset$
    \item $\overset{\circ}{L_2}\cap G_{(+,+)}\neq \emptyset$, $\overset{\circ}{L_2}\cap G_{(-,+)}\neq \emptyset$
    \item $\overset{\circ}{L_3}\cap G_{(-,+)}\neq \emptyset$, $\overset{\circ}{L_3}\cap G_{(-,-)}\neq \emptyset$
    \item $\overset{\circ}{L_4}\cap G_{(-,-)}\neq \emptyset$, $\overset{\circ}{L_4}\cap G_{(+,-)}\neq \emptyset$
\end{itemize} 
\end{minipage}
\hspace{-0.5cm}\begin{minipage}{0.5\linewidth}
\center{\hspace{-0.7cm}\underline{Negative cycle}}
\begin{itemize}
    \item $\overset{\circ}{L_0}\cap G_{(+,-)}\neq \emptyset$, $\overset{\circ}{L_0}\cap G_{(-,-)}\neq \emptyset$ 
    \item $\overset{\circ}{L_1}\cap G_{(-,-)}\neq \emptyset$, $\overset{\circ}{L_1}\cap G_{(-,+)}\neq \emptyset$
    \item $\overset{\circ}{L_2}\cap G_{(-,+)}\neq \emptyset$, $\overset{\circ}{L_2}\cap G_{(+,+)}\neq \emptyset$
    \item $\overset{\circ}{L_3}\cap G_{(+,+)}\neq \emptyset$, $\overset{\circ}{L_3}\cap G_{(+,-)}\neq \emptyset$
    \item $\overset{\circ}{L_4}\cap G_{(+,-)}\neq \emptyset$, $\overset{\circ}{L_4}\cap G_{(-,-)}\neq \emptyset$
\end{itemize}
\end{minipage}

\begin{figure}[h!]
\centering
\includegraphics[scale=0.6]{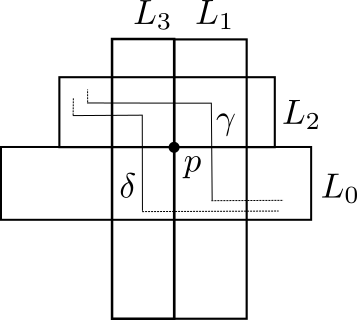}
\caption{}
\label{f.homotopyC}
\end{figure}

\vspace{0.2cm}
In the first case, we will call the cycle $(L_0,L_1,L_2,L_3,L_4)$ \textit{positive} and in the second one \textit{negative}. Notice in particular, that for every cycle $(L_0,L_1,L_2,L_3,L_4)$ around $p$, we have that $\overset{\circ}{L_4}\cap \overset{\circ}{L_0}\neq \emptyset$. 
\begin{lemm}\label{l.twocycles}
Take $p\in \clos{\mathcal{P}}$ a boundary arc point of $\clos{\mathcal{R}}$ and $L_0 \in \clos{\mathcal{R}}$ such that $L_0$ contains $p$ and the germs of two quadrants of $p$. There exist exactly two cycles around $p$ starting from $L_0$: one positive and one negative. 
\end{lemm}
\begin{proof}
Without any loss of generality, we can assume that $p$ is contained in a stable boundary component of $L_0$, say $s$, and that a germ of the $(-,-)$ and $(+,-)$ quadrants of $p$ is contained in $L_0$. By Remark \ref{r.caracterisationarcpoints} and  
Lemmas \ref{l.crossingrectanglesnoperiodicpoints}, \ref{l.crossingrectangleswithperiodicpoints} applied for $\clos{\mathcal{R}}$ (see also Remark \ref{r.equivdefi}), since $p$ is a boundary arc point and any two distinct $s$-crossing predecessors of $L_0$ intersect only along their boundaries, there exist exactly two $s$-crossing predecessors $L_1,L_1'$ of $L_0$ containing $p$. Assume without any loss of generality that $L_1$ contains a germ of the $(+,+)$ and $(+,-)$ quadrants of $p$ and $L_1'$ a germ of the $(-,+)$ and $(-,-)$ quadrants of $p$. Since $L_1$ is an $s$-crossing predecessor of $L_0$, the point $p$ is contained in the interior of some unstable boundary component of $L_1$. Therefore, by our previous argument the rectangle $L_2$ of Definition \ref{d.cyclearcpoint} is uniquely defined. Similarly, the rectangles $L_3$ and $L_4$ of Definition \ref{d.cyclearcpoint} are uniquely defined. The same argument applies for $L_1'$. We thus obtain exactly two cycles around $p$ starting from $L_0$: one positive and one negative. 
\end{proof}
\begin{lemm}
Take $p\in \clos{\mathcal{P}}$ a boundary arc point of $\clos{\mathcal{R}}$ and $L_0\in \clos{\mathcal{R}}$ containing the germs of two quadrants of $p$. If  $(L_0,L_1,L_2,L_3,L_4)$ is a cycle around $p$, then  $(L_0,L_3,L_2,L_1,L_4)$ is the other cycle around $p$ starting from $L_0$. 
\end{lemm}
\begin{proof}
Let us assume without any loss of generality that $p \in \partial^s L_0$. Denote by $s$ the stable boundary component of $L_0$ containing $p$. Assuming that $(L_0,L_1,L_2,L_3,L_4)$ is a cycle around $p$, let us show that $L_3$ is an $s$-crossing predecessor of $L_0$. 

If $s'$ is the stable boundary component of $L_2$ containing $p$, then since $(L_0,L_1,L_2,L_3,L_4)$ is a cycle, we have that $L_3$ is an $s'$-crossing predecessor of $L_2$. By  Remark \ref{r.caracterisationarcpoints} applied for $\clos{\mathcal{R}}$, the unstable boundary of $L_3$ intersects $s'$ along boundary arc points. By the Markovian intersection property, $L_3\cap L_0$ is a vertical subrectangle of $L_0$ (see Figure \ref{f.cycle}). Since $\partial^u L_3\cap s$ consists of two boundary arc points, we deduce from Remark \ref{r.caracterisationarcpoints} that $L_3$ is an $s$-crossing predecessor of $L_0$. We show in the exact same way, that $L_2$ is a crossing successor of $L_3$, that $L_1$ is a crossing predecessor of $L_2$ and that $L_4$ is a crossing successor of $L_1$. Finally, since $(L_0,L_1,L_2,L_3,L_4)$ is a cycle, we have that $\overset{\circ}{L_2}\cap \overset{\circ}{L_0}=\emptyset$, $\overset{\circ}{L_3}\cap \overset{\circ}{L_1}=\emptyset$ and $\overset{\circ}{L_4}\cap \overset{\circ}{L_2}=\emptyset$. By Lemma \ref{l.twocycles}, we conclude that $(L_0,L_3,L_2,L_1,L_4)$ is the second cycle around $p$ starting from $L_0$. 
\end{proof}
\begin{lemm}
    Take $p\in \clos{\mathcal{P}}$ a boundary arc point and $L_0,L'_0\in \clos{\mathcal{R}}$ two rectangles  containing the germs of the same two quadrants of $p$. If  $(L_0,L_1,L_2,L_3,L_4)$ is a positive (resp. negative) cycle around $p$, then  $(L'_0,L_1,L_2,L_3,L_4)$ is also a positive (resp. negative) cycle around $p$.
\end{lemm}
\begin{proof}
    Indeed, assume without any loss of generality that $(L_0,L_1,L_2,L_3,L_4)$ is a positive cycle around $p$ and that $L_0, L'_0$ contain the germs of the $(-,-)$ and $(+,-)$ quadrants of $p$. 
    
    Denote by $r$ and $r'$ the stable boundary components of $L_0$ and $L'_0$ containing $p$. By definition of a positive cycle, $L_1$ contains the germs of the $(+,-)$ and $(+,+)$ quadrants of $p$ and is an $r$-crossing predecessor of $L_0$. By the Markovian intersection axiom, $L_1\cap r=L_1\cap r'$ and $L_1\cap L'_0$ is a non-trivial vertical subrectangle of $L'_0$. Furthermore, thanks to Remark \ref{r.caracterisationarcpoints}, the two extremities of $L_1\cap r$ consist of two boundary arc points and by the same remark applied to $L_1\cap r'$, we get that $L_1$ is an $r'$-crossing predecessor of $L'_0$. By using the same construction as the one that was given prior to Definition \ref{d.cyclearcpoint}, we get that   $(L'_0,L_1,L_2,L_3,L_4)$ is a cycle around $p$ visiting the quadrants of $p$ in the same order as $(L_0,L_1,L_2,L_3,L_4)$, which gives us the desired result.  
\end{proof}
The above lemma implies that the rectangles $L_1,L_2,L_3,L_4$ are canonically associated to $p$ independently of the choice of $L_0$. In particular, this means that $(L_4,L_1,L_2,L_3,L_4)$ defines a cycle around $p$. 
\begin{defi}\label{d.generalisedrectpaths}
Consider $L_0,...,L_n$ a sequence of rectangles in $\clos{\mathcal{R}}$ such that for every $i$, the rectangle $L_{i+1}$ is a successor or predecessor of some generation of $L_i$. We will call $L_0,...,L_n$ a \emph{generalised rectangle path}. By Lemma \ref{l.npredecessor} applied for $\clos{\mathcal{R}}$ (see also Remark \ref{r.equivdefi}), for every $i$ there exists a unique monotonous rectangle path (see Definition \ref{d.rectanglepath})  $L_i=L_{i0},L_{i1},..., L_{is_i},L_{i(s_i+1)}=L_{i+1}$ going from $L_i$ to $L_{i+1}$. 
We are going to call $$L_0,L_{01},...,L_{0s_0},L_1, L_{11},...,L_{1s_1},L_2,...., L_{n-1},L_{(n-1)1},...,L_{(n-1)s_{n-1}}, L_{n}$$
the rectangle path \emph{associated} to  $L_0,...,L_n$.
\end{defi}
\begin{defi}\label{d.homotopyC}
Consider $R_0,...,R_n$ a rectangle path, $p\in \clos{\mathcal{P}}$ a boundary arc point, $L_0\in \clos{\mathcal{R}}$ such that $L_0$ contains $p$ and the germs of two quadrants of $p$. Denote by $(L_0,L_1,L_2,L_3,L_4)$ and $(L_0,L_3,L_2,L_1,L_4)$ the two cycles around $p$ starting from $L_0$. Consider $(L_0,L_1...,L_k)$ the first $k+1$ terms of the first cycle with $k\in \llbracket 1, 4\rrbracket $ and $(L_0,L_3...,L_k)$ the part of the second cycle starting from $L_0$ and ending at $L_k$. 

Following the notations of Definition \ref{d.generalisedrectpaths}, we can associate to $L_0,L_1,...,L_k$ and $L_0,L_3,...,L_k$ the rectangle paths $L_0,L_{01},...,L_{0s_0},L_1, L_{11},...,L_{1s_1},L_2,...,L_k$ and \newline{}$L_0,L'_{01},...,L'_{0s'_0},L_3, L'_{31},...,L'_{3s'_3},...,L_k$ respectively. Assume that the first rectangle path is contained in $R_0,...,R_n$. In other words, assume that $R_0,...,R_n$ is of the form $$ R_0,...,R_m,L_0,L_{01},...,L_{0s_0},L_1, L_{11},...,L_{1s_1},..., L_k,R_l,R_{l+1},...,R_n$$  
Consider now the rectangle path $R'_0,...,R'_{n'}$ defined as follows 
$$R_0,...,R_m,L_0,L'_{01},...,L'_{0s'_0},L_3, L'_{31},...,L'_{3s'_3},...,L_k,R_l,R_{l+1},...,R_n$$

We are going to say that the rectangle paths $R_0,...,R_n$ and $R'_0,...,R'_{n'}$ are \emph{homotopic by a homotopy of type $\mathcal{C}$}.

Take $\gamma$ and $\delta$ two good polygonal curves whose associated rectangle paths starting from $R_0$ are  $R_0,...,R_n$ and $R'_0,...,R'_{n'}$ respectively (see Figure \ref{f.homotopyC}). We will say that $\gamma$ is homotopic to $\delta$ by a \emph{homotopy of type $\mathcal{C}$} (\emph{relatively to $R_0,...,R_n$ and $R'_0,...,R'_{n'}$}). 
\end{defi}


\begin{defi}
Take $P:=r_0,r_1,...,r_n$ and $P':=r_0,r'_1,...,r'_m$ two rectangle paths. We are going to say that $P$ and $P'$ are \emph{homotopic} if there exists $P_0:=P, P_1, ..., P_s:=P'$ a finite sequence of rectangle paths such that  $P_i$ is homotopic to $P_{i+1}$ by a homotopy of type $\mathcal{A}$, $\mathcal{B}$ or $\mathcal{C}$. 
\end{defi}

 Let us remark that thanks to our definition of homotopies of type $\mathcal{A},\mathcal{B} $ and $\mathcal{C}$, any two homotopic rectangle paths in $\clos{\mathcal{R}}$ begin and end by the same rectangles. In particular, a closed rectangle path in $\clos{\mathcal{R}}$ is only homotopic to closed rectangle paths. Therefore, a rectangle path homotopic to the trivial rectangle path is necessarily closed. It turns out that the converse is also true:  

\begin{theorem}\label{t.homotopictotrivialpath}
Any closed rectangle path $r_0,r_1,...,r_n$ in $\clos{\mathcal{P}}$ is homotopic to the trivial rectangle path $r_0$

\end{theorem}
The above theorem constitutes the key argument in the proof of Theorem \ref{t.closedrectanglespathscorrespondtoclosedpaths}, which will be used to prove Theorem B in Section \ref{s.mainresult}. 
\subsection{Theorem \ref{t.homotopictotrivialpath} implies Theorem \ref{t.closedrectanglespathscorrespondtoclosedpaths}}
\label{s.proofoftheoremcorrespondancehom}

Let $\rho_1:G_1 \rightarrow \text{Homeo}(\mathcal{P}_1)$, $\rho_2: G_2 \rightarrow  \text{Homeo}(\mathcal{P}_2)$ be two strong orientation preserving Markovian actions acting on the planes $\mathcal{P}_1$ and $\mathcal{P}_2$, endowed with their stable and unstable foliations $\mathcal{F}^{s,u}_1$ and $\mathcal{F}^{s,u}_2$. Assume that $\rho_1,\rho_2$ preserve two strong Markovian families $\mathcal{R}_1$ and $\mathcal{R}_2$, whose associated classes of geometric types coincide. 

Denote by $\Gamma_1,\Gamma_2$ the set of boundary periodic points of $\mathcal{R}_1, \mathcal{R}_2$, by $\clos{\mathcal{P}_{1}}$, $\clos{\mathcal{P}_{2}}$ the bifoliated planes of $\rho_1$, $\rho_2$ up to surgeries on $\Gamma_1,\Gamma_2$, by $\clos{\mathcal{R}_1}$ and $\clos{\mathcal{R}_2}$ the lifts of $\mathcal{R}_1$ and $\mathcal{R}_2$ on $\clos{\mathcal{P}_{1}}$, $\clos{\mathcal{P}_{2}}$, by $\clos{\Gamma_1},\clos{\Gamma_2}$ the lifts of $\Gamma_1,\Gamma_2$ on $\clos{\mathcal{P}_{1}}$, $\clos{\mathcal{P}_{2}}$ and by $\clos{\rho_1}:\clos{G_1}\rightarrow  \text{Homeo}(\clos{\mathcal{P}_{1}})$, $\clos{\rho_2}:\clos{G_2}\rightarrow  \text{Homeo}(\clos{\mathcal{P}_{2}})$ the lifts of $\rho_1$, $\rho_2$ on $\clos{\mathcal{P}_{1}}$, $\clos{\mathcal{P}_{2}}$. By Proposition  \ref{p.liftedmarkovfamiliessamegeomtype} the Markovian families $\clos{\mathcal{R}_1}$ and $\clos{\mathcal{R}_2}$ are associated to the same class of geometric types and thanks to the realisability lemma, we can choose a set of representatives of every orbit in $\clos{\mathcal{R}_1}$ and $\clos{\mathcal{R}_2}$ and we can orient the stable and unstable foliations in   $\clos{\mathcal{P}_{1}}$, $\clos{\mathcal{P}_{2}}$ so that the geometric types associated to $\clos{\mathcal{R}_1}$ and $\clos{\mathcal{R}_2}$ for the previous choices of orientations and representatives coincide and are equal to $$\mathcal{G}=(n,(h_i)_{i \in \llbracket 1,n \rrbracket}, (v_i)_{i\in \llbracket 1,n \rrbracket}, \mathcal{H}, \mathcal{V},\phi, u)$$

Fix $r_0^1\in \clos{\mathcal{R}_1}$ and $r_0^2\in \clos{\mathcal{R}_2}$ two of the previously chosen representatives such that $r_0^1$ and $r_0^2$ are of the same type. For the rest of this section, $r_0^1$ and $r_0^2$ will be our reference rectangles in $\clos{\mathcal{R}_1}$ and $\clos{\mathcal{R}_2}$.

Assuming that Theorem \ref{t.homotopictotrivialpath} stands, let us show that:
\begin{theorem*}[Theorem \ref{t.closedrectanglespathscorrespondtoclosedpaths}]
Any closed centered rectangle path in $\clos{\mathcal{R}_1}$ is associated (see Definition  \ref{d.associatingrectanglepaths}) to a closed centered rectangle path in $\clos{\mathcal{R}_2}$.
\end{theorem*} 
\begin{proof}
It suffices to prove that two centered and homotopic rectangle paths in $\clos{\mathcal{P}_1}$ correspond to centered and homotopic rectangle paths in $\clos{\mathcal{P}_2}$. Indeed, let $r^1_0,...,r^1_n$ be a closed and centered rectangle path in $\clos{\mathcal{P}_1}$ and $r^2_0,...,r^2_n$ its associated centered rectangle path in $\clos{\mathcal{P}_2}$. By Theorem \ref{t.homotopictotrivialpath}, $r^1_0,...,r^1_n$ is homotopic to the trivial rectangle path $r_0^1$. If homotopic rectangle paths in $\clos{\mathcal{P}_1}$ correspond to homotopic rectangle paths in $\clos{\mathcal{P}_2}$, then $r^2_0,...,r^2_n$ is homotopic to $r^2_0$. By our discussion at the end of the last section, this implies that  $r^2_0,...,r^2_n$ is closed, which gives us the desired result. Let us now show that homotopic rectangle paths in $\clos{\mathcal{P}_1}$ correspond to homotopic rectangle paths in $\clos{\mathcal{P}_2}$.

Since any homotopy of rectangle paths can be described as a sequence of homotopies of type $\mathcal{A},\mathcal{B}$ or $\mathcal{C}$, it suffices to show that homotopies of type $\mathcal{A}$ (resp. $\mathcal{B}$, $\mathcal{C}$) in $\clos{\mathcal{P}_1}$ correspond to homotopies of type $\mathcal{A}$ (resp. $\mathcal{B}$, $\mathcal{C}$) in $\clos{\mathcal{P}_2}$. The result is trivially true for homotopies of type $\mathcal{A}$.  

Consider now $C_1,C_2$ two centered rectangle paths in $\clos{\mathcal{P}_1}$ that are homotopic by a homotopy of type $\mathcal{B}$. By Definition \ref{d.homotopyB}, the paths $C_1,C_2$ are of the form  $r^1_0,...,r^1_k,r^1_{k+1},r^1_k,r^1_{k+2},...r^1_n$ and $r^1_0,...,r^1_k,r^1_{k+2},...,r^1_n$. By our construction in Definition \ref{d.associatingrectanglepaths}, it is not difficult to see that $C_1,C_2$ correspond in $\clos{\mathcal{P}_2}$ to two  rectangle paths of the form $r^2_0,...,r^2_k,r^2_{k+1},r^2_k,r^2_{k+2},...r^2_n$ and  $r^2_0,...,r^2_k,r^2_{k+2},...,r^2_n$. The previous two rectangle paths are homotopic by a homotopy of type $\mathcal{B}$ in $\clos{\mathcal{P}_2}$. 

Finally, consider $C_1,C_2$ two centered rectangle paths in $\clos{\mathcal{P}_1}$ that are homotopic by a homotopy of type $\mathcal{C}$. Following the notations of Definition \ref{d.homotopyC}, we may assume that $C_1$ and $C_2$ are respectively of the form: $$r^1_0,...,r^1_k,L^1_0,L^1_{01},...,L^1_{0s_0},L^1_1, L^1_{11},...,L^1_{1s_1},..., L^1_k,r^1_l,r^1_{l+1},...,r^1_n$$  
 $$r^1_0,...,r^1_k,L^1_0,L'^1_{01},...,L'^1_{0s'_0},L^1_3, L'^1_{31},...,L'^1_{3s'_3},...,L^1_k,r^1_l,r^1_{l+1},...,r^1_n$$
where $(L^1_0,L^1_1,L^1_2,L^1_3,L^1_4)$ and $(L^1_0,L^1_3,L^1_2,L^1_1,L^1_4)$ are the two cycles starting from $L^1_0$ around a boundary arc point $p^1\in L^1_0 - \{\text{corners of $L^1_0$}\} $, $k\in \llbracket 1,4 \rrbracket$,  $L^1_0,L^1_{01},...$ $,L^1_{0s_0},L^1_1, L^1_{11},$ $...,L^1_{1s_1},..., L^1_k$ is the rectangle path associated to the generalized rectangle path $(L^1_0,L^1_1,...,L^1_k)$ and $L^1_0,L'^1_{01},...,L'^1_{0s'_0},L^1_3, L'^1_{31},...,L'^1_{3s'_3},...,L^1_k$ is the rectangle path associated to the generalized rectangle path $(L^1_0,L^1_3,...,L^1_k)$. 

Assume without any loss of generality that $L_0^1$ contains a germ of the $(-,-)$ and $(+,-)$ quadrants of $p^1$ and that $(L^1_0,L^1_1,L^1_2,L^1_3,L^1_4)$ is a positive cycle around $p^1$ (see Figure \ref{f.cycle}). Denote by $$r^2_0,...,r^2_k,L^2_0,L^2_{01},...,L^2_{0s_0},L^2_1, L^2_{11},...,L^2_{1s_1},..., L^2_k,r^2_l,r^2_{l+1},...,r^2_n$$ 
$$r^2_0,...,r^2_k,L^2_0,L'^2_{01},...,L'^2_{0s'_0},L^2_3, L'^2_{31},...,L'^2_{3s'_3},...,L^2_k,r^2_l,r^2_{l+1},...,r^2_n$$
the centered rectangle paths in $\clos{\mathcal{P}_2}$ associated to $C_1$ and $C_2$ respectively. We would like to show that there exists $p^2$ a boundary arc point in $\clos{\mathcal{P}_2}$ such that  $(L^2_0,L^2_1..., L^2_k)$ \big(resp. $(L^2_0,L^2_3..., L^2_k)$\big) is part of the positive (resp. negative) cycle of $p^2$ starting from $L^2_0$. 

\textit{\underline{$(L^2_0,L^2_1..., L^2_k)$ is part of a positive cycle in $\clos{\mathcal{P}_2}$} }

Let $s^1_0$ be the uppermost stable boundary component of $L_0^1$. Notice that by our previous hypotheses $p^1\in s^1_0$ and also by our definition of a positive cycle, $L_1^1$ is an $s^1_0$-crossing predecessor of $L_0^1$ containing $p^1$ in its leftmost unstable boundary. Let us show that $L^2_1$ is an $s^2_0$-crossing predecessor of $L_0^2$, where $s^2_0$ is the uppermost stable boundary component of $L_0^2$ (we use here our choice of orientation of the unstable foliation in $\clos{\mathcal{P}_2}$ fixed in the beginning of this section). 

Indeed, using the geometric interpretation of $\mathcal{G}$, denote by $R^{\mathcal{G}}$ the rectangle in the geometric type $\mathcal{G}$ associated to $R\in \clos{\mathcal{R}_1}$ or $R\in \clos{\mathcal{R}_2}$ (see Definition \ref{d.geometrictypemarkovfamily}). By our definition of an $s_0^1$-crossing predecessor and thanks to Definition \ref{d.generalisedrectpaths}, $L^1_{01}$ is a predecessor of $L^1_0$ that does not cross $s_0^1$. Recall that every successor of $L^1_{01}$ corresponds to a unique horizontal subrectangle of $(L^1_{01})^{\mathcal{G}}$. By our previous arguments and depending on our choice of representatives (which determines the vertical and horizontal orientations of the rectangles of $\mathcal{G}$), $L^1_0$ corresponds to the highermost or lowermost horizontal subrectangle of $(L^1_{01})^G$ and $s_0^1\cap L^1_{01}$ to the highermost or lowermost respectively stable boundary component of $(L^1_{01})^{\mathcal{G}}$. Therefore, by our definition of the associated rectangle path, $L^2_{01}$ corresponds to a predecessor of $L^2_0$ that does not cross $s_0^2$. Similarly, for every $i\in \llbracket 2, s(0)\rrbracket$, we get that $L^2_{0i}$ is  a predecessor of $L^2_{0(i-1)}$ that does not cross $s_0^2\cap L^2_{0(i-1)}$ and that $L_1^2$ is an $s_0^2$-crossing predecessor of $L^2_0$.

 We will associate to $p^1$ the unique point of intersection, say $p^2\in \clos{\mathcal{P}_2}$, of $s^2_0$ with the leftmost unstable boundary component of $L_1^2$. Notice that, since $L_1^1$ is not  the leftmost $s^1_0$-crossing predecessor of $L_0^1$ ($L^1_3$ is an $s^1_0$-crossing predecessor of $L_0^1$ at the left of $L_1^1$, see Figure \ref{f.cycle}), neither is $L_1^2$. Hence, $p^2$ is not a corner point of $L_0^2$ and $L_0^2$ contains a germ of the $(-,-)$ and $(+,-)$ quadrants of $p^2$.

By a similar argument, one can prove that $L_2^2$ is a crossing successor of $L_1^2$ in $\clos{\mathcal{P}_2}$, such that $p^2\in \partial L_1^2 \cap \partial L_2^2$ and  $\inte{L_0^2}\cap \inte{L_2^2}=\emptyset$ and by a finite repetition of this argument that $(L^2_0,L^2_1..., L^2_k)$ is part of the positive cycle of $p^2$ starting from $L^2_0$. Similarly, we also get that $(L^2_0,L^2_3..., L^2_k)$ is part of the negative cycle of $p^2$ starting from $L^2_0$.

\textit{\underline{End of the proof of the theorem}}

By our previous arguments, the rectangle paths $C_1$ and $C_2$ are respectively associated two rectangle paths in $\clos{\mathcal{P}_2}$ of the form: 

$$r^2_0,...,r^2_k,L^2_0,L^2_{01},...,L^2_{0s(0)},L^2_1, L^2_{11},...,L^2_{1s(1)},..., L^2_k,r^2_l,r^2_{l+1},...,r^2_n$$  
 $$r^{2}_0,...,r^{2}_k,L^{2}_0,L'^{2}_{01},...,L'^{2}_{0s'_0},L^{2}_3, L'^2_{31},...,L'^{2}_{3s'_3},...,L^{2}_k,r^{2}_l,r^2_{l+1},...,r^{2}_n$$
where $(L^{2}_0,L^{2}_1,L^{2}_2,L^{2}_3,L^{2}_4)$ and $(L^{2}_0,L^{2}_3,L^{2}_2,L^{2}_1,L^{2}_4)$ are the two cycles starting from $L^{2}_0$ around $p^2\in L^2_0$, the path $L^2_0,L^2_{01},...,L^2_{0s_0},L^2_1, L^2_{11},...,L^2_{1s_1},..., L^2_k$ is the rectangle path associated to the generalized rectangle path $(L^2_0,L^2_1,...,L^2_k)$ and $L^2_0,L'^2_{01},...,L'^2_{0s'(0)},L^2_3, L'^2_{31},$ $...,L'^2_{3s'_3},...,L^2_k$ is the rectangle path associated to the generalized rectangle path $(L^2_0,L^2_3,...,L^2_k)$. By Definition \ref{d.homotopyC}, the previous rectangle paths are homotopic by a homotopy of type $\mathcal{C}$, which gives us the desired result. 
\end{proof}

\subsection{Proof of Theorem \ref{t.homotopictotrivialpath}} 
\label{s.proofoftheoremhomtrivialpath}

Fix $\rho: G\rightarrow \text{Homeo}(\mathcal{P})$ an orientation preserving strong Markovian action, preserving the pair of singular foliations $\mathcal{F}^s$ and $\mathcal{F}^u$ and leaving invariant a strong Markovian family $\mathcal{R}$. We denote by $\Gamma$ the boundary periodic point of $\mathcal{R}$, by $\clos{\mathcal{P}}$ the bifoliated plane of $\rho$ up to surgeries on $\Gamma$, by $\clos{\mathcal{R}}$ the lift of $\mathcal{R}$ on $\clos{\mathcal{P}}$, by $\clos{\Gamma}$ the lift of $\Gamma$ on $\clos{\mathcal{P}}$ and by $\clos{\rho}: \clos{G}\rightarrow \text{Homeo}(\clos{\mathcal{P}})$ the lift of the action $\rho$ on $\clos{\mathcal{P}}$. Consider $\clos{\mathcal{F}^s}$ and $\clos{\mathcal{F}^u}$ the stable and unstable foliations of $\clos{\rho}$ endowed with an orientation.

Our main goal in this section consists in showing Theorem \ref{t.homotopictotrivialpath}, the proof of which will be split in two parts: 
\begin{enumerate}
    \item Proving Theorem \ref{t.homotopictotrivialpath} in the case of rectangle paths that can be associated to  simple closed polygonal curves: 
    \begin{prop}\label{p.homotopictotrivialsimplecase}
    Let $\gamma$ be a simple closed good polygonal curve in $\clos{\mathcal{P}}$ and $R_0,..,R_n$ one of its associated rectangle paths. 
    \begin{itemize}
        \item If $R_0=R_n$, then $R_0,...,R_n$ is homotopic to the trivial rectangle path.
        \item If $R_n$ is a successor or  predecessor of generation $k$ of $R_0$ and $R_0=r_0,r_1,..r_k=R_n$ is the unique monotonous rectangle path from $R_0$ to $R_n$ (see Lemma \ref{l.npredecessor}), then $R_0,...,R_n$ is homotopic to $r_0,...,r_k$. 
    \end{itemize}
    \end{prop}
    \item Proposition \ref{p.homotopictotrivialsimplecase} implies Theorem \ref{t.homotopictotrivialpath}
\end{enumerate}
\addtocontents{toc}{\protect\setcounter{tocdepth}{2}}
\subsubsection*{Proof of Proposition \ref{p.homotopictotrivialsimplecase}}
We will prove  Proposition \ref{p.homotopictotrivialsimplecase} by induction. 
Let $\gamma: [0,1] \rightarrow \clos{\mathcal{P}}$  be a simple good closed polygonal curve, $R_0\in \clos{\mathcal{R}}$ such that $\gamma(0)\in \inte{R_0}$, $R_0,...,R_n$ the rectangle path starting from $R_0$ associated to $\gamma$ and $M$ be the number of boundary arc points of $\clos{\mathcal{R}}$ in the interior of $\gamma$ (see Definition \ref{d.simple}). Notice that since $\gamma$ is good, by Remark \ref{r.caracterisationarcpoints} applied for $\clos{\mathcal{R}}$, it can not contain any boundary arc or boundary periodic points; hence, the number of boundary arc points in the closure of the interior of $\gamma$ is also equal to $M$. 

Let us show that $M<\infty$. Indeed, by compactness, we can cover the closure of the interior of $\gamma$ by a finite number of rectangles of our Markovian family $\clos{\mathcal{R}}$. By Lemmas  \ref{l.crossingrectanglesnoperiodicpoints},  \ref{l.crossingrectangleswithperiodicpoints} and Remark \ref{r.caracterisationarcpoints} applied for $\clos{\mathcal{R}}$, the number of boundary arc points in the closure of the interior of $\gamma$ is finite except if a point of $\clos{\Gamma}$ lies in the interior of $\gamma$. This is impossible, since $\gamma$ is good and no loop in $\clos{\mathcal{P}}-\clos{\Gamma}$ can go around a point in $\clos{\Gamma}$ (see Remark \ref{r.gammabarneverinterior}). 

Let $l(\gamma)$ be the length of the good polygonal curve $\gamma$. We will apply an induction on the well ordered set of couples $(M,l)\in \mathbb{N}\times \mathbb{N}$ endowed with the lexicographic order: $(M_1,l_1)<(M_2,l_2)$ if and only if $M_1<M_2$ or $M_1=M_2$ and $l_1<l_2$.  At every step of our induction, we  will perform a finite sequence of homotopies of type $\mathcal{A}, \mathcal{B}$ or $\mathcal{C}$ on $\gamma$ -therefore also on its associated rectangle path-, that will produce a new simple, good and closed polygonal curve $\gamma'$ such that $(M(\gamma'),l(\gamma'))<(M(\gamma),l(\gamma))$.
 
Notice that any simple and good polygonal curve has length at least equal to 4. Hence, the shortest such curves bound a rectangle in $\clos{\mathcal{P}}$. We will therefore initialize our induction by considering the case where $\gamma$ bounds a rectangle containing no boundary arc points in its interior. 

Let us first fix some notations. Recall that $R_0,...,R_n$ is a rectangle path associated to $\gamma$ and that by Remark \ref{r.polygonalcurverectassociation}, there exists $0=c_0<c_1<...c_{n+1}=1$ a $(n+2)$-uple such that $\gamma([c_i,c_{i+1}])\subset R_i$ and a function $\textit{Rect}_{\gamma,R_0}$ associated to $c_0<c_1<...c_{n+1}$ sending points of $[0,1]$ to rectangles in $\lbrace R_0,...,R_n\rbrace$. Fix such a collection of $c_i$ satisfying the results of Lemma \ref{l.choiceci}. Next, for the sake of simplicity,  for any interval $I\subset [0,1]$ if $A=\gamma(I)$, we will denote $\textit{Rect}_{\gamma,R_0}(I)$ also by  $\textit{Rect}_{\gamma,R_0}(A)$. By following the rectangles associated to the different points in $I$ (see Convention \ref{conv.rectanglepathimage}), we will think of $\textit{Rect}_{\gamma,R_0}(A)=\textit{Rect}_{\gamma,R_0}(I)$ as a rectangle path, instead of just a set of rectangles.

\subsubsection*{Initializing the induction} 
Assume that $(M(\gamma),l(\gamma))=(0,4)$. Otherwise said, assume that the simple and good polygonal curve $\gamma$ bounds a rectangle in $\clos{\mathcal{P}}$ containing no boundary arc points in its interior. We will show that under this hypothesis, the rectangle path $R_0,...,R_n$ is either trivial or monotonous, which gives us the desired result.

 By Remark \ref{r.gamma0},  $\gamma(0)=\gamma(1)$ is a corner of the rectangle $\gamma$. Let us denote by $u',s,u,s'$ the four sides of $\gamma$, following the order in which $\gamma$ visits them.

We will associate to one of the previous  segments of $\gamma$, say $S$, the number $0$ if the rectangle path $\textit{Rect}_{\gamma,R_0} (S)$ is trivial and the number $1$ if not. We write $S\rightarrow 0$ in the first case and $S\rightarrow 1$ in the latter. In this way, we can associate to $\gamma$ a unique element of $\{0,1\}^4$. Our choice of $c_i$ (see Lemma \ref{l.choiceci}) and our method for associating rectangle paths to curves (see Definition \ref{d.curvetorectanglepath}) restrict the elements of $\{0,1\}^4$ that can be associated to $\gamma$. 

\begin{lemm}\label{l.possiblecases}
If $\gamma$ bounds a rectangle, then the only elements of $\{0,1\}^4$ that can be associated to $\gamma$ are $(0,0,0,0),(0,1,1,1)$, $(0,1,0,0), (1,0, 0,0)$ and $(1,1,1,1)$. Furthermore, $\gamma$ is associated to  $(0,0,0,0),(0,1,1,1)$ or $(0,1,0,0)$ if and only if it contains no boundary arc points in its interior. 
\end{lemm}
\begin{proof}
Without any loss of generality, let us assume  that $u,u'$ are unstable segments and $s,s'$ are stable segments. Recall that, $R_0,R_1,...,R_n$, the rectangle path starting from $R_0$ associated to $\gamma$, can be obtained by the juxtaposition of $\textit{Rect}_{\gamma,R_0} (u')$, followed by $\textit{Rect}_{\gamma,R_{0}} (s)$, followed by $\textit{Rect}_{\gamma,R_{0}} (u)$ and finally followed by $\textit{Rect}_{\gamma,R_{0}} (s')$.

If $(u',s) \rightarrow (0,0)$, then by our choice of $c_i$ (see Lemma \ref{l.choiceci}), the rectangle path starting from $R_0$ associated to the juxtaposition of $u'$ followed by $s$ is trivial and therefore, the segments $u'$ and $s$ do not exit $R_0$. Since $R_0$ is trivially bifoliated and $\gamma$ bounds a rectangle, we also have that $s',u\subset R_0$. Therefore, in this case,  $\gamma$ does not exit $R_0$ at all and its associated element in $\{0,1\}^4$ is $(0,0,0,0)$ (recall that 
the rectangle path associated to $\gamma$ starting from $R_0$ is not trivial if and only if $\gamma$ exits $R_0$). 

Suppose now that $(u',s)\rightarrow(1,1)$. Endow $[0,1]$ with the left to right orientation. Endow also $u'$ and $s$ with the orientation given by $\gamma: [0,1]\rightarrow \clos{\mathcal{P}}$ and assume without any loss of generality that the segments $u'$ and $s$ are positively oriented for our choice of orientations of $\clos{\mathcal{F}^{s,u}}$ given in the beginning of this section. In the following lines, for any rectangle $R\in \clos{\mathcal{R}}$, using the orientations of $\clos{\mathcal{F}^{s,u}}$, we will denote by $\partial^s_+R$ (resp. $\partial^s_-R$) its upper (resp. lower) stable boundary. We similarly define $\partial^u_+R$ and $\partial^u_-R$. 

Thanks to the fact that $u'\rightarrow 1$ and $s\rightarrow 1$, to our choice of $c_i$ and to our definition of $\textit{Rect}_{\gamma,R_0}$, we have that 
\begin{itemize}
\item $u'$ exits $R_0$. Since $u'$ is positively oriented, we can suppose that $u'$ will exit $R_0$ in order to visit $R_{u_1}$, a $\partial^s_+R_0$-crossing predecessor of $R_0$, then it will exit $R_{u_1}$ in order to visit $R_{u_2}$, a $\partial^s_+R_{u_1}$-crossing predecessor of $R_{u_1}$, and so on, until it exits $R_{u_{N-1}}$ and visits $R_{u_N}:=\textit{Rect}_{\gamma,R_0}(u'\cap s)$, a predecessor of some generation of $R_0$ that contains $u'$
\item $\textit{Rect}_{\gamma,R_0} (u')= R_0,R_1,...,R_{u_1},...,R_{u_2},... R_{u_N}$ 
\item Similarly, $s$ exits $R_{u_N}$. Since $s$ is positively oriented, we can suppose that $s$ will exit $R_{u_N}$ in order to visit $R_{u_{N+1}}$, a $\partial^u_+R_{u_N}$-crossing successor of $R_{u_N}$, next $s$ will exit $R_{u_{N+1}}$ in order to visit a $\partial^u_+R_{u_{N+1}}$-crossing successor of $R_{u_{N+1}}$, and so on until it exits $R_{u_{m-1}}$ and visits  $R_{u_m}:=\textit{Rect}_{\gamma,R_0}(s\cap u)$ (see Figure \ref{f.impossiblecasesc}), a successor of some generation of $R_{u_N}$ containing $s$. 
\item $\textit{Rect}_{\gamma,R_{0}} (s)= R_{u_N},...,R_{u_N+1},..., R_{u_m}$

\end{itemize}

Notice that since $u'$ exits $R_{u_{N-1}}$, thanks to Lemma \ref{l.npredecessor}, we have that $R_{u_{N+1}}$ and $R_{u_{N-1}}$ correspond two  successors (of some generation) of $R_{u_N}$ with disjoint interiors. Also, $R_{u_m}$ (resp. $R_0$) is a successor of some generation of $R_{u_{N+1}}$ (resp. $R_{u_{N-1}}$).  By the previous facts and Lemma \ref{l.intersectionsuccessorspredecessors} applied for $\clos{\mathcal{R}}$ (see also Remark \ref{r.equivdefi}), $\inte{R_{u_m}}\cap \inte{R_0}=\emptyset$. Consequently, if $u$ does not exit $R_{u_m}$, then the curve $\gamma$ cannot bound a rectangle. Indeed, this would imply that the stable leaf containing $s'$ intersects $R_{u_N}$ along more than one connected components, which would contradict Item (5) of Proposition \ref{p.propertiesoffoliformarkovianactions}. Hence, as before, $u$ exits $R_{u_m}$ in order to visit $R_{u_{m+1}}$, a $\partial^s_-R_{u_{m}}$-crossing predecessor of $R_{u_m}$, next $u$ exits $R_{u_{m+1}}$ in order to visit a  crossing predecessor of $R_{u_{m+1}}$ and so on until it exits $R_{u_{k-1}}$ and visits $R_{u_k}:=\textit{Rect}_{\gamma,R_0}(u\cap s')$, a predecessor of some generation of $R_{u_m}$ that contains $u$. By the same arguments as before, $\inte{R_{u_k}}\cap \inte{ R_{u_N}}=\emptyset$; hence $s'$ has to exit $R_{u_k}$ in order to reach $\gamma(0)\in \inte{R_{u_N}}$. We conclude that if $(u',s) \rightarrow (1,1)$, then the element of $\{0,1\}^4$ associated to $(u',s,u,s')$ is $(1,1,1,1)$. 
\begin{figure}[h]
\centering 
\includegraphics[scale=0.7]{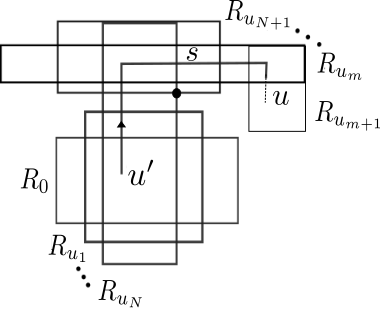}
\caption{}
\label{f.impossiblecasesc}
\end{figure}

If $(u',s) \rightarrow (1,0)$, then as before $u'$ exits $R_0$ in order to visit $R_{u_1}$, a crossing predecessor of $R_0$, next $u'$ exits $R_{u_1}$ in order to visit a crossing predecessor of $R_{u_1}$ and so on until it exits $R_{u_{N-1}}$ and visits $R_{u_N}:=\textit{Rect}_{\gamma,R_0}(u'\cap s)$ a predecessor of some generation of $R_0$ that contains the entire segment $u'$. Since, by hypothesis, $s$ does not exit $R_{u_N}$ and $\gamma$ bounds a rectangle, $u,s'$ do not exit $R_{u_N}$ either (once again we use here Item (5) of Proposition \ref{p.propertiesoffoliformarkovianactions}). Therefore, if $(u',s)\rightarrow (1,0)$, then the element of $\{0,1\}^4$ associated to $(u',s,u,s')$ is $(1,0,0,0)$.

We show in a similar way that if $(u',s)\rightarrow (0,1)$, then the element of $\{0,1\}^4$ associated to $(u',s,u,s')$ is either  $(0,1,0,0)$ or $(0,1,1,1)$.

Finally, notice that when $\gamma$ is associated to $(0,0,0,0), (1,0,0,0)$ or $(0,1,0,0)$, then by our previous arguments (see cases $(u',s)\rightarrow (0,0)$ and $(u',s)\rightarrow (1,0)$), $\gamma$ is contained in the interior of a rectangle of $\clos{\mathcal{R}}$. Hence, its interior cannot contain boundary arc points. Let us now show that if $\gamma$ is associated to either $(1,1,1,1)$ or $(0,1,1,1)$, then the interior of $\gamma$ necessarily contains boundary arc points. Assume that $(u',s,u,s') \rightarrow (1,1,1,1)$ (the case $(0,1,1,1)$ can be treated in a similar way). Following our previous notations (see the case $(u',s)\rightarrow (1,1)$), $u'\subset  \inte{R_{u_{N}}}$, $s$ exits $R_{u_N}$ in order to visit $R_{u_{N+1}}$, a $\partial^u_+R_{u_N}$-crossing successor of $R_{u_N}$ and $\inte{R_{u_{N-1}}}\cap \inte{R_{u_{N+1}}}=\emptyset $. Since $\gamma$ bounds a rectangle, its interior contains the point  $\partial^u_+R_{u_N}\cap \partial^s_-R_{u_{N+1}}$ (see Figure \ref{f.impossiblecasesc}), which is a boundary arc point according to Remark \ref{r.caracterisationarcpoints}. 
\end{proof}
We are now ready to describe the initialization of our induction, in which we suppose $\gamma$ to be the boundary of a rectangle, whose interior contains no boundary arc points of $\clos{\mathcal{R}}$. By the previous lemma, we will consider the following 3 cases:  

\begin{enumerate}
    \item $\gamma$ is associated to $(0,0,0,0)$: in this case, by our proof of Lemma \ref{l.possiblecases}, we have that $\gamma$ never exits $R_0$; its associated rectangle path is therefore  trivial. 
    \item $\gamma$ is associated to $(1,0,0,0)$: in this case, following the notations of Lemma \ref{l.possiblecases} (see the case where $(u',s)\rightarrow (1,0)$), $u'$ exits $R_0$ in order to visit $R_{u_1}$, a crossing predecessor of $R_0$, next $u'$ exits $R_{u_1}$ in order to visit a crossing predecessor of $R_{u_1}$ and so on until it exits $R_{u_{N-1}}$ and visits $R_{u_N}:=\textit{Rect}_{\gamma,R_0}(u'\cap s)$ a predecessor of some generation of $R_0$ that, as we have previously shown, contains completely $u',s,u,s'$. It follows that the rectangle path starting from $R_0$ associated to $\gamma$ coincides with the unique  decreasing rectangle path from $R_0$ to $R_{u_N}$.  
    \item $\gamma$ is associated to $(0,1,0,0)$: in this case, similarly to the previous case, one can show that the rectangle path starting from $R_0$ associated to $\gamma$ is increasing. 
\end{enumerate}
We deduce that if  $(M(\gamma),l(\gamma))=(0,4)$, then Proposition \ref{p.homotopictotrivialsimplecase} is true. Let us now describe the step of our induction. 

\subsubsection*{Induction step: the case where $l(\gamma)=4$ and $M(\gamma)>0$}

Following our previous notations, thanks to Lemma \ref{l.possiblecases}, the element of $\{0,1\}^4$ associated to $\gamma$ is either  $(1,1,1,1)$ or $(0,1,1,1)$. Suppose first that $\gamma$ corresponds to $(1,1,1,1)$. Assume once again, without any loss of generality, that the (oriented by $\gamma$) segments $u'$ and $s$ are positively oriented for our choice of orientation of $\clos{\mathcal{F}^{s,u}}$. 
\begin{figure}
    \centering
    \includegraphics[scale=0.7]{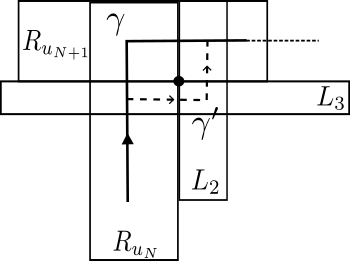}
    \caption{}
    \label{f.caserectanglehomtypec}
\end{figure}

During our proof of Lemma \ref{l.possiblecases}, we established that $u'$ will exit $R_0$ in order to visit $R_{u_1}$, a $\partial^s_+R_0$-crossing predecessor of $R_0$, then it will exit $R_{u_1}$ in order to visit $R_{u_2}$, a $\partial^s_+R_{u_1}$-crossing predecessor of $R_{u_1}$, and so on, until it exits $R_{u_{N-1}}$ and visits $R_{u_N}$, a predecessor of some generation of $R_0$ containing $u'$. We also showed that $s$ exits $R_{u_N}$ in order to visit $R_{u_{N+1}}$ and that $\partial^u_+R_{u_N}\cap \partial^s_-R_{u_{N+1}}$ is a boundary arc point contained in the interior of $\gamma$ (see Figure \ref{f.caserectanglehomtypec}). Consider the negative cycle $(R_{u_N},R_{u_{N+1}},L_2,L_3,L_4)$ around the previous boundary arc point starting from  $R_{u_N}$. Using the previous cycle, we can  deform $\gamma$ by a homotopy of type $\mathcal{C}$ to a good and simple polygonal curve $\gamma'$ containing strictly less boundary arc points in its interior (see Figure \ref{f.caserectanglehomtypec}).

We treat the case where $\gamma$ corresponds to  $(0,1,1,1)$ in a similar way.

\subsubsection*{Induction step in the general case} We will assume here that $l(\gamma)>4$. In this case, we begin by choosing one stable or unstable tangency $s$ such that :
\begin{enumerate}
    \item $s$ is of type (a) or (b) and has a complete domain $Dom$ (see Lemma \ref{l.canonicalneighbourhoodtangency})
    \item $\gamma(0)\notin Dom$ or $\gamma(0)$ is a corner point of the rectangle $Dom$ that is not in $s$
    \item the interior of one side of $\partial Dom$ is disjoint from $\gamma$
\end{enumerate}
The existence of such a tangency is proven in the following lemma. We will call a tangency with the above properties, a tangency with \emph{property $(\star)$}. If furthermore $\gamma(0)\notin Dom$, we will say that the tangency satisfies the \emph{strong  property $(\star)$}. 

\begin{lemm}\label{l.existencetangenciespropstar}
Let $\gamma$ be a good, simple and closed polygonal curve of length strictly greater than 4. There exists $s$ a stable or unstable tangency of $\gamma$ with property $(\star)$. 
\end{lemm}
\begin{proof}Indeed, take $s$ to be a stable tangency of $\gamma$ of type (a) or (b) (the existence of at least 2 of such tangencies is assured by Lemma \ref{l.numbertangencies}). Suppose that the domain of $s$, that we will denote by $Dom$, is not complete. In that case, by Remark \ref{r.incompletedomain}, there exists a stable segment $s'$ of $\gamma$ intersecting the interior of the stable boundary component of $Dom$ that is not $s$. Since $\inte{Dom}\cap \gamma=\emptyset$ and $Dom$ is incomplete, $s'$ is a stable tangency of $\gamma$ of type (c) or (d). Therefore, to every tangency of $\gamma$ of type (a) or (b) with an incomplete domain we can injectively associate a tangency of type (c) or (d). We deduce by Lemma \ref{l.numbertangencies} that there exist at least two stable tangencies of type (a) or (b) for $\gamma$, whose domains are complete. 


Consider now $s$ a stable tangency of $\gamma$ of type (a) or (b) with a complete domain $Dom$. The function $\gamma:\mathbb{S}^1\rightarrow \clos{\mathcal{P}}$ induces a cyclic order on the stable/unstable segments forming $\gamma$. Denote by $u$ (resp. $u'$) the unstable segment of $\gamma$ after (resp. before) $s$  and $s'$ the stable boundary component of $Dom$ that is not $s$ (see Figure \ref{f.reductioncase1bis}). By  definition, the domain of $s$ satisfies,  $\inte{Dom}\cap \gamma=\emptyset$, hence if  $\gamma(0)\in Dom$ then $\gamma(0)\in \partial Dom$. 

If $\gamma(0)\in Dom$, then by Remark \ref{r.gamma0}, $\gamma(0)$ cannot belong to the interior of $s$, $u$ or $u'$. Therefore, $\gamma(0)\in s'$ or $\gamma(0)\in \partial s$. Since there are at least 2 stable tangencies of $\gamma$ of type (a) or (b) with a complete domain, we can  assume without any loss of generality that $\gamma(0) \notin s$. 

\begin{figure}[h]
\centering 
\includegraphics[scale=0.45]{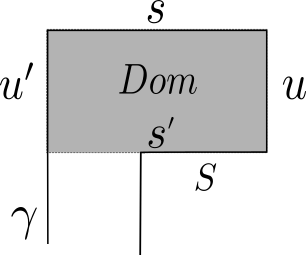}
\caption{}
\label{f.reductioncase1bis}
\end{figure}

Assume now that the interior of $s'$ intersects $\gamma$. Assume without any loss of generality that $u\subset Dom$ (since $Dom$ is complete). Since $\gamma$ is good, any two stable sides of $\gamma$ do not belong to the same stable leaf of $\clos{\mathcal{F}^s}$. Therefore, we have necessarily that the stable segment of $\gamma$ after $u$, say $S$, is contained in $s'$. Since $\gamma$ does not bound a rectangle $S \subsetneq s'$; we are therefore in the situation of Figure  \ref{f.reductioncase1bis}. Hence, $u$ is an unstable  tangency of type (a) or (b) with a complete domain.  If $\gamma(0)\notin u$, then $u$ has property $(\star)$. If $\gamma(0)\in u$, then take $s_{fin}\neq s$ to be another stable tangency of $\gamma$ of type (a) or (b) with a complete domain (we showed previously that there are at least 2 such tangencies) and  notice that $\gamma(0)\notin s_{fin}$. If  $s_{fin}$ does not satisfy $(\star)$, then  by repeating the arguments of this paragraph, either the unstable segment of $\gamma$ before or after $s_{fin}$ is an unstable tangency of $\gamma$ with property $(\star)$. 
\end{proof}
\begin{figure}[h]
    \centering
    \includegraphics[scale=0.8]{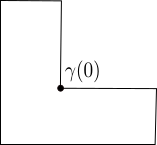}
    \caption{This simple, good, closed polygonal curve does not admit stable or unstable tangencies with strong property $(\star)$}
    \label{f.Lshape}
\end{figure}
Even when $\gamma$ does not bound a rectangle, it is not always possible to find a tangency with strong property $(\star)$ (see for instance Figure \ref{f.Lshape}). However, there exist numerous cases where this is possible: 

\begin{lemm}\label{l.existencetangenciesgood}
Let $\gamma$ be a simple, closed and good polygonal curve that does not bound a rectangle. Assume that $\gamma(0)$ belongs to a stable or unstable tangency of $\gamma$ of type (a) or (b). Then there exists a stable or unstable tangency of $\gamma$ with strong property $(\star)$. 
\end{lemm}
\begin{proof}
Indeed, assume without any loss of generality that $\gamma(0)\in S$, where $S$ is a stable tangency of $\gamma$ of type (a) or (b). By our proof of Lemma \ref{l.existencetangenciespropstar}, there exists $s$ a stable tangency of $\gamma$ of type (a) or (b) with a complete domain $Dom$ such that $\gamma(0)\notin s$. Notice that $Dom$ cannot contain $\gamma(0)$, since the domain of a stable tangency of type (a) or (b) cannot intersect another  stable tangency of type (a) or (b), except when $\gamma$ bounds a rectangle. By our proof Lemma \ref{l.existencetangenciespropstar}, either $s$ has strong property $(\star)$ or we are in the case of Figure \ref{f.reductioncase1bis}. In the latter case, using the notations of Figure \ref{f.reductioncase1bis}, the segment $u$ has property $(\star)$ and its domain $D'\subset Dom$ does not contain $\gamma(0)$; hence $u$ has strong property $(\star)$. 
\end{proof}

We will now describe our induction step in the following two subcases:
\begin{enumerate}
    \item there exists $s$ a tangency in $\gamma$ with strong property $(\star)$ 
    \item no such tangency exists
\end{enumerate}
\subsubsection*{Induction step in the case where a tangency with strong $(\star)$ property exists}

\begin{figure}[h]
\centering 
\includegraphics[scale=0.5]{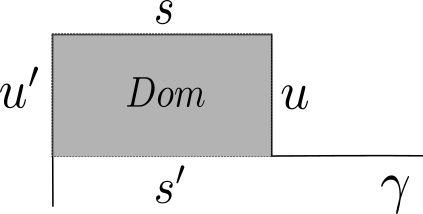}
\caption{}
\label{f.reduction1case}
\end{figure}

Assume without any loss of generality that the stable tangency $s$ of $\gamma$ has strong property $(\star)$. The function $\gamma:[0,1]\rightarrow \clos{\mathcal{P}}$ endows the segments forming $\gamma$ with a total order, for which $s$ is neither the first nor the last segment (since $\gamma(0)\notin s$). Let us denote by $U'$ (resp. $U$) the unstable segment of $\gamma$ before (resp. after) $s$ and by $s'$ the stable side of the domain $Dom$ of $s$ that is not $s$. We define $u':=U'\cap Dom$ and $u:=U\cap Dom$. We are therefore -up to a change of orientation and up to interchanging $u'$ and $u$- in the case of Figure \ref{f.reduction1case}. 


As before, we will associate to a stable or unstable segment of $\gamma$, say $S$, the number $0$ if $\textit{Rect}_{\gamma,R_0} (S)$ consists of one rectangle and the number $1$ if not. We will write $S\rightarrow 0$ in the first case and $S\rightarrow 1$ in the latter. We can therefore associate to $U',s,U$ a unique element of $\{0,1\}^3$. 

We are now ready to describe the induction step in this case: 
\begin{enumerate}[label=(\arabic*{$'$}), leftmargin=0.01cm, itemsep=0.1cm]
    \item If $(U',s,U)\rightarrow(0,0,0)$, then  $\textit{Rect}_{\gamma,R_0}(U')=\textit{Rect}_{\gamma,R_0}(U)=\textit{Rect}_{\gamma,R_0}(s)$. In other words, the segments $U',s,U$ of $\gamma$ are all contained in a rectangle $R$ that is part of the rectangle path starting from $R_0$ associated to $\gamma$. We deduce that $Dom\subset R$ and by pushing $s$ along $Dom\subset R$ in order that $s$ be identified with $s'$ (see Figure \ref{f.reduction1case}), we obtain a new simple, closed, good polygonal curve $\gamma'$ of strictly smaller length and such that $M(\gamma)=M(\gamma')$. It is easy to check that $\gamma$ and $\gamma'$ are homotopic by a homotopy of type $\mathcal{A}$ (relatively to $R_0,...,R_n$). 
    
    \begin{figure}[h]
    
    \centering
  \begin{minipage}[c]{0.4\textwidth}
  \hspace{0.72cm}
    \includegraphics[width=0.6\textwidth]{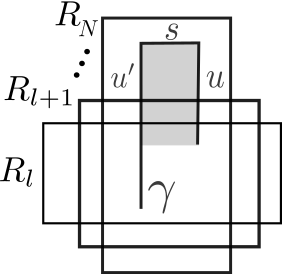}
    \vspace{0.3cm}
  \caption*{(a)}
    \end{minipage}
 \begin{minipage}[c]{0.4\textwidth}
    \hspace{0.8cm}
     \includegraphics[width=0.6\textwidth]{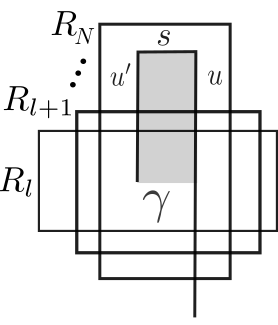}
    \caption*{(b)}
    
  \end{minipage}
\caption{}
\label{f.impossiblecasesab}
\end{figure}
    \item If $(U',s,U)\rightarrow(1,0,1)$, then by Lemma \ref{l.choiceci}, $\textit{Rect}_{\gamma,R_0}(U')$ corresponds to a non-trivial decreasing rectangle path, say $R_l,...,R_N$, that is part of the rectangle path $R_0,...,R_n$. Moreover,  $U'\subset R_N$, $s$ does not exit $R_N$ and $U$ exits $R_N$ in order to visit a crossing predecessor of $R_N$ (see Figure \ref{f.impossiblecasesab}b). Thanks to the above facts, $u'=U'$. By pushing $\gamma$ along $Dom$ so that $s$ be identified with $s'$, we can construct a simple, good and closed polygonal curve $\gamma'$ of strictly smaller length and such that $M(\gamma)=M(\gamma')$.

    Let us show that the previous movement corresponds to a homotopy of type $\mathcal{A}$. Indeed, our  movement pushes the  segments $u',s,u$ of $\gamma$ to the segment $s'$ of $\gamma'$, while keeping $
\gamma-(u'\cup s \cup u)$ fixed. Let $a,b\in [0,1]$ such that $\gamma([a,b])=u'\cup s \cup u$. By eventually reparametrizing $\gamma'$ (the rectangle path associated to $\gamma':[0,1]\rightarrow \clos{\mathcal{P}}$ does not depend on the choice of parametrization), we can assume that $\gamma'([a,b])=s'$. Similarly to $\gamma$, using Lemma \ref{l.choiceci}, we can define the function $\textit{Rect}_{\gamma',R_0}$. Notice that 
\begin{itemize}
    \item  $\textit{Rect}_{\gamma',R_0}([0,a])$ and $\textit{Rect}_{\gamma,R_0}([0,a])$ correspond to the same rectangle paths. In other words, the rectangle paths starting from $R_0$ associated to the good polygonal curves $\gamma_{|[0,a]}$ and $ \gamma'_{|[0,a]} \equiv \gamma_{|[0,a]}$ are both equal to $R_0,...,R_l$.  
    \item since $s'$ does not exit $R_l$, $ \textit{Rect}_{\gamma',R_0}(s')=\textit{Rect}_{\gamma,R_0}(a)=\{R_l\}$
    \item $u_{\text{after}}:=U-u\subset \gamma'$ will exit $R_l$ in order to visit a crossing predecessor of $R_l$ and will eventually exit $R_N$ in order to visit a crossing  predecessor of $R_N$. Recall that by our initial hypothesis and by Lemma \ref{l.npredecessor}, $R_l,R_{l+1},...,R_N$ is the unique decreasing rectangle path from $R_l$ to $R_N$. It follows that the rectangle path $\textit{Rect}_{\gamma',R_0}(u_{\text{after}})$ is of the form $R_l,R_{l+1},...,R_N,...$. Consequently, since $u_{\text{after}}\subset U$, the rectangle paths $\textit{Rect}_{\gamma,R_0}(U'\cup s\cup U)$ and $  \textit{Rect}_{\gamma',R_0}( s'\cup u_{\text{after}})$ are the same. 
\end{itemize} Finally, since $\gamma$ and $\gamma'$ coincide on $[b,1]$, we deduce that the rectangle paths starting from $R_0$ associated to $\gamma'$ and $\gamma$ are the same and thus our movement corresponds indeed to a homotopy of type $\mathcal{A}$.  
    \item If $(U',s,U)\rightarrow(0,0,1)$, $\textit{Rect}_{\gamma,R_0}(U')=\textit{Rect}_{\gamma,R_0}(s)=\{R_l\}$ and $U$ exits $R_l$ in order to visit a crossing predecessor of $R_l$. We are therefore in the case of Figure \ref{f.impossiblecasesab}b, where $R_l=R_{l+1}=R_{l+2}=...=R_N$ and $U'=u'$. As in Case $(1')$, by pushing $\gamma$ along $Dom\subset R_l$ so that we erase $U'$ and thus performing a homotopy of type $\mathcal{A}$, we get  that $\gamma$ is homotopic to a simple,  good and closed polygonal curve $\gamma'$ of strictly smaller length and such that $M(\gamma)=M(\gamma')$. 
    \item If $(U',s,U)\rightarrow(1,1,0)$, then $\textit{Rect}_{\gamma,R_0}(U')$ is a non-trivial decreasing rectangle path of the form $R_l,...,R_N$,  $\textit{Rect}_{\gamma,R_0}(s)$ is a non-trivial rectangle path of the form $R_N,R_{N+1}
,...,$ $R_m$ and $U$ does not exit $R_m$ (see Figure \ref{f.impossiblecasescsecond}). As before, $R_l,R_{l+1},...,R_m$ is part of the rectangle path $R_0,...R_n$. Moreover, thanks to the above facts, we get that $u=U$. By pushing $\gamma$ along $Dom\subset R_m$ so that we erase $U$ and thus performing  a homotopy of type $\mathcal{A}$, we get that $\gamma$ is homotopic to a simple, good and closed polygonal curve $\gamma'$ of strictly smaller length and such that $M(\gamma)=M(\gamma')$. 
\begin{figure}[h]
\centering 
\includegraphics[scale=0.5]{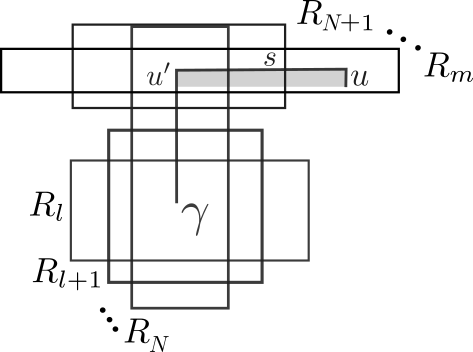}
\caption{}
\label{f.impossiblecasescsecond}
\end{figure}
\item If $(U',s,U)\rightarrow(0,1,0)$, then $\textit{Rect}_{\gamma,R_0}(U')=\{R_l\}$,  $\textit{Rect}_{\gamma,R_0}(s)$ is a non-trivial increasing rectangle path, say $R_l,R_{l+1}
,...,R_m$, that is part of $R_0,...,R_n$, and $U$ does not exit $R_m$. We are therefore in the case of Figure \ref{f.impossiblecasescsecond} with $R_l=R_{l+1}=...=R_N$.  As in Case $(4')$, by pushing $\gamma$ along $Dom\subset R_m$ so that we erase $U$ or $U'$ and thus performing a homotopy of type $\mathcal{A}$, we get that $\gamma$ is homotopic to a simple, good and closed polygonal curve $\gamma'$ of strictly smaller length and such that $M(\gamma)=M(\gamma')$. 
\item Assume now that $(U',s,U)\rightarrow(1,1,1)$. 
Denote by $R_{u'},...,R_s$ the non-trivial rectangle path defined by $\textit{Rect}_{\gamma,R_0}(s)$ (see Figure \ref{f.impossiblecase1}). Notice that 
\begin{itemize}
    \item since $\textit{Rect}_{\gamma,R_0}(s)$ is not trivial, $R_s\neq R_{u'}$ 
    \item since $\textit{Rect}_{\gamma,R_0}(U)$ is not trivial, $U$ will exit $R_s$ in order to enter a crossing predecessor of $R_s$, say $R_u$
    \item since $\textit{Rect}_{\gamma,R_0}(U')$ is not trivial, following $U'$ negatively starting from the unique point in $s\cap u'$, $U'$ must exit $R_s$. 
\end{itemize} We are therefore in the case of Figure \ref{f.impossiblecase1} (in which $U'$ can be eventually  contained in $R_{u'}$, but still exits $R_s$). Hence, there exists a point in $\partial^u{R_u}\cap \partial^s{R_s}$ inside $Dom$. By Remark \ref{r.caracterisationarcpoints}, this point is a boundary arc point in the interior of $\gamma$. Finally, as in the case described in \textit{Induction step: the case where $l(\gamma)=4$ and $M(\gamma)>0$}, by performing a homotopy of type $\mathcal{C}$ inside $Dom$, we get that $\gamma$ is homotopic to a simple, good and closed polygonal curve with strictly less boundary arc points in its interior. 
\begin{figure}[h]
    \centering
    \includegraphics[scale=0.5]{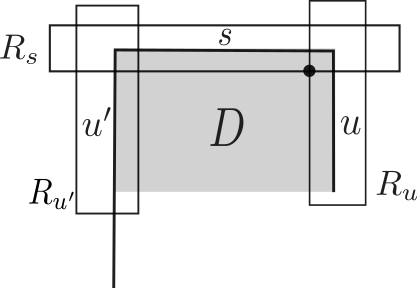}
    \caption{$Dom$ contains boundary arc points}
    \label{f.impossiblecase1}
\end{figure}

\item If $(U',s,U)\rightarrow(0,1,1)$, then denote once again $R_{u'},...,R_s$ the non-trivial rectangle path defined by $\textit{Rect}_{\gamma,R_0}(s)$. In this case, we have that $\textit{Rect}_{\gamma,R_0}(U')= \{R_{u'}\}$ which implies in particular that $U'\subset R_{u'}$. 
\begin{itemize}
    \item If $U'\not\subset R_s$, by the exact same argument as in Case $(6')$, both $U'$ and $U$ exit $R_s$. We are therefore in the case of Figure \ref{f.impossiblecase1} (with $U'\subset R_{u'}$). Once again, by performing a homotopy of type $\mathcal{C}$ inside $Dom$, we get that $\gamma$ is homotopic to a simple, good and closed polygonal curve with strictly less boundary arc points in its interior.  
    \item If $U'\subset R_s$, then we are in the case of Figure \ref{f.case011}, where $\textit{Rect}_{\gamma,R_0}(s)$ corresponds to the non-trivial increasing rectangle path $R_{u'},R_{u'+1},...,R_s$ and $\textit{Rect}_{\gamma,R_0}(U)$ to the non-trivial decreasing path $R_s,R_{s+1}
,...,R_m$. Notice that since $U$ exits $R_s$ and $U'$ does not, we have that $u'=U'$. As in Case $(4')$, by pushing $\gamma$ along $Dom\subset R_s$ so that we erase $U'$ and thus performing a homotopy of type $\mathcal{A}$, we get  that $\gamma$ is homotopic to a simple,  good and closed polygonal curve $\gamma'$ of strictly smaller length and such that $M(\gamma)=M(\gamma')$.
     
\end{itemize}

\begin{figure}[h]
\centering 
\includegraphics[scale=0.55]{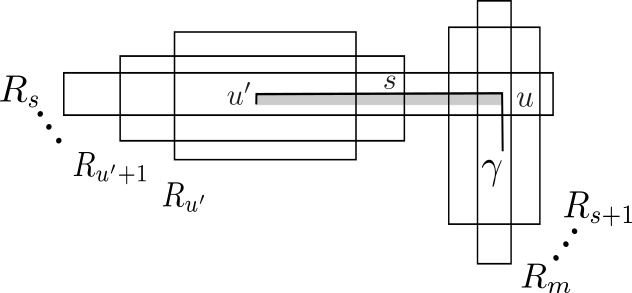}
\caption{}
\label{f.case011}
\end{figure}
    \item If $(U',s,U)\rightarrow (1,0,0)$, then we are in the case of case of Figure \ref{f.impossiblecasesab}a, where $\textit{Rect}_{\gamma,R_0}(U')$ corresponds to the decreasing rectangle path $R_l,...,R_N$ and $U'\cup s \cup U \subset R_N$. Let $a,b \in [0,1]$ such that $\gamma([a,b])=U'\cup s \cup U$.
    Before describing the induction step in this case, let us first prove a very useful lemma: 
    
    \begin{lemm}\label{l.increasingrectanglepaths} 
    Let $B:=L_1,...,L_k$ be a decreasing rectangle path. Consider any rectangle path of the form $L_p,w_B,L_{p'}$, where $p,p'\in \llbracket 1,k\rrbracket$ and $w_B$ is a rectangle path formed by rectangles in $B$. If $p'>p$ (resp. $p>p'$, $p=p'$), the previous rectangle path is homotopic by a sequence of homotopies of type  $\mathcal{B}$ to a decreasing (resp. increasing, trivial) rectangle path. 
    \end{lemm}
    \begin{proof}
        We will show the desired result in the case where $p'>p$. The other cases can be similarly treated. If $L_p,w_B,L_{p'}$ is not decreasing, we can find a non-monotonous rectangle path in $L_p,w_B,L_{p'}$ of length 3. The previous rectangle path is of the form $L_k,L_{k+1},L_k$ and can be ``erased" by a homotopy of type $\mathcal{B}$.  This reduces the length of $w_B$ and thus by a finite number of homotopies of type $\mathcal{B}$, we get the desired result.
    \end{proof}
    
    \begin{enumerate}[leftmargin=0.01cm, label=(8$'$\alph*)]
        
        \item Assume first that $\gamma_{|[a,1]}$ eventually exits $\underset{i=l}{\overset{N}{\cup}} R_i$ in order to visit a crossing predecessor (or successor) $L\notin \{R_l,...,R_N\}$  of $R_{k}\in \{R_l,...,R_N\}$. By eventually changing our choice of $R_k$, we will assume that if $L$ is a crossing predecessor (resp. successor) of $R_k$, then there is no predecessor (resp. successor) of $R_k$ in $\{R_l,...,R_N\}$, say $R_p$, such that $L$ is a crossing predecessor (resp. successor) of $R_p$. In this case, by pushing $\gamma$ along $Dom$ so that $s$ be identified with $s'$ (see Figure \ref{f.impossiblecasesab}a),  we obtain a simple, good and closed polygonal curve $\gamma'$ of strictly smaller length and such that $M(\gamma)=M(\gamma')$. Let us show that the previous movement corresponds to a sequence of homotopies of type  $\mathcal{B}$.

    \par{\quad{} We will assume that $U=u$ (the case $U'=u'$ follows from similar arguments). As in the Case $(2')$, our  movement pushes the  segments $u',s,u$ of $\gamma$ to the segment $s'$ of $\gamma'$, while keeping $
    \gamma-(u'\cup s \cup u)$ fixed. By eventually reparametrizing $\gamma'$, we can assume that $\gamma'([a,b])=u'_{\text{after}}\cup s'$, where  $u'_{\text{after}}:=U'-u'\subset \gamma'$. Therefore, $\gamma\equiv \gamma'$ on $[0,a]\cup [b,1]$ and the rectangle paths starting from $R_0$ associated to $\gamma_{|[0,a]}$ and $\gamma'_{|[0,a]}\equiv \gamma_{|[0,a]} $ are identical. In particular, both rectangle paths end by $R_l$. It therefore suffices to show that the rectangle paths starting from $R_l$ associated to $\gamma_{|[a,1]}$ and $\gamma'_{|[a,1]}$ are homotopic by a sequence of homotopies of type  $\mathcal{B}$.} 

    \par{\quad{} As $u'_{\text{after}}\cup s'\subset \underset{i=l}{\overset{N}{\cup}} R_i$ and $\gamma_{|[b,1]}\equiv \gamma'_{|[b,1]} $, both $\gamma([a,1])$ and $\gamma'([a,1])$ remain in 
    $\underset{i=l}{\overset{N}{\cup}} R_i$ until they exit this region in order to enter the crossing predecessor (or successor) $L$ of $R_k$. Therefore, there exist $R\notin \{R_l,...,R_N\}$ a predecessor (resp. successor) of $R_k$,  two rectangle paths, say $w_A,w'_A$, formed by rectangles in  $\{R_l,...,R_N\}$ and $W,W'$ two rectangle paths such that 
    \begin{itemize}
        \item the rectangle path starting from $R_l$ associated to   $\gamma_{|[a,1]} $ is $R_l,w_A,R_{k},R, W$
        \item the rectangle path starting from $R_l$ associated to $\gamma'_{|[a,1]} $ is $R_l,w'_A,R_{k},R, W'$
    \end{itemize} }
    \par{\quad{} Since $\gamma_{|[b,1]}\equiv \gamma'_{|[b,1]} $, the curves $\gamma$ and $\gamma'$ after exiting $\underset{i=l}{\overset{N}{\cup}} R_i$ are exactly the same; hence $W=W'$. We thus get the desired result by remarking that $R_l,w_A,R_{k}$ and $R_l,w'_A,R_{k}$ are homotopic by a sequence of homotopies of type $\mathcal{B}$, thanks to  Lemma \ref{l.increasingrectanglepaths}.}
        
       \vspace{0.5cm}  
    \item Assume finally that    $\gamma([a,1])\subset \underset{i=l}{\overset{N}{\cup}} R_i$. In this case, we will show directly, using the induction hypothesis, that the rectangle path associated to $\gamma$ is homotopic to a trivial or monotonous rectangle path. 
    
    \quad{}Consider $\gamma'$ the curve obtained by pushing $\gamma$ along $Dom$ so that $s$ be identified with $s'$. The curve $\gamma'$ is simple, good, closed and of strictly smaller length than $\gamma$. By eventually reparametrizing $\gamma'$, assume once again that $\gamma\equiv \gamma'$ on $[0,a]\cup[b,1]$; hence the rectangle paths starting from $R_0$ associated to $\gamma_{|[0,a]}$ and $ \gamma'_{|[0,a]} $ are identical. Recall that the last rectangle of the previous rectangle paths is $R_l$. Since $\gamma([a,1]),\gamma'([a,1]) \subset \underset{i=l}{\overset{N}{\cup}} R_i$, there exist $R_{k},R_{k'}\in \{R_l,...,R_N\}$ and $w_A,w'_A$ rectangle paths formed by rectangles in $\{R_l,...,R_N\}$ such that the rectangle paths starting from $R_l$ associated to  $\gamma'_{|[a,1]} $ and  $\gamma_{|[a,1]}$ are respectively  $R_l,w'_A,R_{k'}$ and $R_l,w_A,R_{k}$. Notice that $R_k$ is not necessarily equal to $R_{k'}$ and thus the rectangle paths starting from $R_0$ associated to $\gamma$ and $\gamma'$ are not homotopic. However, it is true that $\inte{R_k}\cap \inte{R_{k'}}\neq \emptyset$, since $\gamma(1)=\gamma'(1)$.

    \quad{}Assume without any loss of generality that $R_{k'}$ is a predecessor of some generation of $R_0$ (recall that $\gamma(0)=\gamma(1)\in \inte{R_0}$) and apply the induction hypothesis for $\gamma'$, which verifies $(M(\gamma'),l(\gamma'))<(M(\gamma),l(\gamma))$. We have that the rectangle path associated to $\gamma'$ is homotopic to the decreasing path  $R_0=r_0',r'_1,...,r'_p=R_{k'}$. Consider the rectangle path associated to $\gamma$, which is of the form:  $$R_0,R_1,...,R_l,w_A,R_{k}$$ By Lemma \ref{l.increasingrectanglepaths}, the previous rectangle path is homotopic by a sequence of homotopies of type $\mathcal{B}$ to the rectangle path $$R_0,R_1,...,R_l,w'_A,R_{k'},w'^{-1}_A, R_l, w_A,R_{k}$$ where $w'^{-1}_A$ is the rectangle path $w_A'$ followed from its end to its beginning. The first part of this rectangle path is by hypothesis homotopic to the decreasing rectangle path $R_0=r_0',r'_1,...,r'_p=R_{k'}$. Therefore, the previous rectangle path is homotopic to $$R_0=r_0',r'_1,...,r'_p=R_{k'},w'^{-1}_A, R_l, w_A,R_{k}$$ Recall that the second part of this rectangle path consists only of rectangles in $\{R_l,...,R_N\}$ and that $\inte{R_k}$ intersects $\inte{R_{k'}}$. If $R_k=R_{k'}$, then by Lemma \ref{l.increasingrectanglepaths} the path $R_{k'},w'^{-1}_A, R_l, w_A,R_{k}$ is homotopic to a trivial path and therefore 
 the rectangle path associated to $\gamma$ is homotopic to the decreasing path $R_0,r'_1,...,r'_{p}=R_{k}$, which gives us the desired result. If $R_{k}$ is a predecessor of some generation of $R_{k'}$, then by Lemma \ref{l.increasingrectanglepaths}, $R_{k'},w'^{-1}_A, R_l, w_A,R_{k}$ is homotopic to the decreasing rectangle path $r'_p=R_{k'}, r'_{p+1},...,r'_{p+s}=R_k$. We deduce in this case that the rectangle path associated to $\gamma$ is homotopic to the decreasing rectangle path $R_0=r_0',r'_1,...,r'_p,r'_{p+1},...,r'_{p+s}=R_k$, which gives us the desired result. Finally, assume that $R_k$ is a successor of some generation of $R_{k'}$. By the same argument, the rectangle path associated to $\gamma$ is homotopic to the concatenation of the decreasing rectangle path $R_0,r'_1,...,r'_p=R_{k'}$ and of an increasing rectangle path of the form $r'_p=R_{k'},r'_{p+1},...,r'_{p+s}=R_k$. However, by Lemma \ref{l.npredecessor}, $r'_p=R_{k'},r'_{p+1},...,r'_{p+s}=R_k$ (resp. $r'_p=R_{k'},r'_{p-1}...,r'_1,R_0$) is the unique increasing rectangle path from $r'_p$ to $R_k$ (from $r'_p$ to $R_0$). As $\inte{R_k}\cap \inte{R_0}\neq \emptyset$, thanks to Lemma \ref{l.intersectionsuccessorspredecessors} applied for $\clos{\mathcal{R}}$, we get that $r'_{p+1}=r'_{p-1}$. Therefore, by a homotopy of type $\mathcal{B}$, we can decrease the length of the path $R_0=r_0',r'_1,...r'_{p-1},r'_p,r'_{p+1},...,r'_{p+s}=R_k$ and we can keep doing this until we obtain an increasing,  decreasing or trivial rectangle path, which gives us the desired result.  
     \end{enumerate}  
\end{enumerate}


\vspace{0.5cm}
\subsubsection*{Induction step in the case where no tangency with strong $(\star)$ property exists}

Assume that the domains of all stable and unstable tangencies of $\gamma$ with property $(\star)$ contain $\gamma(0)$. Consider $s$ a tangency with property $(\star)$ and $Dom$ its domain (its existence is ensured by Lemma \ref{l.existencetangenciespropstar}). We will assume without any loss of generality that $s$ is a stable tangency. The function $\gamma:[0,1]\rightarrow \clos{\mathcal{P}}$ induces a cyclic order on the stable and unstable segments forming $\gamma$. Let us  denote once again by $U'$ (resp. $U$) the unstable segment of $\gamma$ before (resp. after) $s$, by $s'$ the stable side of $Dom$ that is not $s$, $u':=U'\cap Dom$ and $u:=U\cap Dom$. Recall that by definition of property $(\star)$, $\gamma(0)\notin s$ and $\gamma(0)$ is a corner point of $Dom$, therefore $\gamma(0)\in s' \cap U'$ or $\gamma(0)\in s' \cap U$.

\begin{lemm}\label{l.reductiontofirstsubcase} 
Let $R_0$ and $R_n$ be the first and last rectangles of the rectangle path associated to $\gamma$. If $\gamma(0)\in U'$ (resp. $\gamma(0)\in U$) and $U'\subset R_0\cap R_n$ (resp. $U\subset R_0\cap R_n$), then $\gamma$ is homotopic by a homotopy of type $\mathcal{A}$ to a simple, closed and good polygonal curve having at least one tangency with strong property $(\star)$. 
\end{lemm}
\begin{proof}
Assume that $\gamma(0)\in U'$. Let $\gamma':[0,1]\rightarrow \clos{\mathcal{P}}$ be a simple, closed, good polygonal curve defined by the following properties: $\gamma'(0)$ is the unique point in $U'\cap s$, the topological circles $\gamma([0,1])$ and $\gamma'([0,1])$ coincide, and they are identically oriented by $\gamma$ and $\gamma'$. Notice that since $U'\subset R_0$ and $\gamma$ is good, we have that $\gamma'(0)\in \inte{R_0}$. Moreover, by definition of $\gamma'$ and using the fact that $U'\subset R_n$, one can easily check that the rectangle paths starting from $R_0$ associated to the good polygonal curves $\gamma$ and $\gamma'$ coincide; hence $\gamma$ and $\gamma'$ are homotopic by a homotopy of type $\mathcal{A}$, which gives us the desired result by Lemma \ref{l.existencetangenciesgood}.
\end{proof}

We are now ready to describe the induction step in this final case: 
\begin{enumerate}[label=(\arabic*{$''$}),leftmargin=0.01cm, itemsep=0.1cm]
    \item If $(U',s,U)\rightarrow(0,0,0)$, by our argument in Case $(1')$,  we can push $\gamma$ by a homotopy of type $\mathcal{A}$ along $Dom$ so that $s$ comes arbitrarily close to $s'$. By Lemma \ref{l.reductiontofirstsubcase}, up to a homotopy of type $\mathcal{A}$, $\gamma$  has a tangency with strong property $(\star)$. We deduce the induction step by applying Cases $(1')-(8')$.
    \item If $(U',s,U)\rightarrow(1,0,1)$, as we have showed in Case $(2')$, we are in the case of Figure \ref{f.impossiblecasesab}b, with $u'=U'$ and therefore $\gamma(0)\in u'$. By our argument in (2$'$), we can push $\gamma$ by a homotopy of type $\mathcal{A}$ along $Dom$ so that $s$ comes arbitrarily close to $s'$. By doing so  $\textit{Rect}_{\gamma,R_0}(u')$ becomes trivial and therefore by Lemma \ref{l.reductiontofirstsubcase} up to a homotopy of type $\mathcal{A}$, $\gamma$ has a tangency with strong property $(\star)$. We deduce the induction step by applying Cases $(1')-(8')$.
    \item If $(U',s,U)\rightarrow(0,0,1)$, we are in the case of  Figure \ref{f.impossiblecasesab}b, with $R_l=R_{l+1}=...=R_N$. We apply in this case the same argument as in Case $(2'')$. 
    \item If $(U',s,U)\rightarrow(1,1,0)$, then we are in the case of Figure \ref{f.impossiblecasescsecond}, where $\{R_l,...,R_N\}=\textit{Rect}_{\gamma,R_0}(u')$,  $\{R_N,R_{N+1}
,...,R_m\}=\textit{Rect}_{\gamma,R_0}(s)$, $u=U$ and $\gamma(0)\in U$. By the same arguments as before, up to a homotopy of type $\mathcal{A}$, $\gamma$ has a tangency with strong property $(\star)$. We deduce the induction step by applying Cases $(1')-(8')$.
    \item If $(U',s,U)\rightarrow(0,1,0)$, then we are in the case of Figure \ref{f.impossiblecasescsecond} with $R_l=R_{l+1}=...=R_N$. In this case either $\gamma(0)\in u\cap s'$ or $\gamma(0)\in u'\cap s'$. In both cases, by the same arguments as before, up to a homotopy of type $\mathcal{A}$, $\gamma$ has a tangency with strong property $(\star)$. We deduce the induction step by applying Cases $(1')-(8')$.
\item If $(U',s,U)\rightarrow(1,1,1)$, then we are in the case of Figure \ref{f.impossiblecase1}. In this case, we can apply in the exact same way the homotopy performed in Case (6$'$). 
\item Assume that $(U',s,U)\rightarrow(0,1,1)$. Denote by $R_{u'},...,R_s$ the rectangle path $\textit{Rect}_{\gamma,R_0}(s)$. As we have shown in Case  ($7'$):
\begin{itemize}
    \item If $U'\not\subset R_s$, we are in the case of Figure \ref{f.impossiblecase1}. In this case, we can apply in the exact same way the homotopy performed in Case (7$'$).   
    \item If $U'\subset R_s$, then we are in the case of Figure \ref{f.case011}. In this case, $u'=U'$ and therefore $\gamma(0)\in u'$. Once again, we can push $\gamma$ by a homotopy of type $\mathcal{A}$ along $Dom$ so that $s$ comes arbitrarily close to $s'$. By Lemma \ref{l.reductiontofirstsubcase}, up to a  homotopy of type $\mathcal{A}$, $\gamma$ has a tangency with strong property $(\star)$. We deduce the induction step by applying Cases $(1')-(8')$.
     
\end{itemize}
    \item If $(U',s,U)\rightarrow(1,0,0)$, we are in the case of case of Figure \ref{f.impossiblecasesab}a, where $\textit{Rect}_{\gamma,R_0}(U')$ corresponds to the non-trivial decreasing rectangle path $R_l,...,R_N$.  Let $a,b\in[0,1]$ such that $\gamma([a,b])=U'\cup s \cup U$. \underline{Assume first that $\gamma(0)\in U' \cap s' $}. In this case, we have that $R_l=R_0$.  
    \begin{enumerate}[leftmargin=0.01cm, label=(8$''$\alph*)]
        
        \item If $\gamma_{|[b,1]}$ eventually exits $\underset{i=l}{\overset{N}{\cup}} R_i$ in order to visit a crossing predecessor (or successor) of $R_{k}\in \{R_l,...,R_N\}$, as in Case (8$'$a), we can push $\gamma$ along $Dom$  by a sequence of homotopies of type $\mathcal{B}$ so that $s$ comes arbitrarily close to $s'$. Hence, by Lemma \ref{l.reductiontofirstsubcase} up to a homotopy,  $\gamma$ has a tangency with strong property $(\star)$. We obtain the induction step by applying the Cases $(1')-(8')$.
        
    \item If $\gamma([a,1]) \subset \underset{i=l}{\overset{N}{\cup}} R_i$, since $\gamma(0)\in U'\cap s'$, we get that $\gamma([0,1]) \subset \underset{i=l}{\overset{N}{\cup}} R_i=\underset{i=0}{\overset{N}{\cup}} R_i$. Therefore, the rectangle path associated to $\gamma$ is a 
 formed by rectangles in $\{R_0,...,R_N\}$. The previous rectangle path is homotopic to a monotonous or trivial rectangle path, thanks to Lemma \ref{l.increasingrectanglepaths}. 
    
\item \underline{Assume now that $\gamma(0) \in s' \cap U$}. Notice that in this case, since $s$ does not satisfy the strong $(\star)$ property, we have that $\gamma(0)\in Dom$ and thus $U\subset Dom$ and  $U=u$. Push $\gamma$ along $Dom$ so that $s$ comes very close to $s'$. By Lemma \ref{l.reductiontofirstsubcase}, up to a homotopy of type $\mathcal{A}$, we can assume that this new simple, closed, good polygonal curve, say $\gamma'$,  has a tangency with strong property $(\star)$. By applying our algorithm (Cases $(1')-(8')$), $\gamma'$ is homotopic to a simple, good, closed polygonal curve $\gamma''$ of strictly smaller length or with less boundary arc points in its interior. By the induction hypothesis, the rectangle path associated to $\gamma''$ is homotopic to a monotonous or trivial path. Hence, the same applies to the rectangle path associated to $\gamma'$. By the exact same argument as in Case $(8'b)$, we can now conclude that the rectangle path associated to $\gamma$ is homotopic to a monotonous or trivial rectangle path. 
\end{enumerate}
\end{enumerate}
This finishes the description of our induction step. In all the previous cases, we have either defined a homotopy producing a simple, closed and good polygonal curve $\gamma'$ such that $(M(\gamma'), l(\gamma'))<(M(\gamma), l(\gamma))$ or we have shown that $\gamma$ satisfies the desired result. This finishes the proof by induction of Proposition \ref{p.homotopictotrivialsimplecase}.
\addtocontents{toc}{\protect\setcounter{tocdepth}{3}}
\subsubsection*{Proposition \ref{p.homotopictotrivialsimplecase} implies  Theorem \ref{t.homotopictotrivialpath}} 
Consider $r_0,...,r_n$ a closed rectangle path in $\clos{\mathcal{P}}$ and $\gamma$ a closed and good polygonal curve associated to the previous rectangle path (the existence of such a $\gamma$ is guaranteed by Proposition  \ref{p.rectanglepathtocurve}). We will show that $r_0,...,r_n$ is homotopic to a trivial rectangle path.

By our definition of a good polygonal curve, no two stable (resp. unstable) sides of $\gamma$ belong to the same leaf in $\clos{\mathcal{F}^s}$ (resp. $\clos{\mathcal{F}^u}$). Consequently, every self-intersection of $\gamma$ is contained in the intersection of a stable and an unstable segment of $\gamma$. By Item (4) of Proposition \ref{p.propertiesoffoliformarkovianactions}, this implies that $\gamma$ has only a finite number, say $M$, of self-intersections, all of which are topologically transverse. We will prove Theorem \ref{t.homotopictotrivialpath} by induction on $M$. If $M=0$, then we obtain the desired result from Proposition \ref{p.homotopictotrivialsimplecase}. Suppose now that $M>0$ and that the result holds for all $N\in \llbracket 0, M-1\rrbracket$.

Consider $x_0 \in [0,1]$ the biggest element in $[0,1]$ for which $\gamma_{|[0,x_0)}$ is injective. Take $y\in [0,x_0)$ such that $\gamma(y)=\gamma(x_0)$. Let us remark that $\gamma(x_0)=\gamma(y)$ lies in the intersection of a stable and unstable segment of $\gamma$, therefore it is not contained in the stable or unstable leaf of a boundary periodic point. Consider 
\begin{itemize}
    \item $r_0,...,r_k$ the rectangle path starting from $r_0$ associated to the simple and 
good polygonal curve $\gamma_{|[0,y]}$
\item $r_k,...,r_l$ the rectangle path starting from $r_k$ associated to the simple, closed and 
good polygonal curve $\gamma_{|[y,x_0]}$ (the curve $\gamma_{|[y,x_0]}$ is indeed good, since $\gamma$ is good and the intersection at $\gamma(x_0)$ is transverse)
\item $r_l,...,r_n$ the rectangle path starting from $r_l$ associated to the  
good polygonal curve $\gamma_{|[x_0,1]}$
\end{itemize} 
Notice that by our construction of the rectangle path associated to a good polygonal curve (see Definition  \ref{d.curvetorectanglepath}), the juxtaposition of the previous rectangle paths gives the rectangle path starting from $r_0$ associated to $\gamma$.

Since $\gamma_{|[y,x_0]}$ is simple and closed, by Proposition \ref{p.homotopictotrivialsimplecase} its associated rectangle path $r_k,...,r_l$ is homotopic to a monotonous or trivial rectangle path $$r_k=R_0,R_{1}...,R_s=r_l$$ Similarly, consider $\gamma'$ the good, closed polygonal curve formed by the juxtaposition of $\gamma([x_0,1])$ followed by $\gamma([0,y])$. Notice that $\gamma'$ has at most $M-1$ self-intersections and therefore by our induction hypothesis its associated rectangle path starting from $r_l$  $r_l,...,r_n=r_0,...,r_k$ is homotopic to monotonous or trivial rectangle path $$r_l=R'_0,R'_1...,R'_m=r_k$$

Consider now the rectangle path starting from $r_0$ associated to $\gamma$: $r_0,...,r_n$. By a sequence of homotopies of type $\mathcal{B}$ the previous rectangle path is homotopic to $$r_0,..r_k,...r_l,..,r_n=r_0,r_1,...,r_{k-1},r_k,r_{k-1},...,r_0$$
By our previous discussion the above rectangle path is homotopic to 
\begin{equation}\label{eq.rectpath}
r_0,..r_k=R_0,R_1...,R_s=r_l=R_0',R'_1...,R'_m=r_k,r_{k-1},...,r_0
\end{equation}

If $r_k=R_0,R_1...,R_s=r_l$ is the trivial rectangle path, then $r_k=r_l$ and $r_l=R'_0,R'_1...,R'_m=r_k$ is also the trivial rectangle path. In this case, $r_0,...,r_n$ is homotopic to the rectangle path $r_0,...,r_{k-1},r_k,r_{k-1}...,r_0$, which is homotopic to a trivial rectangle path by a sequence of homotopies of type $\mathcal{B}$. 

Assume now without any loss of generality that $r_k=R_0,R_1...,R_s=r_l$ is a decreasing rectangle path; hence $r_l$ is a predecessor of some generation of $r_k$. We deduce that $r_l=R_0',R'_1...,R'_m=r_k$ is an increasing rectangle path. Furthermore, by Lemma \ref{l.npredecessor}, we have that there exists a unique decreasing (resp. increasing) rectangle path from $r_k$ to $r_l$ (resp. from $r_l$ to $r_k$). We deduce that the rectangle paths $r_k=R_0,R_1...,R_s=r_l$ and $r_k=R_m', R_{m-1}'...,R'_1, R'_0=r_l$ are exactly the same. Therefore, by a sequence of homotopies of type $\mathcal{B}$ the rectangle path (\ref{eq.rectpath}) is homotopic to: 
$$r_0,...r_k,r_{k-1},..,r_0$$
Finally, the above rectangle path is homotopic to the trivial rectangle path by a sequence of homotopies of type $\mathcal{B}$; we thus get 
the desired result. 

\section{Proof of Theorem B}\label{s.mainresult}

Let $\rho_1:G_1 \rightarrow \text{Homeo}(\mathcal{P}_1)$, $\rho_2: G_2 \rightarrow  \text{Homeo}(\mathcal{P}_2)$ be two strong Markovian actions acting on the planes $\mathcal{P}_1$ and $\mathcal{P}_2$ endowed with their stable and unstable singular foliations $\mathcal{F}_{1}^{s,u}$,$\mathcal{F}_{2}^{s,u}$. Assume that $\rho_1$ and $\rho_2$ preserve  two strong Markovian families $\mathcal{R}_1$ and $\mathcal{R}_2$ whose associated classes of geometric types coincide. 

For every $i\in \{1,2\}$, let us denote by $\Gamma_{i}$ the boundary periodic points of $\mathcal{R}_{i}$, by $\clos{\mathcal{P}_{i}}$ the bifoliated plane of $\rho_i$ up to surgeries on $\Gamma_{i}$, by $\clos{\rho_i}:\clos{G_i}\rightarrow  \text{Homeo}(\clos{\mathcal{P}_{i}})$ the lift of $\rho_i$ on $\clos{\mathcal{P}_{i}}$, by $\clos{\mathcal{F}^{s}_{i}}, \clos{\mathcal{F}^{u}_{i}}$ the stable and unstable  foliations of $\clos{\rho_i}$, by $\clos{\mathcal{R}_{i}}$ the lift of $\mathcal{R}_{i}$ on $\clos{\mathcal{P}_{i}}$ and by $\clos{\Gamma_{i}}$ the lift of $\Gamma_{i}$ on $\clos{\mathcal{P}_{i}}$. 

Thanks to Corollary, \ref{c.specialgeomtypes}, by appropriately choosing representatives for every rectangle orbit in $\clos{\mathcal{R}_1}$ and $\clos{\mathcal{R}_2}$, and orientations for the foliations $\clos{\mathcal{F}_{1}^{s,u}}$,$\clos{\mathcal{F}_{2}^{s,u}}$, we may assume that $\clos{\mathcal{R}_1}$ and $\clos{\mathcal{R}_2}$ are associated to the same (special) geometric type $$\mathcal{G}=(R_1,...,R_n,(h_i)_{i \in \llbracket 1,n \rrbracket}, (v_i)_{i\in \llbracket 1,n \rrbracket}, \mathcal{H}, \mathcal{V},\phi, u)$$

We would like to show the following generalization of Theorem B in the setting of strong Markovian actions: 

\begin{theorem}[Theorem B']
    There exists a homeomorphism $h: \clos{\mathcal{P}_1}\rightarrow \clos{\mathcal{P}_2}$ and an isomorphism $\alpha: \clos{G_1} \rightarrow \clos{G_2}$ such that: 

\begin{itemize}
    \item the image by $h$ of any stable (resp. unstable) leaf in $\clos{\mathcal{F}_1^{s}}$ (resp. $\clos{\mathcal{F}_1^{u}}$) is a stable (resp. unstable) leaf in $\clos{\mathcal{F}_2^{s}}$ (resp. $\clos{\mathcal{F}_2^{u}}$)
    \item for every $g\in \clos{G_1}$ and every $x\in \clos{\mathcal{P}_1}$ we have that $$h(\clos{\rho_1}(g)(x))= \clos{\rho_2}(\alpha(g))(h(x))$$ 
\end{itemize}
\end{theorem}

We remark here that Theorem B is an immediate consequence of Theorems B' and C. We will split the proof of the previous theorem in two parts: 
\begin{enumerate}
    \item the construction of an isomorphism $\alpha: \clos{G_1} \rightarrow \clos{G_2}$
    \item the construction of a homeomoprhism $h: \clos{\mathcal{P}_1}\rightarrow \clos{\mathcal{P}_2}$ with the desired properties 
\end{enumerate}

Let us begin with the construction of $\alpha$. 

\begin{conv}
    By a small abuse of notation and in order to simplify our notations, from now on we will denote for any $g\in \clos{G_1}$ (resp. $g\in \clos{G_2}$) the homemomorphism $\clos{\rho_1}(g)$ (resp. $\clos{\rho_2}(g)$) by $g$. 

\end{conv}

Consider $r_0^1\in \clos{\mathcal{R}_1},r_0^2\in \clos{\mathcal{R}_2}$ two rectangles of the same type taken among our previously chosen set of representatives. We will fix $r_0^1$ and $ r_0^2$ as the origin rectangles in $ \clos{\mathcal{P}_1}$ and $ \clos{\mathcal{P}_2}$ respectively. Let $\mathcal{O}(r_0^1)\subset \clos{\mathcal{R}_1}$ be the $\clos{G_1}$-orbit of $r_0^1$ in $ \clos{\mathcal{P}_1}$ and   $g^1\in\clos{G_1} $. Consider  $R_0^1:=r_0^1,...,R_n^1=g^1(R_0^1)$ a centered rectangle path in $\clos{\mathcal{R}_1}$ going from $r_0^1$ to $g^1(r_0^1)$  (such a rectangle path exists thanks to Proposition \ref{p.rectpathsstartingending}) and $R_0^2:=r_0^2,...,R_n^2$ its associated centered rectangle path in $ \clos{\mathcal{R}_2}$. Since $R_0^1, g^1(R_0^1), R_0^2, R_n^2$ are by definition of the associated rectangle path (see Definition  \ref{d.associatingrectanglepaths}) of the same type, there exists $g^2\in \clos{G_2}$ such that $R_n^2=g^2(R_0^2)$. We define 

\begin{align*}
\alpha: \clos{G_1} &\rightarrow  \clos{G_2}\\
 g^1 &\rightarrow g^2
\end{align*}
 
\begin{prop}\label{p.existenceiso}
The map $\alpha$ does not depend on our previous choices of rectangle paths and defines an isomorphism between $\clos{G_1}$ and $ \clos{G_2}$. 
\end{prop}
\begin{proof}

Let us begin by showing that the map $\alpha$ is well defined. Indeed, if $g\in \clos{G_1}$ fixes one rectangle $R$ in $\mathcal{O}(r_0^1)$, then it has a fixed point in $R$ and by Proposition \ref{p.kernelprojectiongroup}  we get that  $g=id$. A similar result is true for the elements in $\clos{G_2}$. Consequently, following our previous notations, $g^1$ (resp. $g^2$) is the unique element in $\clos{G_1}$ (resp. $\clos{G_2}$) that sends $r_0^1$ (resp. $r_0^2$) to $R_n^1$ (resp. $R_n^2$). 

Furthermore, thanks to Theorem \ref{t.endingbysamerectangle}, the rectangle $R_n^2$ does not depend on our initial choice of rectangle path going from $r_0^1$ to $R_n^1=g^1(r_0^1)$. It follows that $g^2$ also does not depend on our initial choice of rectangle path and thus $\alpha$ is well defined. Let us now show that $\alpha$ defines an isomorphism. 

Fix $L^1_1, L_2^1\in \mathcal{O}(r_0^1)$. Consider $G_1^1, G^1_2 \in \clos{G_1}$ the unique elements of $\clos{G_1}$ sending $r_0^1$  to $L^1_1, L_2^1$ respectively. Let $r_0^1, r_1^1,...,r_n^1=L^1_1$  (resp. ${r'}_0^1:=r_0^1, {r'}_1^1,...,{r'}_m^1=L^1_2$) be a centered rectangle path in $\clos{\mathcal{R}_1}$ going from $r_0^1$ to $L^1_1$ (resp. $L_2^1$)  and $r_0^2, r_1^2,...,L^2_1:=r_n^2$ (resp. ${r'}_0^{2}:=r_0^2, {r'}_1^2,...,L^2_2:={r'}_m^2$) its associated centered rectangle path in $\clos{\mathcal{R}_2}$.  By definition of $\alpha$, $\alpha(G_1^1)$ (resp. $\alpha(G_2^1)$) is the unique element of $\clos{G_2}$ sending $r_0^2$ to $L_1^2$ (resp. $L_2^2$). 

Consider now the rectangle path $r_0^1,r_1^1,...,L_1^1=G_1^1(r_0^1),G_1^1({r'}_1^1),G_1^1({r'_2}^1),...,G_1^1(L_2^1)=G_1^1 \circ G_2^1 (r_0^1)$. It is not difficult to see that the associated rectangle path in $\clos{\mathcal{R}_2}$ of the previous previous path is $r_0^2,r_1^2,...,L_1^2=G_1^2(r_0^2),G_1^2({r'}_1^2),...,G_1^2(L_2^2)=G_1^2 \circ G_2^2 (r_0^2)$. It follows by the definition of $\alpha$ that $\alpha(G_1^1 \circ G_2^1)= G_1^2 \circ G_2^2=\alpha(G_1^1) \circ \alpha(G_2^1)  $ and thus $\alpha$ defines a morphism from $\clos{G_1}$ to $\clos{G_2}$. By changing the roles of $\clos{G_1}$ and $\clos{G_2}$, we can produce by the same exact procedure a morphism $\alpha': \clos{G_2} \rightarrow \clos{G_1}$ such that $\alpha$ and $\alpha'$ are inverses of each other. Consequently, $\alpha$ defines an isomorphism from  $\clos{G_1}$ to $\clos{G_2}$ and we get the desired result.

\end{proof}
We are now ready to proceed to the construction of the homeomorphism $h: \clos{\mathcal{P}_1}\rightarrow \clos{\mathcal{P}_2}$. We are first going to define a bijection $H$ from  $\clos{\mathcal{R}_1}$ to $\clos{\mathcal{R}_2}$ that is equivariant with respect to the group actions and then we are going to show that a map $h: \clos{\mathcal{P}_1}\rightarrow \clos{\mathcal{P}_2}$ such that $h(R)=H(R)$ for every rectangle $R\in \clos{\mathcal{R}_1}$ defines a unique homeomorphism from $\clos{\mathcal{P}_1}$ to $\clos{\mathcal{P}_2}$ that has the desired properties. 

Indeed, let $R\in \clos{\mathcal{R}_1}$. Consider a centered rectangle path in $\clos{\mathcal{P}_1}$ ending at $R$ and its corresponding centered rectangle path in $\clos{\mathcal{P}_2}$ ending at $R'\in \clos{\mathcal{R}_2}$. We define $H(R)=R'$. Notice that $R'$ does not depend on our initial choice of rectangle path thanks to Theorem \ref{t.endingbysamerectangle} and that $H$ sends predecessors (resp. successors)  to predecessors (resp. successors). Also, by our definition of $\alpha$, we have that for any $g\in \clos{G_1}$ and $R\in \clos{\mathcal{R}_1}$: \begin{equation}\label{eq.equivarianceh}
    H(g(R))=\alpha(g)(H(R))
\end{equation} Finally, by changing the roles of $\clos{\mathcal{R}_1}$, $\clos{\mathcal{R}_2}$ and by applying the same arguments as before, we can construct a right and left inverse for $H$. It follows that $H: \clos{\mathcal{R}_1}\rightarrow  \clos{\mathcal{R}_2}$ is a bijection that is equivariant with respect to the actions of $\clos{G_1}, \clos{G_2}$ . 

In order to prove Theorem B' it suffices to prove the following proposition: 
\begin{prop}\label{p.existencehomeo}
There exists a unique homeomorphism $h$ from $\clos{\mathcal{P}_1}$ to $\clos{\mathcal{P}_2}$ such that: 

\begin{itemize}
    \item for every $R\in \clos{\mathcal{R}_1}$ we have $h(R)=H(R)$ 
    \item the image by $h$ of any stable (resp. unstable) leaf in $\clos{\mathcal{F}_1^{s}}$ (resp. $\clos{\mathcal{F}_1^{u}}$) is a stable (resp. unstable) leaf in $\clos{\mathcal{F}_2^{s}}$ (resp. $\clos{\mathcal{F}_2^{u}}$)
    \item for every $g\in \clos{G_1}$ and every $x\in \clos{\mathcal{P}_1}$ we have $$h(g(x))= \alpha(g)(h(x))$$ 
\end{itemize}
\end{prop}
\begin{proof}
Recall that in the beginning of this section we fixed a choice of representatives of every rectangle orbit in $\clos{\mathcal{R}_1}$, $\clos{\mathcal{R}_2}$ and a choice of orientations on $\clos{\mathcal{F}^{s,u}_1}$, $\clos{\mathcal{F}^{s,u}_2}$ such that the (special) geometric types  associated to $\clos{\mathcal{R}_1}$, $\clos{\mathcal{R}_2}$ for the previous choices of representatives and orientations are the same. 

Let $x\in\clos{\mathcal{P}_1}$. Recall that a quadrant of $x$ is  the closure of a connected component of $\clos{\mathcal{P}_1}$ minus the union of the stable and unstable leaves crossing $x$. For any quadrant $Q$ of $x$, there exists a unique (up to reindexation) bi-infinite sequence $...R^{Q}_{-1}, R^{Q}_0,R^{Q}_1,...$ of rectangles in $\clos{\mathcal{R}_1}$ such that: 
\begin{itemize} 
\item for every $i$, $R^{Q}_{i+1}$ is a predecessor of $R^{Q}_i$ 
\item for every $i$, $R^{Q}_i$ contains a neighborhood of $x$ in $Q$
\end{itemize}
Indeed, by the definition of a Markovian family there exists one rectangle $R^{Q}_0$ containing a neighborhood of $x$ in the quadrant $Q$. By Lemma \ref{l.existenceofpredecessors} applied to $\clos{\mathcal{R}_1}$ (see also Remark \ref{r.equivdefi}), there exists a unique predecessor and a unique successor of $R^{Q}_0$ containing a neighborhood of $x$ in $Q$. We will denote the previous rectangles by $R^Q_{1}$ and $R^Q_{-1}$ respectively. By a repeated application of this argument, we can construct a bi-infinite sequence of rectangles in $\clos{\mathcal{R}_1}$ with the desired properties. Furthermore, by Lemma \ref{l.npredecessor} applied for $\clos{\mathcal{R}_1}$, any rectangle in $\clos{\mathcal{R}_1}$ containing a neighborhood of $x$ in $Q$ is a predecessor or a successor of some generation of $R^{Q}_0$. In other words, any such rectangle appears in our previously constructed bi-infinite sequence. It follows that $...R^{Q}_{-1}, R^{Q}_0,R^{Q}_1,...$ is uniquely defined by the above properties up to a reindexation of the rectangles. 

By Lemmas \ref{l.infiniteintersectionverticalrectangles} and \ref{l.infiniteintersectionhorizontalrectangles} applied for $\clos{\mathcal{R}_1}$, we have that $\{x\}=\overset{+\infty}{\underset{k=-\infty}{\cap}} R^{Q}_{k}$. We now define $h(x):=\overset{+\infty}{\underset{k=-\infty}{\cap}} H(R^{Q}_{k})\in \clos{\mathcal{P}_2}$. Notice that by definition of $H$, for every $k\in \mathbb{Z}$ the rectangle $H(R^{Q}_k)$ is a predecessor of $H(R^{Q}_{k-1})$. Therefore by the Lemmas \ref{l.infiniteintersectionverticalrectangles} and \ref{l.infiniteintersectionhorizontalrectangles}, the set $\overset{+\infty}{\underset{k=-\infty}{\cap}} H(R^{Q}_{k})$ corresponds to a unique point in $\clos{\mathcal{P}_2}$. We would now like to show that this point does not depend on the initial choice of quadrant for $x$. 

\vspace{0.5cm}
\underline{Independence from the choice of quadrant} 

Assume first that $x$ is not a boundary periodic point of $\clos{\mathcal{R}_1}$ or equivalently that $x\notin \clos{\Gamma_1}$. Denote by $\clos{(\mathcal{F}^{s}_1)_+}(x)$ (resp. $\clos{(\mathcal{F}^{u}_1)_+}(x)$, $\clos{(\mathcal{F}^{s}_1)_-}(x)$,  $\clos{(\mathcal{F}^{u}_1)_-}(x)$) the positive (resp. positive, negative, negative) stable (resp. unstable, stable, unstable) separatrix of $x$ given by our choice of orientation of $\clos{\mathcal{F}^s_1}$ and $\clos{\mathcal{F}^u_1}$. For every $\epsilon,\epsilon' \in \{+,-\}$, denote also by $(\epsilon,\epsilon')$ the quadrant of $x$ bounded by $\clos{(\mathcal{F}^s_1)_\epsilon}(x)$ and $\clos{(\mathcal{F}^u_1)_{\epsilon'}}(x)$. Assume that the quadrant $Q$ was the $(\epsilon,\epsilon')$ quadrant of $x$, where $\epsilon,\epsilon' \in \{+,-\}$. Consider the quadrant $(\epsilon,-\epsilon')$ of $x$ and $$i:= \inf \{k\in\mathbb{Z}| \forall j\geq k ~ R^{(\epsilon, \epsilon')}_j \text{ contains a germ of the $(\epsilon,\epsilon')$ and $(\epsilon,-\epsilon')$ quadrants of $x$ }\}$$ Let us first show that $i$ belongs in $ \mathbb{Z}\cup \{-\infty\}$. Indeed, suppose that $R^{(\epsilon, \epsilon')}_0$ contains a germ of the $(\epsilon,\epsilon')$ quadrant of $x$, but not of the $(\epsilon,-\epsilon')$ quadrant of $x$. This implies that $x$ belongs to a stable boundary component of $R^{(\epsilon, \epsilon')}_0$, say $s$. By Lemmas  \ref{l.crossingrectanglesnoperiodicpoints},  \ref{l.crossingrectangleswithperiodicpoints} there exists a unique $s$-crossing predecessor of $R^{(\epsilon, \epsilon')}_0$ that contains a germ of the $(\epsilon,\epsilon')$ and $(\epsilon,-\epsilon')$ quadrants of $x$. Next, by Lemma \ref{l.npredecessor}, there exists $k_0 \in \mathbb{N}$ such that the previous $s$-crossing predecessor is of the form $R^{(\epsilon, \epsilon')}_{k_0}$. By the Markovian intersection property, for every $k\geq k_0$, $R^{(\epsilon, \epsilon')}_{k}$ also contains a germ of the $(\epsilon,\epsilon')$ and $(\epsilon,-\epsilon')$ quadrants of $x$. We thus get that $i \in \mathbb{Z}\cup \{-\infty\}$. 

If $i=-\infty$, then the sequence associated to the quadrant $(\epsilon,-\epsilon')$ is the exact same sequence as for $(\epsilon,\epsilon')$. This can happen for instance when $x$ does not belong to the stable or unstable leaf of some boundary periodic point. Therefore, in this case we get $\overset{+\infty}{\underset{k=-\infty}{\cap}} H(R^{(\epsilon,\epsilon')}_{k})= \overset{+\infty}{\underset{k=-\infty}{\cap}} H(R^{(\epsilon,-\epsilon')}_{k})$. 

If $i\in \mathbb{Z}$, then for every $j\geq i$ the rectangle $R_j^{(\epsilon,\epsilon')}$ is contained in the bi-infinite sequence associated to the quadrant $(\epsilon,-\epsilon')$. We therefore get up to changing our initial choice of indexes,  $R^{(\epsilon, -\epsilon')}_i=R^{(\epsilon, \epsilon')}_i,R^{(\epsilon, -\epsilon')}_{i+1}=R^{(\epsilon, \epsilon')}_{i+1},...,R^{(\epsilon, -\epsilon')}_n=R^{(\epsilon, \epsilon')}_n,...$.
The rectangle $R^{(\epsilon, -\epsilon')}_{i-1}$ is the unique successor of $R^{(\epsilon, -\epsilon')}_i=R^{(\epsilon, \epsilon')}_i$ that contains $x$ and that is not equal to $R^{(\epsilon, \epsilon')}_{i-1}$. Since any two distinct successors of $R^{(\epsilon, -\epsilon')}_i=R^{(\epsilon, \epsilon')}_i$ have disjoint interiors, we have that $x\in \partial^s{R^{(\epsilon, \epsilon')}_{i-1}}\cap \partial^s{R^{(\epsilon, -\epsilon')}_{i-1}}$. 

Assume without any loss of generality that $R^{(\epsilon, -\epsilon')}_{i-1}$ is the successor of $R^{(\epsilon, \epsilon')}_i$ that is right above $R^{(\epsilon, \epsilon')}_{i-1}$ (see Figure \ref{f.hmap}). In this case, thanks to our choice of representatives and orientations and to our construction of $H$, the map $H$ will send $R^{(\epsilon, -\epsilon')}_{i-1}$ to the successor of $H(R^{(\epsilon, -\epsilon')}_i)$ that is right above $H(R^{(\epsilon, \epsilon')}_{i-1})$. Furthermore, since $R^{(\epsilon, -\epsilon')}_{i-2}$ contains $x$ and a germ of its $(\epsilon,-\epsilon')$ quadrant, $R^{(\epsilon, -\epsilon')}_{i-2}$  is the unique successor of $R^{(\epsilon, -\epsilon')}_{i-1}$ that contains $x$, which corresponds to the lowermost successor of $R^{(\epsilon, -\epsilon')}_{i-1}$. Therefore, $H(R^{(\epsilon, -\epsilon')}_{i-2})$ will correspond to the lowermost successor of $H(R^{(\epsilon, -\epsilon')}_{i-1})$ or in other words the successor of $H(R^{(\epsilon, -\epsilon')}_{i-1})$ that contains $h(x)=\overset{+\infty}{\underset{k=-\infty}{\cap}} H(R^{(\epsilon,\epsilon')}_{k})\in \clos{\mathcal{P}_2}$. Similarly, one can prove that for every $k< i$ the rectangle  $H(R^{(\epsilon, -\epsilon')}_{k})$ contains $h(x)=\overset{+\infty}{\underset{k=-\infty}{\cap}} H(R^{(\epsilon,\epsilon')}_{k})\in \clos{\mathcal{P}_2}$. Thanks to Lemma \ref{l.infiniteintersectionhorizontalrectangles}, we have that  $\overset{i-1}{\underset{k=-\infty}{\cap}} H(R^{(\epsilon, -\epsilon')}_{k})$ corresponds to a stable segment in $\clos{\mathcal{F}_2^s}$ crossing $H(R^{(\epsilon, -\epsilon')}_{i-1})$ and containing $h(x)$. 
\begin{figure}
    \centering
    \includegraphics[scale=0.6]{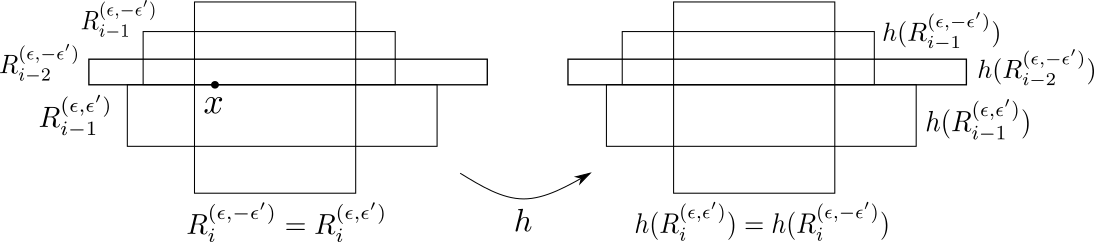}
    \caption{}
    \label{f.hmap}
\end{figure}

Thanks to Lemma \ref{l.infiniteintersectionverticalrectangles}, we also know  that $\overset{+\infty}{\underset{k=i}{\cap}} H(R^{(\epsilon, -\epsilon')}_{k})=\overset{+\infty}{\underset{k=i}{\cap}} H(R^{(\epsilon, \epsilon')}_{k})$ corresponds to an unstable segment in $\clos{\mathcal{F}^u}$ crossing $H(R^{(\epsilon, -\epsilon')}_{i})=H(R^{(\epsilon, \epsilon')}_{i})$ and containing $h(x)$. We conclude that $\overset{\infty}{\underset{k=-\infty}{\cap}} H(R^{(\epsilon, -\epsilon')}_{k})=\{h(x)\}$. By a similar argument we can show that $\overset{\infty}{\underset{k=-\infty}{\cap}} H(R^{(-\epsilon, -\epsilon')}_{k})=\overset{\infty}{\underset{k=-\infty}{\cap}} H(R^{(-\epsilon, \epsilon')}_{k})= \{h(x)\}$, which gives us the desired result. 

Assume now that $x$ is a boundary periodic point or equivalently that $x\in \clos{\Gamma_1}$. In this case, we can define $h(x)$ by the same exact argument for some choice of quadrant $Q$ of $x$ ($x$ has infinitely many quadrants in $\clos{\mathcal{P}_1}$). Consider now $Q'\neq Q$ the quadrant of $x$ that intersects $Q$ along a stable leaf of $x$. By a similar argument, we can show that for every $N$, $\overset{N}{\underset{k=-\infty}{\cap}} H(R^{Q'}_{k})$ corresponds to a stable segment in $\clos{\mathcal{F}_2^s}(h(x))$. Furthermore, if $t\in \text{Stab}(x)$ acts as a contraction on each stable leaf of $x$, then we have that for every $k\in \mathbb{Z}$ and every quadrant $q$ of $x$, $t^k(R^{q}_{0})$ intersects a neighborhood of $x$ in $q$ and therefore  $(t^k(R^{q}_{0}))_{k\in\mathbb{Z}}$ is a subsequence of $ (R^{q}_{k})_{k\in\mathbb{Z}}$. It follows that $\big(\alpha(t)^k(H(R^{q}_{0}))\big)_{k\in\mathbb{Z}}$ is also a subsequence of $ (H(R^{q}_{k}))_{k\in\mathbb{Z}}$. Notice that by Lemmas \ref{l.infiniteintersectionverticalrectangles} and \ref{l.infiniteintersectionhorizontalrectangles}, the sets $\overset{+\infty}{\underset{k=-\infty}{\cap}}\alpha(t)^k(H(R^{Q'}_{0}))$, $\overset{+\infty}{\underset{k=-\infty}{\cap}}\alpha(t)^k(H(R^{Q}_{0}))$ consist each of a unique point fixed by $\alpha(t)$. Our previous arguments imply that both $\overset{+\infty}{\underset{k=-\infty}{\cap}} H(R^{Q'}_{k})$ and $\overset{+\infty}{\underset{k=-\infty}{\cap}} H(R^{Q}_{k})$ are fixed by $\alpha(t)$ and belong to the same stable leaf in $\clos{\mathcal{F}^s}$. By Proposition \ref{p.propertiespbaraction}, we conclude that the previous two points correspond to the same point $h(x)$, which gives us the desired result.

\vspace{0.5cm}
\underline{$h(\clos{\Gamma_1})\subset \clos{\Gamma_2}$ and $h(\clos{\mathcal{P}_1}-\clos{\Gamma_1})\subset \clos{\mathcal{P}_2}-\clos{\Gamma_2}$ }

Indeed, assume that  $x\in\clos{\Gamma_1}$ and $h(x)\in \clos{\mathcal{P}_2}-\clos{\Gamma_2}$. There are infinitely many quadrants around $x$, therefore there exists an infinite  family of rectangles $(R_i)_{i\in \mathbb{N}}$ of $\clos{\mathcal{R}_1}$ intersecting $x$ and having pairwise disjoint interiors. For every $i$ the rectangle $H(R_i)$ intersects $h(x)$ that has only finitely many quadrants. Therefore, there exist $i\neq j$ such that $\inte{H(R_i)}\cap \inte{H(R_j)} \neq \emptyset$. Therefore, $H(R_i)$ is a predecessor or successor of some generation of $H(R_j)$. This contradicts that the fact that $H$ defines a bijection from $\clos{\mathcal{R}_1}$ to $\clos{\mathcal{R}_2}$ that respects the relations of predecessor/successor. 

By a similar argument we can show that $h(\clos{\mathcal{P}_1}-\clos{\Gamma_1})\subset \clos{\mathcal{P}_2}-\clos{\Gamma_2}$. 
\vspace{0.5cm}

\underline{$h$ is a bijection for which $h(R)=H(R)$ for all $R\in \clos{\mathcal{R}_1}$}

By changing the roles of $\clos{\mathcal{P}_1}$ and  $\clos{\mathcal{P}_2}$ and using the fact that $H$ is a bijection between $\clos{\mathcal{R}_1}$ and $\clos{\mathcal{R}_2}$, we can show that $h$ admits a right and left inverse, therefore $h$ is a bijection. 

Next, take $R\in \clos{\mathcal{R}_1}$. For every $x\in R$ there exists a quadrant $Q$ of $x$ such that $R$ contains a neighborhood of $x$ in $Q$. Therefore, $R$ belongs in the sequence of rectangles $(R^{Q}_{k}(x))_{k\in \mathbb{Z}}$ defined in the previous paragraphs and by definition of $h$, we get that $h(x)\in H(R)$ and by extension that $h(R)\subset H(R)$. By a similar argument, we can show that $h^{-1}(H(R))\subset R$. By combining the two previous results, we get that $h(R)=H(R)$.

\underline{$h$ is a homeomorphism}

Recall that by our previous arguments $h$ is a bijection. Consider now $x\in \clos{\mathcal{P}_1}-\clos{\Gamma_1}$. Denote by $(R^{(\epsilon, \epsilon')}_{k}(x))_{k\in \mathbb{Z}}$ the sequence of rectangles defined in the previous paragraphs for the quadrant $(\epsilon, \epsilon')$ of $x$. In order to show that $h$ is continuous on $x$, it suffices to show that for every $(\epsilon, \epsilon')\in \{+,-\}^2$ and for any sequence $x_n$ converging to $x$ in the $(\epsilon, \epsilon')$ quadrant of $x$, $h(x_n)$ converges to $h(x)$. Fix  $\epsilon, \epsilon'\in \{+,-\}$ and consider a sequence $x_n$ in the $(\epsilon, \epsilon')$ quadrant of $x$ converging to $x$. For any $N\in \mathbb{N}$, if $n$ is sufficiently big $x_n\in \overset{N}{\underset{k=-N}{\cap}} R^{(\epsilon, \epsilon')}_{k}(x) $, hence $h(x_n) \in \overset{N}{\underset{k=-N}{\cap}} H(R^{(\epsilon, \epsilon')}_{k}(x))$. By Lemmas \ref{l.infiniteintersectionverticalrectangles} and \ref{l.infiniteintersectionhorizontalrectangles}, the set $ \overset{N}{\underset{k=-N}{\cap}} H(R^{(\epsilon, \epsilon')}_{k}(x))$ converges when $N\rightarrow +\infty$ to $h(x)$ for the Hausdorff topology. We thus get that $h(x_n)$ converges to $h(x)$.

 Consider now $x \in \clos{\Gamma_1}$. Recall that $x$ has infinitely many quadrants. By an argument similar to the one we used in the previous paragraph, if $(x_n)_{n\in\mathbb{N}}$ converges to $x$ and $(x_n)_{n\in\mathbb{N}}$ visits only a finite number of quadrants of $x$, then $(h(x_n))_{n\in\mathbb{N}}$ converges to $h(x)$. In order to show that $h$ is continuous on $x$, it suffices to show that if $(x_n)_{n\in\mathbb{N}}$ visits an infinite number of quadrants of $x$, then $(x_n)_{n\in\mathbb{N}}$ can not converge to $x$. Indeed, take $(x_n)_{n\in\mathbb{N}}$ such a sequence. By eventually considering a subsequence, we can assume that for every $n$, $x_n$ belongs in a different quadrant of $x$, say $q_n$. For every quadrant $q$ of $x$ that belongs in $\{q_n|n\in \mathbb{N}\}$, take $U_q$ a neighborhood of $x$ in $q$ that does not contain $x_n$. For every quadrant $q$ of $x$ such that $q\notin \{q_n|n\in \mathbb{N}\}$, take $U_q$ any neighborhood of $x$ in $q$. By Proposition \ref{p.baseoftopology}, the union of the $U_q$ contains an open set of $x$ that by construction does not intersect $x_n$. Hence, $x_n$ does not converge to $x$. 
   
   The above prove that $h$ is continuous. By applying the same arguments to $h^{-1}$, we get that $h$ is a homeomorphism. 

\vspace{0.5cm}

\underline{$h$ is equivariant with respect to the actions of the group actions}

Consider $x\in \clos{\mathcal{P}_1}$ and $g\in \clos{G_1}$. For any quadrant $q$ of $x$, using our previous notations, we have that: $$h(g(x))=h(g(\overset{+\infty}{\underset{l=-\infty}{\cap}} R^{q}_{l}(x)))= \overset{+\infty}{\underset{l=-\infty}{\cap}} H(g(R^q_l(x)))=\overset{+\infty}{\underset{l=-\infty}{\cap}} \alpha(g)(H(R^q_l(x)))=\alpha(g)(h(x))$$

\vspace{0.5cm}
\underline{$h$ sends oriented stable/unstable leaves to oriented stable/unstable leaves}

Consider $p$ a periodic point in $\clos{\mathcal{P}_1} - \clos{\Gamma_1}$,  $x\in \clos{\mathcal{F}_1^s}(p)$ and $g\in \text{Stab}(p)$ such that $g^n(x)\rightarrow p$ when $n\rightarrow +\infty$. By eventually replacing $g$ by $g^2$, assume that $g$ preserves the orientations of the foliations $\clos{\mathcal{F}_1^{s,u}}$. 

Since $h$ is a homeomorphism that is equivariant with respect to the  actions of $\clos{G_1}$ and $\clos{G_2}$,  we have that $\alpha(g)\in \text{Stab}(h(p))$ and $\alpha(g)^n(h(x))\rightarrow h(p)$ when $n\rightarrow +\infty$. Furthermore, for any rectangle $R\in \clos{\mathcal{R}_1}$ containing $p$, we get by the definition of $g$ and by Lemma \ref{l.npredecessor} that $g(R)$ is a predecessor of some generation of $R$. Since $H$ and thus also $h$ sends predecessors to predecessors, this implies that $\alpha(g)(h(R))$ is a predecessor of some generation of $h(R)$ and therefore $\alpha(g) \in \text{Stab}(h(p))$ acts as a contraction on $\clos{\mathcal{F}_2^s}(h(p))$. Finally, since  $\alpha(g)^n(h(x))\rightarrow h(p)$ when $n\rightarrow +\infty$, we get that $h(x)\in \clos{\mathcal{F}_2^s}(h(p))$. We have thus proven that $h(\clos{\mathcal{F}_1^s}(p))\subset\clos{\mathcal{F}_2^s}(h(p))$. By using the same argument for $h^{-1}$, we deduce that $h(\clos{\mathcal{F}_1^s}(p))=\clos{\mathcal{F}_2^s}(h(p))$. 

By similar arguments one can prove that $h(\clos{\mathcal{F}_1^u}(p))=\clos{\mathcal{F}_2^u}(h(p))$ and also that $h$ sends a  stable (resp. unstable) leaf of a point 
$p'\in \clos{\Gamma_1}$ to a stable (resp. unstable) leaf of the point  $h(p')\in \clos{\Gamma_2}$. Moreover, by Definition \ref{d.markovianaction}, the periodic stable (resp. unstable) leaves form a dense set in the leaf space of $\clos{\mathcal{F}_1^s}$ (resp. $\clos{\mathcal{F}_1^u}$). Using the fact that $h$ is a homeomorphism, our previous arguments imply that  $h(\clos{\mathcal{F}_1^s})=\clos{\mathcal{F}_2^s}$ and $h(\clos{\mathcal{F}_1^u})=\clos{\mathcal{F}_2^u}$. 

Finally, recall that  $\clos{\mathcal{F}_1^{s,u}}$ and  $\clos{\mathcal{F}_2^{s,u}}$ have been endowed with orientations for which the map $H$ sends the predecessors of any rectangle $r\in \clos{\mathcal{R}_1}$ ordered from left to right, to the predecessors of $H(r)$ ordered from left to right. Same for the successors of any rectangle in $\clos{\mathcal{R}_1}$. We deduce that $h$ respects the orientations of $\clos{\mathcal{F}_1^{s,u}}$ and  $\clos{\mathcal{F}_2^{s,u}}$.

\end{proof}

\section{An invariant of a strong Markovian action up to conjugacy}\label{s.theoremsDE}
Fix $\rho: G\rightarrow \text{Homeo}(\mathcal{P})$ an orientation preserving strong Markovian action, preserving the pair of singular foliations $\mathcal{F}^s$ and $\mathcal{F}^u$ and leaving invariant a strong Markovian family $\mathcal{R}$. Denote by $\Gamma$ the set of boundary periodic points of $\mathcal{R}$.


In Section \ref{s.homotopiesofpaths}, we introduced the notion of cycle around a boundary arc point of $\clos{\mathcal{P}}$. We would like in this section to adapt this notion of cycle for boundary periodic points in $\mathcal{P}$, in order to obtain information about the neighborhood of these orbits in $\mathcal{P}$. This will allow us to prove the analogues of Theorem-Definition D and Theorem E for strong Markovian actions. 

According to the following lemma, we can not extend Definition \ref{d.cyclearcpoint} for boundary periodic points in $\mathcal{P}$.
\begin{lemm}\label{l.nocyclearoundperiod}
Let $\gamma\in \Gamma$ and $R,R'\in\mathcal{R}$ two rectangles containing $\gamma$ such that $\inte{R}\cap \inte{R'}\neq \emptyset$. For any  quadrant of $\gamma$, say $Q$, we have that $R$ contains a neighborhood of $\gamma$ in $Q$ if and only if $R'$ contains a neighborhood of $\gamma$ in $Q$. 
\end{lemm}
\begin{proof}
Assume that $R$ contains a neighborhood of $\gamma$ in $Q$, but not $R'$. Since $\inte{R}\cap \inte{R'}\neq \emptyset$ and $\gamma\in R\cap R'$, by the Markovian intersection axiom, both $R$ and $R'$ must contain the a neigbhorhood of $\gamma$ in some quadrant $Q'$ of $\gamma$ that is adjacent to $Q$. The previous facts imply that $R$ must be a predecessor or successor of $R'$ of some generation crossing a stable or unstable boundary of $R$. Assume without any loss of generality, that $R$ is a predecessor of $R'$ of some generation crossing its stable boundary. Since $\gamma\in R\cap R'$ and $R'$ does not contain a neighborhood of $\gamma$ in $Q$, we have that $\gamma \in \partial^s R'\cap \partial^u R$.  Consider now $g\in \text{Stab}(\gamma)$ such that the stable boundaries of $\rho(g)(R')$ become very thin and the unstable ones very long. The intersection of the rectangles $\rho(g)(R')$ and $R$ is not Markovian, which leads to an absurd.
\end{proof}

\begin{lemm}\label{l.finite}
    Take $\gamma\in \Gamma$ and $q$ one of its quadrants in $\mathcal{P}$. The set of rectangles in $\mathcal{R}$ containing a neighborhood of $\gamma$ in $q$ is finite up to the action of $\text{Stab}^+(\gamma)$, the subgroup of $\text{Stab}(\gamma)$ preserving all the stable and unstable separatrices of $\gamma$.
\end{lemm}
\begin{proof}
    Indeed, fix $R$ a rectangle in $\mathcal{R}$ containing a neighborhood of $\gamma$ in $q$. Take $R'$ a rectangle in $\mathcal{R}$ with the same property that is not in the $G$-orbit of $R$. By Lemma \ref{l.npredecessor}, $R'$ is a successor or predecessor of $R$ of some generation. Denote by $t$ the generator of $\text{Stab}^{+}(\gamma)$ that acts as a topological expansion on every unstable separatrix of $\gamma$. Notice that $\rho(t)(R)$ is a predecessor of some generation, say $n$ (with $n\in \mathbb{N}^*$),  of $R$ containing a neighborhood of $\gamma$ in $q$. By eventually replacing $R'$ by $\rho(t^k)(R')$ for some $k\in \mathbb{Z}$, we can assume that $R'$ is a predecessor of some generation of $R$ and a successor of some generation of $t(R)$. It follows that $R'$ is a predecessor of generation $m$ of $R$, where $m<n$ and since the predecessors of $R$ of generation less than $n$ form a finite set (see Lemma \ref{l.existenceofpredecessors}), we get the desired result. 
\end{proof}
Take $\gamma\in \Gamma$ a boundary periodic point of $\mathcal{R}$ and $L_0\in \mathcal{R}$ containing $\gamma$. Denote by $p\geq 2$ the total number of stable separatrices of $\gamma$ and by $Q_0,...,Q_{2p-1}$ the cyclically ordered set of quadrants of $\gamma$. By eventually reindexing the previous set of quadrants, we can assume without any loss of generality that $L_0$ contains a neighborhood of $\gamma$ in the quadrant $Q_0$ and does not intersect the interior of $Q_{2p-1}$. Notice that by Lemma \ref{l.finite}, there exist  finitely many such rectangles up to the action of $\text{Stab}^+(\gamma)$. Assume without any loss of generality that the quadrant $Q_1$ intersects $Q_0$ along an unstable separatrix $u_0$ of $\gamma$. Denote by $R^{Q_1}$ the set of rectangles $r\in \mathcal{R}$ containing a neighborhood of $\gamma$ in $Q_1$ and such that $L_0\cap u_0 \subseteq r\cap u_0$. Consider also $L_1\in R^{Q_1}$ such that for all $r\in R^{Q_1}$ $$L_1\cap \mathcal{F}^u(\gamma)\subseteq r\cap \mathcal{F}^u(\gamma)$$

\textbf{\underline{Claim} :} Such a rectangle $L_1$ exists, $L_1$ is uniquely defined by the above property, that we will call the \emph{minimality condition}, and finally either $L_1=L_0$ or $\inte{L_1}\cap \inte{L_0}=\emptyset$. 

\begin{proof}[Proof of the claim]Indeed, thanks to Lemma \ref{l.nocyclearoundperiod}, we have exactly one of the following three cases: 
\begin{enumerate}[leftmargin=0.7cm,itemsep=0.2cm]
     \item The rectangles of $R^{Q_1}$ contain  neighborhoods of $\gamma$ in the quadrants $Q_0$ and $Q_1$
    \item The rectangles of $R^{Q_1}$ contain neighborhoods of $\gamma$ in the quadrants $Q_1$ and $Q_2$ 
    \item The rectangles of $R^{Q_1}$ contain neighborhoods of $\gamma$ only in the quadrant $Q_1$
\end{enumerate}
Using this fact together with the Markovian intersection axiom, for any $r,r'\in R^{Q_1}$ we have that  $r\cap \mathcal{F}^u(\gamma)\subseteq r'\cap \mathcal{F}^u(\gamma)$ or $r'\cap \mathcal{F}^u(\gamma)\subseteq r\cap \mathcal{F}^u(\gamma)$. Moreover, even though the set $R^{Q_1}$ is not invariant by the action of $\text{Stab}^+(\gamma)$, it is contained by Lemma \ref{l.finite} in the union of a finite number of $\text{Stab}^+(\gamma)$-orbits of rectangles. By combining the previous facts, it follows that there exists at least one rectangle in $R^{Q_1}$ satisfying the minimality condition.

If there exist two rectangles $L_1,L'_1\in R^{Q_1}$ satisfying the minimality condition, then $L_1\cap \mathcal{F}^u(\gamma)= L'_1\cap \mathcal{F}^u(\gamma)$. In any of the above three cases, this implies by the Markovian intersection property that either $\partial^sL_1 \subseteq \partial^sL'_1$ or $\partial^sL'_1 \subseteq \partial^sL_1$. Equivalently, the existence of two distinct rectangles $L_1,L'_1\in R^{Q_1}$ satisfying the minimality condition, would imply that $L_1\cap L_1'= L_1$ or $L_1\cap L_1'= L'_1$. This is impossible, since by the Markovian intersection axiom, $L_1\cap L'_1$ is a non-trivial horizontal or vertical subrectangle of $L_1$ (resp. $L'_1$). We deduce the uniqueness of $L_1$.

Finally, if $\inte{L_1}\cap \inte{L_0}\neq \emptyset$, then by Lemma \ref{l.nocyclearoundperiod}, all the rectangles of $R^{Q_1}$ contain neighborhoods of $\gamma$ in the quadrants $Q_0$ and $Q_1$. This implies by Lemma \ref{l.nocyclearoundperiod} that $L_0\in R^{Q_1}$ and since $L_0$ clearly satisfies the minimality condition, by the uniqueness of $L_1$, we have that $L_1=L_0$.
\end{proof}
Next, denote by $s_0$ the stable separatrix of $\gamma$ separating the quadrants $Q_1$ and $Q_1$. Denote by $R^{Q_2}$ the set of rectangles $r\in \mathcal{R}$ containing a neighborhood of $\gamma$ in $Q_2$ and such that $L_1\cap s_0 \subseteq r\cap s_0$. Again, we define $L_2$ as the unique rectangle in $R^{Q_2}$ for which all $r\in R^{Q_2}$ satisfy $$L_2\cap \mathcal{F}^s(\gamma)\subseteq r\cap \mathcal{F}^s(\gamma)$$ We repeat the previous construction a finite number of times, in order to produce a finite set of rectangles $L_3,...,L_{2p-1}$ and the sets of rectangles $R^{Q_3},...,R^{Q_{2p-1}}$.  
\begin{figure}
    \centering
    \includegraphics[scale=0.4]{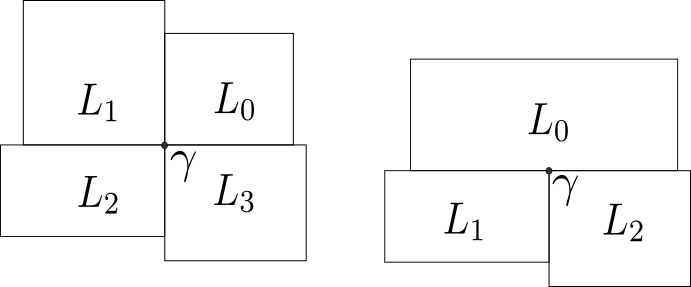}
    \caption{Some examples of pre-cycles starting from $L_0$ around a regular boundary periodic point}
    \label{f.precycles}
\end{figure}
 By Lemma \ref{l.nocyclearoundperiod} and our previous construction, we have that any two rectangles among $L_0,L_1,..., L_{2p-1}$ have either disjoint interiors or are equal (see Figure \ref{f.precycles}).

By eventually removing some rectangles from $L_0,L_1,..., L_{2p-1}$ so that every rectangle appears exactly once and by reindexing the remaining rectangles, we obtain a sequence of $k\in \llbracket 2,2p-1\rrbracket$ rectangles $L_0,...,L_{k-1}$. 
 \begin{defi}\label{d.precycle}
     We will call $L_0,...,L_{k-1}$ \emph{a pre-cycle around $\gamma$ starting from $L_0$}. 
 \end{defi}

\begin{rema}\label{r.precyclesfiniteuptoaction}
\begin{itemize}
    \item To every rectangle $L_0$ containing $\gamma$ and to every choice of cyclic order of the quadrants of $\gamma$, we associated by our previous construction a unique pre-cycle around $\gamma$ starting from $L_0$
    \item By changing the cyclic order of the quadrants of $\gamma$, we can define by the same procedure a second pre-cycle around $\gamma$ starting from $L_0$. It follows that for every rectangle containing $\gamma$, there exist two pre-cycles around $\gamma$ starting from this rectangle
    \item The image of a pre-cycle around $\gamma$ by an element of $\text{Stab}(\gamma)$ is a pre-cycle around $\gamma$. Therefore, by Lemma \ref{l.finite} there exist only a finite number of pre-cycles around $\gamma$ up to the action of $\text{Stab}(\gamma)$ or $\text{Stab}^+(\gamma)$. 
\end{itemize}

 \end{rema}

We will now associate to any pre-cycle around $\gamma$ a closed rectangle path of $\mathcal{P}$ that we will call a \emph{cycle}.

 Consider $L_0\in \mathcal{R}$ containing $\gamma$, $Q_0,Q_1,...,Q_{2p-1}$ the set of cyclically ordered quadrants of $\gamma$ (for some choice of cyclic order). By eventually reindexing the quadrants of $\gamma$, we can assume without any loss of generality that $L_0$ contains a neighborhood of $\gamma$ in $Q_0$ and does not intersect the interior of $Q_{2p-1}$. Assume also for the sake of simplicity that the pre-cycle starting from $L_0$ associated to our choice of cyclic order of the $Q_i$ is of the form  $(L_0,L_1,...,L_{2p-1})$, where the $L_i$ have pairwise disjoint interiors. Finally, assume without any loss of generality that $Q_0$ and $Q_1$ intersect along the unstable separatrix $u_0$ of $\gamma$. 

 Let $l_0=L_0\cap u_0$ and $r_0$ be the $l_0$-crossing successor of $L_0$ containing the endpoint of $l_0$ that is not $\gamma$ (this rectangle exists by Lemma \ref{l.crossingrectangleswithperiodicpoints}). By the  construction of $L_1$ (see our discussion prior to Definition \ref{d.precycle}), $r_0$ is a successor of $L_1$ of some generation (see Figure \ref{f.precycles}). Similarly, if $s_0$ is the separatrix of $\gamma$ lying in the intersection of $Q_1$ and $Q_2$, we define $l_1:=L_1\cap s_0$ and $r_1$ as the $l_1$-crossing predecessor of $L_1$ containing the endpoint of $l_1$ that is not $\gamma$. Again, $r_1$ is a predecessor of $L_2$ of some generation. We repeat this construction in order to construct the rectangles $r_2,...,r_{2p-2}$ which verify the following property: $r_i$ is a crossing successor (resp. predecessor) of $L_i$ when $i$ is even (resp. odd). 

Finally, if $s_p$ is the stable separatrix of $\gamma$ lying in the intersection of $Q_{2p-1}$ and $Q_0$ and $L_0\cap s_p\subset L_3\cap s_p$ we define $l_{2p-1}:=L_0\cap s_p$ and $r_{2p-1}$ as the $l_{2p-1}$-crossing predecessor of $L_0$ intersecting the endpoint of $l_{2p-1}$ that is not $\gamma$. If instead, $L_{2p-1}\cap s_p\subset L_0\cap s_p$ we define $l_{2p-1}:=L_{2p-1}\cap s_p$ and $r_{2p-1}$ as the $l_{2p-1}$-crossing predecessor of $L_{2p-1}$ intersecting the endpoint of $l_{2p-1}$ that is not $\gamma$. In any case, $r_{2p-1}$ is a predecessor of some generation of $L_0$ and also of $L_{2p-1}$.

We have thus constructed a generalized rectangle path in $\mathcal{P}$, namely  $$(L_0,r_0,L_1,r_1,...,L_{2p-1},r_{2p-1},L_0)$$

\begin{defi}\label{d.cycles}
The closed rectangle path associated to the generalized rectangle path $(L_0,r_0,L_1,r_1,...,L_{2p-1},r_{2p-1},L_0)$ (see Definition \ref{d.generalisedrectpaths}) is called the \emph{cycle associated to the pre-cycle $(L_0,L_1,...,L_{2p-1})$}. We define similarly the cycle associated to a pre-cycle with less than $2p-1$ rectangles. 

The set of cycles associated to all pre-cycles around $\gamma$ will be called the set of \emph{cycles} of $\gamma$ and will be denoted by $\text{Cycl}(\gamma)$.
\end{defi}
\begin{rema}\label{r.cyclesfiniteuptoaction}
If $c$ is a pre-cycle around $\gamma$ and $g\in \text{Stab}(\gamma)$. The image by $\rho(g)$ of the cycle  associated to $c$, is the cycle  associated to the pre-cycle $\rho(g)(c)$. Therefore, by Remark  \ref{r.precyclesfiniteuptoaction}, the set of cycles around any boundary periodic orbit $\gamma$  consists of a finite number of orbits by the action of $\text{Stab}(\gamma)$ or $\text{Stab}^+(\gamma)$.
\end{rema}

More generally, 
\begin{rema}\label{r.orbitofcycles}
For any $g\in G$ and any boundary periodic point $\gamma$, $\rho(g)(\text{Cycl}(\gamma))=\text{Cycl}(\rho(g)(\gamma))$. Therefore, the set of all cycles around any boundary periodic point in the $G$-orbit of $\gamma$, consists of a finite number of rectangle paths up to the action of $\rho$.
\end{rema}
  We will now prove that the set of rectangle paths $\underset{g\in G}{\cup}\text{Cycl}(\rho(g)(\gamma))$ can be described as a unique \emph{closed rectangle path} in any geometric type associated to $\mathcal{R}$. 

\begin{defi}
Consider a geometric type $\mathcal{G}=(n,(h_i)_{i \in \llbracket 1,n \rrbracket}, (v_i)_{i\in \llbracket 1,n \rrbracket}, \mathcal{H}, \mathcal{V},\phi, u)$, where $\mathcal{H}=\{H_i^j, i\in \llbracket 1,n \rrbracket, j\in \llbracket 1, h_i \rrbracket\}$ and $  \mathcal{V}=\{V_i^j, i\in \llbracket 1,n \rrbracket, j\in \llbracket 1, v_i \rrbracket\}$. We define a \emph{rectangle path} in $\mathcal{G}$ as a finite sequence of the form $i_0,s_0,i_1,s_1,i_2,...,s_m, i_{m+1}$ where \begin{itemize}     
    \item for all $k\in \llbracket 0, m+1\rrbracket$, $i_k\in \llbracket 1,n\rrbracket$ 
\item for all $k\in \llbracket 0, m \rrbracket$,  $s_k$ is of the form $V_{i_k}^{\bullet}$ or $H_{i_k}^{\bullet}$. 
     \item for all $k\in \llbracket 0, m\rrbracket$, if $s_k\in \mathcal{V}$ (resp. $s_k\in \mathcal{H}$)  then it is the image by $\phi$ (resp. $\phi^{-1}$) of an element of $\mathcal{H}$ (resp. $\mathcal{V}$) of the form $H_{i_{k+1}}^{\bullet}$ (resp. $V_{i_{k+1}}^{\bullet}$). 
\end{itemize}
Furthermore, the rectangle path $i_0,s_0,i_1,s_1,i_2,...,s_m, i_{m+1}$ will be called \emph{closed} when 
$i_0=i_{m+1}$
\end{defi}
Recall that, using its geometric interpretation, a geometric type $\mathcal{G}$ is a set of rectangles $R_1,...,R_n$, endowed each with a collection of horizontal and vertical subrectangles identified by the map $(\phi,u)$. In these terms, by identifying every $i_l\in \llbracket 1,n \rrbracket$ with $R_{i_l}$, a rectangle path $i_0,s_0,i_1,s_1,i_2,...,s_m, i_{m+1}$ in $\mathcal{G}$, can be thought as a sequence of rectangles and subrectangles of the form $R_0,s_0,...,R_m,s_m,R_{m+1}$, where for every $k\in \llbracket 0,m\rrbracket $, $s_k$ corresponds to a vertical (resp. or horizontal) subrectangle of $R_k$ and $\phi^{-1}(s_k)$ (resp. $\phi(s_k)$)  to a horizontal (resp. or vertical) subrectangle of $R_{k+1}$. 
\addtocontents{toc}{\protect\setcounter{tocdepth}{2}}
\subsubsection*{The relation between rectangle paths in $\mathcal{P}$ and rectangle paths in $G$}

Choose  $r_1,...,r_n$ representatives of every rectangle orbit in $\mathcal{R}$, an orientation of $\mathcal{P}$ and choose also an orientation of the stable and unstable foliations inside every $r_i$ such that the positive stable direction followed by the positive unstable direction inside every $r_i$ define an orientation compatible with our choice of orientation of $\mathcal{P}$. Let $\mathcal{G}=(R_1,...,R_n,(h_i)_{i \in \llbracket 1,n \rrbracket}, (v_i)_{i\in \llbracket 1,n \rrbracket}, \mathcal{H}, \mathcal{V},\phi, u)$ be the geometric type associated to $\mathcal{R}$ by Remark \ref{r.canonicalassociationgeometrictype} for the previous choices of representatives and orientations (we use here the geometric interpretation of $\mathcal{G}$).

Denote by $R_i$ the rectangle in $G$ associated to $r_i$. Recall that by Definition \ref{d.geometrictypemarkovfamily}, the vertical (resp. horizontal) subrectangles of the rectangles $R_1,...,R_n$ (in other words the elements of the sets $\mathcal{V}$ and $\mathcal{H}$) correspond to the predecessors (resp. successors) of the rectangles $r_1,...,r_n$. 

Consider now $l_0,l_1,...,l_k$ a rectangle path in $\mathcal{P}$. Using the fact that for every $i\in \llbracket 0, k \rrbracket$ there exists a unique $g_i\in G$ such that $\rho(g_i)(l_i)\in \{ r_1,...,r_n\}$, we get that 
\begin{itemize}
    \item to every $l_i$ we can associate a unique rectangle $R_{l_i}$ of the geometric type $\mathcal{G}$
    \item since for every $i\in \llbracket 1,k\rrbracket $ the rectangle $l_{i}$ is a successor or a predecessor of $l_{i-1}$, it follows that $\rho(g_{i-1})(l_{i})$ is also a successor or a predecessor of $\rho(g_{i-1})(l_{i-1})\in \{ r_1,...,r_n\}$. Therefore, the rectangle $l_i$ corresponds to a unique horizontal or vertical subrectangle $s_{l_i}$ of $R_{l_{i-1}}$.
\end{itemize}
 It is now easy to check that the  sequence $R_{l_0},s_{l_1},....,R_{l_{k-1}},s_{l_k}, R_{l_{k}}$ is a rectangle path in $\mathcal{G}$, which we will call \emph{the projection of the rectangle path $l_0,...,l_k$ in $\mathcal{G}$}. 

Conversely, consider  $R_{0},s_{0},R_1,s_1....,R_{k-1},s_{k-1}, R_{k}$ a rectangle path in $\mathcal{G}$. As before, the rectangle $R_{0}$ corresponds to a unique $\rho$-orbit of rectangles in $\mathcal{R}$. Let $l_0$ be an element of the previous orbit of rectangles and fix $g_0$ the unique element in $G$ that verifies $\rho(g_0)(l_0)\in \{r_1,...,r_n\}$. Notice that  $s_{0}$ corresponds to a unique successor or predecessor of $\rho(g_0)(l_0)$, whose $\rho$-orbit corresponds to the rectangle $R_1$ in $\mathcal{G}$. Thanks to the uniqueness of $g_0$, we get that $s_0$ corresponds to a unique successor or predecessor of $l_0$, say $l_1$. By repeating the previous procedure, we can construct a rectangle path $l_0,...,l_k$ in $\mathcal{P}$, that we will call \emph{the lift on $\mathcal{P}$ of $R_{0},s_{0},R_1,s_1....,R_{k-1},s_{k-1}, R_{k}$ starting from $l_0$}. 

The previous rectangle path  depends on our choice of $l_0$. For any $g\in G$, it is easy to check that the lift on $\mathcal{P}$ of $R_{0},s_{0},R_1,s_1....,R_{k-1},s_{k-1}, R_{k}$ starting from $\rho(g)(l_0)$, corresponds to the image by $\rho(g)$ of $l_0,...,l_k$.  We conclude that the rectangle path  $R_{0},s_{0},R_1,s_1....,R_{k-1},s_{k-1}, R_{k}$ can be lifted to a  $\rho$-orbit of rectangle paths in $\mathcal{P}$. 

\begin{rema}\label{r.associationgeomtypesinPgeomtypesinG}
    One can check, using our previous construction, that the set of lifts on $\mathcal{P}$ of the projection on $\mathcal{G}$ of any rectangle path $l_0,...,l_k$ in $\mathcal{P}$, coincides with the $\rho$-orbit of $l_0,...,l_k$. Conversely, it is also true that the projection on $\mathcal{G}$ of any lift on $\mathcal{P}$ of a rectangle path $R_{0},s_{0},R_1,s_1....,R_{k-1},s_{k-1}, R_{k}$ in $\mathcal{G}$ coincides with $R_{0},s_{0},R_1,s_1....,R_{k-1},s_{k-1}, R_{k}$. 
\end{rema}

\textit{On the proofs of Theorems D and E}

According to Theorem B', an equivalence class of geometric types of $\mathcal{R}$ describes up to conjugacy the lift of $\rho$ on the bifoliated plane up to surgeries along $\Gamma$, the set of boundary periodic points of $\mathcal{R}$. We would like at this point to define a combinatorial object describing a strong Markovian action up to conjugacy. 

\begin{defi}\label{d.geomtypecycles}
Let $\mathcal{G}$ be a geometric type and $\mathcal{A}$ a finite set of closed rectangle paths in $\mathcal{G}$. We will call $(\mathcal{G},\mathcal{A})$ a \emph{geometric type with cycles}.

Two geometric types with cycles $(\mathcal{G},\mathcal{A})$ and $(\mathcal{G}',\mathcal{A}')$ will be called \emph{equivalent} if there exists an equivalence $h$ between $G$ and $G'$ (see Definition \ref{d.equivalentgeomtypes}) that bijectively associates closed rectangle paths in $\mathcal{A}$ to closed rectangle paths in $\mathcal{A}'$. 
\end{defi}
\begin{rema}
Using Remark \ref{r.finitenumbergeometrictypes}, we can show that any equivalence class of geometric types with cycles contains only a finite number of geometric types with cycles. 
\end{rema}
Once again, choose a set of representatives $\{r_1,...,r_n\}$ for every rectangle orbit in $\mathcal{R}$, an orientation on $\mathcal{P}$ and an orientation of the stable and unstable foliations in every $r_i$ such that inside every $r_i$ the positive stable direction followed by the positive unstable direction define an orientation on $\mathcal{P}$ compatible with our initial choice of orientation. By Remark \ref{r.canonicalassociationgeometrictype}, together with this choice of representatives and orientations we can associate to $\mathcal{R}$ a unique geometric type $\mathcal{G}$.

Recall that to every boundary periodic point $\gamma\in \Gamma$ we can associate  the infinite set of cycles around $\gamma$, namely the set of closed rectangle paths $\text{Cycl}(\gamma)$ (see Definition \ref{d.cycles}). Denote by $\mathcal{A}$ the projection in $G$ of the infinite set of closed rectangle paths $\underset{\gamma\in \Gamma}{\cup}\text{Cycl}(\gamma)$.

\begin{defi}\label{d.geomtypecyclesmarkov}
  We will call $(G,\mathcal{A})$ a \emph{geometric type with cycles associated to $\mathcal{R}$} or more simply a \emph{geometric type with cycles of $\mathcal{R}$}. We will also call a geometric type with cycles of some strong Markovian action $\rho$, a \emph{geometric type with cycles associated to $\rho$} or a \emph{geometric type with cycles of $\rho$}.   
\end{defi}

A geometric type with cycles of $\mathcal{R}$ is a geometric type with cycles in the sense of Definition \ref{d.geomtypecycles}. Indeed, using our previous notations, since $\Gamma$ consists of a finite number of $\rho$-orbits of points in $\mathcal{P}$ (see Proposition \ref{p.boundaryperiodicpoints}), by Remark \ref{r.orbitofcycles} we have that $\underset{\gamma\in \Gamma}{\cup}\text{Cycl}(\gamma)$ consists of finitely many closed rectangle paths in $\mathcal{P}$ up to the action of $\rho$. Hence, thanks to Remark \ref{r.associationgeomtypesinPgeomtypesinG}, we get that  $\mathcal{A}$ is a finite set of closed rectangle paths in $G$.

Moreover, as in the case of geometric types: 
\begin{rema}
    By our previous construction, the following choices of representatives and orientations define a unique geometric type with cycles of $\mathcal{R}$:
    
    \begin{enumerate}
        \item a choice of representatives $\{r_1,...,r_n\}$ for every rectangle orbit in $\mathcal{R}$
        \item a choice of orientation on $\mathcal{P}$
        \item a choice of orientation of the stable and unstable foliations in every $r_i$ such that inside every $r_i$ the positive stable direction followed by the positive unstable direction define an orientation on $\mathcal{P}$ compatible with our previous choice of orientation
    \end{enumerate}
    \end{rema}

If $\gamma$ is a boundary periodic point in $\mathcal{R}$, the construction of $\text{Cycl}(\gamma)$ (see Definition \ref{d.cycles}) does not depend on our choice of representatives of rectangle orbits in $\mathcal{R}$, but also it does not depend on our choice of orientations of $\mathcal{P}$ or of the stable and unstable foliations inside the $r_i$. Hence, $\underset{\gamma\in \Gamma}{\cup}\text{Cycl}(\gamma)$ is an infinite set of closed rectangle paths canonically associated to $\mathcal{R}$. It follows that by changing our choices of representatives and orientations, we produce $(G',\mathcal{A}')$ a new geometric type with cycles associated to $\mathcal{R}$, such that the rectangle paths in $\mathcal{A}'$ and in $\mathcal{A}$ lift to the exact same set of rectangle paths in $\mathcal{P}$. By repeating our proof of Theorem \ref{t.associatemarkovfamiliestogeometrictype}, it is now not hard to show that 

\begin{theorem}[Theorem-Definition D']\label{t.geometrictypeswithcyclesclass}
    Take $\rho: G\rightarrow \text{Homeo}(\mathcal{P})$ an orientation preserving strong Markovian action on $\mathcal{P}$ leaving invariant a strong Markovian family $\mathcal{R}$. The set of geometric types with cycles associated to $\mathcal{R}$ is included in a unique class of equivalent geometric types with cycles.  

\end{theorem}

We would like now to prove the following generalization of Theorem E for strong Markovian actions, according to which any equivalence class of geometric types with cycles associated to some strong Markovian action describes the action up to conjugacy. More precisely,

\begin{theorem}[Theorem E']\label{t.theoreme'}
    Let $\rho_1: G_1\rightarrow \text{Homeo}(\mathcal{P}_1)$ $\rho_2: G_2\rightarrow \text{Homeo}(\mathcal{P}_2)$ be two orientation preserving strong Markovian actions on $\mathcal{P}_1$ and $\mathcal{P}_2$ preserving the strong Markovian families   $\mathcal{R}_1,\mathcal{R}_2$. Denote by $\mathcal{F}_1^{s,u}, \mathcal{F}_2^{s,u}$ the stable and unstable foliations of $\rho_1,\rho_2$. If the equivalence classes  geometric types with cycles associated to  $\mathcal{R}_1,\mathcal{R}_2$ are the same, then there exists $\beta$ an isomorphism from $G_1$ to $G_2 $ and $H$ a homeomorphism from $\mathcal{P}_1$ to $\mathcal{P}_2$ such that : 
\begin{itemize}
    \item the image by $H$ of any stable/unstable leaf in $\mathcal{F}_1^{s,u}$ is a  stable/unstable leaf in $\mathcal{F}_2^{s,u}$ 
    \item for every $g\in G_1$ and every $x\in \mathcal{P}_1$ we have $$H(\rho_1(g)(x))= \rho_2(\beta(g))(H(x))$$ 
\end{itemize}
\end{theorem}

We remark here that Theorem E is a direct consequence of Theorem E' and Barbot's theorem (see Theorem \ref{t.barbottheor}). 
\begin{proof}
Let $\Gamma_1,\Gamma_2$ be the boundary periodic points of $\mathcal{R}_1,\mathcal{R}_2$ in $\mathcal{P}_1,\mathcal{P}_2$ and $\mathcal{C}_1, \mathcal{C}_2$ the set of all cycles around the points in $\Gamma_1$ and $\Gamma_2$ respectively. Denote by $\clos{\mathcal{P}_1}$,  $\clos{\mathcal{P}_2}$ the bifoliated planes of $\rho_1$, $\rho_2$ up to surgeries along $\Gamma_1$, $\Gamma_2$ together with $\clos{\rho_1}: \clos{G_1}\rightarrow \text{Homeo}(\clos{\mathcal{P}_1})$ $\clos{\rho_2}: \clos{G_2}\rightarrow \text{Homeo}(\clos{\mathcal{P}_2})$ the lifts of $\rho_1,\rho_2$ on  $\clos{\mathcal{P}_1}$, $\clos{\mathcal{P}_2}$ and $\clos{\pi_1},\clos{\pi_2}$ the projections from $\clos{\mathcal{P}_1}$,  $\clos{\mathcal{P}_2}$ to $\mathcal{P}_1,\mathcal{P}_2$. Denote also by $\clos{\mathcal{F}_1^{s,u}},\clos{\mathcal{F}_2^{s,u}} $ the stable and unstable foliations of $\clos{\rho_1},\clos{\rho_2}$, by $\clos{\Gamma_1}, \clos{\Gamma_2}$ the lifts of $\Gamma_1,\Gamma_2$ on $\clos{\mathcal{P}_1}, \clos{\mathcal{P}_2}$ and by $\clos{\mathcal{R}_1}, \clos{\mathcal{R}_2}$ the lifts of  $\mathcal{R}_1, \mathcal{R}_2$ on $\clos{\mathcal{P}_1}, \clos{\mathcal{P}_2}$. Finally, let $\mathfrak{p_1}, \mathfrak{p_2}$ be the morphisms from $\clos{G_1},\clos{G_2}$ to $G_1,G_2$ and $\clos{\mathcal{C}_1},\clos{\mathcal{C}_2}$ be the sets of rectangle paths in $\clos{\mathcal{R}_1}, \clos{\mathcal{R}_2}$, whose projections by $\clos{\pi_1},\clos{\pi_2}$ belong in $\mathcal{C}_1, \mathcal{C}_2$.

Choose $\{\clos{r^1_1},...,\clos{r^1_n}\}$, $\{\clos{r^2_1},...,\clos{r^2_n}\}$ a set of representatives for every rectangle orbit in $\clos{\mathcal{R}_1}$ and $\clos{\mathcal{R}_2}$, and also orient the foliations $\clos{\mathcal{F}_1^{s,u}}$ and $\clos{\mathcal{F}_2^{s,u}} $ so that the unique geometric types associated to $\clos{\mathcal{R}_1}$ and $\clos{\mathcal{R}_2}$ for the previous choices of representatives and orientations coincide (see Proposition \ref{p.liftedmarkovfamiliessamegeomtype} and Corollary \ref{c.specialgeomtypes}). Denote by $\mathcal{G}$ the previous geometric type and recall that $\mathcal{G}$ belongs to the common equivalence class of geometric types associated to $\mathcal{R}_1, \mathcal{R}_2$.  

Given a rectangle path in $\clos{\mathcal{R}_1}$, by projecting it by $\clos{\pi_1}$ to a rectangle path in $\mathcal{R}_1$ and then to a rectangle path in $\mathcal{G}$, we define similarly to the case of $\mathcal{R}_1$, the \emph{projection in $ \mathcal{G}$ of a rectangle path in $\clos{\mathcal{R}_1}$}. We similarly define the projection in $ \mathcal{G}$ of a rectangle path in $\clos{\mathcal{R}_2}$. Conversely, given a rectangle path in $\mathcal{G}$, by considering the set of all its lifts in $\mathcal{R}_1$ and then the set of rectangle paths in $\clos{\mathcal{R}_1}$ that project to the previous set of rectangle paths in $\clos{\mathcal{R}_1}$, we define the \emph{set of lifts in $\clos{\mathcal{R}_1}$ of a rectangle path in $\mathcal{G}$}. We similarly define the \emph{set of lifts in $\clos{\mathcal{R}_2}$ of a rectangle path in $\mathcal{G}$}.

By our hypothesis, if $\mathcal{C}_1^{\mathcal{G}}, \mathcal{C}_2^{\mathcal{G}}$ denote the projections on $\mathcal{G}$ of $\mathcal{C}_1, \mathcal{C}_2$ (or equivalently of $\clos{\mathcal{C}_1},\clos{\mathcal{C}_2}$), then we know that the geometric types with cycles $(\mathcal{G},\mathcal{C}_1^{\mathcal{G}})$ and $(\mathcal{G},\mathcal{C}_2^{\mathcal{G}})$ are equivalent. By eventually using the previous equivalence, we can assume without any loss of generality that $(\mathcal{G},\mathcal{C}_1^{\mathcal{G}})= (\mathcal{G},\mathcal{C}_2^{\mathcal{G}})$. For the sake of simplicity, we will denote $\mathcal{A}:=\mathcal{C}_1^{\mathcal{G}}=\mathcal{C}_2^{\mathcal{G}}$. Consider a rectangle in $\{\clos{r^1_1},...,\clos{r^1_n}\}$ and another rectangle in $\{\clos{r^2_1},...,\clos{r^2_n}\}$ such that the two rectangles are of the same type. We will fix those rectangles as the origin rectangles in $\clos{\mathcal{R}_1}$ and $\clos{\mathcal{R}_2}$.

Thanks to Theorem B', there exist an isomorphism $\alpha: \clos{G_1}\rightarrow \clos{G_2}$ and a homeomorphism $h:\clos{\mathcal{P}_1}\rightarrow \clos{\mathcal{P}_2}$ such that: 

\begin{itemize}
    \item the image by $h$ of any oriented stable/unstable leaf in $\clos{\mathcal{F}_1^{s,u}}$ is an oriented stable/unstable leaf in $\clos{\mathcal{F}_2^{s,u}}$ 
    \item for every $g\in \clos{G_1}$ and every $x\in \clos{\mathcal{P}_1}$ we have $$h(\clos{\rho_1}(g)(x))= \clos{\rho_2}(\alpha(g))(h(x))$$ 
\end{itemize}

Let us now show that that under the additional hypothesis on the cycles of $\mathcal{R}_1$ and $\mathcal{R}_2$ the map $h$ projects to a conjugacy between $\rho_1$ and $\rho_2$. Recall first that by construction of $h$ (see Proposition \ref{p.existencehomeo}), if $path_1$ is a centered rectangle path in $\clos{\mathcal{R}_1}$, then its image by $h$ is the centered rectangle path in  $\clos{\mathcal{R}_2}$ associated to $path_1$. In other words, the image of $path_1$ by $h$ defines a rectangle path in $\clos{\mathcal{R}_2}$, say $path_2$, such that the projections of $path_1$ and $path_2$ on $\mathcal{G}$ coincide. It follows that $h(\clos{\mathcal{C}_1})$ corresponds to a rectangle path in $\clos{\mathcal{C}_2}$ and vice versa. 

The positive stable direction followed by the positive unstable direction define in $\clos{\mathcal{P}_1}$ and $\clos{\mathcal{P}_2}$ two orientations. Endow $\clos{\mathcal{P}_1}$ and $\clos{\mathcal{P}_2}$ with the previous orientations.  

Let $\clos{\gamma}\in\clos{\Gamma_1}$ and $...,s_{-1},s_0,s_1,s_2,...$ be the stable leaves in $\clos{\mathcal{F}_1^s}$ containing $\clos{\gamma}$ ordered counter-clockwisely. In other words, $s_{i+1}$ is not separated from $s_i$ (see Proposition \ref{p.singularitiesoffoliations}) and for any point $x\in s_i$, the direction in $s_i$ pointing away from $\clos{\gamma}$ followed by the direction pointing inside the connected component of $\clos{\mathcal{P}_1}-s_i$ containing $s_{i+1}$ (see Remark \ref{r.propertiesfolipbar}) defines an orientation compatible with our choice of orientation of $\clos{\mathcal{P}_1}$. 

In Proposition \ref{p.mpformbasis}, we showed that using the previous orientations, one can define a canonical base of $\text{Stab}(\clos{\gamma})$ as follows\footnote{We remark here that the section in which Proposition \ref{p.mpformbasis} was proven dealt only with strong Markovian actions arising from pseudo-Anosov flows. However, Proposition \ref{p.mpformbasis} does not use at all the pseudo-Anosov hypothesis}. Consider $\gamma:=\clos{\pi}(\clos{\gamma})$, $n$ the smallest positive integer such that $\clos{\pi}(s_n)=\clos{\pi}(s_0)$ ($n$ corresponds to the number of stable prongs of $\gamma$) and define  

\begin{itemize}
    \item $m^1_{\clos{\gamma}}$ as the unique element in $\text{ker}(\mathfrak{p_1})\leq \clos{G_1}$ taking $s_0$ to $s_n$
    \item $p_{\clos{\gamma}}^1$ as the unique element in $\clos{G_1}$ taking $s_0$ to a stable leaf in $\{s_1,...,s_{n-1}, s_n\}$ and such that $\mathfrak{p}(p_{\clos{\gamma}}^1)=g_{\clos{\gamma}}^1$, where $g_{\clos{\gamma}}^1$ is the generator of $\text{Stab}(\gamma)$ in $G_1$ acting as an expansion on $\mathcal{F}^s(\gamma)$
\end{itemize}

We similarly define $m^2_{h(\clos{\gamma})}$ and $p^2_{h(\clos{\gamma})}$ for $h(\clos{\gamma})\in \clos{\Gamma_2}$, using our choice of orientation of $\clos{\mathcal{P}_2}$. Let us now prove that $$\alpha(m^1_{\clos{\gamma}})=m^2_{h(\clos{\gamma})}$$

Indeed, take any cycle $l_0,...,l_k=l_0$ in $\mathcal{C}_1$ around $\gamma$. By definition of a cycle, the previous closed rectangle path goes once around $\gamma$ in $\mathcal{P}_1$. It is therefore not hard to see that if $\clos{l_0},\clos{l_1},...,\clos{l_k}$ is a rectangle path in $\clos{\mathcal{C}_1}$ that lifts $l_0,...,l_k$, then $\clos{l_k}= \clos{\rho_1}(m^1_{\clos{\gamma}})(\clos{l_0})$ or $\clos{l_k}= \clos{\rho_1}((m^1_{\clos{\gamma}})^{-1})(\clos{l_0})$. Naturally, the same applies for any element of $\clos{\mathcal{C}_2}$. Assume without any loss of generality that $\clos{l_k}= \clos{\rho_1}(m^1_{\clos{\gamma}})(\clos{l_0})$. Since $h$ sends the rectangle paths in $\clos{\mathcal{C}_1}$ to rectangle paths in $\clos{\mathcal{C}_2}$ and is equivariant with respect to the group actions, we get that $$ \alpha(m^1_{\clos{\gamma}})=(m^2_{h(\clos{\gamma})})^{\pm}$$
Thanks to our proof of Theorem B', $h$ sends positively oriented stable/unstable segments in $\clos{\mathcal{F}_1^{s,u}}$ to positively oriented stable/unstable segments in $\clos{\mathcal{F}_2^{s,u}}$. It follows that if $\clos{\rho_1}(m^1_{\clos{\gamma}})$ acts as a positive translation on the counter-clockwise ordered set of stable separatrices on $\clos{\gamma}$, then $\clos{\rho_2}(\alpha(m^1_{\clos{\gamma}}))$ must also act as a positive translation on the counter-clockwise ordered set of stable separatrices on $h(\clos{\gamma})$. This implies that  $$\alpha(m^1_{\clos{\gamma}})=m^2_{h(\clos{\gamma})}$$ 

Next, let us prove that $$\alpha(p^1_{\clos{\gamma}})=p^2_{\clos{\gamma}}$$ Once again since $h$ is equivariant with respect to the group actions, $\alpha$ must send a basis of $\text{Stab}(\clos{\gamma})$ to a basis of $\text{Stab}(h(\clos{\gamma}))$. It follows that there exist $a\in \mathbb{Z}$ such that $$\alpha(p^1_{\clos{\gamma}})=(m^2_{h(\clos{\gamma})})^a \cdot (p^2_{\clos{\gamma}})^{\pm}$$ 

Using once more the equivariance of $h$, we have that $\alpha(p^1_{\clos{\gamma}})$ must act as an expansion on the set of stable leaves of $h(\clos{\gamma})$. Using the fact that $m^2_{\clos{\gamma}} \in \text{ker}(\mathfrak{p_2})$ and Proposition \ref{p.actionequivariance}, we get that  

$$\alpha(p^1_{\clos{\gamma}})=(m^2_{h(\clos{\gamma})})^a \cdot p^2_{\clos{\gamma}}$$ 

Finally, recall that $\rho_1(p^1_{\clos{\gamma}})$ (resp. $\rho_1(m^1_{\clos{\gamma}})$) acts as a translation on the set of stable leaves of $\clos{\gamma}$ by $+1,+2,...$ or $+n$ (resp. $+n$), where $n$ is the total number of stable prongs of $\gamma$. By the equivariance of $h$, a similar statement should be true for  $\rho_2(\alpha(p^1_{\clos{\gamma}}))$. This forces $a=0$ and implies that 

$$\alpha(p^1_{\clos{\gamma}})=p^2_{\clos{\gamma}}$$

To conclude, thanks to Proposition \ref{p.kernelprojectiongroup}, $\text{ker}(\mathfrak{p_1})$ is normally generated in $\clos{G_1}$ by the set $\{m^1_{\clos{\gamma}}|\clos{\gamma}\in \clos{\Gamma}\}$. Same for $\text{ker}(\mathfrak{p_2})$. Our previous arguments imply that $\alpha(\text{ker}(\mathfrak{p_1}))=\text{ker}(\mathfrak{p_2})$. By Propositions  \ref{p.kernelprojectiongroup}, \ref{p.actionequivariance}, this implies that the isomorphism $\alpha$ projects to an isomorpshism $\beta$ between $G_1$ and $G_2$ and that $h$ projects to a homeomorphism $H$ between $\mathcal{P}_1$ and $\mathcal{P}_2$ such that 

\begin{itemize}
    \item the image by $H$ of any stable/unstable leaf in $\mathcal{F}_1^{s,u}$ is a  stable/unstable leaf in $\mathcal{F}_2^{s,u}$ 
    \item for every $g\in G_1$ and every $x\in \mathcal{P}_1$ we have $$H(\rho_1(g)(x))= \rho_2(\beta(g))(H(x))$$ 
\end{itemize}

We thus get the desired result. 
\end{proof}


\begin{thebibliography}{11}
\bibitem[A]{Anosov} V.Dmitry Anosov. \emph{Geodesic flows on closed Riemannian manifolds of negative curvature}, Trudy Mat. Inst. Steklov., Vol. 90, pp. 3-210 (1967); Proc. Steklov Inst. Math., Vol. 90, pp. 1-235 (1967)
\bibitem[AgTs]{Agol} Ian Agol and Chi Cheuk Tsang, \emph{Dynamics of veering triangulations: infinitesimal components of their flow graphs and applications}, Algebraic $\&$ Geometric Topology, Vol. 24, No. 6, pp. 3401–3453 (2024) 
\bibitem[Ba]{Barbotthese}
Thierry Barbot. \textit{G\'eom\'etrie transverse des flots d’Anosov}, Thesis, Ecole Normale Sup\'erieure de Lyon (1992)
\bibitem[Ba1]{Ba1} Thierry Barbot. \emph{Caract\'erisation des flots d'Anosov en dimension 3 par leurs feuilletages faibles},
Ergodic Theory and Dynamical Systems, Vol. 15, No. 2, pp. 247-270 (1995)
\bibitem[Ba2]{Barbotbola} Thierry Barbot. \emph{Generalizations of the Bonatti-Langevin example of Anosov flow and their
classification up to topological equivalence}, Communications in 
Analysis and Geometry, Vol. 6, No. 4,  pp.749-798 (1998)
\bibitem[Bart]{Introanosov} Thomas Barthelmé. \emph{Anosov flows in dimension 3, A preliminary version}, Unpublished manuscript available online \hyperlink{https://drive.google.com/file/d/1JMfXlAs6i6f8YXAERFWHz-HVhd0mrvPv/view}{here}.
\bibitem[BaBaMa]{nontransitiveanosovlike} Thomas Barthelmé, Christian Bonatti, Kathryn Mann \emph{Non-transitive pseudo-Anosov flows}, arXiv: 2411.03586
\bibitem[BaFe]{BarbotFenley} Thierry Barbot, Sergio Fenley. \emph{Pseudo-Anosov flows in toroidal manifolds}, Geometry \& Topology, Vol. 17, No. 4, pp. 1877–1954 (2013)
\bibitem[BaFe1]{flowsingraphmanifolds} Thierry Barbot, Sergio Fenley. \emph{Classification and rigidity of totally periodic pseudo-Anosov flows in graph manifolds}, Ergodic Theory and Dynamical Systems, Vol. 35, No. 6, pp. 1681-1722 (2015)
\bibitem[BaFePo]{collapsed} Thomas Barthelmé, Sergio Fenley, Rafael Potrie, \emph{Collapsed Anosov flows and self orbit equivalences}, arXiv:2008.06547
\bibitem[BaFrMa]{circleatinfinity} Thomas Barthelmé, Steven Frankel, Kathtyn Mann. \emph{Orbit equivalences for Pseudo-Anosov flows}, Inventiones mathematicae, Vol. 240, pp. 1119–1192 (2025)
\bibitem[BaMa]{planeapproach} Thomas Barthelmé, Kathryn Mann. \emph{Pseudo-Anosov flows: a plane approach}, arXiv:2509.15375
\bibitem[Be]{Beguin} Fran\c{c}ois B\'eguin. \emph{Smale diffeomorphisms of surfaces: a classification algorithm}, Discrete and Continuous Dynamical  Systems, Vol. 11, No. 2$\&$3, pp. 261--310 (2004)
\bibitem[BeBoYu]{plugs} Fran\c{c}ois B\'eguin, Christian Bonatti, Bin Yu. \emph{Building Anosov flows on 3-manifolds}, Geometry and Topology, Mathematical Sciences Publishers, Vol. 21, No. 3, pp. 1837 - 1930 (2017)
\bibitem[BeBoYu1]{plugsdecomp} Fran\c{c}ois B\'eguin, Christian Bonatti, Bin Yu. \emph{A spectral-like decomposition for transitive Anosov flows in dimension three},  Mathematische Zeitschrift, Vol. 282, pp. 889–912 (2016)
\bibitem[BeYu]{uniqueglueing} Fran\c{c}ois B\'eguin, Bin Yu. \emph{A uniqueness theorem for transitive Anosov flows obtained by gluing hyperbolic plugs}, arXiv:1905.08989
\bibitem[Bo]{circleinfinitybonatti} Christian Bonatti. \emph{Action on the circle at infinity of foliations of $\mathbb{R}^2$}, 	arXiv:2301.04530
\bibitem[BoIa]{BonattiIakovoglou}
Christian Bonatti, Ioannis Iakovoglou. \textit{Anosov flows on 3-manifolds: the surgeries and the foliations},  Ergodic Theory and Dynamical Systems, Vol. 43, No. 4, pp. 1129-1188 (2023)
\bibitem[BoIa1]{realisablegeomtypes}Christian Bonatti, Ioannis Iakovoglou. \emph{Geometric types associated to Anosov flows in dimension 3} (in progress) 
\bibitem[BoLa]{Asterisque}
Christian Bonatti, R\'emi Langevin, avec la collaboration de Emmanuelle Jeandenans. \textit{Diff\'eomorphismes de Smale des surfaces}, Ast\'erisque, Tome 250 (1998)
\bibitem[BoLa1]{BonattiLangevin} Christian Bonatti, R\'emi Langevin. \emph{Un exemple de flot d'Anosov transitif transverse à un tore et non conjugué à une suspension}, Ergodic Theory and Dynamical Systems, Vol. 14, No. 4, pp. 633-643 (1994)
\bibitem[BoMa]{Bowdenmann} Jonathan Bowden, Kathryn Mann. \emph{$C^0$ stability of Anosov boundary action and inequivalent Anosov flows}, Annales Scientifiques de l'Ecole Normale Superieure, Série 4, Tome 55, No. 4, pp.  1003-1046 (2022)
\bibitem[Br]{Brown} Morton Brown. \emph{Locally Flat Embeddings of Topological Manifolds}, Annals of Mathematics, Vol. 75, No. 2, pp. 331-341 (1962)
\bibitem[Bru]{Brunella} Marco Brunella. \emph{Expansive flows on three-manifolds}, Thesis, International School for Advanced Studies (1992)
\bibitem[CaCo]{Candel} Alberto Candel, Lawrence Conlon. \emph{Foliations I.}, Graduate Studies in Mathematics, Vol. 23, American Mathematical Society (2000) 
\bibitem[Ch]{Chevallay} Claude Chevalley. \emph{Theory of Lie groups}, Princeton Univ. Press (1946)
\bibitem[ClPi]{ClPi} Adam Clay, Tali Pinsky. 
\emph{Graph manifolds that admit arbitrarily many Anosov flows}, arXiv:2006.09101 
\bibitem[De]{lefthand} Pierre Dehornoy. \textit{Which geodesic flows are left-handed ?, Groups, Geometry, and Dynamics}, Vol. 11, pp. 1347-1376 (2017)
\bibitem[DeSh]{almostequiv} Pierre Dehornoy, Mario Shannon. \emph{Almost equivalence of algebraic Anosov flows}, arXiv:1910.08457 
\bibitem[FaMa]{Primer} Benson Farb, Dan Margalit. \emph{A primer on mapping class
 groups}, Princeton University Press, 2012.
\bibitem[Fe]{Fe} Sergio Fenley. \emph{The structure of branching in Anosov flows of 3-manifolds}, Comment. Math. Helv., Vol. 73, No. 2, pp. 259-297 (1998)
\bibitem[Fe1]{Fe1} Sergio Fenley. \emph{Anosov flows in 3-manifolds}, Ann. of Math., Vol. 139, No. 1, pp. 79-115  (1994) 
\bibitem[Fe2]{Fe2} Sergio Fenley. \textit{Foliations with good geometry}, J. Amer. Math. Soc. Vol. 12, No.3, pp.619–676 (1999)
\bibitem[Fe3]{Fe3} Sergio Fenley. \textit{$\mathbb{R}$-covered foliations and transverse pseudo-Anosov flows in atoroidal pieces}. Comment. Math. Helv. Vol.98, No. 1, pp. 1-39 (2023)
\bibitem[Fe4]{Fe4} Sergio Fenley. \emph{Ideal boundaries of pseudo-Anosov flows and uniform convergence groups with connections and applications to large scale geometry}, Geom. Topol., Vol. 16, No. 1, pp. 1-110 (2012)
\bibitem[FeMo]{Fenleymosher} Sergio Fenley, Lee Mosher. \emph{Quasigeodesic flows in hyperbolic 3-manifolds}, Topology, Vol. 40, No. 3, pp. 503-537 (2001)
\bibitem[FePo]{partial}
Sergio R. Fenley, Rafael Potrie, \textit{Partial hyperbolicity and pseudo-Anosov dynamics}, arXiv:2102.02156
\bibitem[Fl]{noexpansivesurface} Laurence W. Flinn, \emph{Expansive flows}, Phd thesis, Warwick University (1972)	https://wrap.warwick.ac.uk/71294/
\bibitem[Fri]{Fried} David Fried. \emph{Transitive Anosov flows and pseudo-Anosov maps.} Topology, Vol. 22, No. 3, pp. 299-303 (1983)
\bibitem[GK]{GerberKatok} Marlies Gerber, Anatole  Katok. \emph{Smooth models of Thurston's pseudo-Anosov maps}, Annales scientifiques de l'École Normale Supérieure, Série 4, Tome 15, No. 1, pp. 173-204 (1982)
\bibitem[GO]{Gabai} David Gabai and Ulrich Oertel, \textit{Essential Laminations in 3-Manifolds}, Annals of Mathematics, Vol. 130, No. 1, pp. 41-73 (1989)
\bibitem[Hae]{Haefliger} André Haefliger, \textit{Variétés feuilletées}, Annali della Scuola Normale Superiore di Pisa, Classe di Scienze 3e
série, Tome 16, No 4, pp. 367-397 (1962)
\bibitem[Hi]{Hiraide} Koichi Hiraide, \textit{Expansive homeomorphisms of compact surfaces are pseudo-Anosov},
Osaka J. Math., Vol. 27, No. 1, pp.117–162 (1990)
\bibitem[HH]{Hirsh} G. Hector and U. Hirsch, \textit{Introduction to the geometry of foliations, Part B}, Vieweg, Braunschweig (1983)
\bibitem[HiPS]{stablemanifold} Morris Hirsch, Charles Chapman Pugh, Michael Shub. \emph{Invariant manifolds}, Lecture Notes in Mathematics, Vol. 583 (1977), Springer-Verlag
\bibitem[Ia]{preprint}
Ioannis Iakovoglou, \textit{A new combinatorial invariant characterizing Anosov flows on 3-manifolds}, arXiv:2212.13177
\bibitem[Ia1]{fundgroup} Ioannis Iakovoglou, \textit{On the fundamental groups of manifolds supporting Anosov flows}, (in preparation)
\bibitem[Ia2]{totperiodicioannis}
Ioannis Iakovoglou, \textit{On the complete classification of a family of Anosov flows by geometric types} (in progress) 
\bibitem[Ia3]{realisableactions} Ioannis Iakovoglou, \textit{Group actions arising from Anosov flows}, (in progress) 
\bibitem[Ia4]{thesisioannis} Ioannis Iakovoglou. \emph{Classifying Anosov flow in dimension 3 by geometric types}, Thesis, Université de Bourgogne  (2023)
\bibitem[Ia5]{markovpseudoanosov} Ioannis Iakovoglou. \textit{Markov parititions for non-transitive pseudo-Anosov flows}, Comptes Rendus Mathématique (accepted 2024) 
\bibitem[Ia6]{surgeriesnontransitive} Ioannis Iakovoglou. \textit{Surgeries on non-transitive expansive flows}, arXiv:2410.21470 (2024)
\bibitem[InMa]{Inaba} Takashi Inaba and Shigenori Matsumoto. \textit{Nonsingular expansive flows on 3-manifolds and foliations with circle prong singularities}, Japan. J. Math.Vol. 16, No. 2 (1990)
\bibitem[KeSe]{Keynes} H.B. Keynes and M. Sears, \emph{Real-expansive flows and topological dimension}, Ergodic Theory and Dynamical Systems, No. 1, pp. 179-195 (1981)
\bibitem[Le]{Lewowicz} Jorge Lewowicz, \textit{Expansive homeomorphisms of surfaces}, Bol. Soc. Brasil. Mat.
(N.S.), Vol. 20, No. 1, pp.113–133 (1989)
\bibitem[Lu]{Aspherical} Wolfgang Lueck. \emph{Survey on aspherical manifolds}, 	arXiv:0902.2480
\bibitem[Mo]{Mosher1} Lee Mosher, \textit{Dynamical systems and the homology norm of a 3–manifold. I: Efﬁcient intersection of surfaces and ﬂows}, Duke Math. J. Vol. 65, pp.449-500 (1992)
\bibitem[Mo1]{Mosher2} Lee Mosher, \textit{Dynamical systems and the homology norm of a 3–manifold. II}, Invent. Math. Vol. 107, pp.243–281 (1992)
\bibitem[Mo2]{Mosher3} Lee Mosher, \textit{Laminations and flows transverse to finite depth foliations}, manuscript available online \href{https://homepages.warwick.ac.uk/~masgar/Maths/1996lamina-tions\textunderscore and\textunderscore flows\textunderscore transverse\textunderscore to\textunderscore finite\textunderscore depth\textunderscore foliations\textunderscore OCRed.pdf}{here}
\bibitem[No]{Novikov} S. P. Novikov, \emph{The topology of foliations}, Trudy Moskovskogo Matematicheskogo Obshchestva, Vol. 14, pp. 248–278 (1965) 
\bibitem[Oka]{Oka}
Masatoshi Oka, \textit{Singular foliations on cross-sections of expansive flows on 3-manifolds}, Osaka J. Math., Vol 26, pp. 863-883 (1990)
\bibitem[Pa]{Palmeira} Carlos Frederico Borges Palmeira. \emph{Open manifolds foliated by planes}, Annals of Mathematics, Vol. 107, No. 1, pp. 109-131 (1978)
\bibitem[Pat]{Paternain} Miguel Paternain. \emph{Expansive flows and the fundamental group}, Bol. Soc. Bras. Mat., Vol. 24, No. 2, pp. 179-199 (1993)
\bibitem[Pau]{Neige} Neige Paulet. \emph{Flots d’Anosov en dimension trois construits par recollements de blocs}, Thesis, Université Paris-Nord (2023)
\bibitem[Pl]{Plante} J.F.Plante. \emph{Anosov flows}, American Journal of Mathematics,  Vol. 94, No. 3, pp. 729-754 (1972)
\bibitem[Ra]{Ratner}
M. Ratner. \textit{Markov partitions for C-flows on 3-dimensional manifolds}, Mat. Zametki, Vol. 6, No. 6, pp. 693-704 (1969) 
\bibitem[Sc]{Scott}
P. Scott. \emph{The geometries of 3-manifolds}, Bull. London Math. Soc., Vol.15, No.5, pp. 401–487 (1983)
\bibitem[Si1]{Sin1}
Yakov Grigorevich Sinai. 
\textit{Markov partitions and C-diffeomorphisms}, Functional Anal. Appl., Vol. 2, pp. 64-89 (1968)
\bibitem[Si2]{Sin2}
Yakov Grigorevich Sinai. \textit{Construction of Markov partitions}, Functional Anal. Appl., Vol. 2, pp. 70-80 (1968)
\bibitem[Sh]{Mariothese}
Mario Shannon. \emph{Dehn surgeries and smooth structures on 3-dimensional transitive Anosov flows}, Thesis, Université de Bourgogne (2020)
\bibitem[So]{Solodov} V.V.Solodov, \emph{Components of topological foliations}, Mathematics of the USSR-Sbornik, Volume 47, Issue 2, pp. 329–343 (1984)
\bibitem[Whi]{Whitney}
Hassler Whitney, \textit{Regular families of curves}, Annals of Mathematics, No.34 (1933)
\bibitem[YaYu]{plugsclassification}Jiagang Yang, Bin Yu. \emph{Classifying expanding attractors on figure eight knot complement space and non-transitive Anosov flows on Franks-Williams manifold}, arXiv:2004.08921


\end{thebibliography}
\end{document}